\documentclass[%draft,
reqno]{amsbook}
          \usepackage{amssymb}
	  \usepackage{amsmath}
      \usepackage{stmaryrd} % for square double bracket
          \usepackage{amsfonts}
          \usepackage[english]{babel}
          \usepackage[utf8]{inputenc}
\usepackage{bbm}
\usepackage{cancel,comment}

\usepackage[nopostdot,sort=def]{glossaries}
\usepackage{glossary-superragged}
\setglossarystyle{superragged3col}
\makeglossaries 

\usepackage[unicode,colorlinks,plainpages=false,hyperindex=true,bookmarksnumbered=true,bookmarksopen=false,pdfpagelabels]{hyperref}
 \hypersetup{urlcolor=cyan,linkcolor=blue,citecolor=red,colorlinks=true,
 pdftitle={Double Poisson (vertex) algebra cohomology},
 pdfsubject={Double Poisson vertex algebra cohomology (2025)},
 pdfauthor={Maxime Fairon and Daniele Valeri}}
\usepackage{color}

\usepackage{enumerate}

\usepackage{tikz} 
\usetikzlibrary{arrows}  
\usepackage{tikz-cd}

\theoremstyle{plain}
\newtheorem{theorem}{Theorem}[chapter]
\newtheorem{lemma}[theorem]{Lemma}
\newtheorem{proposition}[theorem]{Proposition}
\newtheorem{corollary}[theorem]{Corollary}

\theoremstyle{definition}
\newtheorem{definition}[theorem]{Definition}
\newtheorem{example}[theorem]{Example}

\theoremstyle{remark}
\newtheorem{remark}[theorem]{Remark}
\newtheorem{convention}[theorem]{Convention}

\numberwithin{equation}{chapter}
\numberwithin{section}{chapter}
\numberwithin{subsection}{section}

\allowdisplaybreaks[1] %For breaking long equations

\newcommand{\N}{\ensuremath{\mathbb{N}}}

\newcommand{\Z}{\ensuremath{\mathbb{Z}}}
\newcommand{\kk}{\ensuremath{\Bbbk}}

\newcommand{\mc}[1]{\mathcal{#1}}
\newcommand{\mf}[1]{\mathfrak{#1}}
\newcommand{\mb}[1]{\mathbb{#1}}

\newcommand{\tint}{{\textstyle\int}}

\newcommand{\tr}{\operatorname{tr}}

\newcommand{\Hom}{\operatorname{Hom}}
\newcommand{\End}{\operatorname{End}}

\newcommand{\Der}{\operatorname{Der}}
\newcommand{\Mat}{\operatorname{Mat}}
\newcommand{\Id}{\operatorname{Id}}
\newcommand{\Gl}{\operatorname{GL}}
\newcommand{\Rep}{\operatorname{Rep}}
\newcommand{\Spec}{\operatorname{Spec}}
\newcommand{\res}{\operatorname{res}}
\DeclareMathOperator{\mult}{m}
\DeclareMathOperator{\proj}{P}

\DeclareMathOperator{\Vect}{Vect}
\DeclareMathOperator{\im}{im}
\DeclareMathOperator{\inv}{inv}
\DeclareMathOperator{\sgn}{sgn}

\DeclareMathOperator{\Span}{span}
\DeclareMathOperator{\RCDer}{RCDer}
\DeclareMathOperator{\cent}{Z}
\DeclareMathOperator{\Cas}{Cas}
\DeclareMathOperator{\inn}{Inn}

\DeclareMathOperator{\ev}{ev}
\DeclareMathOperator{\Bil}{Bil}
\newcommand{\cA}{\mathcal{A}}
\newcommand{\cG}{\mathcal{G}}
\newcommand{\XX}{\mathtt{X}}

\newcommand{\gl}{\ensuremath{\mathfrak{gl}}}
\newcommand{\g}{\ensuremath{\mathfrak{g}}}

\newcommand{\del}{\ensuremath{\partial}}
\newcommand{\dd}{\operatorname{d}}
\newcommand{\wdd}{\widehat{\operatorname{d}}}
\newcommand{\BRA}{\operatorname{BR}}
\newcommand{\wBRA}{\widehat{\operatorname{BR}}}
\DeclareMathOperator{\SN}{SN}
\newcommand{\coH}{\operatorname{H}}
\newcommand{\dH}{\operatorname{dH}}
\newcommand{\PH}{\operatorname{PH}}
\newcommand{\dPH}{\operatorname{dPH}}
\newcommand{\gdPH}{\operatorname{gdPH}}
\newcommand{\gPH}{\operatorname{gPH}}
\newcommand{\PVH}{\operatorname{P_vH}}
\newcommand{\dPVH}{\operatorname{dP_vH}}
\newcommand\br[1]{\{ #1 \}}
\newcommand\dgal[1]{  \left\{\!\!\left\{#1\right\}\!\!\right\} }
\newcommand\dsq[1]{  \left\llbracket#1\right\rrbracket }
\newcommand{\ldb}{\{\!\!\{}
\newcommand{\rdb}{\}\!\!\}}
\newcommand{\DDer}{\mathbb{D}\text{er}} 
\newcommand\dSN[1]{\dgal{#1}_{\operatorname{SN}}}
\newcommand\brSN[1]{\br{#1}_{\operatorname{SN}}}

\definecolor{light}{gray}{.9}

\newglossaryentry{Asharp}{name={\ensuremath{\cA_\sharp}},
  description={Vector space $\cA/[\cA,\cA]$ for an algebra $\cA$}}

\newglossaryentry{Arep}{name={\ensuremath{\cA_{\bf n}}},
  description={Representation algebra of $\cA$ of dimension $\bf n$}}

\newglossaryentry{Balg}{name={\ensuremath{B}},
  description={A base algebra, usually  $B=\oplus_{s\in S} \kk e_s$}}

\newglossaryentry{BRA}{name={\ensuremath{\BRA(\cA)}},
  description={Graded vector space $\BRA(\cA)=\oplus_{n\geq 1} \BRA(\cA)_n$ with $n$-brackets in degree $n$; its completion $\wBRA(\cA)$ contains $\cA_\sharp$ in degree $0$}}

\newglossaryentry{PolyLamb}{name={\ensuremath{C(V)}},
  description={The space $C(V)=\bigoplus_{n\in\mb Z_{\geq0}}C^n(V)$ of all poly-$\lambda$-brackets of a commutative differential algebra $V$}}

\newglossaryentry{PolyLambG}{name={\ensuremath{C(V)^G}},
  description={The subspace $C(V)^G \subset C(V)$ of $G$-invariant poly-$\lambda$-brackets of $V$}}

\newglossaryentry{AsPolyLamb}{name={\ensuremath{C(\mc V)}},
  description={The space $C(\mc V)=\bigoplus_{n\in\mb Z_{\geq0}}C^n(\mc V)$ of all $n$-fold $\lambda$-brackets, $n\geq 0$, of an associative differential algebra $\mc V$}}  

\newglossaryentry{shorthand-C}{
name={\ensuremath{C_{u_1 v_1,\ldots,u_k v_k} }},
description={shorthand notation for the element
$c^{(1)}_{u_1 v_1} \cdots  c^{(k)}_{u_k v_k} \in \cA_{\bf n}$, where
$C=c^{(1)}\otimes\dots\otimes c^{(k)}\in\mc A^{\otimes k}$}
}  

\newglossaryentry{casimir}{
name={\ensuremath{\Cas(\mc V)}},
description={The space of Casimir elements of a dPVA $\mc V$: the center of the Lie algebra action of $\mc V_\sharp$ on $\mc V$}
}
  
\newglossaryentry{dd}{name={\ensuremath{\dd}},
  description={The differential on $C(V)$ (or $C(\mc V)$, or $\Sigma(\mc V)$) defined from $[-_\lambda-]$ (or $\dsq{-_\lambda-}$)}}

\newglossaryentry{dPnc}{name={\ensuremath{\dd_P}},
  description={Differential on $(\mb T^\ast \cA)_\sharp$ associated with a noncommutative Poisson bivector $P$}}

\newglossaryentry{dPi}{name={\ensuremath{\dd_\Pi}},
  description={Differential on $\mf X(A)$ associated with a Poisson bivector $\Pi$}}

\newglossaryentry{wdd}{name={\ensuremath{\wdd}},
  description={Differential on $\wBRA_B(\cA)$ associated with a double Poisson bracket, or its gauged version on $\widehat{\mc D}_\cA$}}

\newglossaryentry{dHbasV}{name={\ensuremath{\dH_{\textrm{bas}}(\mc V)}},
  description={Basic dPVA cohomology $\coH(\widetilde \Gamma(\mc V),\tilde{\delta})$}} 

\newglossaryentry{dHredV}{name={\ensuremath{\dH_{\textrm{red}}(\mc V)}},
  description={Reduced dPVA cohomology $\coH(\Gamma(\mc V),\delta)$}}

\newglossaryentry{dPH-PV}{name={\ensuremath{\dPH(\cA)}},
  description={Double Poisson cohomology $\coH((\mb T^\ast \cA)_{\sharp},\dd_P)$ according to Pichereau-Van de Weyer}} 

\newglossaryentry{hatdPH}{name={\ensuremath{\widehat{\dPH}(\cA)}},
  description={Completed double Poisson cohomology $\coH(\wBRA_B(\cA),\wdd)$}} 

\newglossaryentry{dPVHA}{name={\ensuremath{\dPVH(\mc V)}},
  description={Variational dPVA cohomology $\coH(C(\mc V),\dd)$}}

\newglossaryentry{D-i}{name={\ensuremath{D_{(i)}}},description={Linear map $D_{(i)}:\mc A^{\otimes n}\rightarrow \mc A^{\otimes m+n+1}$ obtained from $D:\mc A\rightarrow\mc A^{\otimes m}$ by acting on the $i$-th factor of $\mc A^{\otimes n}$, see \eqref{20240805:eq1}}}

\newglossaryentry{D-LR}{name={\ensuremath{D_L\,,D_R}},description={Special cases: $D_L=D_{(1)}$ and $D_{R}=D_{(n)}$}}

\newglossaryentry{calDA}{name={\ensuremath{\mc D_\cA^k,\, \widehat{\mc D}_\cA^k}}, description={$k$-th vector space used in the complex defining the (completed) gauged double Poisson cohomology}}

\newglossaryentry{DDer}{name={\ensuremath{\DDer(\cA)}},
  description={Double derivations over $\mc A$}}

\newglossaryentry{evmap}{name={\ensuremath{\ev}},
  description={Evaluation map $\ev:\overline{\Sigma}(\mc R_\ell)/\partial \overline{\Sigma}(\mc R_\ell) \stackrel{\simeq}{\longrightarrow} \Sigma(\mc R_\ell)$}}

\newglossaryentry{gdPH}{name={\ensuremath{\gdPH(\cA)}},
  description={Gauged double Poisson cohomology $\coH(\mc D_\cA,\dd_P)$}} 
  
\newglossaryentry{HatgdPH}{name={\ensuremath{\widehat{\gdPH}(\cA)}},
  description={Completed double Poisson cohomology $\coH(\widehat{\mc D}_\cA,\wdd)$}} 
  
\newglossaryentry{gaugedPHM}{name={\ensuremath{\gPH_G(M)}},
  description={$G$-gauged Poisson cohomology $\coH({\mf X}(\kk[M/\!/G]) , \dd_\Xi)$ of the variety $M$ endowed with a gauged Poisson bivector $\Xi$}}  

\newglossaryentry{quivhm}{name={\ensuremath{h_{\mathrm{E},C}}}, description={The homotopy operator $h_{\mathrm{E},C} = (L_{\mathrm{E}})^{-1} \circ \iota_{C}$ on $\wBRA_B(\cA)$ from \eqref{Quiv:hm}}}   

\newglossaryentry{PVAhm}{name={\ensuremath{h_{M}}}, description={The homotopy operator $h_{M} = (L_{\Delta,M})^{-1} \circ \iota_{\Delta,M}$  on $\widetilde{\Omega}(\mc V)$ from \eqref{eq:hm}}} 

\newglossaryentry{HCgW}{name={\ensuremath{\coH(\g;W)}},
  description={Lie algebra cohomology of the Lie algebra $\g$ acting on $W$}} 

\newglossaryentry{HbasV}{name={\ensuremath{\coH_{\textrm{bas}}(V)}},
  description={Basic PVA cohomology $\coH(\widetilde \Gamma(V),\tilde{\delta})$}}

\newglossaryentry{HCEA}{name={\ensuremath{\coH_{CE}(A)}},
  description={Poisson cohomology of the Poisson algebra $(A,\br{-,-})$}} 

\newglossaryentry{HredV}{name={\ensuremath{\coH_{\textrm{red}}(V)}},
  description={Reduced PVA cohomology $\coH(\Gamma(V),\delta)$}} 

\newglossaryentry{jet-comm}{
name={\ensuremath{J_{\infty}A}},
description={Jet algebra of a commutative algebra $A$}
}

\newglossaryentry{JetAss}{
name={\ensuremath{J_\infty \mc A}},
description={Jet algebra of the associative algebra $\mc A$}
}  

\newglossaryentry{quivL}{name={\ensuremath{L_{\mathrm{E}}}}, description={The Lie derivative on $\wBRA_B(\cA)$ of the Euler operator $\mathrm{E}$}}   

\newglossaryentry{PVAL}{name={\ensuremath{L_{\Delta,M}}}, description={The even derivation of $\widetilde{\Omega}(\mc V)$ from \eqref{Eq:Lm}}}   

\newglossaryentry{mult}{name={\ensuremath{\mult}},description={Multiplication map $\mult:\mc A^{\otimes n}\rightarrow\mc A$}}

\newglossaryentry{mult-i}{name={\ensuremath{\mult_{(i,i+1)}}},description={Extension of the multiplication map $\mult:\mc A\otimes\mc A\rightarrow\mc A$ to
the linear map $\mult_{(i,i+1)}:\mc A^{\otimes m}\rightarrow \mc A^{\otimes m-1}$ defined in \eqref{eq:mii+1}}}

\newglossaryentry{ProjOp}{name={\ensuremath{\widetilde{\proj}_n,\proj_n}},
  description={Projection operators $\widetilde{\proj}_n:\widetilde{\Gamma}^n(\mc V)\to C^n(\mc V)$ and $\proj_n:\Gamma^n(\mc V)\rightarrow C^n(\mc V)$}} 

\newglossaryentry{PHA}{name={\ensuremath{\PH(A)}},
  description={Poisson cohomology of $A$ based on a Poisson bivector $\Pi$}}

\newglossaryentry{qPcoh}{name={\ensuremath{\PH_G(\kk[M])}},
  description={Quasi-Poisson cohomology $\coH(\mf X(\kk[M])^{G},\dd_\Xi)$ of the quasi-Poisson $G$-variety $(M,\Xi)$}}

\newglossaryentry{PVHA}{name={\ensuremath{\PVH(V)}},
  description={Variational PVA cohomology $\coH(C(V),\dd)$}}

\newglossaryentry{gPVHA}{name={\ensuremath{\PVH_G(V)}},
  description={$G$-invariant variational PVA cohomology of $V$}}

\newglossaryentry{qV}{
name={\ensuremath{q(V)}},
description={Poisson algebra associated to a PVA $V$}
}

\newglossaryentry{qV-nc}{
name={\ensuremath{q(\mc V)}},
description={Double Poisson algebra associated to a dPVA $\mc V$}
}

\newglossaryentry{Rell}{name={\ensuremath{\mc R_\ell}},
  description={The algebra of non-commutative differential polynomials 
  $\mc R_\ell=\kk\langle u_i^{(n)}\mid i=1,\dots,\ell,n\in\mb Z_+\rangle$}}  
  
\newglossaryentry{RepSp}{name={\ensuremath{\operatorname{Rep}_B(\cA, {\bf n})}},
  description={Representation space of $\cA$ relative to $B$ of dimension $\bf n$}}

\newglossaryentry{TstarA}{name={\ensuremath{\mb T^\ast \cA}},
  description={Graded tensor algebra $T^\ast \cA:=T_{\cA} \DDer_B(\cA)$}}
  
\newglossaryentry{Vsharp}{name={\ensuremath{\mc V_\sharp}},
  description={Vector space $\mc V/([\mc V,\mc V]+\partial \mc V)$ for a differential algebra $\mc V$}}

\newglossaryentry{basic-coch}{
name={\ensuremath{X_{\lambda_1,\dots,\lambda_n}}},
description={Basic $n$-cochain}
}

\newglossaryentry{X-s}{name={\ensuremath{X^{(s)}_{\lambda_1,\dots,\lambda_n}}},description={Extension of the map $X_{\lambda_1,\dots,\lambda_n}:\mc V^{\otimes n}\rightarrow\mc V^{\otimes n+1}[\lambda_1,\dots,\lambda_n]$ to
the linear map defined in \eqref{eq:XL} as $X_{\lambda_1,\dots,\lambda_n}^{(s)}:\mc V^{\otimes n+m}\rightarrow\mc V^{\otimes n+m+1}[\lambda_1,\dots,\lambda_n]$ }}

\newglossaryentry{mfXA}{name={\ensuremath{\mf X^n(A)}}, description={Subspace of $\Hom(\wedge^n A,A)$ made of skewsymmetric multilinear derivations on the commutative algebra $A$}}

\newglossaryentry{vect-partial}{name={\ensuremath{\Vect^\partial(\mc V)}},
description={The space of derivations of the associative product of $\mc V$ commuting with $\partial$}
}

\newglossaryentry{ZPA}{name={\ensuremath{Z_P(\cA_\sharp)}}, description={Centre of the Lie algebra $(\cA_\sharp,\br{-,-}_{P,\sharp})$}}

\newglossaryentry{ZPAA}{name={\ensuremath{Z_P(\cA_\sharp;\cA)}}, description={Elements of $\cA_\sharp$ for which $\br{\bar{a},-}_P\in \Der(\mc A)$ vanishes.}}

\newglossaryentry{center-nc}{
name={\ensuremath{Z(\mc V)}},
description={Center of the dPVA $\mc V$}
}

\newglossaryentry{gammapTil}{name={\ensuremath{\widetilde{\gamma}}},
  description={The map in the identification $\widetilde{\gamma}:\widetilde{\Gamma}(\mc R_\ell)\stackrel{\simeq}{\longrightarrow} \widetilde{\Omega}(\mc R_\ell)$ 
  leading to $\gamma:\Gamma(\mc R_\ell)\stackrel{\simeq}{\longrightarrow} \Omega(\mc R_\ell)$}}
 
\newglossaryentry{Gammatilde}{name={\ensuremath{\widetilde{\Gamma}(V)}},
  description={The superalgebra $\widetilde{\Gamma}(V)
=\bigoplus_{n\in\mb Z_{\geq0}}\widetilde{\Gamma}^n(V)
$ of all basic $n$-cochains, $n\geq0$, of a commutative differential algebra $V$}}

\newglossaryentry{GammaRed}{name={\ensuremath{\Gamma(V)}},
  description={The space $\Gamma(V)=\widetilde{\Gamma}(V)/\partial \widetilde{\Gamma}(V)$ of reduced $n$-cochains of $V$}}

\newglossaryentry{AsGammatilde}{name={\ensuremath{\widetilde{\Gamma}(\mc V)}},
  description={The space $\widetilde{\Gamma}(\mc V)$ of all basic $n$-cochains, $n\geq0$, of an associative differential algebra $\mc V$}}

\newglossaryentry{AsGammaRed}{name={\ensuremath{\Gamma(\mc V)}},
  description={The space $\Gamma(\mc V)
=\widetilde{\Gamma}(\mc V)/(\partial\widetilde{\Gamma}(\mc V)
+[\widetilde{\Gamma}(\mc V),\widetilde{\Gamma}(\mc V)])$ of all reduced $n$-cochains, $n\geq0$, of $\mc V$}}

\newglossaryentry{Tildelta}{name={\ensuremath{\tilde \delta}},
  description={The differential on $\widetilde{\Gamma}(V)$ (or $\widetilde{\Gamma}(\mc V)$, or $\widetilde{\Omega}(\mc V)$)}}
  
\newglossaryentry{deltadiff}{name={\ensuremath{\delta}},
  description={The differential on $\Gamma(V)$ (or $\Gamma(\mc V)$, or $\Omega(\mc V)$) obtained by reduction from $\tilde \delta$}}

\newglossaryentry{DelGauge}{name={\ensuremath{\Delta_s}}, description={The $s$-th gauge element $\Delta_s\in \DDer(\cA)$}}

\newglossaryentry{quiviota}{name={\ensuremath{\iota_{C}}}, description={The contraction operator on $\wBRA_B(\cA)$ from \eqref{Quiv-iota}}}  

\newglossaryentry{PVAiota}{name={\ensuremath{\iota_{\Delta,M}}}, description={The contraction operator on $\widetilde{\Omega}(\mc V)$  from \eqref{Eq:iotam}}}   

\newglossaryentry{piDelta}{name={\ensuremath{\iota^\Delta_k,\,\pi^\Delta_k,\, \widehat{\pi}^\Delta_k}}, description={Morphisms used in the (completed) gauged double Poisson cohomology}} 

\newglossaryentry{sigmaction}{name={\ensuremath{\sigma_\tau}}, description={Action on $\mc A^{\otimes n}$ by permutation of factors according to $\tau\in S_n$. We denote $\sigma:=\sigma_{(12\cdots n)}$}}

\newglossaryentry{Sigmarray}{name={\ensuremath{\Sigma(\mc V)}}, description={The space of skewsymmetric arrays $\Sigma(\mc V)$}} 

\newglossaryentry{Sigmabar}{name={\ensuremath{\overline\Sigma(\mc V)}}, description={Subspace of $\widetilde \Omega(\mc V)$ such that  $\overline\Sigma(\mc V)/\partial \overline\Sigma(\mc V) \simeq \Sigma(\mc V)$}} 

\newglossaryentry{Sigmatilde}{name={\ensuremath{\widetilde\Sigma(\mc V)}}, description={Subspace in the splitting $\widetilde \Omega(\mc V)=\widetilde\Sigma(\mc V)\oplus [\widetilde \Omega(\mc V),\widetilde \Omega(\mc V)]$}}

\newglossaryentry{phimap}{name={\ensuremath{\widetilde{\phi}_n}},
  description={The map  $\widetilde{\phi}_n:\widetilde{\Omega}^n(\mc R_\ell)\to\Sigma^n(\mc R_\ell)$ 
  leading to the identification $\phi:\Omega(\mc R_\ell)\stackrel{\simeq}{\longrightarrow} \Sigma(\mc R_\ell)$}}

\newglossaryentry{psimap}{name={\ensuremath{\psi}},
  description={The map in the identification $\psi:C(\mc R_\ell)\stackrel{\simeq}{\longrightarrow} \Sigma(\mc R_\ell)$}}

\newglossaryentry{Omegarray}{name={\ensuremath{\widetilde \Omega(\mc V)}}, description={The space of arrays $\widetilde \Omega(\mc V)$}} 

\newglossaryentry{Redomegarray}{name={\ensuremath{\Omega(\mc V)}}, description={Reduced space $\Omega(\mc V)=\widetilde \Omega(\mc V)/(\partial \widetilde \Omega(\mc V) + [\widetilde \Omega(\mc V),\widetilde \Omega(\mc V)])$}}

\newglossaryentry{bullet}{name={\ensuremath{\bullet}},
description={The associative product in
$\mc A\otimes\mc A^{\mathrm{op}}$. Also used to denote the left action of $\mc A\otimes\mc A$ on $\mc A^{\otimes n}$}}

\newglossaryentry{inner}{name={\ensuremath{\star}},
description={Inner bimodule structure on $\mc A\otimes \mc A$}}

\newglossaryentry{bullet-i}{name={\ensuremath{\bullet_{(i,i+1)}}},
description={Right action of $\mc A\otimes\mc A$ on $\mc A^{\otimes n}$}}

\newglossaryentry{inner-j}{name={\ensuremath{\star_j}},
description={Left and right $\mc A$-module structure of $\mc A^{\otimes n}$ (it reduces to the outer bimodule structure for $j=0,n$)}}

\newglossaryentry{otimes-i}{name={\ensuremath{\otimes_i}},description={Linear map $-\otimes_i C:\mc A^{\otimes n}\rightarrow\mc A^{\otimes n+m}$, $C\in\mc A^{\otimes m}$, see \eqref{eq:tensor-i-notation}}}

\newglossaryentry{tint}{name={\ensuremath{\tint}},description={Quotient map $\mc V\rightarrow\mc V_{\sharp}$. Also used to denote its extension to the quotient maps $\widetilde{\Gamma}(\mc V)\rightarrow\Gamma(\mc V)$ and $\widetilde{\Omega}(\mc V)\rightarrow\Omega(\mc V)$}}

\newglossaryentry{lambda-b}{name={\ensuremath{\{\cdot_\lambda\cdot\}}},description={$\lambda$-bracket on a commutative differential algebra. It is denoted by $[\cdot_\lambda\cdot]$ when we want to distinguish it from a poly-$\lambda$-bracket}}

\newglossaryentry{n-lambda}{name={\ensuremath{
\{\cdot_{\lambda_1}\!\dots_{\lambda_{n\!-\!1}}\!\cdot\}}},
description={$n$-$\lambda$-bracket on a commutative differential algebra}
}

\newglossaryentry{dPbracket}{name={\ensuremath{\ldb\cdot,\cdot\rdb}},description={Double (Poisson) bracket on an algebra. It is denoted by $\llbracket-,-\rrbracket$ when we want to distinguish it from an
$n$-bracket}}

\newglossaryentry{lambda-b-nc}{name={\ensuremath{\ldb\cdot_\lambda\cdot\rdb}},description={$2$-fold $\lambda$-bracket on a differential algebra. It is denoted by $\llbracket\cdot_\lambda\cdot\rrbracket$ when we want to distinguish it from an
$n$-fold $\lambda$-bracket}}

\newglossaryentry{n-lambda-nc}{name={\ensuremath{
\ldb\cdot_{\lambda_1}\!\dots_{\lambda_{n\!-\!1}}\!\cdot\rdb
}},
description={$n$-fold $\lambda$-bracket on a differential algebra}
}

\newglossaryentry{partialder-nc}{
name={\ensuremath{\frac{\partial}{\partial u_i^{(n)}}}},
description={Partial derivatives of the algebra of differential polynomials $\mc R_\ell$. They are commuting $2$-fold derivations}
}

\begin{document}

\title[dP(v) cohomology]{Double Poisson (vertex) algebra cohomology}
% Other idea: Cohomology theories from double Poisson (vertex) algebras

\author[M.~Fairon]{Maxime Fairon}
\address[Maxime Fairon]{Université Bourgogne Europe, CNRS, IMB UMR 5584, F-21000 Dijon, France}
\email{maxime.fairon@u-bourgogne.fr}
\author[D.~Valeri]{Daniele Valeri}
\address[Daniele Valeri]{Dipartimento di Matematica \& INFN, Sapienza Universit\`a di Roma,
P.le Aldo Moro 5, 00185 Rome, Italy}
\email{daniele.valeri@uniroma1.it}

%\date{today}

\begin{abstract} 
A noncommutative version of Poisson geometry was initiated by Van den Bergh 
by introducing at the level of associative algebras the formalism of \emph{double Poisson brackets}. Their key property is to induce (standard) Poisson brackets under each representation functor. 
Then, Pichereau and Van de Weyer developed and studied the corresponding cohomology theory under the assumption that there exists a noncommutative bivector defining the double Poisson bracket. 
Our first main result is that one can remove this assumption by constructing a completed double Poisson cohomology valid in any situation, hence generalizing the approach of Pichereau and Van de Weyer. As an application, we show that the double Poisson cohomology complex associated to the path algebra of a quiver is acyclic.
Furthermore, we show that this new double Poisson cohomology theory can be adapted to weaker forms of double Poisson brackets (called quasi-Poisson and gauged Poisson), and that it is compatible with each representation functor. 

A second focus of this memoir concerns the formalism of \emph{double Poisson vertex algebras}.
These were introduced 10 years ago by De Sole, Kac and the second author, as noncommutative versions of Poisson vertex algebras, which induce the latter structures under each representation functor. 
Our second main result is the development of cohomology theories for double Poisson vertex algebras. 
These are noncommutative analogues of the basic, reduced and variational Poisson vertex algebra cohomologies. 
More importantly, we prove that under each representation functor these cohomology theories are compatible with their commutative counterparts. As an application, we compute the double Poisson vertex algebra cohomology of the generalized noncommutative de Rham complex
and of the generalized noncommutative variational complex.

Finally, we describe the relation between the double Poisson algebra and double Poisson vertex algebra cohomologies using the jet and quotient functors.
\end{abstract}

\maketitle

\setcounter{tocdepth}{1}

%\addcontentsline{toc}{chapter}{Bibliography}

\tableofcontents

%%%%%%%%%%% NEW CHAPTER %%%%%%%%%%%%%%%
%%%%%%%%%%% NEW CHAPTER %%%%%%%%%%%%%%%
%%%%%%%%%%% NEW CHAPTER %%%%%%%%%%%%%%%
%%%%%%%%%%% NEW CHAPTER %%%%%%%%%%%%%%%

\chapter*{Introduction}
A modern research direction in noncommutative algebraic geometry was established by Kontsevich and Rosenberg \cite{KR00} through the use of the representation functors. 
Namely, one aims at defining a noncommutative structure $\mc S_{\mathrm{nc}}$ on an associative algebra $\cA$ such that it induces a standard structure $\mc S$ on each (commutative) coordinate ring $\cA_N := \kk[\Rep(\cA,N)]$ of the representation scheme $\Rep(\cA,N)$, $N\geq 1$. (For us, the base field $\kk$ is an algebraically closed field of characteristic zero.)
This feature is nowadays called the \emph{Kontsevich-Rosenberg principle} and it was originally illustrated by showing that formally smooth (also called quasi-free) algebras induce smooth representation schemes. The earlier noncommutative symplectic geometry of Kontsevich \cite{Ko} could also be modified to satisfy this principle  \cite{CBEG}. 
This is also the case for Van den Bergh's noncommutative Poisson geometry \cite{VdB1} and the noncommutative Poisson vertex algebras of De Sole, Kac and the second author \cite{DSKV}.  
This memoir aims at introducing various cohomology theories within the last two formalisms, which are all compatible with the Kontsevich-Rosenberg principle.

\subsection*{The Poisson case} Fix an associative algebra $\cA$ over $\kk$. 
The noncommutative Poisson geometry of Van den Bergh \cite{VdB1} is based on the notion of \emph{double bracket}, that is a $\kk$-linear map 
$$\dgal{-,-}:\cA \otimes \cA \to \cA \otimes \cA, \quad 
a\otimes b \mapsto \dgal{a,b},$$
satisfying a skewsymmetry axiom \eqref{Eq:db-skew} and a derivation rule in each argument \eqref{Eq:db-Rleib}-\eqref{Eq:db-Lleib}. 
Such a structure naturally induces a skewsymmetric biderivation on each representation algebra $\cA_N$ (see\footnote{To simplify the exposition in the introduction, we have left the case of a general semisimple base ring (e.g. for quiver algebras) to the main text.} Theorem \ref{Thm:Rep-Dbr}). 
More importantly, Van den Bergh has shown that under a ``Poisson property" valued in $\cA^{\otimes 3}$ given by \eqref{Eq:dJacExpl}, the induced operation is in fact a Poisson bracket. 
Therefore, such a couple $(\cA,\dgal{-,-})$ called a \emph{double Poisson algebra} turns $\cA_N$ into a Poisson algebra, hence it satisfies the Kontsevich-Rosenberg principle. 

It is common in Poisson geometry to work with Poisson bivectors instead of Poisson brackets, and this can usually (though not always!) be done also at the noncommutative level. 
Instead of vector fields, we use the $\cA$-bimodule $\DDer(\cA)$ of \emph{double derivations}, i.e. maps $\cA \to \cA^{\otimes 2}$ satisfying a derivation property, cf. \eqref{DDer}. 
Then, forming the tensor algebra $\mb T^\ast \cA=T_\cA \DDer(\cA)$, 
one can consider noncommutative $n$-vector fields as elements of $(\mb T^\ast \cA)_n$, and 
Van den Bergh has shown that any such element gives rise to an $n$-bracket, i.e. an operation $\cA^{\otimes n}\to \cA^{\otimes n}$ with skewsymmetry and derivation properties, cf. Proposition \ref{Pr:MapMu}. 
In particular, noncommutative bivector fields are elements of $(\mb T^\ast \cA)_2$, and they define double brackets.  
What is important for defining an $n$-bracket is, in fact, the class of a  noncommutative $n$-vector field in $(\mb T^\ast \cA)_\sharp:=\mb T^\ast \cA/[\mb T^\ast \cA,\mb T^\ast \cA]$, i.e. its expression modulo graded commutators. 
Furthermore, the graded vector space $(\mb T^\ast \cA)_\sharp$ inherits a graded Lie bracket $\brSN{-,-}:(\mb T^\ast \cA)_\sharp\times (\mb T^\ast \cA)_\sharp \to (\mb T^\ast \cA)_\sharp$ of degree $-1$ from a ``double analog" of the Schouten-Nijenhuis bracket. 
It is then standard that any element $P\in (\mb T^\ast \cA)_{\sharp,2}$ satisfying 
$\brSN{P,P}=0$ will define a square-zero differential on the complex $(\mb T^\ast \cA)_\sharp$. 
The upshot is that the condition $\brSN{P,P}=0$ for $P\in (\mb T^\ast \cA)_{\sharp,2}$ is precisely the condition that it defines a double \emph{Poisson} bracket.
Thus, the cohomology $\dPH(\cA)$ of the complex $((\mb T^\ast \cA)_{\sharp},\brSN{P,-})$ was naturally called \emph{double Poisson cohomology} by Pichereau and Van de Weyer \cite{PV,VdW}. 
More importantly, this construction is compatible with the Kontsevich-Rosenberg principle, because under each representation functor there is a morphism 
$\dPH(\cA) \to \PH(\cA_N)$ 
to the standard Poisson(-Lichnerowicz) cohomology built from the Poisson bivector $\tr(P)$ acting on skewsymmetric multiderivations $\mc X(\cA_N)$ on $\cA_N$ (and this restricts to $\Gl_N$-invariant elements), cf. \cite{PV,VdW}. 

\medskip 

At this point, the reader may have expected to find more manuscripts studying and computing double Poisson cohomology groups... But there are none! This may indicate that the computations carried out in \cite{PV,VdW} are too technical to be pursued in various examples, and that another perspective is needed. 
Furthermore, the cautious reader has kept in mind that not all double Poisson brackets are defined from an element in $(\mb T^\ast \cA)_{\sharp,2}$, hence they would expect to see a cohomology theory developed for \emph{any} double Poisson brackets. 
Tackling these two remarks is precisely the objective of the first part of this memoir, where we shall introduce and compute the \emph{completed double Poisson cohomology} from Chapter \ref{CH:Gen-dPcoh} onward. 

To motivate this new construction, recall that for a commutative algebra $A$ endowed with a Poisson bracket, one can define its Poisson cohomology $\coH_{CE}(A)$ as the Chevalley-Eilenberg cohomology restricted to the complex of skewsymmetric multiderivations $\mc X(A)$, cf. Subsection~\ref{ss:comPCoh}. 
If the Poisson bracket is obtained from a Poisson bivector $\Pi$ on $\Spec(A)$, then one has a natural identification $\coH_{CE}(A)\simeq \PH(A)$ with the Poisson cohomology defined with the differential $[\Pi,-]$ (with $[-,-]$ being the Schouten-Nijenhuis bracket).  
Thus, for a double Poisson algebra $\cA$, one is led to replace $\mc X(A)$ with 
the complex $\wBRA(\cA)$ having $\cA_\sharp:=\cA/[\cA,\cA]$ in degree zero and 
$n$-brackets in degree $n\geq 1$. The differential on $\wBRA(\cA)$ has to be properly introduced, cf. Definition \ref{def:wdd} and the subsequent theorem. 
This yields the completed double Poisson cohomology $\widehat{\dPH}(\cA)$ of $\cA$, which can be introduced \emph{for any} double Poisson bracket on $\cA$. 
We shall also prove in Section~\ref{ss:ReldPCoh} that, when the double Poisson bracket on $\cA$ is defined from some $P\in (\mb T^\ast \cA)_{2,\sharp}$, there is a relation with the double Poisson cohomology of Pichereau-Van de Weyer given as a morphism $\dPH(\cA)\to \widehat{\dPH}(\cA)$. 
In our opinion, computations for the completed cohomology $\widehat{\dPH}(\cA)$ are easier. 
Thus, when there is an identification of the two theories (which happens when $(\mb T^\ast \cA)_{\sharp}\simeq \wBRA(\cA)$), this provides a new way to compute cohomology classes in $\dPH(\cA)$. 
Finally, we obtain in Subsection~\ref{ss:InducdP} that this new cohomology theory falls within the Kontsevich-Rosenberg principle by constructing a map $\widehat{\dPH}(\cA) \to \coH_{CE}(\cA_N)$ for any $N\geq 1$. 
In the presence of a corresponding $P\in (\mb T^\ast \cA)_{\sharp,2}$, this map is 
compatible with the one $\dPH(\cA) \to \PH(\cA_N)$ of Pichereau and Van de Weyer. 

\medskip 

We mentioned without further comment that a double Poisson bracket not only induces a Poisson structure on $\cA_N$, but also on the invariant algebra $\cA_N^{\Gl_N}$ which is the coordinate ring of the moduli space $\Rep(\cA,N)/\!/\Gl_N$ with $N\geq 1$. 
At a geometric level, it is clear that any invariant element of $\PH(\cA_N)$ descends to $\PH(\cA_N^{\Gl_N})$. However, this map is far from being injective: any invariant multivector field %$\Pi'\in
in $\mc X(\cA_N)$ having a wedge factor given by the infinitesimal action of $\gl_N$ on $\cA_N$ vanishes when restricted to $\cA_N^{\Gl_N}$. 
Therefore, one may need to relax the (completed) double Poisson cohomology on $\cA$ if one is only interested in the Poisson cohomology $\PH(\cA_N^{\Gl_N})$ of $\cA_N^{\Gl_N}$. 
This idea was first considered by Alekseev, Kawazumi, Kuno and Naef~\cite{AKKN}, and we shall pursue it in greater details in Chapter~\ref{CH:qgauge}.
Moreover, we are going to adapt all the previous considerations to the case of double quasi-Poisson algebras, the noncommutative analog of quasi-Poisson spaces~\cite{AKSM}.
For all these instances, we are also going to prove that the new cohomology theories are compatible with the ``$\Gl_N$-invariant" Kontsevich-Rosenberg principle as there is always a morphism to $\PH(\cA_N^{\Gl_N})$, $N\geq 1$, cf. Section \ref{Sec:InducGauge}.

\medskip 

To conclude the Poisson case, let us stress that the cohomology theories that we introduce are not only mathematically pleasing, but that they have a potential of applications for seeking integrable systems on the representation spaces $\Spec(\cA_N)$ or their moduli $\Spec(\cA_N^{\Gl_N})$ for all $N\geq 1$. 
Indeed, recall that the bi-Hamiltonian approach  to integrability (see e.g. \cite{MCFP}) is based on having a pair of Poisson brackets $\br{-,-}_i$, $i=0,1$, that are compatible: any linear combination of them is again a Poisson bracket. 
The compatibility condition can be reformulated as saying that $\br{-,-}_1$ is a cocycle in the Poisson cohomology defined by $\br{-,-}_0$. 
Thus, when $\cA$ is a double Poisson algebra, our previous considerations entail that any double bracket in $\widehat{\dPH}^2(\cA)$ which is Poisson will induce a compatible Poisson bracket on $\Spec(\cA_N)$ and $\Spec(\cA_N^{\Gl_N})$ for all $N \geq 1$. 
As a byproduct, one obtains a pair of compatible double Poisson brackets, as introduced in \cite{ORS}.

\subsection*{The Poisson vertex case}

Recall from \cite{BDSK} that a Poisson vertex algebra (PVA) is a commutative differential algebra $(V,\partial)$ endowed with a $\kk$-linear map 
$$V\otimes V \to V[\lambda], \qquad a\otimes b \mapsto \br{a_\lambda b}\,,$$
satisfying sesquilinearity and skewsymmetry, together with Leibniz rules and a version of Jacobi identity, see the beginning of Section \ref{sec:pva}.
%Sesquilinearity establishes a compatibility of the  $\lambda$-bracket $\br{-_\lambda -}$ with the differential $\partial$, and a polynomial in the indeterminate $\lambda$ should be viewed as the symbol of an operator $p(\partial)\in V[\partial]$.  
There are two natural motivations \cite{KacPisa} for PVA: 
on one hand they appear as quasi-classical versions of vertex algebras, on the other hand
they provide a convenient
algebraic framework for studying Hamiltonian PDEs. 
These two points are analogs of the fact that Poisson algebras appear as quasi-classical limits of (filtered) associative algebras, and that they are used for defining and studying Hamiltonian ODEs.  
In fact, a PVA version of bi-Hamiltonian geometry has been exploited in \cite{DSKV} and several applications to bi-Hamiltonian integrable systems have been later developed. As in the Poisson case, compatibility of a pair  of $\lambda$-brackets $\br{-_\lambda-}_i$, $i=0,1$, on $V$ is a condition that has a cohomological interpretation: 
this is equivalent to require that $\br{-_\lambda-}_1 \in \PVH^2(V)$, where
$\PVH^2(V)$ is the second \emph{variational PVA cohomology} space, see \cite{BDSK20}, associated to the $\lambda$-bracket $\br{-_\lambda-}_0$ (and vice-versa).  
For computations, the variational PVA cohomology can be technical, so this is why it is also important to consider the basic PVA cohomology $\coH_{\mathrm{bas}}(V)$  and its reduced version $\coH_{\mathrm{red}}(V)$ \cite{DSK13}, cf.   Chapter \ref{Ch:PVAcoh} for a brief review. One should think of the variational $\PVH(V)$ and reduced $\coH_{\mathrm{red}}(V)$ cohomologies as counterparts to $\coH_{CE}(A)$ and $\PH(A)$ on a Poisson algebra, although they are not always equivalent and one only has a map $\coH_{\mathrm{red}}(V)\to \PVH(V)$. 
Nevertheless, in many cases of interest such as for algebras of differential polynomials in finitely many variables,  this map becomes an isomorphism $\coH_{\mathrm{red}}(V)\simeq \PVH(V)$, cf. \cite{BDSHKV}.

\medskip 

The importance of PVA and the development of noncommutative Poisson geometry formed an  incentive for investigating noncommutative versions of PVA.   
To do so
De Sole, Kac and the second author \cite{DSKV} defined 
a \emph{double Poisson vertex algebra} (dPVA)
by merging the notions
of a double Poisson algebra and of a PVA. 
A dPVA is an associative differential (noncommutative) algebra $(\mc V,\partial)$ endowed with a $\kk$-linear map $$\dgal{-_\lambda-}:\mc V \otimes \mc V \to \mc V \otimes \mc V, \quad 
a\otimes b \mapsto \dgal{a_\lambda b}\,,$$
called a $2$-fold $\lambda$-bracket, whose precise axioms are given in Section \ref{sec:dPVA}, see Definition \ref{20140606:def-2}. 
Any derivation naturally lifts to each representation algebra, thus 
the pair $(\mc V_N,\partial)$ is a commutative differential algebra for any $N\geq 1$. 
It turns out \cite{DSKV} (cf. Section \ref{sec:dPVAtoPVA-1}) that the $2$-fold $\lambda$-bracket on $\mc V$, induces a PVA structure on
$(\mc V_N,\partial)$
in agreement with the Kontsevich-Rosenberg principle. 

The structure of dPVA can be used to address the study of integrability of non-abelian PDEs through a noncommutative Lenard-Magri scheme \cite{DSKV}. 
Since this scheme is inherent to the bi-Hamiltonian formalism, it is natural to expect that cohomology theories for dPVA should be defined to seek pairs of compatible dPVA structures. 
The introduction of such cohomologies was our second objective. 
In Section \ref{sec:omega1}, we consider the space $\widetilde{\Gamma}(\mc V)$ of basic cochains, 
whose elements of degree $n$ are linear maps $\mc V^{\otimes n}\to \mc V^{\otimes (n+1)}[\lambda_1,\ldots,\lambda_n]$ satisfying sesquilinearity
\eqref{eq:sesquimaps} and Leibniz rules \eqref{eq:Leibnizmaps}. However, differently from the commutative case, \emph{no} skewsymmetry properties are required. 
If $\mc V$ is a dPVA, its $2$-fold $\lambda$-bracket turns $\widetilde{\Gamma}(\mc V)$ into a complex, from which the basic dPVA cohomology $\dH_{\mathrm{bas}}(\mc V)$ is obtained, cf. Definition \ref{def:bas-dPVA}. 
Note that we can extend the derivation on $\mc V$ to an (even) derivation on $\widetilde{\Gamma}(\mc V)$ and make it a graded differential algebra; this structure is compatible with the (square-zero) differential $\tilde{\delta}$ defining the cohomology  $\dH_{\mathrm{bas}}(\mc V)$. 
Hence, one can induce $\tilde{\delta}$ onto the quotient space of $\widetilde{\Gamma}(\mc V)$ modulo  graded commutators and the image of $\partial$, which leads to the reduced dPVA cohomology 
$\dH_{\mathrm{red}}(\mc V)$. 

A slightly different approach consists in starting with the space
of (skewsymmetric) cochains $C(\mc V)$, 
whose elements of degree $n$ are $n$-fold $\lambda$-brackets defined in \cite{DSKV}; these are linear maps $\mc V^{\otimes n}\to \mc V^{\otimes n}[\lambda_1,\ldots,\lambda_{n-1}]$ satisfying sesquilinearity axioms \eqref{20140702:eq4} and \eqref{20140702:eq5}, Leibniz rules \eqref{20140702:eq6}, \emph{together with} the skewsymmetry axiom \eqref{eq:nfold-skew}.
As in the previous cases, if $\mc V$ is a dPVA we can see $C(\mc V)$ as a complex, and we obtain the 
variational dPVA cohomology $\dPVH(\mc V)$, cf. Definition \ref{def:var-dPVA}.
The alert reader will have noticed that $\widetilde{\Gamma}(\mc V)$ and $C(\mc V)$ differ in each degree $n$ by a tensor factor and  $\lambda_n$. 
One can in fact construct a ``projection operator'' $\widetilde{\proj}:\widetilde{\Gamma}(\mc V) \to C(\mc V)$ by cyclically multiplying two consecutive tensor factors, see \eqref{eq:Pn}. 
This map is a morphism of complexes which factors through the quotient 
thus inducing a morphism of complexes $\proj:\Gamma(\mc V)\rightarrow C(\mc V)$. The map $\proj$ induces a map in cohomology 
$\dH_{\mathrm{red}}(\mc V)\to \dPVH(\mc V)$; this is a dPVA analogue of 
the map $\dPH(\cA)\to \widehat{\dPH}(\cA)$ between the double Poisson cohomology and its completed version and of the map $\coH_{\mathrm{red}}(V)\to \PVH(V)$ previously mentioned in the PVA framework. 
The latter can be an isomorphism, hence it is not a surprise that we can get 
$\dH_{\mathrm{red}}(\mc V)\simeq \dPVH(\mc V)$ in some cases, e.g. when $\mc V$ is an algebra of noncommutative differential polynomials
in finitely many variables, cf. Proposition \ref{20240819:prop2}.
Finally, we prove as part of Section \ref{sec:dPVAtoPVA} 
that the variational dPVA cohomology satisfies the Kontsevich-Rosenberg principle: 
under each representation functor, we get a morphism of cochain complexes 
$C(\mc V)\to C(\mc V_N)$, which descends to a linear map $\dPVH(\mc V)\to \PVH(\mc V_N)$ in cohomology for any $N\geq 1$. (An analogous statement holds for the complex $\Gamma(\mc V)$, but we omitted this result since the memoir is already sufficiently long.)

As an application of the cohomology theories developed, in Chapter \ref{CH:PVAexa} we compute the basic and reduced dPVA cohomologies (the latter is isomorphic to the variational dPVA cohomology) for a constant $2$-fold $\lambda$-bracket on an algebra of noncommutative differential polynomials which are noncommutative analogues of the ``generalized de Rham complex"
and ``variational complex" in \cite{DSK13}. The main feature of this computation is that, unlike the commutative case, the cohomology spaces have infinite dimension.

Our approach to dPVA cohomology has been motivated by the paper by De Sole and Kac \cite{DSK13}. This approach allowed us to make a parallel treatment
between the double Poisson and double Poisson vertex cases. However, recently in \cite{BDSHK}, Bakalov, De Sole, Kac and Heluani developed an operadic approach to the theory of Poisson vertex algebras (and vertex algebras) cohomology which has been very useful in the computation of the cohomology of several PVA, see \cite{BDSK20}. We plan to investigate the operadic approach to dPVA in a subsequent work.

Finally, concerning the Poisson vertex case, we have not developed the quasi-Poisson and gauged dPVA cohomology theories. We hope that this gap could be filled in the future.

\subsection*{Relating the two cases} 

Given a commutative algebra $A$, one can define its jet algebra $J_\infty A$ as a universal differential commutative algebra associated with $A$ \cite{AM}. 
It was shown by Arakawa \cite{Ar12} that when $A$ is a Poisson algebra, $J_\infty A$ is naturally equipped with a PVA structure. 
In particular, it is a folklore result (cf. Corollary \ref{Cor:PAPVA1}) that we get a map in cohomology $\coH_{CE}(A)\to \PVH(J_\infty A)$. 
An analogous statement holds in cohomology if we consider the quotient construction 
$V\mapsto V/\langle \partial V \rangle$ that produces a Poisson algebra from a PVA. 

The jet and quotient functors, $\cA \mapsto J_\infty \cA$ and $\mc V \mapsto \mc V /\langle \partial \mc V \rangle$,  were considered in the noncommutative (associative) setting by Bozec, Moreau and the first author \cite{BFM} to relate double Poisson algebras and dPVA. 
In analogy with the commutative setting, one can expect maps in cohomology induced by these constructions. This is precisely what we achieve in  Chapter \ref{Ch:dPA-dPVA}, where we obtain linear maps 
\[
\widehat{\dPH}(\cA)\to \dPVH(J_\infty \cA), \qquad 
\dPVH(\mc V) \to \widehat{\dPH}(\mc V /\langle \partial \mc V \rangle)\,
\]
between the double Poisson algebra and the variational dPVA cohomologies constructed in the previous parts of this memoir. 
Furthermore, in view of the previous discussions, one should also prove that these maps between the noncommutative cohomology theories induce their commutative counterparts under the representation functors. This final aim will be delivered as part of Theorem \ref{Thm:RepJetQuot}.

\section*{Outline of the memoir}

This memoir starts with Chapter \ref{Ch:Setup}, where we gather our conventions and many elementary but useful results in multilinear algebra.  
Then, the first part (Chapters \ref{Ch:StandPoiss} to \ref{CH:rep-dPA}) is concerned with double Poisson algebra cohomologies, 
the second part (Chapters  \ref{Ch:PVAcoh} to \ref{Ch:repVardPVA}) with double Poisson vertex algebra cohomologies, 
and the third part (Chapters \ref{Ch:PA-PVA} and \ref{Ch:dPA-dPVA})  concerns the relation between the previous two. 

%% chapt 2
In Chapter \ref{Ch:StandPoiss}, we recall the standard construction of the Poisson cohomology of a Poisson algebra $A$ using either the Chevalley-Eilenberg cohomology approach for Lie algebras, or the one based on the Schouten-Nijenhuis bracket on skewsymmetric multilinear derivations. 
Then, we generalize the quasi-Poisson cohomology of \cite{AKSM} defined in the presence of a group action by considering \emph{gauged Poisson} invariant bivector (cf. Definition \ref{Def:gaugBiv}), which lead to the cohomology theory in Definition \ref{Def:gPH}. 
To be exhaustive, we also present a ``bivector-free" approach to quasi-Poisson cohomology using Theorem \ref{Thm:qPcohGen}.  

%% chapt 3
In Chapter \ref{Ch:dPA}, we start with the noncommutative setting by introducing the notions of $n$-brackets and noncommutative multivector fields, following Van den Bergh \cite{VdB1}, which span the graded vector spaces\footnote{To ease the presentation, we work over the base field $B=\kk$ for the outline of the first part. In the text, we work over a more general semisimple ring.} $\wBRA(\cA)$ and $\mb T^\ast \cA$.
Double Poisson brackets are elements of $\wBRA(\cA)_2$ introduced in Definition \ref{def:dbr-Poiss}. 
When a double Poisson bracket is defined by an element $P\in (\mb T^\ast \cA)_2$, 
it gives rise to a square zero differential on $(\mb T^\ast \cA)_\sharp$ as noticed by Pichereau and Van de Weyer \cite{PV,VdW}. We recall the corresponding \emph{double Poisson cohomology} $\dPH(\cA)$ in Definition \ref{Def:dPH} and we end with some elementary results.

%% chapt 4
In Chapter \ref{CH:Gen-dPcoh}, we introduce and investigate our new notion of \emph{completed double Poisson cohomology} $\widehat{\dPH}(\cA)$ (cf. Definition \ref{Def:dPH-comp}). 
Given a double Poisson bracket, it is obtained using the differential $\wdd$ on $\wBRA(\cA)$ given in Definition \ref{def:wdd} which is well-defined and squares to zero by Theorem \ref{Thm:g-dPcoh1}. 
Then, we prove in Theorem \ref{Thm:g-dPcoh2} that the map $(\mb T^\ast \cA)_\sharp \to \wBRA(\cA)$ is a morphism of complexes if one considers the differential of Pichereau and Van de Weyer, and the differential $\wdd$, respectively.
We also present in Proposition \ref{Pr:Chemla} an elegant form for $\wdd$ that we call Chemla's formula, and we establish its relation to the double Lie-Rinehart algebras cohomologies of Chemla \cite{Ch}.
We finish by explaining how fusion of double Poisson algebras induces maps in cohomology, see Corollary \ref{Cor:Fus}. 

%% chapt 5
In Chapter \ref{CH:qgauge}, we introduce several cohomologies analogous to the (completed) double Poisson cohomology. 
We start by generalizing the noncommutative Poisson cohomology of \cite{AKKN}. 
What we call \emph{(completed) gauged double Poisson cohomology} ($\gdPH(\mc A)$ and $\widehat{\gdPH}(\mc A)$, in Definition \ref{Def:gdPH}) is obtained by inducing the differentials on $(\mb T^\ast \cA)_\sharp$ and $\wBRA(\cA)$ 
to the quotients $\mc D_{\cA}$ and $\widehat{\mc D}_{\cA}$ by particular subspaces based on the gauge elements $\Delta_s$ in \eqref{Eq:Deltas}.  
Next, we consider double quasi-Poisson brackets \cite{VdB1}, and we prove in Propositions \ref{Pr:dquasiPcoh1} and \ref{Pr:dquasiPcoh2} that the differential of Pichereau and Van de Weyer on $(\mb T^\ast \cA)_\sharp$ and $\wdd$ on $\wBRA(\cA)$ are still squaring to zero in that setting.
Again, the map $(\mb T^\ast \cA)_\sharp \to \wBRA(\cA)$ is a morphism of complexes and therefore it descends in cohomology, see Corollary \ref{Cor:Iso-dquasiPcoh}. 
Finally, we relate these two approaches in Section \ref{ss:RevGDP} through the notion of gauged Poisson element in $(\mb T^\ast \cA)_{\sharp,2}$ (Definition \ref{Def:gaugDBR}), which is the weakest structure leading to a square-zero differential on $\mc D_{\cA}$. 

%% chapt 6
In Chapter \ref{CH:dPA-examples}, we present many examples of double Poisson cohomology groups (and their variants).
The first 3 groups in the completed double Poisson cohomology are determined for \emph{all} double Poisson brackets on $\kk[x]$ in Subsection~\ref{ss:dP-coh-x}, and on $\kk[x]/(x^r)$ in Subsection~\ref{ss:dP-coh-xr}.
We also calculate all the groups in the gauged double Poisson cohomology $\gdPH(\kk[x])$ in Proposition \ref{Pr:gdPH-kx0}. 
Then, these results on $\kk[x]$ are adapted to the quasi-Poisson and the gauged Poisson setting. 
After that, we consider the quartic double Poisson bracket \eqref{Eq:dP-uv} on $\kk\langle u,v\rangle$ and compute in Proposition \ref{Pr:dPH-kuv} the first two groups in $\widehat{\dPH}(\kk\langle u,v\rangle)$.
This double Poisson bracket descends to $\kk\langle u,v\rangle/(u^2,v^2)$ where both the double Poisson cohomology and its completed version can be considered; this provides a first example where these two cohomologies are defined but not isomorphic.
Next, we study in Section \ref{Sec:Quiv} the case of a non-degenerate double Poisson bracket on the path algebra of a quiver (without relations). 
We characterize all the completed double Poisson cohomology groups as part of Theorem \ref{Thm:DPcoh-Quiv}, in particular establishing that they are trivial in positive degree. 
This chapter end with Section \ref{ss:OthExmp} by gathering all previously known examples of double Poisson cohomology (and its variants) while giving a new point of view on some of them.

%% chapt 7
In Chapter \ref{CH:rep-dPA}, we explain how all the constructions are compatible with the Kontsevich-Rosenberg principle. 
The general framework to define a skewsymmetric multilinear bracket on the $N$-th representation algebra\footnote{We work over $B=\kk$ in this outline so the array ${\bf n}$ in the text corresponds to $N$.} 
$\cA_N$ from an $n$-bracket (in particular how to get a Poisson bracket from a double Poisson bracket) is recalled in Section~\ref{sec:7.1}. 
This also guarantees the existence of maps from $\mc D_{\cA}^k$ and $\widehat{\mc D}_{\cA}^k$ to $\mf X^k(\cA_N^{\Gl_N})$. 
Then, we consider all the ways to induce a Poisson cohomology from a noncommutative one using such maps. 
We first recall in Theorem \ref{Thm:dP-rep1} how to obtain the map $\dPH(\cA)\to \PH(\cA_N)$ of Pichereau and Van de Weyer. 
This is adapted to the completed case as part of  Theorem \ref{Thm:dP-rep2} where we get a map 
$\widehat{\dPH(\cA)}\to \coH_{CE}(\cA_N)$; both ways to induce a Poisson cohomology are in agreement in view of Corollary \ref{Cor:dPH-PH}. 
Similarly, we show that we get a map $\gdPH(\cA)\to \PH(\cA_N^{\Gl_N})$ (and an analog in the completed case) in the gauged double Poisson case, cf. Theorem \ref{Thm:gdP-rep1} which is motivated by \cite{AKKN}. 
We also treat the cases where one starts from a double quasi-Poisson bracket in Theorems \ref{Thm:dqP-rep1} and \ref{Thm:dqP-rep2},
or from a gauged double Poisson bracket in Theorem \ref{Thm:gdP-rep2}. 
This ends the first part of the memoir. 

\medskip 

%% chapt 8
In Chapter \ref{Ch:PVAcoh}, we deal with the commutative differential case and we recall  the definition of a Poisson vertex algebra (PVA) and several associated cohomology theories. 
Firstly, we introduce the space $\widetilde{\Gamma}(V)$ of basic cochains. This space can be equipped with a square-zero differential $\tilde{\delta}$ given by \eqref{eq:diff-PVA} (cf. Theorem \ref{thm:DSK1} due to De Sole and Kac \cite{DSK13}). This defines the basic PVA cohomology $\coH_{\mathrm{bas}}(V)$ of $V$. 
Moreover, we can get the reduced PVA cohomology  $\coH_{\mathrm{red}}(V)$ 
by inducing $\tilde{\delta}$ from $\widetilde{\Gamma}(V)$ to 
the complex $\widetilde{\Gamma}(V)/\partial\widetilde{\Gamma}(V)$ where $\partial$ is a natural extension  to $\widetilde{\Gamma}(V)$ of the derivation of $V$. 
Secondly, we consider the space $C(V)$ of $n$-$\lambda$-brackets which is also equipped with a differential, denoted $\dd$ \eqref{eq:dH-poly-n}. It is again a result of De Sole and Kac \cite{DSK13} that $\dd^2=0$, and the corresponding cohomology theory $\PVH(V)$ is called the variational PVA cohomology of $V$.
Finally, we discuss in Section \ref{sec:PVAcoh-inv} the induced cohomology in the presence of a group action by PVA automorphisms. 

%% chapt 9
In Chapter \ref{Ch:BasRed-dPVA}, we go to the noncommutative differential setting. 
We recall the notion of a double Poisson vertex algebra (dPVA) following \cite{DSKV}. 
Then, we define in Section \ref{sec:omega1} the space  $\widetilde{\Gamma}(\mc V)$ of basic cochains, where an element of degree $n$ is a map 
$\mc V^{\otimes n}\to \mc V^{\otimes (n+1)}[\lambda_1,\ldots,\lambda_n]$ satisfying  axioms \eqref{eq:sesquimaps} and \eqref{eq:Leibnizmaps} (cf. Section~\ref{sec:omega1}).  
The space $\widetilde{\Gamma}(\mc V)$ can be turned into a Lie superalgebra with a derivation $\partial$ extending the one on $\mc V$. 
This space admits a differential $\tilde{\delta}$ \eqref{eq:diff}, which squares to zero by Theorem \ref{thm:PVAcoh}. This allows us to introduce a first noncommutative dPVA cohomology 
which we call the \emph{basic dPVA cohomology} $\dH_{\mathrm{bas}}(\mc V)$ (Definition \ref{def:bas-dPVA}).
In Proposition \ref{20230802:prop1}, we show that $\tilde{\delta}$ preserves the (super) Lie bracket of $\widetilde{\Gamma}(\mc V)$ and it commutes with $\partial$. Hence it descends to a differential on
the quotient space
$$\Gamma(\mc V)=\widetilde{\Gamma}(\mc V)/([\widetilde{\Gamma}(\mc V),\widetilde{\Gamma}(\mc V)] + \partial \widetilde{\Gamma}(\mc V))$$
which squares to zero. In that way, we get the \emph{reduced dPVA cohomology} $\dH_{\mathrm{red}}(\mc V)$ (Definition \ref{def:bas-dPVA}). 
These two dPVA cohomologies are related through the long exact sequence \eqref{eq:les2}. 

%% chapt 10
In Chapter \ref{sec:compl-dPVA}, we consider the space $C(\mc V)$ spanned by the $n$-fold $\lambda$-brackets introduced in \cite{DSKV}. 
As in the commutative setting, we can define a linear map $\dd$ \eqref{eq:dP-1} on $C(\mc V)$ and we prove in Theorem \ref{Thm:g-dPVcoh1} that $\dd$ is a square-zero differential.
We call the cohomology of the complex $(C(\mc V),\dd)$ the \emph{variational dPVA cohomology} of $\mc V$ and we denote it by $\dPVH(\mc V)$, see Definition \ref{def:var-dPVA}.
Then, we characterize the first few cohomology groups in Section \ref{sec:explicit}
and we compare the formula \eqref{eq:dP-1} for the differential $\dd$ on a dPVA with the one of $\wdd$ for a dPA, see Proposition \ref{Thm:g-dPVcoh1-compare}.
We also compare equation \eqref{eq:dP-1} with Chemla's formula, see Remark \ref{rem:chemla}.
Finally, we relate the reduced and variational dPVA cohomologies. 
This relation is based on the projection operators $\widetilde{\proj}_n:\widetilde{\Gamma}^n(\mc V)\to C^n(\mc V)$, defined in \eqref{eq:Pn}, which induce a map $\proj:\Gamma(\mc V)\to C(\mc V)$. 
We show in Proposition \ref{20250720.prop1} that this is a morphism of complexes which induces a linear map $\dH_{\mathrm{red}}(\mc V)\to \dPVH(\mc V)$ in cohomology. 

%% chapt 11
In Chapter \ref{sec:PVAdiff}, we restrict our attention to the case of the algebra of noncommutative differential polynomials $\mc R_\ell$ in a finite number of generators $u_1,\dots,u_\ell$ for which we establish that the map $\dH_{\mathrm{red}}(\mc R_\ell)\to \dPVH(\mc R_\ell)$ is an isomorphism. 
To do so, we first remark that a basic $n$-cochain $X\in \widetilde{\Gamma}^n(\mc R_\ell)$ is equivalent to an array 
with entries ($i_1,\ldots,i_n\in I=\{1,\ldots,\ell\}$)
\[X_{i_1,\ldots,i_n}
:=X(u_{i_1}\otimes\dots\otimes u_{i_n})\in \mc R_\ell^{\otimes (n+1)}[\lambda_1,\ldots,\lambda_n]\,,
\]
obtained by evaluating $X$ on any $n$-uple of generators of $\mc R_\ell$. 
This gives an identification 
$\tilde{\gamma}:\widetilde{\Gamma}(\mc R_\ell)\stackrel{\simeq}{\longrightarrow} \widetilde{\Omega}(\mc R_\ell)$.
Furthermore, this identification descends to a map 
$\gamma:\Gamma(\mc R_\ell)\stackrel{\simeq}{\longrightarrow} \Omega(\mc R_\ell)$, 
where $\Omega(\mc R_\ell)$ is the quotient of $\widetilde{\Omega}(\mc R_\ell)$ by graded commutators and the image of the derivation $\partial$, see \eqref{eq:omega}. 
In this way, we can transfer the square-zero differential $\tilde{\delta}$ through $\widetilde{\gamma}$ and compute the basic and reduced PVA cohomologies using 
$\widetilde{\Omega}(\mc R_\ell)$ and $\Omega(\mc R_\ell)$ as in \eqref{20250723:eq2}. 
In a similar way, we can identify through \eqref{eq:psi} the space $C(\mc R_\ell)$ of $n$-fold $\lambda$-brackets with the space $\Sigma(\mc R_\ell)$ of skewsymmetric arrays from \cite{DSKV} and, by inducing the differential $d$ through this identification, we get an isomorphism $\dPVH(\mc R_\ell)\simeq \coH(\Sigma(\mc R_\ell),d)$. 
Building on the identification $\Omega(\mc R_\ell)\simeq \Sigma(\mc R_\ell)$ proved in \cite{DSKV}, we get the key isomorphism $\Gamma(\mc R_\ell)\simeq C(\mc R_\ell)$, which is proved in Proposition \ref{20240819:prop2} to be the projection operator $\proj$ introduced in the previous chapter. 
Since $\proj$ is the morphism of complexes inducing the linear map 
$\dH_{\mathrm{red}}(\mc R_\ell)\to \dPVH(\mc R_\ell)$, this last map is indeed an isomorphism for $\mc V = \mc R_\ell$ as we claimed.

%% chapt 12
In Chapter \ref{CH:PVAexa}, we compute all variational dPVA cohomology groups for a non-degenerate constant $2$-fold $\lambda$-bracket of degree $M\geq0$ (cf. \eqref{Eq:uu-odd-gen}) on $\mc R_\ell$, which is a noncommutative version of the ``generalized variational complex" in \cite{DSK13}. 
Thanks to Theorem \ref{Thm:DPVcoh-Const}, we obtain $\dH_{\mathrm{bas}}(\mc V)$ by determining the equivalent cohomology 
$\coH(\widetilde{\Omega}(\mc R_\ell),\tilde{\delta})$. 
The main technical ingredient is the construction of a homotopy operator $h_M$ \eqref{eq:hm} such that any element $Y\in \widetilde{\Omega}^n(\mc R_\ell)$ of `sufficiently large power in all $\lambda_1,\ldots,\lambda_n$' satisfying $\tilde{\delta}(Y)=0$ can be written as a coboundary $Y=\tilde{\delta}(h_M(Y))$, see Proposition \ref{Pr:Homot}. 
Then, we make use of the long exact sequence \eqref{eq:les2Omega} to characterize   
$\coH(\Omega(\mc R_\ell),\delta)\simeq \dPVH(\mc R_\ell)$ in  Theorem \ref{Thm:HOmega-Const}. In particular, we compute the dimension of the cohomology spaces $\coH^n(\Omega(\mc R_\ell),\delta)$, $n\geq0$, in Proposition \ref{cor:dim}. It follows from Proposition \ref{20250905:prop1} that these cohomology spaces are always nontrivial for $M\geq1$. This is a completely different behavior compared to the analogous result in the commutative case provided by Theorem 11.10 in \cite{DSK13} (see also \cite{CCS}).
As an application of the results of this chapter, we provide a description of the nontrivial first order deformations of the dPVA $\mc R_\ell$ in the case of $M=1$.

%% chapt 13
In Chapter \ref{Ch:repVardPVA}, we induce the variational dPVA cohomology to representation algebras in agreement with the Kontsevich-Rosenberg principle. 
For any $N\geq 1$, we start by building a linear map $\tr: C(\mc V)\to C(\mc V_N)$ assigning an $n$-$\lambda$-bracket on $\mc V_N$ to an $n$-fold $\lambda$-bracket on $\mc V$, see Lemmas  \ref{20240829:lem1} (for $n=0$), \ref{20240902:lem1} (for $n=1$) and Theorem \ref{Thm:IndBr-PVA} (for any $n$). 
Then, it is shown in Theorem \ref{Thm:dP-rep2-PVA} that the map $\tr$ is a morphism of complexes descending to a linear map $\dPVH(\mc V)\to \PVH(\mc V_N)$. 
The $\Gl_N$-invariance of the induced $n$-$\lambda$-bracket is obtained in Proposition \ref{20250627:prop1}, so that the previous map induces a linear map 
$\dPVH(\mc V)\to \PVH(\mc V_N^{\Gl_N})$.
This ends the second part of this memoir. 

\medskip 

%% chapt 14
In Chapter \ref{Ch:PA-PVA}, we prove explicitly the folklore result stating that that the jet and quotient functors induce linear maps in cohomology (in the commutative setting).  
Given a Poisson algebra $A$, its jet algebra $J_\infty A$ is a PVA and by extending skewsymmetric multilinear derivations we get a linear map 
$\mc X(A)\to C(J_\infty A)$, which descends in cohomology as 
$\coH_{CE}(A)\to \PVH(J_\infty A)$, cf. Corollary \ref{Cor:PAPVA1}. 
Similarly, given a PVA $V$, its quotient $q(V)=V/\langle \partial V\rangle$ is a Poisson algebra and we can induce a linear map $\PVH(V)\to \coH_{CE}(q(V))$, see Corollary \ref{Cor:PAPVA2}.  

%% chapt 15
In Chapter \ref{Ch:dPA-dPVA}, we adapt the considerations from the previous chapter to the noncommutative setting based on the (associative) jet and quotient functors from \cite{BFM}. 
Namely, in Corollary \ref{Cor:dPAdPVA1} we obtain a map 
$\widehat{\dPH}(\cA)\to \dPVH(J_\infty\cA)$ starting from a dPA $\cA$, while 
in Corollary \ref{Cor:dPAdPVA2} we get a map 
$\dPVH(\mc V)\to \widehat{\dPH}(q(\mc V))$. 
We show in Corollaries \ref{20250822:cor1} and \ref{20250822:cor2} that one obtains commutative diagrams if we consider these maps in cohomology together with their commutative analogs associated with the representation algebras $A=\cA_N$ or $V=\mc V_N$. 

The final two chapters relate the preceding two parts of this memoir and they end the core of the text.
The reader can find after this an (hopefully helpful) Index of Notation and the Bibliography. 

\subsection*{Acknowledgments}

We are grateful to Anton Alekseev and Sophie Chemla for useful exchanges about the papers \cite{AKKN,Ch}, and to Yong Zheng for sharing the unpublished works \cite{ZT,ZT2}. We also thank Ralph Kaufmann and Sasha Tsymbaliuk for stimulating conversations.
M.F. thanks La Sapienza for hospitality.  
We thank IHES, where part of this work was written, for the excellent working conditions.

M.F. acknowledges financial support from CNRS (projet PEPS JCJC 2024), and from EIPHI Graduate School (contract ANR-17-EURE-0002) for funding IMB. 
D.V. acknowledges the financial support of the project MMNLP (Mathematical Methods in Non Linear Physics) of INFN and of the project PRIN 2022HMBTTL of MUR.

%\pecetta{Should we add the symbols (e.g. $\bullet_{i,i+1}$, $\ast_i$, $\dgal{-,-}$,etc.) to the glossary ???
%
%Daniele: good idea, so far I added only the symbols in the Preliminary Section (feel free to modify or rearrange them). Should we add symbols for all sections? (maybe just the most important ones)
%
%Maxime: yes, let's only put the `most important'}

%%%%%%%%%%% NEW CHAPTER %%%%%%%%%%%%%%%
%%%%%%%%%%% NEW CHAPTER %%%%%%%%%%%%%%%
%%%%%%%%%%% NEW CHAPTER %%%%%%%%%%%%%%%
%%%%%%%%%%% NEW CHAPTER %%%%%%%%%%%%%%%

\chapter{Set up and preliminary results} 
\label{Ch:Setup}

Throughout the book, we fix a field $\kk$, assumed to be of characteristic zero and algebraically closed for simplicity. Unadorned tensor products are over $\kk$, i.e. $\otimes := \otimes_\kk$. A linear map is a $\kk$-linear map.

By \emph{algebra}, we mean an associative unital $\kk$-algebra. 
The multiplication $\mult:\cA\otimes \cA \to \cA$ of an algebra $\cA$ is usually denoted by concatenation:  $ab=\mult(a,b)$ for $a,b\in \cA$. 
We denote by $\Id_\cA$, or simply $\Id$, the identity morphism
on $\cA$. 
We will explicitly state when an algebra is required to be commutative or finitely generated.

Given an algebra $\cA$ and a subalgebra $B\subset \cA$, it will be useful to view $\cA$ as a $B$-algebra. 
For example, if 
$\cA$ admits a complete finite set of orthogonal idempotents $(e_s)_{s\in S}$, i.e. $e_s e_t=\delta_{st} e_s$ and $1=\sum_{s\in S} e_s$, we view $\cA$ as an algebra over $B=\oplus_{s\in S} \kk e_s$\glslink{Balg}{}. 
(We used Kronecker's delta function.)

A \emph{differential algebra} is a pair $(\mc V,\del)$ 
where $\mc V$ is an algebra equipped with a $\kk$-linear map $\del:\mc V\to \mc V$ satisfying the Leibniz rule $\del(ab)=a\del(b)+\del(a)b$ ($a,b\in \mc V$).  
We say that a subset $\{a_i\}_{i\in I}\subset \mc V$ generates $\mc V$ as a differential algebra if
any $b\in \mc V$ can be written (non-uniquely) as a finite sum
\begin{equation*}
    b=\sum_{n\in\mb Z_{\geq 0}} \sum_{\underline i\in I^n}
    \sum_{m_1,\ldots,m_n\in\mb Z_{\geq 0}}
    \gamma_{\underline i}\, \del^{m_1}(a_{i_1}) \cdots \del^{m_n}(a_{i_n}), \quad 
    \underline i=(i_1,\ldots,i_n)\in I^n, \,\, \gamma_{\underline i}\in \kk\,.
\end{equation*}
If the cardinality of $I$ is finite, we say that $\mc V$ is a finitely generated differential algebra.
We will explicitly state when an algebra is required to be commutative.

\begin{convention}
We shall denote by $\cA$ an algebra, and by $A$ a commutative algebra. 
Similarly, we shall write $\mc V$ for a differential algebra, and $V$ for a commutative differential algebra. 
We keep $B$ for a (base) subalgebra. 
\end{convention}

In this chapter, we gather all the elementary operations and notations used in the rest of the text. 
The reader can consult \cite{DSKV,VdB1} for an analogous presentation. 

%%%
\section{Operations from associative algebras}\label{sec:1.1}

\subsection{}
Let $\mc A$ be an algebra. We endow $\cA^{\otimes 2}:=\cA\otimes \cA$ with the natural associative product 
\begin{equation*}
    (a'\otimes a'')\, (b'\otimes b'') = a'b' \otimes a''b'', 
\end{equation*}
where $a',a'',b',b''\in \cA$. 
(We shall freely use Sweedler's notation $c'\otimes c'':=C\in \cA^{\otimes 2}$ throughout the text.)  
We also consider the \emph{bullet product} which is the associative product defined by the formula  
\glslink{bullet}{}
\begin{equation}\label{20140609:eqc1}
 (a'\otimes a'')\bullet (b'\otimes b'') = a'b' \otimes b''a''\,.
\end{equation}
(It is the product in $\cA\otimes \cA^{\mathrm{op}}$). The following are compatible $\mc A$-bimodule structures on $\cA^{\otimes 2}$,  ($a,b\in \cA$, $C\in \cA^{\otimes 2}$)
\begin{subequations}
\begin{align}
  a \,C\, b &:= (a\otimes 1) \,C\, (1\otimes b)
 =ac'\otimes c'' b \,, &&  \text{(outer bimodule)}
 \label{bmodouter}\\
 a\star C \star b &:= (1\otimes a) \,C\, (b\otimes 1)
 = c'b \otimes a c''\,.  && \text{(inner bimodule)}
 \label{bmodinner}
\glslink{inner}{}\end{align}
\end{subequations}
Given $n\geq 2$, the $n$-fold tensor product $\cA^{\otimes n}$ is an algebra for  
\begin{equation*}
  (a_1\otimes a_2 \otimes \ldots \otimes a_n ) \, 
  (b_1\otimes b_2 \otimes \ldots \otimes b_n ) 
  = a_1 b_1\otimes a_2 b_2 \otimes \ldots \otimes a_n b_n , 
\end{equation*}
with all $a_i,b_i\in \cA$. 
For any $1\leq i \leq n-1$, we have a right action (with respect to the bullet product \eqref{20140609:eqc1}) of $\mc A\mc \otimes \mc A$ on $\mc A^{\otimes n}$ obtained by computing the bullet product in position $(i,i+1)$, i.e. ($C\in\mc A\otimes\mc A$) 
\begin{equation} \label{bullet-ii}\glslink{bullet-i}{}
  (a_1\otimes a_2 \otimes \ldots \otimes a_n ) \bullet_{(i,i+1)} C
  = a_1\otimes \ldots \otimes a_i c' \otimes c'' a_{i+1} \otimes \ldots \otimes a_n  .
\end{equation} 
It can be checked that the right actions $\bullet_{(i,i+1)}$ and $\bullet_{(j,j+1)}$ commute if $i\neq j$. 
We also have a left action (with respect to the bullet product \eqref{20140609:eqc1}) of $\mc A\mc \otimes \mc A$ on $\mc A^{\otimes n}$ given by
\begin{equation}\label{eq:bullyaction}\glslink{bullet}{}
C\bullet (a_1\otimes a_2\otimes \dots\otimes a_{n-1}\otimes a_n)
=c' a_1\otimes a_2\otimes \dots\otimes a_{n-1}\otimes a_n c''
\,,
\end{equation}
which reduces to \eqref{20140609:eqc1} for $n=2$. Note that the left action \eqref{eq:bullyaction} commutes with the right actions \eqref{bullet-ii} for every $i=1,\dots,n-1$.

We  introduce for any $0\leq j\leq n-1$ the following left and right $\cA$-module structures on $\mc A^{\otimes n}$:
\begin{equation}
\begin{split}\label{eq:star}\glslink{inner-j}{}
b\star_j (a_1\otimes \ldots \otimes a_n) &= 
a_1\otimes \ldots \otimes a_j\otimes  b a_{j+1}\otimes a_{j+2}\otimes \ldots \otimes a_n  \\  
 (a_1\otimes \ldots \otimes a_n)\star_j b&=
a_1\otimes \ldots \otimes a_{n-j-1}\otimes  a_{n-j}b\otimes a_{n-j+1}\otimes \ldots \otimes a_n \,.   
\end{split}
\end{equation}
The operation $\star_n:=\star_0$ is the outer bimodule multiplication which acts by multiplying the leftmost factor on the left, and the rightmost factor on the right, respectively. 
We may simply write $\star_0$ as concatenation, and consider its generalization for any $m,n \geq2$ as the linear map $\cA^{\otimes m}\times \cA^{\otimes n}\to \cA^{\otimes (m+n-1)}$ given by 
\begin{align} \label{prod-nm}
(a_1\otimes  \ldots \otimes a_m) \cdot (b_1 \otimes \ldots \otimes b_n)
=a_1\otimes \ldots \otimes a_{m-1}\otimes  a_m b_1\otimes b_2\otimes \ldots \otimes b_n\,,
\end{align}
where all $a_i,b_j\in \cA$.

For any $n\geq 2$, an element $\tau \in S_n$ of the symmetric group acts on $\cA^{\otimes n}$ by tensor permutation, denoted $\sigma_\tau:\cA^{\otimes n}\to \cA^{\otimes n}$. 
Explicitly, one has\glslink{sigmaction}{}
\begin{equation*}
    \sigma_\tau(a_1\otimes a_2\otimes  \ldots \otimes a_n)
    =a_{\tau^{-1}(1)}\otimes a_{\tau^{-1}(2)}\otimes  \ldots \otimes a_{\tau^{-1}(n)}\,.
\end{equation*}
In many constructions, it will be important to act by the cyclic permutation $(12 \cdots n)$, and we shall use the shorthand notation $\sigma := \sigma_{(12\cdots n)}$. 
In such a case, we may indicate $\sigma$ as an exponent, e.g. 
$C^\sigma = c'' \otimes c'$, for $C\in \cA^{\otimes 2}$. 
We also introduce, for any $C\in\mc \cA^{\otimes m}$ and $1\leq i \leq n$ the map 
$-\otimes_i C : \cA^{\otimes n}\to \cA^{\otimes (n+m)}$ defined by 
\begin{equation}\label{eq:tensor-i-notation}\glslink{otimes-i}{}
(a_1\otimes \dots\otimes a_n)\otimes_{i}C
=a_1\otimes\dots\otimes a_{n-i}\otimes C\otimes a_{n-i+1}\otimes \dots\otimes a_n\,.
\end{equation} 
The operation in \eqref{eq:tensor-i-notation} is not associative but it satisfies a ``sort" of associativity condition stated in the next result. 
\begin{lemma}\label{lem:PVAcoh-pre6}
Let $X\in\mc A^{\otimes n}$, $Y\in\mc A^{\otimes m}$ and $Z\in\otimes \mc A^{\otimes l}$.
For every $0\leq i\leq n$ and $0\leq j\leq m$, we have 
\begin{equation}\label{20230811:eq1a}
(X\otimes_i Y)\otimes_{i+j} Z=X\otimes_i(Y\otimes_j Z)\,.
\end{equation}
\end{lemma}
\begin{proof}
Straightforward.
\end{proof}
We also describe the behavior of the operation \eqref{eq:tensor-i-notation} with respect to the left and right $\mc A$-module structures \eqref{eq:star}.
\begin{lemma}\label{lem:PVAcoh-pre4}
Let $a\in\mc A$, $X\in\mc A^{\otimes n}$ and $Y\in\mc A^{\otimes m}$. For every $i=0,\dots n-1$ we have
$$
\left(a\star_i X\right)\otimes_j Y
=\left\{
\begin{array}{ll}
a\star_i\left(X\otimes_j Y\right)\,,
& 0\leq j\leq n-1-i
\,,
\\
a\star_{m+i}\left(X\otimes_j Y\right)\,,
& n-i\leq j\leq n
\,,
\end{array}
\right.
$$
and
$$
\left(X\star_i a\right)\otimes_j Y
=\left\{
\begin{array}{ll}
\left(X\otimes_j Y\right)\star_{m+i}a\,,
& 0\leq j\leq i
\,,
\\
\left(X\otimes_j Y\right)\star_{i}a\,,
& i+1\leq j\leq n
\,.
\end{array}
\right.
$$
\end{lemma}
\begin{proof}
This is a straightforward computation using
\eqref{eq:star} and~\eqref{eq:tensor-i-notation}.
\end{proof}
\subsection{ }\label{sec:1.1.2}
Let $D:\cA \to \cA^{\otimes m}$, $m\geq 1$, be an arbitrary linear map. 
For any $n\geq 2$ and $1\leq i \leq n$, we have an extended map  $D_{(i)}:\cA^{\otimes n}\to \cA^{\otimes (m+n-1)}$ obtained by action on the $i$-th factor: 
\begin{equation} \label{20240805:eq1}\glslink{D-i}{}
   D_{(i)}(a_1\otimes  \ldots \otimes a_n) 
   = a_1\otimes \ldots \otimes D(a_i)\otimes \ldots \otimes a_n\,.
\end{equation}
For the action on the left- and right-most factor, we use 
\begin{equation}\label{20240805:eq1b}\glslink{D-LR}{}
    D_L:=D_{(1)}, \qquad D_R:= D_{(n)}\,. 
\end{equation}
The following easily follows from the notation, cf. \cite[Lem.~1.4]{DSKV}.
\begin{lemma}\label{20140606:lem}
Let $D:\, \mc A\to \mc A^{\otimes n}$ be a linear map.
For $C\in \mc A^{\otimes m}$, we have
\begin{equation}\label{20140606:eq4}
(D_{(i)}C)^\sigma=D_{(i+1)}(C^\sigma)
\,\,\text{ for } 1\leq i\leq m-1
\,; \quad  
(D_{(m)}C)^{\sigma^n}=D_{(1)}(C^\sigma)
\,.
\end{equation}
\end{lemma}
A linear map $D:\cA \to \cA^{\otimes m}$ is an \emph{$m$-fold derivation} if it satisfies 
\begin{equation*}
    D(ab)=a\,D(b) + D(a)\, b, \qquad a,b\in \cA , 
\end{equation*}
where we recall that concatenation means the outer bimodule multiplication $\star_0$ \eqref{eq:star}. 
In that case, it is natural to extend $D$ to $D:\cA^{\otimes n}\to \cA^{\otimes (m+n-1)}$ as 
\begin{equation} \label{mfold-ext}
  D :=\sum_{i=1}^n D_{(i)}, \quad  a_1\otimes  \ldots \otimes a_n \mapsto 
   \sum_{i=1}^n a_1\otimes \ldots \otimes D(a_i)\otimes \ldots \otimes a_n\,.
\end{equation} 
If $D'$ is an $M$-fold derivation, we can consider the composition $D\circ D'$ using the extension \eqref{mfold-ext}. Then, the commutator $[D,D']$ is itself a $(m+M-1)$-fold derivation.  
For $m=2$, we talk about \emph{double derivations} and get \glslink{DDer}{}
\begin{equation} \label{DDer}
    \begin{aligned}
        \DDer(\cA):=&\Der(\cA,(\cA^{\otimes 2})_{out}) \\
=&\{D \in \Hom_\kk(\cA,\cA^{\otimes 2}) \mid D(ab)=a\,D(b)+D(a)\,b \}\,.
    \end{aligned}
\end{equation} 
Using the inner bimodule structure on $\cA^{\otimes 2}$ \eqref{bmodinner}, 
we can naturally see $\DDer(\cA)$ as an $\cA$-bimodule. 
If $\cA$ is a $B$-algebra, we require $m$-fold derivations to be $B$-linear, i.e. $D(b)=0$ for all $b\in B$. 
Then, the subgroup $\DDer_B(\cA)=\{D \in \DDer(\cA)\mid D(B)=0\}$ of double derivations relative to $B$ is also an $\cA$-bimodule. 

We write $\cA_\sharp:=\cA/[\cA,\cA]$ 
\glslink{Asharp}{}
for the vector space obtained by identifying all commutators to zero in $\cA$. 
%The quotient map $\tr:\cA\to \cA_\sharp$ is called the trace map. 
We let $\bar{a}%=\tr(a)
\in \cA_\sharp$ for the image of $a\in \cA$ and, when no confusion can arise, we may denote it by $a\in \cA_\sharp$ as well. 
Given an $m$-fold derivation $D$ we get a derivation $\mult\circ D\in\Der(\mc A)$, where 
\begin{equation*}\glslink{mult}{}
\mult:\cA^{\otimes m}\to \cA, \quad \mult(a_1\otimes  \ldots \otimes a_m)
=a_1\cdots a_m\,.
\end{equation*}
%The induced map is a derivation, i.e. $\mult \circ D\in \Der(\cA)$. 
(If we start with a $B$-linear $m$-fold derivation $D$, then $\mult \circ D$ is a $B$-linear derivation.) 
In turn, we obtain an induced linear map $D_\sharp : \cA_\sharp\to \cA_\sharp$.

\subsection{ }
For $m\geq2$ and $i=1,\dots m-1$ we also extend the multiplication map $\mult:\mc A\otimes\mc A\to\mc A$ to a linear map $\mult_{(i,i+1)}:\mc A^{\otimes m}\to\mc A^{\otimes (m-1)}$ by setting
$\mult_{(i,i+1)}=\Id^{\otimes (i-1)}\otimes\mult\otimes \Id^{\otimes (m-1-i)}$.
Namely, we have
\begin{equation}\label{eq:mii+1}\glslink{mult-i}
\mult_{(i,i+1)}(a_1\otimes \dots \otimes a_m)=a_1\otimes a_{i-1}\otimes a_ia_{i+1}\otimes a_{i+2}\dots a_m
\,.
\end{equation}
For $X\in\mc A^{\otimes m}$ and $Y\in\mc A^{\otimes n}$ the product \eqref{prod-nm} can be then rewritten as
$$
XY=\mult_{(m,m+1)}(X\otimes Y)
\,.
$$
The following result expresses some compatibility properties between the multiplication maps
\eqref{eq:mii+1} and the action of the cyclic permutation $\sigma$, the left and right
$\mc A$-module structures \eqref{eq:star} and the products \eqref{prod-nm}
and \eqref{eq:tensor-i-notation}
that will be needed throughout the memoir.
\begin{lemma}\label{20240724:lem1}
\begin{enumerate}[(a)]
\item Let $n\geq1$. For every $k=0,\dots,n-1$, we have the following identity in
$\Hom(\mc A^{\otimes (n+1)},\mc A^{\otimes n})$:
\begin{equation}\label{20240724:eq1}
\sigma^k\circ \mult_{(h,h+1)}=
\left\{
\begin{array}{ll}
\mult_{(h+k,h+k+1)}\circ\sigma^k\,,&
h=1,\dots,n-k\,,
\\
\mult_{(h+k-n,h+k+1-n)}\circ\sigma^{k+1}\,,&
h=n-k+1,\dots,n
\,.
\end{array}
\right.
\end{equation}
\item Let $n\geq1$ and $i\in\{1,\dots,n\}$. For every $X\in\mc V^{\otimes (n+1)}$, $a\in\mc A$ and $h=0,\dots,n-1$
we have
\begin{equation}\label{20240724:eq4}
a\star_i(\mult_{(h+1,h+2)}\circ \sigma^{h+1}(X))
=\mult_{(h+1,h+2)}\circ \sigma^{h+1}(a\star_{\sigma^{-h}(i)}X)
\end{equation}
and
\begin{equation}\label{20240724:eq5}
(\mult_{(h+1,h+2)}\circ \sigma^{h+1}(X)){\star_{n-i}} a
=\mult_{(h+1,h+2)}\circ \sigma^{h+1}(X\star_{n+1-\sigma^{-h}(i)}a)
\,.
\end{equation}
\item Let $(m,n)\in\mb Z_{\geq0}^2\setminus\{(0,0)\}$. For every $X\in \mc A^{\otimes (m+1)},Y\in\mc A^{\otimes (n+1)}$ we have
\begin{equation}\label{20240724:eq2}
\begin{split}
&\mult_{(h+1,h+2)}\circ\sigma^{h+1}(XY)
\\
&=
\left\{
\begin{array}{ll}
\mult_{(m+h+1,m+h+2)}\circ\sigma^{m+h+1}(YX)\,,&
h=0,\dots,n-1\,,
\\
\mult_{(h+1-n,h+2-n)}\circ\sigma^{h+1-n}(YX)\,,&
h=n,\dots,n+m-1
\,.
\end{array}
\right.
\end{split}
\end{equation}
\item
For every $X\in\mc A^{\otimes (n+1)}$, $a\in\mc V$ and $i=1,\dots,n$, we have
\begin{equation}\label{20240823:eq4}
\begin{split}
&(\mult_{(h+1,h+2)}\sigma^{h+1}(X))\otimes_{n-i}a
\\
&=\left\{
\begin{array}{ll}
\mult_{(h+1,h+2)}\sigma^{h+1}\left(X\otimes_{n+1-(i-h)}a\right)\,,
&0\leq h\leq i-1\,,
\\
\mult_{(h+2,h+3)}\sigma^{h+2}\left(X\otimes_{n+1-(n+i-h)}a\right)
&t\leq h\leq n-1\,.
\end{array}\right.
\end{split}
\end{equation}
\end{enumerate}
\end{lemma}
\begin{proof}
Equation \eqref{20240724:eq1} for $k=1$ follows from definitions similarly to the proof of the claim of Lemma \ref{20140606:lem}. For $k>1$ it can be easily proved by induction. This proves part (a).
For part (b) we emphasize that $\sigma$ acts on $i\in \{1,\ldots,n\}$ as the cyclic permutation $(1\dots n)$, so that $\sigma^{-h}(i)=i-h$ for $h+1\leq i \leq n$ while 
$\sigma^{-h}(i)=n+i-h$ for $1\leq i \leq h$. Thus, by a straightforward computation one can check that both sides of \eqref{20240724:eq5} applied to $x_1\otimes\dots \otimes x_{n+1}\in\mc A^{\otimes (n+1)}$ give
\begin{equation*}
    \left\{ 
\begin{array}{ll}
x_{n-h+1}\otimes\dots\otimes x_{n-h+i}a\otimes \dots\otimes x_{n+1}x_1\otimes \dots \otimes x_{n-h}   &  1\leq i \leq h\,, \\
x_{n-h+1}\otimes\dots\otimes x_{n+1}x_1 a\otimes\dots\otimes x_{n-h}   &  i = h+1 \,, \\
x_{n-h+1}\otimes\dots\otimes x_{n+1}x_1\otimes\dots\otimes x_{i-h}a\otimes\dots\otimes x_{n-h})   &  h+1< i \leq n \,. 
\end{array}
\right.
\end{equation*}
Similar (actually simpler) verification for \eqref{20240724:eq4}. 
We are left to prove part (c). Let $X=x_1\otimes \ldots \otimes x_{m+1}\in\mc A^{\otimes (m+1)}$ and
$Y=y_1\otimes \ldots \otimes y_{n+1}\in\mc A^{\otimes (n+1)}$. Then, it is immediate to verify that
\begin{equation}\label{20240823:eqlem}
\mult_{(1,2)}\sigma(XY)=y_{n+1}x_1\otimes\dots\otimes x_{m+1}y_1\otimes\dots\otimes y_n
=\mult_{(m+1,m+2)}\sigma^{m+1}(YX)
\,,
\end{equation}
which gives \eqref{20240724:eq2} for $h=0$. The general case follows by applying $\sigma^h$, $h=1,\dots,n+m+1$, to both sides of
\eqref{20240823:eqlem} and using equation \eqref{20240724:eq1} of part (a). The proof of part (d) is straightforward.
\end{proof}
In the same way, the next result expresses some compatibility properties between the multiplication maps
\eqref{eq:mii+1} and the action of a double derivation by extension \eqref{20240805:eq1}.
\begin{lemma}\label{Lem:Dmult}
Fix a double derivation $D :\cA \to \cA^{\otimes 2}$, $a\in \cA$ and $X\in \cA^{\otimes (n+1)}$ with $n\geq 1$. For any $1\leq t \leq n$, one has: 
\begin{enumerate}[(1)]
\item For any $1\leq s \leq t-1$, 
\begin{equation}
\label{Eq:Dmult1}
   D_{(t)}(\mult_{(s,s+1)}(a\star_sX))=\mult_{(s,s+1)}(a\star_s D_{(t+1)}(X))
\,.
\end{equation}
\item For $s=t$, 
\begin{equation}
  \begin{aligned}
\label{Eq:Dmult2}
&D_{(t)}(\mult_{(t,t+1)}a\star_t X)=\mult_{(t+1,t+2)}(a\star_{t+1}D_{(t)}(X))
\\
%&\qquad\qquad +x_1\otimes \dots \otimes x_tD(a)'\otimes %D(a)''x_{t+1}\otimes\dots\otimes x_{n+1} \\
&+X\bullet_{(t,t+1)}D(a)+\mult_{(t,t+1)}(a\star_{t}D_{(t+1)}(X))
\,.
\end{aligned}  
\end{equation}
\item For any $t+1\leq s \leq n$, 
\begin{equation}
    \label{Eq:Dmult3}
D_{(t)}(\mult_{(s,s+1)}(a\star_sX))=\mult_{(s+1,s+2)}(a\star_{s+1} D_{(t)}(X))\,.
\end{equation}
\end{enumerate}
Moreover, we have
\begin{equation}
\begin{split}\label{Eq:Dmult5}
&\sigma^{-1}\,D_{L}(\mult_{(1,2)}a\star_1 X)=\mult_{(1,2)}(a\star_{1} \sigma^{-1} D_{L}(X))
%\nonumber 
\\
& 
%\qquad\qquad +\sigma^{-1}(x_1 D(a)'\otimes D(a)''x_{2}\otimes x_3\otimes \dots\otimes x_{n+1})  \\
%&\qquad\qquad
+\sigma^{-1}(X\bullet_{(1,2)}D(a))+\mult_{(n+1,n+2)}(a\star_{n+1} \sigma^{-1} D_{L}\sigma^{-1}(X)) \,,
%\nonumber 
\end{split}
\end{equation}
and, for every $2\leq s\leq n$,
\begin{align}
    \label{Eq:Dmult4}
&\sigma^{-1}\, D_{L}(\mult_{(s,s+1)}(a\star_sX))=\mult_{(s,s+1)}(a\star_{s}\sigma^{-1} D_{L}(X))\,.
\end{align}
\end{lemma}
\begin{proof}
Parts (1), (2) and (3) are obtained by direct calculations. Equation \eqref{Eq:Dmult5} (respectively \eqref{Eq:Dmult4}) follows by letting $t=1$ in part (2) (respectively part (3)) and applying
$\sigma^{-1}$ to both sides.
\end{proof}
\subsection{ }
If $\cA$ has an additional $\Z$-grading, all the previous notations can be considered. 
However, one needs to be careful of the sign rule of Loday for multiplying tensors, e.g. 
\begin{align} \label{gradmult}
   (a' \otimes a'') (b'\otimes b'')=(-1)^{|a''|\,|b'|} a'b'\otimes a''b''\,,
\end{align}
and in particular there is a sign $(-1)^{|b||d''|+|a||d'|}$ in \eqref{bmodinner}. One also has 
\begin{align*}
   (a_1\otimes a_2\otimes \ldots \otimes a_n)^\sigma = 
  (-1)^{\rho}\,  (a_2\otimes \ldots \otimes a_n \otimes a_1), \quad 
  \rho:=|a_1| \sum_{j=2}^n |a_n|.
\end{align*}
Moreover, the notation 
$\cA_\sharp$ stands for the quotient of $\cA$ by graded commutators  
\begin{equation} \label{Eq:GradComm}
[\cA,\cA]=\operatorname{span}_\kk\{ab-(-1)^{|a||b|}ba \mid a,b\in \cA\}\,.   
\end{equation}

%%%
\section{Operations from differential algebras}\label{sec:1.2}
For a differential algebra $(\mc V,\del)$, we can complement the operations defined in the previous subsections by new ones involving the derivation $\partial$. 
Hereafter, given an algebra $\cA$ (such as $\cA=\mc V^{\otimes n}$), we let 
$\cA[\lambda_1,\ldots,\lambda_j]:=\cA\otimes \kk[\lambda_1,\ldots,\lambda_j]$whose elements are written as polynomials 
$\sum_{\underline{k}} a_{\underline{k}} \, \lambda_1^{k_1}\cdots \lambda_j^{k_j}$, 
with $\underline{k}=(k_1,\ldots,k_j)\in \Z_{\geq 0}^j$ and finitely many nonzero $a_{\underline{k}} \in \cA$.  
Given $a(\lambda)=\sum_ia_i\lambda^i\in \mc V^{\otimes l}[\lambda]$
and $b\in\mc V^{\otimes m},c\in \mc V^{\otimes n}$, we let 
\begin{equation}\label{eq:notation}
\begin{split}
&a(\lambda+x)\big(\big|_{x=\partial}b)c
=\sum_i a_i(\lambda+\partial)^i(b)\,c
\\
&=\sum_i \sum_{t=0}^i \binom{i}{t} a_i\,\partial^t(b)\,c\,\lambda^{i-t}\, \in \mc V^{\otimes (l+m+n-2)}[\lambda]
\,.
\end{split}
\end{equation}
Fix $n\geq 2$ and a polynomial 
\begin{equation}
    p(\lambda_1,\dots,\lambda_{n-1})
=\sum_{\underline{k}} a_1\otimes\ldots\otimes a_n \,\, 
\lambda_1^{k_1}\cdots\lambda_{n-1}^{k_{n-1}}
\in \mc V^{\otimes n}[\lambda_1,\dots,\lambda_{n-1}]
\,.
\end{equation}
We consider for any $f\in \mc V$ the following extensions of \eqref{eq:star}:  
\begin{equation} \label{20140707:eq6-left}
   \begin{array}{l}
\displaystyle{
\vphantom{\Big)}
(|_{x=\partial}f)\star_{j}
p(\lambda_1,\dots,\lambda_j+x,\dots,\lambda_{n-1})
} \\
\displaystyle{
\vphantom{\Big)}
=\sum_{\underline{k}} \sum_{t=0}^{k_j}\binom{k_j}{t}
a_1\otimes\ldots  \otimes(\partial^t f)a_{j+1}\otimes\ldots\otimes a_n \, 
\lambda_1^{k_1}\cdots \lambda_j^{k_j-t}\cdots\lambda_{n-1}^{k_{n-1}}
,}
\end{array} 
\end{equation} 
and
\begin{equation}\label{20140707:eq6}
\begin{array}{l}
\displaystyle{
\vphantom{\Big)}
p(\lambda_1,\dots,\lambda_i+x,\dots,\lambda_{n-1}) \star_{n-j}
(|_{x=\partial}f)
} \\
\displaystyle{
\vphantom{\Big)}
=\sum_{\underline{k}} \sum_{t=0}^{k_j}\binom{k_j}{t}
a_1\otimes\ldots \otimes a_j(\partial^t f)\otimes\ldots\otimes a_n\, 
\lambda_1^{k_1}\cdots \lambda_j^{k_j-t}\cdots\lambda_{n-1}^{k_{n-1}}
.}
\end{array}
\end{equation} 
There, one replaces
$\lambda_i+x$ by $\lambda_i+\partial$, and the parentheses mean that $\partial$ should be applied to $f$. 
Similarly, for $C\in \mc V^{\otimes m}$ we can adapt the tensor notation \eqref{eq:tensor-i-notation} as 
\begin{equation}\label{20140707:eq6-tens}
\begin{aligned}
    &p(\lambda_1,\dots,\lambda_i+x,\dots,\lambda_{n-1}) \star_{n-j}
(|_{x=\partial} C)
\\
&=\sum_{\underline{k}} \sum_{t=0}^{k_j}\binom{k_j}{t}
a_1\otimes\ldots a_j\otimes \partial^t(C)\otimes a_{j+1} \ldots\otimes a_n\, 
\lambda_1^{k_1}\cdots \lambda_j^{k_j-t}\cdots\lambda_{n-1}^{k_{n-1}}
, 
\end{aligned}
\end{equation}
where $\partial$ acts on $C$ through the extension \eqref{mfold-ext}.  

\medskip 

For a collection of ordered elements $\alpha_1,\dots,\alpha_{n+1}$, introduce 
the notation 
\begin{equation}
    \begin{aligned} \label{Eq:NotCheck}
\alpha_1,\stackrel{s}{\check{\dots}},\alpha_{n+1}
&= \alpha_1,\dots,\alpha_{s-1},\alpha_{s+1},\dots,\alpha_{n+1} \,, \\ 
\alpha_1,\stackrel{s}{\check{\dots}},\stackrel{t}{\check{\dots}},\alpha_{n+1} 
&= \alpha_1,\dots,\alpha_{s-1},\alpha_{s+1},\dots, 
\alpha_{t-1},\alpha_{t+1},\dots,\alpha_{n+1}, 
\end{aligned}
\end{equation}
where $1\leq s\leq n+1$ and $s<t\leq n+1$.  
Namely we remove from the collection the element in position $s$, and the elements in positions $s,t$, respectively.  
For example, given a linear map  
\begin{equation} \label{20250627-IHES}
\begin{aligned}
&X_{\lambda_1,\dots,\lambda_{n}}:\mc V^{\otimes n}\to 
\mc V^{\otimes (n+1)}[\lambda_1,\dots,\lambda_{n}]
\,,
\\
&a_1\otimes\dots\otimes a_{n}\mapsto X_{\lambda_1,\dots,\lambda_{n}}(a_1,\dots,a_{n})
\end{aligned}
\end{equation}
and elements $b_1,\ldots,b_{n+1}\in \mc V$, 
we can consider 
\begin{equation}\label{Eq:NotCheck1}
\begin{split}
&X_{\lambda_1,\stackrel{s}{\check{\dots}},\lambda_{n+1}}
(b_1,\stackrel{s}{\check{\dots}},b_{n+1})
\\
&= X_{\lambda_1,\dots,\lambda_{s-1},\lambda_{s+1},\dots,\lambda_{n+1}}
(b_1,\dots,b_{s-1},b_{s+1},\dots,b_{n+1})\,.     
\end{split}
\end{equation}
Given a map $X_{\lambda_1,\dots,\lambda_{n}}$ as in \eqref{20250627-IHES}, 
we introduce ($a_i\in\mc V$, $C\in\mc V^{\otimes m}$) 
\begin{equation}
\begin{split}\label{eq:XL}\glslink{X-s}{}
&X_{\lambda_1,\dots,\lambda_n}^{(s)}(a_1,\dots,a_{s-1},a_s\otimes C,a_{s+1},\dots,a_{n})
\\
&=
X_{\lambda_1,\dots,\lambda_s+x,\dots,\lambda_n}(a_1,\dots,a_{s-1},a_s,a_{s+1},\dots,a_{n})
\otimes_{n+1-s}(|_{x=\partial}C)
\,,
\end{split}
\end{equation}
with the notation \eqref{eq:tensor-i-notation}. 

We shall consider the vector space $\mc V_\sharp=\mc V/(\partial \mc V+[\mc V,\mc V])$ 
\glslink{Vsharp}{}
and denote by $\tint:\mc V\to\mc V_{\sharp}$ \glslink{tint}{}the quotient map, which is the differential algebra analogue of the quotient map $\mc A\to\mc A_{\sharp}$ introduced in Section \ref{sec:1.1.2}. 
We warn the reader that, for a differential algebra $\mc V$, we shall not consider the vector space $\mc V/[\mc V,\mc V]$ so no confusion can arise from the notation. 
If $\mc V$ is equipped with a $\Z$-grading, we use the same definition for $\mc V_\sharp$ with the vector space of graded commutators defined as in \eqref{Eq:GradComm}.

%%%
\section{Representation algebras}\label{ss:Rep-Not}
Following \cite[\S7]{VdB1}, let $\cA$ be an algebra and $B=\bigoplus_{s \in S} \kk e_s \subset \cA$ be a semisimple algebra made of a complete set of orthogonal idempotents $(e_s)_{s\in S}$, with $S$ a finite index set identified with $\{1,\ldots,|S|\}$.  
For a dimension vector ${\bf n}=(n_s)_{s\in S}$, 
the \emph{representation space} of $\cA$ relative to $B$ of dimension $\bf{n}$, which we denote $\Rep_B(\cA, {\bf n})$\glslink{RepSp}{}, is defined as follows.
For $N=\sum_{s\in S} n_s$ and a fixed decomposition $\kk^N=\oplus_{s\in S}\kk^{n_s}$, a point $\rho \in\Rep_B(\cA, {\bf n})$ is an algebra homomorphism $\cA \to \Hom(\kk^N)$ such that $\rho(e_s)$ is the projection onto $\kk^{n_s}$ for each $s \in S$. 
The representation space is an affine scheme, 
and we focus on its coordinate ring $\cA_{\bf n}:=\kk[\Rep_B(\cA, {\bf n})]$.  
The \emph{representation algebra} $\cA_{\bf n}$\glslink{Arep}{} 
is generated by symbols 
\begin{equation} \label{Eq:genAn}
    a_{uv}, \quad a\in \cA, \,\, 1\leq u,v \leq N, 
\end{equation}
satisfying $(\lambda a+b)_{uv}=\lambda a_{uv}+b_{uv}$, $(ab)_{uv}=\sum_{1\leq r \leq N}a_{ur}b_{rv}$  for $a,b \in \cA$ and $\lambda \in \kk$, as well as 
\begin{equation*}
(e_s)_{uv}= \left\{ 
\begin{array}{cl}
  \delta_{uv}&\text{ if }n_1+\ldots+n_{s-1}+1 \leq u,v \leq n_1+\ldots+n_{s} \,, \\
  0& \text{ else.}
    \end{array}
    \right. 
\end{equation*}
Any $a\in \cA$ thus gives rise to a matrix $\XX(a)\in \Mat_{N\times N}(\cA_{\bf n})$, 
$\XX(a)_{uv}:=a_{uv}$ for $1\leq u,v\leq N$. 
In particular, the matrices $\XX(e_1),\ldots,\XX(e_{|S|})$ correspond to a block decomposition of $\Id_N$, 
with $\XX(e_s)$ having for only nonzero block that of size $n_s$ placed in position $s$. 

We generalize the notation \eqref{Eq:genAn} for the generators of $\cA_{\bf n}$ as follows: 
for $C=c^{(1)} \otimes \ldots \otimes c^{(k)} \in \cA^{\otimes k}$ (or linear combination thereof) and for any $1\leq u_j,v_j\leq N$ with $1\leq j \leq k$, we put 
\begin{equation} \label{Eq:Not-RepIndex}\glslink{shorthand-C}{}
    C_{u_1 v_1,\ldots,u_k v_k} := 
    c^{(1)}_{u_1 v_1} \cdots  c^{(k)}_{u_k v_k} \in \cA_{\bf n}\,.
\end{equation}
In particular, for any $1\leq j\leq k$ 
(with $\sigma=\sigma_{(12\cdots k)}$) 
\begin{equation} \label{Eq:RepIndex1}
C_{u_1 v_1,\ldots,u_k v_k} 
= (\sigma^{1-j} C)_{u_j v_j, \ldots, u_k v_k , u_1 v_1, \ldots, u_{j-1} v_{j-1}} \,.
\end{equation}
Given another $B$-algebra $\widetilde{\cA}$, any morphism (of algebras) $\theta: \cA\to \widetilde{\cA}$ can be extended naturally to representation algebras of the same dimension vector as a morphism (of commutative algebras) by setting 
\begin{equation} \label{ExtendMorphAn}
    \theta: \cA_{\bf n}\to \widetilde{\cA}_{\bf n}, \quad 
    \theta(a_{uv}):=(\theta(a))_{uv} \,,
\end{equation}
and extending \eqref{ExtendMorphAn} to products. 
This defining identity simply reads $\theta(\XX(a))=\XX(\theta(a))$ when written on matrices. 
Similarly, we can consider this extension to turn a derivation $\theta\in \Der_B(\cA)$ of $\cA$ into a derivation $\theta\in \Der_B(\cA_{\bf n})$, and it can be evaluated on \eqref{Eq:Not-RepIndex} using \eqref{20240805:eq1} as  
$$
\theta(C_{u_1 v_1,\ldots,u_k v_k})
= \sum_{1\leq j \leq k} \big(\theta_{(j)}C \big)_{u_1 v_1,\ldots,u_k v_k} \,.
$$
There is a natural  left action of $\Gl_{\bf n}:=\prod_{s\in S}\Gl_{n_s}(\kk)$ on 
$\Rep_B(\cA, {\bf n})$, which is induced by  the conjugation action on 
$\Mat_{N\times N}$. On the coordinate ring $\cA_{\bf n}$, this action becomes an automorphism which, for $g\in \Gl_{\bf n}$, reads as follows on generators 
\begin{equation}
    g\cdot a_{ij} = \sum_{u,v=1}^N (g^{-1})_{iu} \, a_{uv} \, g_{vj}\,.
\end{equation}
This simply becomes $g\cdot \XX(a)=g^{-1}\XX(a) g$ on the matrix-valued elements. 
Furthermore, we obtain an infinitesimal action of the Lie algebra $\gl_{\bf n}:=\prod_{s\in S}\gl_{n_s}(\kk)$, where 
$\xi_{\cA_{\bf n}}\in \Der(\cA_{\bf n})$, $\xi \in \gl_{\bf n}$, is defined in matrix form as 
\begin{equation} \label{Eq:InfAct}
    \xi_{\cA_{\bf n}}\XX(a) := [\XX(a),\xi]\,.
\end{equation}
For any $a\in \cA$, we can define the invariant element
$\tr(a)=\sum_{1\leq u \leq N} a_{uu} \in \cA^{\Gl_{\bf n}}$.
This element is clearly vanishing for $a\in [\cA,\cA]$, hence we get the following compatible $\kk$-linear maps
\begin{equation} \label{Eq:Tr-morph0}
    \tr: \cA \to \cA^{\Gl_{\bf n}}, \,\, a\mapsto \tr(a), \qquad 
    \tr: \cA_\sharp \to \cA^{\Gl_{\bf n}}, \,\, \bar{a}\mapsto \tr(a)\,.
\end{equation}

%%%%%%%%%%% NEW PART %%%%%%%%%%%%%%%
%%%%%%%%%%% NEW PART %%%%%%%%%%%%%%%
%%%%%%%%%%% NEW PART %%%%%%%%%%%%%%%
%%%%%%%%%%% NEW PART %%%%%%%%%%%%%%%
%%%%%%%%%%% NEW PART %%%%%%%%%%%%%%%
%%%%%%%%%%% NEW PART %%%%%%%%%%%%%%%
%%%%%%%%%%% NEW PART %%%%%%%%%%%%%%%

\part{Double Poisson algebra cohomologies}

%%%%%%%%%%% NEW CHAPTER %%%%%%%%%%%%%%%
%%%%%%%%%%% NEW CHAPTER %%%%%%%%%%%%%%%
%%%%%%%%%%% NEW CHAPTER %%%%%%%%%%%%%%%
%%%%%%%%%%% NEW CHAPTER %%%%%%%%%%%%%%%

\chapter{Standard Poisson cohomologies} \label{Ch:StandPoiss}

We recall the classical cohomology theories defined from a Lie algebra and a Poisson algebra. The reader can consult \cite{LGPV} for additional details, as we follow the conventions taken therein. Then, we recall the case of quasi-Poisson cohomology \cite{AKSM}, which we extend to the new \emph{gauged Poisson} setting. 

\section{Poisson cohomology}

\subsection{Definition as a Lie algebra cohomology}  \label{ss:comPCoh}
Fix a Lie algebra $\g$ over $\kk$ with Lie bracket denoted $[-,-]:\g\times \g \to \g$. That is $\g$ is a vector space with a bilinear map $[-,-]$ satisfying $(x,y,z\in \g)$ 
\begin{subequations}
\begin{align}
& \text{(skewsymmetry)}
&
&[x,y] = -[y,x]\,,
\label{LieSkew}\\
& \text{(Jacobi identity)}
&
&[x,[y,z]]-[y,[x,z]]=[[x,y],z]\,.
 \label{LieJac}
\end{align}
\end{subequations}
A Lie algebra representation of $\g$  on a vector space $W$ is given by a $\kk$-linear map 
$\rho: \g\to \End(W)$, $x\mapsto \rho(x)$, satisfying 
$\rho(x)\rho(y)-\rho(y)\rho(x)=\rho([x,y])$ for all $x,y\in \g$.   
For $x\in \g$ and $w\in W$, we simply write $x\cdot w:=\rho(x)(w)$. 
We consider the space of skew-symmetric $n$-linear maps $C^n(\g,W):=\Hom(\wedge^n \g,W)$. We can form the cochain complex $(C(\g,W),\delta_{\g,W})$ for $C(\g,W)=\oplus_{n\geq 0} C^n(\g,W)$ and $\delta_{\g,W}$ denoting the sequence of linear maps  
$\delta_{\g,W}^n:C^n(\g,W)\to C^{n+1}(\g,W)$, $n\geq 0$, defined as follows. 
For any $Q\in C^n(\g,W)$ and $x_1,\ldots,x_{n+1}\in \g$, we have 
\begin{equation}
\begin{aligned}
    \delta_{\g,W}^n(Q)&(x_1,\ldots,x_{n+1})=
\sum_{1\leq i \leq n+1} (-1)^{i+1} x_i\cdot Q(x_1,\stackrel{i}{\check{\dots}},x_{n+1}) \\
&+     \sum_{1\leq i<j \leq n+1} (-1)^{i+j}  Q([x_i,x_j],x_1,\stackrel{i}{\check{\dots}},\stackrel{j}{\check{\dots}}, x_{n+1})\,.
\end{aligned}
\end{equation}
(We use the notation \eqref{Eq:NotCheck}.)
A standard calculation allows to check the condition $\delta_{\g,W}^n\circ \delta_{\g,W}^{n-1}=0$. 
From this complex, we get the Lie algebra cohomology of $\g$ valued in $W$ as  
$\coH(\g;W)=(H^n(\g;W))_{n \geq 0}$ for\glslink{HCgW}{} 
\begin{equation}
    \coH^n(\g;W)=\frac{\ker\, \delta_{\g,W}^n: C^n(\g,W)\to C^{n+1}(\g,W)}{\im \, \delta_{\g,W}^{n-1}:C^{n-1}(\g,W)\to C^n(\g,W)}\,.
\end{equation}
The Chevalley-Eilenberg cohomology consists in picking the adjoint representation $\g^{\mathrm{ad}}$ of $\g$ on itself, and it is denoted by $\coH_{CE}(\g):= \coH(\g;\g^{\mathrm{ad}})$. 

\medskip 

Given a commutative algebra $A$, 
the space of skewsymmetric $n$-linear derivations $\mf X^n(A)$ of $A$ is the subspace of $\Hom(\wedge^n A,A)$\glslink{mfXA}{} spanned by elements $P$ such that 
$(a,b,a_2,\ldots,a_k\in A)$ 
$$P(ab,a_2,\ldots,a_k)=a P(b,a_2,\ldots,a_k)+ b P(a,a_2,\ldots,a_k)\,.$$
In particular, $\mf X^1(A)=\Der(A,A)$. 
We say that $A$ is a Poisson algebra if it is equipped with a Poisson bracket 
$\br{-,-}:A\times A\to A$, which is a Lie bracket satisfying  $(x,y,z\in A)$
\begin{align} \label{PoiLei}
&\text{(Leibniz rules)}
&\br{x,yz} =y\br{x,z}+\br{x,y}z\,, \quad 
 \br{xy,z} = x\br{y,z}+\br{x,z}y\,.
\end{align}
We can restrict the Chevalley-Eilenberg cohomology to the skewsymmetric multilinear derivations $\mf X(A)=\bigoplus_{n\geq 0}\mf X^n(A)$. 
This is the cochain complex considered for the Poisson cohomology of $A$. 
The linear map $\delta^n_{A,\br{-,-}}:\mf X^n(A)\to \mf X^{n+1}(A)$ can be written
explicitly for any $Q\in \mf X^n(A)$ and $a_1,\ldots,a_{n+1}\in A$ as 
\begin{equation}
\begin{aligned} \label{Eq:Diff-Pcoh1}
    \delta^n_{A,\br{-,-}}(Q)&(a_1,\ldots,a_{n+1})=
\sum_{1\leq i \leq n+1} (-1)^{i+1} \br{a_i, Q(a_1,\stackrel{i}{\check{\dots}}, a_{n+1}) } \\
&+  \quad  \sum_{1\leq i<j \leq n+1} (-1)^{i+j}  Q(\br{a_i,a_j},a_1,\stackrel{i}{\check{\dots}},\stackrel{j}{\check{\dots}}, a_{n+1})\,.
\end{aligned}
\end{equation}
From this complex, we get the Poisson cohomology of $A$ with respect to the Poisson bracket $\br{-,-}$ as  
${\mathrm H}_{CE}(A)=({\mathrm H}_{CE}^n(A))_{n \geq 0}$\glslink{HCEA}{} for 
\begin{equation}
    {\mathrm H}_{CE}^n(A)
=\frac{\ker\, \delta_{A,\br{-,-}}^n: \mf X^n(A)\to \mf X^{n+1}(A)}{\im \, \delta_{A,\br{-,-}}^{n-1}:\mf X^{n-1}(A)\to \mf X^n(A)}\,.
\end{equation}

\subsection{Definition from the Schouten-Nijenhuis bracket} \label{sec:standSN}

We first recall a graded version of Lie algebras and Poisson algebras. 
A graded Lie algebra (of degree $-1$) is a graded vector space $\g=\oplus_{i\in \Z} \g_i$ endowed with a graded Lie bracket of degree $-1$ denoted $[-,-]: \g\times \g \to \g$, i.e. it is a $\kk$-bilinear map 
such that $[\g_i,\g_j]\subset \g_{i+j-1}$ for any $i,j\in \Z$, 
and  ($x,y,z\in \g$)
\begin{align}
    &[x,y]= -(-1)^{(|x|-1)(|y|-1)} [y,x] \,, \label{Eq:gLie-skew} \\
    &[x,[y,z]] -(-1)^{(|x|-1)(|y|-1)} [y,[x,z]] =[[x,y],z]\,.  \label{Eq:gLie-Jac}
\end{align} 
A Gerstenhaber algebra $\mf G=\oplus_{i\in \Z}\mf G_i$ is a graded Lie algebra of degree $-1$ endowed with a graded-commutative product (i.e. $ab=(-1)^{|a||b|}ba$ for $a,b\in \mf G$) satisfying on homogeneous elements $x,y,z\in \mf G$  
\begin{equation}\label{Eq:gLie-Lei}
\begin{aligned}
    \br{x,yz} &=(-1)^{(|x|-1)|y|}y\br{x,z}+\br{x,y}z\,, \\ 
     \br{xy,z} &= x\br{y,z}+(-1)^{|y|(|z|-1)}\br{x,z}y\,.  
\end{aligned}
\end{equation}

The Schouten-Nijenhuis bracket $[-,-]_{\SN}$ is the operation on skewsymmetric multilinear derivations which satisfies for any $k,l\geq 0$ that 
$$[-,-]_{\SN}:\mf X^k(A) \times \mf X^l(A)\to \mf X^{k+l-1}(A)$$ 
is given by 
%{\color{red}(I Changed the sign convention compared to other file!!!)}
\begin{align*}
&[P,Q]_{\SN}(a_1,\ldots,a_{k+l-1})  \\
    =& (-1)^{(k-1)(l-1)}
\sum_{\sigma \in S_{l,k-1}} \sgn(\sigma) 
P(Q(a_{\sigma(1)},\ldots,a_{\sigma(l)}),a_{\sigma(l+1)},\ldots,a_{\sigma(k+l-1)}) \\
&-\sum_{\sigma \in S_{k,l-1}} \sgn(\sigma) 
Q(P(a_{\sigma(1)},\ldots,a_{\sigma(k)}),a_{\sigma(k+1)},\ldots,a_{\sigma(k+l-1)})\,,
\end{align*}
for all $a_1,\ldots, a_{k+l-1}\in A$. 
Here, $S_{p,q}\subset S_{p+q}$ denotes the subset of $(p,q)$-shuffles, which are permutations satisfying $\sigma(1)<\ldots <\sigma(p)$ and $\sigma(p+1)<\ldots <\sigma(p+q)$. 

It is well-known that $(\mf X(A),[-,-]_{\SN})$ is a graded Lie algebra of degree $-1$. Endowing $\mf X(A)$ with the wedge product, which is a graded-commutative product, turns  $(\mf X(A),[-,-]_{\SN})$ into a Gerstenhaber algebra due to the graded Leibniz rules \eqref{Eq:gLie-Lei}. 
%\footnote{Note the convention that we adopt  here for the Leibniz rule. Some authors prefer to replace this with $[P\wedge Q,R]=P\wedge [Q,R]+(-1)^{q(r-1)}[P,R]\wedge Q$. One can move from one convention to the other by changing signs appropriately depending on the degrees, as mentioned in Remark \ref{Rem:Ger-CW} for the double case.} 
%$[P\wedge Q,R]_{\SN}=[P,R]_{\SN} \wedge Q + (-1)^{k(r-1)} P\wedge [Q,R]_{\SN}$ if $|P|=k$ and $|R|=r$. 
We note that $[-,-]_{\SN}$ satisfies $[X,a]_{\SN}=X(a)$ and $[X,Y]_{\SN}=[X,Y]=X\circ Y - Y\circ X$, for any $X,Y\in \mf X^1(A)$ and $a\in A$. This means that $[-,-]_{\SN}$ is a natural extension through \eqref{Eq:gLie-Lei} of the Lie bracket of derivations given by the commutator. 
\begin{remark} \label{Rem:Conv-gLei}
    We warn the reader that there is an overall factor $(-1)^{(k-1)(l-1)}$ if one compares our definition of the Schouten-Nijenhuis bracket $[P,Q]_{\SN}$, for  $P\in  \mf X^k(A)$ and $Q\in  \mf X^l(A)$, with \cite[Eq.~(3.36)]{LGPV}. 
    This is because there is another possible convention for the Leibniz rules \eqref{Eq:gLie-Lei} in a Gerstenhaber algebra $\mc A$, which reads 
\begin{equation}\label{Eq:gLie-Lei2} 
\begin{aligned}
    \br{x,yz}' &=y\br{x,z}'+ (-1)^{(|x|-1)|z|} \br{x,y}'z\,, \\ 
     \br{xy,z}' &=(-1)^{|x|(|z|-1)} x\br{y,z}'+\br{x,z}'y\,. 
\end{aligned}
\end{equation}
To pass from \eqref{Eq:gLie-Lei} to \eqref{Eq:gLie-Lei2}, one needs to set on arbitrary homogeneous elements $x,y\in \cA$,  
$\br{x,y}'=(-1)^{(|x|-1)(|y|-1)}\br{x,y}$. 
\end{remark}

Any $\Pi\in \mf X^2(A)$ is a skewsymmetric biderivation, but it may not always define a Poisson bracket through 
\begin{equation} \label{Eq:brFromP}
    \br{-,-}_\Pi:A\times A \to A, \qquad \br{a,b}_\Pi:= \Pi(a,b)\,,
\end{equation}
as it may fail to satisfy the Jacobi identity. In fact, \eqref{Eq:brFromP} defines a Poisson bracket if and only if $[\Pi,\Pi]_{\SN}=0$. 
In that case, denote $\dd_\Pi:=[\Pi,-]_{\SN}$\glslink{dPi}{} to ease notations. It follows from \eqref{Eq:gLie-Jac} that this differential has square zero since $[\Pi,\Pi]_{\SN}=0$. As a consequence, there is another definition of the Poisson cohomology of $A$ with respect to the Poisson bracket $\br{-,-}_\Pi$ as  
$\PH(A)=(\PH^n(A))_{n \geq 0}$\glslink{PHA}{} for 
\begin{equation}
    \PH^n(A)=\frac{\ker\, \dd_\Pi: \mf X^n(A)\to \mf X^{n+1}(A)}{\im \, \dd_\Pi:\mf X^{n-1}(A)\to \mf X^n(A)}\,.
\end{equation}
A direct comparison with \eqref{Eq:Diff-Pcoh1} yields 
%that\footnote{This sign is removed if we use the convention yielding the alternative Leibniz rule, see the previous footnote.} 
$\dd_\Pi=\delta_{A,\br{-,-}_\Pi}$ when restricted to $\mf X(A)$.  
This yields in cohomology $\coH_{CE}(A) \simeq \PH(A)$, so that there is no problem in saying that both constructions yield \emph{the} Poisson cohomology of $A$.

\section{Quasi-Poisson cohomology and a gauged version} \label{ss:classQPCoh}

\subsection{ }
We start by recalling the necessary basics of quasi-Poisson cohomology following Alekseev, Kosmann-Schwarzback and Meinrenken \cite{AKSM} (conveniently adapted to our setting). 
Consider an affine variety $M$ over $\kk$ with an action of 
a reductive algebraic subgroup $G$ of $\Gl_N(\kk)$. 
Assume that $\g=\operatorname{Lie}(G)$ is endowed with an invariant and non-degenerate bilinear form $\langle -,-\rangle$. 
(For our purpose, $G$ will be a product $\Gl_{{\bf n}}=\prod_i \Gl_{n_i}(\kk)$ with the trace pairing on $\gl_{{\bf n}}=\prod_i \gl_{n_i}(\kk)$.)
This leads to an operation 
$\bigwedge^3 \g \to \kk$, $(\zeta_1,\zeta_2,\zeta_3)\mapsto \langle \zeta_1,[\zeta_2,\zeta_3]\rangle$
which we identify with an element $\phi\in \bigwedge^3\g$, the Cartan $3$-tensor. 
The infinitesimal action of $\g$ on $M$, $\zeta\mapsto \zeta_M$, is therefore inducing the $3$-vector field $\phi_M$. 

A bivector $\Xi\in \mf X^2(\kk[M])$ on $M$ is called \emph{quasi-Poisson} if $[\Xi,\Xi]_{\SN}=\phi_M$. 
Letting $\dd_\Xi:=[\Xi,-]_{\SN}$, we compute $d_\Xi^2=\frac12 [\phi_M,-]_{\SN}$. 
This is a square-zero differential when restricted to $G$-invariant multivectors $\mf X^k(\kk[M])^{G}$. 
Indeed, $[\zeta_M,R]_{\SN}=0$ for any $R\in \mf X^k(\kk[M])^{G}$ and $\xi\in \g$. 
The corresponding cohomology of $(\mf X(\kk[M])^{G},\dd_\Xi)$ is called \emph{quasi-Poisson cohomology} \cite[\S4]{AKSM}; 
we denote it $\PH_G(\kk[M])$.\glslink{qPcoh}{} 

\begin{remark} \label{Rem:qPdef}
If one is only interested in $G$-invariant cohomology theories, the quasi-Poisson condition can be replaced by $[\Xi,\Xi]_{\SN}=c\,\phi_M$ for any $c\in \kk$ (thus encompassing the Poisson case $c=0$). 
Similarly, if one can decompose $G=\prod_k G_k$ and write $\phi^{(k)}$ for the Cartan trivector of the $k$-th factor, the condition can be relaxed as  
$[\Xi,\Xi]_{\SN}=\sum_k c_k\,\phi^{(k)}_M$ with $c_k\in \kk$. 
In full generalities, one can in fact take an arbitrary $G$-invariant $\phi \in \bigwedge^3 \g$ and then ask for the condition  $[\Xi,\Xi]_{\SN}=\phi_M$ to hold, cf. \cite{Hue}. 
\end{remark}

Motivated by \cite[Rem.~4.7]{AKKN}, we generalize the previous observations. 
We introduce the linear map 
\begin{equation}
j^G_{M,k} :  \g \times {\mf X}^{k-1}(\kk[M]) \to  {\mf X}^{k}(\kk[M]), \quad 
(\zeta, R) \mapsto \zeta_M \wedge R,
\end{equation}
and we let $\im_{G-\inv}(j^G_{M,k})=\im(j^G_{M,k}) \cap \mf X^k(\kk[M])^G$ be the $G$-invariant part of its image.  

\begin{definition} \label{Def:gaugBiv}
We call a $G$-invariant bivector $\Xi\in \mf X^2(\kk[M])^G$ on $M$ \emph{gauged Poisson} if 
$[\Xi,\Xi]_{\SN} \in \im(j^G_{M,3})$. 
\end{definition}
The condition on the Schouten bracket yields a decomposition  
\begin{equation} \label{Eq:geo-gP}
    [\Xi,\Xi]_{\SN}=\sum_r \, (\zeta_r)_M \wedge \widetilde{\Xi}_r\,, \qquad 
    \zeta_r\in \g,\quad  \widetilde{\Xi}_r\in \mf X^2(\kk[M]).
\end{equation}
Invariance of $\Xi$ guarantees that $[\Xi,\Xi]_{\SN}\in \mf X^3(\kk[M])^G$, and therefore 
$[\Xi,\Xi]_{\SN} \in \im_{G-\inv}(j^G_{M,3})$. 
Thus, both $\Xi$ and $[\Xi,\Xi]_{\SN}$ are multiderivations on $\kk[M/\!/G]$ for the GIT quotient $M/\!/G:= \Spec(\kk[M]^G)$. 
Hence we get that $\dd_\Xi=\brSN{\Xi,-}$ is defined on 
$\mf X(\kk[M/\!/G])$.  
The main condition \eqref{Eq:geo-gP} implies for $R\in \mf X^k(\kk[M])^{G}$, 
\begin{equation*}
    \dd_\Xi^2(R)=\frac12 \sum_r (\xi_r)_M \wedge [\widetilde{\Xi}_r,R]_{\SN}\,,
\end{equation*}
which vanishes on $M/\!/G$ since $G$-invariant functions are in the kernel of each infinitesimal action $\zeta_M$ for any $\zeta\in \g$.

\begin{definition} \label{Def:gPH}
Let $\Xi$ be a gauged Poisson bivector on $M$. 
The cohomology of the complex $({\mf X}(\kk[M/\!/G]) , \dd_\Xi)$ is called 
the $G$-\emph{gauged Poisson cohomology} of $M$ with respect to $\Xi$. 
Explicitly, we have $\gPH_G(M)=(\gPH_G^k(M))_{k \geq 0}$ for\glslink{gaugedPHM}{} 
\begin{equation}
\gPH_G^k(M)=\frac{\ker\, \dd_\Xi:{\mf X}^k(\kk[M/\!/G])\to{\mf X}^{k+1}(\kk[M/\!/G])}{\im \, \dd_\Xi:{\mf X}^{k-1}(\kk[M/\!/G])\to {\mf X}^k(\kk[M/\!/G])}\,.
\end{equation}
\end{definition}
We may write $\gPH_G(\kk[M])$ instead of $\gPH_G(M)$ to emphasize that this is an algebraic theory. 
Note that, if we induce $\Xi$ as a bivector on $M/\!/G$, it becomes Poisson and 
therefore $\gPH_G(M)$ is simply the Poisson cohomology $\PH(M/\!/G)$. 

\subsection{ }
Note that quasi-Poisson cohomology can be defined in the absence of a bivector. Namely, assume that $\kk[M]$ is equipped with a skewsymmetric biderivation 
\begin{equation*}
    \br{-,-} : \kk[M] \times \kk[M] \longrightarrow \kk[M] 
\end{equation*}
subject to the quasi-Poisson property\footnote{Here, $\phi$ can be an arbitrary invariant element, cf. Remark \ref{Rem:qPdef}.} $(f_1,f_2,f_3\in \kk[M]$) 
 \begin{equation} \label{Eq:JacPhi}
\br{f_1,\br{f_2,f_3}} +\br{f_2,\br{f_3,f_1}} + \br{f_3,\br{f_1,f_2}} = \frac12 \phi_M(f_1,f_2,f_3)\,.
 \end{equation}
We call such a $\br{-,-}$ a \emph{quasi-Poisson bracket}. 
A standard computation yields that a quasi-Poisson bivector $\Xi$ defines a quasi-Poisson bracket since $[\Xi,\Xi]_{\SN}(f_1,f_2,f_3)$ is twice the left-hand side of \eqref{Eq:JacPhi}.

\begin{theorem} \label{Thm:qPcohGen}
If $\br{-,-}$ is quasi-Poisson, then the linear map $\delta_{\kk[M],\br{-,-}}$ of degree $+1$ defined on $\mf X(\kk[M])^{G}$ through \eqref{Eq:Diff-Pcoh1} is a square-zero differential. 
Furthermore, it defines a square-zero differential on $\mf X(\kk[M]^G)$.
\end{theorem}
The corresponding cohomology of $(\mf X(\kk[M])^{G},\delta_{\kk[M],\br{-,-}})$ is called \emph{quasi-Poisson cohomology} and is denoted 
$\coH_{CE;G}(\kk[M])$. 
We get a map 
\begin{equation}
    \coH_{CE;G}(\kk[M]) \longrightarrow 
    \coH_{CE}(\kk[M]^G)
\end{equation}
where the latter is the Poisson cohomology defined on $\mf X(\kk[M]^G)$ through \eqref{Eq:Diff-Pcoh1}.

\begin{proof}[Proof of Theorem~\ref{Thm:qPcohGen}] 
We work with $\phi$ being an arbitrary $G$-invariant element in $\bigwedge^3 \g$, cf. Remark~\ref{Rem:qPdef}. 
Since the quasi-Poisson bracket is a Poisson bracket on $\kk[M]^G$ due to \eqref{Eq:JacPhi}, the second part of the statement follows from the usual Poisson case described in Subsection~\ref{ss:comPCoh}.

Let    $R\in \mf X^k(\kk[M])^{G}$ for $k\geq 0$ and $a_0,\ldots,a_k \in \kk[M]$. 
A routine computation using \eqref{Eq:Diff-Pcoh1} shows that 
\begin{subequations}
    \begin{align}
&(\delta_{\kk[M],\br{-,-}}^{k+1} \circ \delta_{\kk[M],\br{-,-}}^k)(R)
(a_0,\ldots,a_{k+1}) \nonumber \\
=&\sum_{i<j} (-1)^{i+j} \big\{a_j,\br{a_i,R(a_0,\stackrel{i+1}{\check{\dots}},\stackrel{j+1}{\check{\dots}},a_{k+1} )}\big\} \label{Eq:qPcohA} \\
&+\sum_{i<j} (-1)^{i+j-1} \big\{a_i,\br{a_j,R(a_0,\stackrel{i+1}{\check{\dots}},\stackrel{j+1}{\check{\dots}},a_{k+1} )}\big\} \label{Eq:qPcohB}\\
&+\sum_{i<j} (-1)^{i+j} \big\{ \br{a_i,a_j},R(a_0,\stackrel{i+1}{\check{\dots}},\stackrel{j+1}{\check{\dots}},a_{k+1} )\big\} \label{Eq:qPcohC} \\
&+\sum_{s<i<j} (-1)^{i+j+s+1} R(\big\{ \br{a_i,a_j}, a_s\big\} ,a_0,\stackrel{s+1}{\check{\dots}},\stackrel{i+1}{\check{\dots}},\stackrel{j+1}{\check{\dots}},a_{k+1} )  \label{Eq:qPcohD} \\
&+\sum_{s<i<j} (-1)^{i+j+s} R(\big\{ \br{a_s,a_j}, a_i\big\} ,a_0,\stackrel{s+1}{\check{\dots}},\stackrel{i+1}{\check{\dots}},\stackrel{j+1}{\check{\dots}},a_{k+1} ) \label{Eq:qPcohE} \\
&+\sum_{s<i<j} (-1)^{i+j+s-1} R(\big\{ \br{a_s,a_i}, a_j\big\} ,a_0,\stackrel{s+1}{\check{\dots}},\stackrel{i+1}{\check{\dots}},\stackrel{j+1}{\check{\dots}},a_{k+1} ) \,. \label{Eq:qPcohF}
    \end{align}
\end{subequations}
For $i<j$, let us introduce the notation 
\begin{equation*}
    R_{ij}:=R(a_0,\stackrel{i+1}{\check{\dots}},\stackrel{j+1}{\check{\dots}},a_{k+1} )\,.
\end{equation*}
We also write $R_{ij}[a_k\mapsto b]$ to denote that we replace the argument $a_k$ by $b\in \kk[M]$ in $R_{ij}$. 
First, remark that thanks to \eqref{Eq:JacPhi}: 
\begin{align*}
 &\eqref{Eq:qPcohA}+   \eqref{Eq:qPcohB}+\eqref{Eq:qPcohC} \\
 =&\sum_{i<j} (-1)^{i+j} \left( \big\{a_j,\br{a_i,R_{ij} }\big\}
+  \big\{a_i,\br{R_{ij},a_j } +  \big\{R_{ij},\br{a_j,a_{i} }\big\} 
 \right) \\
 =&\frac12 \sum_{i<j} (-1)^{i+j} \ \phi_M(a_j,a_i,R_{ij}) \,.
\end{align*}
Second, a similar manipulation yields 
\begin{align*}
 &\eqref{Eq:qPcohD}+   \eqref{Eq:qPcohE} +\eqref{Eq:qPcohF} \\
 =&\sum_{s<i<j} (-1)^{i+j+s} R(\alpha_{ijs} ,a_0,\stackrel{s+1}{\check{\dots}},\stackrel{i+1}{\check{\dots}},\stackrel{j+1}{\check{\dots}},a_{k+1} ) \,, \\
& \text{ for } 
\alpha_{ijs}:= \big\{a_s,\br{a_i,a_j}\big\} + \big\{a_i,\br{a_j,a_s}\big\} 
+ \big\{a_j,\br{a_s,a_i}\big\}= \frac12 \phi_M(a_i,a_j, a_s)\,.
\end{align*}

At this point, it is clear that we get a square-zero differential on $\mf X(\kk[M]^{G})$ since the two reduced expressions vanish by $\g$-invariance of $a_i,a_j,a_s$. But our first claim is stronger: we want the two expressions to cancel out by invariance of $R$. 

Write $\phi=\sum_{t\in T} \zeta^{t,1} \wedge \zeta^{t,2} \wedge \zeta^{t,3}$ where $\zeta^{t,l}\in \g$. We find for twice the first expression, 
\begin{equation} \label{Eq:qPCohGen1}
    \begin{aligned}
& \sum_{i<j} (-1)^{i+j} \ \phi_M(a_j,a_i,R_{ij})    \\ 
&= \sum_{i<j}(-1)^{i+j}\sum_{\substack{t\in T\\ \tau\in S_3}} \sgn(\tau) 
 \zeta^{t,\tau(1)}(a_j) \zeta^{t,\tau(2)}(a_i)  \zeta^{t,\tau(3)}(R_{ij}) \\
&= \sum_{\substack{i<j\\r\neq i,j}}(-1)^{i+j}\sum_{\substack{t\in T\\ \tau\in S_3}} \sgn(\tau) \, 
 \zeta^{t,\tau(1)}(a_j) \zeta^{t,\tau(2)}(a_i)  R_{ij}[a_r \mapsto \zeta^{t,\tau(3)}(a_r)]
    \end{aligned}
\end{equation}
where the second equality holds by $G$-invariance of $R$. 

For the second expression, we start by calculating 
\begin{equation}
    \begin{aligned}
&R(\phi_M(a_i,a_j, a_s),a_0,\stackrel{s+1}{\check{\dots}},\stackrel{i+1}{\check{\dots}},\stackrel{j+1}{\check{\dots}},a_{k+1} ) \\
=&
\sum_{\substack{t\in T\\ \tau\in S_3}} \sgn(\tau)(-1)^{s}\, R_{ij}[a_s \mapsto \zeta^{t,\tau(3)}(a_s)]\, 
\zeta^{t,\tau(1)}(a_i) \zeta^{t,\tau(2)}(a_j) \\
&+ \sum_{\substack{t\in T\\ \tau\in S_3}}\sgn(\tau) (-1)^{i-1}\, R_{sj}[a_i \mapsto \zeta^{t,\tau(1)}(a_i)]\, 
\zeta^{t,\tau(2)}(a_j) \zeta^{t,\tau(3)}(a_s) \\
&+ \sum_{\substack{t\in T\\ \tau\in S_3}} \sgn(\tau)(-1)^{j-2}\, R_{si}[a_j \mapsto \zeta^{t,\tau(2)}(a_j)]\, 
\zeta^{t,\tau(1)}(a_i) \zeta^{t,\tau(3)}(a_s) 
    \end{aligned}
\end{equation}
where the signs come from putting the first argument (i.e. $\zeta^{t,\tau(3)}(a_s)$ for the first line, etc.) appearing in $R$ in the corresponding position, keeping in mind that $s<i<j$. 
Hence we can write  
\begin{equation}
    \begin{aligned}
&\sum_{s<i<j} (-1)^{i+j+s} R(\phi_M(a_i,a_j, a_s),a_0,\stackrel{s+1}{\check{\dots}},\stackrel{i+1}{\check{\dots}},\stackrel{j+1}{\check{\dots}},a_{k+1} ) \\
=& \sum_{r<i<j}
\sum_{\substack{t\in T\\ \tau\in S_3}}\sgn(\tau) (-1)^{i+j}\, R_{ij}[a_r \mapsto \zeta^{t,\tau(3)}(a_r)]\, 
\zeta^{t,\tau(1)}(a_i) \zeta^{t,\tau(2)}(a_j) \\
&+ \sum_{i<r<j} 
\sum_{\substack{t\in T\\ \tau\in S_3}}\sgn(\tau) (-1)^{i+j+1}\, R_{ij}[a_r \mapsto \zeta^{t,\tau(1)}(a_r)]\, 
\zeta^{t,\tau(2)}(a_j) \zeta^{t,\tau(3)}(a_i) \\
&+ \sum_{i<j<r}
\sum_{\substack{t\in T\\ \tau\in S_3}}\sgn(\tau) (-1)^{i+j}\, R_{ij}[a_r \mapsto \zeta^{t,\tau(2)}(a_r)]\, 
\zeta^{t,\tau(1)}(a_j) \zeta^{t,\tau(3)}(a_i) 
    \end{aligned}
\end{equation}
where we relabeled indices. Using the total skewsymmetry of $\phi$, this becomes 
\begin{equation}
    \sum_{\substack{i<j\\r\neq i,j}} \sum_{\substack{t\in T\\ \tau\in S_3}} \sgn(\tau) (-1)^{i+j+1}
 \zeta^{t,\tau(1)}(a_j) \zeta^{t,\tau(2)}(a_i) \ R_{ij}[a_r \mapsto \zeta^{t,\tau(3)}(a_r)]\,,
\end{equation}
which is the opposite of \eqref{Eq:qPCohGen1}. Hence, we can conclude.  
\end{proof}

%%%%%%%%%%% NEW CHAPTER %%%%%%%%%%%%%%%
%%%%%%%%%%% NEW CHAPTER %%%%%%%%%%%%%%%
%%%%%%%%%%% NEW CHAPTER %%%%%%%%%%%%%%%
%%%%%%%%%%% NEW CHAPTER %%%%%%%%%%%%%%%

\chapter{Double Poisson algebras and a double Poisson cohomology} 
\label{Ch:dPA}

We fix an algebra $\cA$ and a subalgebra $B\subset \cA$ such that we can view $\cA$ as a $B$-algebra. 
All results from this chapter are taken from \cite{VdB1} for the basic constructions, and \cite{PV} for the double Poisson cohomology.

\section{Multi-brackets, double Poisson brackets and graded versions}

\subsection{ } 
\label{def:n-dbr}
For $n\geq 1$, an $n$-bracket on $\cA$ is a linear map 
$\dgal{-,\ldots,-}:\cA^{\otimes n} \to \cA^{\otimes n}$ such that ($a_i,b,c\in \cA$)
\begin{subequations}
\begin{align}
& \text{\small (cyclic skewsymmetry)} 
&
&\dgal{a_1,a_2,\ldots,a_n} = (-1)^{n-1} \dgal{a_2,\ldots,a_n,a_1}^\sigma , 
\label{Eq:nbr-Cycl} \\
&
\text{\small (Leibniz rule)} &
&\dgal{a_1,\ldots,a_{n-1},bc}= b\dgal{a_1,\ldots,a_{n-1},c}
\label{Eq:nbr-DerOut} \\
&
&
&
+ \dgal{a_1,\ldots,a_{n-1},b}c\,.
\nonumber 
\end{align}
\end{subequations}
In \eqref{Eq:nbr-DerOut}, one uses the outer bimodule structure corresponding to $j=0$ in \eqref{eq:star}. 
An $n$-bracket $\dgal{-,\ldots,-}$ is $B$-linear if $\dgal{a_1,\ldots,a_{n}}=0$ whenever 
$a_i \in B$ for some $i$.  
A $1$-bracket is an element of $\Der(\cA)$, or $\Der_B(\cA)$ for the $B$-linear case. 

Due to the cyclic skewsymmetry \eqref{Eq:nbr-Cycl}, an $n$-bracket is an $\cA^{\otimes n}$-valued derivation in each of its entries. More precisely, using the star product \eqref{eq:star}, 
\begin{equation} \label{Eq:nbr-DerAll}
\begin{aligned}
   \dgal{a_1,\ldots, a_i \tilde{a}_i,\ldots,a_{n}}
=&\, a_i \star_{i}\dgal{a_1,\ldots, \tilde{a}_i,\ldots,a_{n}} \\
&+\dgal{a_1,\ldots, a_i,\ldots,a_{n}} \star_{n-i} \tilde{a}_i\,, 
\end{aligned}
\end{equation}
for any $a_1,\ldots,a_n,\tilde{a}_i\in \cA$. 
We shall refer to a ($B$-linear) $2$-bracket as a ($B$-linear) double bracket; 
its derivation rule in the first entry is then taken with respect to the inner bimodule structure \eqref{bmodinner}. 

We denote by $\BRA(\cA)_n$\glslink{BRA}{} the $\kk$-vector space of all $n$-brackets, and form the graded vector space $\BRA(\cA)=\oplus_{n\geq 1} \BRA(\cA)_n$ with $n$-brackets in degree $n$. We write $\BRA_B(\cA)=\oplus_{n\geq 1}\BRA_B(\cA)_n$ for the corresponding subspace of $B$-linear $n$-brackets. 
It will be useful to add $\cA_\sharp=\cA/[\cA,\cA]$ as a subspace of elements of degree $0$. Thus, we introduce the completed vector space $\wBRA(\cA)=\oplus_{n\geq 0} \wBRA(\cA)_n$ with $\wBRA(\cA)_0=\cA_\sharp$ and $\wBRA(\cA)_n=\BRA(\cA)_n$ for $n\geq 1$; we define $\wBRA_B(\cA)$ in the same way. 

\medskip 

From now on, assume that $\cA$ is a $B$-algebra and all $n$-brackets are assumed to be $B$-linear. (The standard case consists in taking $B=\kk$.) 
Since the subgroup $\DDer_B(\cA)\subset \DDer(\cA)$ of double derivations relative to $B$ is 
an $\cA$-bimodule, cf. \eqref{DDer} and below,  
we can form $\mb T^\ast \cA:=T_{\cA} \DDer_B(\cA)$\glslink{TstarA}{} as the tensor algebra over $\cA$ having $\cA$ in degree $0$ and $\DDer_B(\cA)$ in degree $+1$. 
Thus, $\mb T^\ast \cA$ is a graded algebra for the tensor multiplication. 
It follows that any $Q\in (\mb T^\ast \cA)_n$ is a linear combination of terms 
$\delta_1\cdots \delta_n$ with $\delta_i\in \DDer_B(\cA)$. 
Elements of $\mb T^\ast \cA$ are called \emph{noncommutative multivector fields}.

\begin{proposition}[\cite{VdB1}, \S4.1] \label{Pr:MapMu} 
For any $n\geq 1$, there is a well-defined map 
$\mu_n: (\mb T^\ast \cA)_n \longrightarrow \wBRA_B(\cA)_n$  
where 
$$\mu_n(Q):=\dgal{-,\ldots,-}_Q=\sum_{0\leq i \leq n-1} (-1)^{(n-1)i}
\sigma^i \circ \dgal{-,\ldots,-}_Q^{\sim}\circ \sigma^{-i}\,,$$
is obtained $\kk$-linearly from 
\begin{equation} \label{Eq:MapMu}
\begin{aligned}
    &\dgal{a_1,\ldots,a_n}^{\sim}_{\delta_1\ldots \delta_n} \\
    &\quad := \, 
\delta_n(a_n)'\delta_1(a_1)''\otimes\delta_1(a_1)'\delta_{2}(a_2)''\otimes \,\cdots\,
\otimes \delta_{n-1}(a_{n-1})' \delta_n(a_n)''\,, 
\end{aligned}   
\end{equation}
for $a_j \in \cA$ and $\delta_j\in \DDer_B(\cA)$ with $1\leq j \leq n$. 
Moreover, the map $\mu_n$ factors through $(\mb T^\ast \cA)_{\sharp,n}$, i.e. $\mu_n(Q)$ only depends on $Q$ modulo graded commutators. 
\end{proposition}
We trivially have a map $\mu_0:(\mb T^\ast \cA)_{\sharp,0}\to \wBRA_B(\cA)_0$ in degree $0$ which is $\Id_{\cA_\sharp}$, hence is an isomorphism. To have isomorphisms in higher degrees, some extra conditions are necessary. 
\begin{proposition}[\cite{VdB1}, \S4.1] \label{Pr:MapMu-Iso}
Assume that $\cA$ is finitely generated over $B$, left and right flat, and that the $\cA$-bimodule of $B$-relative differential forms $\Omega^1_{\cA/B}$ is projective. 
Then $\mu_n: (\mb T^\ast \cA)_{\sharp,n} \to \wBRA_B(\cA)_n$ is an isomorphism for each $n\geq 1$. 
\end{proposition}

\begin{remark}
The fact that $\mu_n$ is not, in general, an isomorphism is the reason why the definition of double Poisson cohomology on $\mb T^\ast\cA$ by Pichereau and Van de Weyer \cite{PV} that we shall give below is not always sufficient. We shall present another definition of double Poisson cohomology on $\wBRA_B(\cA)$ in Chapter \ref{CH:Gen-dPcoh}.   
\end{remark}

\subsection{ } Let us now focus on double brackets, i.e. elements of $\wBRA_B(\cA)_2$. 
Given an element $\dgal{-,-}\in \wBRA_B(\cA)_2$,  
the defining rules \eqref{Eq:nbr-Cycl}-\eqref{Eq:nbr-DerOut} and \eqref{Eq:nbr-DerAll} read  
\begin{align}
\dgal{a,b}&= - \dgal{b,a}^\sigma , 
\label{Eq:db-skew} \\
    \dgal{a,bc}&=b\dgal{a,c} + \dgal{a,b}c \,, 
   \label{Eq:db-Rleib} \\
    \dgal{ab,c}&=a \star \dgal{b,c} + \dgal{a,c} \star b\,,
\label{Eq:db-Lleib}
\end{align}
where we use the bimodule structures \eqref{bmodouter}-\eqref{bmodinner}.  
We form the $B$-linear map 
\begin{equation} \label{Eq:ExtLeft}
    \dgal{-,-}_L: \cA\times \cA^{\otimes 2}\to \cA^{\otimes 3}, \qquad 
\dgal{a,d}_L:=\dgal{a,d'}\otimes d'', 
\end{equation}
cf. \eqref{20240805:eq1}. 
Then, we can define a $B$-linear map $\dgal{-,-,-}:\cA^{\times 3}\to \cA^{\otimes 3}$ by 
\begin{equation} \label{Eq:dJac}
 \dgal{a,b,c}  :=
    \dgal{a,\dgal{b,c}}_L+ (\dgal{b,\dgal{c,a}}_L)^\sigma 
    + (\dgal{c,\dgal{a,b}}_L)^{\sigma^2} \,,
\end{equation}
for any $a,b,c\in \cA$. By \S2.3 in \cite{VdB1}, we have $\dgal{-,-,-}\in \wBRA_B(\cA)_3$.

\begin{definition}
\label{def:dbr-Poiss}\glslink{dPbracket}{}
A $B$-linear double bracket $\dgal{-,-}\in \wBRA_B(\cA)_2$ such that the triple bracket  
$\dgal{-,-,-}\in \wBRA_B(\cA)_3$ defined through \eqref{Eq:dJac} identically vanishes is 
a \emph{double Poisson bracket}. We then say that $(\cA,\dgal{-,-})$ is a \emph{double Poisson algebra}. 
\end{definition}
Let us record that the property of being Poisson is equivalent to 
\begin{equation} \label{Eq:dJacExpl}
\sum_{t\in \Z_3} \sigma^t \circ
(\dgal{-,-}\otimes \Id_{\cA} ) \circ (\Id_{\cA} \otimes \dgal{-,-}) \circ \sigma^{-t}=0\,.
\end{equation}

If $P \in (\mb T^\ast \cA)_2$, we can write $P=\sum_{\ell \in L}\delta_1^{(\ell)} \delta_2^{(\ell)}$ 
with $\delta_i^{(\ell)}\in \DDer_B(\cA)$ and a finite set $L$. 
Then, Proposition \ref{Pr:MapMu} defines the double bracket $\mu_2(P)=\dgal{-,-}_P$  as 
\begin{equation} \label{Eq:dbr-P}
\begin{aligned}
        \dgal{a,b}_P= \sum_{\ell \in L} \Big( &
\delta_2^{(\ell)}(b)' \delta_1^{(\ell)}(a)'' \otimes  \delta_1^{(\ell)}(a)'  \delta_2^{(\ell)}(b)''  \\
&\quad - \delta_1^{(\ell)}(b)' \delta_2^{(\ell)}(a)'' \otimes  \delta_2^{(\ell)}(a)'  \delta_1^{(\ell)}(b)''
\Big)\,, \qquad a,b\in \cA.
\end{aligned}
\end{equation}

\begin{example} \label{Ex:Classifkx}
Any $n$-bracket $\dgal{-}\in \BRA(\kk[x])_n$ is determined by the value 
\begin{equation*}
    \dgal{x,\ldots,x} \in \kk[x]^{\otimes n}\,\, \text{ with }
    \dgal{x,\ldots,x}-(-1)^{n+1}\dgal{x,\ldots,x}^\sigma=0,
\end{equation*}
due to the derivation rules \eqref{Eq:nbr-DerOut}, \eqref{Eq:nbr-DerAll}, 
and the condition \eqref{Eq:nbr-Cycl} of cyclic skewsymmetry. Thus, 
for integers $a,b \geq 0$, there is a double bracket $\dgal{-,-}_{a,b}$ uniquely defined by 
\begin{align}  \label{Eq:dbr-ab} 
\dgal{x,x}_{a,b}=x^a \otimes x^b - x^b \otimes x^a\,. 
\end{align}
Note that $\dgal{-,-}_{a,b}=-\dgal{-,-}_{b,a}$. 
Hence, an arbitrary nonzero double bracket on $\kk[x]$ is of the form $\dgal{-,-}=\sum_{0\leq b < a} c_{a,b} \dgal{-,-}_{a,b}$ for finitely many nonzero $c_{a,b}\in \kk$. 
It is a result of Van den Bergh \cite{VdB1} that the cases $(a,b)=(1,0)$ and $(a,b)=(2,1)$ define double Poisson brackets. Furthermore,  Powell \cite{P16} classified all double Poisson brackets on $\kk[x]$: they are of the form 
\begin{equation} \label{Eq:dbr-xlin}
    \dgal{-,-} := \lambda \dgal{-,-}_{1,0} + \mu \dgal{-,-}_{2,0} + \nu \dgal{-,-}_{2,1}\,,
\end{equation}
for $\lambda,\mu,\nu \in \kk$ satisfying $\lambda\nu-\mu^2=0$. 
All $\mu_n:(\mb T^\ast \kk[x])_{\sharp,n}\to \wBRA(\kk[x])_n$ from Proposition \ref{Pr:MapMu-Iso} are isomorphisms. 
In particular, if we introduce 
\begin{equation} \label{DDer-delx}
 \partial_x\in \DDer(\kk[x]), \qquad \partial_x(x)=1\otimes 1,
\end{equation}
one has $(x^a \partial_x)(x)=1\otimes x^a$ for any $a\geq 1$. We easily see from \eqref{Eq:dbr-P} that \eqref{Eq:dbr-ab} can be written as  
\begin{equation*}
   \dgal{-,-}_{a,b} = \mu_2(x^a \partial_x x^b \partial_x) \,.
\end{equation*}
\end{example}

A double Poisson bracket $\dgal{-,-}$ on $\cA$ gives rise to a map 
$$\br{-,-}=\mult \circ \dgal{-,-}:\cA\times \cA \to \cA\,,$$
after multiplication of the tensor factors. 
Noting that $\br{[\cA,\cA],-}=0$, cf. \cite[\S2.4]{VdB1}, we get a linear map (denoted in the same way)
\begin{equation} \label{Eq:mult-br-1}
    \br{-,-}=\cA_\sharp\times \cA \to \cA\,, \qquad \cA_\sharp:= \cA/[\cA,\cA]\,.
\end{equation}
The operation in \eqref{Eq:mult-br-1} is 
such that $\br{\bar{a},-}:\cA\to \cA$ is a $B$-linear derivation for any $\bar{a}\in \cA_\sharp$.  
Furthermore, we note that \eqref{Eq:mult-br-1} descends to a skewsymmetric map
\begin{equation} \label{Eq:mult-brLie}
    \br{-,-}_\sharp=\cA_\sharp\times \cA_\sharp \to \cA_\sharp\,.
\end{equation}
Furthermore, one can check that \eqref{Eq:mult-brLie} satisfies the Jacobi identity because $\dgal{-,-}$ is a double Poisson bracket. 
Hence, $(\cA_\sharp,\br{-,-}_\sharp)$ is a Lie algebra.

\subsection{ } % graded case 
Let $\cG$ be a $\Z$-graded algebra. 

\begin{definition}[\cite{VdB1}, \S2.7] \label{def:gdPA}
A \emph{graded double bracket} (of degree $-1$) on $\cG$ is a $\kk$-bilinear map 
$$\dgal{-,-}: \cG\otimes \cG \longrightarrow \cG \otimes \cG$$ 
of degree $-1$ satisfying for any homogeneous  $a,b,c\in \cG$, 
\begin{align}
\dgal{a,b}&= -(-1)^{(|a|-1)(|b|-1)} \dgal{b,a}^\sigma , 
\label{Eq:gdb-skew} \\
    \dgal{a,bc}&=(-1)^{(|a|-1)|b|}b\dgal{a,c} + \dgal{a,b}c \,, 
   \label{Eq:gdb-Rleib} \\
    \dgal{ab,c}&=a \star \dgal{b,c} + (-1)^{|b|(|c|-1)} \dgal{a,c} \star b\,. 
\label{Eq:gdb-Lleib}
\end{align}
One needs to be careful of the signs coming from the permutation of tensor factors in those expressions, cf. \eqref{gradmult}. 
Moreover, using the left extension \eqref{Eq:ExtLeft}, if the graded double bracket is such that  
\begin{equation} \label{Eq:gdb-dJac}
   \begin{aligned}
 \dgal{a,\dgal{b,c}}_L&+(-1)^{(|a|-1)(|b|+|c|)} (\dgal{b,\dgal{c,a}}_L)^\sigma \\
&+(-1)^{(|c|-1)(|a|+|b|)} (\dgal{c,\dgal{a,b}}_L)^{\sigma^2} = 0\,,        
   \end{aligned}
\end{equation}
 then we say that $\dgal{-,-}$ is a \emph{graded double Poisson bracket}. 
We refer to the pair $(\cG,\dgal{-,-})$  as a \emph{double Gerstenhaber algebra}. 
\end{definition}
As in the ungraded case, the two Leibniz rules \eqref{Eq:gdb-Rleib} and \eqref{Eq:gdb-Lleib} are equivalent under the cyclic skewsymmetry \eqref{Eq:gdb-skew}. 

\begin{remark}
Definition \ref{def:gdPA} is Van den Bergh's original definition. 
For a double (Poisson) bracket of degree $d \in \Z$, it suffices to replace all the factors of the form $(|v|-1)$ with $v\in \cG$ in the exponents by $(|v|+d)$, see \cite[\S5.1]{BCER}. 
The case $d=0$ is considered by D'Alesio \cite{DA} (on \emph{differential} graded algebras). 
The results presented below have natural analogues in any degree $d \in \Z$.  
\end{remark}
\begin{remark} \label{Rem:Ger-CW}
Casati and Wang \cite{CW} have considered a particular class of double Gerstenhaber algebras. However, the derivation rules that they use are \emph{different} than \eqref{Eq:gdb-Rleib}--\eqref{Eq:gdb-Lleib}. 
The relation can be made as follows: if $\dgal{-,-}$ defines a double Gerstenhaber algebra structure as in Definition \ref{def:gdPA}, 
then the map $\dgal{-,-}_{CW}$ defined on homogeneous elements $a,b\in V$ by $\dgal{a,b}_{CW}=(-1)^{(|a|-1)(|b|-1)}\dgal{a,b}$ satisfies the properties of \cite{CW} (and vice-versa). 
This choice of alternative conventions already exist in the commutative case, cf. Remark \ref{Rem:Conv-gLei}. 
\end{remark} 

Recall the notion of graded Lie algebra of degree $-1$ from  Subsection~\ref{sec:standSN}.

\begin{proposition} \label{Pr:gdPA-gLie}
If $(\cG,\dgal{-,-})$ is a double Gerstenhaber algebra, 
then the $\kk$-bilinear map 
$$\br{-,-}:\cG\times \cG \to \cG\,, \quad \br{a,b}=\mult \dgal{a,b}=\dgal{a,b}'\dgal{a,b}''$$
descends to a graded Lie bracket on $\cG_\sharp$. 
\end{proposition}
\begin{proof}
This result in the non-graded case is shown in  \cite[\S2.4]{VdB1}, while the graded case with operations of degree $d=0$ is in \cite[\S2.3]{DA}.  
Our statement, which corresponds to degree $d=-1$, can be proved in the same way after carefully accounting for the signs due to the grading and the fact that $\cG_\sharp=[\cG,\cG]$ makes use of \eqref{Eq:GradComm}. 
\end{proof} 

The previous setting can be applied to $\cG=\mb T^\ast \cA$ as follows. 

\begin{theorem} [\cite{VdB1},\S3.2] \label{Thm:dSN}
There is a unique structure of double Gerstenhaber algebra on $\mb T^\ast \cA$ whose graded double Poisson bracket $\dSN{-,-}$ of degree $-1$  is determined for 
$a,b\in \cA$ and $\delta,\Delta\in \DDer_B(\cA)$ by  
\begin{subequations}
\begin{align}
    \dSN{a,b}&=0\,,  \label{Eq:dSN-a}\\
    \dSN{\delta,a}&=\delta(a)\,,    \label{Eq:dSN-b} \\
    \dSN{\delta,\Delta}&= 
\tau_{(23)} \Big((\delta \otimes \Id_{\mc A})\circ \Delta 
- (\Id_{\mc A}\otimes \Delta)\circ \delta \Big) \label{Eq:dSN-c} \\
&\quad +\tau_{(12)} \Big((\Id_{\mc A} \otimes \delta)\circ \Delta 
- (\Delta \otimes \Id_{\mc A})\circ \delta \Big)\,.
\nonumber 
\end{align}
\end{subequations}
Consequently, $(\mb T^\ast \cA)_\sharp$ inherits the graded Lie bracket $\brSN{-,-}:=\mult \circ \dSN{-,-}$. 
\end{theorem}

We obtain the following standard result from the fact that $((\mb T^\ast \cA)_\sharp,\brSN{-,-})$ is a graded Lie algebra. 
%, it is standard that an element $P\in (\mb T^\ast \cA)_{2,\sharp}$ satisfying 
%$\brSN{P,P}=0$ induces a square-zero differential $\dd_P=\brSN{P,-}$ of degree $+1$.  
\begin{proposition}[\cite{PV,VdW}] \label{Prop:PVdW-diff}
If $P\in (\mb T^\ast \cA)_2$ satisfies $\brSN{P,P}=0$ modulo  $[\mb T^\ast \cA,\mb T^\ast \cA]$, then the map\glslink{dPnc}{}  
\begin{equation} \label{Eq:dP-PVdW}
    \dd_P=\brSN{P,-}: (\mb T^\ast \cA)_\sharp \to (\mb T^\ast \cA)_\sharp 
\end{equation}
is a square zero differential of degree $+1$. 
\end{proposition}
\begin{proof}
The degree $+1$ of $\dd_P$ comes from the fact that $\dSN{-,-}$ (hence $\brSN{-,-}$) has degree $-1$, while $P$ has degree $+2$. 
Since $\brSN{-,-}$ is a graded Lie bracket, we can write using \eqref{Eq:gLie-Jac}  
\begin{align*}
    \dd_P^2(R)&=\brSN{P,\brSN{P,R}} \\
&=\frac12\brSN{P,\brSN{P,R}}+\frac12\left(\brSN{\brSN{P,P},R} - \brSN{P,\brSN{P,R}} \right), 
\end{align*}
for any $R\in (\mb T^\ast \cA)_\sharp$; 
this is zero by assumption on $P$. 
\end{proof}

\begin{theorem}[\cite{VdB1},\S4.2] \label{Thm:dJac-PP}
Fix $P\in (\mb T^\ast \cA)_2$ and let $\dgal{-,-}:=\dgal{-,-}_P$ be defined as in Proposition \ref{Pr:MapMu}. Then 
\begin{equation} \label{Eq:dJac-PP}
\frac12 \dgal{a,b,c}_{\brSN{P,P}} =
    \dgal{a,\dgal{b,c}}_L+ (\dgal{b,\dgal{c,a}}_L)^\sigma 
    + (\dgal{c,\dgal{a,b}}_L)^{\sigma^2} \,,
\end{equation}
for any $a,b,c \in \cA$. In particular, $\dgal{-,-}$ is a double Poisson bracket 
if $\brSN{P,P}=0\in (T^\ast \cA)_\sharp$. 
\end{theorem}

\subsection{ } 
Fix $P\in(\mb T^\ast \cA)_k$ and $Q\in (\mb T^\ast \cA)_n$. 
One can compute a representative of $\brSN{P,Q}\in (\mb T^\ast \cA)_{k+n-1}$ as follows. 

We can assume without loss of generality that $P=P_1\ldots P_k$ and $Q=Q_1\ldots Q_n$ for $P_i,Q_i\in \DDer(\mc A)$. Since $\dSN{-,-}$ is a graded double Poisson bracket of degree $-1$, we get from Definition \ref{def:gdPA} that $\dSN{P,Q}$ equals 
\begin{align*}
    &\sum_{\substack{1\leq i \leq k\\1\leq j \leq n}}  (-1)^{(k-1)(j-1)} 
Q_1\cdots Q_{j-1}\left(P_1\cdots P_{i-1}\star \dSN{P_i,Q_j}\star P_{i+1}\cdots P_k \right) Q_{j+1}\cdots Q_n\\
&=\sum_{\substack{1\leq i \leq k\\1\leq j \leq n}} (-1)^{(k-1)(j-1)} 
(-1)^{(i-1)|\dSN{P_i,Q_j}'|+(k-i)|\dSN{P_i,Q_j}''|+(i-1)(k-i)} \\
&\qquad Q_1\ldots Q_{j-1} \dSN{P_i,Q_j}' P_{i+1}\ldots P_k \otimes P_1\ldots P_{i-1} \dSN{P_i,Q_j}'' Q_{j+1}\ldots Q_n\,.
\end{align*}
By definition, $\dSN{\delta,\Delta}\in \DDer(\mc A)\otimes \mc A\oplus \mc A\otimes \DDer(\mc A)$ for any $\delta,\Delta\in \DDer(\mc A)$, see \eqref{Eq:dSN-c}, hence we can write\footnote{As is usual with such multilinear operations, we should have a linear combination of such terms, which we do not write down explicitly to ease notation because all the operations are linear.}  
\begin{equation*}
    \dSN{P_i,Q_j}=\delta^i_j \otimes d^i_j + c^i_j \otimes \gamma^i_j
\end{equation*}
for $b^i_j,c^i_j\in \cA$ (of degree $0$), $\delta^i_j,\gamma^i_j\in \DDer(\cA)$ (of degree $1$) such that 
\begin{align*} 
\delta^i_j\otimes d^i_j:=&\tau_{(23)}\left((P_i\otimes \Id_{\cA})\circ Q_j - (\Id_{\cA}\otimes Q_j)\circ P_i\right) \,,     \\
c^i_j\otimes \gamma^i_j:=&\tau_{(12)}\left((\Id_{\cA} \otimes P_i)\circ Q_j - (Q_j \otimes \Id_{\cA})\circ P_i\right) \,.  
\end{align*}
Explicitly, we have by evaluation on an arbitrary $a\in \cA$ the identities 
\begin{subequations}
\begin{align}
\delta^i_j(a)\otimes d^i_j=&
P_i(Q_j(a)')' \otimes Q_j(a)'' \otimes P_i(Q_j(a)')'' \nonumber \\
&- P_i(a)' \otimes Q_j(P_i(a)'')''\otimes Q_j(P_i(a)'')'\,, \label{Eq:dSN-QP1}    \\
c^i_j\otimes \gamma^i_j(a)=& 
P_i(Q_j(a)'')' \otimes Q_j(a)'\otimes P_i(Q_j(a)'')'' \nonumber \\
&- Q_j(P_i(a)')''\otimes  Q_j(P_i(a)')' \otimes P_i(a)'' \,. \label{Eq:dSN-QP2}  
\end{align}
\end{subequations}
We can therefore write 
\begin{equation*}
    \begin{aligned}
   & \mult\circ \dSN{P,Q}  \\
=&\sum_{\substack{1\leq i \leq k\\1\leq j \leq n}} (-1)^{(k-1)(j+i)} 
 Q_1\ldots Q_{j-1} \delta^i_j P_{i+1}\cdots P_k P_1\cdots P_{i-1} d^i_j Q_{j+1}\ldots Q_n \\
&+\sum_{\substack{1\leq i \leq k\\1\leq j \leq n}} (-1)^{(k-1)(j+i-1)} 
 Q_1\ldots Q_{j-1} c^i_j P_{i+1}\cdots P_k P_1\cdots P_{i-1} \gamma^i_j Q_{j+1}\ldots Q_n.
\end{aligned}
\end{equation*}
The class of this element in $(\mb T^\ast \cA)_\sharp$ gives precisely $\brSN{P,Q}$ by definition, see the end of Theorem \ref{Thm:dSN}. 
For later use, note that this expression when $k=2$ and $P=P_1P_2$ becomes 
\begin{equation}
  \begin{aligned} \label{Eq:dSN-QPpoiss}
    \mult\circ \dSN{P,Q} 
=&\sum_{i=1,2} \sum_{j=1}^n (-1)^{j+i} 
 Q_1\ldots Q_{j-1} \delta^i_j P_{i+1} d^i_j Q_{j+1}\ldots Q_n \\
&-\sum_{i=1,2} \sum_{j=1}^n (-1)^{j+i} 
 Q_1\ldots Q_{j-1} c^i_j P_{i+1} \gamma^i_j Q_{j+1}\ldots Q_n\,,
\end{aligned}  
\end{equation}
with the index for $P_{i+1}$ understood modulo $2$ in $\{1,2\}$. 

\section{Double Poisson cohomology, after Pichereau and Van de Weyer} \label{ss:dPcoh}

Recall from Theorem \ref{Thm:dJac-PP} that, if $P\in (\mb T^\ast \cA)_2$ and $\brSN{P,P}=0$ in $(\mb T^\ast \cA)_\sharp$, then the double bracket $\dgal{-,-}_P$ defined through \eqref{Eq:dbr-P} is Poisson.
Furthermore, by Proposition \ref{Prop:PVdW-diff}, $P$ induces the square-zero differential $\dd_P$ \eqref{Eq:dP-PVdW} on $(\mb T^\ast \cA)_\sharp$. 
This motivates the following, which is due to Pichereau and Van de Weyer \cite{PV,VdW}. 

\begin{definition} \label{Def:dPH}
The cohomology of the complex $((\mb T^\ast \cA)_\sharp , \dd_P)$ is called the \emph{double Poisson cohomology} of $\cA$ with respect to $P\in (\mb T^\ast \cA)_2$. 
Explicitly, we have $\dPH(\cA)=(\dPH^k(\cA))_{k \geq 0}$\glslink{dPH-PV}{} for 
\begin{equation}
\dPH^k(\cA)=\frac{\ker\, \dd_P:(\mb T^\ast \cA)_{\sharp,k}\to (\mb T^\ast \cA)_{\sharp,k+1}}{\im \, \dd_P:(\mb T^\ast \cA)_{\sharp,k-1}\to (\mb T^\ast \cA)_{\sharp,k}}\,.
\end{equation}
\end{definition}

In general, to evaluate $\dd_P$ on some $Q\in (\mb T^\ast \cA)_{\sharp,n}$ 
one can use \eqref{Eq:dSN-QPpoiss} and extend this equality linearly. 

%Similarly to Proposition \ref{Pr:gdPA-gLie} in the graded setting, 
Let us now introduce some objects necessary for describing the first two cohomology groups. 
If $\dgal{-,-}_P$ is the double Poisson bracket associated with $P\in (\mb T^\ast \cA)_2$, let
\begin{align}  \glslink{ZPA}{} \glslink{ZPAA}{}
Z_P(\cA_\sharp;\cA):=&\{\bar{a}\in \cA_\sharp \mid \br{\bar{a},b}_P =0 \,\, \forall b\in \cA\}\,, \\
    Z_P(\cA_\sharp):=&\{\bar{a}\in \cA_\sharp \mid \br{\bar{a},\bar{b}}_{P,\sharp}=0 \,\, \forall \bar{b}\in \cA_\sharp\}\,,
\end{align}
where the operations $\br{-,-}_P$ and $\br{-,-}_{P,\sharp}$ are defined from $\dgal{-,-}_P$ using \eqref{Eq:mult-br-1} and \eqref{Eq:mult-brLie}, respectively. 
Hence $Z_P(\cA_\sharp)$ is the center of the Lie algebra $(\cA_\sharp,\br{-,-}_{P,\sharp})$ and $Z_P(\cA_\sharp;\cA)\subset Z_P(\cA_\sharp)$ where the inclusion is proper e.g. if there exists some $\bar{a}\in \cA_\sharp$ such that $\br{\bar{a},-}$ is a nonzero inner derivation on $\cA$. 

A \emph{double Poisson derivation} is an element $\delta\in \DDer_B(\cA)/[\cA,\DDer_B(\cA)]$ such that $\brSN{P,\delta}=0$. 
We say that a \emph{double Hamiltonian derivation} is an element $\delta_f\in \DDer_B(\cA)/[\cA,\DDer_B(\cA)]$ such that $\delta_f=\brSN{P,f}$ for some $f\in \cA$.

\begin{lemma} [\cite{PV}, Section 2]  \label{Lem:PV-1}
The first two cohomology groups of $\dPH(\cA)$ are such that  
\begin{align}
\dPH^0(\cA)& \subseteq Z_P(\cA_\sharp;\cA)\,, \label{Eq:Lem-HP0} \\
\dPH^1(\cA)&=\frac{\{\text{double Poisson derivations}\}}{\{\text{double Hamiltonian derivations}\}}\,.
\label{Eq:Lem-HP1}
\end{align}
Furthermore,  \eqref{Eq:Lem-HP0} is an equality if $\mu_1:(\mb T^\ast \cA)_{\sharp,1}\to \Der_B(\cA)$ is injective. 
\end{lemma}
\begin{proof}
By definition, $\dPH^0(\cA)=\ker\, \dd_P:(\mb T^\ast \cA)_{\sharp,0}\to (\mb T^\ast \cA)_{\sharp,1}$. 
Without loss of generality, we can write $P=\delta_1\delta_2$ and we have for an arbitrary lift $a\in \cA$ of $\bar{a}\in \cA_\sharp$, 
\begin{equation}
\begin{aligned} \label{Eq:dP-0}
    \mult \circ \dSN{P,a}&=
    \mult \circ (\delta_1\ast \dSN{\delta_2,a} - \dSN{\delta_1,a} \ast \delta_2) \\
    &=\delta_2(a)'\delta_1 \delta_2(a)'' - \delta_1(a)' \delta_2 \delta_1(a)''\,,
\end{aligned}    
\end{equation}
where we used \eqref{Eq:gdb-Lleib} and \eqref{Eq:dSN-b}. 
Thus $\dd_P(\bar{a})$ is the equivalence class of \eqref{Eq:dP-0} in $(\mb T^\ast \cA)_{\sharp,1}$. Meanwhile, we observe from \eqref{Eq:dbr-P} that the double derivation \eqref{Eq:dP-0} applied to $b\in \cA$ equals $-\dgal{a,b}_P$. 
If $\bar{a}\in \ker \dd_P$, the double derivation $\dgal{a,-}_P$ vanishes modulo commutators, hence $\mu_1(\dd_P(\bar{a}))=-\br{\bar{a},-}_P=0$.  
Conversely,  we have $\br{\bar{a},-}_P=0$ by definition if $\bar{a}\in Z_P(\cA_\sharp;\cA)$; by the previous computations $\dd_P(\bar{a})=0$ if $\mu_1$ is injective.  

The statement for $\dPH^1(\cA)$ is direct by definition. 
\end{proof}

\begin{remark} \label{Rem:PV}
It is claimed in \cite[\S2]{PV} that $\dPH^0(\cA)= Z_P(\cA_\sharp)$, but the proof contains a typo: they consider $\mu_1(\dd_P(\bar{a}))$ to be the zero derivation modulo commutators, instead of requiring  $\dd_P(\bar{a})$ to be the zero double derivation modulo commutators. 
Let us give explicitly an example where $\dPH^0(\cA)\subsetneq Z_P(\cA_\sharp)$. The double bracket on $\cA:=\kk\langle x,y \rangle$ satisfying
\begin{equation} \label{Eq:dbP-DSKV}
    \dgal{x,x}= y\otimes x - x\otimes y\,, \quad 
\dgal{x,y}=0,\,\, \dgal{y,y}=0,
\end{equation}
is Poisson, see \cite[\S2.1]{ORS} or \cite[\S2.6]{DSKV}. It corresponds to $P=y \partial_x x \partial_x$ where $\partial_x\in \DDer(\cA)$ satisfies   $\partial_x(x)=1\otimes 1$ and $\partial_x(y)=0$. Using \eqref{Eq:gdb-Lleib} and \eqref{Eq:dSN-a}--\eqref{Eq:dSN-b}, we find 
\begin{equation*}
\dSN{P,x}=y \partial_x x \star \dSN{\partial_x,x} 
- y \star \dSN{\partial_x,x} \star x \partial_x
=1\otimes y \partial_x x - x \partial_x \otimes y\,.
\end{equation*}
Modulo commutators, we can write $\dd_P(x)=\brSN{P,x}=(xy-yx)\partial_x$, thus $x \notin \ker \dd_P$. (We also denote by $x$ the class of $x\in \cA$ inside $\cA_\sharp=\kk[x,y]$.) 
Meanwhile, we have $\br{x,a}_P=[y,a]$ for all $a\in \cA$ because this holds on $a=x,y$ by applying the multiplication map to \eqref{Eq:dbP-DSKV}. 
Since $\br{x,-}_P$ is an inner derivation, $x\in Z_P(\cA_\sharp)$. 
In fact, $\cA$ is a free algebra so it satisfies the assumptions of Proposition \ref{Pr:MapMu-Iso} and therefore $\dPH^0(\cA)= Z_P(\cA_\sharp;\cA)$ by Lemma \ref{Lem:PV-1}. 
\end{remark}

It would be interesting to determine if the inclusion in \eqref{Eq:Lem-HP0} is always an equality or if it can be proper. 
For the latter case to happen, one needs $\mu_1$ \emph{not} injective. This is the case e.g. for $\cA=\kk[x]/(x^k)$ with $k\geq 2$, as Van den Bergh's gauge element  $\Delta\in \DDer(\cA)$, $\Delta:x\mapsto x\otimes 1 - 1\otimes x$ \cite{VdB1}, defines a nontrivial class in $(T^\ast \cA)_{\sharp,1}$ while $\mu_1(\Delta)=0$. (This contrasts with the polynomial case $\cA=\kk[x]$ where we can write $\Delta=[\partial_x, x]$ for $\partial_x\in \DDer(\cA)$, $\partial_x:x\mapsto  1\otimes 1$.) However, all the cases of double Poisson brackets on $\kk[x]/(x^k)$ that come from Powell's classification \cite{P16} on $\kk[x]$, cf. Example \ref{Ex:Classifkx}, do \emph{not}  correspond to an element of $(T^\ast \cA)_2$. Hence there is no notion of double Poisson cohomology in those cases according to the definition of Pichereau-Van de Weyer. 
This motivates the study of a more general cohomology theory which is defined in terms of double Poisson brackets directly. We start this investigation in the next chapter and we invite the interested reader to look at Subsection~\ref{ss:dP-coh-xr} where the theory can be applied to $\kk[x]/(x^k)$.

%%%%%%%%%%% NEW CHAPTER %%%%%%%%%%%%%%%
%%%%%%%%%%% NEW CHAPTER %%%%%%%%%%%%%%%
%%%%%%%%%%% NEW CHAPTER %%%%%%%%%%%%%%%
%%%%%%%%%%% NEW CHAPTER %%%%%%%%%%%%%%%

\chapter{Completed double Poisson cohomologies} \label{CH:Gen-dPcoh}

In the previous chapter, we considered an element $P\in (\mb T^\ast \cA)_2$ satisfying  $\brSN{P,P}=0$ in $(\mb T^\ast \cA)_\sharp$, which in turn defines a double Poisson bracket $\dgal{-,-}_P$ on $\cA$ through \eqref{Eq:dbr-P}. 
This allowed for a first construction of a cohomology theory due to Pichereau and Van de Weyer \cite{PV}, cf. Section \ref{ss:dPcoh}. 
However,  there exist $n$-brackets which are not obtained from an element $Q\in (\mb T^\ast \cA)_n$, see e.g. \cite[\S4.4]{VdB1}. 
Below, we introduce a general notion of double Poisson cohomology which does not rely on the existence of a noncommutative bivector defining the double Poisson bracket. 
%\footnote{We shall nevertheless denote our double bracket by $\dgal{-,-}_P$ and an arbitrary $n$-bracket by $\dgal{-}_Q$ \emph{as if} we were in the situation where they are induced by elements $P\in (T^\ast A)_2$ or $Q\in (T^\ast A)_n$. This notation is only chosen to avoid confusion between the two operations.} 

\section[The differential and definition of the cohomology]{The differential \texorpdfstring{$\wdd$}{} and definition of the cohomology}

Hereafter, we fix a $B$-linear double Poisson bracket\footnote{We use square brackets in order to emphasize that this double bracket is not necessarily of the form $\mu_2(P)$ for some $P\in (\mb T^\ast \cA)_2$, and to make it distinct from the $n$-bracket on which we apply the differential $\wdd$ of our complex.} $\dsq{-,-}\in \wBRA_B(\cA)_2$. 

\begin{definition} \label{def:wdd}
The linear map\glslink{wdd}{}  
\begin{equation} 
    \wdd:\wBRA_B(\cA)\to \oplus_{n\geq 0} \Hom_\kk(\cA^{\otimes (n+1)},\cA^{\otimes (n+1)})
\end{equation}
is defined in degree $0$ for any $\bar{a}\in \cA/[\cA,\cA]$ by 
\begin{equation}  \label{Eq:dP-gen-0}
    \wdd(\bar{a}):\cA\to \cA, \quad 
    \wdd(\bar{a}):=-\mult \circ \dsq{a,-}\,, \qquad 
(a\in \cA \text{ is a lift of }\bar{a}) 
\end{equation}
and in degree $n\geq 1$ for any $n$-bracket $\dgal{-}\in \wBRA_B(\cA)_n$ by 
\begin{equation} \label{Eq:dP-gen}
    \begin{aligned}
&\wdd(\dgal{-}):    \cA^{\otimes (n+1)}\to \cA^{\otimes (n+1)}\,, \\
&\wdd(\dgal{-}):=\sum_{s\in \Z_{n+1}} (-1)^{ns}\sigma^s \circ  
(\dgal{-}\otimes \Id_{\cA})\circ (\Id_{\cA}^{\otimes(n-1)}\otimes \dsq{-,-})
\circ \sigma^{-s} \\
&\qquad +(-1)^n\sum_{s\in \Z_{n+1}} (-1)^{ns}\sigma^s \circ  
(\dsq{-,-}\otimes \Id_{\cA}^{\otimes(n-1)}) \circ (\Id_{\cA}\otimes \dgal{-})
\circ \sigma^{-s} \,.
    \end{aligned}
\end{equation}
\end{definition}

\begin{example} \label{Ex:wddn0}
It is easy to see from \eqref{Eq:dP-gen-0} and the right Leibniz rule \eqref{Eq:db-Rleib} that 
$\wdd(\bar{a})\in \Der(\cA)$. One also has $\wdd(\bar{a})(b)=0$ for all $b\in B$ because $\dsq{-,-}$ is $B$-linear, so that $\wdd(\bar{a})\in \wBRA_B(\cA)_1$. 
In particular, we can consider $\wdd^2(\bar{a})$. 
Let us show that it vanishes as a consequence of the vanishing of \eqref{Eq:dJac} (this is the condition for being Poisson). 
Using the cyclic skewsymmetry \eqref{Eq:db-skew} of $\dsq{-,-}$, we have from \eqref{Eq:dP-gen}  
\begin{equation}
    \begin{aligned}
\wdd^2(\bar{a})=&(\wdd(\bar{a})\otimes \Id_{\cA} + \Id_{\cA} \otimes \wdd(\bar{a})) \circ \dsq{-,-} \\ 
&- \dsq{-,-} \circ (\wdd(\bar{a})\otimes \Id_{\cA} + \Id_{\cA} \otimes \wdd(\bar{a})) \,.
    \end{aligned}
\end{equation}
Combining this expression with \eqref{Eq:dP-gen-0} and evaluating on $b\otimes c\in \cA^{\otimes 2}$ gives 
\begin{align*}
    \wdd^2(\bar{a})(b\otimes c) =&
-\mult \circ     \dsq{a,\dsq{b,c}'}\otimes \dsq{b,c}''
-\dsq{b,c}'\otimes \mult \circ \dsq{a,\dsq{b,c}''} \\
&+\dsq{b,\mult \circ \dsq{a,c}} 
+ \dsq{\mult \circ \dsq{a,b},c} \,.
\end{align*}
Using further the cyclic skewsymmetry and the derivation rules for $\dsq{-,-}$ yields  
\begin{align*}
    \wdd^2(\bar{a})(b\otimes c) =&
-(\mult\otimes \Id_{\cA}) \circ     \dsq{a,\dsq{b,c}}_{L}
+(\Id_{\cA}\otimes \mult) \circ \sigma_{(123)} \dsq{a,\dsq{c,b}}_{L} \\
&-(\mult \otimes \Id_{\cA}) \circ \sigma_{(123)} \dsq{b,\dsq{c,a}}_{L} 
+(\Id_{\cA}\otimes \mult) \circ \dsq{b,\dsq{a,c}}_{L} \\
&-(\mult \otimes \Id_{\cA}) \circ \sigma_{(132)} \dsq{c,\dsq{a,b}}_{L}
+(\Id_{\cA}\otimes \mult) \circ \sigma_{(132)} \dsq{c,\dsq{b,a}}_{L} ,
\end{align*}
where we recall $\dsq{a,b\otimes c}_{L}=\dsq{a,b}\otimes c$, while for $\zeta\in S_3$ we denote by $\sigma_\zeta$ the corresponding permutation of factors on ${\cA}^{\otimes 3}$. 
Using the triple bracket defined through \eqref{Eq:dJac}, we can write 
\begin{equation} \label{Eq:dP2-abc}
\wdd^2(\bar{a})(b\otimes c) =
-(\mult\otimes \Id_{\cA}) \circ  \dsq{a,b,c}
+(\Id_{\cA}\otimes \mult) \circ \dsq{b,a,c}\,,
\end{equation}
which is zero since the triple bracket $\dsq{-,-,-}$ identically vanishes by definition of $\dsq{-,-}$ being Poisson. 
\end{example}

\subsection{Main statement}

\begin{theorem} \label{Thm:g-dPcoh1}
Consider the operation $\wdd$ from Definition \ref{def:wdd}. 
\begin{enumerate}[(1)]
    \item For any $n\geq 0$, $\wdd(\wBRA_B(\cA)_n)\subset  \wBRA_B(\cA)_{n+1}$. 

\noindent Hence, we can view $\wdd$ as a degree $+1$ endomorphism of $\wBRA_B(\cA)$. 

\item The linear map $\wdd$ is a square-zero differential on $\wBRA_B(\cA)$. 
\end{enumerate}
%This yields a linear map $\widehat{d}_P:\wBRA_B(A)_n\to \wBRA_B(A)_{n+1}$ for any $n\geq 0$. 
%Furthermore, the operations 
%$\widehat{d}_P:\wBRA_B(A)_\bullet \to \wBRA_B(A)_{\bullet+1}$ define a square-zero differential on the complex $\wBRA_B(A)$ if $\dgal{-,-}_P$ is a double Poisson bracket.  
\end{theorem}

\begin{remark}
The reader can check in the proof given below that the first item of  Theorem \ref{Thm:g-dPcoh1} is true for any double bracket $\dsq{-,-}$. The condition of being Poisson is only used for proving the square-zero property. 
\end{remark}

\subsection{Proof of Theorem \ref{Thm:g-dPcoh1}}
(1) We already proved this for $n=0$ in Example \ref{Ex:wddn0}, so we can take $n\geq 1$. 
The linear map $\wdd(\dgal{-})$ is clearly cyclically skewsymmetric by definition. Note also that if $\dsq{-,-}$ and $\dgal{-}$ are $B$-linear, then so too is $\wdd(\dgal{-})$. To conclude that  
$\wdd(\dgal{-})$ \eqref{Eq:dP-gen} is a $(n+1)$-bracket, it remains to verify the outer derivation property \eqref{Eq:nbr-DerOut}.  
Let us split \eqref{Eq:dP-gen} in terms of the two sums appearing on the right-hand side, which we call $D_1$ and $D_2$. We get from \eqref{Eq:nbr-DerAll} that $D_1(a_0,\ldots,a_{n-1},bc)$ equals 
\begin{align*}
&\sum_{2\leq s\leq n} (-1)^{ns} \sigma^s \circ 
\left[ \Big(b \star_{n+1-s} \dgal{a_0,\ldots,c,\ldots,a_{s-3},\dsq{a_{s-2},a_{s-1}}' }\Big) \otimes \dsq{a_{s-2},a_{s-1}}'' \right] \\
&+\sum_{2\leq s\leq n} (-1)^{ns} \sigma^s \circ 
\left[ \Big(\dgal{a_0,\ldots,b,\ldots,a_{s-3},\dsq{a_{s-2},a_{s-1}}' }\star_{s-1} c\Big) \otimes \dsq{a_{s-2},a_{s-1}}'' \right] \\
 &+ \dgal{a_0,\ldots,a_{n-2},b\dsq{a_{n-1},c}' } \otimes \dsq{a_{n-1},c}'' \\
& +\dgal{a_0,\ldots,a_{n-2},\dsq{a_{n-1},b}' } \otimes \dsq{a_{n-1},b}''c \\
&+(-1)^n\sigma \circ\dgal{a_1,\ldots,a_{n-1},\dsq{b,a_0}' c } \otimes \dsq{b,a_0}'' \\
&+(-1)^n\sigma \circ\dgal{a_1,\ldots,a_{n-1},\dsq{c,a_0}' } \otimes b\dsq{c,a_0}'' \\
&=b D_1(a_0,\ldots,a_{n-1},c)+ D_1(a_0,\ldots,a_{n-1},b)c \\
&\quad + \dgal{a_0,\ldots,a_{n-2},b } \cdot  \dsq{a_{n-1},c}
+(-1)^n   \dsq{b,a_0}^\sigma \cdot  \dgal{a_1,\ldots,a_{n-1},c}\,, 
\end{align*}
where we used for the last two terms the multiplication \eqref{prod-nm}.   
Similarly, we compute that $D_2(a_0,\ldots,a_{n-1},bc)$ equals  
\begin{align*}
&b D_2(a_0,\ldots,a_{n-1},c)+ D_2(a_0,\ldots,a_{n-1},b)c \\
&+(-1)^n (\sigma^{n-1}\dgal{b,a_0,\ldots,a_{n-2}} )\cdot  \dsq{a_{n-1},c}
+(-1)^n   \dsq{a_0,b} \cdot  \dgal{a_1,\ldots,a_{n-1},c} \\
=&b D_2(a_0,\ldots,a_{n-1}c)+ D_2(a_0,\ldots,a_{n-1}b)c \\
&-\dgal{a_0,\ldots,a_{n-2},b } \cdot  \dsq{a_{n-1},c}
+(-1)^n   \dsq{a_0,b} \cdot  \dgal{a_1,\ldots,a_{n-1},c}\,,    
\end{align*}
after using \eqref{Eq:nbr-Cycl}. Summing both terms, we obtain the desired identity: 
$$\wdd(\dgal{-})(a_0,\ldots,a_{n-1},bc)
=b\, \wdd(\dgal{-})(a_0,\ldots,a_{n-1},c)+ \wdd(\dgal{-})(a_0,\ldots,a_{n-1},b)\,c\,.$$ 

\medskip

(2) We already proved the square-zero property for $n=0$ in Example \ref{Ex:wddn0}, so we can take $n\geq 1$. 
We need to prove the vanishing of $\wdd^2(\dgal{-})$, which is a $(n+2)$-bracket.  
We shall use the decomposition 
\begin{equation} \label{Eq:dP-pf1}
\wdd^2(\dgal{-})=
\sum_{t\in \Z_{n+2}}\sum_{s\in \Z_{n+1}} (-1)^{(n+1)t+ns}
\sigma^t \circ (C_1^{(s)}+C_2^{(s)}+C_3^{(s)}+C_4^{(s)})\circ  \sigma^{-t}\,, 
\end{equation}
where we recall $\sigma:=\sigma_{(1,2,\ldots,n+2)}$ and we set 
\begin{subequations}
\begin{align}
C_1^{(s)}=& \sigma^s_{(1,\ldots,n+1)} \circ (\dgal{-} \otimes \Id_{\cA}^{\otimes 2}) \circ 
(\Id_{\cA}^{\otimes (n-1)}\otimes \dsq{-,-}\otimes \Id_{\cA}) \nonumber \\
&\qquad\circ \sigma_{(1,\ldots,n+1)}^{-s} \circ (\Id_{\cA}^{\otimes n}\otimes \dsq{-,-} ) 
\label{Eq:dP-pf2a} \\
C_2^{(s)}=&(-1)^n \,\, \sigma^s_{(1,\ldots,n+1)} \circ (\dsq{-,-} \otimes \Id_{\cA}^{\otimes n}) 
\circ (\Id_{\cA}\otimes \dgal{-} \otimes \Id_{\cA})  \nonumber \\
&\qquad\circ \sigma_{(1,\ldots,n+1)}^{-s} \circ (\Id_{\cA}^{\otimes n}\otimes \dsq{-,-} )
\label{Eq:dP-pf2b} \\
C_3^{(s)} =&(-1)^{n+1}\,\, ( \dsq{-,-} \otimes \Id_{\cA}^{\otimes n})\circ \sigma^s_{(2,\ldots,n+2)} \circ (\Id_{\cA}\otimes \dgal{-} \otimes \Id_{\cA})  \nonumber \\ 
&\qquad\circ (\Id_{\cA}^{\otimes n}\otimes \dsq{-,-})  \circ \sigma^{-s}_{(2,\ldots,n+2)}
\label{Eq:dP-pf2c} \\
C_4^{(s)} =& - \, ( \dsq{-,-}\otimes \Id_{\cA}^{\otimes n})\circ \sigma^s_{(2,\ldots,n+2)} \circ (\Id_{\cA}\otimes \dsq{-,-}\otimes \Id_{\cA}^{\otimes (n-1)}) \nonumber \\
&\qquad\circ (\Id_{\cA}^{\otimes 2}\otimes \dgal{-})    \circ \sigma^{-s}_{(2,\ldots,n+2)}
\label{Eq:dP-pf2d}
\end{align}
\end{subequations}
Let us first handle the case $n=1$. Then, we put $\delta:=\dgal{-}\in \Der_B(\cA)$ and we get 
for the sum over $C_1^{(s)}$ in \eqref{Eq:dP-pf1}: 
\begin{align*}
\sum_{t\in \Z_3} \sigma^t\circ (\delta \otimes \Id_{\cA}^{\otimes 2}+\Id_{\cA} \otimes \delta \otimes \Id_{\cA}) 
(\dsq{-,-} \otimes \Id_{\cA} ) \circ (\Id_{\cA} \otimes \dsq{-,-})\circ \sigma^{-t}, 
\end{align*}
where we used the cyclic skewsymmetry \eqref{Eq:db-skew} for $s=1$. 
Similarly, we get for the sum over $C_2^{(s)}$ : 
\begin{align*}
-\sum_{t\in \Z_3} \sigma^t \circ(\dsq{-,\delta(-)}\otimes \Id_{\cA} + \dsq{\delta(-),-} \otimes \Id_{\cA} ) \circ (\Id_{\cA} \otimes \dsq{-,-})\circ \sigma^{-t} , 
\end{align*}
for the sum over $C_3^{(s)}$ : 
\begin{align*}
+\sum_{t\in \Z_3} \sigma^t \circ(\dsq{-,\delta(-)}\otimes \Id_{\cA} + \dsq{-,-}\otimes \delta) \circ (\Id_{\cA} \otimes \dsq{-,-}) \circ \sigma^{-t} ,
\end{align*}
and for the sum over $C_4^{(s)}$ : 
\begin{align*}
-\sum_{t\in \Z_3} \sigma^t \circ(\dsq{-,-}\otimes \Id_{\cA} ) \circ (\Id_{\cA} \otimes \dsq{-,\delta(-)} + \Id_{\cA} \otimes \dsq{\delta(-),-}) \circ\sigma^{-t}  .
\end{align*}
The first terms in the sums  over $C_2^{(s)}$ and $C_3^{(s)}$ cancel out. 
Using the extension \eqref{mfold-ext} of $\delta$ to $\cA^{\otimes 3}$, we obtain
\begin{align*}
\wdd^2(\delta) =&
\delta\circ \sum_{t\in \Z_3} \sigma^t \circ
(\dsq{-,-}\otimes \Id_{\cA} ) \circ (\Id_{\cA} \otimes \dsq{-,-}) \circ \sigma^{-t} \\
&-\sum_{t\in \Z_3} \sigma^t  \circ
(\dsq{-,-}\otimes \Id_{\cA} ) \circ (\Id_{\cA} \otimes \dsq{-,-})\circ  \sigma^{-t} \circ \delta\,,
\end{align*}
and we conclude because the double bracket is Poisson, cf. \eqref{Eq:dJacExpl}. 

We now have to show the vanishing of \eqref{Eq:dP-pf1} for $n\geq 2$ using a similar reasoning. 
We start by splitting each linear map $\sum_{s\in \Z_{n+1}} (-1)^{ns} C_j^{(s)}$, $1\leq j \leq 4$, into 3 convenient endomorphisms of $\cA^{\otimes (n+2)}$.  
We will repeatedly use the extensions \eqref{20240805:eq1} of the linear maps 
$\dsq{b,-}:\cA\to \cA^{\otimes 2}$ and 
$\dgal{b_1,\ldots,b_{s-1},-,b_{s+1},\ldots,b_n}:\cA\to \cA^{\otimes n}$ for any $b,b_j\in \cA$.

%%%% j=1
For $j=1$, we use \eqref{Eq:dP-pf2a} to get that $\sum_{s\in \Z_{n+1}} (-1)^{ns} C_1^{(s)}(a_1,\ldots,a_{n+2})$ equals 
\begin{subequations}
\begin{align}
&\dgal{a_1,\ldots,a_{n-1},-}_L  ( \dsq{a_n,\dsq{a_{n+1},a_{n+2}}}_{L})
\label{Eq:l-C1-0}    \\
+&(-1)^n \sigma_{(1,\ldots,n+1)} \!\circ\! \dgal{a_2,\ldots,a_{n},-}_L \!
 \left(\dsq{\dsq{a_{n+1},a_{n+2}}',a_1}\!\otimes\! \dsq{a_{n+1},a_{n+2}}'' \right) 
\label{Eq:l-C1-1}    \\
+&\sum_{s=2}^n (-1)^{ns} \sigma_{(1,\ldots,n+1)}^s \circ 
\Big(\dgal{a_{s+1},\ldots,a_n,\dsq{a_{n+1},a_{n+2}}',a_1,\ldots,a_{s-2},\dsq{a_{s-1},a_{s}}'} \nonumber \\
& \hspace{3cm} \otimes  \dsq{a_{s-1},a_{s}}'' \otimes \dsq{a_{n+1},a_{n+2}}'' \Big) \,. \label{Eq:l-C1-out} 
\end{align}
\end{subequations} 
We respectively denote by $C_{1,0}$, $C_{1,1}$ and $C_{1,\geq 2}$ the elements of $\End(\cA^{\otimes (n+2)})$ which return \eqref{Eq:l-C1-0}, \eqref{Eq:l-C1-1} and \eqref{Eq:l-C1-out} when evaluated on $a_1\otimes \ldots \otimes a_{n+2}$. 
Thus, $\sum_{s\in \Z_{n+1}} (-1)^{ns} C_1^{(s)}=C_{1,0}+C_{1,1}+C_{1,\geq2}$. 

%%%% j=2
For $j=2$, we use \eqref{Eq:dP-pf2b} to write $\sum_{s\in \Z_{n+1}} (-1)^{ns} C_2^{(s)}(a_1,\ldots,a_{n+2})$ as 
\begin{subequations}
\begin{align}
&(-1)^n \, \dsq{a_1,-}_L \circ \dgal{a_2,\ldots,a_n,-}_L
( \dsq{a_{n+1},a_{n+2}} )    \label{Eq:l-C2-0}   \\
+&\sum_{s=2}^n (-1)^{ns} \sigma_{(1,\ldots,n+1)}^{s-1} \circ 
\dsq{a_s,-}_L   \nonumber \\
&\hspace{1.2cm}\circ\dgal{a_{s+1},\ldots,a_n,-,a_1,\ldots,a_{s-1}}_L  (\dsq{a_{n+1},a_{n+2}}) 
\label{Eq:l-C2-mid}   \\
+& \sigma_{(1,\ldots,n+1)}^{-1} \circ \dsq{\dsq{a_{n+1},a_{n+2}}',-}_L  
\left(\dgal{a_1,\ldots,a_n}\otimes \dsq{a_{n+1},a_{n+2}}'' \right).
 \label{Eq:l-C2-n}  
\end{align}
\end{subequations}
We respectively denote by $C_{2,0}$, $C_{2,mid}$ and $C_{2,n}$ the elements of $\End(\cA^{\otimes (n+2)})$ which return \eqref{Eq:l-C2-0}, \eqref{Eq:l-C2-mid} and \eqref{Eq:l-C2-n} when evaluated on $a_1\otimes \ldots \otimes a_{n+2}$. 

%%%% j=3
For $j=3$, we use \eqref{Eq:dP-pf2c} to write $\sum_{s\in \Z_{n+1}} (-1)^{ns} C_3^{(s)}(a_1,\ldots,a_{n+2})$ as 
\begin{subequations}
\begin{align} 
-&(-1)^n \dsq{a_1,-}_L \circ 
\dgal{a_2,\ldots,a_n,-}_L  (\dsq{a_{n+1},a_{n+2}})
    \label{Eq:l-C3-0}   \\
-&  \dsq{a_1,-}_L  \circ \sigma \circ   \dgal{a_3,\ldots,a_{n+1},-} (\dsq{a_{n+2},a_{2}})
\label{Eq:l-C3-1}   \\
-&\sum_{s=2}^n (-1)^{n(s+1)}  \dsq{a_1,-}_L  \circ \sigma^s \circ \nonumber \\
&\qquad \qquad \circ  \dgal{a_{s+2},\ldots,a_{n+2},a_2,\ldots,a_{s-1},-}_L 
(\dsq{a_{s},a_{s+1}})\,,
 \label{Eq:l-C3-out}  
\end{align}
\end{subequations}
where  we use the permutation $\sigma=\sigma_{(1,\ldots,n+1)}$ on $\cA^{\otimes (n+1)}$ in \eqref{Eq:l-C3-out}. 
We respectively denote by $C_{3,0}$, $C_{3,1}$ and $C_{3,\geq2}$ the elements of $\End(\cA^{\otimes (n+2)})$ which return \eqref{Eq:l-C3-0}, \eqref{Eq:l-C3-1} and \eqref{Eq:l-C3-out} when evaluated on $a_1\otimes \ldots \otimes a_{n+2}$. 

%%%% j=4
Finally for $j=4$, we use \eqref{Eq:dP-pf2d} to write $\sum_{s\in \Z_{n+1}} (-1)^{ns} C_4^{(s)}(a_1,\ldots,a_{n+2})$ as 
\begin{subequations}
\begin{align}
-&  \dsq{a_1,-}_{L} \circ \dsq{a_2,-}_L (\dgal{a_3,\ldots,a_{n+2}} )
\label{Eq:l-C4-0}   \\
-&\sum_{s=1}^{n-1}(-1)^{ns} 
\dsq{a_1,-}_{L} \circ \sigma^s 
\circ \dsq{a_{s+2},-}_L  \left(\dgal{a_{s+3},\ldots,a_{n+2},a_2,\ldots,a_{s+1}}\right) 
\label{Eq:l-C4-mid}   \\
-&(-1)^{n} \dsq{a_1,-}_{L} \circ \sigma^{-1} \circ  \dsq{a_{n+2},-}_L   
\left(\dgal{a_{2},\ldots,a_{n+1}} \right) \,, \label{Eq:l-C4-n}
\end{align}
\end{subequations}
with $\sigma=\sigma_{(1,\ldots,n+1)}$ acting on $\cA^{\otimes (n+1)}$ in \eqref{Eq:l-C4-mid}-\eqref{Eq:l-C4-n}. 
We respectively denote by $C_{4,0}$, $C_{4,mid}$ and $C_{4,n}$ the elements of $\End(\cA^{\otimes (n+2)})$ which return \eqref{Eq:l-C4-0}, \eqref{Eq:l-C4-mid} and \eqref{Eq:l-C4-n} when evaluated on $a_1\otimes \ldots \otimes a_{n+2}$.

It is direct to see that $C_{2,0}+C_{3,0}=0$ from \eqref{Eq:l-C2-0} and \eqref{Eq:l-C3-0}. 
For the other endomorphisms $C_{j,k}$, one needs to consider them under the sum over $t \in \Z_{n+2}$ \eqref{Eq:dP-pf1} in order to get further cancellations and claim that $\wdd_P^2(\dgal{-})=0$. 

\medskip 

%%%%%%%%%%% Step i 
\underline{Step i:} We prove the equality 
\begin{equation} 
\label{Eq:dP-pf-step1}
\sum_{t\in \Z_{n+2}}(-1)^{(n+1)t} \, \sigma^t \circ [C_{1,0}+C_{1,1}+C_{3,1}] \circ \sigma^{-t}=0\,.
\end{equation} 
\noindent Owing to the sum over $t\in \Z_{n+2}$, it is equivalent to show this identity with $C_{1,1}$ replaced by $\widetilde{C}_{1,1}:=(-1)^{n+1}\sigma^{-1} \circ C_{1,1} \circ \sigma$, and 
$C_{3,1}$ replaced by $\widetilde{C}_{3,1}:=\sigma^{-2} \circ C_{3,1} \circ \sigma^2$. 
Using \eqref{Eq:l-C1-1},  $\widetilde{C}_{1,1}(a_1,\ldots,a_{n+2})$ equals 
\begin{align}
&-\sigma^{-1}\circ \sigma_{(1,\ldots,n+1)}\circ \dgal{a_1,\ldots,a_{n-1},-}_L  
 \left(\dsq{\dsq{a_{n},a_{n+1}}',a_{n+2}}\otimes \dsq{a_{n},a_{n+1}}'' \right) \nonumber  \\
&= \dgal{a_1,\ldots,a_{n-1},-}_L  
\left(\sigma_{(132)} \dsq{a_{n+2} , \dsq{a_{n},a_{n+1} } }_{L}\right)\,,
\label{Eq:l-C1-1-tw} 
    \end{align}
due to $\sigma^{-1}\circ \sigma_{(1,\ldots,n+1)}=\sigma_{(n+1,n+2)}$ and the cyclic skewsymmetry of $\dsq{-,-}$.  
Using \eqref{Eq:l-C3-1}, $\widetilde{C}_{3,1}(a_1,\ldots,a_{n+2})$ equals 
\begin{equation}
\begin{aligned}
&- \dgal{a_1,\ldots,a_{n-1},\dsq{a_{n},a_{n+2}}'} \otimes \dsq{a_{n+1},\dsq{a_{n},a_{n+2}}''} \\
=& \dgal{a_1,\ldots,a_{n-1},-}_L  
\left(\sigma_{(123)} \dsq{a_{n+1} , \dsq{a_{n+2},a_n } }_{L}\right)\,.
\label{Eq:l-C3-1-tw} 
\end{aligned}    
\end{equation}
Combining these expressions with \eqref{Eq:l-C1-0}, 
$[C_{1,0}+\widetilde{C}_{1,1}+\widetilde{C}_{3,1}](a_1,\ldots,a_{n+2})$ 
can be written as the composition of $ \dgal{a_1,\ldots,a_{n-1},-}_L$ with 
\begin{align*}
\dsq{a_n, \dsq{a_{n+1}, a_{n+2}} }_{L}
+\sigma_{(123)} \dsq{a_{n+1} , \dsq{a_{n+2},a_n} }_{L} 
+ \sigma_{(132)} \dsq{a_{n+2} , \dsq{a_{n},a_{n+1}} }_{L} \,,
\end{align*}
which vanishes as $\dsq{-,-}$ is Poisson. Thus $C_{1,0}+\widetilde{C}_{1,1}+\widetilde{C}_{3,1}=0$, and  \eqref{Eq:dP-pf-step1} holds. 

\medskip

%%%%%%%%%%% Step ii 
\underline{Step ii:} We prove the equality 
\begin{equation} 
\label{Eq:dP-pf-step2}
 \sum_{t\in \Z_{n+2}}(-1)^{(n+1)t} \, \sigma^t \circ [C_{2,n}+C_{4,0}+C_{4,n}] \circ \sigma^{-t}=0\,.
\end{equation} 
\noindent Similarly to step i, in the sum we can replace $C_{2,n}$ by $\widetilde{C}_{2,n}:=\sigma^2 \circ C_{2,n} \circ \sigma^{-2}$ and replace $C_{4,n}$ by $\widetilde{C}_{4,n}:=(-1)^{n+1} \sigma \circ C_{4,n} \circ \sigma^{-1}$.  
Using \eqref{Eq:l-C2-n}, $\widetilde{C}_{2,n}(a_1,\ldots,a_{n+2})$ equals 
\begin{equation}
    \begin{aligned}
&\sigma^2 \circ \sigma_{(1,\ldots,n+1)}^{-1} \circ 
\dsq{\dsq{a_{1},a_{2}}',-}_L
\left(\dgal{a_3,\ldots,a_{n+2}} \otimes \dsq{a_{1},a_{2}}'' \right) \\
=&
-\left(\sigma_{(132)}\dsq{-,\dsq{a_{1},a_{2}}}_{L}\otimes \Id_{\cA}^{\otimes (n-1)} \right) 
\left(\dgal{a_3,\ldots,a_{n+2}} \right) \,,
 \label{Eq:l-C2-n-tw}      
    \end{aligned}
\end{equation}
due to $\sigma^2 \circ \sigma_{(1,\ldots,n+1)}^{-1}=\sigma_{(2,\ldots,n+2)}$
and cyclic skewsymmetry.
Next, we remark that for any $b_1,\ldots,b_n\in \cA$ 
 \begin{align*}
&\sigma \circ \dsq{a_2,-}_L  \circ \sigma^{-1} \circ  \dsq{a_{1},-}_L  (b_1\otimes \ldots \otimes b_n) \\
=& \dsq{a_1,b_1}'\otimes \dsq{a_2, \dsq{a_1,b_1}''} \otimes b_2 \otimes \ldots \otimes b_n \\
=&-\left(\sigma_{(123)}\dsq{a_{2},\dsq{-,a_{1}} }_{L}\otimes \Id_{\cA}^{\otimes (n-1)} \right)  (b_1\otimes \ldots \otimes b_n)\,.
 \end{align*}
(In the first line, notice that $\sigma=\sigma_{(1,\ldots,n+2)}$ acts on $\cA^{\otimes (n+2)}$ and $\sigma^{-1}=\sigma^{-1}_{(1,\ldots,n+1)}$ acts on $\cA^{\otimes (n+1)}$.) 
Using this with \eqref{Eq:l-C4-n}, $\widetilde{C}_{4,n}(a_1,\ldots,a_{n+2})$ equals 
\begin{equation}
-\left(\sigma_{(123)}\dsq{a_{2},\dsq{-,a_{1}}}_{L}\otimes \Id_{\cA}^{\otimes (n-1)} \right)  
( \dgal{a_{3},\ldots,a_{n+2}}) \,.
\label{Eq:l-C4-n-tw}
\end{equation}
Combining \eqref{Eq:l-C2-n-tw} and \eqref{Eq:l-C4-n-tw} with \eqref{Eq:l-C4-0},  
$[\widetilde{C}_{2,n}+C_{4,0}+\widetilde{C}_{4,n}](a_1,\ldots,a_{n+2})$ is obtained by applying the linear map 
\begin{equation} \label{Eq:l-Stepii}
-\left[\dsq{a_{1},\dsq{a_{2},-}}_{L}
+ \sigma_{(123)}\dsq{a_{2},\dsq{-,a_{1}}}_{L}
+\sigma_{(132)}\dsq{-,\dsq{a_{1},a_{2}}}_{L} \right] 
\otimes \Id_{\mc A}^{\otimes (n-1)}    
\end{equation}
to $\dgal{a_{3},\ldots,a_{n+2}} \in \cA^{\otimes n}$. 
Since $\dsq{-,-}$ is Poisson, \eqref{Eq:l-Stepii} vanishes identically, hence $\widetilde{C}_{2,n}+C_{4,0}+\widetilde{C}_{4,n}=0$ and \eqref{Eq:dP-pf-step2} holds.

\medskip 

%%%%%%%%%%% Step iii 
\underline{Step iii:} We prove the equality 
\begin{equation} 
\label{Eq:dP-pf-step3}
\sum_{t\in \Z_{n+2}}(-1)^{(n+1)t} \, \sigma^t \circ [C_{2,mid}+C_{3,\geq 2}] \circ \sigma^{-t}=0\,.
\end{equation} 
We first relabel $s$ by $n+2-s$ in \eqref{Eq:l-C3-out}, so that 
$C_{3,\geq 2}(a_1,\ldots,a_{n+2})$ becomes 
\begin{equation}
    \begin{aligned} \label{Eq:dP-pf-3a}
-&\sum_{s=2}^n (-1)^{ns}  \dsq{a_1,-}_L \circ \sigma^{1-s} \circ  \\
&\qquad  \dgal{a_{n+4-s},\ldots,a_{n+2},a_2,\ldots,a_{n+1-s},-}_L (\dsq{a_{n+2-s},a_{n+3-s}})\,.
\end{aligned}
\end{equation}
For $2\leq s \leq n$, we label by $C_{3,\geq 2}^{(s)}$ the endomorphism of $\cA^{\otimes (n+2)}$ which returns the $s$-th summand of \eqref{Eq:dP-pf-3a} when applied to $a_1 \otimes \ldots \otimes a_{n+2}$. 
We can then write 
$$\sum_{t\in \Z_{n+2}}(-1)^{(n+1)t} \sigma^t \circ C_{3,\geq 2} \circ \sigma^{-t}
%=\sum_{t\in \Z_{n+2}} (-1)^{(n+1)t} \sigma^t \circ \sum_{s=2}^n C^{(s)}_{3,\geq 2} \circ \sigma^{-t}
=\sum_{t\in \Z_{n+2}} (-1)^{(n+1)t} \sigma^t \circ \sum_{s=2}^n \widetilde{C}^{(s)}_{3,\geq 2} \circ \sigma^{-t}$$ after introducing  
$\widetilde{C}^{(s)}_{3,\geq 2}=(-1)^{(n+1)(s-1)} \sigma^{s-1} \circ C_{3,\geq 2}^{(s)}\circ \sigma^{-(s-1)}$. 
With this new endomorphism, we see that $\widetilde{C}^{(s)}_{3,\geq 2}(a_1\otimes \ldots \otimes a_{n+2})$ equals
\begin{align*}
&- (-1)^{ns} (-1)^{(n+1)(s-1)} \sigma^{s-1} \circ \dsq{a_s,-}_L \circ \sigma^{1-s} \circ \nonumber \\
&\qquad  \circ\dgal{a_1,\ldots,a_{s-1},a_{s+1},\ldots,a_{n},-}_L ( \dsq{a_{n+1},a_{n+2}}) \\
=&- (-1)^{ns} \sigma^{s-1} \circ  \dsq{a_s,-}_L \circ \sigma^{1-s} \circ \sigma^{s-1}_{(1,\ldots,n)} \circ \nonumber \\
&\qquad  \circ \dgal{a_{s+1},\ldots,a_{n},-,a_1,\ldots,a_{s-1}}_L (\dsq{a_{n+1},a_{n+2}})
\end{align*}
where in the second line we used $\dgal{-}=(-1)^{(n+1)(s-1)} \sigma^{s-1}\circ \dgal{-} \circ \sigma^{-(s-1)}$ by cyclic skewsymmetry \eqref{Eq:nbr-Cycl}. 
Let us emphasize that the permutation composed after $\dsq{a_s,-}_L$ is $\sigma^{s-1}=\sigma^{s-1}_{(1,\ldots,n+2)}$ acting on $\cA^{\otimes (n+2)}$, while the permutation composed before it is $\sigma^{1-s}=\sigma^{1-s}_{(1,\ldots,n+1)}$ acting on $\cA^{\otimes (n+1)}$. 
Since $2\leq s \leq n$, we observe  
\begin{equation*}
\sigma^{s-1} \circ \dsq{a_s,-}_L 
=  
\sigma^{s-1}_{(1,\ldots,n+1)} \circ \dsq{a_s,-}_L \circ \sigma^{-1}_{(n+2-s,\ldots,n+1)}\,,
\end{equation*}
as an equality in $\Hom(\cA^{\otimes (n+1)},\cA^{\otimes (n+2)})$, and 
$\sigma^{1-s} \circ \sigma^{s-1}_{(1,\ldots,n)}=\sigma_{(n+2-s,\ldots,n+1)}$ on $\cA^{\otimes (n+1)}$. 
This finally allows us to write 
$\sum_{s=2}^n \widetilde{C}^{(s)}_{3,\geq 2}(a_1\otimes \ldots \otimes a_{n+2})$ as 
\begin{equation*} 
- \sum_{s=2}^n (-1)^{ns} \sigma^{s-1}_{(1,\ldots,n+1)} \circ   \dsq{a_s,-}_L \circ 
 \dgal{a_{s+1},\ldots,a_{n},-,a_1,\ldots,a_{s-1}}_L(\dsq{a_{n+1},a_{n+2}})
\end{equation*}
which is the opposite of \eqref{Eq:l-C2-mid}. We conclude that 
$ C_{2,mid}+\sum_{s=2}^n \widetilde{C}^{(s)}_{3,\geq 2}=0$, and therefore \eqref{Eq:dP-pf-step3} holds.

\medskip 

%%%%%%%%%%% Step iv
\underline{Step iv:} We prove the equality 
\begin{equation} 
\label{Eq:dP-pf-step4}
\sum_{t\in \Z_{n+2}}(-1)^{(n+1)t} \sigma^t \circ C_{1,\geq 2} \circ \sigma^{-t}=0\,.
\end{equation} 
Due to \eqref{Eq:l-C1-out}, we want the vanishing of the following sum 
\begin{align*}
&\sum_{t=0}^{n+1}(-1)^{(n+1)t+ns} 
\sum_{s=2}^n \sigma^t \circ \sigma_{(1,\ldots,n+1)}^s \circ 
 \nonumber \\
&\,\, \Big(\dgal{a_{s+t+1},\ldots,a_{t+n},\dsq{a_{t+n+1},a_{t+n+2}}',
a_{t+1},\ldots,a_{t+s-2},\dsq{a_{t+s-1},a_{t+s}}'} \\
&\qquad \qquad \qquad \otimes \dsq{a_{t+s-1},a_{t+s}}'' \otimes \dsq{a_{t+n+1},a_{t+n+2}}'' \Big)\,,
\end{align*}
where the index $j$ of $a_j$ is taken modulo $n+2$ in $\{1,\ldots,n+2\}$. All the summands can be gathered pairwise by groups of 2 terms only involving $\dgal{a_i,a_{i+1}}_P$ and $\dgal{a_{i+k},a_{i+k+1}}_P$ with $0\leq i<n$ and $2\leq k \leq n-i+1$. 
For such a fixed pair $(i,k)$, the two summands are obtained as follows. 
First, one takes $t=i+k+1$ and $s=n+2-k$  (so that $t+s-1=i$, $t+n+1=i+k$ as indices), which yields 
\begin{equation} \label{Eq:l-C1-ik1}
    \begin{aligned}
& (-1)^{(n+1)(i+1)+k+n} \, 
 \sigma^{i+k+1} \circ \sigma_{(1,\ldots,n+1)}^{n+2-k} \circ   \\
&\quad \Big(\dgal{a_{i+2},\ldots,a_{i+k-1},\dsq{a_{i+k},a_{i+k+1}}',
a_{i+k+2},\ldots,a_{i-1},\dsq{a_{i},a_{i+1}}'} \otimes \\
&\qquad \qquad \qquad \otimes \dsq{a_{i},a_{i+1}}'' \otimes \dsq{a_{i+k},a_{i+k+1}}'' \Big) 
\end{aligned}
\end{equation}
Second, one takes $t=i+1$ and $s=k$ (for $t+s-1=i+k$, $t+n+1=i$ as indices), which yields 
\begin{equation} \label{Eq:l-C1-ik2}
\begin{aligned}
&(-1)^{(n+1)(i+1)+nk} \, \sigma^{i+1} \circ \sigma_{(1,\ldots,n+1)}^k 
\circ  \\
&\quad  \Big(\dgal{a_{i+k+2},\ldots,a_{i-1},\dsq{a_{i},a_{i+1}}',
a_{i+2},\ldots,a_{i+k-1},\dsq{a_{i+k},a_{i+k+1}}'} \otimes \\
&\qquad \qquad \qquad \otimes \dsq{a_{i+k},a_{i+k+1}}'' \otimes \dsq{a_{i},a_{i+1}}'' \Big)\,,
\end{aligned}
\end{equation}
By cyclic skewsymmetry \eqref{Eq:nbr-Cycl}, we have 
$\dgal{-}=(-1)^{(n+1)(k-1)}\sigma^{-(k-1)}\circ \dgal{-}\circ \sigma^{k-1}$, so that we can rewrite \eqref{Eq:l-C1-ik2} as 
\begin{equation} \label{Eq:l-C1-ik2-b}
\begin{aligned}
&-(-1)^{(n+1)(i+1)+k+n} \, \sigma^{i+1} \circ \sigma_{(1,\ldots,n+1)}^k 
\circ \sigma_{(1,\ldots,n)}^{-(k-1)} \circ  \\
&\quad  \circ \Big(\dgal{a_{i+2},\ldots,a_{i+k-1},\dsq{a_{i+k},a_{i+k+1}}',a_{i+k+2},\ldots,a_{i-1},\dsq{a_{i},a_{i+1}}'} \otimes \\
&\qquad \qquad \qquad \otimes \dsq{a_{i+k},a_{i+k+1}}'' \otimes \dsq{a_{i},a_{i+1}}'' \Big)\,.
\end{aligned}
\end{equation}
We then see that $\eqref{Eq:l-C1-ik1}+\eqref{Eq:l-C1-ik2-b}=0$ provided that 
\begin{equation*}
    \sigma^k \circ \sigma_{(1,\ldots,n+1)}^{1-k}
= \sigma^k_{(1,\ldots,n+1)}\circ \sigma^{1-k}_{(1,\ldots,n)}\circ \sigma_{(n+1,n+2)}\,,
\end{equation*}
which is readily checked for $2\leq k <n+2$ as both side equal $\sigma_{(k,\ldots,n+2)}$. 
Therefore all the pairs for fixed $(i,k)$ cancel out, and \eqref{Eq:dP-pf-step4} is true. 

\medskip

%%%%%%%%%%% Step v
\underline{Step v:} We prove the equality 
\begin{equation} 
\label{Eq:dP-pf-step5}
\sum_{t\in \Z_{n+2}}(-1)^{(n+1)t} \sigma^t \circ C_{4,mid} \circ \sigma^{-t}=0\,.
\end{equation} 
From \eqref{Eq:l-C4-mid}, we are seeking for the vanishing of the following sum 
\begin{align*}
-\sum_{t=0}^{n+1} (-1)^{(n+1)t}\sum_{s=1}^{n-1}(-1)^{ns} & \sigma^t \circ 
\dsq{a_{t+1},-}_{L} \circ \sigma^s \circ  \dsq{a_{s+t+2},-}_L    \\
&  \left(\dgal{a_{s+t+3},\ldots,a_{t},a_{t+2},\ldots,a_{s+t+1}} \right)   \,,
\end{align*}
where the index of each $a_k$ is taken modulo $n+2$. 
We can gather all these summands by groups of 2 terms containing the operators $\dgal{a_i,-}_L$ 
and $\dgal{a_j,-}_L$ for some fixed $1\leq i <j\leq n+2$ such that $i-j \text{ mod }n+2\neq \pm 1$. 
For such a fixed pair $(i,j)$, these two terms are obtained as follows. 
First, one takes $t=j-1$ and $s=n+1+i-j$ (so that $t+1=j$ and $s+t+2=i$ as indices) which yields 
\begin{equation}
    \begin{aligned}  \label{Eq:l-C4-ij1}
& (-1)^{\nu_{ij}} \, \sigma^{j-1} \circ \dsq{a_{j},-}_L
\circ \sigma^{i-j}  \circ \dsq{a_{i},-}_L 
 \dgal{a_{i+1},\ldots,a_{j-1},a_{j+1},\ldots,a_{i-1}}  
    \end{aligned}
\end{equation}
with $\nu_{ij}=n(i-1)+j$. 
Second, one takes $t=i-1$ and $s=j-i-1$ (for $t+1=i$ and $s+t+2=j$ as indices) which yields 
\begin{equation}
    \begin{aligned} \label{Eq:l-C4-ij2}
& (-1)^{\mu_{ij}} \, \sigma^{i-1} \!\circ\! \dsq{a_{i},-}_L 
\!\circ \! \sigma^{j-i-1}  \!\circ \! \dsq{a_{j},-}_L 
\dgal{a_{j+1},\ldots,a_{i-1},a_{i+1},\ldots,a_{j-1}}  
    \end{aligned}
\end{equation}
with $\mu_{ij}=nj+i$. 
Since $i<i+1<j$, we note that 
\begin{align*}
&\sigma^{j-1} \circ \dsq{a_{j},-}_L
\circ \sigma^{i-j}  \circ \dsq{a_{i},-}_L (b_1\otimes \ldots \otimes b_n) \\
=& b_{n-i+2}\otimes \ldots \otimes b_n \otimes \dsq{a_i,b_1}  \otimes 
b_2 \otimes \ldots \otimes b_{j-i-1} \otimes \dsq{a_j,b_{j-i}} \otimes 
b_{j-i+1} \otimes \ldots \otimes b_{n-i+1} \\
=& \dsq{a_{i},-}_{(i)} \circ \dsq{a_{j},-}_{(j-1)} \circ \sigma^{i-1} 
\circ  (b_1\otimes \ldots \otimes b_n) \,,
\end{align*}
using the notation \eqref{20240805:eq1}. 
We can then write \eqref{Eq:l-C4-ij1} as 
\begin{equation}
    \begin{aligned}  \label{Eq:l-C4-ij3}
&(-1)^{\nu_{ij}} \, \dsq{a_{i},-}_{(i)} \circ \dsq{a_{j},-}_{(j-1)} \circ \sigma^{i-1}   \left(\dgal{a_{i+1},\ldots,a_{j-1},a_{j+1},\ldots,a_{i-1}} \right)   \\
=&(-1)^{\nu_{ij}+(n-1)(i-1)} \, \dsq{a_{i},-}_{(i)} \circ \dsq{a_{j},-}_{(j-1)} 
(\ldb{a_1,\stackrel{i}{\check{\dots}},\stackrel{j}{\check{\dots}} ,a_n\rdb} ) 
    \end{aligned}
\end{equation}
after using cyclic skewsymmetry \eqref{Eq:nbr-Cycl} of $\dgal{-}$. 
Similarly, we note that 
\begin{align*}
&\sigma^{i-1} \!\circ\! \dsq{a_{i},-}_L 
\!\circ \! \sigma^{j-i-1}  \!\circ \! \dsq{a_{j},-}_L  
(b_1\otimes \ldots \otimes b_n) \\
=& b_{n-j+3}\otimes \ldots \otimes b_{n-j+i+1} \otimes \dsq{a_i,b_{n-j+i+2}} \otimes  \\
&\qquad \qquad \otimes 
b_{n-j+i+3} \otimes \ldots \otimes b_n \otimes \dsq{a_j,b_1} \otimes 
 \ldots \otimes b_{n-j+2} \\
=& \dsq{a_{i},-}_{(i)} \circ \dsq{a_{j},-}_{(j-1)} \circ \sigma^{j-2} 
(b_1\otimes \ldots \otimes b_n)\,,
\end{align*}
Using this calculation and cyclic skewsymmetry \eqref{Eq:nbr-Cycl}, 
\eqref{Eq:l-C4-ij2} becomes 
\begin{equation}
    \begin{aligned}  \label{Eq:l-C4-ij4} 
&(-1)^{\mu_{ij}}\, \dsq{a_{i},-}_{(i)} \circ \dsq{a_{j},-}_{(j-1)} \circ \sigma^{j-2} 
\dgal{a_{j+1},\ldots,a_{i-1},a_{i+1},\ldots,a_{j-1}}
\\
&(-1)^{\mu_{ij}+(n-1)(j-2)}\, \dsq{a_{i},-}_{(i)} \circ \dsq{a_{j},-}_{(j-1)} 
(\ldb{a_1,\stackrel{i}{\check{\dots}},\stackrel{j}{\check{\dots}} ,a_n\rdb} ) 
    \end{aligned}
\end{equation}
since $i<j$. 
Finally, noting that $\nu_{ij}+(n-1)(i-1)=\mu_{ij}+(n-1)(j-2)+1$ modulo $2$, 
\eqref{Eq:l-C4-ij3} and \eqref{Eq:l-C4-ij4} cancel out. 
This establishes \eqref{Eq:dP-pf-step5}, hence we can conclude.  
\qed

\subsection{The completed double Poisson cohomology}

In view of Theorem \ref{Thm:g-dPcoh1}, we have a square zero differential of degree $+1$ on the graded vector space $\wBRA_B(\mc A)$, denoted $\wdd$. 
This operation is defined thanks to the double Poisson bracket $\dsq{-,-}$ according to \eqref{Eq:dP-gen-0} and \eqref{Eq:dP-gen}. 

\begin{definition}  \label{Def:dPH-comp}
The cohomology of the complex $(\wBRA_B(\cA) , \wdd)$ is called the \emph{completed double Poisson cohomology} of $\cA$ with respect to the double Poisson bracket $\dsq{-,-}$. 
Explicitly, we have $\widehat{\dPH}(\cA)=(\widehat{\dPH}^k(\cA))_{k \geq 0}$\glslink{hatdPH}{} for 
\begin{equation}
\widehat{\dPH}^k(\cA)=\frac{\ker\, \wdd:\wBRA_B(\cA)_k\to \wBRA_B(\cA)_{k+1}}{\im \,  \wdd:\wBRA_B(\cA)_{k-1}\to \wBRA_B(\cA)_k}\,.
\end{equation}
\end{definition}

\subsection{Cohomology for double Lie algebras}

Recall from \cite{DSKV} (see also \cite{ORS,Sch}) that a \emph{double Lie algebra} is a vector space $\cA$ equipped with a linear map $\dsq{-,-}: \cA\otimes \cA \to \cA \otimes \cA$ satisfying cyclic skewsymmetry \eqref{Eq:db-skew} and the Poisson property \eqref{Eq:dJacExpl}.
Put $\BRA_{\operatorname{dLie}}^0(\cA)=\cA$, and for $n\geq 1$ let $\BRA_{\operatorname{dLie}}^n(\cA)$ be the vector space spanned by linear maps
$\dgal{-,\ldots,-} : \cA^{\otimes n} \to \cA^{\otimes n}$
satisfying only the cyclic skewsymmetry \eqref{Eq:nbr-Cycl}.
Set
$\BRA_{\operatorname{dLie}}(\cA)=\bigoplus_{n\geq 0} \BRA_{\operatorname{dLie}}^n(\cA)$.
If one forgets about the derivation properties of $n$-brackets as part of the proof of Theorem \ref{Thm:g-dPcoh1},
one readily sees that the operation
\[
 \wdd : \BRA_{\operatorname{dLie}}(\cA) \to \BRA_{\operatorname{dLie}}(\cA)
\]
given by \eqref{Eq:dP-gen-0} and \eqref{Eq:dP-gen} is well-defined, and it squares to zero.
Hence we can define the \emph{double Lie algebra cohomology} of $(\cA, \dsq{-,-})$
as $\widehat{\operatorname{dLH}}(\cA)= \coH(\BRA_{\operatorname{dLie}}(\cA), \wdd)$.

\section[Relation to the double Poisson cohomology of P-V\MakeLowercase{d}W]{Relation to the double Poisson cohomology of Pichereau and Van de Weyer} 
\label{ss:ReldPCoh}

We are ready to compare the double Poisson cohomology $\dPH(\cA)$ defined on $(\mb T^\ast \cA)_\sharp$ by Pichereau and Van de Weyer \cite{PV}, cf. Definition \ref{Def:dPH}, 
and its completed version $\widehat{\dPH}(\cA)$ defined on $\wBRA_B(\cA)$ as part of Definition \ref{Def:dPH-comp}.
We recall from Proposition \ref{Pr:MapMu} that we have well-defined maps 
$\mu_n:(\mb T^\ast \cA)_{\sharp,n}\to \wBRA_B(\cA)_n$ sending 
(the class of) a noncommutative $n$-vector field to an $n$-bracket. 

\subsection{Statements}

\begin{theorem} \label{Thm:g-dPcoh2}
If $\dsq{-,-}=\dgal{-,-}_P$ is a double Poisson bracket defined by $P\in (T^\ast \mc A)_2$ through \eqref{Eq:dbr-P}, then the sequence of maps $(\mu_n)_{n\geq 0}$ defines a morphism of complexes 
$(\mb T^\ast \cA)_\sharp\longrightarrow \wBRA_B(\cA)$ endowed with the square-zero differentials $\dd_P$ \eqref{Eq:dP-PVdW} and $\wdd$ from Definition \ref{def:wdd} as follows 
%\pecetta{the diagram has been commented to speed up compiling}
%\begin{comment}
\begin{center}
      \begin{tikzpicture}[scale=0.9]
%%%% top row
 \node  (zero) at (-4.5,0) {$0$};
 \node  (T0) at (-2.8,0) {$(\mb T^\ast \cA)_{\sharp,0}$};
 \node  (T1) at (0,0) {$(\mb T^\ast \cA)_{\sharp,1}$};
 \node  (T2) at (3,0) {$(\mb T^\ast \cA)_{\sharp,2}$};
 \node  (T3) at (6,0) {$(\mb T^\ast \cA)_{\sharp,3}$};
 \node  (dots) at (8.5,0) {$\cdots$};
\path[->,>=angle 90,font=\small]  
   (zero) edge  (T0) ;
\path[->,>=angle 90,font=\small]  
   (T0) edge node[above] {$\dd_P$}  (T1) ;
\path[->,>=angle 90,font=\small]  
   (T1) edge node[above] {$\dd_P$}  (T2) ;
\path[->,>=angle 90,font=\small]  
   (T2) edge node[above] {$\dd_P$}  (T3) ;
\path[->,>=angle 90,font=\small]  
   (T3) edge node[above] {$\dd_P$}  (dots) ;
%%%% bottom row 
 \node  (Bzero) at (-4.5,-3) {$0$};
 \node  (B0) at (-2.8,-3) {$\wBRA_B(\cA)_0$};
 \node  (B1) at (0,-3) {$\wBRA_B(\cA)_1$};
 \node  (B2) at (3,-3) {$\wBRA_B(\cA)_2$};
 \node  (B3) at (6,-3) {$\wBRA_B(\cA)_3$};
 \node  (Bdots) at (8.5,-3) {$\cdots$}; 
 \path[->,>=angle 90,font=\small]  
   (Bzero) edge  (B0) ;
\path[->,>=angle 90,font=\small]  
   (B0) edge node[above] {$\wdd$}  (B1) ;
\path[->,>=angle 90,font=\small]  
   (B1) edge node[above] {$-\wdd$}  (B2) ;
\path[->,>=angle 90,font=\small]  
   (B2) edge node[above] {$\wdd$}  (B3) ;
\path[->,>=angle 90,font=\small]  
   (B3) edge node[above] {$-\wdd$}  (Bdots) ;
%%%% vertical arrows 
\path[->,>=angle 90,font=\small] 
(T0) edge node[left] {$\mu_0$}  (B0) ;
\path[->,>=angle 90,font=\small] 
(T1) edge node[left] {$\mu_1$}  (B1) ;
\path[->,>=angle 90,font=\small] 
(T2) edge node[left] {$\mu_2$}   (B2) ;
\path[->,>=angle 90,font=\small] 
(T3) edge node[left] {$\mu_3$}   (B3) ;
   \end{tikzpicture}
\end{center}
%\end{comment}
Thus, $(\mu_n)_{n\geq 0}$ descends to a $\kk$-linear map $\dPH(\cA)\to \widehat{\dPH}(\cA)$ in cohomology.
\end{theorem}

\begin{corollary} \label{Cor:Iso-dPcoh}
Under the assumptions of Proposition \ref{Pr:MapMu-Iso}, the morphism of complexes presented in Theorem \ref{Thm:g-dPcoh2} is an isomorphism. 
In particular, the double Poisson cohomology $\dPH(\cA)$ and its completed version $\widehat{\dPH}(\cA)$ are isomorphic. 
\end{corollary}

\subsection{Proof of Theorem \ref{Thm:g-dPcoh2}}

We start with the easiest case: showing that the first square in the diagram is commutative. This amounts to prove that $\wdd=\mu_1\circ \dd_P$ since $\mu_0=\Id_{\cA_\sharp}$.  
Without loss of generality, assume that $P=P_1 P_2$ with $P_1,P_2\in \DDer_B(\cA)$. 
Then, given $\bar{a}\in \cA_\sharp$ and an arbitrary lift $a\in \cA$, we get from \eqref{Eq:dP-0} modulo commutators 
$$\dd_P(\bar{a})=P_1 P_2(a)'' P_2(a)' -  P_2 P_1(a)'' P_1(a)'\, \in (T^\ast \mc A)_{\sharp,1}\,.$$
Since by definition $\mu_1(\delta)(b)=\delta(b)' \delta(b)''=(\mult \circ \delta)(b)$ for any $\delta\in \DDer_B(\cA)$ and $b\in \cA$, we see that 
\begin{equation}
\mu_1(\dd_P(\bar{a}))(b) = 
P_1(b)' P_2(a)'' P_2(a)'  P_1(b)''
- P_2(b)' P_1(a)''P_1(a)'  P_2(b)''\,.
\end{equation}
(Here, we used that the $\cA$-bimodule structure of $\DDer_B(\cA)$ is based on the inner bimodule structure \eqref{bmodinner} of $\cA^{\otimes 2}$.) 
A direct comparison with \eqref{Eq:dbr-P} yields that 
$\mu_1(\dd_P(\bar{a}))(b)=-\mult \circ \dsq{a,b}$. By \eqref{Eq:dP-gen-0}, this is $\wdd(\bar{a})(b)$ and we can conclude.

Next, we look at the general case: 
$\mu_{n+1}\circ \dd_P=(-1)^n \wdd \circ \mu_\ell$. 
Applied to $Q\in (\mb T^\ast \cA)_{n}$, the left-hand side is the $(n+1)$-bracket $\mu_{n+1}(\brSN{P,Q})$. This is given by, cf. Proposition \ref{Pr:MapMu},  
\begin{equation} \label{Eq:dP2-pf1}
\mu_{n+1}(\brSN{P,Q})=
\sum_{s=0}^{n} (-1)^{n s} \sigma^s \circ \dgal{-,\ldots,-}_{\brSN{P,Q}}^{\sim} \circ \sigma^{-s}\,.
\end{equation}
Explicitly, if we write $P=P_1P_2$ and $Q=Q_1\cdots Q_n$, 
$\brSN{P,Q}$ has the form \eqref{Eq:dSN-QPpoiss} (with \eqref{Eq:dSN-QP1}--\eqref{Eq:dSN-QP2}) and 
we can write down $\dgal{a_0,\ldots,a_\ell}_{\brSN{P,Q}}^{\sim}$ using \eqref{Eq:MapMu} as  
\begin{small}
\begin{align*}
% first group of terms
&\sum_{i=1,2}  (-1)^{n+i} 
 P_{i+1}(a_n)' d^{i}_{n} Q_{1}(a_{0})''\otimes \ldots \otimes Q_{n-1}(a_{n-2})' \delta^{i}_{n} (a_{n-1})'' \otimes  \delta^{i}_{n} (a_{n-1})' P_{i+1}(a_n)'' \\
&+\sum_{\substack{i=1,2 \\ 1\leq j \leq n-1}}  (-1)^{j+i} 
Q_n(a_n)'Q_1(a_0)''\otimes \ldots \otimes 
Q_{j-1}(a_{j-2})'\delta^{i}_{j}(a_{j-1})''\otimes \delta^{i}_{j}(a_{j-1})' P_{i+1}(a_j)'' \\
&\hspace{3cm} \otimes P_{i+1}(a_j)' d^{i}_{j} Q_{j+1}(a_{j+1})''\otimes \ldots \otimes Q_{n-1}(a_{n-1})'Q_n(a_n)'' \\
%& second group of terms 
&-\sum_{i=1,2}  (-1)^{n+i} 
 \gamma^{i}_{n}(a_n)' Q_{1}(a_{0})''\otimes \ldots \otimes Q_{n-1}(a_{n-2})' c^{i}_{n} P_{i+1}(a_{n-1})'' \otimes  P_{i+1}(a_{n-1})' \gamma^{i}_{n}(a_n)'' \\
&-\!\!\sum_{\substack{i=1,2 \\ 1\leq j \leq n-1}}  (-1)^{j+i} 
Q_n(a_n)'Q_1(a_0)''\otimes \ldots \otimes 
Q_{j-1}(a_{j-2})'c^{i}_{j} P_{i+1}(a_{j-1})''\otimes 
P_{i+1}(a_{j-1})' \gamma^{i}_{j}(a_{j})'' \\
&\hspace{3cm} \otimes \gamma^{i}_{j}(a_j)' Q_{j+1}(a_{j+1})''\otimes \ldots \otimes Q_{n-1}(a_{n-1})'Q_n(a_n)''\, .
\end{align*}
\end{small}
Thanks to \eqref{Eq:dSN-QP1}--\eqref{Eq:dSN-QP2}, we get 
\begingroup
\allowdisplaybreaks
\begin{subequations}
    \begin{align}
&\dgal{a_0,\ldots,a_\ell}_{\brSN{P,Q}}^{\sim}   = \nonumber \\
&\sum_{i=1,2}  (-1)^{n+i} 
 P_{i+1}(a_n)' P_i(Q_n(a_{n-1})')'' Q_{1}(a_{0})''\otimes \ldots \nonumber \\
&\hspace{2cm} \ldots\otimes Q_{n-1}(a_{n-2})' Q_n(a_{n-1})'' \otimes  P_i(Q_n(a_{n-1})')' P_{i+1}(a_n)''
\label{Eq:dP2-B1}\\
&-\sum_{i=1,2}  (-1)^{n+i} 
 P_{i+1}(a_n)' Q_n(P_i(a_{n-1})'')' Q_{1}(a_{0})''\otimes \ldots \nonumber \\
&\hspace{2cm} \ldots \otimes Q_{n-1}(a_{n-2})' Q_n(P_i(a_{n-1})'')'' \otimes P_i(a_{n-1})' P_{i+1}(a_n)'' 
\label{Eq:dP2-B2} \\
&+\sum_{\substack{i=1,2 \\ 1\leq j \leq n-1}}  (-1)^{j+i} 
Q_n(a_n)'Q_1(a_0)''\otimes \ldots \otimes 
Q_{j-1}(a_{j-2})' Q_j(a_{j-1})''\otimes \nonumber \\
&\hspace{1cm}  P_i(Q_j(a_{j-1})')' P_{i+1}(a_j)'' \otimes P_{i+1}(a_j)' P_i(Q_j(a_{j-1})')'' Q_{j+1}(a_{j+1})''\otimes\ldots \nonumber \\
&\hspace{4cm} \ldots \otimes Q_{n-1}(a_{n-1})'Q_n(a_n)'' \label{Eq:dP2-B3} \\ 
&-\sum_{\substack{i=1,2 \\ 1\leq j \leq n-1}}  (-1)^{j+i} 
Q_n(a_n)'Q_1(a_0)''\otimes \ldots \otimes 
Q_{j-1}(a_{j-2})' Q_j(P_i(a_{j-1})'')'' \otimes \nonumber \\
&\hspace{1cm} P_i(a_{j-1})' P_{i+1}(a_j)'' \otimes P_{i+1}(a_j)' Q_j(P_i(a_{j-1})'')' Q_{j+1}(a_{j+1})''\otimes \ldots \nonumber \\
&\hspace{4cm} \ldots \otimes Q_{n-1}(a_{n-1})'Q_n(a_n)'' \label{Eq:dP2-B4} \\
&-\sum_{i=1,2}  (-1)^{n+i} 
 Q_n(a_n)'' Q_{1}(a_{0})''\otimes \ldots \nonumber \\
 &\ldots \otimes Q_{n-1}(a_{n-2})' P_i(Q_n(a_n)'')' P_{i+1}(a_{n-1})'' \otimes  P_{i+1}(a_{n-1})' P_i(Q_n(a_n)'')'' \label{Eq:dP2-B5}  \\
 &+\sum_{i=1,2}  (-1)^{n+i}  Q_n(P_i(a_n)')' Q_{1}(a_{0})''\otimes \ldots \nonumber \\
 &\hspace{1cm}\ldots \otimes Q_{n-1}(a_{n-2})' Q_n(P_i(a_n)')'' P_{i+1}(a_{n-1})'' \otimes  P_{i+1}(a_{n-1})' P_i(a_n)'' \label{Eq:dP2-B6}  \\
 &-\!\!\sum_{\substack{i=1,2 \\ 1\leq j \leq n-1}}  (-1)^{j+i} 
Q_n(a_n)'Q_1(a_0)''\otimes \ldots \otimes 
Q_{j-1}(a_{j-2})' P_i(Q_j(a_j)'')' P_{i+1}(a_{j-1})''\otimes \nonumber \\
&\hspace{1.5cm} P_{i+1}(a_{j-1})' P_i(Q_j(a_j)'')''  \otimes Q_j(a_j)' Q_{j+1}(a_{j+1})''\otimes \ldots \nonumber \\
&\hspace{4cm} \ldots \otimes Q_{n-1}(a_{n-1})'Q_n(a_n)'' \label{Eq:dP2-B7}  \\
&+\!\!\sum_{\substack{i=1,2 \\ 1\leq j \leq n-1}}  (-1)^{j+i} 
Q_n(a_n)'Q_1(a_0)''\otimes \ldots \otimes 
Q_{j-1}(a_{j-2})' Q_j(P_i(a_j)')'' P_{i+1}(a_{j-1})''\otimes \nonumber \\
&\hspace{1.5cm} P_{i+1}(a_{j-1})' P_i(a_j)''  \otimes Q_j(P_i(a_j)')' Q_{j+1}(a_{j+1})''\otimes \ldots 
\nonumber \\
&\hspace{4cm} \ldots \otimes Q_{n-1}(a_{n-1})'Q_n(a_n)'' \,. \label{Eq:dP2-B8} 
    \end{align}
\end{subequations}
\endgroup 
%\delta^i_j(a)\otimes d^i_j=&
%P_i(Q_j(a)')' \otimes Q_j(a)'' \otimes P_i(Q_j(a)')'' \nonumber \\
%&- P_i(a)' \otimes Q_j(P_i(a)'')''\otimes Q_j(P_i(a)'')'\,, \label{Eq:dSN-QP1}    \\
%c^i_j\otimes \gamma^i_j(a)=& 
%P_i(Q_j(a)'')' \otimes Q_j(a)'\otimes P_i(Q_j(a)'')'' \nonumber \\
%&- Q_j(P_i(a)')''\otimes  Q_j(P_i(a)')' \otimes P_i(a)'' \,. \label{Eq:dSN-QP2} 
Let us now group these terms together by repeatedly using the defining properties of double derivations: 
$\delta(a)'\otimes \delta(a)'' b + a\delta(b)'\otimes \delta(b)''=\delta(ab)$ for any $a,b\in \cA$, which is satisfied by all $P_i,Q_j$. 
We gather \eqref{Eq:dP2-B3}, with \eqref{Eq:dP2-B5} and the summands $j=2,\ldots,n-1$ of \eqref{Eq:dP2-B7} (where for the latter we replace $j$ by $j+1$) as 
\begin{equation*}
    \begin{aligned}
&\sum_{\substack{i=1,2 \\ 1\leq j \leq n-1}}  (-1)^{j+i} 
Q_n(a_n)'Q_1(a_0)''\otimes \ldots \otimes 
Q_{j-1}(a_{j-2})' Q_j(a_{j-1})''\otimes \nonumber \\
&\qquad  P_i(Q_j(a_{j-1})'Q_{j+1}(a_{j+1})'')' P_{i+1}(a_j)'' \otimes P_{i+1}(a_j)' P_i(Q_j(a_{j-1})'Q_{j+1}(a_{j+1})'')'' \otimes \nonumber \\
&\hspace{4cm} Q_{j+1}(a_{j+1})' Q_{j+2}(a_{j+2})''\otimes \ldots \otimes Q_{n-1}(a_{n-1})'Q_n(a_n)'' 
    \end{aligned}
\end{equation*}
In view of \eqref{Eq:dbr-P}, the double Poisson bracket $\dsq{-,-}=\dgal{-,-}_{P_1P_2}$ can be written as 
$\dsq{a,b}=\sum_{i=1,2}(-1)^{i}P_{i}(b)'P_{i+1}(a)''\otimes P_{i+1}(a)'P_{i}(b)''$ for any $a,b\in \cA$. Thus, the previous expression becomes  
\begin{equation}
    \begin{aligned} \label{Eq:dP2-C1}
&\sum_{1\leq j \leq n-1}  (-1)^{j} 
Q_n(a_n)'Q_1(a_0)''\otimes \ldots \otimes 
Q_{j-1}(a_{j-2})' Q_j(a_{j-1})''\otimes  \\
&\hspace{4cm} \dsq{a_j, Q_j(a_{j-1})'Q_{j+1}(a_{j+1})''} \otimes \\
&\hspace{3cm} Q_{j+1}(a_{j+1})' Q_{j+2}(a_{j+2})''\otimes \ldots \otimes Q_{n-1}(a_{n-1})'Q_n(a_n)'' 
    \end{aligned}
\end{equation}
Noting from Proposition \ref{Pr:MapMu} that for any $b_1,\ldots,b_\ell\in \cA$,  
\begin{align}   \label{Eq:dP2-Q}
\dgal{b_1,\ldots,b_\ell}_Q^{\sim} = 
Q_{n}(b_{n})'Q_{1}(b_1)''\otimes \ldots \otimes   
 Q_{n-1}(b_{n-1})'Q_{n}(b_n)''\,,
\end{align}    
we can finally write \eqref{Eq:dP2-C1} as 
\begin{equation}
    \begin{aligned} \label{Eq:dP2-C1bis}
&\sum_{j=1}^{n-1}  (-1)^{j}  
\dsq{a_j,-}_{(j)}\circ \ldb a_0,\stackrel{j}{\check{\dots}},a_n\rdb_{Q}^{\sim} \,.
    \end{aligned}
\end{equation}
Next, we sum \eqref{Eq:dP2-B4} with \eqref{Eq:dP2-B8} (where we replace $i$ by $i+1$) to obtain 
\begin{align*}
-&\sum_{\substack{i=1,2 \\ 1\leq j \leq n-1}}  (-1)^{j+i} 
Q_n(a_n)'Q_1(a_0)''\otimes \ldots \otimes 
Q_{j-2}(a_{j-3})' Q_{j-1}(a_{j-2})'' \otimes \\
&Q_{j-1}(a_{j-2})' Q_j(P_{i+1}(a_j)' P_i(a_{j-1})'')'' \otimes 
 P_i(a_{j-1})' P_{i+1}(a_j)'' \otimes \\
&\hspace{2cm} Q_j(P_{i+1}(a_j)' P_i(a_{j-1})'')' Q_{j+1}(a_{j+1})''\otimes  \\
&\hspace{2cm} 
Q_{j+1}(a_{j+1})' Q_{j+2}(a_{j+2})'' \ldots \otimes  
\otimes Q_{n-1}(a_{n-1})'Q_n(a_n)'' \,.    
\end{align*}
With the help of our previous observations, this can be written as 
\begin{equation}
    \begin{aligned} \label{Eq:dP2-C2}
&\sum_{j=1}^{n-1}  (-1)^{j}  
\ldb a_0,\ldots,a_{j-2},\dsq{a_{j-1},a_j}',a_{j+1}\ldots,a_n\rdb_{Q}^{\sim} \otimes_{n-j} \dsq{a_{j-1},a_j}'' ,
    \end{aligned}
\end{equation}
where we us the notation \eqref{eq:tensor-i-notation}. 
Similarly, we sum \eqref{Eq:dP2-B2} with \eqref{Eq:dP2-B6} as 
\begin{align*}
 -\sum_{i=1,2}  &(-1)^{n+i} 
  Q_n(P_{i+1}(a_n)'P_i(a_{n-1})'')' Q_{1}(a_{0})''\otimes \ldots \nonumber \\
&\ldots \otimes Q_{n-1}(a_{n-2})' Q_n(P_{i+1}(a_n)'P_i(a_{n-1})'')'' \otimes P_i(a_{n-1})' P_{i+1}(a_n)'' ,    
\end{align*}
and it becomes with the notation \eqref{20240805:eq1}
\begin{equation}
    \begin{aligned} \label{Eq:dP2-C3}
& (-1)^{n}  
(\ldb a_0,\ldots,a_{n-2},-\rdb_{Q}^{\sim})_L\circ\dsq{a_{n-1},a_n}.
    \end{aligned}
\end{equation}
The only two terms that have not yet been considered are \eqref{Eq:dP2-B1} and the summand $j=1$ from \eqref{Eq:dP2-B7}, which we denote \eqref{Eq:dP2-B7}$|_{j=1}$. 
In the end, we are only interested in \eqref{Eq:dP2-pf1}, and therefore we can replace (as we did several times in the proof of Theorem \ref{Thm:g-dPcoh1}) \eqref{Eq:dP2-B7}$|_{j=1}$ 
by the corresponding term occurring in $(-1)^{n s} \sigma^s \dgal{a_s,\ldots,a_n,a_0,\ldots,a_{s-1}}_{\brSN{P,Q}}^{\sim}$ with some $1\leq s\leq n$. Thus, 
\begin{align*}
\eqref{Eq:dP2-B7}|_{j=1} =
 \sum_{i=1,2} & (-1)^{i}\,
Q_n(a_n)' P_i(Q_1(a_1)'')' P_{i+1}(a_{0})''\otimes \\
&P_{i+1}(a_{0})' P_i(Q_1(a_{1})'')''\otimes \ldots \otimes Q_{n-1}(a_{n-1})' Q_n(a_n)''   
\end{align*}
can be replaced with the following term coming from the summand $s=n$ in \eqref{Eq:dP2-pf1}: 
\begin{align*} 
 \sum_{i=1,2} & (-1)^{i+n}\,
P_{i+1}(a_{n})' P_i(Q_1(a_{0})'')''\otimes \ldots \otimes Q_{n-1}(a_{n-2})' Q_n(a_{n-1})'' \otimes \\
&\hspace{2cm} Q_n(a_{n-1})' P_i(Q_1(a_0)'')' P_{i+1}(a_{n})''  \,.
\end{align*}
Summing this expression with \eqref{Eq:dP2-B1} yields 
\begin{align*} 
 \sum_{i=1,2} & (-1)^{i+n}\,
P_{i+1}(a_{n})' P_i(Q_n(a_{n-1})' Q_1(a_0)'')''\otimes \ldots \otimes Q_{n-1}(a_{n-2})' Q_n(a_{n-1})'' \otimes \\
&\hspace{2cm} P_i(Q_n(a_{n-1})' Q_1(a_0)'')' P_{i+1}(a_{n})''  \,,
\end{align*}
which can be written as 
\begin{equation}
    \begin{aligned} \label{Eq:dP2-C4}
& (-1)^{n}\, \sigma^{-1}  \circ 
\dsq{a_n,-}_{L}\circ \ldb a_0,\ldots,a_{n-1}\rdb_{Q}^{\sim} \,.
    \end{aligned}
\end{equation}
We are now in position to derive from \eqref{Eq:dP2-pf1},  \eqref{Eq:dP2-C1bis}, \eqref{Eq:dP2-C2}, \eqref{Eq:dP2-C3} and \eqref{Eq:dP2-C4} that 
$\mu_{n+1}(\brSN{P,Q}) =
\sum_{s=0}^{n} (-1)^{n s} \sigma^s \circ \Xi \circ \sigma^{-s}$ for 
\begin{align*}
    \Xi(a_0,&\ldots,a_n)=
    \sum_{j=1}^{n-1}  (-1)^{j}  
\dsq{a_j,-}_{(j)}\circ \ldb a_0,\stackrel{j}{\check{\dots}},a_n\rdb_{Q}^{\sim} \\%C1
&+(-1)^{n}\, \sigma^{-1}  \circ 
\dsq{a_n,-}_{L}\circ \ldb a_0,\ldots,a_{n-1}\rdb_{Q}^{\sim}\\ %C4
&+\sum_{j=1}^{n-1}  (-1)^{j}  
\ldb a_0,\ldots,a_{j-2},\dsq{a_{j-1},a_j}',a_{j+1}\ldots,a_n\rdb_{Q}^{\sim} \otimes_{n-j} \dsq{a_{j-1},a_j}'' \\ %C2
&+(-1)^{n}  
\, \ldb a_0,\ldots,a_{n-2},\dsq{a_{n-1},a_n}'\rdb_{Q}^{\sim}
\otimes \dsq{a_{n-1},a_n}''\,.  % C3
\end{align*}
As we already noticed, we can make use of the $\Z_{n+1}$ cyclic skewsymmetry of the $(n+1)$-bracket $\mu_{n+1}(\brSN{P,Q})$ to replace $\Xi$ with $\widehat{\Xi}$ defined as  
\begin{align*}
    &\widehat\Xi(a_0,\ldots,a_n)=
    \sum_{j=0}^{n-1}  (-1)^{j(n-1)}  
\dsq{a_0,-}_{L}\circ (\sigma^{-j}\circ\ldb -\rdb_{Q}^{\sim}\circ \sigma^j)(a_1,\ldots,a_n) \\
&+  \sum_{j=0}^{n-1}  (-1)^{n+j(n-1)}
(\sigma^{-j}\circ\ldb -\rdb_{Q}^{\sim}\circ \sigma^j)(a_0,\ldots,a_{n-2},\dsq{a_{n-1},a_n}') \otimes \dsq{a_{n-1},a_n}''
\\
&=\dsq{a_0,-}_{L}( \ldb a_1,\ldots,a_n \rdb_{Q} )
+(-1)^{n}  
\, \ldb a_0,\ldots,a_{n-2},-\rdb_{Q,L} (\dsq{a_{n-1},a_n})
\end{align*}
Combined with Proposition \ref{Pr:MapMu}, we can write 
\begin{equation} 
    \begin{aligned}
\mu_{n+1}(&\brSN{P,Q}) =
\sum_{s=0}^n (-1)^{ns}\sigma^s \circ  
(\dsq{-,-}\otimes \Id_{\cA}^{\otimes(n-1)}) \circ (\Id_{\cA}\otimes \dgal{-}_Q) \circ \sigma^{-s}\\
& +(-1)^n\sum_{s=0}^n (-1)^{ns}\sigma^s \circ  
(\dgal{-}_Q\otimes \Id_{\cA})\circ (\Id_{\cA}^{\otimes(n-1)}\otimes \dsq{-,-})
\circ \sigma^{-s}  \,.
    \end{aligned}
\end{equation}
This is precisely $(-1)^n\,\wdd(\dgal{-}_Q)=(-1)^n\,\wdd(\mu_n(Q))$ in view of \eqref{Eq:dP-gen}.  
\qed

\begin{remark} \label{Rem:ZTdPA}
After the first version of this memoir appeared on arXiv, the authors of \cite{ZT2} have shared with us their work in which they define a dPA cohomology.
In that work, the authors construct a graded Lie bracket $[-,-]_{\operatorname{DP}}$ on $\wBRA_B(\cA)$ and they show that the map
 $$\mu : ((\mb T^\ast \cA)_\sharp,\br{-,-}_{\SN})\longrightarrow (\wBRA_B(\cA),[-,-]_{\operatorname{DP}})$$
 is a morphism of graded Lie algebras, see \cite[Thm.~4.1]{ZT2}.
If  $\mc P := \dsq{-,-}$ is a double Poisson bracket, then the adjoint action
$\operatorname{ad}_{\mc P}:=[\mc P,-]_{\operatorname{DP}}$ defines a square-zero differential on $\wBRA_B(\cA)$ which is related to the differential $\wdd$ from Theorem \ref{Thm:g-dPcoh1} through $\operatorname{ad}_{\mc P}(Q)=(-1)^n\wdd(Q)$ for $Q\in \wBRA_B(\cA)_n$.
This provides an alternative proof of Theorem \ref{Thm:g-dPcoh2}.
\end{remark}

\section{Chemla's formula}

\begin{proposition}  \label{Pr:Chemla}
Fix a double bracket $\dsq{-,-} \in \wBRA_B(\mc A)_2$ and an $n$-bracket  $\dgal{-}\in \wBRA_B(\mc A)_n$, $n\geq 1$. 
Then, the $(n+1)$-bracket $\wdd(\dgal{-})$  \eqref{Eq:dP-gen} satisfies 
\begin{equation} \label{Eq:Chemla}
\begin{aligned}
&\wdd(\dgal{-})(a_1 \otimes \ldots \otimes a_{n+1}) \\
=&\sum_{i=1}^{n+1} (-1)^{ni}\, \sigma^i \, 
\dgal{a_{i+1},\ldots,a_{n+1},a_1,\ldots,a_{i-2}, \dsq{a_{i-1},a_i}}_L   \\
&+\sum_{i=1}^{n+1} (-1)^{ni}\, \sigma^{i-1} \,   
\dsq{a_i,\dgal{a_{i+1},\ldots,a_{n+1},a_1,\ldots,a_{i-1}} }_{L}  \,,
\end{aligned}
\end{equation}
for any $a_1,\ldots,a_{n+1}\in \cA$. We call \eqref{Eq:Chemla} ``Chemla's formula''.
\end{proposition}
\begin{proof}
Using \eqref{Eq:dP-gen} and the left extensions of $\dsq{-,-}$ and $\dgal{-}$, cf. \eqref{20240805:eq1}, we can write the left-hand side of \eqref{Eq:Chemla} as 
\begin{align*}
&\sum_{s=1}^{n+1} (-1)^{ns}\sigma^s \circ  
\dgal{-}_L \circ (\Id_{\cA}^{\otimes(n-1)}\otimes \dsq{-,-})
(a_{s+1} \otimes \ldots \otimes a_{n+1} \otimes a_1 \otimes \ldots \otimes a_{s})  \\
+&\sum_{s=0}^{n} (-1)^{n(s+1)}\sigma^s \circ  
\dsq{-,-}_L \circ (\Id_{\cA}\otimes \dgal{-})
(a_{s+1} \otimes \ldots \otimes a_{n+1} \otimes a_1 \otimes \ldots \otimes a_{s})  \,,
\end{align*}
with indices taken in $\{1,\ldots,n+1\}$ modulo $n+1$. 
We get the required formula after relabeling indices. 
\end{proof}

\subsection{A differential from double Lie-Rinehart algebras}

Recall that the operation  $\wdd$ from Theorem \ref{Thm:g-dPcoh1} is a square-zero differential when $\dsq{-,-}$ is a double Poisson bracket. 
A few months after defining this operation around the end of 2022, we noticed that our differential looked similar to a formula obtained by S.~Chemla in her study \cite{Ch} of differentials on double Lie-Rinehart algebras\footnote{They are called \emph{double Lie algebroids} in \cite{Ch,VdB2}.}.
We shall describe the precise relation with \cite[Prop.~5.8]{Ch} below. Our choice of name for Chemla's formula \eqref{Eq:Chemla} will then become transparent.

In full generalities, fix an $\cA$-bimodule $M$ and let  
\begin{align*}
 M^\ast&=\{ q \in \Hom_{\kk} (M , \cA^{\otimes 2}) \mid q(a D b)=q(D)' b \otimes a q(D)''\quad \forall a,b\in \cA,\,D\in M\}, \\
 M_\ast&=\{ q \in \Hom_{\kk} (M , \cA^{\otimes 2}) \mid q(a D b)=aq(D)' \otimes  q(D)'' b\quad \forall a,b\in \cA,\,D\in M\} .
\end{align*}
Both $M^\ast$ and $M_\ast$ are $\cA$-bimodules when endowed with the remaining $\cA$-bimodule structure on $\cA^{\otimes 2}$. Furthermore, the map $q \mapsto \tau_{(12)}\circ q$ is an isomorphism of bimodules between $M^\ast$ and $M_\ast$. 
For any $q_1,\ldots,q_n\in M_\ast$, consider their product $q_1\cdots q_n:=q_1 \otimes_\cA \ldots\otimes_\cA q_n$  in $T_\cA(M_\ast)$. There is a well-defined operation 
\begin{equation}
\begin{aligned}
 &   (T_\cA(M_\ast))_n \to \Hom(M^{\otimes n},\cA^{\otimes n}), \\
&\dgal{-,\ldots,-}_{q_1\cdots q_n}  
=\sum_{i=1}^n \tau_{(1\ldots n)}^i \circ \Phi_M(q_1 \otimes \ldots \otimes  q_n) \circ \tau_{(1\ldots n)}^{-i} \,, 
\end{aligned}
\end{equation}
where for any $D_1,\ldots,D_n\in M$, 
\begin{equation}
    \begin{aligned}
&\Phi_M(q_1 \otimes \ldots\otimes  q_n)(D_1\otimes \ldots \otimes D_n) \\ 
=&q_n(D_n)'q_1(D_1)'' \otimes q_1(D_1)'q_2(D_2)'' \otimes \ldots \otimes q_{n-1}(D_{n-1})'q_n(D_n)''\,.
    \end{aligned}
\end{equation}
Furthermore, this operation factors through $(T_\cA(M_\ast))_\sharp$ which is $T_\cA(M_\ast)$ modulo graded commutators (with $M_\ast$ in degree $+1$ and $\cA$ in degree $0$). 
Van den Bergh \cite{VdB1} has shown that this is an isomorphism (of left $(\cA^e)^{\otimes n}$-modules) onto the signed invariants in $\Hom(M^{\otimes n},\cA^{\otimes n})$ provided that $M$ is finitely generated and projective. 
Note that for the $\cA$-bimodule of $B$-relative $1$-forms $M=\Omega^1_{\cA/B}$, we have $M_{\ast}=\DDer_B(\cA)$ and we obtain Proposition \ref{Pr:MapMu}, see \cite[\S4.1]{VdB1}.

\begin{definition}[\cite{Ch,VdB2}]
Let $(\mb T^\ast \cA,\dSN{-,-})$ be the double Gerstenhaber algebra  from Theorem \ref{Thm:dSN}. 
A \emph{double Lie-Rinehart algebra} over $\cA$ is an $\cA$-bimodule $M$ endowed with 
a morphism of $\cA$-bimodules $\rho:M\to \DDer_B(\cA)$ (the \emph{anchor}) and 
a $\kk$-linear map $\dgal{-,-}_M:M\otimes M \to \cA\otimes M \oplus M \otimes \cA$ such that 
\begin{enumerate}[(1)]
\item $\dgal{D_1,D_2}_{M}=-\dgal{D_2,D_1}_{M}^\sigma$;
\item $\dSN{\rho(D_1),\rho(D_2)}=\rho(\dgal{D_1,D_2}_{M})$ for any $D_1,D_2\in M$, where we use the natural extension $\rho:  \cA\otimes M \oplus M \otimes \cA \to \cA\otimes \DDer_B(\cA) \oplus \DDer_B(\cA) \otimes \cA$;
\item $\dgal{D_1,a D_2}_{M}=D_1(a) D_2 + a \dgal{D_1,D_2}_{M}$ and 
$\dgal{D_1,D_2 a}_{M}= D_2 \, D_1(a) + \dgal{D_1,D_2}_{M} a$, where we use the products  
\begin{align*}
D_1(a) D_2:=&\rho(D_1)(a)'\otimes \rho(D_1)(a)'' D_2 \in \cA\otimes M \,, \\
D_2 \, D_1(a):=&D_2 \rho(D_1)(a)'\otimes \rho(D_1)(a)''\in M\otimes \cA\,.
\end{align*} 
\item $\dgal{-,-}_M$ satisfies the double Jacobi identity, i.e. the operation \eqref{Eq:dJac} written for $\dgal{-,-}_M$ is identically zero (where we set $\dgal{D,a}_M:=\rho(D)(a)$ when we apply $\dgal{D,-}_{M,L}$ to $\cA\otimes M$). 
\end{enumerate}
\end{definition}
The most relevant result for us from \cite{Ch} is stated as follows\footnote{Let us emphasize that the results of Chemla \cite{Ch} are originally written in terms of the $A$-bimodule $M^\ast$. We have translated them to the $A$-bimodule $M_\ast$ by a simple application of the permutation map $\tau_{(12)}$.}.

\begin{proposition}[\cite{Ch}, Prop.~5.8] \label{Pr:Chem-5.8}
Assume that $M$ is a double Lie-Rinehart algebra over $\cA$ which is finitely generated and projective as an $\cA$-bimodule. Then there is a square-zero differential $d_M$ on $(T_\cA(M_\ast))_\sharp$ which satisfies 
\begin{equation} \label{Eq:DiffCh}
\begin{aligned}
&\dd_M(\dgal{-,\ldots,-})(D_1 \otimes \ldots \otimes D_{n+1}) \\
=&\sum_{i=1}^{n+1} (-1)^{n(i+1)}\, \sigma^i \, 
\dgal{D_{i+1},\ldots,D_{n+1},D_1,\ldots,D_{i-2}, \dgal{D_{i-1},D_i}_M}_{L}   \\
&+\sum_{i=1}^{n+1} (-1)^{n(i-1)}\, \sigma^{i-1} \,   
\dgal{D_i,\dgal{D_{i+1},\ldots,D_{n+1},D_1,\ldots,D_{i-1}} }_{M,L}  \,,
\end{aligned}
\end{equation}
for any $D_1,\ldots,D_{n+1}\in M$ and $\dgal{-}:=\sum_\ell \dgal{-}_{q_1^{(\ell)}\cdots q_n^{(\ell)}}$ with $q_i^{(\ell)}\in M_\ast$. 
\end{proposition}

\begin{remark} \label{Rem:Chemla}
In the first line of the formula \eqref{Eq:DiffCh}, we use 
$$\dgal{D_{1},\ldots,D_{n-1}, \Delta_1 \otimes a_1 + a_2 \otimes \Delta_2 }_{q_1\cdots q_n,L}
:=\dgal{D_{1},\ldots,D_{n-1}, \Delta_1}_{q_1\cdots q_n} \otimes a_1\,,$$
for $a_1,a_2\in \mc A$ and $D_j,\Delta_1,\Delta_2\in M$. 
In the second line of the formula \eqref{Eq:DiffCh}, we use 
$$\dgal{D,a_1 \otimes \ldots \otimes a_n}_{M, L}
:=\dgal{D,a_1}_{M}\otimes a_2 \otimes \ldots \otimes a_n
:=\rho(D)(a_1)\otimes a_2 \otimes \ldots \otimes a_n\,,$$
for $a_j\in \cA$ and $D\in M$. 
\end{remark}

To connect this result to Proposition \ref{Pr:Chemla}, assume that $(\cA,\dsq{-,-})$ is a $B$-linear double Poisson algebra. Then $M=\Omega^1_{\cA/B}$ is a double Lie-Rinehart algebra for 
\begin{align*}
    \rho:\Omega^1_{\cA/B} \to \DDer_B(\cA), \quad \rho(\dd a):=\dsq{a,-} \,, \\
    \dgal{\dd a,\dd b}_{\Omega^1_{\mc A}}:=(\dd \otimes \Id_\cA + \Id_\cA \otimes \dd\,) \dsq{a,b}\,,
\end{align*}
with $a,b\in \mc A$ and $\dd:\mc A\to \Omega^1_{\mc A/B}$, $a\mapsto \dd a$. 
We observed that $(\Omega^1_{\mc A/B})_\ast = \DDer_B(\mc A)$, and a signed invariant in $\Hom((\Omega^1_{\mc A/B})^{\otimes n},\mc A^{\otimes n})$ is nothing else than a $B$-linear $n$-bracket according to Van den Bergh \cite{VdB1}. 
Indeed, denoting by $\dgal{-}_{Q}$ the $n$-bracket associated with $Q=Q_1\cdots Q_n\in (\mb T^\ast \cA)_n$, the correspondence is given by 
\begin{equation}
\dgal{a_1,\ldots,a_n}_{Q}  = 
    \dgal{\dd a_1,\ldots,\dd a_n}_{Q_1\cdots Q_n} ,
\end{equation}
for any $a_1,\ldots,a_n \in \cA$. It is then easy to check that \eqref{Eq:DiffCh} becomes \eqref{Eq:Chemla} multiplied by $(-1)^n$ under this correspondence after using the conventions spelled out in Remark \ref{Rem:Chemla}. 
Furthermore, this is consistent with Corollary \ref{Cor:Iso-dPcoh} since Chemla proved that Proposition \ref{Pr:Chem-5.8} applied to $\Omega^1_{\cA/B}$ computes the double Poisson cohomology of Pichereau-Van de Weyer, see \cite[Prop.~5.12]{Ch}.

\subsection{Fusion as an operation in cohomology} \label{ss:fusion}

We give a direct application of Chemla's formula from Proposition \ref{Pr:Chemla}. 
For $j=1,2$, let $\cA_j$ be endowed with a $B_j$-linear double Poisson bracket $\dsq{-,-}_j$. Assume that $B_j$ decomposes as $B_j=B_0 \oplus C_j$ (possibly with $B_j=B_0$), and form the \emph{fusion algebra} $\cA:=\cA_1 \ast_{B_0} \cA_2$ which is an algebra over $B:=B_0 \oplus C_1 \oplus C_2$. 
Then Van den Bergh noticed \cite[\S2.5]{VdB1} that $\cA$ is endowed with a unique $B$-linear double Poisson bracket $\dsq{-,-}$, called the \emph{fusion double bracket}, such that 
for $a_1,b_1\in \cA_1$, $a_2,b_2\in \cA_2$
$$\dsq{a_1,b_1}=\dsq{a_1,b_1}_1, \quad 
\dsq{a_2,b_2}=\dsq{a_2,b_2}_2, \quad 
\dsq{a_1,a_2}=0.$$ 
Let us also make an observation: any $n$-bracket on $\cA_1$ can be extended to a $B$-linear $n$-bracket on $\cA$ which vanishes when an argument belongs to $\cA_2$ (and similarly if one swaps $\cA_1$ and $\cA_2$). 
This simple extension yields, for $j=1,2$, a morphism 
\begin{equation} \label{Eq:Mor-Fus}
    \mathrm{ext}_{j,n} : \wBRA_{B_j}(\cA_j)_n \to \wBRA_B(\cA)_n\,, \quad n\geq1.
\end{equation} 
For $n=0$, we also have a natural morphism $\mathrm{ext}_{j,0}:\cA_{j,\sharp} \to \cA_\sharp$.

\begin{lemma} 
Consider the complexes $(\wBRA_B(\cA),\wdd)$ and $(\wBRA_{B_j}(\cA_j),\wdd_j)$, for $j\in \{1,2\}$, equipped with the differentials defined from \eqref{Eq:dP-gen-0}-\eqref{Eq:dP-gen} with $\dsq{-,-}_j$ and $\dsq{-,-}$, respectively. 
Then there is a morphism of complexes 
$$\mathrm{ext}_{j,\bullet}:\wBRA_{B_j}(\cA_j) \to \wBRA_B(\cA)$$
given in each degree by \eqref{Eq:Mor-Fus}. 
\end{lemma}
\begin{proof}
Without loss of generality, we take $j=1$. 
We start in degree $n=0$ with $\bar{a}\in \cA_{1,\sharp}$, which can also be seen as an element of $\cA_\sharp$. Then for any $b_1\in \cA_1$, 
\begin{equation}
    \wdd(\bar{a})(b_1)=-\mult \circ \dsq{a,b_1} = - \mult \circ \dsq{a,b_1}_1 =  \wdd_1(\bar{a})(b_1)\,,
\end{equation}
due to \eqref{Eq:dP-gen-0} and the definition of the fusion bracket. 
For any $b_2\in \cA_2$, $\wdd(\bar{a})(b_2)=-\mult \circ \dsq{a,b_2}=0$ by definition of the fusion bracket; this equals $\mathrm{ext}_{1,1}(\wdd_1(\bar{a}))(b_2)$ since $\wdd_1(\bar{a})\in \wBRA_{B_1}(\cA_1)_1$ is extended to be zero on $\cA_2$. 
Thus, $\wdd\circ \mathrm{ext}_{1,0} = \mathrm{ext}_{1,1}\circ \wdd_1$. 

In degree $n\geq 1$, given $\dgal{-}\in \wBRA_{B_1}(\cA_1)_n$, the extension 
$\mathrm{ext}_{1,n+1}\circ \widehat{d}_1(\dgal{-}) \in \wBRA_{B}(\cA)_{n+1}$ is defined to be only nonzero when evaluated on generators $a_1,\ldots,a_{n+1}$ of $\cA_1$, where it takes the form 
\begin{equation}  
\begin{aligned} 
&\sum_{i=1}^{n+1} (-1)^{ni}\, \sigma^i \, 
\dgal{a_{i+1},\ldots,a_{n+1},a_1,\ldots,a_{i-2}, \dsq{a_{i-1},a_i}_1}_{L}   \\
+&\sum_{i=1}^{n+1} (-1)^{ni}\, \sigma^{i-1} \,   
\dsq{a_i,\dgal{a_{i+1},\ldots,a_{n+1},a_1,\ldots,a_{i-1}} }_{1,L}  \,,
\end{aligned}
\end{equation}
by Chemla's formula \eqref{Eq:Chemla}.  Meanwhile, applying $\wdd$ to the extension of $\dgal{-}$ in $\wBRA_{B}(\cA)_{n+1}$ which we also denote $\dgal{-}$, we get again by Chemla's formula 
\begin{equation}  
\begin{aligned} \label{Eq:Pf-Fus}
&\sum_{i=1}^{n+1} (-1)^{ni}\, \sigma^i \, 
\dgal{a_{i+1},\ldots,a_{n+1},a_1,\ldots,a_{i-2}, \dsq{a_{i-1},a_i} }_{L}   \\
+&\sum_{i=1}^{n+1} (-1)^{ni}\, \sigma^{i-1} \,   
\dsq{a_i,\dgal{a_{i+1},\ldots,a_{n+1},a_1,\ldots,a_{i-1}} }_{L}  \,,
\end{aligned}
\end{equation}
for any $a_1,\ldots,a_{n+1}\in \cA$. If all the $a_i$ belong to $\cA_1$, this is the same formula as above. Assume now that  all the $a_i$ belong to either $\cA_1$ or $\cA_2$, and in addition $a_k\in \cA_2$ for some $k$. Then by definition of the extension of $\dgal{-}$, it must vanish when $a_k$ is in its arguments, so \eqref{Eq:Pf-Fus} becomes 
\begin{equation}  
\begin{aligned} \label{Eq:Pf-Fus2}
&(-1)^{nk}\, \sigma^k \, 
\dgal{a_{k+1},\ldots,a_{n+1},a_1,\ldots,a_{k-2}, \dsq{a_{k-1},a_k} }_{L}  \\
+& (-1)^{n(k+1)}\, \sigma^{k+1} \, 
\dgal{a_{k+2},\ldots,a_{n+1},a_1,\ldots,a_{k-1}, \dsq{a_{k},a_{k+1}} }_{L}   \\
+&(-1)^{nk}\, \sigma^{k-1} \,   
\dsq{a_k,\dgal{a_{k+1},\ldots,a_{n+1},a_1,\ldots,a_{k-1}} }_{L}  \,,
\end{aligned}
\end{equation}
For $\dsq{a_{k-1},a_k}$ in the first term of \eqref{Eq:Pf-Fus2} to be nonzero, we need $a_{k-1}\in \cA_2$. But then  $\dsq{a_{k-1},a_k}=\dsq{a_{k-1},a_k}_2\in \cA_2^{\otimes 2}$ and the first term must be zero. Similarly, the second term is always zero. 
For the third term of \eqref{Eq:Pf-Fus2}, we get by extension that $\dgal{a_{k+1},\ldots,a_{k-1}}$ is either zero or in $\cA_1^{\otimes n}$; in the second case we get zero by applying $\dsq{a_k,-}_L$ by definition of the fusion bracket because $a_k\in \cA_2$. 
We conclude that 
$\dd\circ \mathrm{ext}_{1,n} = \mathrm{ext}_{1,n+1}\circ \wdd_1$, as desired. 
\end{proof}

\begin{corollary} \label{Cor:Fus}
Let $\cA_1$ and $\cA_2$ be equipped with a double Poisson bracket. 
If $\cA=\cA_1\ast_{B_0} \cA_2$  is the fusion algebra equipped with the fusion double Poisson bracket, there are linear maps 
\begin{equation} \label{Eq:Fus-dPH}
  \mathrm{ext}_{1} :  \widehat{\dPH}(\cA_1)\to \widehat{\dPH}(\cA), 
  \quad \mathrm{ext}_{2} :  \widehat{\dPH}(\cA_2)\to \widehat{\dPH}(\cA).
\end{equation}
\end{corollary}
By a straightforward adaptation of these arguments, the same conclusion holds in the non-completed case where there are linear maps 
\begin{equation} \label{Eq:Fus-dPH-b}
  \mathrm{ext}_{1} :  \dPH(\cA_1)\to \dPH(\cA), 
  \quad \mathrm{ext}_{2} :  \dPH(\cA_2)\to \dPH(\cA).
\end{equation}

%%%%%%%%%%% NEW CHAPTER %%%%%%%%%%%%%%%
%%%%%%%%%%% NEW CHAPTER %%%%%%%%%%%%%%%
%%%%%%%%%%% NEW CHAPTER %%%%%%%%%%%%%%%
%%%%%%%%%%% NEW CHAPTER %%%%%%%%%%%%%%%

\chapter{The double quasi-Poisson and gauged double Poisson cases}
\label{CH:qgauge}

We start with Section \ref{ss:gdPCoh} where we reformulate the noncommutative Poisson cohomology introduced by Alekseev, Kawazumi, Kuno and Naef \cite[\S4.2]{AKKN} to make it compatible with the completed double Poisson cohomology that we defined in Chapter \ref{CH:Gen-dPcoh}. 
Then, in Section \ref{ss:dqPCoh}, we introduce a cohomology theory based on Van den Bergh's double quasi-Poisson algebras \cite[\S5]{VdB1}. 
Finally, we define in Section~\ref{ss:RevGDP} the most general gauged double Poisson cohomology, which encompasses all previous constructions and describes a noncommutative analogue of Section \ref{ss:classQPCoh}. 

\medskip 

Throughout the chapter, $\cA$ is an algebra over the semisimple base $B=\oplus_{s\in S} \kk e_s$ made of a complete finite set of orthogonal idempotents. 
For each $s\in S$, we can introduce the double derivation $\Delta_s\in (\mb T^\ast \cA)_1$ satisfying \glslink{DelGauge}{}
\begin{equation} \label{Eq:Deltas}
    \Delta_s: \cA\to \cA^{\otimes 2}, \quad 
    \Delta_s(a)= a e_s\otimes e_s - e_s \otimes e_s a\,, 
\end{equation}
which is called a \emph{gauge element}. 
Using the double Gerstenhaber bracket $\dSN{-,-}$ from Theorem \ref{Thm:dSN}, 
it follows (cf. \cite[\S3.3]{VdB1}) that for any $s\in S$ and $Q\in (\mb T^\ast \cA)_k$, 
\begin{equation} \label{Eq:dSN-gauge}
    \dSN{\Delta_s,Q}=Q e_s \otimes e_s - e_s \otimes e_s Q\,, 
\quad \text{ hence } 
\brSN{Q,\Delta_s}=0\,.
\end{equation} 

\section{Gauged double Poisson cohomology} \label{ss:gdPCoh}

We adapt \cite[\S4.2]{AKKN} by working in full generalities  over the base $B=\oplus_{s\in S} \kk e_s$.  
Let $\cA$ be equipped with a $B$-linear double Poisson bracket $\dsq{-,-}=\dgal{-,-}_P$ defined from some $P\in (\mb T^\ast \cA)_2$ according to Proposition \ref{Pr:MapMu}.  
We introduce for any $k\geq 0$, the $\kk$-linear map%\glslink{iotaDelta}{} 
\begin{equation} \label{Eq:iota-k}
\iota^\Delta_k : (\mb T^\ast \cA)_{k-1}^{|S|}\to (\mb T^\ast \cA)_{\sharp,k}, \quad  
(\alpha_s)_{s\in S} \mapsto \sum_{s\in S} (\alpha_s \Delta_s)_\sharp \, , 
\end{equation} 
and let $\mc D_\cA^k = \operatorname{coker} \iota^\Delta_k$ 
with projection map $\pi^\Delta_k:(\mb T^\ast \cA)_{\sharp,k}\to \mc D_\cA^k$.\glslink{calDA}{}\glslink{piDelta}{} 
When $k=0$, $\iota^\Delta_0:0\to \cA_\sharp$ and $\pi^\Delta_0=\Id_{\cA_\sharp}$. 
We get from $\dd_P$ \eqref{Eq:dP-PVdW} and \eqref{Eq:dSN-gauge} that  
\begin{equation} \label{Eq:dP-gauge}
    \dd_P((Q \Delta_s)_\sharp)=(\brSN{P,Q} \Delta_s)_\sharp \in \im \iota^\Delta_{k+1}\,, \quad Q\in \mb T^\ast \cA, \,\, s\in S.
\end{equation}
Therefore we obtain a well-defined square zero differential on $\mc D_\cA=\bigoplus_{k\geq 0}\mc D_\cA^k$: 
\begin{equation} \label{Eq:dP-AKKN}
    \dd_P : \mc D_\cA^k \to \mc D_\cA^{k+1}, \quad 
    \dd_P( \pi^\Delta_k(Q)) :=\pi^\Delta_{k+1}(\dd_P(Q))\,, 
    \quad k \geq 0.
\end{equation}

Composing each map $\iota^\Delta_k$ with $\mu_k:(\mb T^\ast \cA)_{\sharp,k}\to \wBRA_B(\cA)_k$ from Proposition \ref{Pr:MapMu}, we can also introduce 
$\widehat{\mc D}_\cA^k = \operatorname{coker} \mu_k \circ\iota^\Delta_k$ 
with the projection $\widehat{\pi}^\Delta_k:\wBRA(\cA)_k \to \widehat{\mc D}_\cA^k$. 
By Theorem \ref{Thm:g-dPcoh2}, the maps $\mu$ define a morphism of complexes\footnote{To be precise, there is a morphism if one takes alternating signs for the differential of the second complex. We omit to write those signs since we solely focus on the resulting cohomologies.} 
$$
((\mb T^\ast \cA)_{\sharp},\dd_P)\longrightarrow (\wBRA_B(\cA),\wdd)$$
which implies together with \eqref{Eq:dP-gauge}: 
$$\wdd(\mu_k\circ \iota^\Delta_k( (\mb T^\ast \cA)_{k-1}^{|S|})) \subset \im (\mu_{k+1}\circ \iota^\Delta_{k+1})\,.$$ 
Thus $\wdd$ descends to $\widehat{\mc D}_\cA$ as
\begin{equation}
    \wdd : \widehat{\mc D}_\cA^k \to \widehat{\mc D}_\cA^{k+1}, \quad
    \wdd( \widehat{\pi}^\Delta_k(\dgal{-})) :=\widehat{\pi}^\Delta_{k+1}(\wdd(\dgal{-}))\,,
    \quad k \geq 0,
\end{equation}
and we deduce that the following diagram is commutative:\glslink{wdd}{}
%\pecetta{The diagram has been commented to speed up compilation}
%\begin{comment}
%%% THE NICE FIGURE %%%
\begin{center}
      \begin{tikzpicture}[scale=0.8]
%%%% top row BACK
 \node  (zero) at (-4.7,0) {$0$};
 \node   (T0) at (-2.8,0) {$(\mb T^\ast \cA)_{\sharp,0}$};
 \node   (T1) at (0,0) {$(\mb T^\ast \cA)_{\sharp,1}$};
 \node  (T2) at (3,0) {$(\mb T^\ast \cA)_{\sharp,2}$};
 \node  (T3) at (6,0) {$(\mb T^\ast \cA)_{\sharp,3}$};
 \node  (dots) at (8.5,0) {$\cdots$};
\path[->,>=angle 90,font=\small]  
   (zero) edge  (T0) ;
\path[->,>=angle 90,font=\small]  
   (T0) edge  (T1) ;
\path[->,>=angle 90,font=\small]  
   (T1) edge  (T2) ;
\path[->,>=angle 90,font=\small]  
   (T2) edge  (T3) ;
\path[->,>=angle 90,font=\small]  
   (T3) edge  (dots) ;
%%%% bottom row BACK
 \node[gray]   (Bzero) at (-4.7,-3) {$0$};
 \node[gray]   (B0) at (-2.8,-3) {$\wBRA(\cA)_0$};
 \node[gray]   (B1) at (0,-3) {$\wBRA(\cA)_1$};
 \node[gray]   (B2) at (3,-3) {$\wBRA(\cA)_2$};
 \node[gray]   (B3) at (6,-3) {$\wBRA(\cA)_3$};
 \node[gray]   (Bdots) at (8.5,-3) {$\cdots$}; 
 \path[->,>=angle 90,font=\small,gray]  
   (Bzero) edge  (B0) ;
\path[->,>=angle 90,font=\small,gray]  
   (B0) edge  (B1) ;
\path[->,>=angle 90,font=\small,gray]  
   (B1) edge   (B2) ;
\path[->,>=angle 90,font=\small,gray]  
   (B2) edge   (B3) ;
\path[->,>=angle 90,font=\small,gray]  
   (B3) edge   (Bdots) ;
%%%% vertical arrows BACK
\path[->,>=angle 90,font=\small,gray] 
(T0) edge   (B0) ;
\path[->,>=angle 90,font=\small,gray] 
(T1) edge   (B1) ;
\path[->,>=angle 90,font=\small,gray] 
(T2) edge   (B2) ;
\path[->,>=angle 90,font=\small,gray] 
(T3) edge  (B3) ;
%%%%%%%%%%%%%%% FRONT %%%%%%%%%%%%%%%%%%
%%%% top row Front
 \node  (zeroF) at (-6.4,-1.5) {$0$};
 \node   (T0F) at (-5,-1.5) {$\mc D_\cA^0$};
 \node  (T1F) at (-2,-1.5) {$\mc D_\cA^1$};
 \node  (T2F) at (1,-1.5) {$\mc D_\cA^2$};
 \node  (T3F) at (4,-1.5) {$\mc D_\cA^3$};
 \node  (dotsF) at (7,-1.5) {$\cdots$};
\path[->,>=angle 90,font=\small]  
   (zeroF) edge  (T0F) ;
\path[->,>=angle 90,font=\small]  
   (T0F) edge  (T1F) ;
\path[->,>=angle 90,font=\small]  
   (T1F) edge  (T2F) ;
\path[->,>=angle 90,font=\small]  
   (T2F) edge  (T3F) ;
\path[->,>=angle 90,font=\small]  
   (T3F) edge  (dotsF) ;
\path[->,>=angle 90,font=\small]  
   (T0) edge  (T0F) ;
\path[->,>=angle 90,font=\small]  
   (T1) edge  (T1F) ;
\path[->,>=angle 90,font=\small]  
   (T2) edge  (T2F) ;
\path[->,>=angle 90,font=\small]  
   (T3) edge  (T3F) ;
%%%% bottom row FRONT
 \node  (BzeroF) at (-6.4,-4.5) {$0$};
 \node   (B0F) at (-5,-4.5) {$\widehat{\mc D}_\cA^0$};
 \node  (B1F) at (-2,-4.5) {$\widehat{\mc D}_\cA^1$};
 \node  (B2F) at (1,-4.5) {$\widehat{\mc D}_\cA^2$};
 \node  (B3F) at (4,-4.5) {$\widehat{\mc D}_\cA^3$};
 \node  (BdotsF) at (7,-4.5) {$\cdots$};
\path[->,>=angle 90,font=\small]  
   (BzeroF) edge  (B0F) ;
\path[->,>=angle 90,font=\small]  
   (B0F) edge  (B1F) ;
\path[->,>=angle 90,font=\small]  
   (B1F) edge  (B2F) ;
\path[->,>=angle 90,font=\small]  
   (B2F) edge  (B3F) ;
\path[->,>=angle 90,font=\small]  
   (B3F) edge  (BdotsF) ;
\path[->,>=angle 90,font=\small,gray]  
   (B0) edge  (B0F) ;
\path[->,>=angle 90,font=\small,gray]  
   (B1) edge  (B1F) ;
\path[->,>=angle 90,font=\small,gray]  
   (B2) edge  (B2F) ;
\path[->,>=angle 90,font=\small,gray]  
   (B3) edge  (B3F) ;
%%%% vertical arrows front 
\path[->,>=angle 90,font=\small] 
(T0F) edge   (B0F) ;
\path[->,>=angle 90,font=\small] 
(T1F) edge   (B1F) ;
\path[->,>=angle 90,font=\small] 
(T2F) edge   (B2F) ;
\path[->,>=angle 90,font=\small] 
(T3F) edge  (B3F) ;
%%%%%%% back part kernels
\node   (ker0) at (-0.8,1.5) {$0$};
\node   (ker1) at (2,1.5) {$(\mb T^\ast \cA)_{0}^{|S|}$};
\node  (ker2) at (5,1.5) {$(\mb T^\ast \cA)_{1}^{|S|}$};
\node  (ker3) at (8,1.5) {$(\mb T^\ast \cA)_{2}^{|S|}$};
\path[->,>=angle 90,font=\small]  
   (ker0) edge  node[right] {$\iota^\Delta_0$}  (T0) ;
\path[->,>=angle 90,font=\small]  
   (ker1) edge node[right] {$\iota^\Delta_1$} (T1) ;
\path[->,>=angle 90,font=\small]  
   (ker2) edge node[right] {$\iota^\Delta_2$}  (T2) ;
\path[->,>=angle 90,font=\small]  
   (ker3) edge node[right] {$\iota^\Delta_3$}  (T3) ;
   \end{tikzpicture}
\end{center}
%%% THE NICE FIGURE %%%
%\end{comment}

\begin{definition} \label{Def:gdPH}
The cohomology of the complex $(\mc D_\cA, \dd_P)$ is called the \emph{gauged double Poisson cohomology} of $\cA$ with respect to the noncommutative Poisson bivector $P\in (\mb T^\ast \cA)_2$. 
Explicitly, $\gdPH(\cA)=(\gdPH^k(\cA))_{k \geq 0}$ for\glslink{gdPH}{}  
\begin{equation}
\gdPH^k(\cA):=\frac{\ker\, \dd_P: \mc D_\cA^k \to \mc D_\cA^{k+1}}
{\im \, \dd_P:\mc D_\cA^{k-1}\to \mc D_\cA^k}\,.
\end{equation}
Similarly, the cohomology of the complex $(\widehat{\mc D}_\cA, \wdd)$ is called the \emph{completed gauged double Poisson cohomology} of $\cA$ and is denoted $\widehat{\gdPH}(\cA)=(\widehat{\gdPH}^k(\cA))_{k \geq 0}$ for\glslink{HatgdPH}{}
\begin{equation}
\widehat{\gdPH}^k(\cA):=\frac{\ker\, \wdd: \widehat{\mc D}_\cA^k \to \widehat{\mc D}_\cA^{k+1}}
{\im \, \wdd:\widehat{\mc D}_\cA^{k-1}\to \widehat{\mc D}_\cA^k}\,.
\end{equation}
\end{definition}

The previous commutative diagram yields the following result in cohomology.

\begin{corollary} \label{Cor:dPH-gdPH}
    The following diagram is commutative:  
%\pecetta{this diagram has been commented to speed up compiling}
%\begin{comment}
\begin{center}
      \begin{tikzpicture}
%%%% top row BACK
 \node   (T0) at (-2,1) {$\dPH(\cA)$};
 \node   (T1) at (2,1) {$\widehat{\dPH}(\cA)$};
 \node  (B0) at (-2,-1) {$\gdPH(\cA)$};
 \node  (B1) at (2,-1) {$\widehat{\gdPH}(\cA)$};
\path[->,>=angle 90,font=\small]  
   (T0) edge node[above] {$\mu$} (T1) ;
\path[->,>=angle 90,font=\small]  
   (B0) edge node[below] {$\mu$} (B1) ;
\path[->,>=angle 90,font=\small]  
   (T0) edge node[left] {$\pi^\Delta$} (B0) ;
\path[->,>=angle 90,font=\small]  
   (T1) edge node[right] {$\widehat{\pi}^\Delta$} (B1) ;
   \end{tikzpicture}
\end{center}
%\end{comment}
\end{corollary}

\begin{remark}
When the maps $\mu_k$ are all isomorphisms, the cohomologies $\gdPH(\cA)$ and $\widehat{\gdPH}(\cA)$ are equivalent. 
In general, the maps $\pi^\Delta,\widehat{\pi}^\Delta$ may fail to be injective (see Remark \ref{Rem:dPH-gdPH-inj}) or surjective (see Remark \ref{Rem:dPH-gdPH-surj}).
\end{remark}

\medskip 

We can state a gauged analogue of Lemma \ref{Lem:PV-1}. 
We say that a \emph{gauged double Poisson derivation} is an element $\delta\in {\mc D}^1_\cA$ such that $\dd_P(\delta)=0$, where $\dd_P$ is the differential on ${\mc D}_\cA$ induced by \eqref{Eq:dP-PVdW}.   
\begin{lemma}  \label{Lem:PV-gauge}
The first two cohomology groups of $\gdPH(\cA)$ are such that  
\begin{align}
\gdPH^0(\cA)& \subseteq Z_P(\cA_\sharp)\,, \label{Eq:Lem-gHP0} \\
\gdPH^1(\cA)&=\frac{\{\text{gauged double Poisson derivations}\}}{\pi_1^\Delta(\{\text{double Hamiltonian derivations}\})}\,.
\label{Eq:Lem-gHP1}
\end{align} 
\end{lemma}
\begin{proof}
Using the proof of Lemma \ref{Lem:PV-1}, if $\bar{a} \in \ker \dd_P:{\mc D}^0_\cA\to {\mc D}^1_\cA$, the projection  $\pi^\Delta_1(\dsq{a,-})$ is in the class of $0\in {\mc D}^1_\cA$.
This means that $\dsq{a,-}=\sum_s f_s\Delta_s$ for some $f_s\in \cA$, $s\in S$. 
Therefore, for any $\bar{b}\in \cA_\sharp$,  
$$(\mult \circ \dsq{\bar{a},\bar{b}})_{\sharp}=\sum_s (\mult \circ (f_s\Delta_s)(b))_\sharp 
=\sum_s (b e_s f e_s  - e_s f e_s b)_\sharp =0, $$
showing that $\bar{a}$ is central in $(\cA_\sharp,\mult \circ \dsq{-,-})$. 
For $\gdPH^1(\cA)$, it is by definition given by gauged double Poisson derivations modulo elements $\delta_f\in {\mc D}^1_\cA$ such that $\delta_f=\dd_P(f)$ for some $f\in \cA_\sharp={\mc D}^0_{\mc A}$. Since the differential $\dd_P$ on ${\mc D}_\cA$ is obtained from the one on $(\mb T^\ast \cA)_\sharp$, an element of the later form is an element belonging to the class of a double Hamiltonian derivation (as in Section~\ref{ss:dPcoh}) projected to ${\mc D}^1_\cA$.  
\end{proof}

We finish with a simple criterion for finding possible candidates for cocycles. To state it, note from Proposition \ref{Pr:MapMu} that there is a well-defined map 
\begin{equation}
\begin{aligned} \label{Eq:AKKN-map}
 &\mu_{k,\sharp}:(\mb T^\ast \cA)_{k,\sharp} \to 
 \Hom_\kk(\cA_{\sharp}^k,\cA_\sharp)\,, \\
&\mu_{k,\sharp}(Q_\sharp)(a_{1,\sharp},\ldots,a_{k,\sharp}) 
:=(\mult \circ \dgal{a_1,\ldots,a_k}_Q)_\sharp\,,
\end{aligned}
\end{equation} 
where we take arbitrary lifts $Q\in (\mb T^\ast \cA)_k$ and $a_1,\ldots,a_k\in \cA$, with $k\geq 1$.

\begin{lemma}[Criterion for gauged cocycles] \label{Lem:AKKN}
Let $Q_\sharp \in (\mb T^\ast \cA)_{k,\sharp}$ be such that 
$$\mu_{k+1,\sharp}(\dd_P(Q_\sharp))\neq 0.$$ 
Then the projection $\pi_k^\Delta(Q_\sharp)$ of $Q_\sharp$ to ${\mc D}^k_\cA$ is not a cocycle in $(\mc D_\cA, \dd_P)$.  
\end{lemma}
\begin{proof}
We see from \eqref{Eq:AKKN-map} and \eqref{Eq:MapMu} that 
$\mu_{k+1,\sharp}((\alpha_s \Delta_s)_\sharp))\equiv 0$
for any $s\in S$ and $\alpha_s\in (\mb T^\ast \cA)_{k}$, cf. \cite[Prop.~4.3]{AKKN}. 
If $\pi_k^\Delta(Q_\sharp)$ is a cocycle, then $\dd_P(Q_\sharp)=\sum_s \alpha_s \Delta_s$ for some $\alpha_s\in (\mb T^\ast \cA)_{k}$, leading to a contradiction. 
\end{proof}

\begin{remark}
We do not know under which assumptions the converse of Lemma \ref{Lem:AKKN} holds as the map \eqref{Eq:AKKN-map} may vanish outside $\im \iota_k^\Delta$. 
For example, consider the free algebra $\cA=\kk\langle x_1,\ldots, x_\ell \rangle$ with $\ell>1$ generators. Then, $\mb T^\ast \cA$ 
is generated by the $\ell$ double derivations  
\begin{equation*}
    \del_j\in (\mb T^\ast \cA)_1, \quad \del_j(x_i)=\delta_{ij}\, 1\otimes 1\,, \quad 1\leq i,j\leq \ell. 
\end{equation*}
We can write in this case $\Delta=\sum_{j=1}^\ell [\del_j,x_j]$. For a fixed pair $(i,j)$ of indices, one can check that $\pi_2^\Delta((\del_i\, [\del_j,x_j])_\sharp) \in {\mc D}^2_\cA$ is a nontrivial class. However, plugging \eqref{Eq:dbr-P} in \eqref{Eq:AKKN-map}  yields  for any $1\leq a,b\leq \ell$ 
\begin{align*}
   \mu_{2,\sharp}((\del_i\, [\del_j,x_j])_\sharp)(x_a,x_b)
=&(\delta_{jb}\delta_{ia}+\delta_{ja}\delta_{ib})\, (\mult \circ (x_j\otimes 1 - 1 \otimes x_j))_\sharp  =0\,.
\end{align*}
Hence $(\del_i\, [\del_j,x_j])_\sharp\in \ker \mu_{2,\sharp}$, 
proving that $\im \iota_2^\Delta \subsetneq \ker \mu_{2,\sharp}$. 
\end{remark}

\section{Double quasi-Poisson cohomology}   \label{ss:dqPCoh}

%We introduce a minor variant of Van den Bergh's double quasi-Poisson brackets \cite[\S5]{VdB1}.  
\begin{definition}\label{Def:qPoiss}
An element $P\in (\mb T^\ast \cA)_{\sharp,2}$ is called \emph{quasi-Poisson} if there exist constants $(q_s)\in \kk^S$ such that  
\begin{equation} \label{Eq:qP-PP}
    \brSN{P,P}=\sum_{s\in S} \frac16 q_s (\Delta_s^3)_\sharp \in (\mb T^\ast \cA)_{\sharp,3}\,.
\end{equation}
We then say that the $B$-linear double bracket $\mu_2(P)$ is \emph{quasi-Poisson}.

\noindent More generally, a $B$-linear double bracket $\dsq{-,-}$ is a \emph{double quasi-Poisson bracket} if the associated triple bracket $\dsq{-,-,-}$ \eqref{Eq:dJac} satisfies for any $a,b,c\in \cA$
\begin{equation}
   \begin{aligned} \label{qPabc}
    \dsq{a,b,c}=\sum_{s\in S} \frac{q_s}{4} \Big(
&\, c e_s a \otimes e_s b \otimes e_s  - c e_s a \otimes e_s \otimes b e_s 
- c e_s \otimes a e_s b \otimes e_s
\\
&+ c e_s \otimes a e_s \otimes b e_s - e_s a \otimes e_s b \otimes e_s c 
+ e_s a \otimes e_s \otimes b e_s 
c \\
& + e_s \otimes a e_s b \otimes e_s c - e_s \otimes a e_s \otimes b e_s c \Big)\,.
  \end{aligned}
\end{equation} 
\end{definition} 
\begin{remark} \label{Rem:qP}
    \begin{enumerate}
    \item Van den Bergh's original condition \cite[\S5]{VdB1} consists in taking all $q_s=1$. 
The explicit form \eqref{qPabc} of the triple bracket can be found in \cite{MT14}.
\item If $\dsq{-,-}=\mu_2(P)$, the two conditions \eqref{Eq:qP-PP} and \eqref{qPabc} are equivalent. 
This can be seen by explicitly computing $\mu_3(\Delta_s^3)$ and using \eqref{Eq:dJac-PP}. 
%that $$\dgal{-,-,-}=\frac12\mu_2(\brSN{P,P})\,.$$ cf. \cite[\S4.2]{VdB1}. 
\item By allowing to take $q_s=0$ for each $s$, we recover the theory of double Poisson bracket as well. In that way, we shall introduce a generalization of (completed) double Poisson cohomology.  
    \end{enumerate}
\end{remark}
\begin{lemma} \label{lem:multQP}
    If $\dsq{-,-}$ is quasi-Poisson, then the triple bracket $\dsq{-,-,-}$ defined by  \eqref{Eq:dJac} satisfies 
\begin{equation}
 (\mult \otimes \Id_\cA)\circ \dsq{-,-,-}=0, \,\, \text{ and } 
 (\Id_\cA \otimes \mult)\circ \dsq{-,-,-}=0   \,.
\end{equation}
\end{lemma}
\begin{proof}
    Direct calculation using \eqref{qPabc}.  
\end{proof}

\begin{proposition}[\cite{CW}, \S5.2] \label{Pr:dquasiPcoh1} 
Assume that $P$ is quasi-Poisson. 
Define $\dd_P=\brSN{P,-}:(T^\ast \cA)_{\sharp} \to (T^\ast \cA)_{\sharp}$. 
Then $\dd_P$ is a square-zero differential of degree $+1$. 
\end{proposition}
\begin{proof}
    This is a simple computation which we prove in greater generality as part of Remark \ref{Rem:Delta3} below. 
\end{proof}

We deduce that we can form the complex $((\mb T^\ast \cA)_\sharp , \dd_P)$ also in the quasi-Poisson case. More generally, the quasi-Poisson property is also suitable for the completed case as follows. 

\begin{proposition} \label{Pr:dquasiPcoh2} 
Assume that $\dsq{-,-}$ is a $B$-linear double quasi-Poisson bracket. 
Then,  the conclusion of Theorem \ref{Thm:g-dPcoh1} holds, i.e., the (degree $+1$) operations 
$$\wdd:\wBRA_B(\cA) \to \wBRA_B(\cA)$$ given by \eqref{Eq:dP-gen-0}--\eqref{Eq:dP-gen} define a square-zero differential on the complex $\wBRA_B(\cA)$. 
\end{proposition}

\subsection{Proof of Proposition \ref{Pr:dquasiPcoh2}}

By the first part of Theorem \ref{Thm:g-dPcoh1}, we only need to prove the square-zero property of $\wdd$. 
Given $\bar{a}\in \cA_\sharp$, $\wdd^2(\bar{a})$ evaluated on $b\otimes c \in \cA^{\otimes 2}$ yields \eqref{Eq:dP2-abc}, where the two terms vanish by Lemma~\ref{lem:multQP}. 

Given an $n$-bracket $\dgal{-} \in \wBRA_B(\cA)_n$, we can reproduce most of the computations made in the proof of Theorem \ref{Thm:g-dPcoh1}. There, only Steps i and ii rely on the double Jacobi identity $\dsq{-,-,-}=0$ for the triple bracket \eqref{Eq:dJac} associated with $\dsq{-,-}$. 
In the present case, we only have the weaker rule \eqref{qPabc}. To be precise, what we get from Step i is: 
\begin{align*}
    \mathrm{S}_{\mathrm{i}}:=& 
\sum_{t\in \Z_{n+2}}(-1)^{(n+1)t} \sigma^t \circ D_{\mathrm{i}} \circ \sigma^{-t} \,, \\
D_{\mathrm{i}}(a_1,\ldots,a_{n+2}):=&
 \dgal{a_1,\ldots,a_{n-1},-}_L (\dsq{a_n, a_{n+1}, a_{n+2}} );
\end{align*}
and what we get from Step ii is: 
\begin{align*}
    \mathrm{S}_{\mathrm{ii}}:=& 
\sum_{t\in \Z_{n+2}}(-1)^{(n+1)t} \sigma^t \circ D_{\mathrm{ii}} \circ \sigma^{-t} \,, \\
D_{\mathrm{ii}}(a_1,\ldots,a_{n+2}):=& 
-   \dsq{a_1,a_2,-}_L (\dgal{a_{3},\ldots,a_{n+2}} )\,.
\end{align*}
We shall verify that $\mathrm{S}_{\mathrm{i}}+\mathrm{S}_{\mathrm{ii}}=0$. 
Firstly, we use \eqref{qPabc} and the derivation rule \eqref{Eq:nbr-DerOut} to write 
$D_{\mathrm{i}}=\sum_{s\in S}\frac{q_s}{4} D_{\mathrm{i}}^{(s)}$, where 
\begin{equation} \label{Eq:qP-pf1}
    \begin{aligned}
    D_{\mathrm{i}}^{(s)}(a_1,\ldots,a_{n+2}) 
=&\dgal{a_1,\ldots,a_{n-1},a_{n+2}} e_s a_n \otimes e_s a_{n+1} \otimes e_s    \\
&+ a_{n+2}e_s \dgal{a_1,\ldots,a_{n-1},a_n} \otimes e_s a_{n+1} \otimes e_s \\
%%% 2 lines above: from first term in \eqref{qPabc}
&-\dgal{a_1,\ldots,a_{n-1},a_{n+2}} e_s a_n \otimes e_s \otimes a_{n+1} e_s  \\  
&- a_{n+2}e_s \dgal{a_1,\ldots,a_{n-1},a_n} \otimes e_s \otimes a_{n+1} e_s \\
%%% 2 lines above: from second term in \eqref{qPabc}
&-\dgal{a_1,\ldots,a_{n-1},a_{n+2}} e_s  \otimes a_n e_s a_{n+1} \otimes e_s  \\ 
%%% line above: from third term in \eqref{qPabc}
&+\dgal{a_1,\ldots,a_{n-1},a_{n+2}} e_s  \otimes a_n e_s \otimes a_{n+1} e_s  \\ 
%%% line above: from fourth term in \eqref{qPabc}
&- e_s \dgal{a_1,\ldots,a_{n-1},a_n} \otimes e_s a_{n+1} \otimes e_s a_{n+2} \\
%%% line above: from fifth term in \eqref{qPabc}
&+ e_s \dgal{a_1,\ldots,a_{n-1},a_n} \otimes e_s \otimes a_{n+1} e_s a_{n+2} \,.
%%% line above: from sixth term in \eqref{qPabc}
    \end{aligned}
\end{equation}
Let us look at the 4 terms in $(-1)^{n+1}\,(\sigma \circ D_{\mathrm{i}}^{(s)} \circ \sigma^{-1}) (a_1,\ldots,a_{n+2})$ that correspond to the first, third, fifth and sixth terms in \eqref{Eq:qP-pf1} (these are the 4 terms containing the factor $\dgal{a_1,\ldots,a_{n-1},a_{n+2}}$). They are given by 
\begin{align*}
(-1)^{n+1}\,\sigma \circ& \Big(  \dgal{a_2,\ldots,a_{n},a_{1}} e_s a_{n+1} \otimes e_s a_{n+2} \otimes e_s    \\
&-\dgal{a_2,\ldots,a_{n},a_{1}} e_s a_{n+1} \otimes e_s \otimes a_{n+2} e_s  \\ 
&-\dgal{a_2,\ldots,a_{n},a_{1}} e_s  \otimes a_{n+1} e_s a_{n+2} \otimes e_s  \\ 
&+\dgal{a_2,\ldots,a_{n},a_{1}} e_s  \otimes a_{n+1} e_s \otimes a_{n+2} e_s \Big)  \,, 
  \end{align*}
or, after using the cyclic skewsymmetry \eqref{Eq:nbr-Cycl} of $\dgal{-}_Q^{(n)}$, 
\begin{align*}
&e_s \otimes  (\sigma^{-1} \dgal{a_1,\ldots,a_{n}}) e_s a_{n+1} \otimes e_s a_{n+2}    \\
-& a_{n+2} e_s \otimes (\sigma^{-1} \dgal{a_1,\ldots,a_{n}}) e_s a_{n+1} \otimes e_s \\ 
-&e_s \otimes (\sigma^{-1} \dgal{a_1,\ldots,a_{n}}) e_s  \otimes a_{n+1} e_s a_{n+2}   \\ 
+&a_{n+2} e_s \otimes (\sigma^{-1} \dgal{a_1,\ldots,a_{n}}) e_s  \otimes a_{n+1} e_s \,.
\end{align*}
(Note that in those expressions, $\sigma^{-1} =\sigma_{(1,\ldots,n)}^{-1} $.)
Since we are only interested in what happens for $\mathrm{S}_{\mathrm{i}}$, we can replace the first, third, fifth and sixth terms in \eqref{Eq:qP-pf1} by the terms that we have just computed, i.e., we can replace above $D_{\mathrm{i}}^{(s)}$ by $\widetilde{D}_{\mathrm{i}}^{(s)}$ defined through 
\begin{subequations} \label{Eq:qP-pf2}
 \begin{align} 
\widetilde{D}_{\mathrm{i}}^{(s)}(a_1,\ldots,a_{n+2}) 
=& a_{n+2}e_s \dgal{a_1,\ldots,a_n} \otimes e_s a_{n+1} \otimes e_s \label{Eq:qP-pf2a}\\
&- a_{n+2}e_s \dgal{a_1,\ldots,a_n} \otimes e_s \otimes a_{n+1} e_s \label{Eq:qP-pf2b} \\
&- e_s \dgal{a_1,\ldots,a_n} \otimes e_s a_{n+1} \otimes e_s a_{n+2} \label{Eq:qP-pf2c} \\
&+ e_s \dgal{a_1,\ldots,a_n} \otimes e_s \otimes a_{n+1} e_s a_{n+2} \label{Eq:qP-pf2d} \\ 
&+e_s \otimes  (\sigma^{-1} \dgal{a_1,\ldots,a_{n}}) e_s a_{n+1} \otimes e_s a_{n+2}   \label{Eq:qP-pf2e}  \\
&- a_{n+2} e_s \otimes (\sigma^{-1} \dgal{a_1,\ldots,a_{n}}) e_s a_{n+1} \otimes e_s  \label{Eq:qP-pf2f} \\ 
&-e_s \otimes (\sigma^{-1} \dgal{a_1,\ldots,a_{n}}) e_s  \otimes a_{n+1} e_s a_{n+2}   \label{Eq:qP-pf2g} \\ 
&+a_{n+2} e_s \otimes (\sigma^{-1} \dgal{a_1,\ldots,a_{n}}) e_s  \otimes a_{n+1} e_s \,.  \label{Eq:qP-pf2h}
 \end{align}
\end{subequations}
Secondly, we shall replace $D_{\mathrm{ii}}$ by $\widetilde{D}_{\mathrm{ii}}:=\sigma^{-2}\circ D_{\mathrm{ii}}\circ \sigma^2$ as it leaves $ \mathrm{S}_{\mathrm{ii}}$ unchanged. 
Denoting $\dgal{a_1,\ldots,a_{n-1},a_n}=:c\otimes C_{[2\cdots n]} \in \cA \otimes \cA^{\otimes (n-1)}$ for short, we use \eqref{qPabc} and the derivation rule \eqref{Eq:nbr-DerOut} to write $\widetilde{D}_{\mathrm{ii}}=\sum_{s\in S}\frac{q_s}{4} \widetilde{D}_{\mathrm{ii}}^{(s)}$ with 
\begin{subequations} \label{Eq:qP-pf3}
    \begin{align}
   \widetilde{D}_{\mathrm{ii}}^{(s)}(a_1,\ldots,a_{n+2}) 
= - \sigma^{-2} \circ (& \dsq{a_{n+1},a_{n+2},c} \otimes C_{[2\cdots n]} ) \nonumber  \\
=- \sigma^{-2} \circ \Big( 
&  c e_s a_{n+1} \otimes e_s a_{n+2} \otimes e_s \otimes C_{[2\cdots n]} \label{Eq:qP-pf3a}  \\ %%1
& - c e_s a_{n+1} \otimes e_s \otimes a_{n+2} e_s \otimes C_{[2\cdots n]} \label{Eq:qP-pf3b} \\  %% 2
&- c e_s \otimes a_{n+1} e_s a_{n+2} \otimes e_s \otimes C_{[2\cdots n]} \label{Eq:qP-pf3c} \\  %% 3
&+ c e_s \otimes a_{n+1} e_s \otimes a_{n+2} e_s \otimes C_{[2\cdots n]} \label{Eq:qP-pf3d} \\  %% 4
&- e_s a_{n+1} \otimes e_s a_{n+2} \otimes e_s c \otimes C_{[2\cdots n]} \label{Eq:qP-pf3e} \\  %% 5
&+ e_s a_{n+1} \otimes e_s \otimes a_{n+2} e_s c \otimes C_{[2\cdots n]} \label{Eq:qP-pf3f} \\  %% 6
&+ e_s \otimes a_{n+1} e_s a_{n+2} \otimes e_s c \otimes C_{[2\cdots n]} \label{Eq:qP-pf3g} \\  %% 7
&- e_s \otimes a_{n+1} e_s \otimes a_{n+2} e_s c \otimes C_{[2\cdots n]} \Big)\,. \label{Eq:qP-pf3h}%% 8 
    \end{align}
\end{subequations}
Now, we remark the cancellations 
\begin{align*}
 &\eqref{Eq:qP-pf2a}+\eqref{Eq:qP-pf3f}=0, \quad 
 \eqref{Eq:qP-pf2b}+\eqref{Eq:qP-pf3h}=0, \quad
 \eqref{Eq:qP-pf2c}+\eqref{Eq:qP-pf3e}=0, \\
 & \eqref{Eq:qP-pf2d}+\eqref{Eq:qP-pf3g}=0, \quad 
 \eqref{Eq:qP-pf2e}+\eqref{Eq:qP-pf3a}=0, \quad 
 \eqref{Eq:qP-pf2f}+\eqref{Eq:qP-pf3b}=0, \\
 & \eqref{Eq:qP-pf2g}+\eqref{Eq:qP-pf3c}=0, \quad
 \eqref{Eq:qP-pf2h}+\eqref{Eq:qP-pf3d}=0,
\end{align*}
which imply $\widetilde{D}_{\mathrm{i}}^{(s)}+\widetilde{D}_{\mathrm{ii}}^{(s)}=0$. 
This yield the desired identity $\mathrm{S}_{\mathrm{i}}+\mathrm{S}_{\mathrm{ii}}=0$.  \qed 

\subsection{The cohomology}

From Propositions \ref{Pr:dquasiPcoh1} and \ref{Pr:dquasiPcoh2}, we get complexes with square-zero differential and we can therefore adapt Definitions \ref{Def:dPH} and \ref{Def:dPH-comp}. 

\begin{definition} \label{Def:quasi-dPH} 
If $P\in (\mb T^\ast \cA)_{\sharp,2}$ is quasi-Poisson, the cohomology of the complex $((\mb T^\ast \cA)_\sharp , \dd_P)$ is called the\footnote{Since our definition of double quasi-Poisson bracket encompasses that of a double Poisson bracket, we keep the same name for the cohomology rather than using ``double \emph{quasi-}Poisson cohomology''.} \emph{double Poisson cohomology} of $\cA$ with respect to $P$, and it is denoted $\dPH(\cA)$. 

If $\dsq{-,-} \in \wBRA_B(\cA)_2$ is quasi-Poisson, 
the cohomology of the complex $(\wBRA_B(\cA),\wdd)$ is called the \emph{completed double Poisson cohomology} of $\cA$ with respect to $\dsq{-,-}$, and it is denoted $\widehat{\dPH}(\cA)$. 
\end{definition}

The first part of Theorem \ref{Thm:g-dPcoh2} is independent of the fact that the element $P\in (\mb T^\ast \cA)_2$ is Poisson, as can be seen by inspecting its proof. Thus we deduce the following result in the quasi-Poisson case. 

\begin{corollary} \label{Cor:Iso-dquasiPcoh}
Assume that $P$ is quasi-Poisson. 
Then $(\mu_\ell)_{\ell\geq 0}$ descends to a $\kk$-linear map $\dPH(\cA)\to \widehat{\dPH}(\cA)$ in cohomology.
In particular, under the assumptions of Proposition \ref{Pr:MapMu-Iso}, the double Poisson cohomology $\dPH(\cA)$ and its completed version $\widehat{\dPH}(\cA)$ are isomorphic. 
\end{corollary}

It is also possible to define gauged double Poisson cohomology as in Section~\ref{ss:gdPCoh} when $P$ is quasi-Poisson. We shall proceed to do so in even greater generalities as part of the next section.

\section{Gauged double Poisson cohomology revisited}  \label{ss:RevGDP}

\begin{definition}%[\cite{VdB1},\S5]
\label{Def:gaugDBR} 
An element $P\in (\mb T^\ast \cA)_{\sharp,2}$ is called \emph{gauged Poisson} if there exists $R_s\in (\mb T^\ast \cA)_2$ for each $s\in S$ such that 
\begin{equation} \label{Eq:gaugDBR}
    \brSN{P,P}=\sum_{s\in S} (\Delta_s R_s)_\sharp \in (\mb T^\ast \cA)_{\sharp,3}.
\end{equation}
The corresponding $B$-linear double bracket $\mu_2(P)$ is then said to be a \emph{double gauged Poisson bracket}. 
\end{definition} 

Recall the maps 
$\iota^\Delta_k : (\mb T^\ast \cA)_{k-1}^{|S|}\to (\mb T^\ast \cA)_{\sharp,k}$
from \eqref{Eq:iota-k}. 

\begin{lemma} \label{lem:dqP2iota}
Assume that $P\in (\mb T^\ast \cA)_{\sharp,2}$ is a gauged Poisson element. 
Define the degree $+1$ maps $\dd_P:=\brSN{P,-}:(\mb T^\ast \cA)_{\sharp} \to (\mb T^\ast \cA)_{\sharp}$.  
Then $\im \dd_P^2 \subseteq \im \iota^\Delta$. 
\end{lemma}
\begin{proof}
From the proof of Proposition \ref{Prop:PVdW-diff}, $\dd_P^2(Q)=\frac12 \brSN{\brSN{P,P},Q}$ for any $Q\in (\mb T^\ast \cA)_{\sharp,n}$. 
By assumption, $\brSN{P,P}$ has the form \eqref{Eq:gaugDBR} for some elements $R_s$, therefore the statement holds if we can show  
\begin{align*}
    \brSN{\Delta_s R,Q} \in \im \iota_{n+2}^\Delta, \quad \forall R_s \in (\mb T^\ast \cA)_2, \,\, Q\in (\mb T^\ast \cA)_n.
\end{align*}
We first compute with the graded double bracket $\dSN{-,-}$ that 
\begin{align*}
\dSN{\Delta_s R_s,Q}&=
\Delta_s\star \dSN{R_s,Q}
+ \dSN{\Delta_s,Q}\star R_s \\
&=\Delta_s\star \dSN{R_s,Q} 
+ (Q e_s R_s  \otimes e_s - e_s R_s \otimes e_s Q) \,, 
\end{align*}
using \eqref{Eq:gdb-Lleib} then \eqref{Eq:dSN-gauge} (recalling the sign rules of Subsection~\ref{sec:1.1}). 
Taking the product, we get modulo graded commutators in $ (\mb T^\ast \cA)_{\sharp,n+2}$, 
\begin{equation} \label{Eq:gaugDSN-1}
    \brSN{\Delta_s R_s , Q} =(-1)^{|\dSN{R_s,Q}'|} (\dSN{R_s,Q}' \Delta_s \dSN{R_s,Q}'')_\sharp . 
\end{equation}
This is clearly in the image of 
$\iota^\Delta_{n+1}:(\alpha_s)_{s} \mapsto \sum_{s\in S} (\alpha_s \Delta_s)_\sharp$. 
\end{proof}

\begin{remark} \label{Rem:Delta3}
If $R_s=\Delta_s \delta_s$ with $\delta_s\in \DDer_B(\cA)$, we get the following from \eqref{Eq:gdb-Lleib}, \eqref{Eq:dSN-gauge} and  \eqref{Eq:gaugDSN-1} 
with $Q\in (\mb T^\ast \cA)_n$
    \begin{equation} \label{Eq:gaugDSN-2}
    \brSN{\Delta_s^2 \delta_s , Q} = (\dSN{\delta,Q}' \Delta_s^2 \dSN{\delta,Q}'')_\sharp \in (\mb T^\ast \cA)_{\sharp,n+2}.
\end{equation} 
If furthermore $\delta_s=\Delta_s$ is the gauge element, \eqref{Eq:dSN-gauge} yields $ \brSN{\Delta_s^3 , Q}=0$.
\end{remark}

As a consequence of Lemma \ref{lem:dqP2iota}, the differential of the complex $((\mb T^\ast \cA)_{\sharp},\dd_P)$ squares to zero modulo $\im \iota^\Delta$. We deduce the following from the upper part of the diagram in Section \ref{ss:gdPCoh}. 

\begin{corollary}
Assume that $P\in (\mb T^\ast \cA)_{\sharp,2}$ is a gauged Poisson element. 
Then the linear operation (of degree $+1$) $\dd_P$ on $(\mb T^\ast \cA)_{\sharp}$ descends to a square-zero differential on $\mc D_\cA$.      
\end{corollary}

In view of this corollary, we can repeat Definition \ref{Def:gdPH} in the current setting and introduce the gauged double Poisson cohomology $\gdPH(\cA)$ of $(\mc D_\cA,\dd_P)$ for $P$ a gauged Poisson element. 

\begin{remark} \label{Rem:Compl-gdP}
We can draw the full diagram depicted in Section \ref{ss:gdPCoh} for a gauged Poisson element $P\in (\mb T^\ast \cA)_{\sharp,2}$. Indeed, $\mu$ defines a morphism of complexes $\mb T^\ast \cA \to \wBRA_B(\cA)$ as in Theorem \ref{Thm:g-dPcoh2} even if the differentials are not square-zero. This is because the proof does not require for $P\in (\mb T^\ast \cA)_{\sharp,2}$ to satisfy any assumption (except for inducing a map in cohomology). 
The induced differential on $\widehat{\mc{D}}_\cA$ is square-zero under mild assumptions, e.g. if each $\mu_\ell$ is an isomorphism since the complex becomes isomorphic to $\mc D_\cA$ (and $\dd_P^2=0$ there by the corollary). However, we do not think that the square-zero property is true in full generalities: we conjecture that one could build an $n$-bracket $\dgal{-} \in \wBRA_B(\cA)_n$, $k\geq1$, with nonzero $\wdd^2(\dgal{-})\in \wBRA_B(\cA)_{n+2}$ not in the image of $\mu_{n+2}$ (for a suitably chosen pair $(\cA,P)$). This would guarantee that the class of $\wdd^2(\dgal{-})$ in $\widehat{\mc{D}}_\cA^{n+2}$ is nonzero. 
\end{remark}

%%%%%%%%%%% NEW CHAPTER %%%%%%%%%%%%%%%
%%%%%%%%%%% NEW CHAPTER %%%%%%%%%%%%%%%
%%%%%%%%%%% NEW CHAPTER %%%%%%%%%%%%%%%
%%%%%%%%%%% NEW CHAPTER %%%%%%%%%%%%%%%

\chapter{Computations of double Poisson cohomologies}
\label{CH:dPA-examples}

We study the various cohomology theories defined previously in specific cases: 
for one generator where we can study them all, 
for two generators where we can highlight differences between the double Poisson cohomology of Pichereau and Van de Weyer \cite{PV} and the completed version introduced earlier, 
and for the constant non-degenerate case in an arbitrary number of variables. 
We finish by commenting all the cases that were previously known.

\section{Double Poisson algebras with one generator}

\subsection{Double Poisson brackets on \texorpdfstring{$\kk[x]$}{k[x]}} 
\label{ss:dP-coh-x}

We denote by  $\dsq{-,-}$ the double bracket on $\kk[x]$ from \eqref{Eq:dbr-xlin}, which is determined by 
\begin{equation} \label{Eq:dbr-xlin2}
   \dsq{x,x}= 
   \nu (x^2 \otimes x - x \otimes x^2)
   +\mu (x^2\otimes 1 - 1 \otimes x^2)
   + \lambda (x\otimes 1 - 1 \otimes x)\,,
\end{equation}
for $\lambda,\mu,\nu \in \kk$. It is Poisson when $\lambda\nu-\mu^2=0$ (cf. Example \ref{Ex:Classifkx}), which we assume from now on. As $\kk[x]$ is smooth, $\dsq{-,-} =\mu_2(P)$ for some $P\in (T^\ast \kk[x])_2$. Using $\partial_x\in \DDer(\kk[x])$ defined in \eqref{DDer-delx}, we easily see from \eqref{Eq:dbr-P} that 
\begin{equation} \label{Eq:dbr-kxP}
P=\lambda \, x\partial_x\partial_x
+\mu\, x^2 \partial_x\partial_x
+\nu \, x^2 \partial_x x \partial_x\,,
\end{equation}
which is unique modulo graded commutators. 

We are going to compute the first few double Poisson cohomology groups of $\dsq{-,-}$, noting that $\dPH(\kk[x])\simeq \widehat{\dPH}(\kk[x])$ by Corollary \ref{Cor:Iso-dPcoh}.  Hence we can use the complex introduced in Theorem \ref{Thm:g-dPcoh1} for all our computations. 

\begin{remark}
    Let us emphasize that calculations in the complex $\wBRA(\cA)$ are in general more cumbersome than in $(\mb T^* \cA)_\sharp$. However performing computations in $\wBRA(\kk[x])$ will allow us to deduce many results in the truncated case $\cA=\kk[x]/(x^k)$ presented in the next subsection, where Pichereau-Van de Weyer's double Poisson cohomology is \emph{not} defined. 
\end{remark}

\begin{remark} \label{Rem:Condkx}
Apart from the trivial case of the zero double Poisson bracket, 
there are essentially two cases to consider when computing the higher cohomology groups. 
Indeed, up to performing a linear transformation $x\mapsto x+\alpha$ for $\alpha \in \kk$ and a scaling $x\mapsto \alpha x$ for $\alpha\in \kk^\times$, the double Poisson bracket can be chosen to be either $\dgal{-,-}_3$ (if $\nu$ was nonzero) or $\dgal{-,-}_1$ (if $\nu$ was zero). 
To keep track of these transformations, we either assume that $\nu\neq0$, $\mu\in \kk$ and $\lambda=\mu^2/\nu$, or that $\mu=\nu=0$ and $\lambda\neq 0$. 
\end{remark}

\subsubsection{The group $\widehat{\dPH}^0(\kk[x])$}

\begin{proposition} \label{Pr:dPH-kx0}
For the double Poisson bracket \eqref{Eq:dbr-xlin2}, 
$\widehat{\dPH}^0(\kk[x])=\kk[x]$.
\end{proposition}
\begin{proof}
Since $\kk[x]$ is commutative, it coincides with its abelianization $\kk[x]_\sharp$. For any $r\geq 1$, we get from \eqref{Eq:dP-gen-0}  
$$\wdd(x^r)(x)=-\mult \circ \dsq{x^r,x}=-r \,\dsq{x,x}' x^{r-1} \dsq{x,x}''\,, $$
which vanishes since $\mult \circ \dsq{x,x}_{a,b}=0$ for any $a,b \geq 0$, cf. \eqref{Eq:dbr-ab}.  
Any $f\in \kk[x]$ is of the form $\sum_{r\geq 0}\alpha_r x^r$, $\alpha_r\in \kk$, so it belongs to $\ker \wdd$. 
\end{proof}
\begin{remark}
This result agrees with \eqref{Eq:Lem-HP0} in Lemma \ref{Lem:PV-1}: given the double Poisson bracket $\dsq{-,-}$ of the form \eqref{Eq:dbr-xlin2}, we have seen that $\mult \circ \dsq{f,-}=0$ for any $f$, and therefore $\dsq{-,-}$ descends to the trivial Lie bracket on $\kk[x]_\sharp=\kk[x]$. Thus $Z_P(\kk[x]_\sharp;\kk[x])=Z_P(\kk[x])=\kk[x]$. 
\end{remark}

\subsubsection{The group $\widehat{\dPH}^1(\kk[x])$}

For any integer $\ell\geq 0$, we introduce $\theta_\ell\in \Der(\kk[x])$ as the derivation satisfying $\theta_\ell(x)=x^\ell$. We simply write $\Id_{\kk[x]}$ as $\Id$. 
Recalling $\partial_x$ defined in \eqref{DDer-delx}, 
we note that $\theta_\ell=\mu_1(x^\ell \partial_x)$, where we see $x^\ell \partial_x$ as an element in $(T^\ast \kk[x])_{\sharp,1}$.

\begin{proposition}  \label{Pr:dPH-kx1} 
For the double Poisson bracket \eqref{Eq:dbr-xlin2}, we have 
\begin{align*}
 \widehat{\dPH}^1(\kk[x])&=\kk\left(\nu\theta_2+2\mu \theta_1 +\lambda \theta_0\right) \,,\\
    \dPH^1(\kk[x])&=\kk \left(\nu x^2 \partial_x+2 \mu x \partial_x +\lambda \partial_x  \right)  \,.
\end{align*}
\end{proposition}
\begin{proof}  
The results for $\dPH^1(\kk[x])$ directly follow from the ones for $\widehat{\dPH}^1(\kk[x])$ by applying the inverse of $\mu_1$. Thus, recalling Remark \ref{Rem:Condkx},  it suffices to show: 
\begin{itemize}
\item $\widehat{\dPH}^1(\kk[x])=\kk\left(\theta_2+2\frac{\mu}{\nu} \theta_1 + \frac{\lambda}{\nu} \theta_0  \right)$ 
if $\nu\neq 0$. 
\item $\widehat{\dPH}^1(\kk[x])=\kk \theta_0$ if $\nu=\mu =0$, $\lambda\neq 0$.  
\end{itemize}
First, we note that $\widehat{\dPH}^1(\kk[x])=\ker \{\wdd:\Der(\kk[x])\to \wBRA(\kk[x])_2\}$ because the image of $\wBRA(\kk[x])_0$ under $\wdd$ is zero by Proposition \ref{Pr:dPH-kx0}. 
If $\theta\in  \Der(\kk[x])$, the double bracket $\wdd(\theta)$ is completely determined by its value $\wdd(\theta)(x,x)$ due to the derivation rules enjoyed by a double bracket. For $\ell \geq 1$, we compute using \eqref{Eq:dP-gen} and the cyclic skewsymmetry of $\dsq{-,-}$ that 
\begin{align*}
    \wdd(\theta_\ell)(x,x)=&
(\theta_\ell \otimes \Id + \Id \otimes \theta_\ell) (\dsq{x,x}) 
-\dsq{x,\theta_\ell(x)}  -\dsq{\theta_\ell(x),x} \, .
\end{align*}
Since $\kk[x]$ is commutative, 
$$\dsq{x^\ell,x}=\sum_{\lambda=0}^{\ell-1} (x^\lambda \otimes x^{\ell-\lambda-1}) \dsq{x,x}\,,$$ 
and the same expansion holds for $\dsq{x,x^\ell}$. 
As $\theta_\ell(x^a)=a x^{a+\ell-1}$, a short computation then yields 
\begin{equation} \label{Eq:dPtheta}
    \wdd(\theta_\ell)(x,x)=
\nu (x^2 \otimes x^\ell - x^\ell \otimes x^2) 
+2 \mu (x \otimes x^\ell - x^\ell \otimes x) 
+ \lambda (1 \otimes x^\ell - x^\ell \otimes 1) \, .
\end{equation}
This identity also holds for $\ell=0$.

An arbitrary $\theta\in  \Der(\kk[x])$ can be written as $\theta=\sum_{\ell=0}^m c_\ell \theta_\ell$ with $c_m\neq 0$ for some $m \geq 0$. 
In the case $\nu =\mu = 0$, $\wdd(\theta)(x,x)$ contains the nonzero term $\lambda c_m 1\otimes x^m$ if $m\geq 1$ so $\wdd(\theta)(x,x)\neq 0$ (the $x^a \otimes x^b$, $a,b\geq 0$, form a basis) and thus  $\theta\notin \ker \wdd$. We directly have that  $\wdd(\theta_0)(x,x)=0$ by \eqref{Eq:dPtheta}  and we find in this way $\widehat{\dPH}^1(\kk[x])=\kk \theta_0$. 
In the case $\nu \neq 0$, $\wdd(\theta)(x,x)$ contains the nonzero term $\nu c_m x^2 \otimes x^m$ if $m\geq 3$ from which $\theta\notin \ker \wdd$. If $m\leq 2$, we have 
\begin{align*}
\wdd(\theta)(x,x)=& 
+(2 \mu c_2-\nu c_1) (x \otimes x^2 - x^2 \otimes x) 
+ (\lambda c_2- \nu c_0) (1 \otimes x^2 - x^2 \otimes 1) \\
&+(2 \mu c_0-\lambda c_1) (x \otimes 1 - 1\otimes x) \,.
\end{align*}
(We set $c_2=0$ if $m=0,1$ and $c_1=0$ if $m=0$.) We deduce that elements of $\widehat{\dPH}^1(\kk[x])=\ker \wdd$ must be of the form $c_2 (\theta_2+2\frac{\mu}{\nu} \theta_1 + \frac{\lambda}{\nu} \theta_0)$. 
\end{proof}

\subsubsection{The group $\widehat{\dPH}^2(\kk[x])$}

Recall that on $\kk[x]$, any double bracket (i.e. element of $\wBRA(\kk[x])$) is a linear combination of $\dgal{-,-}_{a,b}$ \eqref{Eq:dbr-ab}. 

\begin{proposition}  \label{Pr:dPH-kx2} 
For the double Poisson bracket \eqref{Eq:dbr-xlin2}, we have 
$$\widehat{\dPH}^2(\kk[x])= \{0\}, \qquad \dPH^2(\kk[x])= \{0\}\,.$$
\end{proposition} 
\begin{proof}
Let us start by noting from \eqref{Eq:dPtheta} that $\wdd(\Der(\kk[x]))$ is the $\kk$-linear span of the following double brackets: 
\begin{enumerate}
    \item[($\nu=0$)] \quad $\dgal{-,-}_{\ell,0}$ for $\ell\geq 1$ if $\mu=\nu=0$;  
    \item[($\nu\neq0$)] \quad  $\dgal{-,-}_{\ell,2}+2\frac{\mu}{\nu}\dgal{-,-}_{\ell,1}+\frac{\lambda}{\nu}\dgal{-,-}_{\ell,0}$ for $\ell\geq 0$ if $\nu\neq 0$.
\end{enumerate} 
Furthermore, in the case $\nu \neq 0$, we note using $\dgal{-,-}_{a,b}=-\dgal{-,-}_{b,a}$ that the elements corresponding to $\ell=0$ and $\ell=1$ are respectively given by
$-\dgal{-,-}_{2,0} -\frac{2\mu}{\nu} \dgal{-,-}_{1,0}$ and 
$-\dgal{-,-}_{2,1} +\frac{\lambda}{\nu} \dgal{-,-}_{1,0}$; the element corresponding to $\ell=2$ is  $\frac{2\mu}{\nu}\dgal{-,-}_{2,1} +\frac{\lambda}{\nu} \dgal{-,-}_{2,0}$ hence it is a linear combination of the elements with $\ell=0$ and $\ell=1$.

Next, we compute $\wdd(\dgal{-,-}_{a,b})$ for any $a>b \geq 0$, using that $\dgal{-,-}_{a,b}$ is defined by \eqref{Eq:dbr-ab}. This triple bracket is completely determined by its value on $(x,x,x)$. 
We have from \eqref{Eq:dP-gen} that 
\begin{align*}
    \wdd(\dgal{-,-}_{a,b})(x,x,x)=&
(1+\sigma+\sigma^2)\, \dgal{x,\dsq{x,x}}_{a,b;L} \\
&+(1+\sigma+\sigma^2)\, \dsq{x,\dgal{x,x}_{a,b}}_{L}\,,
\end{align*}
with $\sigma=\sigma_{(123)}$ acting on $\kk[x]^{\otimes 3}$. 
Thus, we are interested in the $\Z_3$-orbit of 
\begin{align*}
&  \dgal{x,\dsq{x,x}}_{a,b;L}+ \dsq{x,\dgal{x,x}_{a,b}}_{L} \\
=&\nu (\dgal{x,x^2}_{a,b}\otimes x - \dgal{x,x}_{a,b} \otimes x^2) 
+\mu \dgal{x,x^2}_{a,b} \otimes 1 \\
&+\lambda \dgal{x,x}_{a,b} \otimes 1 
+ \dsq{x,x^a} \otimes x^b - \dsq{x,x^b} \otimes x^a\,.
\end{align*}
This last expression can be explicitly expanded to yield 
\begin{align*}
&\quad \nu (x^{a+1}\otimes x \otimes x^b -x^b \otimes x^{a+1} \otimes x)
+\nu (x^{a+1} \otimes x^b \otimes x -x \otimes  x^{a+1} \otimes x^b )\\
&+\nu (x^{a}\otimes x^{b+1} \otimes x -x^{b+1} \otimes x \otimes x^a)
+\nu (x \otimes x^{b+1}\otimes x^{a}  -x^{b+1} \otimes x^a \otimes x)\\
&+\nu (x^b \otimes x^{a} \otimes x^2 - x^a \otimes x^b \otimes x^2) \\
&+\mu (x^{a+1} \otimes 1 \otimes x^b - x^b \otimes x^{a+1}\otimes 1)
+\mu (x^{a+1}  \otimes x^b \otimes 1 - 1 \otimes x^{a+1} \otimes x^b) \\
&+\mu (x^{a}  \otimes x^{b+1} \otimes 1 -  x^{b+1}\otimes 1 \otimes x^a)
+\mu (1  \otimes x^{b+1} \otimes x^{a} -  x^{b+1}\otimes x^a \otimes 1) \\
&+\mu (x^a \otimes x \otimes x^b + x \otimes x^b \otimes x^a)
-\mu (x \otimes x^a \otimes x^b + x^b \otimes x \otimes x^a) \\
&+ \lambda (x^a \otimes x^b \otimes 1 - 1 \otimes x^a \otimes x^b)
+\lambda (x^a \otimes 1\otimes x^b - x^b \otimes x^a \otimes 1) \\
&+\lambda  ( 1\otimes x^b \otimes x^a - x^b \otimes 1 \otimes x^a ) \,.
\end{align*}
Taking the $\Z_3$-orbit, we find after canceling out terms that 
\begin{equation}
    \begin{aligned} \label{Eq:dP-dbrx}
\wdd(\dgal{-,-}_{a,b})(x,x,x)=&
\nu \, (1+\sigma+\sigma^2) (x^b \otimes x^{a} \otimes x^2 - x^a \otimes x^b \otimes x^2) \\
&+2\mu (1+\sigma+\sigma^2)  (x^b \otimes x^a\otimes x - x^a \otimes x^b \otimes x) \\
&+\lambda (1+\sigma+\sigma^2)  (x^b \otimes x^a\otimes 1 - x^a \otimes x^b \otimes 1) 
    \end{aligned}
\end{equation}
Since a double bracket on $\kk[x]$ is of the form $\dgal{-,-}=\sum_{0\leq b < a} c_{a,b} \dgal{-,-}_{a,b}$ for some $c_{a,b}\in \kk$, its image under $\wdd$ evaluated on $(x,x,x)$ is simply obtained by linearity from \eqref{Eq:dP-dbrx}. 

In the case $\nu=\mu=0$, $\wdd(\dgal{-,-})(x,x,x)$ contains a nonzero term $-\lambda c_{a,b} x^a \otimes x^b \otimes 1$ if $c_{a,b}\neq 0$ with $b>0$. 
Therefore, $\ker \wdd \subseteq \oplus_{a>0} \kk \dgal{-,-}_{a,0}$. 
This is easily seen to be an equality as $\wdd(\dgal{-,-}_{a,0})(x,x,x)=0$ using \eqref{Eq:dP-dbrx} with $a>0$ in this case $\nu=\mu=0$. 
From our remark at the beginning of the proof, the kernel of $\wdd:\wBRA(\kk[x])_2 \to \wBRA(\kk[x])_3$ coincides with $\wdd(\Der(\kk[x]))$ and the second cohomology group is zero. 

In the case $\nu\neq 0$, we see from \eqref{Eq:dP-dbrx} that $\wdd(\dgal{-,-})(x,x,x)$ contains a nonzero term $-\nu c_{a,b} x^a \otimes x^b \otimes x^2$ of highest degree in the first component if $c_{a,b}\neq 0$ with $a,b>2$. 
This observation ensures that a double bracket $\dgal{-,-}$ such that $\wdd(\dgal{-,-})=0$ must be of the form 
$$\dgal{-,-}=\sum_{a> 2} c_{a,2} \dgal{-,-}_{a,2} +\sum_{a> 1} c_{a,1} \dgal{-,-}_{a,1} 
+ \sum_{a> 0} c_{a,0} \dgal{-,-}_{a,0} \,.$$
Using \eqref{Eq:dP-dbrx} again, we get for such a double bracket that 
\begin{align*}
    &\wdd(\dgal{-,-})(x,x,x) \\
=&\sum_{a>2} (\nu c_{a,1}-2\mu c_{a,2}) \, (1+\sigma+\sigma^2)\, (x^a \otimes x^2 \otimes x - x^a \otimes x \otimes x^2) \\
&+\sum_{a>2} (\nu c_{a,0}-\lambda c_{a,2}) \, (1+\sigma+\sigma^2)\, (x^a \otimes x^2 \otimes 1 - x^a \otimes 1 \otimes x^2) \\
&+\sum_{a>2} (2\mu c_{a,0}-\lambda c_{a,1}) \, (1+\sigma+\sigma^2)\, (x^a \otimes x \otimes 1 - x^a \otimes 1 \otimes x) \\
&+(\nu c_{1,0} - 2 \mu c_{2,0} + \lambda c_{2,1})\, (1+\sigma+\sigma^2) \,
(x^2 \otimes 1 \otimes x - x^2 \otimes x \otimes 1)\,.
\end{align*}
The vanishing of these terms imposes that an element  $\dgal{-,-}$ in $\ker \wdd:\wBRA(\kk[x])_2 \to \wBRA(\kk[x])_3$ must take the form 
\begin{equation}
\begin{aligned} \label{Eq:dbr-cab}
    \dgal{-,-}=&\sum_{a> 2} c_{a,2} \left( \dgal{-,-}_{a,2} +\frac{2\mu}{\nu} \dgal{-,-}_{a,1} 
    +\frac{\lambda}{\nu} \dgal{-,-}_{a,0} \right) \\
&+ c_{2,1} \left( \dgal{-,-}_{2,1}     -\frac{\lambda}{\nu} \dgal{-,-}_{1,0} \right) \\
&+ c_{2,0} \left( \dgal{-,-}_{2,0}     +\frac{2\mu}{\nu} \dgal{-,-}_{1,0} \right)\,,
\end{aligned}
\end{equation}
for arbitrary constants $c_{a,b}$. 
This expression clearly belongs to $\wdd(\Der(\kk[x]))$ if we look at its basis given at the beginning of the proof. Thus, $\widehat{\dPH}^2(\kk[x])$ is zero. 

Finally, the result for $\dPH^2(\kk[x])$ is direct if we apply the inverse of $\mu_2:(T^\ast \kk[x])_{\sharp,2}\to \wBRA(\kk[x])_2$ in cohomology. 
\end{proof}

\begin{corollary} \label{cor:dPh-kx}
The double Poisson bracket \eqref{Eq:dbr-xlin2} is exact in the completed double Poisson cohomology that it defines.
\end{corollary}

Alternatively, \eqref{Eq:dPtheta} provides an explicit formula proving Corollary \ref{cor:dPh-kx}: 
\begin{align*}
    \dsq{-,-}=&\wdd(-\theta_1) \,, && \text{if }\mu=\nu=0\,, \\
    \dsq{-,-}=&\wdd(\theta_1 + \frac{\mu}{\nu} \theta_0)\,, && \text{if }\nu\neq 0\, .
\end{align*}

\subsection{Double Poisson brackets on  \texorpdfstring{$\kk[x]/(x^r)$}{k[x]/(xr)}} 
\label{ss:dP-coh-xr}

Fix $r\geq 2$. 
Let us note that each double bracket $\dgal{-,-}_{a,b}$ defined by \eqref{Eq:dbr-ab}, $a>b\geq 0$, satisfies 
\begin{equation}
    \dgal{x^r,x}_{a,b},\dgal{x,x^r}_{a,b}\in (x^r) \otimes \kk[x] + \kk[x] \otimes (x^r)\, . 
\end{equation} 
Hence, $\dgal{-,-}_{a,b}$ descends to a well-defined double bracket on $\kk[x]/(x^r)$. 
Furthermore, the double bracket $\dsq{-,-}$ given by \eqref{Eq:dbr-xlin2} stays Poisson if $\lambda\nu-\mu^2=0$. However, it is not hard to check that this double Poisson bracket is \emph{no longer} an element of $\mu_1(T^\ast (\kk[x]/(x^r)))_{\sharp,2})$ as the bivector \eqref{Eq:dbr-kxP} is \emph{not} defined on $\kk[x]/(x^r)$; indeed, the double derivation $\partial_x$ \eqref{DDer-delx} does not descend to the quotient $\kk[x]/(x^r)$.  
In particular, the completed double Poisson cohomology $\widehat{\dPH}(\kk[x]/(x^r))$ (see Definition \ref{Def:dPH-comp}) of this double Poisson bracket exists, but we can \emph{not} define its double Poisson cohomology (see Definition \ref{Def:dPH}) following Pichereau-Van de Weyer \cite{PV,VdW}. 

We compute the first 3 groups of the completed double Poisson cohomology on $\kk[x]/(x^r)$ in this subsection. 
By comparing these results with those from Subsection~\ref{ss:dP-coh-x}, we can observe that passing from $\kk[x]$ to $\kk[x]/ (x^r)$ may shrink a cohomology group (for $\widehat{\dPH}^1$ with $\nu,\mu\neq 0$) or it may enlarge a cohomology group (for $\widehat{\dPH}^2$ when $r\geq 3$ with $\nu\neq 0$ and $\mu=\lambda=0$).

\begin{proposition}  \label{Pr:dPH-kxTrun} Fix $r \geq 3$. 
For the double Poisson bracket defined on $\kk[x]/(x^r)$ by \eqref{Eq:dbr-xlin2}, we have 
\begin{enumerate}[(1)]
\item $\widehat{\dPH}^0(\kk[x]/(x^r))=\kk[x]/(x^r)$.
\item If $\nu\neq 0$ and $\mu=\lambda=0$, then $\widehat{\dPH}^1(\kk[x]/(x^r))=\kk\,\theta_2$. 
Otherwise, $\widehat{\dPH}^1(\kk[x]/(x^r))=\{0\}$.
\item If $\nu\neq 0$ and $\mu=\lambda=0$, then $\widehat{\dPH}^2(\kk[x]/(x^r))=\kk  \dgal{-,-}_{2,0}$. Otherwise,
$\widehat{\dPH}^2(\kk[x]/(x^r))=\{0\}$.
\end{enumerate}
\end{proposition}
\begin{proof}
For (1), we can repeat the proof of Proposition \ref{Pr:dPH-kx0}. 

For (2), note that only the derivation $\theta_0$ does not descend from $\kk[x]$ to $\kk[x]/(x^r)$. Thus any derivation is of the form $\theta=\sum_{\ell=1}^{r-1} c_\ell \theta_\ell$. Repeating the proof of Proposition \ref{Pr:dPH-kx0}, we can derive \eqref{Eq:dPtheta} and then see that if $c_\ell\neq 0$ for some $\ell>2$ the double bracket $\wdd(\theta)$ is nonzero. Otherwise we have
\begin{align*}
\wdd(\theta)(x,x)=& 
(2 \mu c_2-\nu c_1) (x \otimes x^2 - x^2 \otimes x) 
+ \lambda c_2 (1 \otimes x^2 - x^2 \otimes 1) \\
&-2\lambda c_1(x \otimes 1 - 1\otimes x) \,,
\end{align*}
for $\theta=c_1 \theta_1 + c_2 \theta_2$. 
Requiring the vanishing of the last equation and remembering the condition $\mu^2=\lambda\nu$, we get the desired statement.

For (3), recall that each double bracket $\dgal{-,-}_{a,b}$ defined by \eqref{Eq:dbr-ab} descends to $\kk[x]/(x^r)$. 
Meanwhile, the non-existence of $\theta_0$ on $\kk[x]/(x^r)$ gives that the basis of 
$\wdd(\Der(\kk[x]/(x^r)))$ is the same as the one given in the proof of Proposition \ref{Pr:dPH-kx2}  if $\nu=0$. Thus, in that case, we can repeat the proof of Proposition \ref{Pr:dPH-kx2} verbatim. 
If $\nu\neq 0$ instead, $\wdd(\Der(\kk[x]/(x^r)))$  
is the $\kk$-linear span of: 
\begin{enumerate}   
    \item[($\nu\neq0$)] \quad  $\dgal{-,-}_{\ell,2}+2\frac{\mu}{\nu}\dgal{-,-}_{\ell,1}+\frac{\lambda}{\nu}\dgal{-,-}_{\ell,0}$ for $\ell\geq 1$.
\end{enumerate} 
(I.e. we can not consider $\ell=0$.) The cases $\ell=1,2$ give the elements 
$-\dgal{-,-}_{2,1} +\frac{\lambda}{\nu} \dgal{-,-}_{1,0}$ and  $\frac{2\mu}{\nu}\dgal{-,-}_{2,1} +\frac{\lambda}{\nu} \dgal{-,-}_{2,0}$. Note that the $\ell=2$ term identically vanishes if $\mu$ (thus $\lambda$) is zero.  
The proof of Proposition \ref{Pr:dPH-kx2} can then be used verbatim to derive that $\wdd(\dgal{-,-})=0$ provided $\dgal{-,-}$ is of the form \eqref{Eq:dbr-cab}.  
In that expression, the term $\dgal{-,-}_{2,0}+\frac{2\mu}{\nu} \dgal{-,-}_{1,0}$ is in $\wdd(\Der(\kk[x]/(x^r)))$ only if $\lambda$ is invertible, which yields the statement.
\end{proof}

In the case $r=2$ with $\kk[x]/(x^2)$, we only need to consider the double Poisson bracket $\dsq{-,-}=\dgal{-,-}_{1,0}$. 
\begin{proposition}  \label{Pr:dPH-kxTr2} 
For the double Poisson bracket $\dgal{-,-}_{1,0}$ on $\kk[x]/(x^2)$,  
\begin{enumerate}[(1)]
\item $\widehat{\dPH}^0(\kk[x]/(x^2))=\kk[x]/(x^2)$.
\item $\widehat{\dPH}^1(\kk[x]/(x^2))=\{0\}$.
\item $\widehat{\dPH}^2(\kk[x]/(x^2))=\{0\}$.  
\end{enumerate}
\end{proposition}
\begin{proof}
 This is easily adapted from the case $\mu=\nu=0$ on $\kk[x]$. 
\end{proof}

\section{Gauged double Poisson cohomology for one generator} 

We present computations of the (completed) gauged double Poisson cohomology introduced in Definition \ref{Def:gdPH}. For $\kk$-algebras in one generator $x$, the gauge element is $\Delta:x \mapsto x\otimes 1 - 1\otimes x$. 

\subsection{The case of $\kk[x]$}
\label{ss:gdP-coh-x}

Given any double Poisson bracket, the maps $(\mu_\ell)_{\ell\geq 0}$ yield an isomorphism $\widehat{\gdPH}(\kk[x])\simeq \gdPH(\kk[x])$. 
Due to Corollary \ref{Cor:dPH-gdPH}, the cohomology classes in $\dPH(\kk[x])$ are sent to classes in $\gdPH(\kk[x])$. The aim is to determine if some of these classes coincide, and if new classes appear, i.e. if there exists $Q\in (\mb T^\ast \kk[x])_\sharp$ with $\dd_P(Q)=(\alpha \Delta)_\sharp \neq 0$ for some $\alpha \in \mb T^\ast \kk[x]$. 
As we shall see, computations can be trivialized  in $({\mc D}_{\kk[x]},\dd_P)$ than in $((\mb T^\ast \kk[x])_\sharp,\dd_P)$, and we can compute \emph{all} cohomology groups.

We denote by $\alpha_\Delta \in {\mc D}_{\kk[x]}^k$ the projection of an element $\alpha_\sharp \in (\mb T^\ast \kk[x])_{k,\sharp}$ under $\pi_k^\Delta: (\mb T^\ast \kk[x])_{k,\sharp}\to {\mc D}_{\kk[x]}^k$, cf. Section~\ref{ss:gdPCoh}. 
Recall the  double derivation $\partial_x$ \eqref{DDer-delx}. 
\begin{lemma} \label{Lem:Dkx-span}
For $k\geq 1$, ${\mc D}_{\kk[x]}^k=\operatorname{span}_\kk \{(x^j \del_x^k)_\Delta \mid j\geq 0\}$.   In particular, 
\begin{enumerate}[(1)]
    \item ${\mc D}_{\kk[x]}^k=\oplus_{j\geq 0} \kk (x^j \del_x^k)_\Delta$ if $k\geq 1$ is odd;
    \item ${\mc D}_{\kk[x]}^k=0$ if $k\geq 2$ is even.
\end{enumerate}
\end{lemma}
\begin{proof}
By definition, the following identity holds in ${\mc D}_{\kk[x]}^k$ for any $k\geq 1$: 
\begin{equation} \label{Eq:gdP-pf-3}
    (\alpha x \partial_x)_\Delta = (x \alpha \partial_x)_\Delta \,, \quad \forall \alpha \in (\mb T^\ast \kk[x])_{k-1}\,,
\end{equation}
since $\Delta=\del_x x - x \del_x$. 
Hence an arbitrary element $(x^{a_1}\del_x \ldots x^{a_k}\del_x)_\Delta \in {\mc D}_{\kk[x]}^k$ is in the class $(x^{a_1+\ldots+a_k}\del_x^k)_\Delta$. 

Assume that $k$ is odd. Using the map \eqref{Eq:AKKN-map}, we get for $j\geq 0$
\begin{equation}
    \mu_{k,\sharp}((x^j \del_x^k)_\sharp)(x,\ldots,x)
=\sum_{0\leq i\leq k-1} \sigma^i\circ (x^j\otimes 1 \otimes \ldots \otimes 1)\,.
\end{equation}
In particular, no nontrivial linear combination of the elements $((x^j \del_x^k)_\sharp)_{j\geq 0}$ is sent to the zero $k$-bracket under $\mu_{k,\sharp}$. 
Similarly to the proof of Lemma \ref{Lem:AKKN}, we deduce that the classes $(x^j \del_x^k)_\Delta \in {\mc D}_{\kk[x]}^k$ are nonzero and linearly independent. 

Assume that $k$ is even; we claim that $(x^b \del_x^{k})_\Delta =0$ for any $b\geq 0$. This is direct if $b=0$ as $(\del_x^k)_\sharp = (-1)^{k-1}(\del_x^k)_\sharp$ in $(\mb T^\ast \kk[x])_{k,\sharp}$.  
For $b\geq 1$, we compute 
\begin{equation*}
   (x^b \del_x^{k})_\sharp =(x^{b-1} \del_x^{k-1}\Delta)_\sharp + (x^{b-1} \del_x^{k-1}x \del_x)_\sharp 
   =(x^{b-1} \del_x^{k-1}\Delta)_\sharp - (x \del_xx^{b-1} \del_x^{k-1})_\sharp \,.
\end{equation*}
In ${\mc D}^{k}_{\kk[x]}$, where $x$ and $\del_x$ commute, this entails  
$(x^b \del_x^{k})_\Delta=-(x^b \del_x^{k})_\Delta$, and the claim follows. 
\end{proof}

\begin{proposition}  \label{Pr:gdPH-kx0}   
For the (nonzero) double Poisson bracket defined on $\kk[x]$ by \eqref{Eq:dbr-xlin2},  
each map $\dd_P:{\mc D}_{\kk[x]}^{k}\to {\mc D}_{\kk[x]}^{k+1}$ is the zero map. 

\noindent Therefore $\gdPH^k(\kk[x])={\mc D}_{\kk[x]}^k$ for all $k\geq 0$, and 
\begin{enumerate}[(1)]
\item $\gdPH^0(\kk[x])=(\mb T^\ast \kk[x])_{0,\sharp}=\kk[x]$; 
\item $\gdPH^1(\kk[x])=(\mb T^\ast \kk[x])_{1,\sharp}=\oplus_{j \geq 0}\, \kk (x^j \partial_x)_\sharp$; 
\item $\gdPH^{2\ell}(\kk[x])=0$ and $\gdPH^{2\ell+1}(\kk[x])=\oplus_{j\geq 0} \kk (x^j \del_x^{2\ell+1})_\Delta$  for $\ell\geq 1$. 
\end{enumerate}
\end{proposition}
\begin{proof}
From Proposition \ref{Pr:dPH-kx0}, we have $\dPH^0(\kk[x])=\kk[x]=(\mb T^\ast \kk[x])_{0,\sharp}$. Since ${\mc D}_{\mc A}^0=(\mb T^\ast \kk[x])_{0,\sharp}$, we deduce part (1) from the diagram depicted in Section~\ref{ss:gdPCoh}.    

We get that $\dd_P:{\mc D}_{\kk[x]}^{k}\to {\mc D}_{\kk[x]}^{k+1}$ is zero from part (1) if $k=0$ and from Lemma \ref{Lem:Dkx-span} if $k\geq 1$. Therefore $\gdPH^k(\kk[x])={\mc D}_{\kk[x]}^k$ by Definition \ref{Def:gdPH}. 

We deduce part (3) from Lemma \ref{Lem:Dkx-span}.  

For part (2), note that $\im \iota_1^\Delta: (\mb T^\ast \kk[x])_{0}\to (\mb T^\ast \kk[x])_{1,\sharp}$ reduces to $\{0\}$, because $(f\Delta)_\sharp = ([x,f] \del_x)_\sharp=0$ for any $f\in \kk[x]$. Hence ${\mc D}_{\kk[x]}^1=(\mb T^\ast \kk[x])_{1,\sharp}$.   
\end{proof}

\begin{remark} \label{Rem:dPH-gdPH-surj}
The natural morphism $\pi^\Delta:\dPH(\kk[x]) \to \gdPH(\kk[x])$ is clearly not surjective by comparing the first cohomology groups in Propositions \ref{Pr:dPH-kx1} and \ref{Pr:gdPH-kx0}.
\end{remark}

\subsection{Comment on the case of $\kk[x]/(x^r)$}
\label{ss:gdP-coh-xr}

Each double Poisson bracket on $\cA:=\kk[x]/(x^r)$ of the form \eqref{Eq:dbr-xlin2} can \emph{not} be defined by a bivector $P\in (\mb T^\ast \cA)_2$, as we mentioned in Subsection~\ref{ss:dP-coh-xr}. 
Therefore we can not consider their (completed) gauged double Poisson cohomology since this last theory is defined using the corresponding bivector, cf. Definition \ref{Def:gdPH}.

\subsection{The quasi-Poisson setting for $\kk[x]$}
We start by recalling the following result, cf. \cite[Prop.~4.1]{F21}. 

\begin{lemma}
Any double quasi-Poisson bracket on $\kk[x]$ is of the form \eqref{Eq:dbr-xlin2} for arbitrary $\lambda,\mu,\nu \in \kk$. 
\end{lemma}
\begin{proof}
An explicit computation yields for the double bracket    \eqref{Eq:dbr-xlin2} 
\begin{equation} \label{Eq:qP-xxx}
   \dgal{x,x,x}=(\mu^2-\lambda \nu) \ (\Id+\sigma +\sigma^2) \circ 
   (t^2\otimes t \otimes 1 - 1\otimes t \otimes t^2)\,,
\end{equation} 
which yields \eqref{qPabc} (with $|S|=1$) for $q=4(\mu^2-\lambda \nu)$ and $a=b=c=x$. 
The result follows for arbitrary $a,b,c\in \kk[x]$ since $\dgal{-,-,-}$ is a triple bracket. 
\end{proof}

\begin{remark}
If one uses the original definition of Van den Bergh \cite{VdB1}, the quasi-Poisson condition becomes $\mu^2-\lambda \nu=\frac14$. 
\end{remark}

\begin{proposition}  \label{Pr:dqPH} 
For the double quasi-Poisson bracket \eqref{Eq:dbr-xlin2}, we have 
\begin{enumerate}[(1)]
\item $\widehat{\dPH}^0(\kk[x])=\dPH^0(\kk[x])=\kk[x]$;
\item  $\widehat{\dPH}^1(\kk[x])= \kk (\nu \theta_2+2\mu \theta_1+\lambda \theta_0)$, $\dPH^1(\kk[x])= \kk (\nu x^2\partial_x +2\mu \partial_x+\lambda \partial_x)$;
    \item $\widehat{\dPH}^2(\kk[x])= \{0\}$ and $\dPH^2(\kk[x])= \{0\}$.
\end{enumerate} 
\end{proposition} 
\begin{proof}
This basically follows from the proofs of Propositions \ref{Pr:dPH-kx0}, \ref{Pr:dPH-kx1} and \ref{Pr:dPH-kx2}. Indeed, the case $\mu=\nu=0$ is treated therein, and the case $\nu\neq 0$ does not use the Poisson condition $\mu^2-\lambda \nu=0$ so that it holds for arbitrary $\lambda,\mu \in \kk$. 
Thus, one only needs to adapt the end of the proofs to handle the case $\nu=0,\mu\neq 0$, $\lambda\in \kk$; this is left to the reader. 
\end{proof}

\begin{remark}
Recall from Corollary \ref{cor:dPh-kx} that each double Poisson bracket is exact in its completed double Poisson cohomology. 
This is \emph{not} true for a generic double quasi-Poisson bracket with $\mu^2-\lambda \nu\neq 0$.  Indeed, the condition $\wdd(\dgal{-,-})=0$ is equivalent to $\dd_P(P)=0$ (where $P$ is given by \eqref{Eq:dbr-kxP}) while the quasi-Poisson property entails $\dd_P(P)=\frac23(\mu^2-\lambda \nu)(\Delta^3)_\sharp$.
\end{remark}

\subsection{More on the gauged double Poisson cohomology on $\kk[x]$} 

By Definition \ref{Def:gaugDBR}, $P\in (\mb T^\ast \kk[x])_{\sharp,2}$ is gauged Poisson if 
$\brSN{P,P}=(\Delta R)_\sharp$ for some $R\in (\mb T^\ast \kk[x])_2$. 
Equivalently, this means that 
$\brSN{P,P}\in \im \iota_3^\Delta$, cf. Section \ref{ss:gdPCoh}. 
\begin{lemma}
For any    $P\in (\mb T^\ast \kk[x])_{\sharp,2}$,  $\brSN{P,P}\in \im \iota_3^\Delta$. 
In particular, $P$ is always a gauged Poisson element.
\end{lemma}
\begin{proof}
Any $P$ is a linear combination of the elements $(x^a \partial_x x^b \partial_x)_\sharp$, so it suffices to prove the claim for  $P=(x^a \partial_x x^b \partial_x)_\sharp$. 

By Theorem \ref{Thm:dSN} and \eqref{Eq:gdb-skew}, $\dSN{x,x}=0$, $\dSN{\partial_x,\partial_x}=0$ and $\dSN{\partial_x,x}=1\otimes 1$, $\dSN{x,\partial_x}=-1\otimes 1$. 
Using the graded Leibniz rule \eqref{Eq:gdb-Lleib}, we get 
\begin{align*}
    \dSN{P,x^c}&=\sum_{\gamma=0}^{c-1} (x^\gamma \otimes x^a \partial_x x^{b+c-\gamma-1}
    - x^{b+\gamma} \partial_x \otimes x^{a+c-\gamma-1}) \,,\\
    \dSN{P,\partial_x}&=\sum_{\beta=0}^{b-1} x^\beta \partial_x \otimes x^a \partial_x x^{b-\beta-1}
 - \sum_{\alpha=0}^{a-1}  x^{\alpha}\partial_x x^b \partial_x \otimes x^{a-\alpha-1}) \,. \\
\end{align*}
Similarly, we can obtain  
\begin{align*}
 \brSN{P,P}=&
-x^a\partial_x x^b \brSN{P,\partial_x}
-x^a \partial_x \brSN{P,x^b} \partial_x \\
&+x^a \brSN{P,\partial_x} x^b \partial_x
+ \brSN{P,x^a} \partial_x x^b \partial_x. 
\end{align*}
Plugging the first two equalities in the last one, we find  modulo graded commutators 
%(where $\brSN{-,-}$ now stands for the image of $\mult \circ \dSN{-,-}$ in $(T^\ast \kk[x])_\sharp$) 
\begin{align*}
 \brSN{P,P}=&
2 \sum_{\alpha=0}^{a-1}  (\partial_x x^b \partial_x x^{2a-\alpha-1} \partial_x x^{b+\alpha})_\sharp
-2 \sum_{\alpha=0}^{a-1}  (\partial_x x^b \partial_x x^{a+b-\alpha-1} \partial_x x^{a+\alpha})_\sharp \\
&+2 \sum_{\beta=0}^{b-1}  (\partial_x x^a \partial_x x^{2b-\beta-1} \partial_x x^{a+\beta})_\sharp
-2 \sum_{\beta=0}^{b-1}  (\partial_x x^a \partial_x x^{a+b-\beta-1} \partial_x x^{b+\beta})_\sharp\,.
\end{align*}
If we write $Q_\Delta$ for the class of $Q_\sharp$ modulo $\im \iota^\Delta$, we get from \eqref{Eq:gdP-pf-3} that 
$$(\partial_x x^{c_1}\partial_x x^{c_2}\partial_x x^{c_3})_\Delta=(\partial_x^3 x^{c_1+c_2+c_3})_\Delta\,, \qquad c_1,c_2,c_3 \geq 0.$$
Therefore, the previous expression for $\brSN{P,P}$ vanishes modulo $\im \iota_3^\Delta$.
\end{proof}

We obtain the following result by repeating the proof of Proposition \ref{Pr:gdPH-kx0}. 
\begin{proposition}  \label{Pr:gdPH-kxGen}   
Fix  $P\in (\mb T^\ast \kk[x])_{\sharp,2}$.   
Each map $\dd_P:{\mc D}_{\kk[x]}^{k}\to {\mc D}_{\kk[x]}^{k+1}$ is the zero map, and therefore $\gdPH^k(\kk[x])={\mc D}_{\kk[x]}^k$ for all $k\geq 0$. In particular, 
\begin{enumerate}[(1)]
\item $\gdPH^0(\kk[x])=(\mb T^\ast \kk[x])_{0,\sharp}=\kk[x]$; 
\item $\gdPH^1(\kk[x])=(\mb T^\ast \kk[x])_{1,\sharp}=\oplus_{j \geq 0}\, \kk (x^j \partial_x)_\sharp$; 
\item $\gdPH^{2\ell}(\kk[x])=0$ and $\gdPH^{2\ell+1}(\kk[x])=\oplus_{j\geq 0} \kk (x^j \del_x^{2\ell+1})_\Delta$  for $\ell\geq 1$. 
\end{enumerate}
\end{proposition}

\section[\MakeLowercase{d}PA with two generators and its truncation]{A double Poisson algebra with two generators and its truncation}

Consider the double bracket $\dsq{-,-}$ on $\kk\langle u,v \rangle$ satisfying 
\begin{equation}  \label{Eq:dP-uv}
    \dsq{u,u}=0,\quad \dsq{v,v}=0,\quad \dsq{u,v}=vu \otimes uv\,.
\end{equation}
(In particular, $\dsq{v,u}= - uv\otimes vu$ by cyclic skewsymmetry.) 
It is shown in \cite[Rem.~4.15]{FV} that this double bracket is Poisson. 
Introduce 
\begin{equation}
\begin{aligned}  \label{DDer-deluv}
 \partial_u\in \DDer(\kk\langle u,v\rangle), \qquad 
 \partial_u(u)=1\otimes 1, \,\, \partial_u(v)=0,   \\
 \partial_v\in \DDer(\kk\langle u,v\rangle), \qquad 
 \partial_v(u)=0, \,\, \partial_v(v)=1\otimes 1 .
\end{aligned}
\end{equation}
Then, the bivector $P$ defining $\dsq{-,-}$ \eqref{Eq:dP-uv} through Proposition \ref{Pr:MapMu} can be easily seen to be given by 
\begin{equation} \label{Eq:dPbiv-uv}
    P=u \partial_u u v \partial_v v\,. 
\end{equation} 
Our aim is to compute the first two groups of its (completed) double Poisson cohomology, which yield the first such result for a quartic double Poisson bracket. Indeed, all the cases considered over $\kk\langle u,v \rangle$ by Pichereau and Van de Weyer \cite{PV} are either linear or quadratic, cf. Section \ref{ss:OthExmp}. 

Introduce the truncated algebra $\cA^{\tr}:=\kk\langle u,v \rangle/(u^2,v^2)$ as the quotient algebra of $\kk\langle u,v\rangle$ by the relations $u^2=v^2=0$. 
The double bracket \eqref{Eq:dP-uv} descends to  $\cA^{\tr}$, where it remains a double Poisson bracket. In fact, it is still defined by the element 
$P\in (\mb T^\ast \cA^{\tr})_2$ given by \eqref{Eq:dPbiv-uv} because $u \partial_u u$ and $v \partial_v v$ descend to double derivations on $\cA^{\tr}$ (though $\partial_u$ and $\partial_v$ are \emph{not} defined over $A^{\tr}$). We can thus speak both of the double Poisson cohomology according to Definition \ref{Def:dPH} and the completed version according to Definition \ref{Def:dPH-comp} for this example. Since the maps $(\mu_n)_{n\geq 1}$ may fail to be isomorphisms (e.g. the derivation of $\cA^{\tr}$ satisfying $u\mapsto u$, $v\mapsto 0$ is not in the image of $\mu_1$) it is thus important to compare these two cohomology theories for the double bracket \eqref{Eq:dP-uv}, which is the second main aim of this section. 

\subsection{The case of $\kk\langle u,v \rangle$}

We set $\cA:=\kk\langle u,v \rangle$ from now on. 
Our aim is to compute $\widehat{\dPH}^0$ and $\widehat{\dPH}^1$ for the double Poisson bracket given by \eqref{Eq:dP-uv}. Since $\cA$ is a free algebra, we will also obtain $\dPH^0$ and $\dPH^1$ by Corollary \ref{Cor:Iso-dPcoh}. 
For $f,g\in \cA$, we consider the derivation $\theta_{f,g}\in \Der(\cA)$  which is uniquely determined by  
\begin{equation} \label{Eq:theta-fg}
\theta_{f,g}(u)=f\,, \quad \theta_{f,g}(v)=g\,.
\end{equation}
Any $D \in \Der(\cA)$ is of this form since  $D=\theta_{d(u),d(v)}$. 
It will be convenient to see $\cA$ as a non-negatively graded algebra in 3 ways: using the grading $\deg_u$ where $\deg_u(u)=1$ and $\deg_u(v)=0$; using the grading $\deg_v$ where $\deg_v(u)=0$ and $\deg_v(v)=1$; or using the total grading $\deg=\deg_u+\deg_v$. 

\subsubsection{Preliminary results}

Following \cite{PV}, introduce the Euler derivation $\mathrm{E}\in \Der(\cA)$ given by 
\begin{equation} \label{Eq:Euler}
a\mapsto \mathrm{E}(a):= \mult \circ \left(u\star \partial_u(a) + v \ast \partial_v(a) \right)\,.  
\end{equation}
The noncommutative Euler formula from \cite[Prop.~12]{PV} states that, if $a\in \cA$ is homogeneous for $\deg$, then 
\begin{equation} \label{NCEulerFor}
   \mathrm{E}(a)=\deg(a)\, a\,. 
\end{equation} 
The following is easily obtained, cf. \cite[Lem.~2.6]{DSKV}.
\begin{lemma} \label{Lem:dLdR}
If $\delta,\Delta \in \{\partial_u,\partial_v\}$, 
$\delta_L \circ \Delta = \Delta_R \circ \delta$
for the extensions \eqref{20240805:eq1b}. 
\end{lemma}

\begin{lemma} \label{Lem:dvadva}
Let $f\in \cA$. The following two conditions are equivalent: 
\begin{itemize}
    \item[(i)] $\partial_v(f)=(\partial_v (f))^\sigma$; 
    \item[(ii)] $f=\partial_v(a)'' \partial_v(a)'$ for some $a\in \cA$. 
\end{itemize}    
The same result holds if we replace $v$ by $u$. 
\end{lemma}
\begin{proof}
The proof presented below can be repeated verbatim if one replaces $v$ by $u$, so we only prove the equivalence in the former case. 

Assuming that (ii) holds, we can write 
\begin{align*}
    \partial_v(f)=\partial_v(\partial_v(a)'')' \otimes \partial_v(\partial_v(a)'')''  \partial_v(a)' 
+\partial_v(a)'' \partial_v(\partial_v(a)')' \otimes \partial_v(\partial_v(a)')''\,.
\end{align*}
By Lemma \ref{Lem:dLdR}, we have  
\begin{equation*}
   \partial_v(\partial_v(a)')' \otimes \partial_v(\partial_v(a)')'' \otimes \partial_v(a)''
   = \partial_v(a)' \otimes \partial_v(\partial_v(a)'')' \otimes \partial_v(\partial_v(a)'')''\,,
\end{equation*}
which directly yields 
\begin{align*}
\partial_v(f)&=\partial_v(\partial_v(a)')'' \otimes \partial_v(a)''  \partial_v(\partial_v(a)')'   
+ \partial_v(\partial_v(a)'')''  \partial_v(a)' 
\otimes \partial_v(\partial_v(a)'')' \\
&=(\partial_v(f))^\sigma\,.
\end{align*}

When (i) holds, we can assume without loss of generalities that $f$ has degrees $\deg(f)=d$ and $\deg_v(f)=\ell$ for $0\leq \ell \leq d$ because the map $a\mapsto \partial_v(a)''\partial_v(a)'$ is of degree $-1$ with respect to $\deg$ and $\deg_v$. 
Furthermore, if $\ell=0$ or $\ell=d$, this is straightforward as 
\begin{equation*}
    f=u^d = \mult \circ (\partial_v(u^dv))^\sigma \,, \\ \quad 
    f=v^d = \mult \circ (\partial_v(v^{d+1}/(d+1)))^\sigma\,. 
\end{equation*}
Therefore, we can assume that $f$ is of degree $d\geq 2$ and each of its terms has exactly $1\leq \ell <d$ factors $v$, that is 
\begin{equation}
  f=\sum_{I \in \Sigma_{d,\ell} } 
  \alpha_I \, u^{i_0} v u^{i_1} v u^{i_2} \ldots v u^{i_\ell} \,, \qquad \alpha_I\in \kk \,,
\end{equation}
where we set $\Sigma_{d,\ell}=\{I=(i_0,\ldots,i_\ell)\mid i_j\geq 0;\, \sum_j i_j=d-\ell \}$. We directly compute 
\begin{align*}
   \partial_v(f)= \sum_{I \in \Sigma_{d,\ell} } \sum_{1\leq s \leq \ell} 
  \alpha_I \, \underbrace{u^{i_0}\ldots  v u^{i_{s-1}}}_{s-1 \text{ factors }v} 
\otimes \underbrace{u^{i_s}v \ldots u^{i_\ell} }_{\ell-s \text{ factors }v}\,,
\end{align*}
from which assumption (i) implies the following equalities with $s=1,\ldots,\ell$: 
\begin{align*}
\sum_{I \in \Sigma_{d,\ell} } 
  \alpha_I \, u^{i_0}\ldots  v u^{i_{s-1}} \otimes u^{i_s}v \ldots u^{i_\ell} 
=& 
\sum_{I \in \Sigma_{d,\ell} } 
  \alpha_I \, u^{i_{\ell+1-s}}\ldots  v u^{i_{\ell}} \otimes u^{i_0}v \ldots u^{i_{\ell-s}} \,.
\end{align*}
Equating multiples of the same term, we require that $\alpha_I=\alpha_{\sigma^s(I)}$, with 
$\sigma:\Sigma_{d,\ell}\to \Sigma_{d,\ell}$ given by  $\sigma(i_0,\ldots,i_\ell)=(i_{\sigma^s(0)},\ldots,i_{\sigma^s(\ell)})$ where on the indices $\sigma=\sigma_{(0,\ldots,\ell)}$ is the cyclic permutation of order $\ell+1$. 
It now suffices to check the claim when $f$ has the form 
\begin{equation*}
    f=\sum_{s=0}^\ell u^{i_{\sigma^s(0)}} v u^{i_{\sigma^s(1)}}  \ldots v u^{i_{\sigma^s(\ell)}} \,.
\end{equation*}
Taking $a=u^{i_0}v u^{i_1}v \ldots v u^{i_\ell} v$, we verify that 
\begin{align*}
\partial_v(a)'' \partial_v(a)'=&
\sum_{r=0}^\ell \left( u^{i_{r+1}} v \ldots  u^{i_{\ell}}v  \right) 
\left( u^{i_{0}} v \ldots  v u^{i_{r}} \right) = f\,. \qedhere 
\end{align*}
\end{proof}
\begin{remark}
The proof of Lemma \ref{Lem:dvadva} provides an explicit method to find $a\in \cA$ when we are given an element $f$ satisfying assumption (i). 
For example, $f=v^2u+uv^2+vuv+u^k v u^k$ satisfies (i) for any $k\geq 1$. We can decompose $f$ as the $\Z_3$-orbit of $u v^2$ (for the action of $\sigma$ on $(1,0,0)\in\Sigma_{3,2}$) and the $\Z_2$-orbit of $\frac12 u^k v u^k$ (for the action of $\sigma$ on $(k,k)\in\Sigma_{2k+1,1}$). Therefore we can take $a=uv^3+\frac12 u^kvu^k v$.   
\end{remark}

\begin{lemma} \label{Lem:dudv}
For fixed $a\in \cA$, set $a_u:=(\partial_u a)''(\partial_u a)'$ and $a_v:=(\partial_v a)''(\partial_v a)'$. 
Then $\partial_v(a_u)=(\partial_u(a_v))^\sigma$. 
\end{lemma}
\begin{proof}
    With the notation \eqref{20240805:eq1b}, we can directly write 
\begin{align*}
 \partial_v(a_u)&=(\Id_\cA \otimes \mult)  \big((\partial_v)_R(\partial_u a)\big)^{\sigma^2}
 +(\mult \otimes \Id_\cA)  \big((\partial_v)_L(\partial_u a)\big)^{\sigma}\,, \\
 \partial_u(a_v)^\sigma 
 &=(\mult\otimes \Id_\cA)  \big((\partial_u)_R(\partial_v a)\big)^{\sigma}
 + (\Id_\cA \otimes \mult) \big((\partial_u)_L(\partial_v a)\big)^{\sigma^2}\,.
\end{align*} 
We get from these two expressions 
\begin{align*}
\partial_v(a_u)  -  \partial_u(a_v)^\sigma 
=&(\Id_\cA \otimes \mult)  \big[(\partial_v)_R\circ \partial_u (a) - (\partial_u)_L\circ \partial_v (a)\big]^{\sigma^2} \\
&+(\mult \otimes \Id_\cA)  \big[(\partial_v)_L\otimes \partial_u (a) - (\partial_u)_R\circ \partial_v (a)\big]^{\sigma}\,,
\end{align*}
which vanishes by Lemma \ref{Lem:dLdR}. 
\end{proof}

\subsubsection{The cohomology groups}

Hereafter, $\wdd$ is the differential of Definition \ref{def:wdd} associated with $\dsq{-,-}$ given by \eqref{Eq:dP-uv}. 

\begin{lemma}  \label{Lem:kxy-1}
The image of the $\kk$-linear map $\wdd:\cA_\sharp\to \Der(\cA)$ is spanned by 
the derivations $\theta_{f,g}\in \Der(\cA)$ with 
\begin{equation} 
f=uv \partial_v(a)''\partial_v(a)' vu,\quad 
g=-vu \partial_u(a)''\partial_u(a)' uv, \quad  \text{ for some }a\in \cA\,. 
\end{equation}
Moreover, $\ker (\wdd:\cA_\sharp\to \Der(\cA))=\kk$. 
\end{lemma}
\begin{proof}
 Let $a=a(u,v)\in \cA$ be a lift of an element $\bar{a}\in \cA_\sharp$. 
 The derivation $\wdd(\bar{a})=-\mult \circ \dsq{a,-}$, see \eqref{Eq:dP-gen-0}, is completely determined by its value on the generators $u$ and $v$. We compute from \eqref{Eq:dP-uv} that 
 \begin{equation} \label{Eq:kxy-1}
\begin{aligned} 
\wdd(\bar{a})(u)=&-\mult (\partial_v(a)'\star \dsq{v,u} \star \partial_v(a)'') = uv \partial_v(a)''\partial_v(a)' vu\,, \\
\wdd(\bar{a})(v)=&-\mult (\partial_u(a)'\star \dsq{u,v} \star \partial_u(a)'') = -vu \partial_u(a)''\partial_u(a)' uv\,.  
 \end{aligned} 
 \end{equation} 
This directly gives the image of $\wdd$. If $\bar{a}\in \ker \wdd$,  \eqref{Eq:kxy-1}  entails $\partial_u(a)''\partial_u(a)' =0$ and 
$\partial_v(a)''\partial_v(a)' =0$. 
(We stress that these identities are independent of the chosen lift $a\in \cA$ since $\mult \circ \delta([\cA,\cA])^\sigma=0$ for any $\delta\in \DDer(\cA)$.) 
In particular, we have 
\begin{equation*}
 \mathrm{E}(a)=   \mult \circ \left(u\star \partial_u(a) + v \star \partial_v(a) \right) \in [\cA,\cA]\,,
\end{equation*}
for the Euler derivation \eqref{Eq:Euler}. 
Thus, requiring $\mathrm{E}(a)\in [\cA,\cA]$ means that all terms of $a$ of positive degree are themselves commutator, cf. \eqref{NCEulerFor}, hence its image $\bar{a}\in \cA_\sharp$ must be a constant.  
\end{proof}

\begin{lemma}     \label{Lem:kxy-2}
The kernel of the $\kk$-linear map $\wdd:\Der(\cA)\to \wBRA(\cA)_2$ is spanned by all 
$\theta_{f,g}\in \Der(\cA)$ where 
\begin{equation} \label{Eq:kxy-2-fg}
\begin{aligned}
f&=uv \partial_v(a)''\partial_v(a)' vu + \zeta uvu + \alpha u^2 + \gamma u,\\
g&=-vu \partial_u(a)''\partial_u(a)' uv +\xi vuv + \beta v^2 + \gamma v,
\end{aligned}
\end{equation}
for some $a\in \cA$ and $\alpha,\beta,\gamma,\zeta,\xi\in \kk$. 
\end{lemma}
\begin{proof}
An arbitrary element of $\Der(\cA)$ is of the form $\theta_{f,g}$ as given in \eqref{Eq:theta-fg} for fixed $f,g\in \cA$. Such a derivation belongs to $\ker \wdd$ provided that the double bracket $\wdd(\theta_{f,g})$ vanishes on the three couples $(u,u)$, $(v,v)$ and $(u,v)$. 
We first compute from \eqref{Eq:dP-gen} and \eqref{Eq:dP-uv} that 
\begin{equation*}
  \wdd(\theta_{f,g})(u,u)=
-\dsq{u,f} - \dsq{f,u}
=-(\partial_v f)' vu \otimes uv (\partial_v f)'' 
+uv (\partial_v f)''  \otimes (\partial_v f)' vu .
\end{equation*}
By looking at the left- and right-most factors in each tensor product of this expression, we see that its vanishing imposes that $\partial_v f\in u \cA \otimes \cA u$. 
Hence we can assume that $f\in u \cA u$ can be decomposed as $f=u \hat{f} u + \tilde{f}$ for $\hat{f}=\hat{f}(u,v)\in \cA$ satisfying $\partial_v(\hat{f})\neq 0$ (if nonzero) and $\tilde{f}=\tilde{f}(u)\in \kk[u]$. 
Repeating this argument, we see that we need 
$\partial_v \hat{f}\in \kk \otimes \kk \oplus v \cA \otimes \cA v$, i.e. we can in fact assume that 
$f=uv \hat{f} vu +\zeta uvu+ \tilde{f}$ for $\hat{f}=\hat{f}(u,v)\in \cA$, 
$\zeta\in \kk$, and $\tilde{f}=\tilde{f}(u)\in \kk[u]$. We then find that 
\begin{equation*}
  \wdd(\theta_{f,g})(u,u)= 
-uv (\partial_v \hat{f})' vu \otimes uv (\partial_v \hat{f})'' vu  
+uv (\partial_v \hat{f})''  vu \otimes uv (\partial_v \hat{f})' vu \,,
\end{equation*}
which vanishes if and only if $\partial_v \hat{f}=(\partial_v \hat{f})^\sigma$. By Lemma \ref{Lem:dvadva}, this means that 
\begin{equation} \label{Eq:kxy-2}
    \begin{aligned}
f= uv (\partial_v a)'' (\partial_v a)' vu + \zeta uvu + \tilde{f}\,,& \quad 
\text{ for }a\in \cA,\,\, \zeta\in \kk,\,\, \tilde{f}=\tilde{f}(u)\in \kk[u]\,; \\
g= vu (\partial_u b)'' (\partial_u b)' uv + \xi vuv + \tilde{g}\,,& \quad 
\text{ for }b\in \cA,\,\,\xi\in \kk,\,\, \tilde{g}=\tilde{g}(v)\in \kk[v]\,. 
    \end{aligned}
\end{equation}
Indeed, the form for $g$ is obtained from the vanishing of $\wdd(\theta_{f,g})(v,v)$ after repeating the same argument where $u$ and $v$ that are swapped. 

Finally, assume that $f,g$ are  given by \eqref{Eq:kxy-2} and we require the vanishing of 
\begin{align*}
\wdd(\theta_{f,g})(u,v) 
=&
(\theta_{f,g} \otimes \Id_{\mc A} + \Id_{\mc A} \otimes \theta_{f,g}) (vu \otimes uv) 
- \dsq{u, g} 
- \dsq{f , v}\\
=& \tilde{g}u \otimes uv + vu \otimes u \tilde{g}
+ v \tilde{f} \otimes uv + vu \otimes \tilde{f} v  \\
&- vu (\partial_v \hat{g})' vu \otimes uv (\partial_v \hat{g})'' uv
- (\partial_v \tilde{g})' vu \otimes uv (\partial_v \tilde{g})''  \\
&- vu (\partial_u \hat{f})'' vu \otimes uv (\partial_u \hat{f})' uv
- vu (\partial_u \tilde{f})'' \otimes  (\partial_u \tilde{f})' uv\,, 
\end{align*}
where we put $\hat{f}=(\partial_v a)'' (\partial_v a)'$ and $\hat{g}=(\partial_u b)'' (\partial_u b)'$ to ease notation. 
If $\deg_u(\tilde{f})=d>2$, there would be nonzero elements in $vu^2 \kk[u]\otimes \kk[u]u^2v$ coming from the last term of this expansion. If $\tilde{f}$ contains a constant term, we would have nonzero elements in $\kk (v\otimes uv+vu \otimes v)$ coming from the third and fourth terms of this expansion. 
Thus, we can assume that $\tilde{f}=\alpha u^2+ \gamma u$ for some $\alpha,\gamma \in \kk$. 
By an analogous argument, $\tilde{g}=\beta v^2+ \epsilon v$ for some $\beta,\epsilon \in \kk$. We can thus write 
\begin{equation*}
\wdd(\theta_{f,g})(u,v) 
=(\gamma+\epsilon) \, vu \otimes uv  
- vu (\partial_v \hat{g})' vu \otimes uv (\partial_v \hat{g})'' uv
- vu (\partial_u \hat{f})'' vu \otimes uv (\partial_u \hat{f})' uv \,. 
\end{equation*}
The vanishing of the first term imposes $\gamma+\epsilon=0$, while the other two terms cancel out if and only if 
$\partial_v \hat{g}=-(\partial_u \hat{f})^\sigma$. 
For $b_0=b+a\in \cA$, we can write $\hat{g}=(\partial_u b_0)'' (\partial_u b_0)' - (\partial_u a)'' (\partial_u a)'$, so 
\begin{align*}
   \partial_v \hat{g}=-(\partial_u \hat{f})^\sigma \quad 
   \Leftrightarrow \quad 
   \partial_v ((\partial_u b_0)'' (\partial_u b_0)')=0\,,
\end{align*}
by Lemma \ref{Lem:dudv}. 
This means that $(\partial_u b_0)'' (\partial_u b_0)'\in \kk[u]$. 
%As the kernel of the linear map $A\to A$, $b_0\mapsto (\partial_u b_0)'' (\partial_u b_0)'$, is the subspace $[A,A]\oplus \kk[v]$ 
Therefore, $u(\partial_u b_0)'' (\partial_u b_0)'\in u\kk[u]$, and by the noncommutative Euler formula \eqref{NCEulerFor}, 
$$\deg_u(b_0) b_0=(\partial_u b_0)'u(\partial_u b_0)'' \in u\kk[u]\oplus [\cA,\cA]\,.$$
This means that $b_0\in u\kk[u]\oplus \kk[v]$ modulo commutators, hence we require 
$$a(u,v)+b(u,v) = b_1(u)+b_2(v)\,\, \text{ mod }[\cA,\cA], \quad 
\text{ for } b_1\in u\kk[u],\,\, b_2\in\kk[v].$$
Since we only care about $\hat{f}=(\partial_v a)'' (\partial_v a)'$ and $\hat{g}=(\partial_u b)'' (\partial_u b)'$, we can replace $a$ by $a-b_1(u)$ then replace $b$ by $b-b_2(v)$ without changing $\hat{f}$ and $\hat{g}$. In that case $a=-b$. 
Thus, we obtain that $\wdd(\theta_{f,g})=0$ provided that $f$ and $g$ are of the form \eqref{Eq:kxy-2-fg}. 
\end{proof}

We can summarize all our previous results in the following form. 

\begin{proposition}  \label{Pr:dPH-kuv} 
For the double Poisson bracket $\dsq{-,-}$ on $\kk\langle u,v\rangle$ given by \eqref{Eq:dP-uv}, we have 
\begin{enumerate}[(1)]
\item $\widehat{\dPH}^0(\kk\langle u,v\rangle)=\dPH^0(\kk\langle u,v\rangle)=\kk$.
\item $\widehat{\dPH}^1(\kk\langle u,v\rangle)=\kk \theta_{uvu,0}\oplus \kk \theta_{0,vuv}\oplus \kk \theta_{u^2,0}\oplus \kk \theta_{0,v^2}\oplus \kk \theta_{u,-v}$;  

\noindent $\dPH^1(\kk\langle u,v\rangle)=\kk uvu \partial_u \oplus \kk vuv \partial_v \oplus
\kk u^2 \partial_u \oplus \kk v^2 \partial_v \oplus \kk (u \partial_u -v \partial_v)$.
\end{enumerate}
\end{proposition}
\begin{proof}
We directly get $\widehat{\dPH}^0(\kk\langle u,v\rangle)$ as the kernel appearing in Lemma \ref{Lem:kxy-1}. 
For $\widehat{\dPH}^1(\kk\langle u,v\rangle)$, we mod out Lemma \ref{Lem:kxy-2} by the image given in Lemma \ref{Lem:kxy-1}. 
We get $\dPH^1(\kk\langle u,v\rangle)$ from the isomorphism $\mu_1$ noting that $\mu_1 (f\partial_u+g \partial_v)=\theta_{f,g}$. 
\end{proof}

\subsection{The (completed) double Poisson cohomology on $\kk\langle u,v \rangle/(u^2,v^2)$}

Our aim is to compute $\widehat{\dPH}^0$ and $\widehat{\dPH}^1$ for the double Poisson bracket $\dsq{-,-}$ given by \eqref{Eq:dP-uv} on the truncated algebra 
$\cA^{\tr}:=\kk\langle u,v\rangle/(u^2,v^2)$. 
%Set $\kk\langle u,v \rangle/(u^2,v^2)$ from now on. 
To have a well-defined derivation $\theta_{f,g}\in \Der(\cA^{\tr})$, $f,g\in \cA^{\tr}$,  uniquely determined by  
\begin{equation}
\theta_{f,g}(u)=f\,, \quad \theta_{f,g}(v)=g\,,
\end{equation}
it is required to impose $uf+fu=0$ and $vg+gv=0$.
Any $D\in \Der(\cA^{\tr})$ is then of this form. We use the non-negative gradings $\deg_x$, $\deg_y$ and $\deg$ that descend from $\kk\langle u,v\rangle$ to $\cA^{\tr}$.

\subsubsection{ }

Let us recall that the double derivations $\partial_u,\partial_v$ \eqref{DDer-deluv} do \emph{not} exist on $\cA^{\tr}$. However, we can consider the double derivations 
\begin{equation}
 \begin{aligned} \label{DDer-dtr-uv}
 \delta_u\in \DDer(\cA^{\tr}), \qquad 
 \delta_u(u)=u\otimes u, \,\, \delta_u(v)=0,   \\
 \delta_v\in \DDer(\cA^{\tr}), \qquad 
 \delta_v(u)=0, \,\, \delta_v(v)=v\otimes v ,
\end{aligned}   
\end{equation}
which could be written as $u\partial_u u$ and $v\partial_v v$, by abuse of notation.

\begin{lemma}  \label{Lem:kxytr-1}
The image of the $\kk$-linear map $\wdd:\cA^{\tr}_\sharp\to \Der(\cA^{\tr})$ is spanned by 
the derivations $\theta_{f,g}\in \Der(\cA^{\tr})$ with 
\begin{equation} 
f=uv \partial_v(a)''\partial_v(a)' vu,\quad 
g=-vu \partial_u(a)''\partial_u(a)' uv, \,\text{ where }a\in \cA^{\tr}\,. 
\end{equation}
Moreover, $\widehat{\dPH}^0(\cA^{\tr})=\ker (\wdd:\cA^{\tr}_\sharp\to \Der(\cA^{\tr}))=\kk \oplus \kk u \oplus \kk v$. 
\end{lemma}
\begin{proof}
We can derive \eqref{Eq:kxy-1} as in Lemma \ref{Lem:kxy-1} for any $\bar{a}\in \cA^{\tr}_\sharp$ with lift $a\in \cA^{\tr}$, and it becomes 
 \begin{equation} \label{Eq:kxy-1tr} 
\wdd(\bar{a})(u)= u \delta_v(a)''\delta_v(a)' u\,, \qquad 
\wdd(\bar{a})(v)= -v \delta_u(a)'' \delta_u(a)' v\,.   
 \end{equation} 
This yields the first part of the claim. 

If $\bar{a}\in \ker \wdd$, it suffices to check the vanishing of the identities in \eqref{Eq:kxy-1tr} after assuming that the lift $a\in \cA^{\tr}$ is homogeneous with respect to $\deg$. Indeed, the map $a\mapsto u \delta_v(a)''\delta_v(a)'u$ is of degree $+3$, and the same holds by swapping $u$ and $v$. 
Seeing $\cA^{\tr}$ as a vector space, any $a\in \cA^{\tr}$ can be decomposed in the basis 
\begin{equation}
1,\,\,u,\,\,v,\,\, (uv)^{k},\,\, (vu)^{k},\,\, v(uv)^{k},\,\, u(vu)^{k},\quad k\geq 1\,,
\end{equation}
with elements of respective degrees $0$, $1$, $2k$ and $2k+1$.  
It is clear that the vanishing of \eqref{Eq:kxy-1tr} holds in $\cA^{\tr}$ for $a=1,u,v$. 
In degree $2k$, we can compute that for $a=\alpha (uv)^k-\beta (vu)^k$ where $\alpha,\beta\in \kk$,  
\begin{equation*}
u \delta_v(a)''\delta_v(a)'u=
\sum_{\ell=0}^{k-1} (\alpha-\beta) u(vu)^{k+1}\,.
\end{equation*}
Hence, the first equality in \eqref{Eq:kxy-1tr} vanishes for an element of degree $2k$ whenever it is a multiple of $(uv)^k-(vu)^k$. We also easily see that the second equality in \eqref{Eq:kxy-1tr} vanishes for $a=(uv)^k-(vu)^k$.  

In degree $2k+1$, note that both $a=v(uv)^{k}$ and $a=u(vu)^{k}$ satisfy
\begin{equation} \label{Eq:kxytr-1}
    \begin{aligned}
\delta_u(a)&%\partial_u(a)'u\otimes u\partial_u(a)''
\in\, u \cA^{\tr}\otimes \cA^{\tr}u \oplus v \cA^{\tr}\otimes \cA^{\tr}v        \,, \\
\delta_v(a)&%\partial_v(a)'v\otimes v\partial_v(a)''
\in\, u \cA^{\tr}\otimes \cA^{\tr}u \oplus v \cA^{\tr}\otimes \cA^{\tr}v \,,
    \end{aligned}
\end{equation} 
from which we deduce the vanishing of \eqref{Eq:kxy-1tr} in $\cA^{\tr}$. 
Finally, we note that the elements $(uv)^k-(vu)^k$, $v(uv)^{k}$ and $u(vu)^{k}$ can be omitted as their projection to $\cA^{\tr}_\sharp$ is just $0\in \kk$. 
Hence only $\{1,u,v\}$ forms a basis of $\ker \wdd\subset \cA^{\tr}_\sharp$.  
\end{proof}

As we are only interested in highlighting some differences between the double Poisson cohomology of $\cA^{\tr}$ and its completed version, we will not run the tedious argument of computing $\widehat{\dPH}^1(\cA^{\tr})$. We can nevertheless get  nontrivial cohomology classes based on Proposition \ref{Pr:dPH-kuv}.

\begin{lemma} \label{Lem:kxytr-2}
The cohomology class of $\theta_{u,-v}$ in $\widehat{\dPH}^1(\cA^{\tr})$ is non-trivial.   
\end{lemma}
\begin{proof}
It follows from the computations in the proof of Lemma \ref{Lem:kxy-2} that the double bracket $\wdd(\theta_{u,-v})$ vanishes on any pair of generators $\{(u,u),(v,v),(u,v)\}$, hence it is identically zero. 
We conclude as $\theta_{u,-v}$ is not in the image $\wdd(\cA^{\tr}_\sharp)$ given in Lemma \ref{Lem:kxytr-1}. 
\end{proof}

\subsubsection{ }
We want to compare the previous two lemmas with the non-completed case by computing $\dPH^0(\cA^{\tr})$ and some instances of $\dPH^1(\cA^{\tr})$. Indeed, the double Poisson bracket $\dsq{-,-}$ given by \eqref{Eq:dP-uv} on $\cA^{\tr}$ is defined by 
$P=\delta_u \delta_v \in (\mb T^\ast \cA^{\tr})_2$ where $\delta_u,\delta_v\in \DDer(\cA^{\tr})$ were introduced in \eqref{DDer-dtr-uv}. 
Hence, we can consider the complex $(\mb T^\ast \cA)_\sharp$ with the differential $\dd_P=\brSN{P,-}$. 
%As for our computations of the completed cohomology, we recall that while $\partial_u,\partial_v$ are not defined on $A^{\tr}$, it makes sense to write down $P_u(a)=(\partial_u a)'u\otimes u(\partial_u a)''$ and $P_v(a)=(\partial_v a)'v\otimes v(\partial_v a)''$ for any $a\in A$.
\begin{lemma}  \label{Lem:kxytr-3}
The image of the linear map $\dd_P:\cA^{\tr}_\sharp\to (\mb T^\ast \cA)_{\sharp,1}$ is spanned by 
\begin{equation}
\delta_v(a)''\delta_v(a)'  \delta_u - \delta_u(a)''\delta_u(a)' \delta_v, \,\text{ where }a\in \cA^{\tr}\,. 
\end{equation}
Moreover, $\dPH^0(\cA^{\tr})=\ker (\dd_P:\cA^{\tr}_\sharp\to (\mb T^\ast \cA)_{\sharp,1})=\kk \oplus \kk u \oplus \kk v$. 
\end{lemma}
\begin{proof}
For any $\bar{a}\in \cA^{\tr}_\sharp$, we have by definition $\dd_P(a)=\mult \circ \dSN{P,a}$, where on the right-hand side $a\in \cA^{\tr}$ is any lift. 
By the left Leibniz rule \eqref{Eq:gdb-Lleib} and by \eqref{Eq:dSN-b}, we can write 
\begin{align*}
\dSN{P,a}=& 
\delta_v(a)'\otimes \delta_u \delta_v(a)'' - \delta_u(a)' \delta_v\otimes \delta_v(a)''\,.
\end{align*}
We get $\dd_P(\bar{a})=\delta_v(a)'' \delta_v(a)' \delta_u -\delta_u(a)''\delta_u(a)' \delta_v$ after multiplication and  modulo commutators. We directly find the stated image. 

As for the kernel, it is not hard to see that the vanishing of $\dd_P(\bar a)\in (\mb T^\ast \cA^{\tr})_{\sharp,1}$ imposes that $\delta_v(a)'' \delta_v(a)'=0$ and 
$\delta_u(a)''\delta_u(a)'=0$. We can run an argument similar to the proof of Lemma \ref{Lem:kxytr-1} to conclude. 
\end{proof}

For the next result, we endow $(\mb T^\ast \cA^{\tr})_{1,\sharp}$ with a grading constructed on $\DDer(\cA^{\tr})$ as follows:  we put $\deg(\delta)=d\geq -1$ if $\deg(\delta(a)')+\deg(\delta(a)'')=a+d$. For example, $\deg(\delta_u)=\deg(\delta_v)=+1$. 
%We note that the gauge elements $\Delta_u,\Delta_v \in \DDer(A^{\tr})$ 
\begin{lemma} \label{Lem:kxytr-4}
There exists two non-zero cohomology classes in $\dPH^1(\cA^{\tr})$ with representatives of degree $d\leq 0$ which are spanned by the partial gauge elements $\Delta_u$ and $\Delta_v$ defined by 
$$\Delta_u(u)=u\otimes 1-1\otimes u,\,\, \Delta_u(v)=0; \quad 
\Delta_v(u)=0,\,\, \Delta_v(v)=v\otimes 1-1\otimes v\,.$$
\end{lemma}
\begin{proof}
Any $\delta\in \DDer(\cA^{\tr})$ of degree $d\leq 0$ must be such that 
$$\delta(u)=\alpha_1 u\otimes 1 + \alpha_2 1\otimes u + \alpha_3 v\otimes 1 + \alpha_4 1\otimes v + \alpha_5 1\otimes 1\,,$$
for some $\alpha_j\in \kk$. As $\delta(u^2)=0$, this imposes that only $\alpha_1=-\alpha_2$ could be nonzero. Thus, the elements of $\DDer(\cA^{\tr})$ of degree $d\leq 0$ are of the form 
\begin{equation} \label{Eq:dderAtr}
\delta=\alpha \Delta_u + \beta \Delta_v \,, \qquad \alpha,\beta \in \kk.   
\end{equation}
 By adapting the proof of \cite[Prop.~3.3.1]{VdB1}, we compute for $\delta$ as in \eqref{Eq:dderAtr},   
\begin{equation}
    \dSN{\delta_u,\delta}=\alpha (\delta_u \otimes 1 - 1 \otimes \delta_u), \quad 
    \dSN{\delta_v,\delta}=\beta (\delta_v \otimes 1 - 1 \otimes \delta_v). 
\end{equation}
Therefore, the left Leibniz rule \eqref{Eq:gdb-Lleib}  yields 
\begin{align*}
\dSN{P,\delta}=& \delta_u\star \dSN{\delta_v,\delta} + \dSN{\delta_u,\delta} \star \delta_v \\
=&-\beta (\delta_v \otimes \delta_u + 1 \otimes \delta_u \delta_v) + \alpha (\delta_u \delta_v\otimes 1 + \delta_v \otimes \delta_u)\,.
\end{align*}
Hence, $\dd_P(\delta)=(\mult \circ \dSN{P,\delta})_\sharp$ vanishes modulo graded commutators. 
\end{proof}

\begin{remark}  \label{Rem:dPH-gdPH-inj}
The sum $\Delta_u+\Delta_v \in \DDer(\cA^{\tr})$ is the gauge element 
$$\Delta\in \DDer(\cA^{\tr})\,, \qquad \Delta(a)=a\otimes 1 - 1 \otimes a, \quad \forall a\in \cA^{\tr}\,.$$
Thus, 
under the projection $\pi_1^\Delta:(\mb T^\ast \cA^{\tr})_{1,\sharp}\to {\mc D}^1_{\cA^{\tr}}$ (cf. Section~\ref{ss:gdPCoh}), the element $(\Delta_u+\Delta_v)_\sharp$ is sent to zero. 
This implies that the natural map $\dPH(\cA^{\tr})\to \gdPH(\cA^{\tr})$ in cohomology is not injective. 
\end{remark}

\subsubsection{ } 
We have computed some groups of the double Poisson cohomology $\dPH(\cA^{\tr})$ on $\cA^{\tr}:=\kk\langle u,v \rangle/(u^2,v^2)$ and its completed version $\widehat{\dPH}(\cA^{\tr})$ for the double Poisson bracket given by \eqref{Eq:dP-uv}. 
This raises several remarks. 

First, by comparing Lemmas \ref{Lem:kxytr-1} and \ref{Lem:kxytr-3}, we get $\widehat{\dPH}^0(\cA^{\tr})= \dPH^0(\cA^{\tr})$. In view of Theorem \ref{Thm:g-dPcoh2}, we have always $\dPH^0(\cA^{\tr})\subseteq \widehat{\dPH}^0(\cA^{\tr})$  since the map $\mu_0:=\Id_{\cA^{\tr}_\sharp}$ is an isomorphism, but the equality is not automatic in this case since $\mu_1$ is neither injective ($\mu_1(\Delta_{u})=\mu_1(\Delta_{v})=0$) nor surjective (nothing is mapped onto $\theta_{u,-v}$). 
Then, this last observation has serious consequences when it comes to comparing $\widehat{\dPH}^1(\cA^{\tr})$ and $\dPH^1(\cA^{\tr})$: the nontrivial cohomology classes exhibited in Lemmas \ref{Lem:kxytr-2} and \ref{Lem:kxytr-4} do not exist in the other cohomology theory. 
It is our opinion that the completed double Poisson cohomology is the most suitable theory to study in this context. Namely, 
\begin{itemize}
    \item Only $\widehat{\dPH}^1(\cA^{\tr})$ detects a nontrivial cohomology class from the non-truncated case of $\widehat{\dPH}^1(\kk\langle u,v\rangle)$ considered in Proposition \ref{Pr:dPH-kuv}. In fact, one could argue that it detects them all, since the elements $\theta_{u^2,0},\theta_{0,v^2}\in \Der(\kk\langle u,v\rangle)$ trivially descend to the zero derivation on $\cA^{\tr}$. 
    \item The double derivations $\Delta_u,\Delta_v$ detected by $\dPH^1(\cA^{\tr})$ are only nontrivial because $\DDer(\cA^{\tr})$ is smaller than $\DDer(\kk\langle u,v\rangle)$. In the latter case, we can write that $\Delta_u=\partial_u u-u\partial_u$ and $\Delta_v=\partial_v v-v\partial_v$ which are therefore trivially zero in $(\mb T^\ast \kk\langle u,v\rangle)_{\sharp,1}$. Hence these two classes are not interesting: this is precisely what we observe in $\widehat{\dPH}^1(\cA^{\tr})$ since they are sent to the zero class under $\mu_1$. 
\end{itemize}

\section{Constant non-degenerate cases for quivers} \label{Sec:Quiv}

A quiver $(Q,S,h,t)$ is given by an arrow set $Q$, a vertex set $S$, and the head and tail maps $h,t:Q\to S$. By abuse of notation, we simply label a quiver as $Q$, and we assume that the numbers of arrows $\ell:=|Q|$ and of vertices $|S|$ are finite. 
We write the arrows as $u_1,\ldots,u_\ell$ and let $I=\{1,\ldots,\ell\}$. 
The path algebra $\kk Q$ of $Q$ is generated by 
$\{u_i,e_s \mid i \in I,\, s\in S\}$, with multiplication written as concatenation, subject to the relations 
\begin{equation}
  e_s e_t=\delta_{st}, \quad  u_i=e_{t(u_i)}u_i e_{h(u_i)}, \quad i\in I, \,\, s,t\in S,
\end{equation}
and the completeness condition $ \sum_{s\in S} e_s=1$. 
Thus, $\cA:= \kk Q$ is a $B$-algebra for $B=\oplus_{s\in S} \kk e_s$. 

We fix a skewsymmetric matrix $C=(C_{ij})\in \Gl_{\ell}(\kk)$ such that $C_{ij}=0$ if either $t(u_i)\neq h(u_j)$ or $h(u_j)\neq t(u_i)$. 
Such a matrix exists provided that, for any $s,t\in I$, the number of arrows $s\to t$ is the same as the number of arrows $t\to s$; in particular $\ell$ is even. 
Then, we introduce $\dsq{-,-}:\cA^{\otimes 2}\to \cA^{\otimes 2}$ by setting 
\begin{equation} \label{Eq:Quiv-gen}
\begin{aligned}
    \dsq{e_s,u_j}&=0=\dsq{u_j,e_s}\,, \quad 
    \dsq{e_s,e_t}=0\,,&& \\
    \dsq{u_i , u_j}&= C_{ij}\, (e_{t(j)}\otimes e_{h(j)})\,, && i,j\in I,\,\, s,t\in S,
\end{aligned}
\end{equation} 
which we extend to $\cA^{\otimes 2}$ by the Leibniz rules \eqref{Eq:db-Rleib}-\eqref{Eq:db-Lleib}. 
This operation is $B$-linear, and it satisfies cyclic skewsymmetry \eqref{Eq:db-skew} because it holds on generators by skewsymmetry of $C$. Furthermore, the vanishing of the associated triple bracket \eqref{Eq:dJac} is trivial on generators since all terms are zero. 
Then, it follows that $\dsq{-,-}$ determined by \eqref{Eq:Quiv-gen} defines a double Poisson bracket, see e.g. \cite[Thm.~2.8]{DSKV}. 

\begin{remark}
  The condition   $C_{ij}=0$ unless $t(u_i)=h(u_j)$ and $h(u_j)=t(u_i)$ comes from the following calculation due to the defining relations of $\kk Q$ and \eqref{Eq:db-Lleib}:
\begin{align*}
0&=
\dsq{u_i , u_j}-\dsq{e_{t(i)} u_i e_{h(i)}, u_j} \\
&=C_{ij}\,(e_{t(j)}\otimes e_{h(j)} - e_{t(j)}e_{h(i)}\otimes e_{t(i)}e_{h(j)}) \\
&=(1-\delta_{t(j),h(i)}\delta_{t(i),h(j)})\, C_{ij}\, (e_{t(j)}\otimes e_{h(j)}). 
\end{align*}
\end{remark}

\begin{example} \label{Exa:loops}
 If $|S|=1$,   
 $\kk Q=\kk\langle u_1,\ldots,u_\ell\rangle$ is the free algebra on $\ell$ generators because the quiver $Q$ is made of $\ell$ loops. Then, the condition on $C$ is simply to be a skewsymmetric invertible matrix. 
\end{example}
\begin{example}[\cite{VdB1}] \label{Exa:QuivVdB}
Pick an arbitrary quiver $\mc Q$ with $\ell_0:=|\mc Q|$ arrows and vertex set $S$, and let $Q=\overline{\mc Q}$ be its double; i.e. we add an arrow $u_{\ell_0+i}:h(u_i)\to t(u_i)$ with the orientation opposite to $u_i$ for any $1\leq i\leq \ell_0$. 
Then, the path algebra of the double quiver $\overline{\mc Q}$ has a double Poisson bracket \cite[\S6]{VdB1}, which is of the form \eqref{Eq:Quiv-gen} if one takes  
$C\in \Gl_{2\ell_0}(\kk)$ to be the canonical symplectic matrix, i.e.  
$$C_{ij}=\left\{
\begin{array}{cl}
    \delta_{j, i+\ell_0} &\text{if }1\leq i \leq \ell_0,\\
   -\delta_{j, i-\ell_0}  &\text{if }\ell_0+1\leq i \leq 2\ell_0. 
\end{array}
\right.$$
\end{example}
\noindent The next result is analogous to  \cite[Thm.~2.18]{DSKV} and it will be proved in Subsection~\ref{ss:pfQuiver}. 
\begin{theorem}  \label{Thm:DPcoh-Quiv}
Consider $\kk Q$ with the dPA structure defined by \eqref{Eq:uu-odd-full}.  Then, 
\begin{equation*}
    \dim\left(\widehat{\dPH}^{n}(\kk Q) \right) =  \delta_{n0} \,|S|
    \,,
    \quad n\in\mb Z_{\geq0}
    \,.
\end{equation*}
\end{theorem}
If $|S|=1$ as in Example \ref{Exa:loops}, Theorem \ref{Thm:DPcoh-Quiv} is a Poisson version of the computation that will be carried out in Chapter \ref{CH:PVAexa} for $M=0$.

\subsection{Preparation} On $\cA=\kk Q$, we have the standard double derivations 
\begin{equation} \label{Eq:QuivD1}
    \frac{\partial}{\partial u_i}\in \DDer(\cA), \quad 
    \frac{\partial u_j}{\partial u_i}=\delta_{ij} e_{t(i)}\otimes e_{h(i)}, \quad i,j\in I\,.
\end{equation}
In particular, for any $f\in \cA$, 
\begin{equation} \label{Eq:QuivD2}
    \frac{\partial f}{\partial u_i}\in    \cA e_{t(i)}\otimes e_{h(i)} \cA \,.
\end{equation}
As in the previous section on the case of $2$ generators, introduce the Euler derivation $\mathrm{E}\in \Der(\cA)$ given by 
\begin{equation} \label{Eq:QuivEuler}
f \mapsto \mathrm{E}(f):= \sum_{i\in I} 
\mult \circ \left(u_i\star \frac{\partial f}{\partial u_i} \right)\,,  
\end{equation}
which satisfies for any homogeneous $f\in \cA$,
\begin{equation} \label{QuivNCEulerFor}
   \mathrm{E}(f)=\deg(f)\, f\,, 
\end{equation} 
where the degree is defined using $\deg(u_i)=1$, $\deg(e_s)=0$, for $i\in I$, $s\in S$. 

Consider the double Poisson bracket $\dsq{-,-}$ from \eqref{Eq:Quiv-gen}. 
Since $\dsq{u_i,u_j}\in B^{\otimes 2}$ for any $i,j\in I$, 
the associated differential $\wdd$ of Definition \ref{def:wdd} can be evaluated on any  $\dgal{-}\in \wBRA_B(\cA)_n$, $n\geq 1$, so that 
\begin{align}
&\wdd(\dgal{-})(u_{i_1} \otimes \ldots \otimes u_{i_{n+1}}) \nonumber  \\
=&\sum_{t=1}^{n+1} (-1)^{nt}\, \sigma^{t-1} \,   
\dsq{u_{i_t},\dgal{u_{i_{t+1}},\ldots,u_{i_{n+1}},u_{i_1},\ldots,u_{i_{t-1}}} }_{L}  \,, \label{Eq:QuivChemla} \\
=&\sum_{h\in I}\sum_{t=1}^{n+1} (-1)^{nt} C_{i_t h}\, \sigma^{t-1} \,  
\bigg( \frac{\partial }{\partial u_h}\bigg)_L
\dgal{u_{i_{t+1}},\ldots,u_{i_{n+1}},u_{i_1},\ldots,u_{i_{t-1}}}  \,,
\label{Eq:QuivChemla2}
\end{align} 
by Chemla's formula \eqref{Eq:Chemla}. 
For $n=0$, we have the next result. 

\begin{lemma} \label{Lem:DPCohQ}
If an element $\bar{f}\in \cA_\sharp$ satisfies $\wdd(\bar{f})=0$, then $\bar{f}\in B$.  
In particular, $\widehat{\dPH}^{0}(\kk Q)\simeq \kk^{|S|}$.  
\end{lemma}
\begin{proof}
Using \eqref{Eq:dP-gen-0} which we evaluate on $u_j$, $j\in I$, yields 
\begin{equation*}
  0=\wdd(\bar{f})(u_j)=
  -\sum_{i\in I} C_{ij} \, e_{t(j)}\bigg(\frac{\partial f}{\partial u_i}\bigg)'' \bigg(\frac{\partial f}{\partial u_i}\bigg)' e_{h(j)}\,.
\end{equation*}
By \eqref{Eq:QuivD2} and the defining condition on $C$, we can remove the idempotents appearing in this expression, and then by invertibility of $C$ we get  \begin{equation*}
 \mult \circ \sigma \frac{\partial f}{\partial u_i}= 0\,.
\end{equation*}
Multiplying by $u_i$ and summing over $i\in I$ yields the vanishing of 
$\mathrm{E}(f)$, cf. \eqref{Eq:QuivEuler}, modulo commutators. 
Thus, assuming without loss of generality that $f$ is homogeneous, 
$\mathrm{E}(f)= \deg(f) f=0$ modulo commutators, and therefore $\bar{f}\in B$. 
\end{proof}

The Euler derivation \eqref{Eq:QuivEuler} induces a graded decomposition 
$\cA=\oplus_{k\in \Z_{\geq 0}} \cA_k$ with 
$\cA_k=\{f\in\cA \mid \mathrm{E}(f)=k f\}$. 
Extending $\mathrm{E}$ to a derivation $\mathrm{E}:\cA^{\otimes n}\to\cA^{\otimes n}$ as in \eqref{mfold-ext}, we still have the direct sum decomposition in $\mathrm{E}$-eigenspaces:
\begin{equation}\label{QuivDelta-dec}
\cA^{\otimes n}=\bigoplus_{k\in\mb Z_{\geq0}}(\cA^{\otimes n})_k
\,,
\qquad
(\cA^{\otimes n})_k=\{f\in\mc V^{\otimes n}\mid \Delta(f)=k f\}
\,.
\end{equation}
We then obtain for any $n\geq 1$ a graded decomposition 
\begin{equation}
  \wBRA_B(\cA)_n =\bigoplus_{k\in\mb Z_{\geq0}} \wBRA_B(\cA)_{n; k}, 
\end{equation}
where $\wBRA_B(\cA)_{n; k}$ is spanned by $n$-brackets $\dgal{-}$ satisfying 
$\dgal{u_{i_1},\ldots,u_{i_n}}\in (\cA^{\otimes n})_k$ for any 
$i_1,\ldots,i_n\in I$. 
On $\wBRA_B(\cA)_0=\cA_\sharp$, we get a graded decomposition by considering the degree induced by $\cA$. 
Then, we obtain an operator $L_{\mathrm{E}}:\wBRA_B(\cA)\to \wBRA_B(\cA)$ by setting\glslink{quivL}{} 
\begin{equation}\label{QuivLDelta}
L_{\mathrm{E}}(\dgal{-})=(k+n) \, \dgal{-}\,, \quad \dgal{-}\in \wBRA_B(\cA)_{n; k}. 
\end{equation}
In particular, $L_{\mathrm{E}}$ is invertible on 
$\oplus_{(n,k)\neq (0,0)}\wBRA_B(\cA)_{n; k}$. 
Motivated by \cite[\S2.7]{DSKV}, we can introduce for any $P=(P_1,\ldots,P_\ell)$, $P_j\in \cA$,  the contraction operator 
$\iota_{P}:\wBRA_B(\cA)\to \wBRA_B(\cA)$, being identically zero on $\wBRA_B(\cA)_0$, 
with $\iota_{P}:\wBRA_B(\cA)_n\to \wBRA_B(\cA)_{n-1}$ given  by 
\begin{equation} \label{QuivIotaP}
    \begin{aligned}
&(\iota_P \dgal{-})(u_{i_1},\dots,u_{i_{n-1}}) \\
 =& \sum_{s=1}^{n-1}\sum_{j\in I}(-1)^{s+1} 
\mult_{(s,s+1)}
\big(
P_{j}
\star_s
A_{i_1,\dots,i_{s-1},j,i_s,\dots,i_{n-1}} \big)
\,, 
\end{aligned}
\end{equation}
where, for an $n$-bracket, we let 
\begin{equation} \label{Quiv:NotArray}
A_{j_1,\dots,j_{n}}  :=
\dgal{ u_{j_1},\dots,u_{j_n} }, \quad j_1,\ldots,j_n\in I\,.
\end{equation}
With that notation, the cyclic skewsymmetry \eqref{Eq:nbr-Cycl} entails 
\begin{equation}\label{Quiv:nbr-Cycl}
 A_{j_1,\dots,j_{n}} = (-1)^{(n-1)t} \sigma^t A_{j_{t+1},\dots,j_{n},j_1,\dots,j_t}\,.   
\end{equation}
\begin{lemma}
    The operator $\iota_P$ is well-defined and valued in $\wBRA_B(\cA)$. 
\end{lemma}
\begin{proof}
Fix an $n$-bracket $\dgal{-}$. 
By defining  $\iota_P \dgal{-}$ on any array of $n-1$ generators $(u_{i_1},\dots,u_{i_{n-1}})$ as in \eqref{QuivIotaP},  
we obtain a $B$-linear map $\iota_P \dgal{-}: \cA^{\otimes (n-1)}\to \cA^{\otimes (n-1)}$ by requiring its vanishing when an entry belongs to $B$ and extending through the derivation rules \eqref{Eq:nbr-DerAll} (with $n-1$ instead of $n$) for all $i=1,\ldots,n-1$. 
Thus, one only needs to check the cyclic skewsymmetry \eqref{Eq:nbr-Cycl} of a $(n-1)$-bracket, and it suffices to do so on generators. 

We start by rewriting \eqref{QuivIotaP} as 
\begin{equation} \label{QuivIotaP-bis}
    \begin{aligned}
&(\iota_P \dgal{-})(u_{i_1},\dots,u_{i_{n-1}}) \\
 =& \sum_{s=1}^{n-1}\sum_{j\in I}(-1)^{s+1} 
\sigma^{s-1}\circ \mult_{(1,2)}
\big( P_{j}\star_1 \,\sigma^{-(s-1)}
A_{i_1,\dots,i_{s-1},j,i_s,\dots,i_{n-1}} \big)
\,, 
\end{aligned}
\end{equation}
noting that $\sigma$ denotes $\sigma_{(1,\ldots,n-1)}$ after applying the multiplication, and $\sigma_{(1,\ldots,n)}$ before applying it. 
Thus, we can write 
    \begin{align*}
&(\iota_P \dgal{-})(u_{i_2},\dots,u_{i_{n-1}},u_{i_1}) \\
 =& 
% s=1 
\sum_{j\in I} 
\mult_{(1,2)}
\big( P_{j} \star_1
A_{j,i_2,\dots,i_{n-1},i_1} \big) \\
% s=2, ..., n-2
&+ \sum_{s=2}^{n-2}\sum_{j\in I}(-1)^{s+1} 
\sigma^{s-1}\circ \mult_{(1,2)}
\big( P_{j}\star_1 \,\sigma^{-(s-1)}
A_{i_2,\dots,i_{s},j,i_{s+1},\dots,i_{n-1},i_1} \big) \\
% s = n-1
&+\sum_{j\in I}(-1)^{n} 
\sigma^{n-2}\circ \mult_{(1,2)}
\big( P_{j}\star_1 \,\sigma^{-(n-2)}
A_{i_2,\dots,i_{n-1},j,i_1} \big)
\,.
\end{align*}
Thanks to \eqref{Quiv:nbr-Cycl}, we can write, 
\begin{align*}
    A_{i_2,\dots,i_{s},j,i_{s+1},\dots,i_{n-1},i_1} 
&=(-1)^{n-1} \sigma^{-1} A_{i_1,\dots,i_{s},j,i_{s+1},\dots,i_{n-1}} , \quad s=1,\ldots,n-2, \\
A_{i_2,\dots,i_{n-1},j,i_1}&=\sigma^{-2} A_{j,i_1,\dots,i_{n-1}},  
\end{align*}
which entails 
    \begin{align*}
&(\iota_P \dgal{-})(u_{i_2},\dots,u_{i_{n-1}},u_{i_1}) \\
 =& 
% s=1 
(-1)^{n}\sum_{j\in I}  (-1)
\mult_{(1,2)} 
\big( P_{j} \star_1
 \sigma^{-1}\, A_{i_1,j,i_2,\dots,i_{n-1}} \big) \\
% s=2, ..., n-2
&+(-1)^n \sum_{s=3}^{n-1}\sum_{j\in I}(-1)^{s+1} 
\sigma^{s-2}\circ \mult_{(1,2)}
\big( P_{j}\star_1 \,\sigma^{-(s-1)}
A_{i_1,i_2,\dots,i_{s-1},j,i_{s},\dots,i_{n-1}} \big) \\
% s = n-1
&+(-1)^{n} \sum_{j\in I} 
\sigma^{-1}\circ \mult_{(1,2)}
\big( P_{j}\star_1 \, 
A_{j,i_1,i_2,\dots,i_{n-1}} \big)
\,.
\end{align*}
after using $s+1$ as the summation index in the second line. 
The first and third line correspond, respectively, to the summands $s=2$ and $s=1$. We can then write 
    \begin{align*}
&(\iota_P \dgal{-})(u_{i_2},\dots,u_{i_{n-1}},u_{i_1}) \\
 =&  (-1)^n \sigma^{-1}\bigg[ \sum_{s=1}^{n-1}\sum_{j\in I}(-1)^{s+1} 
\sigma^{s-1}\circ \mult_{(1,2)}
\big( P_{j}\star_1 \,\sigma^{-(s-1)}
A_{i_1,i_2,\dots,i_{s-1},j,i_{s},\dots,i_{n-1}} \big) \bigg] \\
=&(-1)^n \sigma^{-1} (\iota_P \dgal{-})(u_{i_1},\dots,u_{i_{n-1}})\,,
\end{align*}
which is the desired equality. 
\end{proof}

For our purpose, we consider the contraction operator associated with the elements $P_j=\sum_{k\in I} (C^{-1})_{kj} u_k$. 
We write it as  
$\iota_{C}:\wBRA_B(\cA)\to \wBRA_B(\cA)$, with\glslink{quiviota}{} 
\begin{align}
&(\iota_{C} \dgal{-})(u_{i_1},\dots,u_{i_{n-1}})
\nonumber
\\ 
\label{Quiv-iota}
=&\sum_{s=1}^{n-1}\sum_{j,k\in I}(-1)^{s+1}
(C^{-1})_{kj} \mult_{(s,s+1)}
\big( u_{k} \star_s
A_{i_1,\dots,i_{s-1},j,i_s,\dots,i_{n-1}} \big) \,, 
\end{align}
using the notation \eqref{Quiv:NotArray}, and being also identically zero on $\wBRA_B(\cA)_0$. 
In particular, $\iota_{\mathrm{E},C}(\wBRA_B(\cA)_{n,k})\subset \wBRA_B(\cA)_{n-1,k+1}$. 
We can finally define the \emph{homotopy operator}\glslink{quivhm}{} 
\begin{equation}\label{Quiv:hm}
 h_{\mathrm{E},C} = (L_{\mathrm{E}})^{-1} \circ \iota_{C}:
 \bigoplus_{n>0} \bigoplus_{k\geq 0}\wBRA_B(\cA)_{n; k} \to \wBRA_B(\cA)\,.
\end{equation}
The following result is inspired by \cite[Thm.~2.18]{DSKV}, and it will be generalized as\footnote{The statements are analogous but they are not exactly the same, e.g. no sign appears in \eqref{eq:homotopy}.} Proposition \ref{Pr:Homot} in the differential setting. 
\begin{proposition} \label{Pr:QuivHomot}
For any $\dgal{-}\in\bigoplus_{n>0} \bigoplus_{k\geq 0}\wBRA_B(\cA)_{n; k}$, we have  
\begin{equation}\label{Quivhomotopy}
    (\wdd \circ h_{\mathrm{E},C} - h_{\mathrm{E},C} \circ \wdd)(\dgal{-}) = (-1)^{n+1}\dgal{-}\,.
\end{equation}
\end{proposition}
\begin{proof}
Since $\wdd(\wBRA_B(\cA)_{n; \kappa})\subset \wBRA_B(\cA)_{n+1; \kappa-1}$ by inspecting \eqref{Eq:QuivChemla}, the LHS is well-defined. 
By linearity of the operators, we can assume that $\dgal{-}\in \wBRA_B(\cA)_{n;\kappa}$ for fixed $n>0$ and $\kappa\geq 0$. 
Note also that, by \eqref{QuivLDelta}, $L_{\mathrm{E}}^{-1}$ acts by multiplication by $(n+\kappa)^{-1}$ on $\wBRA_B(\cA)_{n; \kappa}$.  
In that case, we have to check \eqref{Quivhomotopy} when evaluated on 
$u_{i_1},\dots,u_{i_n}$ for any $i_1,\ldots,i_n\in I$ which is equivalent to 
\begin{equation} \label{Quivhomotopy2}
 (\wdd \circ \iota_{C} - \iota_{C} \circ \wdd)(\dgal{-})
 (u_{i_1},\dots,u_{i_n})
 =(-1)^{n+1}(\kappa+n)\, A_{i_1,\dots,i_{n}}\,,
\end{equation}
with the notation \eqref{Quiv:NotArray}.
We start by rewriting 
\begin{align*}
& \iota_{C}(\dgal{-})(u_{i_{t+1}},\ldots,u_{i_{n}},u_{i_1},\ldots,u_{i_{t-1}}) \\   
=&(-1)^{n(t-1)} \sigma^{-(t-1)}\, \iota_{C}(\dgal{-})(u_{i_1},\stackrel{t}{\check{\ldots}},u_{i_{n}}) \\
=&\sum_{s=1}^{t-1}\sum_{j,k\in I}(-1)^{n(t-1)+s+1} \sigma^{-(t-1)}\, 
(C^{-1})_{kj} \mult_{(s,s+1)}
\bigg( u_{k} \star_s
A_{i_1,\dots,i_{s-1},j,i_s,\stackrel{t+1}{\check{\ldots}},i_{n}} \bigg)  \\
&+\sum_{j,k\in I}(-1)^{n(t-1)+t+1} \sigma^{-(t-1)}\, 
(C^{-1})_{kj} \mult_{(t,t+1)}
\big( u_{k} \star_t
A_{i_1,\dots,i_{t-1},j,i_{t+1},\dots,i_{n}} \big)  \\
&+\sum_{s=t+1}^{n-1}\sum_{j,k\in I}(-1)^{n(t-1)+s+1} \sigma^{-(t-1)}\,  
(C^{-1})_{kj} \mult_{(s,s+1)}
\bigg( u_{k} \star_s
A_{i_1,\stackrel{t}{\check{\ldots}},i_{s},j,i_{s+1},\dots,i_{n}} \bigg) 
\end{align*}
where we used the cyclic skewsymmetry of $\iota_{C}(\dgal{-})\in \wBRA_B(\cA)_{n-1}$ for the first equality, and \eqref{Quiv-iota} for the second equality. 
Due to \eqref{Eq:QuivChemla2}, the first term in the LHS of \eqref{Quivhomotopy2} becomes  
\begin{align}
%&(\wdd\circ \iota_{C}(\dgal{-})) (u_{i_1},\dots,u_{i_n})    \\
&\sum_{h\in I}\sum_{t=1}^{n} (-1)^{(n-1)t} C_{i_t h}\, \sigma^{t-1} \,  
\bigg( \frac{\partial }{\partial u_h}\bigg)_L
\iota_{C}(\dgal{-})(u_{i_{t+1}},\ldots,u_{i_{n}},u_{i_1},\ldots,u_{i_{t-1}})
\nonumber \\
=&
\sum_{s=1}^{n-1}\sum_{h,j,k\in I}(-1)^{s+1}  C_{i_n h} (C^{-1})_{kj} 
\nonumber \\
&\qquad \qquad \sigma^{-1}\, \bigg( \frac{\partial }{\partial u_h}\bigg)_L \mult_{(s,s+1)}
\big( u_{k} \star_s
A_{i_1,\dots,i_{s-1},j,i_s,\ldots,i_{n-1}} \big) 
\label{pfQuiv1n} \\
&+\sum_{t=1}^{n-1}\sum_{s=1}^{t-1}\sum_{h,j,k\in I}(-1)^{n+t+s+1}C_{i_t h} (C^{-1})_{kj} 
\nonumber \\
&\qquad \qquad \sigma^{t-1}\, \bigg( \frac{\partial }{\partial u_h}\bigg)_L  \sigma^{-(t-1)}
\mult_{(s,s+1)} \bigg( u_{k} \star_s
A_{i_1,\dots,i_{s-1},j,i_s,\stackrel{t+1}{\check{\ldots}},i_{n}} \bigg) 
\label{pfQuiv1a} \\
&+\sum_{t=1}^{n-1}\sum_{h,j,k\in I}(-1)^{n+1}C_{i_t h} (C^{-1})_{kj} 
\nonumber \\
&\qquad \qquad \sigma^{t-1}\, \bigg( \frac{\partial }{\partial u_h}\bigg)_L  \sigma^{-(t-1)}
\mult_{(t,t+1)} \bigg( u_{k} \star_t
A_{i_1,\dots,i_{t-1},j,i_{t+1},\ldots,i_{n}} \bigg) 
\label{pfQuiv1b} \\
&+\sum_{t=1}^{n-1}\sum_{s=t+1}^{n-1}\sum_{h,j,k\in I}(-1)^{n+t+s+1}C_{i_t h} (C^{-1})_{kj} 
\nonumber \\
&\qquad \qquad \sigma^{t-1}\, \bigg( \frac{\partial }{\partial u_h}\bigg)_L  \sigma^{-(t-1)}
\mult_{(s,s+1)} \bigg( u_{k} \star_s
A_{i_1,\stackrel{t}{\check{\ldots}},i_{s},j,i_{s+1},\dots,i_{n}} \bigg) 
\label{pfQuiv1c}
\end{align}
noting that $\sigma$ denotes $\sigma_{(1,\ldots,n)}$ after applying $(\partial/\partial u_h)_L$, and $\sigma_{(1,\ldots,n-1)}$ before it.  
In particular, in \eqref{pfQuiv1a}--\eqref{pfQuiv1c}, we can simply replace 
$\sigma^{t-1}\, (\partial /\partial u_h)_L  \sigma^{-(t-1)}$ by 
$(\partial /\partial u_h)_{(t)}$ due to Lemma \ref{20140606:lem}.  

Let us rewrite \eqref{pfQuiv1n} after isolating the term $s=1$ as 
\begin{align*}
\eqref{pfQuiv1n}&=
\sum_{s=2}^{n-1}\sum_{h,j,k\in I}(-1)^{s+1}  C_{i_n h} (C^{-1})_{kj} 
 \\
&\qquad \qquad \mult_{(s,s+1)} \big( u_{k} \star_s
\sigma^{-1}\, \bigg( \frac{\partial }{\partial u_h}\bigg)_L A_{i_1,\dots,i_{s-1},j,i_s,\ldots,i_{n-1}} \big)   \\
&+ \sum_{h,j,k\in I}   C_{i_n h} (C^{-1})_{kj} 
\bigg[ \delta_{hk}\,\sigma^{-1} A_{j,i_1,\ldots,i_{n-1}}  \\
&\qquad \qquad  
+\mult_{(1,2)} \big( u_{k} \star_1  \sigma^{-1} \bigg( \frac{\partial }{\partial u_h}\bigg)_L A_{j,i_1,\ldots,i_{n-1}} \big)   \\ 
&\qquad \qquad  +\mult_{(n,n+1)} \big( u_{k} \star_n  \sigma^{-1} \bigg( \frac{\partial }{\partial u_h}\bigg)_L \sigma^{-1} A_{j,i_1,\ldots,i_{n-1}} \big) 
\bigg] ,
\end{align*}
where we used \eqref{Eq:Dmult4} and \eqref{Eq:Dmult5} for the first equality (with $n$ instead of $n+1$). This becomes 
\begin{align}
&\eqref{pfQuiv1n}= (-1)^{n-1} A_{i_1,\ldots,i_{n}} \label{pfQuiv1n-B} \\
&+\sum_{s=1}^{n}\sum_{h,j,k\in I}(-1)^{s+1}  C_{i_n h} (C^{-1})_{kj}   \nonumber \\
&\qquad \qquad \mult_{(s,s+1)} \big( u_{k} \star_s
\sigma^{-1}\, \bigg( \frac{\partial }{\partial u_h}\bigg)_L A_{i_1,\dots,i_{s-1},j,i_s,\ldots,i_{n-1}} \big). \label{pfQuiv1n-A} 
\end{align}
where we used cyclic skewsymmetry \eqref{Quiv:nbr-Cycl} to get the first term and the $s=n$ summand. 
Next, we use \eqref{Eq:Dmult2} to rewrite \eqref{pfQuiv1b}  as
\begin{align}
&\sum_{t=1}^{n-1}\sum_{h,j,k\in I}(-1)^{n+1}C_{i_t h} (C^{-1})_{kj} 
\bigg( \frac{\partial }{\partial u_h}\bigg)_{(t)} \mult_{(t,t+1)} \big( u_{k} \star_t
A_{i_1,\dots,i_{t-1},j,i_{t+1},\ldots,i_{n}} \big) \nonumber \\
 =&\sum_{t=1}^{n-1}\sum_{h,j,k\in I}(-1)^{n+1}C_{i_t h} (C^{-1})_{kj}  
\nonumber \\
&\qquad \qquad \mult_{(t+1,t+2)} \bigg( u_{k} \star_{t+1}
\bigg( \frac{\partial }{\partial u_h}\bigg)_{(t)}  A_{i_1,\dots,i_{t-1},j,i_{t+1},\ldots,i_{n}} \bigg)  \label{pfQuiv1b-c}   \\  
 &+ (n-1) (-1)^{n+1}\, A_{i_1,\dots,i_{t-1},j,i_{t+1},\ldots,i_{n}} \label{pfQuiv1b-B} \\
&\sum_{t=1}^{n-1}\sum_{h,j,k\in I}(-1)^{n+1}C_{i_t h} (C^{-1})_{kj}  
\nonumber \\
&\qquad \qquad  \mult_{(t,t+1)} \bigg( u_{k} \star_t
 \bigg( \frac{\partial }{\partial u_h}\bigg)_{(t+1)} A_{i_1,\dots,i_{t-1},j,i_{t+1},\ldots,i_{n}} \bigg)   .  \label{pfQuiv1b-a}
\end{align} 
Using \eqref{Eq:Dmult1} in \eqref{pfQuiv1a} and adding \eqref{pfQuiv1b-a} as the summand $s=t$, we get after swapping the ordering of the sums  
\begin{align}
\eqref{pfQuiv1a} +\eqref{pfQuiv1b-a} =
&\sum_{s=1}^{n-1} \sum_{t=s}^{n-1}\sum_{h,j,k\in I}(-1)^{n+t+s+1}C_{i_t h} (C^{-1})_{kj} 
\nonumber \\
&\quad
\mult_{(s,s+1)} \bigg( u_{k} \star_s
\bigg( \frac{\partial }{\partial u_h}\bigg)_{(t+1)}  A_{i_1,\dots,i_{s-1},j,i_s,\stackrel{t+1}{\check{\ldots}},i_{n}} \bigg) .
\label{pfQuiv1a-X}    
\end{align}
Similarly, using \eqref{Eq:Dmult3} in \eqref{pfQuiv1c} and adding \eqref{pfQuiv1b-c} as the summand $s=t$ yields 
\begin{align}
\eqref{pfQuiv1c} +\eqref{pfQuiv1b-c} =
&\sum_{t=1}^{n-1}\sum_{s=t}^{n-1}\sum_{h,j,k\in I}(-1)^{n+t+s+1}C_{i_t h} (C^{-1})_{kj} 
\nonumber \\
&\quad \mult_{(s+1,s+2)} \bigg( u_{k} \star_{s+1} \bigg( \frac{\partial }{\partial u_h}\bigg)_{(t)} 
A_{i_1,\stackrel{t}{\check{\ldots}},i_{s},j,i_{s+1},\dots,i_{n}} \bigg) .
\label{pfQuiv1c-X}    
\end{align} 

For the second term in the LHS of \eqref{Quivhomotopy2}, we 
use \eqref{Quiv-iota} and \eqref{Eq:QuivChemla2} to get 
\begin{equation*}
    \begin{aligned}
 &(\iota_{C} \circ \wdd(\dgal{-}))
 (u_{i_1},\dots,u_{i_n})  \\
 =&\sum_{s=1}^{n} \sum_{j,k\in I} (-1)^{s+1} (C^{-1})_{kj} 
\mult_{(s,s+1)}\bigg(u_k \star_s 
\wdd(\dgal{-}) (u_{i_1},\dots,u_{i_{s-1}},u_j,u_{i_s},\dots,u_{i_n})\bigg) 
 \\
=&\sum_{s=1}^{n} \sum_{j,k,h\in I} (-1)^{s+1} (C^{-1})_{kj} 
\mult_{(s,s+1)}\bigg(u_k \star_s \bigg[ 
\\
&\qquad \quad \sum_{t=1}^{s-1} (-1)^{nt} C_{i_th} \sigma^{t-1} 
 \bigg(\frac{\partial}{\partial u_h}\bigg)_L
A_{i_{t+1},\dots,i_{s-1},j,i_{s},\dots,i_{n},i_1,\ldots,i_{t-1}} 
\\
&\qquad + (-1)^{ns} C_{jh} \sigma^{s-1} 
\bigg(\frac{\partial}{\partial u_h}\bigg)_L
A_{i_{s},\dots,i_{n},i_1,\ldots,i_{s-1}}   \\
&\qquad + \sum_{t=s+1}^{n+1} (-1)^{nt} C_{i_{t-1}h} \sigma^{t-1} 
 \bigg(\frac{\partial}{\partial u_h}\bigg)_L
A_{i_{t},\dots,i_{n},i_1,\dots,i_{s-1},j,i_{s},\ldots,i_{t-2}} 
\bigg]\bigg)  \,.
\end{aligned}
\end{equation*}
We can swap the permutation $\sigma$ and the extended double derivation $(\partial/\partial u_h)_L$ according to Lemma \ref{20140606:lem}. 
Together with cyclic skewsymmetry written as in \eqref{Quiv:nbr-Cycl}, this can be written as  
    \begin{align}
    &\sum_{s=1}^{n}\sum_{t=1}^{s-1}  \sum_{j,k,h\in I} (-1)^{n+t+s} C_{i_th} (C^{-1})_{kj} 
\mult_{(s,s+1)}\bigg(u_k \star_s  
 \bigg(\frac{\partial}{\partial u_h}\bigg)_{(t)}
A_{i_1,\stackrel{t}{\check{\ldots}},i_{s-1},j,i_{s},\dots,i_{n}} \bigg) 
\nonumber \\
&+\sum_{s=1}^{n}  \sum_{h\in I} (-1)^{n}  
\mult_{(s,s+1)}\bigg(u_h \star_s  
 \bigg(\frac{\partial}{\partial u_h}\bigg)_{(s)}
A_{i_1,\dots,i_{n}} \bigg) 
\nonumber \\
&+\sum_{s=1}^{n} (C^{-1})_{kj}  \sum_{j,k,h\in I} \mult_{(s,s+1)}\bigg(u_k \star_s  \bigg[ \nonumber \\
&\qquad \qquad \sum_{t=s+1}^{n} (-1)^{n+t+s} C_{i_{t-1}h} 
 \bigg(\frac{\partial}{\partial u_h}\bigg)_{(t)}
A_{i_1,\dots,i_{s-1},j,i_{s},\stackrel{t}{\check{\ldots}},i_{n}} \nonumber \\
&\qquad \qquad +C_{i_{n}h} \sigma^{-1}  \bigg(\frac{\partial}{\partial u_h}\bigg)_{(L)}
A_{i_1,\dots,i_{s-1},j,i_{s},\dots,i_{n-1}} 
\bigg]\bigg) 
\nonumber \\
=&\sum_{t=1}^{n-1}\sum_{s=t}^{n-1}  \sum_{j,k,h\in I} (-1)^{n+t+s+1} C_{i_th} (C^{-1})_{kj} \nonumber \\
&\qquad \qquad \mult_{(s+1,s+2)}\bigg( u_k \star_{s+1}  
 \bigg(\frac{\partial}{\partial u_h}\bigg)_{(t)}
A_{i_1,\stackrel{t}{\check{\ldots}},i_{s},j,i_{s+1},\dots,i_{n}} \bigg)  \label{pfQuiv2a} \\
&+(-1)^n \kappa \, A_{i_1,\ldots,i_n} \label{pfQuiv2b} \\
&+\sum_{s=1}^{n-1}\sum_{t=s}^{n-1}  \sum_{j,k,h\in I} (-1)^{n+t+s+1} C_{i_th} (C^{-1})_{kj} \nonumber \\
&\qquad \qquad  \mult_{(s,s+1)}\bigg( u_k \star_{s}  
 \bigg(\frac{\partial}{\partial u_h}\bigg)_{(t+1)}
A_{i_1,\dots,i_{s-1},j,i_{s},\stackrel{t+1}{\check{\ldots}},i_{n}} \bigg)  \label{pfQuiv2c} \\
&+\sum_{s=1}^{n}  \sum_{j,k,h\in I} (-1)^{s+1} C_{i_n h} (C^{-1})_{kj} \nonumber \\
&\qquad \qquad \mult_{(s,s+1)}\bigg( u_k \star_{s}  
 \sigma^{-1} \, \bigg(\frac{\partial}{\partial u_h}\bigg)_{L}
A_{i_1,\dots,i_{s-1},j,i_{s},\dots,i_{n-1}} \bigg)  \label{pfQuiv2n} 
\end{align}
To derive the second equality above, 
we changed the order of summation in the first sum and shifted the index of summation $s$, 
we used in the second sum the definition of the extended Euler operator and the fact that $\dgal{-}\in \wBRA_B(\cA)_{n,\kappa}$, 
and we shifted the index of summation $t$ in the third sum. 

Plugging the various terms in the LHS of \eqref{Quivhomotopy2}, we get the cancellations
$$ \eqref{pfQuiv1c-X} - \eqref{pfQuiv2a} = 0, \quad 
    \eqref{pfQuiv1a-X} - \eqref{pfQuiv2c} = 0, \quad 
    \eqref{pfQuiv1n-A} - \eqref{pfQuiv2n} = 0 \,,$$
which leaves us with 
$$ \eqref{pfQuiv1n-B} +  \eqref{pfQuiv1b-B} - \eqref{pfQuiv2b} = (-1)^{n+1} (\kappa+n) \, A_{i_1,\ldots,i_n}.$$
This gives the claimed identity. 
\end{proof}

\subsection{Proof of Theorem \ref{Thm:DPcoh-Quiv}} \label{ss:pfQuiver}

Since we have computed the $n=0$ case in Lemma \ref{Lem:DPCohQ}, we only have to show that 
$$\dim\left(\widehat{\dPH}^{n}(\kk Q) \right)=0, \qquad n\geq 1.$$ 
We denote $B^{n}=\ker \wdd|_{\wBRA_B(\cA)_n}$ and 
$Z^n=\im \wdd|_{\wBRA_B(\cA)_{n-1}}$, $n\geq1$, so that we have  $\widehat{\dPH}^{n}(\kk Q) =B^n/Z^n$, for $n\geq1$. 
Hence, we want the inclusion $B^n \subset Z^n$. 
 
Let $\dgal{-}\in B^n$. By Proposition \eqref{Pr:QuivHomot}, we can write 
$$\dgal{-}=(-1)^{n-1}\, \wdd \circ h_{\mathrm{E},C} (\dgal{-}).$$
This means that $\dgal{-}\in Z^n$, proving the claim. \qed

\section{Previously studied examples} \label{ss:OthExmp}

\subsection{Double Poisson cohomology}
Besides the calculations provided in this section, only few examples of (gauged) double Poisson cohomology groups are known. 
The first cases considered by Van de Weyer \cite{VdW} concern the matrix algebras $S=\Mat_{d_1}(\kk)\oplus \ldots \oplus \Mat_{d_k}(\kk)$ where $d_j\geq 1$ for each $1\leq j \leq k$. It is shown that $\dPH^0(S)$ is always trivial, and  $\dPH^1(S)$ can be explicitly computed, cf. \cite[Prop.~6]{VdW}.
In \cite[\S6]{PV}, Pichereau and Van de Weyer compute the groups $\dPH^{0,1}(\kk\langle u,v\rangle)$ for the free algebra on 2 generators equipped with 5 non-equivalent double Poisson bivectors. These results can be stated as follows, with $\del_u,\del_v$ defined as in \eqref{Eq:dPbiv-uv}. 
\begin{proposition}[\cite{PV}]
We have on $\cA=\kk\langle u,v\rangle$, 
\begin{itemize}
    \item $\dPH^{0}(\cA)= \kk[u]\oplus \kk[v]$, $\dPH^{1}(\cA)= \kk \del_u \oplus \kk[v] \del_v$ for $P_1=u \del_u^2$; 
    \item $\dPH^{0}(\cA)= \kk[u]\oplus \kk[v]$, $\dPH^{1}(\cA)= \kk \del_u \oplus \kk \del_v$ for $P_2= u\del_u^2+v\del_v^2$; 
    \item $\dPH^{0}(\cA)= \kk$, $\dPH^{1}(\cA)= \kk \del_u$ for $P_3=u\del_u^2+v\del_u \del_v$, $P_4=u\del_u^2+v\del_v \del_u$; 
    \item $\dPH^{0}(\cA)= \kk$,  
    $$\dPH^{1}(\cA)= \kk \big[(uv+vu)\del_u + v^2 \del_v \big] \oplus \kk (u\del_u+v \del_v) \oplus \kk u\del_v \oplus \del_v$$ for $P_5=u \del_u u \del_v$. 
\end{itemize}
\end{proposition} 

Let us comment on the first two cases using the fusion operation in cohomology from Section \ref{ss:fusion} since $\cA=\kk[u]\ast_\kk \kk[v]$. 
Firstly, one can see that $P_1=u \del_u^2$ is obtained by fusion of the case $(\lambda,\mu,\nu)=(1,0,0)$ in \eqref{Eq:dbr-kxP} on $\kk[u]$ with the zero bivector on $\kk[v]$. 
By comparison with Propositions \ref{Pr:dPH-kx0} and \ref{Pr:dPH-kx1}, the first summand in $\dPH^{0}(\cA)$ and in $\dPH^{1}(\cA)$ is therefore obtained from the extension of Corollary \ref{Cor:Fus}. 
In fact, since for the zero bivector  on $\kk[v]$ the cohomology groups will trivially be 
$\dPH^{0}(\kk[v])=\kk[v]$ and in $\dPH^{1}(\kk[v])=\kk[v] \del_v$, one gets that the images of the two extension maps from Corollary \ref{Cor:Fus} are entirely spanning $\dPH^{0}(\cA)$ and $\dPH^{1}(\cA)$. 
Secondly, for $P_2= u\del_u^2+v\del_v^2$, one can see that we consider the fusion of twice the case $(\lambda,\mu,\nu)=(1,0,0)$ in \eqref{Eq:dbr-kxP}. Again, the two cohomology groups $\dPH^{0}(\cA)$ and $\dPH^{1}(\cA)$ are engendered by the extension maps through Corollary \ref{Cor:Fus}.

\subsection{Gauged double Poisson cohomology}
Alekseev, Kawazumi, Kuno and Naef used the Goldman Lie bracket in \cite[\S4.2]{AKKN} to derive some gauged double Poisson cohomology groups. To state their result, let $\cA_{g,n}$ denote (the completion of) the free algebra 
$\kk\langle z_1,\ldots,z_n,x_1,y_1,\ldots,x_g,y_g \rangle$ on $n+2g$ generators. 
Introduce the double derivations $\del_{z_j},\del_{x_i},\del_{y_i}$ which are uniquely defined by 
$$\del_a(a)= 1\otimes 1,\quad \del_a(b)=0 \text{ if }b\neq a\,, \quad \forall a,b\in \{z_1,\ldots,z_n,x_1,y_1,\ldots,x_g,y_g \}. $$ 
The double Poisson bracket under consideration in \cite{AKKN} is defined from the double Poisson bivector 
\begin{equation} \label{Eq:Pi-AKKN}
    P = \sum_{j=1}^n z_j \del_{z_j} \del_{z_j} +  \sum_{i=1}^g \del_{x_i} \del_{y_i} \in (\mb T^\ast \cA_{g,n})_2 \,.
\end{equation}
\begin{theorem}[\cite{AKKN}] \label{Thm:AKKN} We have: 
\begin{itemize}
    \item $\gdPH^0(\cA_{g,n}) = Z_P((\cA_{g,n})_\sharp)$;  
\item $\gdPH^1(\cA_{g,n}) = \oplus_{j=1}^n \kk \del_{z_j}$ if $(g,n)\neq(0,1)$.
\end{itemize}
\end{theorem} 
\begin{remark}
 The assumption $(g,n)\neq(0,1)$ for $\gdPH^1(\cA_{g,n})$ does not appear in \cite[Thm.~4.6]{AKKN} but the proof explicitly makes use of this condition\footnote{We are grateful to the authors of \cite{AKKN} for confirming this fact.}. 
 In the present work, $\gdPH^1(\cA_{0,1})$ has been computed in Proposition \ref{Pr:gdPH-kx0}, as well as all the other cohomology groups for the case $(g,n)=(0,1)$.   
\end{remark}
For $n=0$, the formalism of Section \ref{Sec:Quiv} allows to compute all the groups in the (completed) double Poisson cohomology. 
\begin{proposition}
    We have $\dPH^m(\cA_{g,0})\simeq \widehat{\dPH}^m(\cA_{g,0})\simeq \delta_{m,0}\kk$.
\end{proposition}
\begin{proof}
The first identification holds by Corollary \ref{Cor:Iso-dPcoh} because $\cA_{g,0}$ is a free algebra, so each $\mu_k$ is an isomorphism.

The double Poisson bracket is of the form given by Van den Bergh \cite{VdB1} for the double of a $g$-loop quiver based at a single vertex. 
By Example \ref{Exa:QuivVdB}, such double Poisson brackets are of the form \eqref{Eq:Quiv-gen}, 
and therefore the result follows from Theorem \ref{Thm:DPcoh-Quiv} with $|S|=1$. 
\end{proof}

%%%%%%%%%%% NEW CHAPTER %%%%%%%%%%%%%%%
%%%%%%%%%%% NEW CHAPTER %%%%%%%%%%%%%%%
%%%%%%%%%%% NEW CHAPTER %%%%%%%%%%%%%%%
%%%%%%%%%%% NEW CHAPTER %%%%%%%%%%%%%%%

\chapter[Induced cohomologies on representation algebras]{Induced cohomologies on (invariant) representation algebras}
\label{CH:rep-dPA}

We fix an algebra $\cA$ over the semisimple base $B=\oplus_{s\in S} \kk e_s$ made of a complete finite set of orthogonal idempotents. 
We work with the representation algebra  $\cA_{\bf n}$ of dimension vector ${\bf n}=(n_s)_{s\in S}$, $N=\sum_{s\in S} n_s$, introduced in Section \ref{ss:Rep-Not}. 

Our aim is to explain how the `double Poisson' cohomology theories defined in the previous chapters can be induced as `standard Poisson' cohomologies on $\cA_{\bf n}$, or on the invariant algebra $\cA_{\bf n}^{\Gl_{\bf n}}$. 

\section{General results}\label{sec:7.1}
As in Subsection~\ref{ss:comPCoh}, introduce the vector space of (totally) skewsymmetric multilinear derivations of $\cA_{\bf n}$ as 
\begin{equation*}
 \mf X(\cA_{\bf n})\simeq \bigoplus_{k\geq 0}\bigwedge^k_{\cA_{\bf n}} \Der(\cA_{\bf n})\,.   
\end{equation*}  
For $1\leq u,v\leq N$, any $Q=Q_1\cdots Q_k \in (\mb T^\ast \cA)_k$ gives rise to a $k$-linear derivation 
$$Q_{uv}\in \mf X^k(\cA_{\bf n})  \,, \quad 
Q_{uv}=\sum_{1\leq u_2,\ldots,u_{k}\leq N} (Q_{1})_{uu_2}\wedge (Q_2)_{u_2u_3}\wedge \ldots \wedge (Q_k)_{u_kv}\,,$$
where for $k=1$ the derivation $Q_{uv}\in \Der(\cA_{\bf n})$ is defined by setting 
\begin{equation}
Q_{uv}(a_{pq}):=Q(a)'_{pv} \, Q(a)''_{uq}, \qquad 
\forall a\in \cA, \,\, 1\leq p,q\leq N.
\end{equation} 
We form the associated matrix $\XX(Q):=(Q_{uv})_{1\leq u,v\leq N}$.  
The action of $\Gl_{\bf n}$ on $\cA_{\bf n}$ introduced in Section \ref{ss:Rep-Not} extends to $\mf X^k(\cA_{\bf n})$, and for  $Q \in (\mb T^\ast \cA)_k$, the element $\tr(Q)=\sum_{1\leq u\leq N} Q_{uu}$ is $\Gl_{\bf n}$-invariant. 
Then, the trace map \eqref{Eq:Tr-morph0} extends as 
\begin{equation} \label{Eq:Tr-morph}
    \tr : (\mb T^\ast \cA)_{\sharp,k} \longrightarrow 
    \mf X^k(\cA_{\bf n})^{\Gl_{\bf n}},\,\, Q \mapsto \tr(Q), \quad \forall k\geq 0\,.
\end{equation}
The central result for the constructions presented in this chapter is the following. 

\begin{theorem}[\cite{VdB1},\S7.5] \label{Thm:Rep-Dbr}
If $\dsq{-,-}$ is a $B$-linear double bracket on $\cA$, then there is a unique skewsymmetric biderivation $\br{-,-}$ on $\cA_{\bf n}$ satisfying 
\begin{equation} \label{Eq:Rep-dbr}
    \br{a_{uv},b_{pq}}=\dsq{a,b}_{pv,uq}=\dsq{a,b}_{pv}'\, \dsq{a,b}''_{uq}\,. 
\end{equation}
Furthermore, 
\begin{enumerate}[(1)]
    \item if $\dsq{-,-}$ is Poisson, then  $\br{-,-}$ is a Poisson bracket;
    \item if $\dsq{-,-}=\mu_2(P)$ for some $P\in (\mb T^\ast \cA)_{\sharp,2}$, then 
    $\br{-,-}$ is defined by $\tr(P)$. 
\end{enumerate}
\end{theorem}

The graded double Poisson bracket $\dSN{-,-}$ on $\mb T^\ast \cA$ of Theorem \ref{Thm:dSN} is related to its analogue $[-,-]_{\SN}$ on $\mf X(\cA_{\bf n})$ from Subsection \ref{sec:standSN} in a similar way.

\begin{proposition}[\cite{VdB1},\S7.6-7.7] \label{Pr:Tr-Lie}
We have 
\begin{equation}
[P_{uv},Q_{pq}]_{\SN}= (\dSN{P,Q}')_{pv}\, (\dSN{P,Q}'')_{uq}\,
    \in  \mf X^{r+s-1}(\cA_{\bf n})\,, 
\end{equation}
for any $P\in (\mb T^\ast \cA)_r$, $Q\in (\mb T^\ast \cA)_s$,  with  $r,s\in \N$.
In particular 
$$[\tr(P),\tr(Q)]_{\SN}= \tr(\brSN{P,Q})\,,$$
and the trace map \eqref{Eq:Tr-morph} is a graded Lie algebra homomorphism 
$$((\mb T^\ast \cA)_\sharp,\brSN{-,-}) \longrightarrow 
(\mf X(\cA_{\bf n}),[-,-]_{\SN})$$
factoring through the graded Lie subalgebra $\mf X(\cA_{\bf n})^{\Gl_{\bf n}}$ of $\Gl_{\bf n}$-invariant skewsymmetric multiderivations. 
\end{proposition}

Similarly to \eqref{Eq:Rep-dbr}, we can uniquely  define for any $k\geq 1$ a skewsymmetric multiderivation from a $k$-bracket. 
To state the result, denote by $S^{(1)}_k$ the subgroup of the permutation group $S_k$ on $k$ elements which leave $\{1\}$ invariant. %We let $\epsilon(\tau)$ denote the signature of a permutation $\tau \in S_k$. 

\begin{theorem} \label{Thm:IndBr}
Given $\dgal{-} \in \wBRA_B(\cA)_k$, $k\geq 1$, there exists a unique skewsymmetric multiderivation  
$\tr(\dgal{-}) \in \mf X^k(\cA_{\bf n})^{\Gl_{\bf n}}$ satisfying 
\begin{equation} \label{Eq:TrBrRep}
\tr(\dgal{-})(a^1_{u_1 v_1}, \ldots, a^k_{u_k v_k})
=\sum_{\tilde{\sigma} \in S_{k}^{(1)}}\sgn(\tilde{\sigma})\dgal{a^{1},a^{\tilde{\sigma}(2)},\ldots,a^{\tilde{\sigma}(k)}}_{\tilde{\sigma}(u,v)} 
\end{equation}
for any $a^j\in \cA$, $1\leq u_j,v_j\leq N$ with $1\leq j \leq k$, 
using the notation \eqref{Eq:Not-RepIndex} with 
$$\tilde{\sigma}(u,v):=(u_{\tilde{\sigma}(k)}v_1,u_1 v_{\tilde{\sigma}(2)},\ldots,u_{\tilde{\sigma}(k-1)}v_{\tilde{\sigma}(k)}).$$
This defines a linear map  
\begin{equation} \label{Eq:Tr-morph-2}
    \tr : \wBRA_B(\cA)_k \longrightarrow 
    \mf X^k(\cA_{\bf n})^{\Gl_{\bf n}}, \quad \forall k\geq 0\,,
\end{equation}
whose restriction to $\mu_k((\mb T^\ast \cA)_k)$ coincides with \eqref{Eq:Tr-morph}.
\end{theorem}
\begin{proof}
The operation defined through \eqref{Eq:TrBrRep} is linear in each argument. 
It is clear that exchanging $a^j_{u_j v_j}$ and $a^{j+1}_{u_{j+1} v_{j+1}}$ results in changing sign for $2\leq j \leq k$. 
Cyclically permuting all the elements result in a sign $(-1)^{k-1}$ by combining \eqref{Eq:RepIndex1} and \eqref{Eq:nbr-Cycl}. 
Thus this operation is totally skewsymmetric as well. To get a well-defined multiderivation, it suffices to check  
\begin{align*}
&\tr(\dgal{-})((bc)_{u_1 v_1}, a^2_{u_2 v_2}, \ldots, a^k_{u_k v_k}) \\
=& \sum_{w} 
b_{u_1 w} \tr(\dgal{-})(c_{w v_1}, a^2_{u_2 v_2}, \ldots, a^k_{u_k v_k}) \\
& +\sum_{w} 
c_{w v_1} \tr(\dgal{-})(b_{u_1 w}, a^2_{u_2 v_2}, \ldots, a^k_{u_k v_k})
\end{align*}
which follows by combining \eqref{Eq:TrBrRep} and \eqref{Eq:nbr-DerAll} with $i=1$. 

We get the map \eqref{Eq:Tr-morph-2} by construction where, for $k=0$, $\tr: \bar{a} \mapsto \tr(a)$ as in \eqref{Eq:Tr-morph}.  
The last part of the statement follows from \cite[Lem.~5.2]{F21}, where the formula \eqref{Eq:TrBrRep} is obtained by evaluating the multivector $\tr(Q)$ assuming that $\dgal{-}=\mu_k(Q)$.
\end{proof}

\begin{remark}
For $k=1$, the formula \eqref{Eq:TrBrRep} assigns to a derivation 
$\theta \in \Der(\cA)$ another derivation $\tr(\theta)\in \Der(\cA_{\bf n})$. 
The induced derivation is simply the one obtained by the extension \eqref{ExtendMorphAn} to representation algebras. 
\end{remark}

Note that when evaluated on $\Gl_{\bf n}$-invariant elements, \eqref{Eq:TrBrRep} reduces to 
\begin{equation} \label{Eq:TrBrRep-2}
\tr(\dgal{-})(\tr(a^1), \ldots, \tr(a^k))
=\sum_{\tilde{\sigma} \in S_{k}^{(1)}} \sgn(\tilde{\sigma}) \tr\left( \mult  \dgal{a^{1},a^{\tilde{\sigma}(2)},\ldots,a^{\tilde{\sigma}(k)}}\right) . 
\end{equation} 
In particular, we see that the trace maps \eqref{Eq:Tr-morph} and \eqref{Eq:Tr-morph-2} can be considered to take value in $\mf X^k(\cA_{\bf n}^{\Gl_{\bf n}})$, the vector space of skewsymmetric $k$-linear derivation on invariant elements. 
Following \cite[\S4.1]{AKKN}, let us notice the next result. 

\begin{corollary} \label{Cor:Tr-morphD}
    The linear maps \eqref{Eq:Tr-morph} and \eqref{Eq:Tr-morph-2} induce linear maps 
\begin{equation} \label{Eq:Tr-morph-3}
    \tr : {\mc D}_\cA^k \longrightarrow 
    \mf X^k(\cA_{\bf n}^{\Gl_{\bf n}}), \quad 
    \tr : \widehat{\mc D}_\cA^k \longrightarrow 
    \mf X^k(\cA_{\bf n}^{\Gl_{\bf n}}), \quad
    \forall k\geq 0\,.
\end{equation}    
\end{corollary}
\begin{proof}
If $\dgal{-}=\mu_k((\Delta_s \alpha)_\sharp)$ for $\Delta_s$ a gauge element \eqref{Eq:Deltas} and $\alpha \in (\mb T^\ast \cA)_{k-1}$, the right-hand side of \eqref{Eq:TrBrRep-2} is easily seen to vanish using \eqref{Eq:MapMu}.  
We can conclude from the definition of ${\mc D}_\cA$ and $\widehat{\mc D}_\cA$,  cf. Section \ref{ss:gdPCoh}.    
\end{proof}

\begin{remark}  \label{Rem:DeltaRep}
If we consider the derivation $(\Delta_s)_{uv}\in \Der(\cA_{\bf n})$, $1\leq u,v\leq N$, induced by the gauge element $\Delta_s$ of \eqref{Eq:Deltas}, it acts as 
$$(\Delta_s)_{uv}(a_{pq})=(ae_s)_{pv} (e_s)_{uq} - (e_s)_{pv} (e_sa)_{uq}\,.$$
This derivation is non-trivial only for 
$n_1+\ldots+n_{s-1}+1 \leq u,v \leq n_1+\ldots+n_{s}$, 
and in that case
\begin{equation} \label{Eq:DelRep}
    (\Delta_s)_{uv} = (E^{(s)}_{vu})_{\cA_{\bf n}}, 
\end{equation}
for the infinitesimal action \eqref{Eq:InfAct} of the elementary matrix $E_{vu}^{(s)}\in \gl_{n_s}\hookrightarrow \gl_{\bf n}$ (with only nonzero entry $+1$ in position $(v,u)$), cf. \cite[Prop.~7.9.1]{VdB1}. 
Thus, one should clearly expect Corollary \ref{Cor:Tr-morphD} since a multivector field $\tr(\Delta_s \alpha)$ will act with one factor according to \eqref{Eq:InfAct}, and this action vanishes on invariant elements. 
\end{remark}

\section{From the (completed) double Poisson cohomology} \label{sec:dPH-PH}
 
\subsection{The construction of Pichereau and Van de Weyer}

Recall that an element $P\in (\mb T^\ast \cA)_2$ satisfying 
$\brSN{P,P}=0\in (\mb T^\ast \cA)_{\sharp,2}$ defines a double Poisson bracket on $\cA$, see Theorem \ref{Thm:dJac-PP}. 
In turn, $P$ induces the skewsymmetric biderivation $\tr(P) \in \mf X^2(\cA_{\bf n})$ which is Poisson. This either follows from Theorem \ref{Thm:Rep-Dbr}, or as a consequence of Proposition \ref{Pr:Tr-Lie} because the condition  $\brSN{P,P}=0$ yields $[\tr(P),\tr(P)]_{\SN}=0$. In particular, we can form the square-zero differential $\dd_{\tr(P)}:=[\tr(P),-]_{\SN}$ on  $\mf X(\cA_{\bf n})$, cf. Subsection \ref{ss:comPCoh}.

\begin{theorem}[\cite{PV,VdW}] \label{Thm:dP-rep1}
Given $P\in (\mb T^\ast \cA)_{\sharp,2}$ satisfying $\brSN{P,P}=0$, the Pichereau-Van de Weyer differential $\dd_P$ \eqref{Eq:dP-PVdW} and the differential $\dd_{\tr(P)}$ induce the morphism of complexes 
\begin{equation} \label{Eq:dP-rep1}
    \tr : ((\mb T^\ast \cA)_{\sharp}, \dd_P) \longrightarrow 
    (\mf X(\cA_{\bf n}), \dd_{\tr(P)}),  
\end{equation}
which descends to a linear map $\dPH(\cA) \to \PH(\cA_{\bf n})$ in cohomology. 

Furthermore, by $\Gl_{\bf n}$-invariance of $\tr(P)$, the above morphism of complexes factors through $\mf X(\cA_{\bf n})^{\Gl_{\bf n}}$ where it can be restricted to $\mf X(\cA_{\bf n}^{\Gl_{\bf n}})$, i.e. we get the commutative diagram 
%\pecetta{this diagram has been commented to speed up compiling}
%\begin{comment}
\begin{center}
      \begin{tikzpicture}
%%%% top row
 \node  (zero) at (-4.2,0) {$0$};
 \node  (T0) at (-2.8,0) {$(\mb T^\ast \cA)_{\sharp,0}$};
 \node  (T1) at (0,0) {$(\mb T^\ast \cA)_{\sharp,1}$};
 \node  (T2) at (3,0) {$(\mb T^\ast \cA)_{\sharp,2}$};
 \node  (T3) at (6,0) {$(\mb T^\ast \cA)_{\sharp,3}$};
 \node  (dots) at (7.8,0) {$\cdots$};
\path[->,>=angle 90,font=\small]  
   (zero) edge  (T0) ;
\path[->,>=angle 90,font=\small]  
   (T0) edge node[above] {$\dd_P$}  (T1) ;
\path[->,>=angle 90,font=\small]  
   (T1) edge node[above] {$\dd_P$}  (T2) ;
\path[->,>=angle 90,font=\small]  
   (T2) edge node[above] {$\dd_P$}  (T3) ;
\path[->,>=angle 90,font=\small]  
   (T3) edge   (dots) ;
%%%% middle
 \node  (Bzero) at (-4.2,-2) {$0$};
 \node  (B0) at (-2.8,-2) {$\cA_{\bf n}^{\Gl_{\bf n}}$};
 \node  (B1) at (0,-2) {$\mf X^1(\cA_{\bf n})^{\Gl_{\bf n}}$};
 \node  (B2) at (3,-2) {$\mf X^2(\cA_{\bf n})^{\Gl_{\bf n}}$};
 \node  (B3) at (6,-2) {$\mf X^3(\cA_{\bf n})^{\Gl_{\bf n}}$};
 \node  (Bdots) at (7.8,-2) {$\cdots$}; 
 \path[->,>=angle 90,font=\small]  
   (Bzero) edge  (B0) ;
\path[->,>=angle 90,font=\small]  
   (B0) edge node[above] {$\dd_{\tr(P)}$}  (B1) ;
\path[->,>=angle 90,font=\small]  
   (B1) edge node[above] {$\dd_{\tr(P)}$}  (B2) ;
\path[->,>=angle 90,font=\small]  
   (B2) edge node[above] {$\dd_{\tr(P)}$}  (B3) ;
\path[->,>=angle 90,font=\small]  
   (B3) edge   (Bdots) ;
%%%% vertical arrows top
\path[->,>=angle 90,font=\small] 
(T0) edge node[left] {$\tr$}  (B0) ;
\path[->,>=angle 90,font=\small] 
(T1) edge node[left] {$\tr$}  (B1) ;
\path[->,>=angle 90,font=\small] 
(T2) edge node[left] {$\tr$}   (B2) ;
\path[->,>=angle 90,font=\small] 
(T3) edge node[left] {$\tr$}   (B3) ;
%%%
%%%% bottom row 
 \node  (Czero) at (-4.2,-4) {$0$};
 \node  (C0) at (-2.8,-4) {$\cA_{\bf n}^{\Gl_{\bf n}}$};
 \node  (C1) at (0,-4) {$\mf X^1(\cA_{\bf n}^{\Gl_{\bf n}})$};
 \node  (C2) at (3,-4) {$\mf X^2(\cA_{\bf n}^{\Gl_{\bf n}})$};
 \node  (C3) at (6,-4) {$\mf X^3(\cA_{\bf n}^{\Gl_{\bf n}})$};
 \node  (Cdots) at (7.8,-4) {$\cdots$}; 
 \path[->,>=angle 90,font=\small]  
   (Czero) edge  (C0) ;
\path[->,>=angle 90,font=\small]  
   (C0) edge node[above] {$\dd_{\tr(P)}$}  (C1) ;
\path[->,>=angle 90,font=\small]  
   (C1) edge node[above] {$\dd_{\tr(P)}$}  (C2) ;
\path[->,>=angle 90,font=\small]  
   (C2) edge node[above] {$\dd_{\tr(P)}$}  (C3) ;
\path[->,>=angle 90,font=\small]  
   (C3) edge  (Cdots) ;
%%%% vertical arrows 
\path[->,>=angle 90,font=\small] 
(B0) edge   (C0) ;
\path[->,>=angle 90,font=\small] 
(B1) edge  (C1) ;
\path[->,>=angle 90,font=\small] 
(B2) edge   (C2) ;
\path[->,>=angle 90,font=\small] 
(B3) edge  (C3) ;
   \end{tikzpicture}
\end{center}
%\end{comment}
In particular, this descends to a linear map $\dPH(\cA) \to \PH(\cA_{\bf n}^{\Gl_{\bf n}})$ in cohomology.  
\end{theorem} 

\subsection{Inducing the completed double Poisson cohomology}
\label{ss:InducdP}

\subsubsection{Main statements}

\begin{theorem} \label{Thm:dP-rep2} 
Assume that $\dsq{-,-} \in \wBRA_B(\cA)_2$ is a double Poisson bracket, and let $\br{-,-}$ denote the associated Poisson bracket on $\cA_{\bf n}$ obtained through Theorem \ref{Thm:Rep-Dbr}. 
The differential $\wdd$ of Definition \ref{def:wdd} and the differential $\delta_{\cA_{\bf n},\br{-,-}}$ \eqref{Eq:Diff-Pcoh1} induce the morphism of complexes 
\begin{equation} \label{Eq:dP-rep2}
    \tr : (\wBRA_B(\cA), \wdd) \longrightarrow 
    (\mf X(\cA_{\bf n}), (-1)^{\bullet} \, \delta_{\cA_{\bf n},\br{-,-}} ),  
\end{equation}
which descends to a linear map $\widehat{\dPH}(\cA) \to {\mathrm H}_{CE}(\cA_{\bf n})$ in cohomology. 

Furthermore, by $\Gl_{\bf n}$-invariance of $\br{-,-}$, the morphism of complexes \eqref{Eq:dP-rep2} factors through $\mf X(\cA_{\bf n})^{\Gl_{\bf n}}$ where it can be restricted to $\mf X(\cA_{\bf n}^{\Gl_{\bf n}})$. 
In particular, this descends to a linear map $\widehat{\dPH}(\cA) \to {\mathrm H}_{CE}(\cA_{\bf n}^{\Gl_{\bf n}})$ in cohomology.  
\end{theorem} 

\begin{proof}
It is proved in \S\ref{ss:Proof-dP-rep2} that the map \eqref{Eq:dP-rep2} is a morphism of complexes. The other parts of the statement directly follow by definition. 
\end{proof}

\begin{remark}
The graded Lie bracket $[-,-]_{\operatorname{DP}}$ on $\wBRA_B(\cA)$ mentioned in Remark \ref{Rem:ZTdPA}
allows to adapt the first part of Theorem \ref{Thm:dP-rep2} to get a morphism of graded Lie algebras if we equip $\mf X(\cA_{\bf n})$ with $[-,-]_{\SN}$, see \cite[Thm.~4.2]{ZT2}.
\end{remark}

When the double Poisson bracket is of the form $\mu_2(P)$ with $P\in (\mb T^\ast \cA)_{\sharp,2}$, we can relate Theorems \ref{Thm:dP-rep1} and \ref{Thm:dP-rep2} as follows. 

\begin{corollary} \label{Cor:dPH-PH} 
Assume that $P\in (\mb T^\ast \cA)_{\sharp,2}$ satisfies $\brSN{P,P}=0$. 
Then, the following diagrams are commutative: 
%\pecetta{this diagram has been commented to speed up compiling}
%\begin{comment}
\begin{center}
      \begin{tikzpicture}
%%%% top row BACK
 \node   (T0) at (-2,1) {$\dPH(\cA)$};
 \node   (T1) at (1,1) {$\widehat{\dPH}(\cA)$};
 \node  (B0) at (-2,-1) {$\PH(\cA_{\bf n})$};
 \node  (B1) at (1,-1) {${\mathrm H}_{CE}(\cA_{\bf n})$};
\path[->,>=angle 90,font=\small]  
   (T0) edge node[above] {$\mu$} (T1) ;
\path[->,>=angle 90,font=\small]  
   (B0) edge node[above] {$\sim$} (B1) ;
\path[->,>=angle 90,font=\small]  
   (T0) edge node[left] {$\tr$} (B0) ;
\path[->,>=angle 90,font=\small]  
   (T1) edge node[right] {$\tr$} (B1) ;
%%% right part 
 \node   (T2) at (4,1) {$\dPH(\cA)$};
 \node   (T3) at (7.5,1) {$\widehat{\dPH}(\cA)$};
 \node  (B2) at (4,-1) {$\PH(\cA_{\bf n}^{\Gl_{\bf n}})$};
 \node  (B3) at (7.5,-1) {${\mathrm H}_{CE}(\cA_{\bf n}^{\Gl_{\bf n}})$};
\path[->,>=angle 90,font=\small]  
   (T2) edge node[above] {$\mu$} (T3) ;
\path[->,>=angle 90,font=\small]  
   (B2) edge node[above] {$\sim$} (B3) ;
\path[->,>=angle 90,font=\small]  
   (T2) edge node[left] {$\tr$} (B2) ;
\path[->,>=angle 90,font=\small]  
   (T3) edge node[right] {$\tr$} (B3) ;
   \end{tikzpicture}
\end{center}
%\end{comment}
where we use the differentials $\dd_P$ \eqref{Eq:dP-PVdW}, $\wdd$ \eqref{Eq:dP-gen}, $\dd_{\tr(P)}=[\tr(P),-]_{\SN}$ and $\delta_{\mc A_{\bf n},\{-,-\}}$ \eqref{Eq:Diff-Pcoh1} on the top left, top right, bottom left and bottom right, respectively.   
\end{corollary}
\begin{proof}
It suffices to combine the conclusion of Subsection~\ref{ss:comPCoh} with Theorems \ref{Thm:g-dPcoh2}, \ref{Thm:dP-rep1} and \ref{Thm:dP-rep2} to get the morphisms. 
We conclude since $\tr(P)$ defines the Poisson bracket obtained through Theorem \ref{Thm:Rep-Dbr} which is used at the bottom right corner. 
\end{proof}

\subsubsection{Proof: the map \eqref{Eq:dP-rep2} is a morphism of complexes} \label{ss:Proof-dP-rep2}

We work under the assumptions of Theorem \ref{Thm:dP-rep2}. 
We write $\delta := \delta_{\cA_{\bf n},\{-,-\}}$ for the Chevalley-Eilenberg differential \eqref{Eq:Diff-Pcoh1} on $\cA_{\bf n}$, where $\{-,-\}$ is the Poisson bracket induced by $\dsq{-,-}$  by Theorem \ref{Thm:Rep-Dbr}.  
Our aim is to check that the following diagram is commutative for any $k\geq 1$: 
%\pecetta{this diagram has been commented to speed up compiling}
%\begin{comment}
\begin{center}
      \begin{tikzpicture}
%%%% top row BACK
 \node   (T0) at (-2,1) {$\wBRA_B(\cA)_{k-1}$};
 \node   (T1) at (3,1) {$\wBRA_B(\cA)_{k}$};
 \node  (B0) at (-2,-1) {$\mf X^{k-1}(\cA_{\bf n})$};
 \node  (B1) at (3,-1) {$\mf X^{k}(\cA_{\bf n})$};
\path[->,>=angle 90,font=\small]  
   (T0) edge node[above] {$\wdd$} (T1) ;
\path[->,>=angle 90,font=\small]  
   (B0) edge node[above] {$(-1)^{k-1}\, \delta^{k-1}$} (B1) ;
\path[->,>=angle 90,font=\small]  
   (T0) edge node[left] {$\tr$} (B0) ;
\path[->,>=angle 90,font=\small]  
   (T1) edge node[right] {$\tr$} (B1) ; 
   \end{tikzpicture}
\end{center} 
%\end{comment}
We shall check that the two multilinear derivations obtained in this way are identical when evaluated on generators of $A_{\bf n}$, which will be sufficient to conclude. 

\medskip 

For $k=1$, take $\bar{a} \in \wBRA_B(\cA)_{0}=\cA_\sharp$ and let $a\in \cA$ be a lift of $\bar{a}$. Then,  
\begin{equation}
   \tr(\wdd(\bar{a}))(b_{uv}) 
   =-\tr(\mult \circ \dsq{a,-})(b_{uv})
   =-(\mult \circ \dsq{a,b})_{uv}\,,
\end{equation}
after using \eqref{Eq:dP-gen-0} and \eqref{Eq:TrBrRep} for any $b\in \cA$ and $1\leq u,v\leq N$. Meanwhile, 
\begin{equation}
 \delta^{0}\big(\tr(\bar{a})\big)(b_{uv}) 
   =\br{b_{uv} , \tr(\bar{a})}
   =-\sum_{1\leq j \leq N} \br{a_{jj},b_{uv}}
   =- \sum_{1\leq j \leq N} (\dsq{a,b})_{uj;jv}
\end{equation} 
after using \eqref{Eq:Diff-Pcoh1} and \eqref{Eq:Rep-dbr}. 
Due to the defining relation of $\cA_{\bf n}$, both expressions coincide.

For $k=2$, take $\theta\in \wBRA_B(\cA)_1=\Der_B(\cA)$. For any $a,b\in \cA$ and $1\leq u,v,p,q\leq N$, we compute 
\begin{align*}
&\tr(\wdd(\theta))(a_{uv},b_{pq}) 
   =(\wdd(\theta)(a,b))_{pv,uq} \\
   &=\big(-\sigma \circ \theta_L(\dsq{b,a}) + \theta_L(\dsq{a,b}) 
   -\dsq{a,\theta(b)} + \sigma \circ \dsq{b,\theta(a)} \big)_{pv,uq}\,,
\end{align*}
after using \eqref{Eq:TrBrRep} and Chemla's formula \eqref{Eq:Chemla}. 
On the other hand, 
\begin{align*}
-&\delta^{1}(\tr\theta )(a_{uv},b_{pq})  
=- \br{a_{uv} , \tr(\theta)(b_{pq})} 
+ \br{b_{pq} , \tr(\theta)(a_{uv})} 
+\tr(\theta)(\br{a_{uv} , b_{pq}})  \\
&=- (\dsq{a,\theta(b)})_{pv,uq}  + (\dsq{b,\theta(a)})_{uq,pv} 
+((\theta_L+\theta_R) \dsq{a,b})_{pv,uq} \\
&=(-\dsq{a,\theta(b)} + \sigma \circ \dsq{b,\theta(a)} 
+ \theta_L(\dsq{a,b}) - \sigma \circ \theta_L(\dsq{b,a}) )_{pv,uq} \,,
\end{align*}
after using \eqref{Eq:Diff-Pcoh1}, \eqref{Eq:Rep-dbr}, \eqref{Eq:TrBrRep}, and finally the cyclic skewsymmetry \eqref{Eq:nbr-Cycl}.  
Again, both expressions coincide.

To consider the general case $k\geq 3$, let 
$\dgal{-} \in \wBRA_B(\cA)_{k-1}$, $a^1,\ldots,a^k\in \cA$ and $1\leq u_j,v_j\leq N$ for each $j=1,\ldots,k$. 
Combining \eqref{Eq:TrBrRep} and Chemla's formula \eqref{Eq:Chemla}, we have on the one hand, 
%\begin{small}
\begin{align}
 &\tr(\wdd(\dgal{-}))(a^1_{u_1v_1}, \ldots, a^k_{u_kv_k}) 
= \sum_{i=1}^k (-1)^{(k-1)i} \sum_{\tilde{\sigma} \in S_{k}^{(1)}}\sgn(\tilde{\sigma}) 
(\mathtt{T}_1+\mathtt{T}_2)_{\tilde{\sigma}(u,v)}  \label{Eq:dPrep-1tr}
\\
&\mathtt{T}_1:=  \sigma^i \,\dgal{a^{\tilde{\sigma}(i+1)},\ldots,a^{\tilde{\sigma}(k)}, a^{1},a^{\tilde{\sigma}(2)},\ldots,a^{\tilde{\sigma}(i-2)} , \dsq{a^{\tilde{\sigma}(i-1)},a^{\tilde{\sigma}(i)}} }_L \,,  \label{Eq:dPrep-1A} \\
&\mathtt{T}_2:= \sigma^{i-1} \,\dsq{a^{\tilde{\sigma}(i)}, \dgal{ a^{\tilde{\sigma}(i+1)},\ldots,a^{\tilde{\sigma}(k)}, a^{1},a^{\tilde{\sigma}(2)},\ldots , a^{\tilde{\sigma}(i-1)} } }_{L}   \,, \label{Eq:dPrep-1B} 
\end{align}    
%\end{small}
where $\tilde{\sigma}(u,v):=(u_{\tilde{\sigma}(k)}v_1,u_1 v_{\tilde{\sigma}(2)},\ldots,u_{\tilde{\sigma}(k-1)}v_{\tilde{\sigma}(k)})$ for any $\tilde{\sigma}\in S_k^{(1)}$, and we recall that $\sigma:=\sigma_{(1 \cdots k)}$.  
On the other hand, \eqref{Eq:Diff-Pcoh1} yields 
\begin{align}
&\delta^{k-1}(\tr \dgal{-})(a^1_{u_1v_1}, \ldots, a^k_{u_k v_k}) 
\nonumber \\
=&\sum_{j=1}^k (-1)^{j-1} \, \Big\{a^j_{u_jv_j} , (\tr \dgal{-})(a^1_{u_1v_1},\stackrel{j}{\check{\dots}}, a^k_{u_k v_k})\Big\}  \label{Eq:dPrep-2A} \\
&+\sum_{1\leq j<l \leq k} (-1)^{j+l} \, (\tr \dgal{-})(\br{a^j_{u_jv_j},a^l_{u_lv_l}} ,a^1_{u_1v_1},\stackrel{j}{\check{\dots}},\stackrel{l}{\check{\dots}}, a^k_{u_k v_k}) \,. \label{Eq:dPrep-2B}
\end{align}    

We start by establishing the following equality: 
\begin{equation} \label{Eq:dPrep-Goal1}
    \eqref{Eq:dPrep-2A} = (-1)^{k-1}  \sum_{i=1}^k (-1)^{(k-1)i} \sum_{\tilde{\sigma} \in S_{k}^{(1)}}\sgn(\tilde{\sigma}) (\mathtt{T}_2)_{\tilde{\sigma}(u,v)} \,.
\end{equation}
For $1\leq i<j\leq k$, we let $S_k^{(i,j)}\subset S_k$ denote the subgroup of permutations that fix the 2 elements $i$ and $j$ separately (i.e. it does not contain $\sigma_{(ij)}$).  
We find from \eqref{Eq:TrBrRep} that 
\begin{align}
\tr(\dgal{-})(a^2_{u_2v_2},\ldots, a^k_{u_k v_k})
&=\sum_{\hat{\sigma} \in S_{k}^{(1,2)}}\sgn(\hat{\sigma})\dgal{a^{2},a^{\hat{\sigma}(3)},\ldots,a^{\hat{\sigma}(k)}}_{\hat{\sigma}(u,v)^{\stackrel{1}{\vee}}} ,
\label{Eq:dPrep-3a} \\
\tr(\dgal{-})(a^1_{u_1v_1},\stackrel{j}{\check{\dots}}, a^k_{u_k v_k})
&=\sum_{\hat{\sigma} \in S_{k}^{(1,j)}}\sgn(\hat{\sigma})\dgal{a^{1},a^{\hat{\sigma}(2)} ,\stackrel{j}{\check{\dots}}, a^{\hat{\sigma}(k)}}_{\hat{\sigma}(u,v)^{\stackrel{j}{\vee}}} ,   
\label{Eq:dPrep-3j}
\end{align}
where $2\leq j \leq k$ in \eqref{Eq:dPrep-3j}, and we use the index notation \eqref{Eq:Not-RepIndex} with 
\begin{align*}
  \hat{\sigma}(u,v)^{\stackrel{1}{\vee}}:=&(u_{\hat{\sigma}(k)}v_2,u_2 v_{\hat{\sigma}(3)},\ldots,u_{\hat{\sigma}(k-1)}v_{\hat{\sigma}(k)})  \,, \\
  \hat{\sigma}(u,v)^{\stackrel{j}{\vee}}:=&(u_{\hat{\sigma}(k)}v_1,u_1 v_{\hat{\sigma}(2)},\ldots, u_{\hat{\sigma}(j-1)} v_{\hat{\sigma}(j+1)}, \ldots,u_{\hat{\sigma}(k-1)}v_{\hat{\sigma}(k)}) \,, \quad 2\leq j \leq k\,.
\end{align*}
For the term $j=1$ in \eqref{Eq:dPrep-2A}, we get from \eqref{Eq:dPrep-3a} that 
\begin{align*}
 \eqref{Eq:dPrep-2A}_{j=1}=&
\sum_{\hat{\sigma} \in S_{k}^{(1,2)}} \sgn(\hat{\sigma})
\sum_{r=2}^k 
\left(\dsq{ a^1, \dgal{a^{2},a^{\hat{\sigma}(3)},\ldots,a^{\hat{\sigma}(k)}} }_{(r-1)}
\right)_{ \hat{\sigma}(u,v)^{\stackrel{1}{\vee}}_{r-1} } \,, \\
\hat{\sigma}(u,v)^{\stackrel{1}{\vee}}_{1}:=&(u_{\hat{\sigma}(k)}v_1,u_1 v_2, u_2 v_{\hat{\sigma}(3)}, \ldots,u_{\hat{\sigma}(k-1)}v_{\hat{\sigma}(k)}) \,, \\
\hat{\sigma}(u,v)^{\stackrel{1}{\vee}}_{r-1}:=&(u_{\hat{\sigma}(k)}v_2,\ldots, u_{\hat{\sigma}(r-1)} v_1, u_1 v_{\hat{\sigma}(r)}, \ldots,u_{\hat{\sigma}(k-1)}v_{\hat{\sigma}(k)})\,, \quad 2<r\leq k\,,
\end{align*}
after using \eqref{Eq:Rep-dbr} and the notation \eqref{20240805:eq1}. 
Using cyclic skewsymmetry \eqref{Eq:nbr-Cycl} $(k+1-r)$ times on the $r$-th summand to make $\dsq{-,-}$ act on the first tensor factor and following the corresponding rule \eqref{Eq:RepIndex1} to permute indices, we can write 
\begin{align*}
 \eqref{Eq:dPrep-2A}_{j=1}=&
\sum_{\hat{\sigma} \in S_{k}^{(1,2)}} \sgn(\hat{\sigma})
\sum_{r=2}^k (-1)^{(k-2)(k+1-r)} \\
&\quad \left(\dsq{ a^1, \dgal{a^{\hat{\sigma}(r)},\ldots,a^{\hat{\sigma}(k)},a^{2},\ldots , a^{\hat{\sigma}(r-1)} } }_{L}
\right)_{ I_r }, &
\end{align*}
for $I_2=\hat{\sigma}(u,v)$ and 
$$I_r=(u_{\hat{\sigma}(r-1)} v_1, u_1 v_{\hat{\sigma}(r)}, \ldots,u_{\hat{\sigma}(k-1)}v_{\hat{\sigma}(k)}, u_{\hat{\sigma}(k)}v_2,\ldots, u_{\hat{\sigma}(r-2)} v_{\hat{\sigma}(r-1)} ),\,\, r>2.$$
Hence $I_r=\hat{\sigma}\circ \sigma_{(2\cdots k)}^{r-2}(u,v)$ and it follows by summing over $\tilde{\sigma}=\hat{\sigma}\circ \sigma_{(2\cdots k)}^{r-2}$ that 
\begin{equation} \label{Eq:dPrep-4} 
   \eqref{Eq:dPrep-2A}_{j=1}=  \sum_{\tilde{\sigma} \in S_{k}^{(1)}} \sgn(\tilde{\sigma})
\left(\dsq{ a^1, \dgal{a^{\tilde{\sigma}(2)},\ldots,a^{\tilde{\sigma}(k)} } }_{L}
\right)_{\tilde{\sigma}(u,v)}\,.
\end{equation}
For $2\leq j \leq k$, the corresponding term in \eqref{Eq:dPrep-2A} can be written from \eqref{Eq:Rep-dbr} and \eqref{Eq:dPrep-3j} as 
\begin{align*}
 \eqref{Eq:dPrep-2A}_{j\neq1}=& (-1)^{j-1}
\sum_{\hat{\sigma} \in S_{k}^{(1,j)}} \sgn(\hat{\sigma})
\sum_{r=2}^k 
\left(\dsq{ a^j, \dgal{a^{1},a^{\hat{\sigma}(2)} ,\stackrel{j}{\check{\dots}}, a^{\hat{\sigma}(k)} } }_{(r-1)}
\right)_{ \hat{\sigma}(u,v)^{\stackrel{j}{\vee}}_{r-1} }
\end{align*}
where the index sequence $\hat{\sigma}(u,v)^{\stackrel{j}{\vee}}_{r-1}$ is given for $2 < r \leq j$ by 
$$%%% 2 < r < j+1
(u_{\hat{\sigma}(k)}v_1,u_1 v_{\hat{\sigma}(2)},\ldots, u_{\hat{\sigma}(r-2)}v_j,u_jv_{\hat{\sigma}(r-1)} , \ldots, u_{\hat{\sigma}(j-1)} v_{\hat{\sigma}(j+1)}, \ldots,u_{\hat{\sigma}(k-1)}v_{\hat{\sigma}(k)}) $$
or it is given for $j+1<r$ by 
$$%%% j+1 < r 
(u_{\hat{\sigma}(k)}v_1,u_1 v_{\hat{\sigma}(2)},\ldots, u_{\hat{\sigma}(j-1)} v_{\hat{\sigma}(j+1)}, \ldots,u_{\hat{\sigma}(r-1)}v_j,u_jv_{\hat{\sigma}(r)} , \ldots , u_{\hat{\sigma}(k-1)}v_{\hat{\sigma}(k)}) $$
or it is given for $r=2$ and $r=j+1$ by 
\begin{align*} 
%% r=2
\hat{\sigma}(u,v)^{\stackrel{j}{\vee}}_1:=&(u_{\hat{\sigma}(k)}v_j,u_jv_1,u_1 v_{\hat{\sigma}(2)},\ldots, u_{\hat{\sigma}(j-1)} v_{\hat{\sigma}(j+1)}, \ldots,u_{\hat{\sigma}(k-1)}v_{\hat{\sigma}(k)})  \,, \\
%%% 2 < r < j+1
%\hat{\sigma}(u,v)^{\stackrel{j}{\vee}}_{r-1}:=&(u_{\hat{\sigma}(k)}v_1,u_1 v_{\hat{\sigma}(2)},\ldots, u_{\hat{\sigma}(r-2)}v_j,u_jv_{\hat{\sigma}(r-1)} , \ldots, u_{\hat{\sigma}(j-1)} v_{\hat{\sigma}(j+1)}, \ldots,u_{\hat{\sigma}(k-1)}v_{\hat{\sigma}(k)}) \,, \, \, 2\!<\!r\!\leq\! j \\
%%% r = j+1
\hat{\sigma}(u,v)^{\stackrel{j}{\vee}}_{j}:=&(u_{\hat{\sigma}(k)}v_1,u_1 v_{\hat{\sigma}(2)},\ldots, u_{\hat{\sigma}(j-1)}v_j,u_jv_{\hat{\sigma}(j+1)} , \ldots,u_{\hat{\sigma}(k-1)}v_{\hat{\sigma}(k)}) \,.
%%% j+1 < r 
%\hat{\sigma}(u,v)^{\stackrel{j}{\vee}}_{r-1}:=&(u_{\hat{\sigma}(k)}v_1,u_1 v_{\hat{\sigma}(2)},\ldots, u_{\hat{\sigma}(j-1)} v_{\hat{\sigma}(j+1)}, \ldots,u_{\hat{\sigma}(r-1)}v_j,u_jv_{\hat{\sigma}(r)} , \ldots , u_{\hat{\sigma}(k-1)}v_{\hat{\sigma}(k)})  \,, \,\, j+1< r.
\end{align*} 
(For $j=k$, note that the indices are of the different form $(u_{\hat{\sigma}(k-1)}v_\ast, \ldots,u_\ast v_{\hat{\sigma}(k-1)})$; we omit to comment on that case hereafter.) 
Let us analyze the summands appearing in $\eqref{Eq:dPrep-2A}_{j\neq1}$. 
For the summand $r=2$, we can write by \eqref{Eq:RepIndex1}  
\begin{equation}  \label{Eq:dPrep-5a}
\begin{aligned}
  &\eqref{Eq:dPrep-2A}_{j\neq1}^{r=2}= (-1)^{k-1} \\
&\quad \sum_{\substack{\hat{\sigma} \in S_{k}^{(1,j)}\\ \tilde{\sigma}=\hat{\sigma}\circ \sigma_{(j \cdots k)}}}\sgn(\tilde{\sigma}) 
\left(\sigma^{-1} \dsq{ a^{\tilde{\sigma}(k)}, \dgal{a^{1},a^{\tilde{\sigma}(2)} ,\ldots, a^{\hat{\sigma}(k-1)}} }_{L}
\right)_{\tilde{\sigma}(u,v)} .     
\end{aligned} 
\end{equation}
For $2<r\leq j$, we use cyclic skewsymmetry \eqref{Eq:nbr-Cycl} $(r-2)$ times on the corresponding summand and \eqref{Eq:RepIndex1} to get 
\begin{align*}   
& \eqref{Eq:dPrep-2A}_{j\neq1}^{2<r\leq j}= (-1)^{j-1+kr} \\
&\qquad \sum_{ \hat{\sigma} \in S_{k}^{(1,j)} } \sgn(\hat{\sigma}) 
\left( \dsq{ a^{j}, \dgal{ a^{\hat{\sigma}(r-1)} ,\stackrel{j}{\check{\dots}}, a^{\hat{\sigma}(k)} , a^{1},a^{\hat{\sigma}(2)},\ldots, a^{\hat{\sigma}(r-2)} } }_{L}
\right)_{I^{r,j}} \,,
\end{align*} 
where $I^{r,j}$ is given by the index sequence  
\begin{align*}
I^{r,j}:= (u_{\hat{\sigma}(r-2)}v_j,u_jv_{\hat{\sigma}(r-1)}&, \ldots, u_{\hat{\sigma}(j-1)} v_{\hat{\sigma}(j+1)}, \ldots \\
\ldots&,u_{\hat{\sigma}(k-1)}v_{\hat{\sigma}(k)}, 
u_{\hat{\sigma}(k)}v_1,u_1 v_{\hat{\sigma}(2)},\ldots, u_{\hat{\sigma}(r-3)} v_{\hat{\sigma}(r-2)})
\,.    
\end{align*}
Notice that $I^{r,j}$ is a cyclic permutation of the sequence $\hat{\sigma} \circ \sigma_{(r-1 \cdots j)}^{-1}(u,v)$, so that using \eqref{Eq:RepIndex1} once more yields  
\begin{equation}  \label{Eq:dPrep-5b}
\begin{aligned}
&  \eqref{Eq:dPrep-2A}_{j\neq1}^{2<r\leq j}= (-1)^{(k-1)r}  \\ 
&\sum_{ \substack{\hat{\sigma} \in S_{k}^{(1,j)}\\ \tilde{\sigma}=\hat{\sigma}\circ \sigma_{(r-1 \cdots j)}^{-1}} } \sgn(\tilde{\sigma}) 
\left(\sigma^{r-2} \dsq{ a^{\tilde{\sigma}(r-1)}, \dgal{a^{\tilde{\sigma}(r)}, \ldots, a^{\tilde{\sigma}(k)}, a^{1} ,\ldots, a^{\hat{\sigma}(r-2)}} }_{L}
\right)_{\tilde{\sigma}(u,v)} .     
\end{aligned} 
\end{equation}

For $r=j+1$, we use cyclic skewsymmetry \eqref{Eq:nbr-Cycl} $(j-1)$ times and \eqref{Eq:RepIndex1}, 
\begin{align*}   
&\eqref{Eq:dPrep-2A}_{j\neq1}^{r=j+1} = (-1)^{(k-1)(j-1)} \\
&\sum_{ \hat{\sigma} \in S_{k}^{(1,j)} } \sgn(\hat{\sigma}) 
\left( \dsq{ a^{j}, \dgal{ a^{\hat{\sigma}(j+1)} ,\ldots , a^{\hat{\sigma}(k)} , a^{1},a^{\hat{\sigma}(2)},\ldots, a^{\hat{\sigma}(j-1)} } }_{L}
\right)_{I^{j}} \\
&I^{j}= 
(u_{\hat{\sigma}(j-1)}v_j,u_jv_{\hat{\sigma}(j+1)} , \ldots,u_{\hat{\sigma}(k-1)}v_{\hat{\sigma}(k)}, 
u_{\hat{\sigma}(k)}v_1,u_1 v_{\hat{\sigma}(2)},\ldots, u_{\hat{\sigma}(j-2)} v_{\hat{\sigma}(j-1)})
\,.   
\end{align*} 
Here $I^{j}$ is simply the sequence $\hat{\sigma}(u,v)$ up to a cyclic permutation, and by \eqref{Eq:RepIndex1},  
\begin{equation}  \label{Eq:dPrep-5c}
\begin{aligned}
& \eqref{Eq:dPrep-2A}_{j\neq1}^{r=j+1}= (-1)^{(k-1)(j-1)} \\
&\sum_{ \tilde{\sigma} \in S_{k}^{(1,j)} } \!\! \sgn(\tilde{\sigma}) 
\left(\sigma^{j-1} \!\dsq{ a^{\tilde{\sigma}(j)}, \dgal{ a^{\tilde{\sigma}(j+1)} ,\ldots , a^{\tilde{\sigma}(k)} , a^{1},\ldots, a^{\tilde{\sigma}(j-1)} } }_{L} 
\right)_{\tilde{\sigma}(u,v)}      
\end{aligned} 
\end{equation} 
For $j+1<r\leq k$, we use cyclic skewsymmetry \eqref{Eq:nbr-Cycl} $(r-2)$ times and \eqref{Eq:RepIndex1}, 
\begin{align*}   
&\eqref{Eq:dPrep-2A}_{j\neq1}^{r>j+1} = (-1)^{j-1+kr}  \\
&\sum_{ \hat{\sigma} \in S_{k}^{(1,j)} } \sgn(\hat{\sigma}) 
\left( \dgal{ a^{j}, \dsq{ a^{\hat{\sigma}(r)} ,\ldots , a^{\hat{\sigma}(k)} , a^{1},a^{\hat{\sigma}(2)},\stackrel{j}{\check \ldots}, a^{\hat{\sigma}(r-1)} } }_{L}
\right)_{I^{j,r}} ,  
\end{align*} 
for the multi-index 
\begin{align*}
I^{j,r}= 
(u_{\hat{\sigma}(r-1)}v_j,u_jv_{\hat{\sigma}(r)} , \ldots &,u_{\hat{\sigma}(k-1)}v_{\hat{\sigma}(k)}, 
u_{\hat{\sigma}(k)}v_1,u_1 v_{\hat{\sigma}(2)},\ldots \\
&\ldots, u_{\hat{\sigma}(j-1)} v_{\hat{\sigma}(j+1)}, \ldots, u_{\hat{\sigma}(r-2)} v_{\hat{\sigma}(r-1)})
\,.     
\end{align*}
Notice that $I^{j,r}$ is a cyclic permutation of the sequence $\hat{\sigma} \circ \sigma_{(j \cdots r-1)}(u,v)$, so \eqref{Eq:RepIndex1} implies   
\begin{equation}  \label{Eq:dPrep-5d}
\begin{aligned}
& \eqref{Eq:dPrep-2A}_{j\neq1}^{r>j+1}= (-1)^{(k-1)r} \\
&\sum_{ \substack{\hat{\sigma} \in S_{k}^{(1,j)}\\ \tilde{\sigma}=\hat{\sigma}\circ \sigma_{(j \cdots r-1)}} } \sgn(\tilde{\sigma}) 
\left(\sigma^{r-2} \dsq{ a^{\tilde{\sigma}(r-1)}, \dgal{a^{\tilde{\sigma}(r)}, \ldots, a^{\tilde{\sigma}(k)}, a^{1} ,\ldots, a^{\tilde{\sigma}(r-2)}} }_{L}
\right)_{\tilde{\sigma}(u,v)} .  
\end{aligned} 
\end{equation} 
Now, remark that summing \eqref{Eq:dPrep-5a} over $j>1$ gives 
\begin{equation}  \label{Eq:dPrep-5Y}
(-1)^{k-1}
\sum_{\hat{\sigma} \in S_{k}^{(1)}}  \sgn(\tilde{\sigma}) 
\left(\sigma^{-1} \dsq{ a^{\tilde{\sigma}(k)}, \dgal{a^{1},a^{\tilde{\sigma}(2)} ,\ldots, a^{\tilde{\sigma}(k-1)}} }_{L}
\right)_{\tilde{\sigma}(u,v)} \,. 
\end{equation} 
Making the substitutions $i=r-1 \in \{2,\ldots,j-1\}$, $i=j$, and $i=r-1 \in \{j+1,\ldots,k-1\}$ in \eqref{Eq:dPrep-5b}, \eqref{Eq:dPrep-5c} and \eqref{Eq:dPrep-5d}, respectively, then summing everything over $j>1$ yields 
\begin{equation}  \label{Eq:dPrep-5Z}
\begin{aligned}
 \sum_{i=2}^{k-1} &(-1)^{(k-1)(i-1)} 
\sum_{ \hat{\sigma} \in S_{k}^{(1)} } \sgn(\tilde{\sigma}) \\
&\left(\sigma^{i-1} \dsq{ a^{\tilde{\sigma}(i)}, \dgal{a^{\tilde{\sigma}(i+1)}, \ldots, a^{\tilde{\sigma}(k)}, a^{1} ,\ldots, a^{\tilde{\sigma}(i-1)}} }_{L}
\right)_{\tilde{\sigma}(u,v)} .     
\end{aligned} 
\end{equation}  
It remains to note that, up to a factor $(-1)^{k-1}$, one has: 
\begin{itemize}
    \item \eqref{Eq:dPrep-4} is the summand $i=1$ of $\mathtt{T}_2$ \eqref{Eq:dPrep-1B} in \eqref{Eq:dPrep-1tr}; 
    \item \eqref{Eq:dPrep-5Y} is the summand $i=k$ of $\mathtt{T}_2$ \eqref{Eq:dPrep-1B} in \eqref{Eq:dPrep-1tr}; 
    \item \eqref{Eq:dPrep-5Z} is all the remaining summands $\mathtt{T}_2$ \eqref{Eq:dPrep-1B} in \eqref{Eq:dPrep-1tr}. 
\end{itemize}
This establishes the claimed equality \eqref{Eq:dPrep-Goal1}. 

\medskip 

Our second goal is to prove that  
\begin{equation} \label{Eq:dPrep-Goal2}
    \eqref{Eq:dPrep-2B} =  (-1)^{k-1}  \sum_{i=1}^k (-1)^{(k-1)i} \sum_{\tilde{\sigma} \in S_{k}^{(1)}}\sgn(\tilde{\sigma}) (\mathtt{T}_1)_{\tilde{\sigma}(u,v)}\,.
\end{equation}

We rewrite \eqref{Eq:dPrep-2B} thanks to \eqref{Eq:Rep-dbr} and \eqref{Eq:TrBrRep} as 
\begin{align}
&\eqref{Eq:dPrep-2B}=\sum_{1\leq j<l \leq k} (-1)^{j+l} \, (\tr \dgal{-})
\left((\dsq{a^j,a^l})_{u_lv_j,u_j v_l}  ,a^1_{u_1v_1},\stackrel{j}{\check{\dots}},\stackrel{l}{\check{\dots}}, a^k_{u_k v_k}\right)  \nonumber 
\\
&= \sum_{j<l} (-1)^{j+l} \sum_{\hat{\sigma} \in S_k^{(j,l)}} \sgn(\hat{\sigma}) 
\left( \dgal{ -,  a^{\hat{\sigma}(1)},\stackrel{j}{\check{\dots}},\stackrel{l}{\check{\dots}}, a^{\hat{\sigma}(k)}}_L (\dsq{ a^j,a^l})
\right)_{I^{j,l}_<}
\label{Eq:dPrep-6} \\
&\quad- \sum_{j<l} (-1)^{j+l} \sum_{\hat{\sigma} \in S_k^{(j,l)}} \sgn(\hat{\sigma}) 
\left( \dgal{ -,  a^{\hat{\sigma}(1)},\stackrel{j}{\check{\dots}},\stackrel{l}{\check{\dots}}, a^{\hat{\sigma}(k)}}_L (\dsq{ a^l,a^j})
\right)_{I^{j,l}_>}
\label{Eq:dPrep-7}
\end{align}    
where we used cyclic skewsymmetry of $\dsq{-,-}$ for the second sum, with the index sequences  
\begin{align*}
 I^{j,l}_<=&(u_{\hat{\sigma}(k)} v_j , u_l v_{\hat{\sigma}(1)}, \ldots, 
 u_{\hat{\sigma}(j-1)}v_{\hat{\sigma}(j+1)}, \ldots, 
 u_{\hat{\sigma}(l-1)} v_{\hat{\sigma}(l+1)}, \ldots, 
 u_{\hat{\sigma}(k-1)} v_{\hat{\sigma}(k)} , u_j v_l ) , \\
I^{j,l}_>=&(u_{\hat{\sigma}(k)} v_l , u_j v_{\hat{\sigma}(1)}, \ldots, 
 u_{\hat{\sigma}(j-1)}v_{\hat{\sigma}(j+1)}, \ldots, 
 u_{\hat{\sigma}(l-1)} v_{\hat{\sigma}(l+1)}, \ldots, 
 u_{\hat{\sigma}(k-1)} v_{\hat{\sigma}(k)} , u_l v_j ) .
\end{align*} 
All the terms in \eqref{Eq:dPrep-6}-\eqref{Eq:dPrep-7} with $j=1$ can be collected as 
\begin{align*}
\eqref{Eq:dPrep-6}_{j=1}&=
\sum_{l>1} (-1)^{l+k+1} \sum_{\hat{\sigma} \in S_k^{(1,l)}} \sgn(\hat{\sigma}) 
\left( \dgal{  a^{\hat{\sigma}(2)},\stackrel{l}{\check{\dots}}, a^{\hat{\sigma}(k)} , \dsq{a^1,a^l} }_{L} 
\right)_{J^{1,l}_<} \\
\text{for}\,\, J^{1,l}_< &= (u_l v_{\hat{\sigma}(2)}, \ldots,  
 u_{\hat{\sigma}(l-1)} v_{\hat{\sigma}(l+1)}, \ldots, 
 u_{\hat{\sigma}(k-1)} v_{\hat{\sigma}(k)}, u_{\hat{\sigma}(k)} v_1  , u_1 v_l ) \,, \\
\eqref{Eq:dPrep-7}_{j=1}&= 
-\sum_{l>1} (-1)^{l+k+1} \sum_{\hat{\sigma} \in S_k^{(1,l)}} \sgn(\hat{\sigma}) 
\left( \dgal{  a^{\hat{\sigma}(2)},\stackrel{l}{\check{\dots}}, a^{\hat{\sigma}(k)} , \dsq{a^l,a^1} }_{L} 
\right)_{J^{1,l}_>} \\
\text{for}\,\, J^{1,l}_> &= (u_1 v_{\hat{\sigma}(2)}, \ldots,  
 u_{\hat{\sigma}(l-1)} v_{\hat{\sigma}(l+1)}, \ldots, 
 u_{\hat{\sigma}(k-1)} v_{\hat{\sigma}(k)}, u_{\hat{\sigma}(k)} v_l  , u_l v_1 ) \,,
\end{align*}
after using cyclic skewsymmetry \eqref{Eq:nbr-Cycl} of $\dgal{-}$ and \eqref{Eq:RepIndex1} (for the first $k-1$ tensor factors). 
Notice that $J^{1,l}_<$ is a cyclic permutation of the sequence $\hat{\sigma}\circ \sigma_{(2 \cdots l)}^{-1} (u,v)$, while $J^{1,l}_>$ is a cyclic permutation of the sequence $\hat{\sigma}\circ \sigma_{(l \cdots k)} (u,v)$. Therefore one last application of \eqref{Eq:RepIndex1} allows us to write 
\begin{align}
  &  \eqref{Eq:dPrep-6}_{j=1} + \eqref{Eq:dPrep-7}_{j=1} \nonumber \\
=& \sum_{l>1} (-1)^{k+1} \sum_{ \substack{\hat{\sigma} \in S_k^{(1,l)} \\ \tilde{\sigma} = \hat{\sigma}\circ \sigma_{(2 \cdots l)}^{-1}} } \sgn(\tilde{\sigma}) 
\left( \sigma^2 \dgal{  a^{\tilde{\sigma}(3)},\ldots, a^{\tilde{\sigma}(k)} , \dsq{a^1,a^{\tilde{\sigma}(2)}} }_{L} 
\right)_{\tilde{\sigma}(u,v)} \nonumber \\
& + \sum_{l>1}  \sum_{ \substack{\hat{\sigma} \in S_k^{(1,l)} \\ \tilde{\sigma} = \hat{\sigma}\circ \sigma_{(l \cdots k)}} } \sgn(\tilde{\sigma}) 
\left( \sigma \dgal{  a^{\tilde{\sigma}(2)},\ldots, a^{\tilde{\sigma}(k-1)} , \dsq{a^{\tilde{\sigma}(k)},a^1 } }_{L} 
\right)_{\tilde{\sigma}(u,v)} \nonumber  \\
=&(-1)^{k-1} 
\sum_{ \tilde{\sigma} \in S_k^{(1)} } \sgn(\tilde{\sigma}) 
\left( \sigma^2 \dgal{  a^{\tilde{\sigma}(3)},\ldots, a^{\tilde{\sigma}(k)} , \dsq{a^1,a^{\tilde{\sigma}(2)}} }_{L} 
\right)_{\tilde{\sigma}(u,v)} \label{Eq:dPrep-8a} \\
&+ \sum_{ \tilde{\sigma} \in S_k^{(1)} } \sgn(\tilde{\sigma}) 
 \left( \sigma \dgal{  a^{\tilde{\sigma}(2)},\ldots, a^{\tilde{\sigma}(k-1)} , \dsq{a^{\tilde{\sigma}(k)},a^1 } }_{L} 
\right)_{\tilde{\sigma}(u,v)} . \label{Eq:dPrep-8b}
\end{align}

Next, we fix $j,l$ such that $1<j<l\leq k$. (This is why we dealt with the cases $k=1,2$ separately from $k\geq 3$.) 
By definition, for any $\hat{\sigma}\in S_k^{(j,l)}$, there exists an element $r\neq j,l$ such that $\hat{\sigma}(r)=1$. We shall inspect the three cases where $1\leq r<j$, $j<r<l$, or $r>l$. 
We denote by $S^{(1,j,l)}_k$ the subgroup of $S_k$ made of the permutations fixing the 3 elements $1$, $j$ and $l$ separately. %(i.e. $\sigma_{(jl)}$ or $\sigma_{(1jl)}$ are not in that group). 

%%% first type of terms: \eqref{Eq:dPrep-6} with  1<j<l
We start with \eqref{Eq:dPrep-6}. 
A term with $1<j<l$ fixed can be written as 
\begin{align}
&\eqref{Eq:dPrep-6}_{j,l}
= 
(-1)^{j+l+k} \sum_{\hat{\sigma} \in S_k^{(j,l)}} \sgn(\hat{\sigma}) 
\left( \ldb a^{\hat{\sigma}(1)},\stackrel{j}{\check{\dots}},\stackrel{l}{\check{\dots}}, a^{\hat{\sigma}(k)}, \dsq{ a^j,a^l} \rdb_L   
\right)_{J^{j,l}_<} ,
\label{Eq:dPrep-6jl}  \\
&\hspace{-0.2cm}J^{j,l}_< =
(u_l v_{\hat{\sigma}(1)}, \ldots, 
 u_{\hat{\sigma}(j-1)}v_{\hat{\sigma}(j+1)}, \ldots, 
 u_{\hat{\sigma}(l-1)} v_{\hat{\sigma}(l+1)}, \ldots, 
 u_{\hat{\sigma}(k-1)} v_{\hat{\sigma}(k)} , 
 u_{\hat{\sigma}(k)} v_j , u_j v_l ),  \nonumber 
\end{align}    
after using cyclic skewsymmetry \eqref{Eq:nbr-Cycl} of $\dgal{-}$ and \eqref{Eq:RepIndex1}. 

% First case for \eqref{Eq:dPrep-6} with  1<j<l
In the first case, consider $\hat{\sigma} \in S_k^{(j,l)}$ for which there exists $1\leq r<j$ with $\hat{\sigma}(r)=1$. To understand the corresponding summand of \eqref{Eq:dPrep-6jl}, note that 
\begin{align*} 
&\sgn(\hat{\sigma}) \,
\ldb a^{\hat{\sigma}(1)},\stackrel{j}{\check{\dots}},\stackrel{l}{\check{\dots}}, a^{\hat{\sigma}(k)}, \dsq{ a^j,a^l} \rdb_L \\
=&(-1)^{r-1} \sgn(\dot{\sigma})  \, 
\ldb a^{\dot{\sigma}(2)},\ldots, a^{\dot{\sigma}(r)},a^1, a^{\dot{\sigma}(r+1)},\stackrel{j}{\check{\dots}},\stackrel{l}{\check{\dots}}, a^{\dot{\sigma}(k)}, \dsq{a^j,a^l} \rdb_L\,, 
\end{align*} 
where $\hat{\sigma}=\dot{\sigma} \circ \sigma_{(1,\ldots,r)}$ with $\dot{\sigma}\in S_k^{(1,j,l)}$. 
This can be equivalently written as 
\begin{align*} 
(-1)^{l+j+r} \sgn(\ddot{\sigma})  
\ldb a^{\ddot{\sigma}(2)},\ldots, a^{\ddot{\sigma}(r)},a^1, a^{\ddot{\sigma}(r+1)},\stackrel{j,j+1}{\check{\dots}}, a^{\ddot{\sigma}(k)}, \llbracket a^j,a^{\ddot{\sigma}(j+1)} \rrbracket \rdb_L\,, 
\end{align*} 
after letting $\dot{\sigma}=\ddot{\sigma} \circ \sigma_{(j+1,\ldots,l)}$, then 
\begin{align*} 
(-1)^{l+j+r} \sgn(\dddot{\sigma})  
\ldb a^{\dddot{\sigma}(2)},\ldots, a^{\dddot{\sigma}(r)},a^1, a^{\dddot{\sigma}(r+1)},\ldots, a^{\dddot{\sigma}(k-2)}, \llbracket a^{\dddot{\sigma}(k-1)},a^{\dddot{\sigma}(k)} \rrbracket \rdb_L\,, 
\end{align*} 
after further letting $\ddot{\sigma}=\dddot{\sigma} \circ \sigma_{(j,\ldots,k)}^{-2}$. A final relabeling yields 
\begin{align*} 
&(-1)^{l+j+r}(-1)^{k(r-1)} \sgn(\tilde{\sigma})   \\
&\ldb a^{\tilde{\sigma}(k-r+2)},\ldots, a^{\tilde{\sigma}(k)},a^1, a^{\tilde{\sigma}(2)},\ldots, a^{\tilde{\sigma}(k-r-1)}, 
\llbracket a^{\tilde{\sigma}(k-r)},a^{\tilde{\sigma}(k-r+1)} \rrbracket \rdb_L\,, 
\end{align*} 
by setting $\dddot{\sigma}=\tilde{\sigma} \circ \sigma_{(2,\ldots,k)}^{k-r}$. 
In summary, we introduced $\tilde{\sigma}\in S_k^{(1)}$ according to 
$$\hat{\sigma}=\tilde{\sigma} \circ \sigma_{(2,\ldots,k)}^{k-r} \circ \sigma_{(j,\ldots,k)}^{-2} \circ \sigma_{(j+1,\ldots,l)} \circ \sigma_{(1,\ldots,r)},$$
so that $\tilde{\sigma}(k-r)=j$ and $\tilde{\sigma}(k-r+1)=l$. 
Applying this series of transformations to the corresponding term appearing in \eqref{Eq:dPrep-6jl} turns the index sequence $J_<^{j,l}$ into 
$$(u_{\tilde{\sigma}(k-r+1)} v_{\tilde{\sigma}(k-r+2)}, \ldots , 
u_{\tilde{\sigma}(k)} v_1, \ldots , 
u_{\tilde{\sigma}(k-r)} v_{\tilde{\sigma}(k-r+1)})$$
which is a cyclic permutation of the sequence $\tilde{\sigma}(u,v)$.  
Gathering these manipulations,  we can write the terms under scrutiny that appear in \eqref{Eq:dPrep-6jl} as 
\begin{align} 
&\sum_{r\geq 1} \sum_{\tilde{\sigma} \in T_r}
(-1)^{(k-1)r} \sgn(\tilde{\sigma})  \nonumber \\ 
&\left( \sigma^{k-r+1} 
\dgal{ a^{\tilde{\sigma}(k-r+2)},\ldots, a^{\tilde{\sigma}(k)},a^1,\ldots, a^{\tilde{\sigma}(k-r-1)}, 
\llbracket a^{\tilde{\sigma}(k-r)},a^{\tilde{\sigma}(k-r+1)} \rrbracket }_L
\right)_{\tilde{\sigma}(u,v)}  \nonumber \\
&\text{ for}\quad T_r:=\{\tilde{\sigma} \in S^{(1)}_k  \mid r<\tilde{\sigma}(k-r) < \tilde{\sigma}(k-r+1)\}\,. \label{Eq:dPrep-9a} 
\end{align}

% Second case for \eqref{Eq:dPrep-6} with  1<j<l
In the second case, consider $\hat{\sigma} \in S_k^{(j,l)}$ for which there exists $j< r<l$ with $\hat{\sigma}(r)=1$. To understand the corresponding summand of \eqref{Eq:dPrep-6jl}, note that 
\begin{align*} 
&\sgn(\hat{\sigma}) \,
\ldb a^{\hat{\sigma}(1)},\stackrel{j}{\check{\dots}},\stackrel{l}{\check{\dots}}, a^{\hat{\sigma}(k)}, \dsq{a^j,a^l} \rdb_L \\
=&(-1)^{j+l-1} \sgn(\dot{\sigma}) \, 
\ldb a^{\dot{\sigma}(1)},\ldots, a^{\dot{\sigma}(l-2)}, a^{\dot{\sigma}(l+1)},\ldots, a^{\dot{\sigma}(k)}, \llbracket a^{\dot{\sigma}(l-1)},a^l\rrbracket \rdb_L\,, 
\end{align*} 
where $\hat{\sigma}=\dot{\sigma} \circ \sigma_{(j \cdots l-1)}^{-1}$ with $\dot{\sigma}\in S_k^{(l)}$ satisfying $\dot{\sigma}(r-1)=1$. 
This can be equivalently written as 
\begin{align*} 
(-1&)^{l+j+r-1} \sgn(\ddot{\sigma}) \\  
&\ldb a^{\ddot{\sigma}(2)},\ldots, a^{\ddot{\sigma}(r-1)},a^1, a^{\ddot{\sigma}(r)}, \ldots, a^{\ddot{\sigma}(l-2)}, a^{\ddot{\sigma}(l+1)},\ldots, a^{\ddot{\sigma}(k)}, 
\llbracket a^{\ddot{\sigma}(l-1)},a^{l} \rrbracket \rdb_L\,, 
\end{align*} 
after letting $\dot{\sigma}=\ddot{\sigma} \circ \sigma_{(1 \cdots r-1)}$ where $\ddot{\sigma}\in S_k^{(1,l)}$, then 
\begin{align*} 
(-1&)^{l+j+r-1} \sgn(\dddot{\sigma}) \\
&\ldb a^{\dddot{\sigma}(2)},\ldots, a^{\dddot{\sigma}(r-1)},a^1, a^{\dddot{\sigma}(r)},\ldots, a^{\dddot{\sigma}(k-2)}, \llbracket a^{\dddot{\sigma}(k-1)},a^{\dddot{\sigma}(k)} \rrbracket \rdb_L\,, 
\end{align*} 
after further letting $\ddot{\sigma}=\dddot{\sigma} \circ \sigma_{(l-1 \cdots k)}^{-2}$.
A final relabelling yields 
\begin{align*} 
(-1&)^{l+j-1}(-1)^{(k-1)r} \sgn(\tilde{\sigma})  \\
&\ldb a^{\tilde{\sigma}(k-r+3)},\ldots, a^{\tilde{\sigma}(k)},a^1, a^{\tilde{\sigma}(2)},\ldots, a^{\tilde{\sigma}(k-r)}, 
\llbracket a^{\tilde{\sigma}(k-r+1)},a^{\tilde{\sigma}(k-r+2)} \rrbracket \rdb_L 
\end{align*} 
by setting $\dddot{\sigma}=\tilde{\sigma} \circ \sigma_{(2,\ldots,k)}^{k-(r-1)}$. 
In summary, we introduced $\tilde{\sigma}\in S_k^{(1)}$ according to 
$$\hat{\sigma}=\tilde{\sigma} \circ \sigma_{(2,\ldots,k)}^{k-(r-1)} \circ \sigma_{(l-1,\ldots,k)}^{-2} \circ \sigma_{(1,\ldots,r-1)} \circ \sigma_{(j,\ldots,l-1)}^{-1},$$
so that $\tilde{\sigma}(k-r+1)=j$ and $\tilde{\sigma}(k-r+2)=l$.
Applying this series of transformations to the corresponding term appearing in \eqref{Eq:dPrep-6jl} turns the index sequence $J_<^{j,l}$ into 
$$(u_{\tilde{\sigma}(k-r+2)} v_{\tilde{\sigma}(k-r+3)}, \ldots , 
u_{\tilde{\sigma}(k)} v_1, \ldots , 
u_{\tilde{\sigma}(k-r+1)} v_{\tilde{\sigma}(k-r+2)})$$
which is a cyclic permutation of the sequence $\tilde{\sigma}(u,v)$.  
From \eqref{Eq:RepIndex1}, we can write the corresponding terms in \eqref{Eq:dPrep-6jl} as 
\begin{align} 
&\sum_{r\geq 1} \sum_{\tilde{\sigma} \in U_r}
(-1)^{(k-1)(r-1)} \sgn(\tilde{\sigma})  \nonumber \\ 
&\left( \sigma^{k-r+2} 
\ldb a^{\tilde{\sigma}(k-r+3)},\ldots, a^{\tilde{\sigma}(k)},a^1,\ldots, a^{\tilde{\sigma}(k-r)}, 
\llbracket a^{\tilde{\sigma}(k-r+1)},a^{\tilde{\sigma}(k-r+2)}\rrbracket \rdb_L
\right)_{\tilde{\sigma}(u,v)} \nonumber \\
&\text{ for}\quad U_r:=\{\tilde{\sigma} \in S^{(1)}_k  \mid \tilde{\sigma}(k-r+1) < r < \tilde{\sigma}(k-r+2)\}\,. \label{Eq:dPrep-9b} 
\end{align} 

% Third case for \eqref{Eq:dPrep-6} with  1<j<l
In the third case, consider $\hat{\sigma} \in S_k^{(j,l)}$ for which there exists $r>l$ with $\hat{\sigma}(r)=1$. To understand the corresponding summand of \eqref{Eq:dPrep-6jl}, note that  
\begin{align*} 
&\sgn(\hat{\sigma}) 
\ldb a^{\hat{\sigma}(1)},\stackrel{j}{\check{\dots}},\stackrel{l}{\check{\dots}}, a^{\hat{\sigma}(k)}, \dsq{a^j,a^l}\rdb_L \\
=&(-1)^{k+r+1} \sgn(\dot{\sigma})  
\ldb a^{\dot{\sigma}(k)}, a^{\dot{\sigma}(2)}, 
\stackrel{j}{\check{\dots}},\stackrel{l}{\check{\dots}}, 
a^{\dot{\sigma}(r-1)},a^1,a^{\dot{\sigma}(r)},\ldots, a^{\dot{\sigma}(k-1)}, 
\llbracket a^{j},a^l\rrbracket \rdb_L\,, 
\end{align*} 
where $\hat{\sigma}=\dot{\sigma} \circ \sigma_{(r\cdots k,1)}^{-1}$ with $\dot{\sigma}\in S_k^{(1,j,l)}$. 
This can be equivalently written as 
\begin{align*} 
(-1)^{l+j+k+r} \sgn(\ddot{\sigma})  
\ldb a^{\ddot{\sigma}(k)},a^{\ddot{\sigma}(2)},  \stackrel{j,j+1}{\check{\dots}}, a^{\ddot{\sigma}(r-1)},a^1, a^{\ddot{\sigma}(r)}, \ldots , a^{\ddot{\sigma}(k-1)}, 
\llbracket a^{j},a^{\ddot{\sigma}(j+1)}\rrbracket \rdb_L\,, 
\end{align*} 
after letting $\dot{\sigma}=\ddot{\sigma} \circ \sigma_{(j+1 \cdots l)}$ with  $\ddot{\sigma}\in S_k^{(1,j)}$, then 
\begin{align*} 
(-1&)^{l+j+k+r} \sgn(\dddot{\sigma}) \\ 
&\ldb a^{\ddot{\sigma}(k)},a^{\ddot{\sigma}(2)},  \ldots , a^{\ddot{\sigma}(r-3)},a^1, a^{\ddot{\sigma}(r-2)}, \ldots , a^{\ddot{\sigma}(k-3)}, 
\llbracket a^{\ddot{\sigma}(k-2)},a^{\ddot{\sigma}(k-1)} \rrbracket \rdb_L\,, 
\end{align*} 
after further letting $\ddot{\sigma}=\dddot{\sigma} \circ \sigma_{(j \cdots k-1)}^{-2}$.
A final relabelling yields 
\begin{align*} 
(-1)^{l+j+k}&(-1)^{(k-1)r} \sgn(\tilde{\sigma}) \\
&\ldb a^{\tilde{\sigma}(k-r+4)},\ldots, a^{\tilde{\sigma}(k)},a^1, a^{\tilde{\sigma}(2)},\ldots, a^{\tilde{\sigma}(k-r+1)}, 
\llbracket a^{\tilde{\sigma}(k-r+2)},a^{\tilde{\sigma}(k-r+3)} \rrbracket \rdb_L, 
\end{align*} 
by setting $\dddot{\sigma}=\tilde{\sigma} \circ \sigma_{(2,\ldots,k)}^{k-(r-3)}$. 
In summary, we introduced $\tilde{\sigma}\in S_k^{(1)}$ according to 
$$\hat{\sigma}=\tilde{\sigma} \circ \sigma_{(2,\ldots,k)}^{k-(r-3)} \circ \sigma_{(j,\ldots,k-1)}^{-2} \circ \sigma_{(j+1,\ldots,l)} \circ \sigma_{(r,\ldots,k,1)}^{-1},$$
so that $\tilde{\sigma}(k-r+2)=j$ and $\tilde{\sigma}(k-r+3)=l$. 
Applying this series of transformations to the corresponding term appearing in \eqref{Eq:dPrep-6jl} turns the index sequence $J_<^{j,l}$ into 
$$(u_{\tilde{\sigma}(k-r+3)} v_{\tilde{\sigma}(k-r+4)}, \ldots , 
u_{\tilde{\sigma}(k)} v_1, \ldots , 
u_{\tilde{\sigma}(k-r+2)} v_{\tilde{\sigma}(k-r+3)})$$
which is a cyclic permutation of the sequence $\tilde{\sigma}(u,v)$.  
Thus, using \eqref{Eq:RepIndex1}, we can write the corresponding terms in \eqref{Eq:dPrep-6jl} as 
\begin{align} 
&\sum_{r\geq 1} \sum_{\tilde{\sigma} \in V_r}
(-1)^{(k-1)r} \sgn(\tilde{\sigma})  \nonumber \\ 
&\left( \sigma^{k-r+3} 
\ldb a^{\tilde{\sigma}(k-r+4)},\ldots, a^{\tilde{\sigma}(k)},a^1,\ldots, a^{\tilde{\sigma}(k-r+1)}, 
\llbracket a^{\tilde{\sigma}(k-r+2)},a^{\tilde{\sigma}(k-r+3)} \rrbracket \rdb_L
\right)_{\tilde{\sigma}(u,v)} \nonumber  \\
&\text{ for}\quad V_r:=\{\tilde{\sigma} \in S^{(1)}_k  \mid \tilde{\sigma}(k-r+2) < \tilde{\sigma}(k-r+3) <r\}\,. \label{Eq:dPrep-9c}
\end{align} 
% Gathering the previous calculations 
By setting $i=k-r+1$ in \eqref{Eq:dPrep-9a}, $i=k-r+2$ in \eqref{Eq:dPrep-9b}, and $i=k-r+3$ in \eqref{Eq:dPrep-9c}, respectively, we have obtained 
\begin{align}
&\eqref{Eq:dPrep-6}_{j\neq 1} =  \sum_{1<j<l} \eqref{Eq:dPrep-6jl} 
=\eqref{Eq:dPrep-9a}+\eqref{Eq:dPrep-9b}+\eqref{Eq:dPrep-9c} \nonumber  \\
&=(-1)^{k-1}
\sum_{i=3}^k (-1)^{(k-1)i} \sum_{\tilde{\sigma} \in W_{i,<} }\sgn(\tilde{\sigma}) 
\label{Eq:dPrep-10}
 \\
&\qquad \left(\sigma^i \,\ldb a^{\tilde{\sigma}(i+1)},\ldots,a^{\tilde{\sigma}(k)}, a^{1},a^{\tilde{\sigma}(2)},\ldots,a^{\tilde{\sigma}(i-2)} , \llbracket a^{\tilde{\sigma}(i-1)},a^{\tilde{\sigma}(i)}\rrbracket \rdb_L \right)_{\tilde{\sigma}(u,v)} \nonumber 
\end{align}
where $W_{i,<}:= \{\tilde{\sigma} \in S_k^{(1)} \mid \tilde{\sigma}(i-1)<\tilde{\sigma}(i)\}$. 

\medskip 

%%% second type of terms 
We now proceed with \eqref{Eq:dPrep-7}. 
A term with $1<j<l$ fixed can be written as 
\begin{align}
&\hspace{-0.1cm}\eqref{Eq:dPrep-7}_{j,l}
=
(-1)^{j+l+k} \! \sum_{\hat{\sigma} \in S_k^{(j,l)}}\! \sgn(\hat{\sigma}) 
\left(\! - \ldb a^{\hat{\sigma}(1)},\stackrel{j}{\check{\dots}},\stackrel{l}{\check{\dots}}, a^{\hat{\sigma}(k)}, \llbracket a^l,a^j\rrbracket \rdb_L   
\right)_{J^{j,l}_>},
\label{Eq:dPrep-7jl}  \\
&\hspace{-0.2cm} J^{j,l}_> =
(u_j v_{\hat{\sigma}(1)}, \ldots, 
 u_{\hat{\sigma}(j-1)}v_{\hat{\sigma}(j+1)}, \ldots, 
 u_{\hat{\sigma}(l-1)} v_{\hat{\sigma}(l+1)}, \ldots, 
 u_{\hat{\sigma}(k-1)} v_{\hat{\sigma}(k)} , 
 u_{\hat{\sigma}(k)} v_l , u_l v_j ),  \nonumber 
\end{align}    
after using cyclic skewsymmetry \eqref{Eq:nbr-Cycl} of $\dgal{-}$ and \eqref{Eq:RepIndex1}. 

% First case for \eqref{Eq:dPrep-7} with  1<j<l
In the first case, consider $\hat{\sigma} \in S_k^{(j,l)}$ for which there exists $1\leq r<j$ with $\hat{\sigma}(r)=1$. To understand the corresponding summand of \eqref{Eq:dPrep-7jl}, note that 
\begin{align*} 
&- \sgn(\hat{\sigma}) 
\ldb a^{\hat{\sigma}(1)},\stackrel{j}{\check{\dots}},\stackrel{l}{\check{\dots}}, a^{\hat{\sigma}(k)}, \llbracket a^l,a^j\rrbracket \rdb_L \\
=&(-1)^{r} \sgn(\dot{\sigma})  
\ldb a^{\dot{\sigma}(2)},\ldots, a^{\dot{\sigma}(r)},a^1, a^{\dot{\sigma}(r+1)},\stackrel{j}{\check{\dots}},\stackrel{l}{\check{\dots}}, a^{\dot{\sigma}(k)}, \llbracket a^l,a^j\rrbracket \rdb_L\,, 
\end{align*} 
where $\hat{\sigma}=\dot{\sigma} \circ \sigma_{(1 \cdots r)}$ with $\dot{\sigma}\in S_k^{(1,j,l)}$. 
This can be equivalently written as 
\begin{align*} 
(-1)^{l+j+r} \sgn(\ddot{\sigma})  
\ldb  a^{\ddot{\sigma}(2)},\ldots, a^{\ddot{\sigma}(r)},a^1, a^{\ddot{\sigma}(r+1)},\stackrel{l-1,l}{\check{\dots}}, a^{\ddot{\sigma}(k)}, 
\llbracket a^{\ddot{\sigma}(l-1)},a^{\ddot{\sigma}(l)}\rrbracket \rdb_L\,, 
\end{align*} 
after letting $\dot{\sigma}=\ddot{\sigma} \circ \sigma_{(j \cdots l)}^{-1}$, and then 
\begin{align*} 
(-1)^{l+j+r} \sgn(\dddot{\sigma})  
\ldb a^{\dddot{\sigma}(2)},\ldots, a^{\dddot{\sigma}(r)},a^1, a^{\dddot{\sigma}(r+1)},\ldots, a^{\dddot{\sigma}(k-2)}, 
\llbracket a^{\ddot{\sigma}(k-1)},a^{\dddot{\sigma}(k)} \rrbracket \rdb_L
\end{align*} 
after further letting $\ddot{\sigma}=\dddot{\sigma} \circ \sigma_{(l-1 \cdots k)}^{-2}$.
A final relabeling yields 
\begin{align*} 
(-1&)^{l+j+k}(-1)^{(k-1)r} \sgn(\tilde{\sigma})  \\
&\ldb a^{\tilde{\sigma}(k-r+2)},\ldots, a^{\tilde{\sigma}(k)},a^1, a^{\tilde{\sigma}(2)},\ldots, a^{\tilde{\sigma}(k-r-1)}, 
\llbracket a^{\tilde{\sigma}(k-r)},a^{\tilde{\sigma}(k-r+1)} \rrbracket\rdb_L, 
\end{align*} 
by setting $\dddot{\sigma}=\tilde{\sigma} \circ \sigma_{(2,\ldots,k)}^{k-r}$. 
In summary, we introduced $\tilde{\sigma}\in S_k^{(1)}$ according to 
$$\hat{\sigma}=\tilde{\sigma} \circ \sigma_{(2,\ldots,k)}^{k-r} \circ \sigma_{(l-1,\ldots,k)}^{-2} \circ \sigma_{(j,\ldots,l)}^{-1} \circ \sigma_{(1,\ldots,r)} ,$$
so that $\tilde{\sigma}(k-r)=l$ and $\tilde{\sigma}(k-r+1)=j$.
Applying this series of transformations to the corresponding term appearing in \eqref{Eq:dPrep-7jl} turns the index sequence $J_>^{j,l}$ into 
$$(u_{\tilde{\sigma}(k-r+1)} v_{\tilde{\sigma}(k-r+2)}, \ldots , 
u_{\tilde{\sigma}(k)} v_1, \ldots , 
u_{\tilde{\sigma}(k-r)} v_{\tilde{\sigma}(k-r+1)})$$
which is a cyclic permutation of the sequence $\tilde{\sigma}(u,v)$.  
Thus, using \eqref{Eq:RepIndex1}, we can write the corresponding terms in \eqref{Eq:dPrep-7jl} as 
\begin{align} 
&\sum_{r\geq 1} \sum_{\tilde{\sigma} \in T_r'}
(-1)^{(k-1)r} \sgn(\tilde{\sigma})  \nonumber \\ 
&\left( \sigma^{k-r+1} 
\ldb a^{\tilde{\sigma}(k-r+2)},\ldots, a^{\tilde{\sigma}(k)},a^1,\ldots, a^{\tilde{\sigma}(k-r-1)}, 
\llbracket a^{\tilde{\sigma}(k-r)},a^{\tilde{\sigma}(k-r+1)} \rrbracket \rdb_L
\right)_{\tilde{\sigma}(u,v)} , \nonumber  \\
&\text{ for}\quad T_r':=\{\tilde{\sigma} \in S^{(1)}_k  \mid r< \tilde{\sigma}(k-r+1) < \tilde{\sigma}(k-r) \}\,. \label{Eq:dPrep-11a}
\end{align} 

% Second case for \eqref{Eq:dPrep-7} with  1<j<l
In the second case, consider $\hat{\sigma} \in S_k^{(j,l)}$ for which there exists $j<r<l$ with $\hat{\sigma}(r)=1$. To understand the corresponding summand of \eqref{Eq:dPrep-7jl}, remark that 
\begin{align*} 
&- \sgn(\hat{\sigma}) 
\ldb a^{\hat{\sigma}(1)},\stackrel{j}{\check{\dots}},\stackrel{l}{\check{\dots}}, a^{\hat{\sigma}(k)}, \dsq{a^l,a^j} \rdb_L \\
=&(-1)^{j+l+1} \sgn(\dot{\sigma})  
\ldb a^{\dot{\sigma}(1)},\stackrel{l-1,l}{\check{\dots}}, a^{\dot{\sigma}(k)}, 
\llbracket a^{\dot{\sigma}(l-1)},a^{\dot{\sigma}(l)}\rrbracket \rdb_L\,, 
\end{align*} 
where $\hat{\sigma}=\dot{\sigma} \circ \sigma_{(j \cdots l)}^{-1}$ with $\dot{\sigma}$ such that $\dot{\sigma}(r-1)=1$.  
This is equivalent to  
\begin{align*} 
(-1)^{l+j+r+1} \sgn(\ddot{\sigma})  
\ldb  a^{\ddot{\sigma}(2)},\ldots, a^{\ddot{\sigma}(r-1)},a^1, a^{\ddot{\sigma}(r)},\stackrel{l-1,l}{\check{\dots}}, a^{\ddot{\sigma}(k)}, 
\llbracket a^{\ddot{\sigma}(l-1)},a^{\ddot{\sigma}(l)} \rrbracket \rdb_L , 
\end{align*} 
after letting $\dot{\sigma}=\ddot{\sigma} \circ \sigma_{(1 \cdots r-1)}$, and then 
\begin{align*} 
(-1&)^{l+j+r+1} \sgn(\dddot{\sigma})  \\
&\ldb a^{\dddot{\sigma}(2)},\ldots, a^{\dddot{\sigma}(r-1)},a^1, a^{\dddot{\sigma}(r)},\ldots, a^{\dddot{\sigma}(k-2)}, 
\llbracket a^{\ddot{\sigma}(k-1)},a^{\dddot{\sigma}(k)} \rrbracket \rdb_L\,, 
\end{align*} 
after further letting $\ddot{\sigma}=\dddot{\sigma} \circ \sigma_{(l-1 \cdots k)}^{-2}$.
Finally, this is 
\begin{align*} 
(-1&)^{l+j+r+1}(-1)^{kr} \sgn(\tilde{\sigma})  \\
&\ldb a^{\tilde{\sigma}(k-r+3)},\ldots, a^{\tilde{\sigma}(k)},a^1, a^{\tilde{\sigma}(2)},\ldots, a^{\tilde{\sigma}(k-r)}, 
\llbracket a^{\tilde{\sigma}(k-r+1)},a^{\tilde{\sigma}(k-r+2)}\rrbracket \rdb_L, 
\end{align*} 
by setting $\dddot{\sigma}=\tilde{\sigma} \circ \sigma_{(2,\ldots,k)}^{k-(r-1)}$. 
In summary, we introduced $\tilde{\sigma}\in S_k^{(1)}$ according to 
$$\hat{\sigma}=\tilde{\sigma} \circ \sigma_{(2,\ldots,k)}^{k-(r-1)} \circ \sigma_{(l-1,\ldots,k)}^{-2} \circ \sigma_{(1,\ldots,r-1)} \circ \sigma_{(j,\ldots,l)}^{-1} ,$$
so that $\tilde{\sigma}(k-r+1)=l$ and $\tilde{\sigma}(k-r+2)=j$. 
Applying this series of transformations to the corresponding term appearing in \eqref{Eq:dPrep-7jl} turns the index sequence $J_>^{j,l}$ into 
$$(u_{\tilde{\sigma}(k-r+2)} v_{\tilde{\sigma}(k-r+3)}, \ldots , 
u_{\tilde{\sigma}(k)} v_1, \ldots , 
u_{\tilde{\sigma}(k-r+1)} v_{\tilde{\sigma}(k-r+2)})$$
which is a cyclic permutation of $\tilde{\sigma}(u,v)$.   
Thus, using \eqref{Eq:RepIndex1}, we can write the corresponding terms in \eqref{Eq:dPrep-7jl} as 
\begin{align} 
&\sum_{r\geq 1} \sum_{\tilde{\sigma} \in U_r'}
(-1)^{(k-1)(r-1)} \sgn(\tilde{\sigma})  \nonumber \\ 
&\left( \sigma^{k-r+2} 
\ldb a^{\tilde{\sigma}(k-r+3)},\ldots, a^{\tilde{\sigma}(k)},a^1,\ldots, a^{\tilde{\sigma}(k-r)}, 
\llbracket a^{\tilde{\sigma}(k-r+1)},a^{\tilde{\sigma}(k-r+2)} \rrbracket\rdb_L
\right)_{\tilde{\sigma}(u,v)}, \nonumber  \\
&\text{ for}\quad U_r':=\{\tilde{\sigma} \in S^{(1)}_k  \mid \tilde{\sigma}(k-r+2) <r < \tilde{\sigma}(k-r+1) \}\,. \label{Eq:dPrep-11b}
\end{align} 

% Third case for \eqref{Eq:dPrep-7} with  1<j<l
In the third and final case, consider $\hat{\sigma} \in S_k^{(j,l)}$ for which there exists $r>l$ with $\hat{\sigma}(r)=1$. To understand the corresponding summand of \eqref{Eq:dPrep-7jl}, remark that 
\begin{align*} 
&- \sgn(\hat{\sigma}) 
\ldb a^{\hat{\sigma}(1)},\stackrel{j}{\check{\dots}},\stackrel{l}{\check{\dots}}, a^{\hat{\sigma}(k)}, \dsq{a^l,a^j} \rdb_L \\
=&(-1)^{k+r} \sgn(\dot{\sigma})  
\ldb a^{\dot{\sigma}(k)},a^{\dot{\sigma}(2)},
\stackrel{j}{\check{\dots}},\stackrel{l}{\check{\dots}},
a^{\dot{\sigma}(r-1)},a^1,a^{\dot{\sigma}(r)},\ldots,a^{\dot{\sigma}(k-1)}, \dsq{a^l,a^j} \rdb_L\,, 
\end{align*} 
where $\hat{\sigma}=\dot{\sigma} \circ \sigma_{(r\cdots k,1)}^{-1}$ with $\dot{\sigma}\in S_k^{(1,j,l)}$. 
This is equivalent to  
\begin{align*} 
(-1&)^{l+j+k+r} \sgn(\ddot{\sigma})  \\
&\ldb  a^{\ddot{\sigma}(k)},a^{\ddot{\sigma}(2)},
\stackrel{j,j+1}{\check{\dots}},
a^{\ddot{\sigma}(r-1)},a^1,a^{\ddot{\sigma}(r)},\ldots,a^{\ddot{\sigma}(k-1)},  
\llbracket a^{\ddot{\sigma}(j)},a^{\ddot{\sigma}(j+1)} \rrbracket \rdb_L\,, 
\end{align*} 
after letting $\dot{\sigma}=\ddot{\sigma} \circ \sigma_{(j \cdots l)}$, and then 
\begin{align*} 
(-&1)^{l+j+k+r} \sgn(\dddot{\sigma})  \\
&\ldb a^{\dddot{\sigma}(k)}, a^{\dddot{\sigma}(2)},\ldots, a^{\dddot{\sigma}(r-3)},a^1, a^{\dddot{\sigma}(r-2)},\ldots, a^{\dddot{\sigma}(k-3)}, 
\llbracket a^{\ddot{\sigma}(k-2)},a^{\dddot{\sigma}(k-1)} \rrbracket \rdb_L
\end{align*} 
after further letting $\ddot{\sigma}=\dddot{\sigma} \circ \sigma_{(j \cdots k-1)}^{-2}$.
Finally, this is 
\begin{align*} 
(-1&)^{l+j+k}(-1)^{(k-1)r} \sgn(\tilde{\sigma})  \\
&\ldb a^{\tilde{\sigma}(k-r+4)},\ldots, a^{\tilde{\sigma}(k)},a^1,\ldots, a^{\tilde{\sigma}(k-r+1)}, 
\llbracket a^{\tilde{\sigma}(k-r+2)},a^{\tilde{\sigma}(k-r+3)} \rrbracket \rdb_L, 
\end{align*} 
by setting $\dddot{\sigma}=\tilde{\sigma} \circ \sigma_{(2,\ldots,k)}^{k-(r-3)}$. 
In summary, we introduced $\tilde{\sigma}\in S_k^{(1)}$ according to 
$$\hat{\sigma}=\tilde{\sigma} \circ \sigma_{(2,\ldots,k)}^{k-(r-3)} \circ \sigma_{(j,\ldots,k-1)}^{-2} \circ \sigma_{(j,\ldots,l)} \circ \sigma_{(r,\ldots,k,1)}^{-1},$$
so that $\tilde{\sigma}(k-r+2)=l$ and $\tilde{\sigma}(k-r+3)=j$. 
Applying this series of transformations to the corresponding term appearing in \eqref{Eq:dPrep-7jl} turns the index sequence $J_>^{j,l}$ into 
$$(u_{\tilde{\sigma}(k-r+3)} v_{\tilde{\sigma}(k-r+4)}, \ldots , 
u_{\tilde{\sigma}(k)} v_1, \ldots , 
u_{\tilde{\sigma}(k-r+2)} v_{\tilde{\sigma}(k-r+3)})$$
which is a cyclic permutation of the sequence $\tilde{\sigma}(u,v)$.  
Thus, using \eqref{Eq:RepIndex1}, we can write the corresponding terms in \eqref{Eq:dPrep-7jl} as 
\begin{align} 
&\sum_{r\geq 1} \sum_{\tilde{\sigma} \in V_r'}
(-1)^{(k-1)r} \sgn(\tilde{\sigma})  \nonumber \\ 
&\hspace{-0.1cm}\left( \sigma^{k-r+3} 
\ldb a^{\tilde{\sigma}(k-r+4)},\ldots, a^{\tilde{\sigma}(k)},a^1,\ldots, a^{\tilde{\sigma}(k-r+1)}, 
\llbracket a^{\tilde{\sigma}(k-r+2)},a^{\tilde{\sigma}(k-r+3)} \rrbracket\rdb_L
\right)_{\tilde{\sigma}(u,v)} \!, \nonumber \\
&\text{ for}\quad V_r':=\{\tilde{\sigma} \in S^{(1)}_k  \mid \tilde{\sigma}(k-r+3) < \tilde{\sigma}(k-r+2) < r \}\,. \label{Eq:dPrep-11c}
\end{align} 
By setting $i=k-r+1$ in \eqref{Eq:dPrep-11a}, $i=k-r+2$ in \eqref{Eq:dPrep-11b}, and $i=k-r+3$ in \eqref{Eq:dPrep-11c}, respectively, we have obtained 
\begin{align}
&\eqref{Eq:dPrep-7}_{j\neq 1} =  \sum_{1<j<l} \eqref{Eq:dPrep-7jl} 
=\eqref{Eq:dPrep-11a}+\eqref{Eq:dPrep-11b}+\eqref{Eq:dPrep-11c} \nonumber  \\
&=(-1)^{k-1}
\sum_{i=3}^k (-1)^{(k-1)i} \sum_{\tilde{\sigma} \in W_{i,>} }\sgn(\tilde{\sigma}) 
\label{Eq:dPrep-12}
 \\
&\qquad \left(\sigma^i \,\ldb a^{\tilde{\sigma}(i+1)},\ldots,a^{\tilde{\sigma}(k)}, a^{1},a^{\tilde{\sigma}(2)},\ldots,a^{\tilde{\sigma}(i-2)} , 
\llbracket a^{\tilde{\sigma}(i-1)},a^{\tilde{\sigma}(i)} \rrbracket \rdb_L 
\right)_{\tilde{\sigma}(u,v)} \nonumber 
\end{align}
where $W_{i,>}:= \{\tilde{\sigma} \in S_k^{(1)} \mid \tilde{\sigma}(i-1)>\tilde{\sigma}(i)\}$. 
Noting that $S_k^{(1)}=W_{i,<}\cup W_{i,>}$, 
\begin{align}
  \eqref{Eq:dPrep-10} +& \eqref{Eq:dPrep-12}  
=(-1)^{k-1}
\sum_{i=3}^k (-1)^{(k-1)i} \sum_{\tilde{\sigma} \in S_k^{(1)} }\sgn(\tilde{\sigma}) 
\label{Eq:dPrep-13} \\
&\quad \left(\sigma^i \,\ldb a^{\tilde{\sigma}(i+1)},\ldots,a^{\tilde{\sigma}(k)}, a^{1},a^{\tilde{\sigma}(2)},\ldots,a^{\tilde{\sigma}(i-2)} , \llbracket a^{\tilde{\sigma}(i-1)},a^{\tilde{\sigma}(i)}\rrbracket \rdb_L \right)_{\tilde{\sigma}(u,v)} . \nonumber   
\end{align}

Up to the factor $(-1)^{k-1}$, one has: 
\begin{itemize}
    \item \eqref{Eq:dPrep-8a} is the summand $i=2$ of $\mathtt{T}_1$ \eqref{Eq:dPrep-1A} in \eqref{Eq:dPrep-1tr}; 
    \item \eqref{Eq:dPrep-8b} is the summand $i=1$  of $\mathtt{T}_1$ \eqref{Eq:dPrep-1A} in \eqref{Eq:dPrep-1tr};  
    \item \eqref{Eq:dPrep-13} is all the remaining summands $3\leq i \leq k$ of $\mathtt{T}_1$ \eqref{Eq:dPrep-1A} in \eqref{Eq:dPrep-1tr}. 
\end{itemize}
This establishes the claimed equality \eqref{Eq:dPrep-Goal2}. 
Combining \eqref{Eq:dPrep-Goal1} and \eqref{Eq:dPrep-Goal2} entails 
$\delta^{k-1}(\tr \dgal{-}) = (-1)^{k-1} \tr(\wdd(\dgal{-}))$, which is the commutativity of the diagram presented at the beginning of the proof. 
\qed

\section[From gauged double Poisson \& double quasi-Poisson]{From gauged double Poisson and double quasi-Poisson cohomologies}
\label{Sec:InducGauge}

\subsection{The gauged double Poisson case}

Recall from Proposition \ref{Pr:Tr-Lie} (resp.  Theorem \ref{Thm:IndBr})
that the trace map \eqref{Eq:Tr-morph} (resp. \eqref{Eq:Tr-morph-2}) is a linear map from the complex defining the (resp. completed) double Poisson cohomology to $\mc X(\cA_{\bf n}^{\Gl_{\bf n}})$. 
The main results of Section \ref{sec:dPH-PH} can then be adapted to the gauged case because the trace map can be restricted to $\mc D_{\cA}$ and $\widehat{\mc D}_{\cA}$, see Corollary \ref{Cor:Tr-morphD} and \eqref{Eq:Tr-morph-3}.

\begin{theorem}  \label{Thm:gdP-rep1}
If $P\in (\mb T^\ast \cA)_{\sharp,2}$ satisfies $\brSN{P,P}=0$,
the differential $\dd_P$ of \eqref{Eq:dP-AKKN} and the differential $\dd_{\tr(P)}$ induce the morphism of complexes 
\begin{equation} \label{Eq:gdP-rep1}
    \tr : ({\mc D}_\cA, \dd_P) \longrightarrow 
    (\mf X(\cA_{\bf n}^{\Gl_{\bf n}}), \dd_{\tr(P)}),  
\end{equation}
which descends to a linear map $\gdPH(\cA) \to \PH(\cA_{\bf n}^{\Gl_{\bf n}})$ in cohomology. 
Similarly, we get a linear map $\widehat{\gdPH}(\cA) \to {\mathrm H}_{CE}(\cA_{\bf n}^{\Gl_{\bf n}}) \stackrel{\sim}{\longrightarrow} \PH(\cA_{\bf n}^{\Gl_{\bf n}})$. 
\end{theorem}
\begin{proof}
The first part of the statement over $B=\kk$ is \cite[Rem.~4.7]{AKKN}. 
We get the morphism of  complexes \eqref{Eq:gdP-rep1} from the first part of Theorem \ref{Thm:dP-rep1} which we induce to ${\mc D}_\cA$ thanks to \eqref{Eq:TrBrRep-2}-\eqref{Eq:Tr-morph-3}.  The completed case uses Theorem \ref{Thm:dP-rep2}.   
\end{proof}

\begin{corollary} \label{Cor:gdPH-PH} 
Assume that $P\in (\mb T^\ast \cA)_{2,\sharp}$ satisfies $\brSN{P,P}=0$. 
Then, the following diagram is commutative: 
%\pecetta{this diagram has been commented to speed up compiling}
%\begin{comment}
\begin{center}
      \begin{tikzpicture}
%%%% BACK
 \node   (T0) at (0,2.5) {$\dPH(\cA)$};
 \node   (T1) at (4,2.5) {$\widehat{\dPH}(\cA)$};
 \node  (B0) at (0,-0.5) {$\PH(\cA_{\bf n})$};
 \node  (B1) at (4,-0.5) {${\mathrm H}_{CE}(\cA_{\bf n})$};
\path[->,>=angle 90,font=\small]  
   (T0) edge (T1) ;
\path[->,>=angle 90,font=\small,dashed]  
   (B0) edge  (B1) ;
\path[->,>=angle 90,font=\small,dashed]  
   (T0) edge (B0) ;
\path[->,>=angle 90,font=\small]  
   (T1) edge (B1) ;
%%% FRONT 
 \node   (T2) at (-2.5,1) {$\gdPH(\cA)$};
 \node   (T3) at (1.5,1) {$\widehat{\gdPH}(\cA)$};
 \node  (B2) at (-2.5,-2) {$\PH(\cA_{\bf n}^{\Gl_{\bf n}})$};
 \node  (B3) at (1.5,-2) {${\mathrm H}_{CE}(\cA_{\bf n}^{\Gl_{\bf n}})$};
\path[->,>=angle 90,font=\small]  
   (T2) edge (T3) ;
\path[->,>=angle 90,font=\small]  
   (B2) edge (B3) ;
\path[->,>=angle 90,font=\small]  
   (T2) edge (B2) ;
\path[->,>=angle 90,font=\small]  
   (T3) edge (B3) ;
%%% other maps 
\path[->,>=angle 90,font=\small]  
   (T0) edge (T2) ;
   \path[->,>=angle 90,font=\small]  
   (T1) edge (T3) ;
\path[->,>=angle 90,font=\small,dashed]  
   (B0) edge (B2) ;
\path[->,>=angle 90,font=\small]  
   (B1) edge (B3) ;
   \end{tikzpicture}
\end{center}
%\end{comment}
where the differentials and the linear maps on the back face are defined as in Corollary \ref{Cor:dPH-PH}, from which they are naturally induced on the front face.   
\end{corollary}
\begin{proof}
Commutativity of the top face comes from Corollary \ref{Cor:dPH-gdPH}. 
For the bottom and back faces, this is Corollary \ref{Cor:dPH-PH}. 
Commutativity of the left, right and front faces follows by construction, cf.  the proof of Theorem \ref{Thm:gdP-rep1}.     
\end{proof}

\begin{example} \label{Exmp:Rep-gdP-kx}
For $\cA=\kk[x]$ with the double Poisson bracket \eqref{Eq:dbr-xlin} (where $\lambda,\mu,\nu$ satisfy $\lambda\nu-\mu^2=0$ and are not all zero), the maps $(\mu_\ell)_{\ell\geq 0}$ of Proposition \ref{Pr:MapMu} are isomorphisms so all horizontal arrows in the cube of Corollary \ref{Cor:gdPH-PH} are isomorphisms. Thus only the left face of that cube requires attention. 
For $N\geq 1$, $\cA_N\simeq \kk[\gl_N]$ is generated by the $N^2$ elements $x_{ij}$, $1\leq i,j\leq N$, and the Poisson bracket induced by Theorem \ref{Thm:Rep-Dbr} satisfies 
\begin{equation}
 \begin{aligned} \label{Eq:PBkx-ind}
    \br{x_{ij},x_{kl}} =&
    \lambda (x_{kj} \delta_{il} - \delta_{kj} x_{il})  
+  \mu  \sum_{1\leq r \leq N} (x_{kr}x_{rj} \delta_{il} - \delta_{kj} x_{ir}x_{rl})  \\
&+ \nu  \sum_{1\leq r \leq N} (x_{kr}x_{rj} x_{il} - x_{kj} x_{ir}x_{rl})  \,.     
 \end{aligned}
\end{equation}
(Here, we use the Kronecker delta notation.)
We easily deduce that this descends to the zero Poisson bracket on $\cA_N^{\Gl_N}$, 
so that $\PH(\cA_N^{\Gl_N})\simeq {\mc X}(\cA_N^{\Gl_N})$. 
This is consistent with Proposition \ref{Pr:gdPH-kx0} where we deduce that $\gdPH(\kk[x])$ is given by the full complex ${\mc D}_{\kk[x]}$.  
(For completeness, note that $\Spec(\cA_N^{\Gl_N})$ is the $N$-dimensional affine space $\mb A^N$.)

The first 3 cohomology groups of $\dPH(\kk[x])$ (see Subsection~\ref{ss:dP-coh-x}) are 
\begin{equation}
\begin{aligned}
 \dPH^0(\kk[x])&=\kk[x], \quad  \quad    \dPH^2(\kk[x])= \{0\},   \\ 
\dPH^1(\kk[x])&=\kk\, (\nu x^2 \del_x + 2\mu x \del_x + \lambda \del_x)_\sharp,
\end{aligned}
\end{equation}
and they induce the following nontrivial classes 
\begin{equation}
    \sum_{k \geq 0}\kk \tr(x^k) \subset \PH^0(\cA_N), \quad 
    \kk \tr(\nu x^2 \del_x + 2\mu x \del_x + \lambda \del_x) \subset \PH^1(\cA_N), 
\end{equation}
where the derivation appearing in $\PH^1(\cA_N)$ is given by (cf. Section~\ref{ss:Rep-Not})
\begin{equation*}
 \tr(\nu x^2 \del_x + 2\mu x \del_x + \lambda \del_x) : x_{ij} 
 \mapsto \nu \sum_{1\leq r \leq N} x_{ir}x_{rj} + 2 \mu x_{ij} + \lambda \delta_{ij}\,.
\end{equation*}
It is clear that all the classes induced in $\PH^0(\cA_N)$ are not independent due to the Cayley-Hamilton theorem. 
\end{example}

\begin{remark}
For $\mu=\nu=0$ and $\lambda=1$,  \eqref{Eq:PBkx-ind} gives the Lie-Poisson structure on $\gl_N$ and we recover from this example elements in its first $3$ cohomology classes. 
\end{remark}

%\pecetta{Is there a ref where these groups appear?
%
%For $\mu=\nu=0$, these are invariants in the cohomology 
%of $\gl_N$ with its Lie-Poisson structure}

\begin{example} 
We can adapt the previous example to $\cA:=\kk[x]/(x^r)$, $r\geq 3$, based on Proposition \ref{Pr:dPH-kxTrun}. In that case, $\cA_N$ is generated by the $N^2$ elements $x_{ij}$ subject to the vanishing of the $N^2$ functions $(x^r)_{ij}=0$, $1\leq i,j\leq N$. 
Geometrically, 
\begin{equation*}
   \Spec(\cA_N)\simeq \{Y \in \gl_N \mid Y^r = 0_N\}\,,
\end{equation*}
and for $N=1$, this is simply the fat point $\Spec(\kk[x]/(x^r))$. 
The induced Poisson bracket is still defined through \eqref{Eq:PBkx-ind}. 
The completed double Poisson cohomology induces
cohomology classes in ${\mathrm H}_{CE}(\cA_N)$ through Theorem \ref{Thm:dP-rep2}. 
For the group $\widehat{\dPH}^0(\cA)=\cA$, we obtain 
$\sum_{k \geq 0}\kk \tr(x^k)\subset {\mathrm H}_{CE}^0(\cA_N)$. 
All these classes are in the zero class. 
Indeed, $\cA_N$ is the coordinate ring of the $r$-step nilpotent matrices, whose eigenvalues are all zero. 
The only other nonzero classes appearing in Proposition \ref{Pr:dPH-kxTrun} occur in the case of the Poisson bracket \eqref{Eq:PBkx-ind} with $\nu\neq 0$ and $\lambda=\mu=0$; then 
$\widehat{\dPH}^1(\cA)=\kk\,\theta_2$ and $\widehat{\dPH}^2(\cA)=\kk  \dgal{-,-}_{2,0}$ yield 
\begin{align*}
 \kk \tr(\theta_2)\subset {\mathrm H}_{CE}^1(\cA_N), \quad 
 \tr(\theta_2):x_{ij}\mapsto \sum_{1\leq r \leq N} x_{ir}x_{rj}\,,  
\end{align*}
and $\kk \tr(\dgal{-,-}_{2,0})\subset {\mathrm H}_{CE}^2(\cA_N)$, where
\begin{align*}
 \tr(\dgal{-,-}_{2,0}):(x_{ij},x_{pq})\mapsto 
 \sum_{1\leq r \leq N} (x_{pr}x_{rj} \delta_{iq} - \delta_{pj} x_{ir}x_{rq}) \,.
\end{align*}
%and the class of the Poisson bracket \eqref{Eq:PBkx-ind} in ${\mathrm H}_{CE}^2(\cA_{\bf n})$. 
\end{example}

\subsection{The quasi-Poisson case}

Let us start by recalling the relation between a noncommutative quasi-Poisson bivector, in the sense of Definition \ref{Def:qPoiss}, and the geometric notion of a quasi-Poisson bivector from Section~\ref{ss:classQPCoh}. 

For any $s\in S$, we let $G_{(s)}$ denote the copy of $\Gl_{n_s}(\kk)$ in $s$-th position in the product $\Gl_{\bf n}:=\prod_{i\in I}\Gl_{n_i}(\kk)$. 
Under the embedding $G_{(s)}\hookrightarrow \Gl_{\bf n}$, 
we get an infinitesimal action of $\g_{(s)}=\operatorname{Lie}(G_{(s)})$, cf. \eqref{Eq:InfAct}, and a corresponding $3$-derivation $\phi^{(s)}_{\cA_{\bf n}}$ induced by the Cartan $3$-tensor $\phi^{(s)}\in \bigwedge^3 \g_s$. 

\begin{proposition}[\cite{VdB1}, 7.12.1]  \label{Eq:PhiRep}
We have: $\phi^{(s)}_{\cA_{\bf n}} = \frac16 \tr(\Delta_s^3)$ for any $s\in S$.     
\end{proposition} 
\begin{proof}
 This is a simple computation using that the derivation $(\Delta_s)_{uv}$ is the infinitesimal action of $E^{(s)}_{vu}$, cf. Remark \ref{Rem:DeltaRep}.  
\end{proof}

Combining this result with the definition relation \eqref{Eq:qP-PP} and Proposition \ref{Pr:Tr-Lie}, we directly get the following. 
\begin{proposition}[\cite{VdB1}, 7.12.2] \label{Pr:IndqP}
If $P\in (\mb T^\ast \cA)_{\sharp,2}$ defines a double quasi-Poisson bracket, then 
$\tr(P)$ is quasi-Poisson on $\cA_{\bf n}$ 
(in the sense of Remark \ref{Rem:qPdef}). 
\end{proposition}

The next result is the quasi-Poisson analogue of Theorem \ref{Thm:dP-rep1} due to \cite{PV,VdW}, and it can be summarized in terms of the commutative diagram presented therein.  
\begin{theorem} \label{Thm:dqP-rep1}
Assume that $P\in (\mb T^\ast \cA)_{\sharp,2}$ is quasi-Poisson. 
Then, the differential $\dd_P=\brSN{P,-}$ of Proposition \ref{Pr:dquasiPcoh1} and 
the differential $\dd_{\tr(P)}=[\tr(P),-]_{\SN}$ induce the morphism of complexes 
\begin{equation} \label{Eq:dqP-rep1}
    \tr : ((\mb T^\ast \cA)_{\sharp}, \dd_P) \longrightarrow 
    (\mf X(\cA_{\bf n})^{\Gl_{\bf n}}, \dd_{\tr(P)}),  
\end{equation}
which descends to a linear map $\dPH(\cA) \to \PH_{\Gl_{\bf n}}(\cA_{\bf n})$ in cohomology. 
Furthermore, this morphism can be restricted to a map 
$\dPH(\cA) \to \PH(\cA_{\bf n}^{\Gl_{\bf n}})$. 
\end{theorem}
\begin{proof}   
We know that $P$ induces the quasi-Poisson bivector $\tr(P)$ on representation spaces by Proposition \ref{Pr:IndqP}, so the map of complexes is obtained thanks to the graded Lie algebra homomorphism of Proposition \ref{Pr:Tr-Lie} valued in $\mf X(\cA_{\bf n})^{\Gl_{\bf n}}$. 
We thus get the first map. 

For the second map, it suffices to restrict invariant multiderivations to the invariant ring $\cA_{\bf n}^{\Gl_{\bf n}}$, and then consider the corresponding cohomology theories. 
\end{proof}

Let us now assume that $\dsq{-,-}$ is a ($B$-linear) double quasi-Poisson bracket in the sense that it satisfies \eqref{qPabc}. (We shall comment below on the case where it is defined by a quasi-Poisson $P\in (\mb T^\ast \cA)_{\sharp,2}$, cf. Remark \ref{Rem:Comp-qCoh}.)

\begin{proposition}
The skewsymmetric biderivation $\tr(\dsq{-,-})$ defined on $\cA_{\bf n}$ through Theorem \ref{Thm:IndBr} is a quasi-Poisson bracket.
\end{proposition}
\begin{proof}
It suffices to adapt \cite[Thm.~5.6]{F21} (based on \cite[Rem.~7.12.3]{VdB1}). 
We put $\br{-,-}:=\tr(\dsq{-,-})$.
We have on generators $a_{ij}, b_{kl}, c_{uv}\in \cA_{\bf n}$
\begin{equation}
\begin{aligned}
     & \br{ a_{ij}, \br{b_{kl}, c_{uv}} }  +  \br{b_{kl}, \br{ c_{uv},a_{ij}} }  + 
  \br{c_{uv}, \br{a_{ij},b_{kl}} } \\
=&\dgal{a,b,c}_{uj,il,kv}- \dgal{b,a,c}_{kj,iv,ul} \\
=&\sum_s\frac{1}{12} q_s \left((\dgal{a,b,c}_{\Delta_s^3})_{uj,il,kv} 
- (\dgal{b,a,c}_{\Delta_s^3})_{kj,iv,ul}\right)\\ 
=&\sum_s\frac{1}{12} q_s \tr(\dgal{-,-,-}_{\Delta_s^3})  (a_{ij}, b_{kl}, c_{uv} ) \\
=&\frac12 \sum_s q_s \phi^{(s)}_{\cA_{\bf n}}(a_{ij}, b_{kl}, c_{uv}) 
\end{aligned}
\end{equation}
where the first equality is \cite[Eq.~(5.3)]{F21}, the second follows from Remark~\ref{Rem:qP} and \eqref{qPabc}, the third is obtained by \eqref{Eq:TrBrRep} and the fourth is a consequence of \eqref{Eq:PhiRep}. 
Since this is an equality of multiderivations, the same is true if we replace $a_{ij}, b_{kl}, c_{uv}$ by arbitrary functions $f_1,f_2,f_3\in \cA_{\bf n}$. 
In view of \eqref{Eq:JacPhi}, $\br{-,-}$ is a quasi-Poisson bracket for $\phi:=\sum_s q_s \phi^{(s)}$. 
\end{proof}

\begin{theorem} \label{Thm:dqP-rep2} 
Assume that  $\dsq{-,-} \in \wBRA_B(\cA)_2$ is a double quasi-Poisson bracket, and let $\br{-,-}$ denote the associated quasi-Poisson bracket on $\cA_{\bf n}$ obtained through Theorem \ref{Thm:Rep-Dbr}. 
The differential $\wdd$ \eqref{Eq:dP-gen-0}-\eqref{Eq:dP-gen} (cf. Proposition \ref{Pr:dquasiPcoh2}) and the differential $\delta_{\cA_{\bf n},\br{-,-}}$ \eqref{Eq:Diff-Pcoh1} (cf. Theorem \ref{Thm:qPcohGen}) induce the morphism of complexes 
\begin{equation} \label{Eq:dqP-rep2}
    \tr : (\wBRA_B(\cA), \wdd) \longrightarrow 
    (\mf X(\cA_{\bf n})^{\Gl_{\bf n}}, (-1)^{\bullet} \, \delta_{\cA_{\bf n},\br{-,-}} ),  
\end{equation}
which descends to a linear map $\widehat{\dPH}(\cA) \to {\mathrm H}_{CE;\Gl_{\bf n}}(\cA_{\bf n})$ in cohomology. 
Furthermore, the above morphism of complexes factors can be restricted to $\mf X(\cA_{\bf n}^{\Gl_{\bf n}})$. 
In particular, this descends to a linear map $\widehat{\dPH}(\cA) \to {\mathrm H}_{CE}(\cA_{\bf n}^{\Gl_{\bf n}})$ in cohomology.  
\end{theorem} 

\begin{proof}
We already proved in \S\ref{ss:Proof-dP-rep2} that the map \eqref{Eq:dP-rep2} is a morphism of complexes because it did not require to know that $\dsq{-,-}$ or $\br{-,-}$ was Poisson. The rest directly follows by definition. 
\end{proof}

\begin{remark} \label{Rem:Comp-qCoh}
If the double quasi-Poisson bracket is of the form $\mu_2(P)$, $P\in (\mb T^\ast \cA)_{\sharp,2}$, we can also obtain the commutativity of the diagram depicted in Corollary \ref{Cor:dPH-PH} for the quasi-Poisson case; one should only replace the bottom part of the left diagram by 
$\PH_{\Gl_{\bf n}}(\cA_{\bf n}) \to {\mathrm H}_{CE;\Gl_{\bf n}}(\cA_{\bf n})$ since the corresponding complexes are made of $\Gl_{\bf n}$-invariant multiderivations.
\end{remark}

\subsection{The gauged double Poisson case revisited}

Recall Definition~\ref{Def:gaugDBR} for a gauged Poisson element in $(\mb T^\ast \cA)_{\sharp,2}$, and Definition \ref{Def:gaugBiv} for a gauged Poisson bivector in ${\mf X}^2(\cA_{\bf n})^{\Gl_{\bf n}}$.

\begin{proposition}
If $P\in (\mb T^\ast \cA)_{\sharp,2}$ is gauged Poisson, 
then $\tr(P)\in {\mf X}^2(\cA_{\bf n})^{\Gl_{\bf n}}$ is a gauged Poisson bivector. 
\end{proposition}
\begin{proof}
By definition, there exist some $R_s\in (\mb T^\ast \cA)_2$ such that 
$\brSN{P,P}=\sum_s \Delta_s R_s$ modulo graded commutators. Then, 
as a consequence of \eqref{Eq:DelRep}, we can write 
\begin{equation}
\tr(\brSN{P,P})=\sum_{s\in I} \tr(\Delta_s R_s)
= \sum_{s\in I}\sum_{i,j} (E^{(s)}_{ji})_{\cA_{\bf n}} \wedge (R_s)_{ji}\,.
\end{equation}
This expression is $[\tr(P),\tr(P)]_{\SN}$ due to Proposition~\ref{Pr:Tr-Lie}. 
Thus, 
$[\tr(P),\tr(P)]_{\SN}$ admits a decomposition of the form \eqref{Eq:geo-gP} and it belongs to 
$\im(j^{\Gl_{\bf n}}_{\cA_{\bf n},3})$. 
Since $\tr(P)$ is $\Gl_{\bf n}$-invariant, it is a gauged Poisson bivector by Definition~\ref{Def:gaugBiv}.
\end{proof}

\begin{theorem}  \label{Thm:gdP-rep2} 
If $P\in (\mb T^\ast \cA)_{\sharp,2}$ is a gauged double Poisson element, then 
the differential $\dd_P$ on ${\mc D}_\cA$ from \eqref{Eq:dP-AKKN} and the differential $\dd_{\tr(P)}$ on ${\mf X}(\cA_{\bf n}^{\Gl_{\bf n}})$ (cf. Definition~\ref{Def:gPH}) induce the morphism of complexes 
\begin{equation} \label{Eq:gdP-rep2}
    \tr : ({\mc D}_\cA, \dd_P) \longrightarrow 
    (\mf X(\cA_{\bf n}^{\Gl_{\bf n}}),  \dd_{\tr(P)}),  
\end{equation}
which descends to a linear map $\gdPH(\cA) \to \gPH_{\Gl_{\bf n}}(\cA_{\bf n})$ in cohomology.  
\end{theorem}
\begin{proof}
This is just Theorem \ref{Thm:gdP-rep1} in the present setting.
\end{proof}

\begin{remark}
We can not state the second part of  Theorem \ref{Thm:gdP-rep1} as well as Corollary \ref{Cor:gdPH-PH} in the gauged Poisson setting. Indeed, we were not able to show that the double bracket associated with a gauged double Poisson element defines a square-zero differential on the completed complex $\widehat{{\mc D}}_\cA$ (cf. Section \ref{ss:gdPCoh}) through \eqref{Eq:dP-gen-0}--\eqref{Eq:dP-gen}, as we emphasized in Remark~\ref{Rem:Compl-gdP}.   
\end{remark}

\begin{example}
Recall from Example~\ref{Exmp:Rep-gdP-kx} that $\kk[x]_n\simeq \kk[\gl_n]$. 
In particular, the invariant ring $\kk[\gl_n]^{\Gl_n}$ is a polynomial ring in $n$ variables. 
Its Poisson structure induced by an arbitrary $P\in (\mb T^\ast \kk[x])_{\sharp,2}$ is the zero Poisson structure (because the corresponding double Poisson bracket descends to the zero Lie bracket on the abelianization $\kk[x]_\sharp = \kk[x]$). 
Hence we always get $\PH(\cA_n^{\Gl_n})\simeq {\mc X}(\cA_n^{\Gl_n})$; in the noncommutative setting, we also proved that we get the full complex: $\gdPH(\kk[x])={\mc D}_{\kk[x]}$, cf. Proposition \ref{Pr:gdPH-kxGen}.  
\end{example}

%%%%%%%%%%% NEW PART %%%%%%%%%%%%%%%
%%%%%%%%%%% NEW PART %%%%%%%%%%%%%%%
%%%%%%%%%%% NEW PART %%%%%%%%%%%%%%%
%%%%%%%%%%% NEW PART %%%%%%%%%%%%%%%
%%%%%%%%%%% NEW PART %%%%%%%%%%%%%%%
%%%%%%%%%%% NEW PART %%%%%%%%%%%%%%%
%%%%%%%%%%% NEW PART %%%%%%%%%%%%%%%

\part{Double Poisson vertex algebra cohomologies}

%%%%%%%%%%% NEW CHAPTER %%%%%%%%%%%%%%%
%%%%%%%%%%% NEW CHAPTER %%%%%%%%%%%%%%%
%%%%%%%%%%% NEW CHAPTER %%%%%%%%%%%%%%%
%%%%%%%%%%% NEW CHAPTER %%%%%%%%%%%%%%%

\chapter{Poisson vertex algebra cohomologies}
\label{Ch:PVAcoh}

We review in this chapter the definition of the basic 
and reduced cohomology complexes for Poisson vertex 
algebras \cite{DSK13} and the variational Poisson 
vertex algebra cohomology complex \cite{BDSK20} (the 
latter is called the Poisson vertex algebra cohomology 
complex in \cite{DSK13}). The main reference is 
\cite{DSK13}, see also \cite{BKV,DSK09,BDSK,BDSK20,BDSHKV}.

%%%
\section{Poisson vertex algebras}\label{sec:pva}
Let $V$ be a commutative differential algebra.
A \emph{$\lambda$-bracket}\glslink{lambda-b}{} on $V$ is a linear
map $\{-_{\lambda}-\}:V\otimes V\to V[\lambda]$
satisfying ($a,b,c\in V$)
\begin{subequations}
\begin{align}
\label{sesquiLCA}
&  \text{(sesquilinearity)} 
&
&\{\partial a_\lambda b\}=-\lambda \{a_\lambda b\}\,,
\qquad
\{a_\lambda\partial b\}=(\lambda+\partial)\{a_\lambda b\}
\,;
\\
\label{skewLCA}
&  \text{(skewsymmetry)} 
&
&\{a_\lambda b\}=-(|_{x=\partial}\{b_{-\lambda-x}a\})
\,;
\\
\label{lleibnizPVA}
& \text{(left Leibniz rule)} 
&
&\{a_\lambda bc\}
=\{a_\lambda b\}c+\{a_\lambda c\}b
\,.  
\end{align}
\end{subequations}
%
%\begin{enumerate}[(i)]
%\item sesquilinearity:
%\begin{equation}\label{sesquiLCA}
%\{\partial a_\lambda b\}=-\lambda \{a_\lambda b\}\,,
%\qquad
%\{a_\lambda\partial b\}=(\lambda+\partial)\{a_\lambda b\}
%\,;
%\end{equation}
%\item skewsymmetry:
%\begin{equation}\label{skewLCA}
%\{a_\lambda b\}=-(|_{x=\partial}\{b_{-\lambda-x}a\})
%\,;
%\end{equation}
%\item Leibniz rule:
%\begin{equation}\label{lleibnizPVA}
%\{a_\lambda bc\}
%=\{a_\lambda b\}c+\{a_\lambda c\}b
%\,.
%\end{equation}
%\end{enumerate}
In \eqref{skewLCA} %, \eqref{rleibnizPVA}
and in the sequel
we are using the (commutative version of) notation \eqref{eq:notation}. 
Note that skewsymmetry \eqref{skewLCA} and the Leibniz rule \eqref{lleibnizPVA} imply 
\begin{equation}\label{rleibnizPVA}
\{ab_\lambda c\}
=\{a_{\lambda+x} c\}(|_{x=\partial}b)
+\{b_{\lambda+x} c\}(|_{x=\partial}a)
\,.  \qquad  \text{(right Leibniz rule)} 
\end{equation}
A \emph{Poisson vertex algebra} (PVA) is a commutative differential algebra $V$
endowed with a $\lambda$-bracket $\{-_\lambda-\}:V\otimes V\to V[\lambda]$ satisfying  the Jacobi identity ($a,b,c\in V$)
\begin{equation}\label{JacobiLCA}
\{ a_\lambda \{b_\mu c\}\}-\{b_\mu[a_\lambda c\}\}
=\{\{a_\lambda b\}_{\lambda+\mu} c\}
\,.
\end{equation}
\begin{example}\label{exa:SR}
Recall from \cite{K} that a \emph{Lie conformal algebra} (LCA) is a $\kk[\partial]$-module $R$ endowed with a $\lambda$-bracket
$[-_\lambda -]:R\otimes R\to R[\lambda]$ satisfying
\eqref{sesquiLCA}, \eqref{skewLCA} and \eqref{JacobiLCA}.
The symmetric algebra $S(R)$ of $R$, has a natural structure of PVA with the $\lambda$-bracket on $R$ extended to $S(R)$ by the Leibniz rule
\eqref{lleibnizPVA}.
\end{example}
\begin{example}\label{exa:PVAtoPA}
Let $V$ be a PVA and let $I= \langle\partial V\rangle\subset V$ be the ideal (with respect to the commutative associative product of $V$) generated by the elements $\partial a$, $a\in V$. Then, the quotient commutative algebra $q(V)=V/I$\glslink{qV}{} has the structure of a Poisson algebra with Poisson bracket defined by $(a,b\in V$)
\begin{equation} \label{Eq:br-PVAquotient}
    \{\pi(a),\pi(b)\}=\pi\left(\{a_\lambda b\}|_{\lambda=0}\right) \,,
\end{equation}
where $\pi:V\to V/I=q(V)$ denotes the canonical quotient map.
\end{example}
\begin{example}\label{exa:jetPVA}
Let $A$ be a commutative algebra. 
Recall that the \emph{jet algebra} $J_\infty A$\glslink{jet-comm}{} of $A$
is the unique (up to isomorphism) commutative differential algebra
endowed with an (injective) commutative algebra homomorphism $\iota:A\hookrightarrow J_\infty A$,
satisfying the following universal property:
for every commutative algebra homomorphism $f:A\to B$ from $A$ to a commutative differential algebra $B$,
there exists a unique homomorphism of commutative differential algebras $\tilde f:J_\infty A\to B$
making the following diagram commute:
\begin{equation}\label{eq:univ-jet}
% \text{this diagram has been commented to speed up compiling}
% \begin{comment}
\begin{tikzcd}
A \arrow[r,"f"]\arrow[d,hook,"\iota"']& B
\\
J_{\infty}A\arrow[ur,bend right,"\exists!\tilde f"']
\end{tikzcd}
% \end{comment}
\end{equation}
It is shown in \cite{Ar12} that any Poisson algebra $A$ can be naturally extended to yield a structure of a PVA on its jet algebra $J_{\infty}A$ such that
$$
\{a_\lambda b\}=\{a,b\}\,,
$$
for every $a,b\in A\subset J_{\infty}A$.
Moreover, we have that $q(J_{\infty}A)=A$, where $q(V)$ denotes the Poisson algebra associated to a PVA $V$ in Example~\ref{exa:PVAtoPA}.
\end{example}
Let $V$ be a PVA and let $\tint V\to V_\sharp=V/\partial V$ denote the quotient map. It is shown in \cite{BDSK} that $V_\sharp$ is a Lie algebra with Lie bracket defined by ($\tint a,\tint b\in V_\sharp$)
$$
\{\tint f,\tint g\}=\tint \{a_\lambda b\}|_{\lambda=0}
\,.
$$
Moreover, the Lie algebra $V_\sharp$ acts on $V$ by derivations
(of the commutative associative product of $V$) commuting with $\partial$ via the formula ($\tint a \in V_\sharp$, $b\in V$)
\begin{equation}\label{20250702:eq3}
\{\tint a,b\}=\{a_\lambda b\}|_{\lambda=0}
\,.
\end{equation}

%%%
\section{Basic Poisson vertex algebra cohomology}\label{sec:basic}
We introduce the space of \emph{basic cochains} $\widetilde{\Gamma}(V)$.
For $n=0$ we let $\widetilde{\Gamma}^0(V)=V$, and, for $n\geq1$, $\widetilde{\Gamma}^n(V)$ consists of all linear maps\glslink{basic-coch}{}
$$
X:V^{\otimes n}\to V[\lambda_1,\dots,\lambda_n]\,,
\quad
a_1\otimes \dots\otimes a_n\mapsto X_{\lambda_1,\dots,\lambda_n}(a_1,\dots,a_n)
\,,
$$
which satisfy
\begin{enumerate}[(a)]
\item sesquilinearity:
\begin{equation}\label{eq:sesquimaps-comm}
\begin{split}
&X_{\lambda_1,\dots,\lambda_n}(a_1,\dots,a_{i-1},\partial a_i,a_{i+1},\dots,a_n)
\\
&=-\lambda_i X_{\lambda_1,\dots,\lambda_n}(a_1,\dots,a_{i-1},a_i,a_{i+1},\dots,a_n)\,,
\end{split}
\end{equation}
for all $i=1,\dots,n$ and $a_1,\dots,a_n\in V$;
\item skewsymmetry:
\begin{equation}\label{eq:skewmaps-comm}
X_{\lambda_1,\dots,\lambda_n}(a_1,\dots,a_n)
=-\sgn(\tau)X_{\lambda_{\tau(1)},\dots,\lambda_{\tau(n)}}(a_{\tau(1)},\dots,a_{\tau(n)})\,,
\end{equation}
for every $\tau\in S_n$ and $a_1,\dots,a_n\in V$;
\item  Leibniz rules:
\begin{equation}
\begin{split}\label{eq:Leibnizmaps-comm}
&X_{\lambda_1,\dots,\lambda_n}(a_1,\dots,a_{i-1},b c,a_{i+1},\dots,a_n)
\\
&=X_{\lambda_1,\dots,\lambda_{i-1},\lambda_i+x,\lambda_{i+1},\dots\lambda_n}
(a_1,\dots,a_{i-1},b,a_{i+1},\dots,a_n)(|_{x=\partial}c)\\
&+X_{\lambda_1,\dots,\lambda_{i-1},\lambda_i+x,\lambda_{i+1},\dots\lambda_n}
(a_1,\dots,a_{i-1},c,a_{i+1},\dots,a_n)(|_{x=\partial}b)
\,,
\end{split}
\end{equation}
for all $i=1,\dots, n$ and $a_i,b,c\in V$.
\end{enumerate}
\begin{remark}
The skewsymmetry \eqref{eq:skewmaps-comm} and the Leibniz rule
\eqref{eq:Leibnizmaps-comm} for $i=n$ imply the Leibniz rules
for every $i=1,\dots,n-1$.
\end{remark}
\begin{example}\label{exa:rcder}
We have that $\widetilde{\Gamma}^1(V)=\RCDer(V)$ is the space of \emph{right conformal derivations} \cite{DSK13}, namely linear maps $X:V\to V[\lambda]$ such that ($a\in V$)
$$
X_\lambda(\partial a)=-\lambda X_\lambda(a)\,,
$$
and ($a,b\in V$)
\begin{equation}\label{eq:rcder}
X_\lambda(ab)
=X_{\lambda+x}(a)(|_{x=\partial}b)
+X_{\lambda+x}(b)(|_{x=\partial}a)
\,.
\end{equation}
This space has the structure of a $\kk[\partial]$-module given by ($a\in V$)
$$
(\partial X)_\lambda(a)=(\lambda+\partial)X_\lambda(a)
\,.
$$
Furthermore, for $X,Y\in\RCDer(V)$ we can define a
$\lambda$-bracket with values in $\RCDer(V)[[\lambda]]$
using the formula ($a\in V$)
$$
[X_\lambda Y]_\mu(a)=\left(|_{x=\partial}X_{-\lambda-x}(Y_\mu(a))\right)
-Y_{\lambda+\mu}(X_\mu(a))
\,.
$$
If $V$ is a finitely generated commutative differential algebra then the above $\lambda$-bracket is polynomial
in $\lambda$ and endows $\RCDer(V)$ with the structure of an LCA.
\end{example}
Elements of $\widetilde{\Gamma}^n(V)$ are called basic $n$-cochains. We define the $\mb Z_{\geq0}$-graded space of basic cochains as
$$
\widetilde{\Gamma}(V)
=\bigoplus_{n\in\mb Z_{\geq0}}\widetilde{\Gamma}^n(V)
\,.
\glslink{Gammatilde}{}
$$
On $\widetilde{\Gamma}(V)$ we introduce the concatenation product
$\widetilde{\Gamma}^m(V)\otimes \widetilde{\Gamma}^n(V)\to \widetilde{\Gamma}^{m+n}(V)$ defined for 
$X\in\widetilde{\Gamma}^m(V)$, $Y\in\widetilde{\Gamma}^n(V)$ by 
\begin{equation}\label{eq:conc-prod}
\begin{split}
&(XY)_{\lambda_1,\dots\lambda_{m+n}}(a_1,\dots,a_{m+n})
\\
&=\sum_{\substack{1\leq i_1<\dots<i_m\leq m+n\\1\leq i_{m+1}<\dots <i_{m+n}\leq m+n}}
X_{\lambda_{i_1},\dots,\lambda_{i_m}}(a_{i_1},\dots,a_{i_m})
Y_{\lambda_{i_{m+1}},\dots,\lambda_{i_{m+n}}}(a_{i_{m+1}},\dots,a_{i_{m+n}})
\,.
\end{split}
\end{equation}
The action of $\partial$ on $V$ is extended to a derivation on $\widetilde{\Gamma}(V)$ by letting
($X\in\widetilde{\Gamma}^m(V)$)
$$
(\partial X)_{\lambda_1,\dots,\lambda_m}(a_1,\dots,a_m)
=(\lambda_1+\dots+\lambda_m+\partial)X_{\lambda_1,\dots,\lambda_m}(a_1,\dots, a_m)
\,.
$$
This makes $\widetilde{\Gamma}(V)$ a differential algebra. We also make $\widetilde{\Gamma}(V)$ a superspace, by assigning $|X|=n\bmod 2$, for $X\in\widetilde{\Gamma}^n(V)$ (this is compatible with the $\mb Z_{\geq0}$-grading of $\widetilde{\Gamma}(V)$). 

Let us assume that $V$ is a PVA.
For every $X\in\widetilde{\Gamma}^n(V)$ let
$\tilde\delta(X):V^{\otimes(n+1)}\to V[\lambda_1,\dots,\lambda_{n+1}]$
be the linear map defined by\glslink{Tildelta}{}
\begin{equation}
\begin{split}\label{eq:diff-PVA}
&\tilde\delta(X)_{\lambda_1,\dots,\lambda_{n+1}}(a_1,\dots,a_{n+1})
=
\sum_{s=1}^{n+1}(-1)^{s+1} \{{a_s}_{\lambda_s}
X_{\lambda_1,\stackrel{s}{\check{\dots}},\lambda_{n+1}}
(a_1,\stackrel{s}{\check{\dots}},a_{n+1})\}
\\
&
+\sum_{1\leq s<t\leq n+1}(-1)^{s+t}X_{\lambda_s+\lambda_t,\lambda_1,\stackrel{s}{\check{\dots}} \stackrel{t}{\check{\dots}},\lambda_{n+1}}(\{a_s{}_{\lambda_s} a_t\},a_1,\stackrel{s}{\check{\dots}} \stackrel{t}{\check{\dots}},a_{n+1})
\,.
\end{split}
\end{equation}
\begin{example}\label{exa:deltaP-PVA}
Let $c\in V=\widetilde{\Gamma}^0(V)$. Then, from \eqref{eq:diff-PVA} we have
\begin{equation}\label{exa1-PVA}
\tilde\delta(c)=\{-_\lambda c\}
\,.
\end{equation}
It follows from sesquilinearity \eqref{sesquiLCA}
and the right Leibniz rule \eqref{rleibnizPVA} that
$\tilde\delta(c)\in\RCDer(V)=\widetilde{\Gamma}^1(V)$.
Furthermore, 
for a right conformal derivation
$X\in\widetilde{\Gamma}^1(V)$ we have, using again \eqref{eq:diff-PVA},
\begin{equation}\label{exa2-PVA}
\tilde \delta(X)_{\lambda,\mu}(a,b)
=\{ a_\lambda X_\mu(b)\}-\{ b_\mu X_\lambda(a)\}
-X_{\lambda+\mu}(\{a_\lambda b\})
\,.
\end{equation}
It can be checked directly that $\tilde\delta(X)\in\widetilde{\Gamma}^2(V)$ provided that skewsymmetry \eqref{skewLCA} holds. Moreover, combining equations \eqref{exa1-PVA} and \eqref{exa2-PVA} we get
$$
\tilde\delta^2(c)_{\lambda,\mu}(a,b)
=\{a_\lambda \{ b_\mu c\}\}-\{ b_\mu \{a_\lambda c\}\}
-\{\{a_\lambda b\}_{\lambda+\mu} c\}
$$
which vanishes due to the Jacobi identity \eqref{JacobiLCA}.
\end{example}
By Example~\ref{exa:deltaP-PVA} we have that,
using the axioms of PVA,
equation \eqref{eq:diff-PVA} gives a well defined map
$\tilde\delta:\widetilde{\Gamma}^n(V)\to\widetilde{\Gamma}^{n+1}(V)$, for $n=0,1$, such that
$\tilde\delta^2:\widetilde{\Gamma}^0(V)\to\widetilde{\Gamma}^{2}(V)$ is the zero map. In fact, we have the following result.
\begin{theorem}[\cite{DSK13}]\label{thm:DSK1}
Let $V$ be a PVA and let $\widetilde{\Gamma}(V)$ be the superspace of basic cochains.
\begin{enumerate}[(a)]
\item 
The linear map $\tilde \delta$ defined in \eqref{eq:diff-PVA} gives a well defined map
$\tilde \delta:\widetilde{\Gamma}^n(V)\to\widetilde{\Gamma}^{n+1}(V)$, for every
$n\in\mb Z_{\geq0}$.
\item The linear map $\tilde\delta:\widetilde{\Gamma}(V)\to\widetilde{\Gamma}(V)$ is an odd derivation of the concatenation product \eqref{eq:conc-prod} such that $\tilde \delta^2=0$.
\item 
The linear map $\tilde\delta$ commutes with the action of $\partial$ on $\widetilde{\Gamma}(V)$.
\end{enumerate}
\end{theorem}
By Theorem~\ref{thm:DSK1}(a) and (b) we get a complex $(\widetilde{\Gamma}(V),\tilde\delta)$ which
is called the \emph{basic PVA complex} of the PVA $V$.
The cohomology\glslink{HbasV}{}
$$
\coH_{\textrm{bas}}(V)=\coH(\widetilde{\Gamma}(V),\tilde\delta)
=\bigoplus_{n\in\mb Z_{\geq0}} \coH_{\textrm{bas}}^n(V)\,,
\quad \coH_{\textrm{bas}}^n(V)=\ker(\tilde\delta|_{\widetilde{\Gamma}^n(V)})/\tilde\delta(\widetilde{\Gamma}^{n-1}(V))\,,
$$
of this complex is called the \emph{basic PVA cohomology} of $V$.
\begin{example}\label{exa:basiclow}
From equation \eqref{exa1-PVA} we get
that the zeroth basic cohomology space is
the center $\cent(V)$ of the PVA $V$:
$$
\coH^0_{\textrm{bas}}(V)
=\{a\in V\mid \{a_{\lambda}b\}=0\text{ for every }b\in V\}=\cent(V)
\,.
$$
We call $\tilde\delta(\widetilde{\Gamma}^0(V))$
the space of \emph{inner right conformal derivations}, namely the right conformal derivations of the form $X=\{-_\lambda a\}$ for some $a\in V$ (cf \eqref{exa1-PVA}), and we call $\ker(\tilde\delta|_{\widetilde{\Gamma}^n(V)})$
the space of \emph{Poisson right conformal derivations},
namely the right conformal derivations $X$ which satisfy the following ``derivation" property with respect to the $\lambda$-bracket of $V$ (cf. \eqref{exa2-PVA}):
$$
X_{\lambda+\mu}(\{a_\lambda b\})=\{a_\lambda X_\mu(b)\}
-\{b_\mu X_{\lambda}(a)\}\,,\qquad a,b\in V\,.
$$
Hence, we have
$$
\coH^1_{\textrm{bas}}(V)
=\frac{\{\text{Poisson right conformal derivations}\}}{\{\text{inner right conformal derivations}\}}
\,.
$$
\end{example}

%%%
\section{Reduced Poisson vertex algebra cohomology}\label{sec:red}
By Theorem~\ref{thm:DSK1}(c) we have that $(\partial\widetilde{\Gamma}(V),\tilde\delta)
\subset(\widetilde{\Gamma}(V),\tilde\delta)$ is a subcomplex. We can thus consider the
complex
$(\Gamma(V),\delta)$ where\glslink{GammaRed}{}
$$
\Gamma(V)=\widetilde{\Gamma}(V)/\partial\widetilde{\Gamma}(V)
=\bigoplus_{n\in\mb Z_{\geq0}}\Gamma^n(V)
\,,
\quad
\Gamma^n(V)=\widetilde{\Gamma}^n(V)/\partial\widetilde{\Gamma}^n(V)
\,,
$$
and $\delta:\Gamma(V)\to\Gamma(V)$ is the differential induced by $\tilde\delta$: for $X\in\widetilde{\Gamma}(V)$, the action of $\delta$ on the coset $[X]=X+\partial\widetilde{\Gamma}(V)\in\Gamma (V)$ is given by\glslink{deltadiff}{}
$$
\delta([X])=[\tilde\delta(X)]
\,.
$$
For example, we have $\Gamma^0(V)=V/\partial V=V_\sharp$. Let us denote by $\tint: \widetilde{\Gamma}(V)\to\Gamma(V)$
the canonical quotient map.
Let $X\in\widetilde{\Gamma}^1(V)$. Then $\tint X:V\to V$ is the linear map defined by ($a\in V$)
$$
(\tint X)(a)=(|_{x=\partial}X_{-x}(a))\in V\,.
$$
It follows from the Leibniz rule \eqref{eq:Leibnizmaps-comm} that $\tint X$ is a derivation for the commutative associative product of $V$ which, by the sesquilinearity \eqref{eq:sesquimaps-comm}, commutes with $\partial$. Hence, we get a map 
$\tint:\widetilde{\Gamma}^1(V)\to \Vect(V)^{\partial}$,
where $\Vect(V)^{\partial}$ denotes the space of derivations of the product of $V$ commuting with $\partial$. It is proved in \cite{DSK13} that
$\ker \tint=\partial\widetilde{\Gamma}^1(V)$, so we get an induced injective linear map
\begin{equation}\label{eq:20250702:eq1}
\Gamma^1(V)\hookrightarrow \Vect(V)^{\partial}\,.
\end{equation}
The complex $(\Gamma(V),\delta)$
is called the \emph{reduced PVA complex}
of $V$ and its cohomology
$$
\coH_{\textrm{red}}(V)=\coH(\Gamma(V),\delta)
=\bigoplus_{n\in\mb Z_{\geq0}}\coH_{\textrm{red}}^n(V)
\,,
\quad \coH_{\textrm{red}}^n(V)=\ker(\delta|_{\Gamma^n(V)})/\delta(\Gamma^{n-1}(V))\,,
$$
is called the \emph{reduced PVA cohomology} of $V$.\glslink{HredV}{}

By Theorem~\ref{thm:DSK1} we have a short exact sequence of complexes
$$
0\rightarrow \partial\widetilde{\Gamma}(V)
\rightarrow \widetilde{\Gamma}(V)
\rightarrow \Gamma(V)\rightarrow 0
$$
which leads to the following long exact sequence in cohomology
\begin{equation}\label{eq:les}
% \text{this diagram has been commented to speed up compiling}
% \begin{comment}
\hspace{-1.3cm}
\begin{tikzcd}
0 \arrow[r]& \coH^0(\partial\widetilde{\Gamma}(V),\tilde\delta)\arrow[r]
& \coH^0_{\textrm{bas}}(V)\arrow[r] &
\coH^0_{\textrm{red}}(V) \arrow[dll,in=175,out=-5]&
\\
&\coH^1(\partial\widetilde{\Gamma}(V),\tilde\delta)\arrow[r]&
\coH^1_{\textrm{bas}}(V) \arrow[r]
&\coH^1_{\textrm{red}}(V) \arrow[dll,in=175,out=-5]&
\\
&\coH^2(\partial\widetilde{\Gamma}(V),\tilde\delta)\arrow[r]&
\coH^2_{\textrm{bas}}(V) \arrow[r]
&\coH^2_{\textrm{red}}(V)\arrow[r] &\dots
\end{tikzcd}
% \end{comment}
\end{equation}

\section{The variational Poisson vertex algebra cohomology}\label{sec:var-PVA}
In this section we review the definition 
of poly-$\lambda$-brackets (which are called \emph{cochains} in \cite{DSK13}) and the definition 
of the variational PVA cohomology.

Let $V$ be a commutative differential algebra. For $n\geq1$, an \emph{$n$-$\lambda$-bracket} on $V$ is a linear map\glslink{n-lambda}{}
$$
\{-{}_{\lambda_1}-\dots-{}_{\lambda_{n-1}}-\}:V^{\otimes n}\to V[\lambda_1,\ldots,\lambda_{n-1}]
$$
which maps $a_1\otimes \dots\otimes a_n\mapsto \{ a_1{}_{\lambda_1}a_2\dots a_{n-1}{}_{\lambda_{n-1}}a_n\}$, satisfying
\begin{enumerate}[(a)]
\item sesquilinearity:
\begin{equation}\label{eq:sesquipoly1}
\{a_1{}_{\lambda_1}\cdots{}_{\lambda_{i-1}}(\partial a_i)_{\lambda_i}\cdots a_{n-1}{}_{\lambda_{n-1}}a_n\}
=-\lambda_i\{ a_1{}_{\lambda_1}\cdots a_{n-1}{}_{\lambda_{n-1}}a_n\}\,,
\end{equation}
for all $i=1,\ldots,n-1$, and
\begin{equation}\label{eq:sesquipoly2}
\{a_1{}_{\lambda_1}\cdots a_{n-1}{}_{\lambda_{n-1}}(\partial a_n)\}
=(\lambda_1\cdots+\lambda_{n-1}+\partial)
\{ a_1{}_{\lambda_1}\cdots a_{n-1}{}_{\lambda_{n-1}}a_n\}\,;
\end{equation}
\item skewsymmetry:
\begin{equation}\label{eq:skewpoly}
\{a_1{}_{\lambda_1}\cdots a_{n-1}{}_{\lambda_{n-1}}a_n\}
=\sgn(\tau)
|_{\lambda_n=\lambda_n^\dagger}\{ a_{\tau(1)}{}_{\lambda_{\tau(1)}}\cdots a_{\tau(n-1)}{}_{\lambda_{\tau(n-1)}}a_{\tau(n)}\}\,,
\end{equation}
for every $\tau\in S_n$, where we are denoting
\begin{equation}\label{eq:dagger}
\lambda_n^\dagger=-\lambda_1-\dots-\lambda_{n-1}-\partial
\,.
\end{equation}
\item Leibniz rules:
\begin{equation}\label{eq:leibnizpoly2}
\begin{aligned} 
\{ a_1{}_{\lambda_1}\cdots bc_{\lambda_i}\dots a_{n-1}{}_{\lambda_{n-1}}a_n\}
&=
\{ a_1{}_{\lambda_1}\cdots b_{\lambda_i+x}\dots a_{n-1}{}_{\lambda_{n-1}}a_n\}(|_{x=\partial}c) \\ 
&+\{a_1{}_{\lambda_1}\cdots c_{\lambda_i+x}\dots a_{n-1}{}_{\lambda_{n-1}}a_n\}(|_{x=\partial}b)
\,,
\end{aligned}
\end{equation}
for all $i=1,\ldots,n-1$, and
\begin{equation}\label{eq:leibnizpoly}
\{a_1{}_{\lambda_1}\cdots a_{n-1}{}_{\lambda_{n-1}}bc\}
=\{ a_1{}_{\lambda_1}\cdots a_{n-1}{}_{\lambda_{n-1}}b\}c
+\{a_1{}_{\lambda_1}\cdots a_{n-1}{}_{\lambda_{n-1}}c\}b
\,.
\end{equation}
\end{enumerate}
\begin{remark}\label{rem:skewPVA}
The skewsymmetry \eqref{eq:skewpoly} and the Leibniz rule \eqref{eq:leibnizpoly} imply the Leibniz rules \eqref{eq:leibnizpoly2} for all $i=1,\dots,n-1$.
\end{remark}
We set $C^0(V)=V_\sharp$. For $n\geq1$, we denote by
$C^n(V)$ the space of $n$-$\lambda$-brackets and we define\glslink{PolyLamb}{}
$$
C(V)=\bigoplus_{n\in\mb Z_{\geq0}}C^n(V)
$$
to be the space of poly-$\lambda$-brackets.
For example, we have that
$C^1(V)=\Vect(V)^{\partial}$, while $C^2(V)$ is the 
space of $\lambda$-brackets on $V$ (cf. 
Section~\ref{sec:pva}).
Hence, a Poisson vertex algebra structure on $V$ amounts to the
choice of an element $[-_\lambda-]\in C^2(V)$ satisfying the Jacobi 
identity \eqref{JacobiLCA} (the use of the square brackets will be 
clear in the upcoming equation \eqref{eq:dH-poly-n}).

Let then $V$ be a Poisson vertex algebra with respect to a $\lambda$-bracket $[-_\lambda -]$\glslink{lambda-b}{}.
We define a linear map $d$ on the space of poly-$\lambda$-brackets $ V$ as follows.
For $\tint f\in V_\sharp=C^0(V)$ we let (cf. \eqref{20250702:eq3})\glslink{dd}{}
\begin{equation}\label{eq:dH-poly-0}
\dd(\tint f)(a)=-[\tint f,a]=-[f_{\lambda}a]|_{\lambda=0}\,,
\qquad a\in V
\,.
\end{equation}
For $n\geq1$, given $c=\{-{}_{\lambda_1}-\dots-{}_{\lambda_{n-1}}-\}\in C^n( V)$, we let
\begin{equation}\label{eq:dH-poly-n}
\begin{split}
&\dd(c)_{\lambda_1,\dots,\lambda_{n}}(a_1,\dots,a_{n+1})
=\sum_{s=1}^{n+1}(-1)^{s+1}|_{\lambda_{n+1}=\lambda_{n+1}^\dagger}[a_s{}_{\lambda_s} \{a_1{}_{\lambda_1}\stackrel{s}{\check{\dots}}a_n{}_{\lambda_{n}}a_{n+1}\}]
\\
&+\sum_{1\leq s< t\leq n+1}(-1)^{s+t}
|_{\lambda_{n+1}=\lambda_{n+1}^\dagger}\{ [a_s{}_{\lambda_s}a_t]{}_{\lambda_s+\lambda_t} a_1{}_{\lambda_1}\stackrel{s}{\check{\dots}} \stackrel{t}{\check{\dots}} a_{n}{}_{\lambda_{n}}a_{n+1}\}
\,.
\end{split}
\end{equation}
\begin{remark}
Using skewsymmetry \eqref{eq:skewpoly} we rewrite \eqref{eq:dH-poly-n} as
\begin{equation}\label{20250720:eq1}
\begin{split}
&\dd(c)_{\lambda_1,\dots,\lambda_{n}}(a_1,\dots,a_{n+1})
\\
&=\sum_{s=1}^{n}(-1)^{s+1}[a_s{}_{\lambda_s} \{a_1{}_{\lambda_1}\stackrel{s}{\check{\dots}}a_n{}_{\lambda_{n}}a_{n+1}\}]
\\
&+(-1)^{n+1}[\{a_1{}_{\lambda_1}\dots_{\lambda_{n-1}}a_n\}_{\lambda_1+\dots+\lambda_{n}}a_{n+1}]
\\
&+\sum_{1\leq s< t\leq n}(-1)^{s+t}
\{ [a_s{}_{\lambda_s}a_t]_{\lambda_s+\lambda_t} a_1{}_{\lambda_1}\stackrel{s}{\check{\dots}} \stackrel{t}{\check{\dots}} a_{n}{}_{\lambda_{n}}a_{n+1}\}
\\
&+\sum_{s=1}^n(-1)^{s}
\{ a_1{}_{\lambda_1}\stackrel{s}{\check{\dots}}a_n{}_{\lambda_{n}}[a_s{}_{\lambda_s}a_{n+1}]\}
\,.
\end{split}
\end{equation}
\end{remark}
\begin{example}\label{exa:dP-PVA}
Using the skewsymmetry \eqref{skewLCA} of the $\lambda$-bracket
$[-_\lambda-]$ on $V$ we can rewrite equation \eqref{eq:dH-poly-0} as
(cf. \eqref{20250702:eq3})
\begin{equation}\label{20240828:eq0}
\dd(\tint a)
=-[a_\lambda -]|_{\lambda=0}=-[\tint a,-]
\,.
\end{equation}
As mentioned in Section~\ref{sec:pva}, the RHS of \eqref{20240828:eq0} is a derivation commuting with $\partial$, hence it lies in $C^1(V)$. Let $D\in C^1( V)$. Then, we have that
\begin{equation}\label{20240828:eq1}
\dd(D)_{\lambda}(b,c)
=[D(b)_\lambda c]+[b_\lambda D(c)]-D([b_\lambda c])
\,.
\end{equation}
It is straightforward to verify that the RHS of \eqref{20240828:eq1} lies in $C^2( V)$.
Combining equations \eqref{20240828:eq0} and \eqref{20240828:eq1} we get
$$
\dd^2(\tint a)_{\lambda}(b,c)
=\left([a_\mu [b_\lambda c]]-[b_\lambda [a_\mu c]]
-[[a_\mu b]_{\lambda+\mu} c]\right)|_{\mu=0}\,,
$$
which vanishes by the Jacobi identity \eqref{JacobiLCA}.
\end{example}
The following result has been proven in \cite{DSK13} (see also \cite{DSK09}).
\begin{theorem}\label{thm:DSK2}
For $c\in C^n( V)$ we have that $\dd(c)\in C^{n+1}(V)$ and $\dd^2(c)=0$.
\end{theorem}
By Theorem~\ref{thm:DSK2} we have a complex $(C( V),d)$, called the \emph{variational PVA complex}
of $V$, and its cohomology\glslink{PVHA}{}
$$
\PVH(V)=H(C(V),\dd)
=\bigoplus_{n\in\mb Z_{\geq0}}\PVH^n(V)
\,,
$$
where $\PVH^n(V)=\ker(\dd|_{C^n(V)})/\dd (C^{n-1}(V))$, 
is called the \emph{variational PVA cohomology} of $V$ \cite{BDSK20} (it is called PVA cohomology in \cite{DSK13}).
\begin{example}\label{exa:varlow}
As a consequence of \eqref{20240828:eq0} we have that
$$
\PVH^0(V)=\{\tint a\in V_\sharp\mid [\tint a,b]=0
\text{ for every }b\in V\}=\Cas(V)
$$
is the space of \emph{Casimir elements} of $V$, namely it is the center of the representation of the Lie algebra
$V_\sharp$ on $V$ defined by \eqref{20250702:eq3}.
Moreover, from equation \eqref{20240828:eq1}, we have that
$\ker \dd|_{C^1(V)}=\Der(V)\subset\Vect(V)^\partial$ is the subspace consisting 
of the derivations of the $\lambda$-bracket of $V$. We denote by $\inn(V)\subset\Der(V)$ the subspace of elements $D\in\Der(V)$
of the form $D=[\tint a,-]$ for some $a\in V$. Then, we have
$$
\PVH^1(V)=\Der(V)/\inn(V)
\,.
$$
Finally, it is shown in \cite{DSK13}, that $\PVH^2(V)$
parametrizes the
equivalence classes of first-order deformations of $V$ that preserve the product and the $\kk[\partial]$-module structure
(cf. \cite{NR})
\end{example}
Recall from Section~\ref{sec:basic} that we have $\Gamma^0(V)=V_\sharp=C^0(V)$, and that we have an injective linear map
$\Gamma^1(V)\hookrightarrow C^1(V)$ (cf. \eqref{eq:20250702:eq1}).
In fact, it is shown in \cite{DSK13} that we have an injective linear map
\begin{equation}\label{20250702:eq2}
\Gamma(V)\hookrightarrow C(V)\,,
\end{equation}
which restricts to an injective map $\Gamma^n(V)\hookrightarrow C^n(V)$ on any graded component.
Moreover, it is shown in \cite{BDSHKV}, that if $V$ is an algebra of differential polynomials in finitely many variables (see \cite{BDSK} for the definition), then the map in \eqref{20250702:eq2} is also surjective and we have an isomorphism of cohomologies
$$
\coH_{\textrm{red}}(V)\simeq \PVH(V)
\,.
$$

%%%
\section{Group actions and cohomology}\label{sec:PVAcoh-inv}
Let $G$ be a group acting on the Poisson vertex algebra $ V$, with $\lambda$-bracket $[-_\lambda-]$, by
Poisson vertex algebra automorphisms, namely, we have a group action
$G\times  V\to V$, $(g,a)\mapsto g\cdot a$, such that
$$
g\cdot (ab)=(g\cdot a)(g\cdot b)\,,
\quad
g\cdot\partial a=\partial(g\cdot a)\,,
\quad
g\cdot[a_\lambda b]=[g\cdot a_\lambda g\cdot b]
\,,
\quad g\in G\,,a,b\in V
\,,
$$
where the action of $G$ is extended to an action on $ V[\lambda]$
in the obvious way.

Since the action of $G$ commutes with $\partial$ we have an induced action of $G$ on $C^0( V)= V_\sharp$. We can also extend the action of $G$ on $ V$ to an action on $C^n( V)$, $n\geq1$, by letting, for $c=\{-_{\lambda_1}-\dots-_{\lambda_{n-1}}-\}\in C^n( V)$,
$g\cdot c$ be the element in
$C^{n+1}( V)$ defined by ($a_1,\dots,a_n\in V$)
\begin{equation}\label{eq:G-action-C}
(g\cdot c)(a_1,\dots,a_n)
=g\cdot\{g^{-1}\cdot a_1{}_{\lambda_1}-\dots-_{\lambda_{n-1}}g^{-1}\cdot a_n\}
\,.
\end{equation}
For $n\in\mb Z_{\geq0}$, we let 
$C^n( V)^G=\{c\in C^n( V)\mid g\cdot c=c\text{ for every }g\in G\}$ be the space of \emph{$G$-invariant $n$-$\lambda$-brackets}.\glslink{PolyLambG}{} 
Let $\tint f\in C^0( V)^G=\left( V/\partial V\right)^G$. Then, from \eqref{20240828:eq0} we have ($a\in V$)
\begin{align*}
&(g\cdot \dd(\tint f))(a)=
-g\cdot[f_\lambda g^{-1}\cdot a]|_{\lambda=0}
\\
&=-[g\cdot f_\lambda a]|_{\lambda=0}
=-[f_\lambda a]|_{\lambda=0}
=\dd(\tint f)(a)
\,,
\end{align*}
where in the third equality we used the fact that $g\cdot f-f\in\partial V$
combined with the sesquilinearity \eqref{sesquiLCA}. Hence,
$\dd(C^0( V)^G)\subset C^0( V)^G$. Furthermore, it is straightforward to check, using \eqref{eq:dH-poly-n} and the $G$-invariance
of the $\lambda$-bracket of $ V$ that
$\dd(C^n( V)^G)\subset C^n( V)^G$, for every $n\geq1$.
We then have a subcomplex $(C( V)^G=\oplus_n C^n(V)^G,\dd)\subset(C(V),\dd)$ consisting of $G$-invariant poly-$\lambda$-brackets and we can consider the corresponding cohomology\glslink{gPVHA}{}
$$
\PVH_G( V)=\coH(C(V)^G,\dd)
=\bigoplus_{n\in\mb Z_{\geq0}} \PVH_G^n(V)
\,,
$$
where $\PVH_G^n(V)=\ker(\dd|_{C^n(V)^G})/\dd (C^{n-1}(V)^G)$, 
which we call the \emph{$G$-invariant variational PVA cohomology} of $V$. Clearly, we have linear maps in cohomology
$$
\PVH_G^n(V)\rightarrow \PVH^n(V)
\,,
$$
for every $n\in\mb Z_{\geq0}$, induced by the inclusion $C^n(V)^G\subset C^n(V)$.

Since $ V^G\subset V$ is a PVA subalgebra we can also consider the variational PVA complex $(C( V^G),\dd)$ and its variational PVA cohomology $\coH(C( V^G),\dd)$.
For $n\geq1$, the restriction of $c\in C^n( V)^G$ to $( V^G)^{\otimes n}$ gives an element of $C^n( V^G)$ and we have the following commutative diagram
\begin{equation}\label{diag0}
% \text{this diagram has been commented to speed up compiling}
% \begin{comment}
\begin{tikzcd}
C^n(V) \arrow[r,"\dd"]& C^{n+1}(V)
\\
C^{n}(V)^G \arrow[r,"\dd"]\arrow[d] \arrow[u,hook]& C^{n+1}(V)^G\arrow[u,hook]\arrow[d]
\\
C^n(V^G) \arrow[r,"\dd"] & C^{n+1}(V^G)
\end{tikzcd}
% \end{comment}
\end{equation}
Hence, we have induced linear maps in cohomology
$$
\PVH_G^n(V)\rightarrow \PVH^n(V^G)\,,
$$
for every $n\geq1$.
%

%%%%%%%%%%% NEW CHAPTER %%%%%%%%%%%%%%%
%%%%%%%%%%% NEW CHAPTER %%%%%%%%%%%%%%%
%%%%%%%%%%% NEW CHAPTER %%%%%%%%%%%%%%%
%%%%%%%%%%% NEW CHAPTER %%%%%%%%%%%%%%%

\chapter[Basic and reduced \MakeLowercase{d}PVA cohomologies]{Basic and reduced double Poisson vertex algebra cohomologies}
\label{Ch:BasRed-dPVA}

In this chapter we review the definition and properties of $2$-fold $\lambda$-brackets and double Poisson vertex algebra following \cite{DSKV}. Then, we introduce the complex of basic cochains and we define the basic
and reduced double Poisson vertex algebra cohomology as a ``double" version
of the construction described in Sections \ref{sec:basic} and \ref{sec:red}.

%%%
\section{Double Poisson vertex algebras}
\label{sec:dPVA}
Let $\mc V$ be a differential algebra with derivation $\partial:\mc V\to\mc V$. A \emph{2-fold $\lambda$-bracket} on $\mc V$ is a linear map $\ldb-_\lambda-\rdb:\mc V\otimes\mc V\to(\mc V\otimes\mc V)[\lambda]$\glslink{lambda-b-nc}{} satisfying
($a,b,c\in\mc V$)
\begin{subequations}
\begin{align}
\label{eq:sesqui}
&
\text{(sesquilinearity)} 
&
&\ldb \partial a_\lambda b\rdb=-\lambda \ldb a_\lambda b\rdb
\,\,,\qquad
\ldb a_\lambda \partial b\rdb=(\lambda+\partial) \ldb a_\lambda b\rdb
\,;
\\
\label{eq:skew2}
&  \text{(skewsymmetry)} 
&
&\ldb a_{\lambda}b\rdb=-\left(|_{x=\partial}\ldb b_{-\lambda-x}a\rdb^\sigma\right)\,;
\\
\label{eq:lleibniz}
& \text{(left Leibniz rule)} 
&
&\ldb a_\lambda bc\rdb
=
\ldb a_\lambda b\rdb c+b\ldb a_\lambda c\rdb
\,;
\\
\label{eq:rleibniz}
& \text{(right Leibniz rule)} 
&
&\ldb ab_\lambda c\rdb
=
\ldb a_{\lambda+x} c\rdb \star_1 (|_{x=\partial}b)
+(|_{x=\partial}a)\star_1 \ldb b_{\lambda+x} c\rdb
\,.
\end{align}
\end{subequations}
In \eqref{eq:skew2} we are using the action of the cyclic permutation $\sigma$ on a tensor product defined in Section~\ref{sec:1.1}, and in \eqref{eq:skew2} and
\eqref{eq:rleibniz} we are using the notations \eqref{eq:notation}, \eqref{20140707:eq6-left}
and \eqref{20140707:eq6}.
We also let, for $a,b,c\in\mc V$,
\begin{equation}\label{notation}
\begin{array}{l}
\displaystyle{
\vphantom{\Big(}
\ldb a_\lambda (b\otimes c)\rdb_L:=\ldb a_{\lambda} b\rdb\otimes c
\,\,,\,\,\,\,
\ldb a_\lambda (b\otimes c)\rdb_R:=b\otimes \ldb a_{\lambda} c\rdb
\,,} \\
\displaystyle{
\vphantom{\Big(}
\ldb (a\otimes b)_\lambda c\rdb_L
:=
\ldb a_{\lambda+x} c\rdb \otimes_1 (|_{x=\partial}b)
\,.}
%\\
%\displaystyle{
%\vphantom{\Big(}
%\ldb (a\otimes b)_\lambda c\rdb_R
%:=
%(|_{x=\partial}a)\otimes_1\ldb b_{\lambda+x} c\rdb
%\,.}
\end{array}
\end{equation}

\begin{definition}\label{20140606:def-2}
A \emph{double Poisson vertex algebra} (dPVA) is a differential algebra $\mc V$,
with derivation $\partial:\,\mc V\to\mc V$, 
endowed with a $2$-fold $\lambda$-bracket $\ldb-_\lambda-\rdb:\,\mc V\times \mc V\to \mc V\otimes \mc V[\lambda]$
satisfying the Jacobi identity ($a,b,c\in\mc V$)
\begin{equation}\label{eq:jacobi2}
\ldb a_{\lambda}\ldb b{}_\mu c\rdb\rdb_L
-\ldb b_{\mu}\ldb a{}_\lambda c\rdb\rdb_R
=\ldb\ldb a_{\lambda} b\rdb_{\lambda+\mu} c\rdb_L
\,.
\end{equation}
\end{definition}
\begin{example}\label{exa:dPVAtodPA}
Let $\mc V$ be a dPVA and let $\mc I= \langle\partial \mc V\rangle\subset \mc V$ be the ideal (with respect to the associative product of $\mc V$) generated by the elements $\partial a$, $a\in \mc V$. Then, the algebra $q(\mc V)=\mc V/\mc I$\glslink{qV-nc}{} has the structure of a dPA with double Poisson bracket defined by $(a,b\in \mc V$)
\begin{equation}\label{20250822:eq1}
\{\pi(a),\pi(b)\}=(\pi\otimes\pi)\left(\ldb a_\lambda b\rdb|_{\lambda=0}\right)
\,,
\end{equation}
where $\pi:\mc V\to \mc V/\mc I=q(\mc V)$ denotes the canonical quotient map.
\end{example}
\begin{example}\label{exa:jetdPVA}
Let $\mc A$ be an algebra. 
The \emph{jet algebra}\glslink{JetAss}{} $J_\infty \mc A$ of $\mc A$
is the unique (up to isomorphism)  differential algebra
endowed with an (injective) algebra homomorphism $\iota:\mc A\hookrightarrow J_\infty \mc A$,
satisfying the following universal property:
for every algebra homomorphism $f:\mc A\to \mc B$ to a differential algebra $\mc B$,
there exists a unique homomorphism of differential algebras $\tilde f:J_\infty \mc A\to \mc B$
making the diagram \eqref{eq:univ-jet} (where we replace $A$ by $\mc A$ and $B$ by $\mc B$) commute.
It is shown in \cite{BFM} that any dPA $\mc A$ can be naturally extended to yield a structure of a dPVA on its jet algebra $J_{\infty}\mc A$ such that
$$
\ldb a_\lambda b\rdb=\ldb a,b\rdb\,,
$$
for every $a,b\in \mc A\subset J_{\infty}\mc A$.
Moreover, we have that $q(J_{\infty}\mc A)=\mc A$, where $q(\mc V)$ denotes the double Poisson algebra associated to a dPVA $\mc V$ in Example~\ref{exa:dPVAtodPA}.
\end{example}
Recall the notation $\mc V_\sharp=\mc V/(\partial \mc V+[\mc V,\mc V])$ and the canonical quotient map $\tint:\mc V\to\mc V_{\sharp}$ introduced in Section~\ref{sec:1.2}.
\begin{theorem}[\cite{DSKV}]\label{20140707:thm}
Let $\mc V$ be a dPVA, with
$2$-fold $\lambda$-bracket $\ldb-_{\lambda}-\rdb$, and let
$\{-,-\}$ be defined as in \eqref{20140707:eq5}.
\begin{enumerate}[(a)]
\item
We have a well defined Lie algebra bracket
$\{-,-\}:\,\mc V_\sharp\times\mc V_\sharp
\to\mc V_\sharp$
on $\mc V_\sharp$ given by
\begin{equation}\label{20140707:eq3c-lie}
\{\tint f,\tint g\}:=\tint \mult\ldb f_\lambda g\rdb|_{\lambda=0}
\,.
\end{equation}
\item
We have a well defined
Lie algebra action
$\{-,-\}:\,\mc V_\sharp\times\mc V\to\mc V$
of $\mc V_\sharp$ on $\mc V$
given by
\begin{equation}\label{20140707:eq3b-lie}
\{\tint f, g\}:=\mult\ldb f_\lambda g\rdb|_{\lambda=0}
\,.
\end{equation}
The action 
\eqref{20140707:eq3b-lie} defines a representation of the Lie  algebra
$\mc V_\sharp$ by derivations of $\mc V$
commuting with $\partial$.
\end{enumerate}
\end{theorem}
%

%%%
\section{Properties of 2-fold \texorpdfstring{$\lambda$}{lambda}-brackets}
In this section we let $\mc V$ be a differential algebra and we let $\ldb-_\lambda-\rdb$ be a $2$-fold $\lambda$-bracket on $\mc V$. We list below some results that will be used in the sequel.

Let us introduce the linear map
$\ldb-\,_\lambda-_\mu-\rdb:
\mc V^{\otimes 3}\to\mc V^{\otimes 3}[\lambda,\mu]$
defined by ($a,b,c\in\mc V$):
\begin{equation}
\label{eq:triple1}
\ldb a_\lambda b_\mu c\rdb
:=\ldb a_\lambda \ldb b_\mu c\rdb\rdb_L
-\ldb b_\mu\ldb a_\lambda c\rdb\rdb_R
-\ldb\ldb a_\lambda b\rdb_{\lambda+\mu}c\rdb_L\,,
\end{equation}
and the linear map $\{-,-\}:\mc V\otimes \mc V\to \mc V$
defined by ($a,b\in\mc V$)
\begin{equation}\label{20140707:eq5}
\{a,b\}=\mult\ldb a_{\lambda}b\rdb|_{\lambda=0}
\,.
\end{equation}
The next result has been proved
in \cite{DSKV}.
\begin{proposition}
\label{prop:1bis}
Let $\ldb-_\lambda-\rdb$ be a $2$-fold $\lambda$-bracket on $\mc V$.
Then
the following identity
holds in $\mc V^{\otimes2}[\lambda]$ ($a,b,c\in\mc V$):
$$%\begin{align}
%\begin{split}\label{eq:quasijacobi}
%&
\{a,\ldb b_\lambda c\rdb\}
-\ldb b_\lambda\{a, c\}\rdb
-\ldb\{a, b\}_{\lambda} c\rdb
=(\mult\otimes1)\ldb a_\mu b_\lambda c\rdb|_{\mu=0}
-(1\otimes \mult)\ldb b_\lambda a_\mu c\rdb|_{\mu=0}\,,
%%\end{split}
%%\end{align}
$$
where we set $\{a, b\otimes c\}=\{a, b\}\otimes c+b\otimes\{a, c\}$.
In particular, if $\ldb-_\lambda-\rdb$ defines a dPVA
structure on $\mc V$, then
\begin{equation}
\begin{split}\label{eq:quasijacobi2}
&\{a,\ldb b_\lambda c\rdb\}
=\ldb\{a,b\}_{\lambda} c\rdb
+\ldb b_\lambda\{a,c\}\rdb
\,,
\end{split}
\end{equation}
for every $a,b,c\in\mc V$.
\end{proposition}
Let $a\in\mc V$. Using the notation \eqref{20240805:eq1} for $D=\ldb a_\lambda-\rdb$ we get linear maps
$\ldb a_\lambda-\rdb_{(s)}:\mc V^{\otimes m}\to\mc V^{\otimes (m+1)}[\lambda]$, for ever $s=1,\dots,m$ (extending the notation in the first line of \eqref{notation} which corresponds to the case of $m=2$.).
The next result describes the compatibility properties between the linear maps $\ldb a_\lambda-\rdb_{(s)}$ and the multiplication maps \eqref{eq:mii+1}.
\begin{lemma}\label{20240821:lem1}
\begin{enumerate}[(a)]
\item
For every
$a\in\mc V$,
$B=b_1\otimes \dots\otimes b_n\in\mc V^{\otimes n}$, $i=0,\dots,n-1$, and $s=1,\dots,n$ we have
\begin{equation}\label{20240821:eq1}
\begin{split}
&\mult_{(i+1,i+2)}\ldb a_\lambda B\rdb_{(s)}
\\
&=\left\{
\begin{array}{ll}
\ldb a_\lambda\mult_{(i+1,i+2)}B\rdb_{(s-1)}\,,
&
i\leq s-3\,,
\\
b_1\otimes\dots \otimes b_{s-2}\otimes b_{s-1}\ldb a_\lambda b_s\rdb\otimes b_{s+1}\otimes\dots\otimes b_n\,,
&
i=s-2\,,
\\
b_1\otimes\dots \otimes b_{s-1}\otimes\mult\ldb a_\lambda b_s\rdb\otimes b_{s+1}\otimes\dots\otimes b_n\,,
&
i=s-1\,,
\\
b_1\otimes\dots \otimes b_{s-1}\otimes\ldb a_\lambda b_s\rdb b_{s+1}\otimes b_{s+2}\otimes \dots\otimes b_n\,,
&
i=s\,,
\\
\ldb a_\lambda\mult_{(i,i+1)}B\rdb_{(s)}\,,
&
i\geq s+1\,.
\end{array}
\right.
\end{split}
\end{equation}
\item
For every
$a\in\mc V$,
$B=b_1\otimes \dots\otimes b_n\in\mc V^{\otimes n}$ and $s=1,\dots,n$ we have
\begin{equation}\label{20240821:eq2}
\mult_{(s,s+1)}\ldb a_\lambda B\rdb_{(s+1)}
+\mult_{(s+1,s+2)}\ldb a_\lambda B\rdb_{(s)}
=\ldb a_\lambda \mult_{(s,s+1)}B\rdb_{(s)}
\,.
\end{equation}
\end{enumerate}
\end{lemma}
\begin{proof}
Equation \eqref{20240821:eq1} can be easily checked directly and equation \eqref{20240821:eq2} follows
from \eqref{20240821:eq1} and the Leibniz rule \eqref{eq:lleibniz}.
\end{proof}
The next result describes the compatibility properties between the linear maps $\ldb a_\lambda-\rdb_{(s)}$ and the
left and right $\mc V$-module structures given in \eqref{eq:star}.
\begin{lemma}\label{lem:PVAcoh-pre2}
Let $\mc V$ be a differential algebra and let $\ldb-_\lambda-\rdb$ be a $2$-fold $\lambda$-bracket on $\mc V$.
For every $a,b\in\mc V$, $B\in\mc V^{\otimes m}$ and $i=0,\dots m-1$ we have
$$
\ldb a_\lambda b\star_i B\rdb_{(s)}
=\left\{
\begin{array}{ll}
b\star_{i+1}\ldb a_{\lambda} B\rdb_{(s)}\,,
& 1\leq s\leq i
\,,
\\
b\star_i\ldb a_\lambda B\rdb_{(i+1)}
&s=i+1\,,
\\
+\mult_{(i+2,i+3)}\left(B\otimes_{m-i}\ldb a_\lambda b\rdb\right)
\,,
&
\\
b\star_i\ldb a_\lambda B\rdb_{(s)}\,,
&i+2\leq s\leq m\,,
\end{array}
\right.
$$
and
$$
\ldb a_\lambda B \star_i b\rdb_{(s)}
=\left\{
\begin{array}{ll}
\ldb a_{\lambda} B\rdb_{(s)}\star_i b\,,
& 1\leq s\leq m-1-i
\,,
\\
\ldb a_\lambda B\rdb_{(m-i)}\star_i b
&s=m-i\,,
\\
+\mult_{(m-i,m+1-i)}\left(B\otimes_{i}\ldb a_\lambda b\rdb\right)
\,,
&
\\
\ldb a_\lambda B\rdb_{(s)}\star_{i+1}b\,,
&m-i+1\leq s\leq m\,.
\end{array}
\right.
$$
\end{lemma}
\begin{proof}
The claims follow immediately from
\eqref{eq:star},~\eqref{20240805:eq1} and the Leibniz rule
\eqref{eq:lleibniz}. 
\end{proof}
The next result describes the compatibility properties between the linear maps $\ldb a_\lambda-\rdb_{(s)}$ and the
operation given in \eqref{eq:tensor-i-notation}.
\begin{lemma}\label{lem:PVAcoh-pre5}
Let $a\in\mc V$, $A\in\mc V^{\otimes n}$ and $B\in\mc V^{\otimes m}$. For $0\leq i\leq n$ and 
$1\leq s\leq m+n$ we have
\begin{equation}\label{20230811:eq1}
\ldb a_{\lambda} A\otimes_i B\rdb_{(s)}
=\left\{
\begin{array}{ll}
\ldb a_\lambda A\rdb_{(s)}\otimes_i B\,,
&
1\leq s\leq n-i\,,
\\
A\otimes_i \ldb a_\lambda B\rdb_{(s+i-n)}
\,,
&
n+1-i\leq s\leq n+m-i\,,
\\
\ldb a_\lambda A\rdb_{(s-m)}\otimes_{i+1}B\,,
&
n+m+1-i\leq s\leq n+m\,.
\end{array}
\right.
\end{equation}
\end{lemma}
\begin{proof}
It follows from a straightforward computation using~\eqref{20240805:eq1} and~\eqref{eq:tensor-i-notation}.
\end{proof}
Next, we generalize the right Leibniz rule given in \eqref{eq:rleibniz}.
\begin{lemma}\label{lem:PVAcoh-pre3}
Let $a,b\in\mc V$ and $C\in\mc V^{\otimes m}$. For every $s=1,\dots m$ we have
$$
\ldb ab_\lambda C\rdb_{(s)}=
(|_{x=\partial} a)\star_s\ldb b_{\lambda+x} C\rdb_{(s)}+
\ldb a_{\lambda+x} C\rdb_{(s)}\star_{m+1-s}(|_{x=\partial}b)
\,.
$$
\end{lemma}
\begin{proof}
It follows immediately from the Leibniz rule~\eqref{eq:rleibniz}, and equations~\eqref{eq:star}
and~\eqref{20240805:eq1}.
\end{proof}
We conclude the section by giving a generalization of the Jacobi identity \eqref{eq:jacobi2}.
\begin{lemma}\label{lem:PVAcoh-pre7}
Let $a,b\in\mc V$ and $C=c_1\otimes\dots\otimes c_m \in\mc V^{\otimes m}$.
\begin{enumerate}[(i)]
\item For every $1\leq t\leq m$ and $1\leq s\leq m+1$ we have
\begin{equation}\label{20230814:eq1}
\begin{split}
&\ldb a_{\lambda} \ldb b_\mu C\rdb_{(t)}\rdb_{(s)}
\\
&=\left\{
\begin{array}{ll}
c_1\otimes\dots\otimes \ldb a_\lambda c_s\rdb\otimes\dots\otimes \ldb b_\mu c_t\rdb\otimes\dots\otimes c_m
\,,
&
1\leq s\leq t-1\,,
\\
c_1\otimes\dots\otimes \ldb a_\lambda \ldb b_\mu c_t\rdb \rdb_L\otimes\dots\otimes c_m\,,
&
s=t\,,
\\
c_1\otimes\dots\otimes \ldb a_\lambda \ldb b_\mu c_t\rdb \rdb_R\otimes\dots\otimes c_m\,,
&
s=t+1\,,
\\
c_1\otimes\dots\otimes \ldb b_\mu c_t\rdb\otimes\dots\otimes\ldb a_\lambda c_{s-1}\rdb \otimes\dots\otimes c_m\,,
&
t+2\leq s\leq m+1
\,.
\end{array}\right.
\end{split}
\end{equation}
\item Let us assume that the Jacobi identity~\eqref{eq:jacobi2} holds. For every $1\leq t\leq m$ we have
$$
\ldb a_{\lambda} \ldb b_\mu C\rdb_{(t)}\rdb_{(t)}
-\ldb b_{\mu} \ldb a_\lambda C\rdb_{(t)}\rdb_{(t+1)}
=\ldb \ldb a_\lambda b\rdb_{\lambda+\mu} C\rdb_{L,(t)}
\,.
$$
\end{enumerate}
\end{lemma}
\begin{proof}
Part (i) follows immediately from the definition~\eqref{20240805:eq1}. Part (ii) follows from~\eqref{20240805:eq1} and the assumption using part (i).
\end{proof}

%%%
\section{The space of basic cochains}\label{sec:omega1}
Let $\mc V$ be a differential algebra.
We introduce the space of \emph{basic cochains}
$\widetilde{\Gamma}(\mc V)$ extending the construction 
of Section~\ref{sec:basic} to the noncommutative case.\glslink{AsGammatilde}{}

For $n=0$ we let $\widetilde{\Gamma}^0(\mc V)=\mc V$. 
For every $n\geq1$ we let $\widetilde{\Gamma}^n(\mc V)$ be the space of linear maps\glslink{basic-coch}{}
$$
X:\mc V^{\otimes n}\to \mc V^{\otimes (n+1)}[\lambda_1,\dots,\lambda_n]
\,,
\quad
a_1\otimes\dots\otimes a_n\mapsto X_{\lambda_1,\dots,\lambda_n}(a_1,\dots,a_n)
$$
which satisfy
\begin{enumerate}[(a)]
\item
sesquilinearity:
\begin{equation}\label{eq:sesquimaps}
\begin{split}
&X_{\lambda_1,\dots,\lambda_n}(a_1,\dots,a_{i-1},\partial a_i,a_{i+1},\dots,a_n)
\\
&=-\lambda_i X_{\lambda_1,\dots,\lambda_n}(a_1,\dots,a_{i-1},a_i,a_{i+1},\dots,a_n)\,,
\end{split}
\end{equation}
for all $i=1,\dots n$ and for every $a_1,\dots, a_n\in\mc V$;
\item
Leibniz rules:
\begin{equation}
\begin{split}\label{eq:Leibnizmaps}
&X_{\lambda_1,\dots,\lambda_n}(a_1,\dots,a_{i-1},b c,a_{i+1},\dots,a_n)
\\
&=(|_{x=\partial} b)\star_iX_{\lambda_1,\dots,\lambda_{i-1},\lambda_i+x,\lambda_{i+1},\dots\lambda_n}
(a_1,\dots,a_{i-1},c,a_{i+1},\dots,a_n)
\\
&+X_{\lambda_1,\dots,\lambda_{i-1},\lambda_i+x,\lambda_{i+1},\dots\lambda_n}
(a_1,\dots,a_{i-1},b,a_{i+1},\dots,a_n)\star_{n+1-i}(|_{x=\partial} c)\,,
\end{split}
\end{equation}
for all $i=1,\dots, n$ and $a_i,b,c\in\mc V$. 
\end{enumerate}
An element of $\widetilde{\Gamma}^n(\mc V)$
is called a \emph{basic $n$-cochain}
and we set 
\begin{equation}\label{eq:Omegatilde}
\widetilde{\Gamma}(\mc V)=\bigoplus_{n\in\mb Z_{\geq0}}\widetilde{\Gamma}^n(\mc V)\,.
\end{equation}
to be the space of all basic cochains.
\begin{example}\label{exa:rdcd}
We have that $\widetilde{\Gamma}^1(\mc V)$ consists of 
linear maps $X:\mc V\to (\mc V\otimes\mc V)[\lambda]$ such that ($a\in \mc V$)
$$
X_\lambda(\partial a)=-\lambda X_\lambda(a)\,,
$$
and ($a,b\in \mc V$)
$$
X_\lambda(ab)
=X_{\lambda+x}(a)\star(|_{x=\partial}b)
+(|_{x=\partial}a)\star X_{\lambda+x}(b)
\,.
$$
We call such a linear map a conformal $2$-fold derivation (cf. \eqref{eq:rcder}) with respect to the inner bimodule structure \eqref{bmodinner} of $\mc V\otimes\mc V$.
This space has the structure of a $\kk[\partial]$-module given by ($a\in\mc V$)
$$
(\partial X)_\lambda(a)=(\lambda+\partial)X_\lambda(a)
\,.
$$
\end{example}
\begin{remark}\label{20250722:rem1}
Let $R\subset \mc V$ be a subset that generates $\mc V$ as a differential algebra. From sesquilinearity \eqref{eq:sesquimaps} and the Leibniz rules \eqref{eq:Leibnizmaps} it follows that any basic $n$-cochain
$X\in \widetilde{\Gamma}(\mc V)$ is completely determined by its restriction on $R^{\otimes n}$.
\end{remark}
We extend the usual associative product $\mc V^{\otimes m}\otimes\mc V^{\otimes n}\to\mc V^{\otimes(m+n-1)}$ defined by \eqref{prod-nm}
to an associative product
\begin{equation}\label{20230802:eq2}
\widetilde{\Gamma}^m(\mc V)\otimes\widetilde{\Gamma}^n(\mc V)\to\widetilde{\Gamma}^{m+n}(\mc V)
\end{equation}
by letting, for $X\in\widetilde{\Gamma}^m(\mc V)$ and $Y\in\widetilde{\Gamma}^n(\mc V)$,
$X\cdot Y\in\widetilde\Gamma^{m+n}(\mc V)$
be defined by
\begin{align*}
&(XY)_{\lambda_1,\dots,\lambda_{m+n}}
(a_1,\dots,a_{m+n})
\\
&=X_{\lambda_1,\dots,\lambda_m}(a_1,\dots,a_m)
Y_{\lambda_{m+1},\dots,\lambda_{m+n}}(a_{m+1},\dots,a_{m+n})
\,,
\end{align*}
for every $a_1,\dots,a_{m+n}\in\mc V$.
Furthermore, we consider $\widetilde{\Gamma}(\mc V)$ as a superspace, with superstructure compatible with
the $\mb Z_{\geq0}$-grading \eqref{eq:Omegatilde}, by assigning $|X|=n\mod 2$,
for $X \in\widetilde{\Gamma}^n(\mc V)$.
Hence, $\widetilde{\Gamma}(\mc V)$ becomes a Lie superalgebra, with Lie bracket given 
on homogeneous elements $X\in\widetilde{\Gamma}^m(\mc V)$ and $Y\in\widetilde{\Gamma}^n(\mc V)$
by (cf.~\eqref{Eq:GradComm})
$$
[X,Y]=XY-(-1)^{m n}YX
\,.
$$
Finally, we extend the action of $\partial$ on $\mc V$ to an even derivation on $\widetilde{\Gamma}(\mc V)$ 
by letting,
for $X\in\widetilde{\Gamma}^m(\mc V)$,
the basic cochain $\partial X\in\widetilde{\Gamma}^m(\mc V)$ be defined by
\begin{equation}\label{20230802:eq1}
(\partial X)_{\lambda_1,\dots,\lambda_m}(a_1,\dots,a_m)
=(\lambda_1+\dots+\lambda_m+\partial)X_{\lambda_1,\dots,\lambda_m}(a_1,\dots, a_m)
\,,
\end{equation}
for every $a_1,\dots,a_m\in\mc V$.
In the RHS of equation \eqref{20230802:eq1} the map
$\partial:\mc V\to\mc V$ is extended to a map $\partial:\mc V^{\otimes m}\to\mc V^{\otimes m}$ using \eqref{mfold-ext}.

Let $X\in\widetilde{\Gamma}^n(\mc V)$ be a basic $n$-cochain. Recall from \eqref{eq:XL} that for every $A\in\mc V^{\otimes m}$ and $s=1,\dots,n$ we can define a map
$X^{(s)}:\mc V^{\otimes n}\to\mc V^{\otimes(m+n+1)}[\lambda_1,\dots,\lambda_n]$. The next result describes the Leibniz rules properties satisfied by $X^{(s)}$.
\begin{lemma}\label{lem:PVAcoh-pre4b}
Let $a_1,\dots,a_n,b,c\in\mc V$, $A\in\mc V^{\otimes m}$
and $i=1,\dots,n$.
\begin{enumerate}[(a)]
\item
For $1\leq s\leq i-1$ we have
\begin{equation}
\begin{split}\label{20230808:eq1}
&X_{\lambda_1,\dots,\lambda_n}^{(s)}
(a_1,\dots,a_s\otimes A,\dots,a_{i-1},bc,a_{i+1},\dots, a_n)
\\
&=
(|_{x=\partial}b)\star_{m+i}
X_{\lambda_1,\dots,\lambda_i+x,\dots\lambda_n}^{(s)}
(a_1,\dots,a_s\otimes A,\dots,a_{i-1},c,a_{i+1},\dots, a_n)
\\
&+
X_{\lambda_1,\dots,\lambda_i+x,\dots\lambda_n}^{(s)}
(a_1,\dots,a_s\otimes A,\dots,a_{i-1},b,a_{i+1},\dots, a_n)
\star_{n+1-i}(|_{x=\partial}c)
\,.
\end{split}
\end{equation}
\item
We have
\begin{equation}
\begin{split}\label{20230808:eq2}
&X_{\lambda_1,\dots,\lambda_n}^{(i)}(a_1,\dots,a_{i-1},bc\otimes A,a_{i+1},\dots, a_n)
\\
&=
(|_{x=\partial}b)\star_{m+i}
X_{\lambda_1,\dots,\lambda_i+x,\dots\lambda_n}^{(i)}
(a_1,\dots,a_{i-1},c\otimes A,a_{i+1},\dots, a_n)
\\
&+
X_{\lambda_1,\dots,\lambda_i+x,\dots\lambda_n}^{(i)}
(a_1,\dots,a_{i-1},b\otimes A,a_{i+1},\dots, a_n)
\star_{m+n+1-i}(|_{x=\partial}c)
\,.
\end{split}
\end{equation}
\item
For $i+1\leq s\leq n$ we have
\begin{equation}
\begin{split}\label{20230808:eq3}
&X_{\lambda_1,\dots,\lambda_n}^{(s)}(a_1,\dots,a_{i-1},bc,a_{i+1},\dots,a_s\otimes A,\dots a_n)
\\
&=
(|_{x=\partial}b)\star_i
X_{\lambda_1,\dots,\lambda_i+x,\dots\lambda_n}^{(s)}
(a_1,\dots,a_{i-1},c,a_{i+1},\dots,a_s\otimes A,\dots a_n)
\\
&+
X_{\lambda_1,\dots,\lambda_i+x,\dots\lambda_n}^{(s)}
(a_1,\dots,a_{i-1},b,a_{i+1},\dots,a_s\otimes A,\dots, a_n)
\star_{m+n+1-i}(|_{x=\partial}c)
\,.
\end{split}
\end{equation}
\end{enumerate}
\end{lemma}
\begin{proof}
The proof follows from a straightforward computation using the definition of $X^{(s)}$ given by~\eqref{eq:XL}, the Leibniz rules~\eqref{eq:Leibnizmaps} and Lemma \ref{lem:PVAcoh-pre4}.
\end{proof}
The next results will be needed in the sequel.
\begin{lemma}\label{lem:PVAcoh-pre1}
Let $X\in\widetilde{\Gamma}^m(\mc V)$ and $Y\in\widetilde{\Gamma}^n(\mc V)$.
\begin{enumerate}[(i)]
\item For $1\leq s\leq m+n+1$ we have
\begin{align*}
&(XY)_{\lambda_1,\stackrel{s}{\check{\dots}},\lambda_{m+n+1}}
(a_1,\stackrel{s}{\check{\dots}},a_{m+n+1})
\\
&=\left\{
\begin{array}{ll}
X_{\lambda_1,\stackrel{s}{\check{\dots}},\lambda_{m+1}}
(a_1,\stackrel{s}{\check{\dots}},a_{m+1})
Y_{\lambda_{m+2},\dots,\lambda_{m+n+1}}
(a_{m+2},\dots,a_{m+n+1})
\,, & 1\leq s\leq m\,,
\\
X_{\lambda_1,\dots,\lambda_{m}}
(a_1,\dots,a_{m})
Y_{\lambda_{m+2},\dots,\lambda_{m+n+1}}
(a_{m+2},\dots,a_{m+n+1})\,,
& s=m+1\,,
\\
X_{\lambda_{1},\dots,\lambda_{m}}(a_{1},\dots,a_{m})
Y_{\lambda_{m+1},\stackrel{s}{\check{\dots}},\lambda_{m+n+1}}
(a_{m+1},\stackrel{s}{\check{\dots}},a_{m+n+1})
\,,
&
m+2\leq s\leq m+n+1\,.
\end{array}
\right.
\end{align*}
\item For $1\leq s\leq m+n$ we have
\begin{align*}
&(XY)_{\lambda_1,\dots,\lambda_{s-1},\mu_s,\lambda_{s+2},\dots,\lambda_{m+n+1}}^{(s)}
\\
&=\left\{
\begin{array}{ll}
X_{\lambda_1,\dots,\lambda_{s-1},\mu_s,\lambda_{s+2},\dots,\lambda_{m+1}}^{(s)}
Y_{\lambda_{m+2},\dots,\lambda_{m+n+1}}
\,, & 1\leq s\leq m\,,
\\
X_{\lambda_{1},\dots,\lambda_{m}}
Y_{\lambda_{m+1},\dots,\lambda_{s-1},\mu_s,\lambda_{s+2},\dots,\lambda_{m+n+1}}^{(s-m)}\,,
&
m+1\leq s\leq m+n\,.
\end{array}
\right.
\end{align*}
\end{enumerate}
\end{lemma}
\begin{proof}
Straightforward using the product of basic cochains
\eqref{20230802:eq2}
and equations \eqref{Eq:NotCheck1} and \eqref{eq:XL}.
\end{proof}

\begin{lemma}\label{20240724:lem1-b}
Let $(m,n)\in\mb Z_{\geq0}^2\setminus\{(0,0)\}$. For every $X\in\widetilde{\Gamma}^m(\mc V)$ and
$Y\in\widetilde{\Gamma}^n(\mc V)$ we have
\begin{equation}\label{20240724:eq3}
\begin{split}
&\mult_{(h+1,h+2)}\circ\sigma^{h+1}(XY)_{\lambda_1,\dots,\lambda_{m+n}}\circ\sigma^{-h}
\\
&=
\left\{
\begin{array}{ll}
\mult_{(m+h+1,m+h+2)}\circ\sigma^{m+h+1}(YX)_{\lambda_{\sigma^m(1)},\dots,\lambda_{\sigma^{m}(m+n)}}
\circ\sigma^{-h-m}\,,&
h=0,\dots,n-1\,,
\\
\mult_{(h+1-n,h+2-n)}\circ\sigma^{h+1-n}(YX)_{\lambda_{\sigma^{-n}(1)},\dots,\lambda_{\sigma^{-n}(m+n)}}\circ
\sigma^{-h+n}\,,&
h=n,\dots,n+m-1
\,.
\end{array}
\right.
\end{split}
\end{equation}
\end{lemma}
\begin{proof}
Equation \eqref{20240724:eq3} follows from \eqref{20240724:eq2} and the definition of the product
of $n$-cochains \eqref{20230802:eq2}.
\end{proof}

%%%
\section[Basic \MakeLowercase{d}PVA cohomology complex]{The basic double Poisson vertex algebra cohomology complex}
\label{sec:omega2}

Let $\mc V$ be a dPVA with $2$-fold $\lambda$-bracket
$\ldb-_\lambda-\rdb$
and let $\widetilde{\Gamma}(\mc V)$ be the space of basic cochains constructed in Section~\ref{sec:omega1}.
For every $X\in\widetilde{\Gamma}^n(\mc V)$ we introduce the linear map
$\delta(X):\mc V^{\otimes(n+1)}\to \mc V^{\otimes(n+2)}[\lambda_1,\dots,\lambda_{n+1}]
$
be the linear map defined by\glslink{Tildelta}{}
\begin{equation}
\begin{split}\label{eq:diff}
&\tilde\delta(X)_{\lambda_1,\dots,\lambda_{n+1}}(a_1,\dots,a_{n+1})
=
\sum_{s=1}^{n+1}(-1)^{s+1}\ldb {a_s}_{\lambda_s}
X_{\lambda_1,\stackrel{s}{\check{\dots}},\lambda_{n+1}}
(a_1,\stackrel{s}{\check{\dots}},a_{n+1})\rdb_{(s)}
\\
&
+\sum_{s=1}^{n}(-1)^{s}X_{\lambda_1,\dots,\lambda_{s-1},\lambda_s+\lambda_{s+1},\lambda_{s+2},\dots,\lambda_{n+1}}^{(s)}(a_1,\dots,a_{s-1},\ldb {a_{s}}_{\lambda_s} a_{s+1}\rdb,a_{s+2},\dots,a_{n+1})
\,.
\end{split}
\end{equation}
In the RHS of \eqref{eq:diff} we are using the 
notation~\eqref{20240805:eq1} for $D=\ldb a_s{}_{\lambda_s}-\rdb:\mc V\to\mc V^{\otimes2}[\lambda_s]$, $s=1,\dots,n+1$,
and the notations~\eqref{Eq:NotCheck},~\eqref{Eq:NotCheck1} and~\eqref{eq:XL}.
\begin{example}\label{exa:deltaP}
The linear map $\tilde\delta$ in~\eqref{eq:diff} is the ``double'' analogue of the linear map~\eqref{eq:diff-PVA}. We can check this explicitly for small values of $n$.
Let $c\in\mc V=\widetilde{\Gamma}^0(\mc V)$.
By~\eqref{eq:diff} we have
\begin{equation}\label{exa1}
\tilde\delta(c)=\ldb-_\lambda c\rdb:\mc V\to(\mc V\otimes\mc V)[\lambda]
\,.
\end{equation}
It follows from sesquilinearity~\eqref{eq:sesqui}
and the right Leibniz rule~\eqref{eq:rleibniz}
that $\tilde\delta(c)$ satisfies~\eqref{eq:sesquimaps} and~\eqref{eq:Leibnizmaps} for $n=1$, hence it lies in $\widetilde{\Gamma}^1(\mc V)$.
Furthermore, for $X\in\widetilde{\Gamma}^1(\mc V)$, 
equation~\eqref{eq:diff} reads
\begin{equation}\label{exa2}
\tilde\delta(X)_{\lambda,\mu}(a,b)
=\ldb a_\lambda X_\mu(b)\rdb_{(1)}-\ldb b_\mu X_\lambda(a)\rdb_{(2)}
-X_{\lambda+\mu}^{(1)}(\ldb a_\lambda b\rdb)
\,.
\end{equation}
We leave as an exercise to check that 
the sesquilinearity~\eqref{eq:sesquimaps} and the skewsymmetry~\eqref{eq:skew2} 
axioms imply that $\tilde\delta(X)\in\widetilde{\Gamma}^2(\mc V)$ . 
Finally, combining equations~\eqref{exa1} 
and~\eqref{exa2} we get
$$
\tilde\delta^2(c)_{\lambda,\mu}(a,b)
=\ldb a_\lambda \ldb b_\mu c\rdb\rdb_{L}
-\ldb b_\mu \ldb a_\lambda c\rdb\rdb_{R}
-\ldb \ldb a_\lambda c\rdb_{\lambda+\mu} c\rdb_L=0
\,,
$$
which vanishes by the Jacobi identity~\eqref{eq:jacobi2}.
These computations should be compared with the analogous in Example~\ref{exa:deltaP-PVA}.
\end{example}
\begin{theorem}\label{thm:PVAcoh}
Let $\mc V$ be a dPVA and let $\widetilde{\Gamma}(V)$ be the space of basic cochains.
Then~\eqref{eq:diff} gives a well defined map
$\tilde\delta:\widetilde{\Gamma}^n(\mc V)\to\widetilde{\Gamma}^{n+1}(\mc V)$, for every
$n\in\mb Z_{\geq0}$, such that $\tilde\delta^2=0$.
\end{theorem}
\begin{proof}
To prove the first part of the claim we need to show that $\tilde\delta(X)$ satisfies~\eqref{eq:sesquimaps}
and~\eqref{eq:Leibnizmaps} for every $X\in\widetilde{\Gamma}^n(\mc V)$.
The fact that $\tilde\delta(X)$ satisfies 
the sesquilinearity conditions~\eqref{eq:sesquimaps} can be checked
directly from the definition~\eqref{eq:diff} using the fact that $X$
satisfies~\eqref{eq:sesquimaps} and using the sesquilinearity properties
\eqref{eq:sesqui} of $2$-fold $\lambda$-bracket of $\mc V$. We omit the details.
Let us now prove that
$\tilde\delta(X)$ satisfies the Leibniz rule~\eqref{eq:Leibnizmaps} for every $i=1,\dots,n+1$,
and $a_i,b,c\in\mc V$.
By the definition~\eqref{eq:diff}
of $\tilde\delta(X)$, the fact that $X$ satisfies~\eqref{eq:Leibnizmaps}, and that
$X^{(s)}$ satisfies~\eqref{20230808:eq1} and~\eqref{20230808:eq3} we have
\begin{subequations}
\begin{align}
&\tilde\delta(X)_{\lambda_1,\dots,\lambda_{n+1}}(a_1,\dots,a_{i-1},bc,a_{i+1},\dots,a_{n+1})
\notag\\
\begin{split}
&=
\sum_{s=1}^{i-1}(-1)^{s+1}\ldb {a_s}_{\lambda_s}
(|_{x=\partial}b)\star_{i-1}X_{\lambda_1,\stackrel{s}{\check{\dots}},\lambda_{i-1},\lambda_i+x,\lambda_{i+1},\dots,\lambda_{n+1}}
(a_1,\stackrel{s}{\check{\dots}}
\\
&\quad
\dots,a_{i-1},c,a_{i+1},\dots a_{n+1})\rdb_{(s)}
\label{20230808:eqproof1}
\end{split}
\\
\begin{split}\label{20230808:eqproof2}
&+
\sum_{s=1}^{i-1}(-1)^{s+1}\ldb {a_s}_{\lambda_s}
X_{\lambda_1,\stackrel{s}{\check{\dots}},\lambda_{i-1},\lambda_i+x,\lambda_{i+1},\dots,\lambda_{n+1}}
(a_1,\stackrel{s}{\check{\dots}}
\\
&\quad
\dots,a_{i-1},b,a_{i+1},\dots a_{n+1})\star_{n+2-i}(|_{x=\partial}c)\rdb_{(s)}
\end{split}
\\
\begin{split}\label{20230808:eqproof3}
&+
(-1)^{i+1}\ldb {bc}_{\lambda_i}
X_{\lambda_1,\dots,\lambda_{i-1},\lambda_{i+1}\dots,\lambda_{n+1}}
(a_1,\dots
\\
&\quad
\dots,a_{i-1},a_{i+1}\dots,a_{n+1})\rdb_{(i)}
\end{split}
\\
\begin{split}\label{20230808:eqproof4}
&+
\sum_{s=i+1}^{n+1}(-1)^{s+1}\ldb {a_s}_{\lambda_s}
(|_{x=\partial}b)\star_i
X_{\lambda_1,\dots,\lambda_{i-1},\lambda_i+x,\lambda_{i+1},\stackrel{s}{\check{\dots}},\lambda_{n+1}}
(a_1,\dots
\\
&\quad
\dots,a_{i-1},c,a_{i+1},\stackrel{s}{\check{\dots}},a_{n+1})\rdb_{(s)}
\end{split}
\\
\begin{split}\label{20230808:eqproof5}
&+
\sum_{s=i+1}^{n+1}(-1)^{s+1}\ldb {a_s}_{\lambda_s}
X_{\lambda_1,\dots,\lambda_{i-1},\lambda_i+x,\lambda_{i+1},\stackrel{s}{\check{\dots}},\lambda_{n+1}}
(a_1,\dots
\\
&\quad
\dots,a_{i-1},b,a_{i+1},\stackrel{s}{\check{\dots}},a_{n+1})\star_{n+1-i}(|_{x=\partial}c)\rdb_{(s)}
\end{split}
\\
\begin{split}\label{20230808:eqproof6}
&+\sum_{s=1}^{i-2}(-1)^{s}(|_{x=\partial} b)\star_i
X_{\lambda_1,\dots,\lambda_{s-1},\lambda_s+\lambda_{s+1},\lambda_{s+2},\dots,
\lambda_i+x,\dots,\lambda_{n+1}}^{(s)}(a_1,\dots
\\
&\quad
\dots,a_{s-1},\ldb {a_{s}}_{\lambda_s} a_{s+1}\rdb,a_{s+2},\dots,a_{i-1},c,a_{i+1},\dots,a_{n+1})
\end{split}
\\
\begin{split}\label{20230808:eqproof7}
&
+\sum_{s=1}^{i-2}(-1)^{s}X_{\lambda_1,\dots,\lambda_{s-1},\lambda_s+\lambda_{s+1},\lambda_{s+2},\dots,\lambda_i+x,\dots,\lambda_{n+1}}^{(s)}(a_1,\dots
\\
&\quad
\dots,a_{s-1},\ldb {a_{s}}_{\lambda_s} a_{s+1}\rdb,a_{s+2},\dots,a_{i-1},b,a_{i+1},\dots,a_{n+1})
\star_{n+2-i}(|_{x=\partial}c)
\end{split}
\\
\begin{split}\label{20230808:eqproof8}
&
+(-1)^{i+1}X_{\lambda_1,\dots,\lambda_{i-2},\lambda_{i-1}+\lambda_{i},\lambda_{i+1},\dots,\lambda_{n+1}}^{(i-1)}(a_1,\dots
\\
&\quad
\dots,a_{i-2},\ldb {a_{i-1}}_{\lambda_{i-1}} bc\rdb,a_{i+1},\dots,a_{n+1})
\end{split}
\\
\begin{split}\label{20230808:eqproof9}
&
+(-1)^{i}X_{\lambda_1,\dots,\lambda_{i-1},\lambda_i+\lambda_{i+1},\lambda_{i+2},\dots,\lambda_{n+1}}^{(i)}(a_1,\dots
\\
&\quad
\dots,a_{i-1},\ldb bc_{\lambda_i} a_{i+1}\rdb,a_{i+2},\dots,a_{n+1})
\end{split}
\\
\begin{split}\label{20230808:eqproof10}
&
+\sum_{s=i+1}^{n}(-1)^{s}
(|_{x=\partial}b)\star_i
X_{\lambda_1,\dots,\lambda_i+x,\dots,\lambda_{s-1},\lambda_s+\lambda_{s+1},\lambda_{s+2},\dots,\lambda_{n+1}}^{(s)}
(a_1,\dots
\\
&\quad
\dots,a_{i-1},c,a_{i+1},\dots,a_{s-1},\ldb {a_{s}}_{\lambda_s} a_{s+1}\rdb,a_{s+2},\dots,a_{n+1})
\end{split}
\\
\begin{split}\label{20230808:eqproof11}
&
+\sum_{s=i+1}^{n}(-1)^{s}
X_{\lambda_1,\dots,\lambda_i+x,\dots,\lambda_{s-1},\lambda_s+\lambda_{s+1},\lambda_{s+2},\dots,\lambda_{n+1}}^{(s)}
(a_1,\dots
\\
&\quad
\dots,a_{i-1},b,a_{i+1},\dots,a_{s-1},\ldb {a_{s}}_{\lambda_s} a_{s+1}\rdb,a_{s+2},\dots,a_{n+1})
\star_{n+2-i}(|_{x=\partial}c)\,.
\end{split}
\end{align}
\end{subequations}
Using Lemma~\ref{lem:PVAcoh-pre2} we have
\begin{equation}
\begin{split}
\label{20230808:eqproof1a}
&\eqref{20230808:eqproof1}+\eqref{20230808:eqproof4}
\\
&=
\sum_{s=1}^{i-1}(-1)^{s+1}(|_{x=\partial}b)\star_i\ldb {a_s}_{\lambda_s}
X_{\lambda_1,\stackrel{s}{\check{\dots}},\lambda_i+x,\dots,\lambda_{n+1}}
(a_1,\stackrel{s}{\check{\dots}}
\\
&\quad\quad\quad\quad\quad\quad\quad\quad\quad\quad\quad\quad\quad
\dots,a_{i-1},c,a_{i+1},\dots a_{n+1})\rdb_{(s)}
\\
&+\sum_{s=i+1}^{n+1}(-1)^{s+1}(|_{x=\partial}b)\star_i\ldb {a_s}_{\lambda_s}
X_{\lambda_1,\dots,\lambda_i+x,\stackrel{s}{\check{\dots}},\lambda_{n+1}}
(a_1,\dots
\\
&\quad\quad\quad\quad\quad\quad\quad\quad\quad\quad\quad\quad\quad
\dots,a_{i-1},c,a_{i+1},\stackrel{s}{\check{\dots}},a_{n+1})\rdb_{(s)}
\\
&
+(-1)^i\mult_{(i+2,i+3)}
X_{\lambda_1,\dots,\lambda_{i-1},\lambda_{i}+\lambda_{i+1},\lambda_{i+2},\dots,\lambda_{n+1}}^{(i)}
(a_1,\dots
\\
&\quad\quad\quad\quad\quad\quad\quad\quad\quad\quad\quad\quad\quad
\dots,a_{i-1}, c\otimes \ldb {a_{i+1}}_{\lambda_{i+1}} b\rdb,a_{i+2},\dots,a_{n+1})
\,.
\end{split}
\end{equation}
Similarly, using again Lemma~\ref{lem:PVAcoh-pre2}, we have
\begin{equation}
\begin{split}
\label{20230808:eqproof2a}
&\eqref{20230808:eqproof2}+\eqref{20230808:eqproof5}
\\
&=
\sum_{s=1}^{i-1}(-1)^{s+1}\ldb {a_s}_{\lambda_s}
X_{\lambda_1,\stackrel{s}{\check{\dots}},\lambda_i+x,\dots,\lambda_{n+1}}
(a_1,\stackrel{s}{\check{\dots}}
\\
&\quad\quad\quad\quad\quad\quad\quad\quad\quad\quad\quad\quad
\dots,a_{i-1},b,a_{i+1},\dots a_{n+1})\rdb_{(s)}
\star_{n+2-i}(|_{x=\partial}c)
\\
&+\sum_{s=i+1}^{n+1}(-1)^{s+1}\ldb {a_s}_{\lambda_s}
X_{\lambda_1,\dots,\lambda_i+x,\stackrel{s}{\check{\dots}},\lambda_{n+1}}
(a_1,\dots
\\
&\quad\quad\quad\quad\quad\quad\quad\quad\quad\quad\quad\quad
\dots,a_{i-1},b,a_{i+1},\stackrel{s}{\check{\dots}},a_{n+1})\rdb_{(s)}
\star_{n+2-i}(|_{x=\partial}c)
\\
&+(-1)^i\mult_{(i-1,i)}
X_{\lambda_1,\dots,\lambda_{i-2},\lambda_{i-1}+\lambda_i,\lambda_{i+1},\dots,\lambda_{n+1}}^{(i-1)}
(a_1,\dots
\\
&\quad\quad\quad\quad\quad\quad\quad\quad\quad\quad\quad\quad
\dots,a_{i-2}, b\otimes \ldb {a_{i-1}}_{\lambda_{i-1}} c\rdb,a_{i+1},\dots,a_{n+1})
\,.
\end{split}
\end{equation}
Furthermore, using Lemma~\ref{lem:PVAcoh-pre3} we have
\begin{align}
\begin{split}\label{20230808:eqproof3a}
&\eqref{20230808:eqproof3}=
(-1)^{i+1}(|_{x=\partial}b)\star_i\ldb {c}_{\lambda_i+x}
X_{\lambda_1,\dots,\lambda_{i-1},\lambda_{i+1}\dots,\lambda_{n+1}}
(a_1,\dots
\\
&\quad\quad\quad\quad\quad\quad\quad\quad\quad\quad\quad\quad\quad
\dots,a_{i-1},a_{i+1}\dots,a_{n+1})\rdb_{(i)}
\end{split}
\\
\begin{split}\label{20230808:eqproof3b}
&+(-1)^{i+1}\ldb {b}_{\lambda_i+x}
X_{\lambda_1,\dots,\lambda_{i-1},\lambda_{i+1}\dots,\lambda_{n+1}}
(a_1,\dots
\\
&\quad\quad\quad\quad\quad\quad\quad\quad\quad\quad\quad\quad\quad
\dots,a_{i-1},a_{i+1}\dots,a_{n+1})\rdb_{(i)}
\star_{n+2-i}(|_{x=\partial}c)\,.
\end{split}
\end{align}
Next, we rewrite~\eqref{20230808:eqproof8} as
\begin{align}
&\eqref{20230808:eqproof8}
=(-1)^{i+1}X_{\lambda_1,\dots,\lambda_{i-2},\lambda_{i-1}+\lambda_{i},\lambda_{i+1},\dots,\lambda_{n+1}}^{(i-1)}(a_1,\dots
\notag
\\
&\quad\quad\quad\quad\quad\quad\quad\quad\quad
\dots,a_{i-2},\ldb {a_{i-1}}_{\lambda_{i-1}} b\rdb c,a_{i+1},\dots,a_{n+1})
\notag
\\
&+(-1)^{i+1}X_{\lambda_1,\dots,\lambda_{i-2},\lambda_{i-1}+\lambda_{i},\lambda_{i+1},\dots,\lambda_{n+1}}^{(i-1)}(a_1,\dots
\notag
\\
&\quad\quad\quad\quad\quad\quad\quad\quad\quad
\dots,a_{i-2},b \ldb {a_{i-1}}_{\lambda_{i-1}} c\rdb,a_{i+1},\dots,a_{n+1})
\notag
\\
\begin{split}\label{20230808:eqproof8a}
&
=(-1)^{i+1}X_{\lambda_1,\dots,\lambda_{i-2},\lambda_{i-1}+\lambda_{i}+x,\lambda_{i+1},\dots,\lambda_{n+1}}^{(i-1)}(a_1,\dots
\\
&\quad\quad\quad\quad\quad\quad\quad\quad\quad
\dots,a_{i-2},\ldb {a_{i-1}}_{\lambda_{i-1}} b\rdb,a_{i+1},\dots,a_{n+1})
\star_{n+2-i}(|_{x=\partial}c)
\end{split}
\\
\begin{split}\label{20230808:eqproof8b}
&
+(-1)^{i+1}(|_{x=\partial}b)\star_iX_{\lambda_1,\dots,\lambda_{i-2},\lambda_{i-1}+\lambda_{i}+x,\lambda_{i+1},\dots,\lambda_{n+1}}^{(i-1)}(a_1,\dots
\\
&\quad\quad\quad\quad\quad\quad\quad\quad\quad
\dots,a_{i-2}, \ldb {a_{i-1}}_{\lambda_{i-1}} c\rdb,a_{i+1},\dots,a_{n+1})
\end{split}
\\
\begin{split}\label{20230808:eqproof8c}
&-(-1)^i\mult_{(i-1,i)}
X_{\lambda_1,\dots,\lambda_{i-2},\lambda_{i-1}+\lambda_i,\lambda_{i+1},\dots,\lambda_{n+1}}^{(i-1)}
(a_1,\dots
\\
&\quad\quad\quad\quad\quad\quad\quad\quad\quad
\dots,a_{i-2}, b\otimes \ldb {a_{i-1}}_{\lambda_{i-1}} c\rdb,a_{i+1},\dots,a_{n+1})
\,.
\end{split}
\end{align}
In the first equality we used the Leibniz rule~\eqref{eq:lleibniz} and the linearity of $X^{(s)}$
and in the second equality we used equation~\eqref{20230808:eq2} and the identity
$$
(A\otimes_{n+2-i}\beta)\star_{n+3-i}\alpha
=\mult_{(i-1,i)}(A\otimes_{n+2-i}(\alpha\otimes \beta))
\,,
$$
which can be easily verified to hold for every $A\in\mc V^{\otimes (n+1)}$ and $\alpha,\beta\in\mc V$.
Finally, we rewrite~\eqref{20230808:eqproof9} as
\begin{align*}
&\eqref{20230808:eqproof9}
=
(-1)^{i}X_{\lambda_1,\dots,\lambda_{i-1},\lambda_i+\lambda_{i+1},\lambda_{i+2},\dots,\lambda_{n+1}}^{(i)}(a_1,\dots
\\
&\quad\quad\quad\quad\quad\quad\quad
\dots,a_{i-1},(|_{x=\partial}b)\star\ldb c_{\lambda_i+x} a_{i+1}\rdb,a_{i+2},\dots,a_{n+1})
\\
&
+(-1)^{i}X_{\lambda_1,\dots,\lambda_{i-1},\lambda_i+\lambda_{i+1},\lambda_{i+2},\dots,\lambda_{n+1}}^{(i)}(a_1,\dots
\\
&\quad\quad\quad\quad\quad\quad\quad
\dots,a_{i-1},\ldb b_{\lambda_i+x} a_{i+1}\rdb\star(|_{x=\partial}c),a_{i+2},\dots,a_{n+1})
\\
&=(-1)^{i}(|_{x=\partial}b)\star_iX_{\lambda_1,\dots,\lambda_{i-1},\lambda_i+x+\lambda_{i+1},\lambda_{i+2},\dots,\lambda_{n+1}}^{(i)}(a_1,\dots
\\
&\quad\quad\quad\quad\quad\quad\quad
\dots,a_{i-1},\ldb c_{\lambda_i+x} a_{i+1}\rdb,a_{i+2},\dots,a_{n+1})
\\
&+(-1)^{i}X_{\lambda_1,\dots,\lambda_{i-1},\lambda_i+x+\lambda_{i+1},\lambda_{i+2},\dots,\lambda_{n+1}}^{(i)}(a_1,\dots
\\
&\quad\quad\quad\quad\quad\quad\quad
\dots,a_{i-1},\ldb b_{\lambda_i+x} a_{i+1}\rdb,a_{i+2},\dots,a_{n+1})\star_{n+2-i}(|_{x=\partial}c)
\\
&+(-1)^{i}\mult_{(i+2,i+3)}\Big(
X_{\lambda_1,\dots,\lambda_{i-1},\lambda_i+y+\lambda_{i+1},\lambda_{i+2},\dots,\lambda_{n+1}}(a_1,\dots
\\
&\quad\quad\quad\quad\quad\quad\quad
\dots,a_{i-1},(|_{x=\partial}c),a_{i+2},\dots,a_{n+1})
\otimes_{n+1-i}
\left(|_{y=\partial}\ldb b_{\lambda_i+x} a_{i+1}\rdb^\sigma\right)
\Big)
\,.
\end{align*}
In the first equality we used the Leibniz rule~\eqref{eq:rleibniz} and the linearity of $X^{(s)}$ and in the
second equality we used equation~\eqref{20230808:eq2} and the identities
$$
A\otimes_{n+1-i}bB=b\star_i(A\otimes_{n+1-i}B)
\,,
\,\,
\alpha\star_{i+1}(A\otimes_{n+1-i}\beta)=\mult_{(i+2,i+3)}(A\otimes_{n+1-i}(\beta\otimes\alpha))
\,$$
which can be easily checked to hold for every
$A\in\mc V^{\otimes(n+1)}$, $B\in\mc V^{\otimes m}$ and $b,\alpha,\beta\in\mc V$.
Using the sesquilinearity property~\eqref{eq:sesquimaps} of $X$ and the skew-symmetry assumption
\eqref{eq:skew2} on the 2-fold $\lambda$-bracket of the dPVA $\mc V$, we have
\begin{align*}
&X_{\lambda_1,\dots,\lambda_{i-1},\lambda_i+y+\lambda_{i+1},\lambda_{i+2},\dots,\lambda_{n+1}}(a_1,\dots
\\
&\quad\quad\quad\quad\quad\quad
\dots,a_{i-1},(|_{x=\partial}c),a_{i+2},\dots,a_{n+1})
\otimes_{n+1-i}
\left(|_{y=\partial}\ldb b_{\lambda_i+x} a_{i+1}\rdb^\sigma\right)
\\
&=
X_{\lambda_1,\dots,\lambda_{i-1},\lambda_i+y+\lambda_{i+1},\lambda_{i+2},\dots,\lambda_{n+1}}(a_1,\dots
\\
&\quad\quad\quad\quad\quad\quad
\dots,a_{i-1},c,a_{i+2},\dots,a_{n+1})
\otimes_{n+1-i}
\left(|_{y=\partial}\ldb b_{-\lambda_{i+1}-y} a_{i+1}\rdb^\sigma\right)
\\
&=
-X_{\lambda_1,\dots,\lambda_{i-1},\lambda_i+y+\lambda_{i+1},\lambda_{i+2},\dots,\lambda_{n+1}}(a_1,\dots
\\
&\quad\quad\quad\quad\quad\quad
\dots,a_{i-1},c,a_{i+2},\dots,a_{n+1})
\otimes_{n+1-i}
\left(|_{y=\partial}\ldb {a_{i+1}}_{\lambda_{i+1}} b\rdb\right)
\\
&=
-X_{\lambda_1,\dots,\lambda_{i-1},\lambda_i+\lambda_{i+1},\lambda_{i+2},\dots,\lambda_{n+1}}^{(i)}
(a_1,\dots
\\
&\quad\quad\quad\quad\quad\quad
\dots,a_{i-1},c\otimes \ldb {a_{i+1}}_{\lambda_{i+1}} b\rdb,a_{i+2},\dots,a_{n+1})
\,,
\end{align*}
where in the last identity we used the definition of $X^{(s)}$ given by~\eqref{eq:XL}. Hence, we get
\begin{align}
\begin{split}\label{20230808:eqproof9a}
&\eqref{20230808:eqproof9}
=(-1)^{i}(|_{x=\partial}b)\star_iX_{\lambda_1,\dots,\lambda_{i-1},\lambda_i+x+\lambda_{i+1},\lambda_{i+2},\dots,\lambda_{n+1}}^{(i)}(a_1,\dots
\\
&\quad\quad\quad\quad\quad\quad\quad\quad\quad
\dots,a_{i-1},\ldb c_{\lambda_i+x} a_{i+1}\rdb,a_{i+2},\dots,a_{n+1})
\end{split}
\\
\begin{split}\label{20230808:eqproof9b}
&+(-1)^{i}X_{\lambda_1,\dots,\lambda_{i-1},\lambda_i+x+\lambda_{i+1},\lambda_{i+2},\dots,\lambda_{n+1}}^{(i)}(a_1,\dots
\\
&\quad\quad\quad\quad\quad\quad\quad\quad\quad
\dots,a_{i-1},\ldb b_{\lambda_i+x} a_{i+1}\rdb,a_{i+2},\dots,a_{n+1})\star_{n+2-i}(|_{x=\partial}c)
\end{split}
\\
\begin{split}\label{20230808:eqproof9c}
&
-(-1)^{i}\mult_{(i+2,i+3)}
X_{\lambda_1,\dots,\lambda_{i-1},\lambda_i+\lambda_{i+1},\lambda_{i+2},\dots,\lambda_{n+1}}^{(i)}
(a_1,\dots
\\
&\quad\quad\quad\quad\quad\quad\quad\quad\quad
\dots,a_{i-1},c\otimes \ldb {a_{i+1}}_{\lambda_{i+1}} b\rdb,a_{i+2},\dots,a_{n+1})
\,.
\end{split}
\end{align}
Recalling the definition of $\tilde\delta$ given by~\eqref{eq:diff}, we get
\begin{align*}
&
\eqref{20230808:eqproof1a}+\eqref{20230808:eqproof9c}
+\eqref{20230808:eqproof3a}+\eqref{20230808:eqproof6}
+\eqref{20230808:eqproof8b}+\eqref{20230808:eqproof9a}
+\eqref{20230808:eqproof10}
\\
&=(|_{x=\partial}b)\star_i\tilde\delta(X)_{\lambda_1,\dots,\lambda_i+x,\dots,\lambda_{n+1}}(a_1,\dots,a_{i-1},c,a_{i+1},\dots,a_{n+1})
\,,
\end{align*}
and
\begin{align*}
&
\eqref{20230808:eqproof2a}+\eqref{20230808:eqproof8c}
+\eqref{20230808:eqproof3b}+\eqref{20230808:eqproof7}
+\eqref{20230808:eqproof8a}+\eqref{20230808:eqproof9b}
+\eqref{20230808:eqproof11}
\\
&=\tilde\delta(X)_{\lambda_1,\dots,\lambda_i+x,\dots,\lambda_{n+1}}(a_1,\dots,a_{i-1},b,a_{i+1},\dots,a_{n+1})
\star_{n+2-i}(|_{x=\partial}c)\,,
\end{align*}
thus proving that $\tilde\delta(X)$ satisfies the Leibniz rule
\begin{align*}
&\tilde\delta(X)_{\lambda_1,\dots,\lambda_{n+1}}(a_1,\dots,a_{i-1},bc,a_{i+1},\dots,a_{n+1})
\\
&=(|_{x=\partial}b)\star_i\tilde\delta(X)_{\lambda_1,\dots,\lambda_i+x,\dots,\lambda_{n+1}}(a_1,\dots,a_{i-1},c,a_{i+1},\dots,a_{n+1})
\\
&+\tilde\delta(X)_{\lambda_1,\dots,\lambda_i+x,\dots,\lambda_{n+1}}(a_1,\dots,a_{i-1},b,a_{i+1},\dots,a_{n+1})
\star_{n+2-i}(|_{x=\partial}c)\,,
\end{align*}
for every $i=1,\dots,n+1$ and $a_i,b,c\in\mc V$. This shows that $\tilde\delta_H(X)\in\widetilde{\Gamma}^{n+1}(\mc V)$.

Next, let us prove that $\tilde\delta^2(X)=0$, for every $X\in\widetilde{\Gamma}^n(\mc V)$. From the definition of 
$\tilde\delta$ given by~\eqref{eq:diff} we have
\begin{equation}
\begin{split}\label{eq:dPsquare}
&\tilde\delta^2 (X)_{\lambda_1,\dots,\lambda_{n+2}}(a_1,\dots,a_{n+2})
\\
&=
\sum_{t=1}^{n+2}(-1)^{t+1}\ldb {a_t}_{\lambda_t}
\tilde\delta (X)_{\lambda_1,\stackrel{t}{\check{\dots}},\lambda_{n+2}}
(a_1,\stackrel{t}{\check{\dots}},a_{n+2})\rdb_{(t)}
\\
&
+\sum_{t=1}^{n+1}(-1)^{t}\tilde\delta(X)_{\lambda_1,\dots,\lambda_{t-1},\lambda_t+\lambda_{t+1},\lambda_{t+2},\dots,\lambda_{n+2}}^{(t)}
(a_1,\dots
\\
&\quad\quad\quad\quad\quad\quad\quad\quad\dots,a_{t-1},\ldb {a_{t}}_{\lambda_t} a_{t+1}\rdb,a_{t+2},\dots,a_{n+2})
\,.
\end{split}
\end{equation}
From~\eqref{eq:diff} we have
\begin{equation}
\begin{split}\label{eq:dPsquare1}
&\tilde\delta(X)_{\lambda_1,\stackrel{t}{\check{\dots}},\lambda_{n+2}}
(a_1,\stackrel{t}{\check{\dots}},a_{n+2})
\\
&=\sum_{s=1}^{t-1}(-1)^{s+1}\ldb {a_{s}}_{\lambda_s}
X_{\lambda_1,\stackrel{s}{\check{\dots}}\stackrel{t}{\check{\dots}},\lambda_{n+2}}
(a_1,\stackrel{s}{\check{\dots}}\stackrel{t}{\check{\dots}},a_{n+2})\rdb_{(s)}
\\
&
+\sum_{s=t+1}^{n+2}(-1)^{s}\ldb {a_s}_{\lambda_s}
X_{\lambda_1,\stackrel{t}{\check{\dots}}\stackrel{s}{\check{\dots}},\lambda_{n+2}}
(a_1,\stackrel{t}{\check{\dots}}\stackrel{s}{\check{\dots}},a_{n+2})\rdb_{(s-1)}
\\
&
+\sum_{s=1}^{t-2}(-1)^{s}X_{\lambda_1,\dots,\lambda_{s-1},\lambda_s+\lambda_{s+1},\lambda_{s+2},
\stackrel{t-1}{\check{\dots}},\lambda_{n+2}}^{(s)}(a_1,\dots
\\
&\quad\quad\quad\quad\quad\quad\quad\quad\quad\quad\quad
\dots,a_{s-1},\ldb {a_{s}}_{\lambda_s} a_{s+1}\rdb,a_{s+2},\stackrel{t-1}{\check{\dots}},a_{n+2})
\\
&
+(-1)^{t-1}X_{\lambda_1,\dots,\lambda_{t-2},\lambda_{t-1}+\lambda_{t+1},\lambda_{t+2},
\dots,\lambda_{n+2}}^{(t-1)}(a_1,\dots
\\
&\quad\quad\quad\quad\quad\quad\quad\quad\quad\quad\quad
\dots,a_{t-2},\ldb {a_{t-1}}_{\lambda_{t-1}} a_{t+1}\rdb,a_{t+2},\dots,a_{n+2})
\\
&
+\sum_{s=t+1}^{n+1}(-1)^{s-1}X_{\lambda_1,\stackrel{t}{\check{\dots}},\lambda_{s-1},\lambda_s+\lambda_{s+1},\lambda_{s+2},\dots,\lambda_{n+2}}^{(s-1)}
(a_1,\stackrel{t}{\check{\dots}}
\\
&\quad\quad\quad\quad\quad\quad\quad\quad\quad\quad\quad
\dots,a_{s-1},\ldb {a_{s}}_{\lambda_s} a_{s+1}\rdb,a_{s+2},\dots,a_{n+2})
\,.
\end{split}
\end{equation}
From~\eqref{eq:diff}, using~\eqref{eq:XL}, we have
\begin{equation}
\begin{split}\label{eq:dPsquare2pre}
&\tilde\delta(X)_{\lambda_1,\dots,\lambda_{t-1},\lambda_t+\lambda_{t+1},\lambda_{t+2},\dots,\lambda_{n+2}}^{(t)}
(a_1,\dots,a_{t-1},\ldb {a_{t}}_{\lambda_t} a_{t+1}\rdb,a_{t+2},\dots,a_{n+2})
\\
&=
\tilde\delta(X)_{\lambda_1,\dots,\lambda_{t-1},\lambda_t+\lambda_{t+1}+x,\lambda_{t+2},\dots,\lambda_{n+2}}
(a_1,\dots
\\
&\quad\quad\quad
\dots,a_{t-1},\ldb {a_{t}}_{\lambda_t} a_{t+1}\rdb',a_{t+2},\dots,a_{n+2})
\otimes_{n+2-t}(|_{x=\partial}\ldb {a_{t}}_{\lambda_t} a_{t+1}\rdb'')
\\
&=
\bigg[
\sum_{s=1}^{t-1}(-1)^{s+1}
\ldb {a_s}_{\lambda_s}X_{\lambda_1,\stackrel{s}{\check{\dots}},\lambda_{t-1},\lambda_t+\lambda_{t+1}+x,\lambda_{t+2},\dots,\lambda_{n+2}}
(a_1,\stackrel{s}{\check{\dots}}
\\
&\quad\quad\quad
\dots,a_{t-1},\ldb {a_{t}}_{\lambda_t} a_{t+1}\rdb',a_{t+2},\dots,a_{n+2})
\rdb_{(s)}
\\
&
+(-1)^{t+1}\ldb \ldb {a_{t}}_{\lambda_t} a_{t+1}\rdb'_{\lambda_t+\lambda_{t+1}+x}
X_{\lambda_1,\dots,\lambda_{t-1},\lambda_{t+2},\dots\lambda_{n+2}}
(a_1,\dots
\\
&\quad\quad\quad
\dots,a_{t-1},a_{t+2},\dots,a_{n+2})
\rdb_{(t)}
\\
&+\sum_{s=t+2}^{n+2}(-1)^{s}
\ldb {a_s}_{\lambda_s}X_{\lambda_1,\dots,\lambda_{t-1},\lambda_t+\lambda_{t+1}+x,\lambda_{t+2},
\stackrel{s-1}{\check{\dots}},\lambda_{n+2}}
(a_1,\dots
\\
&\quad\quad\quad
\dots,a_{t-1},\ldb {a_{t}}_{\lambda_t} a_{t+1}\rdb',a_{t+2},\stackrel{s-1}{\check{\dots}},a_{n+2})
\rdb_{(s-1)}
\\
&
+\sum_{s=1}^{t-2}(-1)^{s}
X_{\lambda_1,\dots,\lambda_s+\lambda_{s+1},\dots,
\lambda_t+\lambda_{t+1}+x,
\dots,\lambda_{n+2}}^{(s)}
(a_1,\dots
\\
&\quad\quad\quad
\dots,a_{s-1}\ldb {a_{s}}_{\lambda_s} a_{s+1}\rdb,a_{s+1},\dots,a_{t-1},\ldb {a_{t}}_{\lambda_t} a_{t+1}\rdb',a_{t+2},\dots,a_{n+2})
\\
&
+(-1)^{t-1} X_{\lambda_1,\dots,\lambda_{t-2},
\lambda_{t-1}+\lambda_t+\lambda_{t+1}+x,\lambda_{t+2},
\dots,\lambda_{n+2}}^{(t-1)}
(a_1,\dots
\\
&\quad\quad\quad
\dots,a_{t-2}, \ldb {a_{t-1}}_{\lambda_{t-1}}\ldb {a_{t}}_{\lambda_t} a_{t+1}\rdb'\rdb,a_{t+2},\dots,a_{n+2})
\\
&
+(-1)^{t} X_{\lambda_1,\dots,\lambda_{t-1},
\lambda_{t}+\lambda_{t+1}+\lambda_{t+2}+x,\lambda_{t+3},
\dots,\lambda_{n+2}}^{(t)}
(a_1,\dots
\\
&\quad\quad\quad
\dots,a_{t-1}, \ldb \ldb {a_{t}}_{\lambda_{t}} a_{t+1}\rdb'_{\lambda_t+\lambda_{t+1}+x} a_{t+2}\rdb,a_{t+3},\dots,a_{n+2})
\\
&
+\sum_{s=t+2}^{n+1}(-1)^{s-1}
X_{\lambda_1,\dots,\lambda_t+\lambda_{t+1}+x,\dots,
\lambda_s+\lambda_{s+1},
\dots,\lambda_{n+2}}^{(s-1)}
(a_1,\dots
\\
&\quad\quad\quad
\dots,a_{t-1},\ldb {a_{t}}_{\lambda_t} a_{t+1}\rdb',a_{t+2},\dots,a_{s-1},\ldb {a_{s}}_{\lambda_s} a_{s+1}\rdb,a_{s+2},\dots,a_{n+2})\bigg]
\\
&
\otimes_{n+2-t}(|_{x=\partial}\ldb {a_{t}}_{\lambda_t} a_{t+1}\rdb'')
\,.
\end{split}
\end{equation}
Using the first equality in~\eqref{20230811:eq1} with $i=n+2-t$ and the definition of $X^{(s)}$ given by
\eqref{eq:XL} we get the identity for $1\leq s \leq t-1$
\begin{equation}
\begin{split}\label{20230811:eq2}
&\ldb {a_s}_{\lambda_s}X_{\lambda_1,\stackrel{s}{\check{\dots}},\lambda_{t-1},\lambda_t+\lambda_{t+1}+x,\lambda_{t+2},\dots,\lambda_{n+2}}
(a_1,\stackrel{s}{\check{\dots}}
\\
&\quad\quad
\dots,a_{t-1},\ldb {a_{t}}_{\lambda_t} a_{t+1}\rdb',a_{t+2},\dots,a_{n+2})
\rdb_{(s)}
\otimes_{n+2-t}(|_{x=\partial}\ldb {a_{t}}_{\lambda_t} a_{t+1}\rdb'')
\\
&=
\ldb {a_s}_{\lambda_s}X_{\lambda_1,\stackrel{s}{\check{\dots}},\lambda_{t-1},\lambda_t+\lambda_{t+1},\lambda_{t+2},\dots,\lambda_{n+2}}^{(t-1)}
(a_1,\stackrel{s}{\check{\dots}}
\\
&\quad\quad
\dots,a_{t-1},\ldb {a_{t}}_{\lambda_t} a_{t+1}\rdb,a_{t+2},\dots,a_{n+2})
\rdb_{(s)}
\,.
\end{split}
\end{equation}
Similarly, using the third equality in~\eqref{20230811:eq1} with $i=n+1-t$
and the definition of $X^{(s)}$
we get the identity  for $t+2\leq s \leq n+2$
\begin{equation}
\begin{split}\label{20230811:eq3}
&\ldb {a_s}_{\lambda_s}X_{\lambda_1,\dots,\lambda_{t-1},\lambda_t+\lambda_{t+1}+x,\lambda_{t+2},
\stackrel{s-1}{\check{\dots}},\lambda_{n+2}}
(a_1,\dots
\\
&\quad
\dots,a_{t-1},\ldb {a_{t}}_{\lambda_t} a_{t+1}\rdb',a_{t+2},\stackrel{s-1}{\check{\dots}},a_{n+2})
\rdb_{(s-1)}
\otimes_{n+2-t}(|_{x=\partial}\ldb {a_{t}}_{\lambda_t} a_{t+1}\rdb'')
\\
&
=
\ldb {a_s}_{\lambda_s}X_{\lambda_1,\dots,\lambda_{t-1},\lambda_t+\lambda_{t+1},\lambda_{t+2},
\stackrel{s-1}{\check{\dots}},\lambda_{n+2}}^{(t)}
(a_1,\dots
\\
&\quad
\dots,a_{t-1},\ldb {a_{t}}_{\lambda_t} a_{t+1}\rdb,a_{t+2},\stackrel{s-1}{\check{\dots}},a_{n+2})
\rdb_{(s)}
\,.
\end{split}
\end{equation}
Next, using the identity
$$
\ldb a_{\lambda +x} A\rdb_{(t)}\otimes_{n+1-t} (|_{x=\partial}b)
=\ldb a\otimes b_\lambda A\rdb_{L,(t)}
\,,
$$
which holds for every $a,b\in\mc V$ and $A\in\mc V^{\otimes n}$, we get
\begin{equation}
\begin{split}\label{20230811:eq4}
&\ldb \ldb {a_{t}}_{\lambda_t} a_{t+1}\rdb'_{\lambda_t+\lambda_{t+1}+x}
X_{\lambda_1,\dots,\lambda_{t-1},\lambda_{t+2},\dots\lambda_{n+2}}
(a_1,\dots
\\
&\quad\quad\quad\quad\quad\quad\quad\quad
\dots,a_{t-1},a_{t+2},\dots,a_{n+2})
\rdb_{(t)}
\otimes_{n+2-t}(|_{x=\partial}\ldb {a_{t}}_{\lambda_t} a_{t+1}\rdb'')
\\
&
=\ldb \ldb {a_{t}}_{\lambda_t} a_{t+1}\rdb_{\lambda_t+\lambda_{t+1}}
X_{\lambda_1,\dots,\lambda_{t-1},\lambda_{t+2},\dots\lambda_{n+2}}
(a_1,\dots
\\
&\quad\quad\quad\quad\quad\quad\quad\quad
\dots,a_{t-1},a_{t+2},\dots,a_{n+2})
\rdb_{L,(t)}
\,.
\end{split}
\end{equation}
By equation~\eqref{20230811:eq1a} we get
\begin{equation}
\begin{split}\label{20230811:eq5}
&X_{\lambda_1,\dots,\lambda_{t-2},
\lambda_{t-1}+\lambda_t+\lambda_{t+1}+x,\lambda_{t+2},
\dots,\lambda_{n+2}}^{(t-1)}
(a_1,\dots
\\
&
\dots,a_{t-2}, \ldb {a_{t-1}}_{\lambda_{t-1}}\ldb {a_{t}}_{\lambda_t} a_{t+1}\rdb'\rdb,a_{t+2},\dots,a_{n+2})
\otimes_{n+2-t}(|_{x=\partial}\ldb {a_{t}}_{\lambda_t} a_{t+1}\rdb'')
\\
&=X_{\lambda_1,\dots,\lambda_{t-2},
\lambda_{t-1}+\lambda_t+\lambda_{t+1},\lambda_{t+2},
\dots,\lambda_{n+2}}^{(t-1)}
(a_1,\dots
\\
&
\dots,a_{t-2}, \ldb {a_{t-1}}_{\lambda_{t-1}}\ldb {a_{t}}_{\lambda_t} a_{t+1}\rdb\rdb_L,a_{t+2},\dots,a_{n+2})
\end{split}
\end{equation}
and
%by equation~\eqref{20230811:eq1a} we get
\begin{equation}
\begin{split}\label{20230811:eq6}
&X_{\lambda_1,\dots,\lambda_{t-1},
\lambda_{t}+\lambda_{t+1}+\lambda_{t+2}+x,\lambda_{t+3},
\dots,\lambda_{n+2}}^{(t)}
(a_1,\dots
\\
&
\dots,a_{t-1}, \ldb \ldb {a_{t}}_{\lambda_{t}} a_{t+1}\rdb'_{\lambda_t+\lambda_{t+1}+x} a_{t+2}\rdb,a_{t+3},\dots,a_{n+2})
\otimes_{n+2-t}(|_{x=\partial}\ldb {a_{t}}_{\lambda_t} a_{t+1}\rdb'')
\\
&
=
X_{\lambda_1,\dots,\lambda_{t-1},
\lambda_{t}+\lambda_{t+1}+\lambda_{t+2},\lambda_{t+3},
\dots,\lambda_{n+2}}^{(t)}
(a_1,\dots
\\
&
\dots,a_{t-1}, \ldb \ldb {a_{t}}_{\lambda_{t}} a_{t+1}\rdb_{\lambda_t+\lambda_{t+1}} a_{t+2}\rdb_L,a_{t+3},\dots,a_{n+2})
\,.
\end{split}
\end{equation}
Let $A\in\mc V^{\otimes m}$ and $B\in\mc V^{\otimes l}$. We introduce the following notation
for $i<j$:
\begin{equation}
\begin{split}\label{20230811:eq7}
&X_{\lambda_1,\dots,\lambda_n}^{(i),(j)}(a_1,\dots,a_{i-1},a_i\otimes A,a_{i+1},\dots,a_{j-1},a_j\otimes B,a_{j+1},\dots,a_n)
\\
&=
X_{\lambda_1,\dots,\lambda_j+y,\dots,\lambda_n}^{(i)}
(a_1,\dots,a_i\otimes A,\dots,a_n)\otimes_{n+1-j}(|_{y=\partial}B)
\\
&=
\left(X_{\lambda_1,\dots,\lambda_i+x,\dots,\lambda_j+y,\dots,\lambda_n}
(a_1,\dots,a_n)\otimes_{n+1-i}(|_{x=\partial}A)\right)\otimes_{n+1-j}(|_{y=\partial}B)
\,;
\end{split}
\end{equation}
and for $i>j$:
\begin{equation}
\begin{split}\label{20230811:eq8}
&
X_{\lambda_1,\dots,\lambda_n}^{(i),(j)}
(a_1,\dots,a_{j-1},a_j\otimes B,a_{j+1},\dots,a_{i-1},a_i\otimes A,a_{i+1},\dots,a_n)
\\
&=
X_{\lambda_1,\dots,\lambda_j+y,\dots,\lambda_n}^{(i)}
(a_1,\dots,a_i\otimes A,\dots,a_n)
\otimes_{n+2-j}(|_{y=\partial}B)
\\
&=
\left(X_{\lambda_1,\dots,\lambda_j+y,\dots,\lambda_i+x,\dots,\lambda_n}
(a_1,\dots,a_n)\otimes_{n+1-i}(|_{x=\partial}A)\right)\otimes_{n+2-j}(|_{x=\partial}B)
\,.
\end{split}
\end{equation}
Then, using equations~\eqref{20230811:eq2}-\eqref{20230811:eq6} and the notations~\eqref{20230811:eq7}-\eqref{20230811:eq8}, we rewrite~\eqref{eq:dPsquare2pre} as
\begin{equation}
\begin{split}\label{eq:dPsquare2}
&\tilde\delta(X)_{\lambda_1,\dots,\lambda_{t-1},\lambda_t+\lambda_{t+1},\lambda_{t+2},\dots,\lambda_{n+2}}^{(t)}
(a_1,\dots,a_{t-1},\ldb {a_{t}}_{\lambda_t} a_{t+1}\rdb,a_{t+2},\dots,a_{n+2})
\\
&=
\sum_{s=1}^{t-1}(-1)^{s+1}
\ldb {a_s}_{\lambda_s}X_{\lambda_1,\stackrel{s}{\check{\dots}},\lambda_{t-1},\lambda_t+\lambda_{t+1},\lambda_{t+2},\dots,\lambda_{n+2}}^{(t-1)}
(a_1,\stackrel{s}{\check{\dots}}
\\
&\quad\quad\quad
\dots,a_{t-1},\ldb {a_{t}}_{\lambda_t} a_{t+1}\rdb,a_{t+2},\dots,a_{n+2})
\rdb_{(s)}
\\
&
+(-1)^{t+1}\ldb \ldb {a_{t}}_{\lambda_t} a_{t+1}\rdb_{\lambda_t+\lambda_{t+1}}
X_{\lambda_1,\dots,\lambda_{t-1},\lambda_{t+2},\dots\lambda_{n+2}}
(a_1,\dots
\\
&\quad\quad\quad
\dots,a_{t-1},a_{t+2},\dots,a_{n+2})
\rdb_{L,(t)}
\\
&+\sum_{s=t+2}^{n+2}(-1)^{s}
\ldb {a_s}_{\lambda_s}X_{\lambda_1,\dots,\lambda_{t-1},\lambda_t+\lambda_{t+1},\lambda_{t+2},
\stackrel{s-1}{\check{\dots}},\lambda_{n+2}}^{(t)}
(a_1,\dots
\\
&\quad\quad\quad
\dots,a_{t-1},\ldb {a_{t}}_{\lambda_t} a_{t+1}\rdb,a_{t+2},\stackrel{s-1}{\check{\dots}},a_{n+2})
\rdb_{(s)}
\\
&
+\sum_{s=1}^{t-2}(-1)^{s}
X_{\lambda_1,\dots,\lambda_s+\lambda_{s+1},\dots,
\lambda_t+\lambda_{t+1},
\dots,\lambda_{n+2}}^{(s),(t-1)}
(a_1,\dots
\\
&\quad\quad\quad
\dots,a_{s-1},\ldb {a_{s}}_{\lambda_s} a_{s+1}\rdb,a_{s+2},\dots,a_{t-1},\ldb {a_{t}}_{\lambda_t} a_{t+1}\rdb,a_{t+2},\dots,a_{n+2})
\\
&
+(-1)^{t-1} X_{\lambda_1,\dots,\lambda_{t-2},
\lambda_{t-1}+\lambda_t+\lambda_{t+1},\lambda_{t+2},
\dots,\lambda_{n+2}}^{(t-1)}
(a_1,\dots
\\
&\quad\quad\quad
\dots,a_{t-2}, \ldb {a_{t-1}}_{\lambda_{t-1}}\ldb {a_{t}}_{\lambda_t} a_{t+1}\rdb\rdb_L,a_{t+2},\dots,a_{n+2})
\\
&
+(-1)^{t} X_{\lambda_1,\dots,\lambda_{t-1},
\lambda_{t}+\lambda_{t+1}+\lambda_{t+2},\lambda_{t+3},
\dots,\lambda_{n+2}}^{(t)}
(a_1,\dots
\\
&\quad\quad\quad
\dots,a_{t-1}, \ldb \ldb {a_{t}}_{\lambda_{t}} a_{t+1}\rdb_{\lambda_t+\lambda_{t+1}} a_{t+2}\rdb_L,a_{t+3},\dots,a_{n+2})
\\
&
+\sum_{s=t+2}^{n+1}(-1)^{s-1}
X_{\lambda_1,\dots,\lambda_t+\lambda_{t+1},\dots,
\lambda_s+\lambda_{s+1},
\dots,\lambda_{n+2}}^{(s-1),(t)}
(a_1,\dots
\\
&\quad\quad\quad
\dots,a_{t-1},\ldb {a_{t}}_{\lambda_t} a_{t+1}\rdb,a_{t+2},\dots,a_{s-1},\ldb {a_{s}}_{\lambda_s} a_{s+1}\rdb,a_{s+2},\dots,a_{n+2})
\,.
\end{split}
\end{equation}
Hence, we replace~\eqref{eq:dPsquare1} and~\eqref{eq:dPsquare2} in equation~\eqref{eq:dPsquare} to get 
\begin{subequations}
\begin{align}
&\tilde\delta^2 (X)_{\lambda_1,\dots,\lambda_{n+2}}(a_1,\dots,a_{n+2})
\notag
\\
&=\sum_{t=2}^{n+2}\sum_{s=1}^{t-1}(-1)^{t+s}
\ldb {a_t}_{\lambda_t}
\ldb {a_{s}}_{\lambda_s}
X_{\lambda_1,\stackrel{s}{\check{\dots}}\stackrel{t}{\check{\dots}},\lambda_{n+2}}
(a_1,\stackrel{s}{\check{\dots}}\stackrel{t}{\check{\dots}},a_{n+2})\rdb_{(s)}
\rdb_{(t)}
\label{eq:dPsquareA1}
\\
&
+\sum_{t=1}^{n+1}\sum_{s=t+1}^{n+2}(-1)^{t+s+1}
\ldb {a_t}_{\lambda_t}
\ldb {a_s}_{\lambda_s}
X_{\lambda_1,\stackrel{t}{\check{\dots}}\stackrel{s}{\check{\dots}},\lambda_{n+2}}
(a_1,\stackrel{t}{\check{\dots}}\stackrel{s}{\check{\dots}},a_{n+2})\rdb_{(s-1)}
\rdb_{(t)}
\label{eq:dPsquareA2}
\\
\begin{split}\label{eq:dPsquareB1}
&-\sum_{t=3}^{n+2}\sum_{s=1}^{t-2}(-1)^{t+s}
\ldb {a_t}_{\lambda_t}
X_{\lambda_1,\dots,\lambda_{s-1},\lambda_s+\lambda_{s+1},\lambda_{s+2},
\stackrel{t-1}{\check{\dots}},\lambda_{n+2}}^{(s)}(a_1,\dots
\\
&\quad\quad\quad
\dots,a_{s-1},\ldb {a_{s}}_{\lambda_s} a_{s+1}\rdb,a_{s+2},\stackrel{t-1}{\check{\dots}},a_{n+2})
\rdb_{(t)}
\end{split}
\\
\begin{split}\label{eq:dPsquareD1}
&
+\sum_{t=2}^{n+1}
\ldb {a_t}_{\lambda_t}
X_{\lambda_1,\dots,\lambda_{t-2},\lambda_{t-1}+\lambda_{t+1},\lambda_{t+2},
\dots,\lambda_{n+2}}^{(t-1)}
(a_1,\dots
\\
&\quad\quad\quad
\dots,a_{t-2},\ldb {a_{t-1}}_{\lambda_{t-1}} a_{t+1}\rdb,a_{t+2},\dots,a_{n+2})
\rdb_{(t)}
\end{split}
\\
\begin{split}\label{eq:dPsquareC1}
&
+\sum_{t=1}^{n}\sum_{s=t+1}^{n+1}(-1)^{t+s}
\ldb {a_t}_{\lambda_t}
X_{\lambda_1,\stackrel{t}{\check{\dots}},\lambda_{s-1},\lambda_s+\lambda_{s+1},\lambda_{s+2},\dots,\lambda_{n+2}}^{(s-1)}
(a_1,\stackrel{t}{\check{\dots}}
\\
&\quad\quad\quad
\dots,a_{s-1},\ldb {a_{s}}_{\lambda_s} a_{s+1}\rdb,a_{s+2},\dots,a_{n+2})
\rdb_{(t)}
\end{split}
\\
\begin{split}\label{eq:dPsquareC2}
&
-\sum_{t=2}^{n+1}\sum_{s=1}^{t-1}(-1)^{t+s}
\ldb {a_s}_{\lambda_s}X_{\lambda_1,\stackrel{s}{\check{\dots}},\lambda_{t-1},\lambda_t+\lambda_{t+1},\lambda_{t+2},\dots,\lambda_{n+2}}^{(t-1)}
(a_1,\stackrel{s}{\check{\dots}}
\\
&\quad\quad\quad
\dots,a_{t-1},\ldb {a_{t}}_{\lambda_t} a_{t+1}\rdb,a_{t+2}\dots,a_{n+2})
\rdb_{(s)}
\end{split}
\\
\begin{split}\label{eq:dPsquareA3}
&-\sum_{t=1}^{n+1}\ldb \ldb {a_{t}}_{\lambda_t} a_{t+1}\rdb_{\lambda_t+\lambda_{t+1}}
X_{\lambda_1,\dots,\lambda_{t-1},\lambda_{t+2},\dots\lambda_{n+2}}
(a_1,\dots
\\
&\quad\quad\quad
\dots,a_{t-1},a_{t+2},\dots,a_{n+2})
\rdb_{L,(t)}
\end{split}
\\
\begin{split}\label{eq:dPsquareB2}
&+\sum_{t=1}^{n}\sum_{s=t+2}^{n+2}(-1)^{t+s}
\ldb {a_s}_{\lambda_s}X_{\lambda_1,\dots,\lambda_{t-1},\lambda_t+\lambda_{t+1},\lambda_{t+2},
\stackrel{s-1}{\check{\dots}},\lambda_{n+2}}^{(t)}
(a_1,\dots
\\
&\quad\quad\quad
\dots,a_{t-1},\ldb {a_{t}}_{\lambda_t} a_{t+1}\rdb,a_{t+2},\stackrel{s-1}{\check{\dots}},a_{n+2})
\rdb_{(s)}
\end{split}
\\
\begin{split}\label{eq:dPsquareE1}
&
+\sum_{t=3}^{n+1}\sum_{s=1}^{t-2}(-1)^{t+s}
X_{\lambda_1,\dots,\lambda_s+\lambda_{s+1},\dots,
\lambda_t+\lambda_{t+1},
\dots,\lambda_{n+2}}^{(s),(t-1)}
(a_1,\dots
\\
&\quad\quad\quad
\dots,a_{s-1},\ldb {a_{s}}_{\lambda_s} a_{s+1}\rdb,a_{s+2},\dots,a_{t-1},\ldb {a_{t}}_{\lambda_t} a_{t+1}\rdb,a_{t-2},\dots,a_{n+2})
\end{split}
\\
\begin{split}\label{eq:dPsquareD2}
&
-\sum_{t=2}^{n+1} X_{\lambda_1,\dots,
\lambda_{t-1}+\lambda_t+\lambda_{t+1},
\dots,\lambda_{n+2}}^{(t-1)}
(a_1,\dots
\\
&\quad\quad\quad
\dots,a_{t-2}, \ldb {a_{t-1}}_{\lambda_{t-1}}\ldb {a_{t}}_{\lambda_t} a_{t+1}\rdb\rdb_L,a_{t+2},\dots,a_{n+2})
\end{split}
\\
\begin{split}\label{eq:dPsquareD3}
&
+\sum_{t=1}^{n}X_{\lambda_1,\dots,
\lambda_{t}+\lambda_{t+1}+\lambda_{t+2},
\dots,\lambda_{n+2}}^{(t)}
(a_1,\dots
\\
&\quad\quad\quad
\dots,a_{t-1},\ldb \ldb {a_{t}}_{\lambda_{t}} a_{t+1}\rdb_{\lambda_t+\lambda_{t+1}} a_{t+2}\rdb_L,a_{t+3},\dots,a_{n+2})
\end{split}
\\
\begin{split}\label{eq:dPsquareE2}
&
-\sum_{t=1}^{n-1}\sum_{s=t+2}^{n+1}(-1)^{t+s}
X_{\lambda_1,\dots,\lambda_t+\lambda_{t+1},\dots,
\lambda_s+\lambda_{s+1},
\dots,\lambda_{n+2}}^{(s-1),(t)}
(a_1,\dots
\\
&\quad\quad\quad
\dots,a_{t-1},\ldb {a_{t}}_{\lambda_t} a_{t+1}\rdb,a_{t+2},\dots,a_{s-1},\ldb {a_{s}}_{\lambda_s} a_{s+1}\rdb,a_{s+2}\dots,a_{n+2})
\,.
\end{split}
\end{align}
\end{subequations}
To conclude the proof we show that the RHS above vanishes. Swapping $s$ and $t$ in~\eqref{eq:dPsquareA1}
and changing the order of summation we get
\begin{align*}
&\eqref{eq:dPsquareA2} +~\eqref{eq:dPsquareA1}
\\
&=
\sum_{t=1}^{n+1}\sum_{s=t+1}^{n+2}(-1)^{t+s+1}
\Big(
\ldb {a_t}_{\lambda_t}
\ldb {a_s}_{\lambda_s}
X_{\lambda_1,\stackrel{t}{\check{\dots}}\stackrel{s}{\check{\dots}},\lambda_{n+2}}
(a_1,\stackrel{t}{\check{\dots}}\stackrel{s}{\check{\dots}},a_{n+2})\rdb_{(s-1)}
\rdb_{(t)}
\\
&-
\ldb {a_s}_{\lambda_s}
\ldb {a_{t}}_{\lambda_t}
X_{\lambda_1,\stackrel{t}{\check{\dots}}\stackrel{s}{\check{\dots}},\lambda_{n+2}}
(a_1,\stackrel{t}{\check{\dots}}\stackrel{s}{\check{\dots}},a_{n+2})\rdb_{(t)}
\rdb_{(s)}
\Big)
\\
&
=\sum_{t=1}^{n+1}
\Big(
\ldb {a_t}_{\lambda_t}
\ldb {a_{t+1}}_{\lambda_{t+1}}
X_{\lambda_1,\dots,\lambda_{t-1},\lambda_{t+2},\dots,\lambda_{n+2}}
(a_1,\dots,a_{t-1},a_{t+2},\dots,a_{n+2})\rdb_{(t)}
\rdb_{(t)}
\\
&-
\ldb {a_{t+1}}_{\lambda_{t+1}}
\ldb {a_{t}}_{\lambda_t}
X_{\lambda_1,\dots,\lambda_{t-1},\lambda_{t+2},\dots,\lambda_{n+2}}
(a_1,\dots,a_{t-1},a_{t+2},\dots,a_{n+2})\rdb_{(t)}
\rdb_{(t+1)}
\Big)
\\
&
=\sum_{t=1}^{n+1}
\ldb \ldb {a_t}_{\lambda_t} {a_{t+1}}\rdb _{\lambda_t+\lambda_{t+1}}
X_{\lambda_1,\dots,\lambda_{t-1},\lambda_{t+2},\dots,\lambda_{n+2}}
(a_1,\dots,a_{t-1},a_{t+2},\dots,a_{n+2})\rdb_{L,(t)}
\\
&=-\eqref{eq:dPsquareA3}
\,.
\end{align*}
In the second equality above we used Lemma ~\ref{lem:PVAcoh-pre7}(i) and in the third equality we used Lemma~\ref{lem:PVAcoh-pre7}(ii).
Hence, in the expression for $\tilde\delta^2( X)$ we have the cancellation
$$
\eqref{eq:dPsquareA1}+\eqref{eq:dPsquareA2}+\eqref{eq:dPsquareA3}=0\,.
$$
Next, by swapping $s$ and $t$ in both equations~\eqref{eq:dPsquareB2} and~\eqref{eq:dPsquareC2} and changing the order of summation, it follows that we get the cancellations
$$
\eqref{eq:dPsquareB1}+\eqref{eq:dPsquareB2}=0
\,,
\qquad
\eqref{eq:dPsquareC1}+\eqref{eq:dPsquareC2}=0
\,.
$$
Furthermore, we use the definition of $X^{(s)}$ and the second identity in~\eqref{20230811:eq1}
(with $s=t$ and $i=n+2-t$) to get 
\begin{equation}
\begin{split}\label{eq:dPsquareD1b}
&\eqref{eq:dPsquareD1}=
\\
&=\sum_{t=2}^{n+1}
X_{\lambda_1,\dots,\lambda_{t-1}+\lambda_{t+1}+x,
\dots,\lambda_{n+2}}(a_1,\dots,\ldb {a_{t-1}}_{\lambda_{t-1}} a_{t+1}\rdb',\dots,a_{n+2})
\\
&\quad\quad\quad\quad\quad\quad\quad\quad\quad\quad\quad\quad\quad
\otimes_{n+2-t}
\ldb {a_t}_{\lambda_t}(|_{x=\partial}\ldb {a_{t-1}}_{\lambda_{t-1}} a_{t+1}\rdb'')\rdb
\\
&=\sum_{t=2}^{n+1}
X_{\lambda_1,\dots,\lambda_{t-1}+\lambda_t+\lambda_{t+1}+x,
\dots,\lambda_{n+2}}(a_1,\dots,\ldb {a_{t-1}}_{\lambda_{t-1}} a_{t+1}\rdb',\dots,a_{n+2})
\\
&\quad\quad\quad\quad\quad\quad\quad\quad\quad\quad\quad\quad\quad
\otimes_{n+2-t}
(|_{x=\partial}\ldb {a_t}_{\lambda_t}\ldb {a_{t-1}}_{\lambda_{t-1}} a_{t+1}\rdb''\rdb)
\\
&=\sum_{t=2}^{n+1}
X_{\lambda_1,\dots,\lambda_{t-1}+\lambda_t+\lambda_{t+1},
\dots,\lambda_{n+2}}^{(t-1)}(a_1,\dots
\\
&\quad\quad\quad\quad\quad\quad\quad\quad\quad\quad\quad\quad\quad
\dots,\ldb {a_t}_{\lambda_t}\ldb {a_{t-1}}_{\lambda_{t-1}} a_{t+1}\rdb\rdb_R,\dots,a_{n+2})
\,,
\end{split}
\end{equation}
where in the second equality we used the second sesquilinearity axiom~\eqref{eq:sesqui}, and in the third 
equality we used the definition of $X^{(s)}$ and the notation~\eqref{notation}. Hence,
using~\eqref{eq:dPsquareD1b}, the linearity of $X^{(s)}$, and changing $t$ in $t-1$ in~\eqref{eq:dPsquareD3}
we get 
\begin{align*}
&\eqref{eq:dPsquareD2}+\eqref{eq:dPsquareD1}+\eqref{eq:dPsquareD3}
\\
&=-\sum_{t=2}^{n+1}
X_{\lambda_1,\dots,\lambda_{t-1}+\lambda_t+\lambda_{t+1},
\dots,\lambda_{n+2}}^{(t-1)}(a_1,\dots,a_{t-2},J_t,a_{t+2},\dots,a_{n+2})\,,  \\
&\text{for }J_t:=\ldb {a_{t-1}}_{\lambda_{t-1}}\ldb {a_{t}}_{\lambda_t} a_{t+1}\rdb\rdb_L
-\ldb {a_t}_{\lambda_t}\ldb {a_{t-1}}_{\lambda_{t-1}} a_{t+1}\rdb\rdb_R
\\
&\quad\quad\quad
-\ldb \ldb {a_{t-1}}_{\lambda_{t-1}} a_{t}\rdb_{\lambda_{t-1}+\lambda_{t}} a_{t+1}\rdb_L\,,
\end{align*}
which vanishes by linearity as $J_t=0$ since we are assuming that the Jacobi identity~\eqref{eq:jacobi2} holds.
Finally, let us swap $s$ and $t$ and change the order of summation in~\eqref{eq:dPsquareE2} to get
\begin{align*}
&\eqref{eq:dPsquareE2}
=-\sum_{t=3}^{n+1}\sum_{s=1}^{t-2}(-1)^{t+s}
X_{\lambda_1,\dots,\lambda_s+\lambda_{s+1},\dots,
\lambda_t+\lambda_{t+1},
\dots,\lambda_{n+2}}^{(t-1),(s)}
(a_1,\dots
\\
&\quad
\dots,a_{s-1},\ldb {a_{s}}_{\lambda_s} a_{s+1}\rdb,a_{s+2},\dots,a_{t-1},\ldb {a_{t}}_{\lambda_t} a_{t+1}\rdb,a_{t-2},\dots,a_{n+2})
\\
&
=-\sum_{t=3}^{n+1}\sum_{s=1}^{t-2}(-1)^{t+s}
\Big(
X_{\lambda_1,\dots,\lambda_s+\lambda_{s+1}+y,\dots,
\lambda_t+\lambda_{t+1}+x,
\dots,\lambda_{n+2}}
(a_1,\dots
\\
&\quad
\dots,a_{s-1},\ldb {a_{s}}_{\lambda_s} a_{s+1}\rdb',a_{s+2},\dots,a_{t-1},\ldb {a_{t}}_{\lambda_t} a_{t+1}\rdb',a_{t+2},\dots,a_{n+2})
\\
&\quad
\otimes_{n+2-t}(|_{x=\partial}\ldb {a_{t}}_{\lambda_t} a_{t+1}\rdb'')\Big)
\otimes_{n+2-s}(|_{y=\partial}\ldb {a_{s}}_{\lambda_s} a_{s+1}\rdb'')
\\
&
=-\sum_{t=3}^{n+1}\sum_{s=1}^{t-2}(-1)^{t+s}
\Big(
X_{\lambda_1,\dots,\lambda_s+\lambda_{s+1}+y,\dots,
\lambda_t+\lambda_{t+1}+x,
\dots,\lambda_{n+2}}
(a_1,\dots
\\
&\quad
\dots,a_{s-1},\ldb {a_{s}}_{\lambda_s} a_{s+1}\rdb',a_{s+2},\dots,a_{t-1},\ldb {a_{t}}_{\lambda_t} a_{t+1}\rdb',a_{t+2},\dots,a_{n+2})
\\
&\quad
\otimes_{n+1-s}(|_{y=\partial}\ldb {a_{s}}_{\lambda_s} a_{s+1}\rdb'')\Big)
\otimes_{n+2-t}(|_{x=\partial}\ldb {a_{t}}_{\lambda_t} a_{t+1}\rdb'')
\\
&=-\sum_{t=3}^{n+1}\sum_{s=1}^{t-2}(-1)^{t+s}
X_{\lambda_1,\dots,\lambda_s+\lambda_{s+1},\dots,
\lambda_t+\lambda_{t+1},
\dots,\lambda_{n+2}}^{(s),(t-1)}
(a_1,\dots
\\
&\quad
\dots,a_{a-1},\ldb {a_{s}}_{\lambda_s} a_{s+1}\rdb,a_{s+2},\dots,a_{t-1},\ldb {a_{t}}_{\lambda_t} a_{t+1}\rdb,a_{t+2},\dots,a_{n+2})
=-\eqref{eq:dPsquareE1}
\,,
\end{align*}
where in the second equality we used equation~\eqref{20230811:eq8}, in the third equality we used the identity
$$
(A\otimes_{n-i} a)\otimes_{n+1-j}b
=(A\otimes_{n-j}b)\otimes_{n-i}a\,,
$$
which can be easily verified to hold for every $a,b\in\mc V$, $A\in\mc V^{\otimes n}$ and $i>j$,
and in the fourth equality we used equation~\eqref{20230811:eq7}. Hence, in the expression for
$\delta_H^2X$ we have the further cancellation 
$$
\eqref{eq:dPsquareE1}+\eqref{eq:dPsquareE2}=0
$$
thus concluding the proof that $\tilde\delta^2(X)=0$.
\end{proof}
By Theorem~\ref{thm:PVAcoh}
we have a complex $(\widetilde{\Gamma}(\mc V),\tilde\delta)$.
\begin{definition}\label{def:bas-dPVA}
Let $\mc V$ be a dPVA.
The complex $(\widetilde{\Gamma}(\mc V),\tilde\delta)$ is called the \emph{basic dPVA complex} of $\mc V$. 
The cohomology\glslink{dHbasV}{}
$$
\dH_{\textrm{bas}}(\mc V)=\coH(\widetilde{\Gamma}(\mc V),\tilde\delta)
=\bigoplus_{n\in \mb Z_{\geq0}}\dH_{\textrm{bas}}^n(\mc V)\,,
\quad \dH_{\textrm{bas}}^n(\mc V)=\ker(\tilde\delta|_{\widetilde{\Gamma}^n(\mc V)})/\tilde\delta(\widetilde{\Gamma}^{n-1}(\mc V))\,,
$$
of this complex is called the \emph{basic dPVA cohomology} of $\mc V$.
\end{definition}
\begin{example}
As in Example~\ref{exa:basiclow} we can describe explicitly the first few cohomology spaces of the basic dPVA cohomology of $\mc V$.
By \eqref{exa1}, we have
$$
\dH^0_{\textrm{bas}}(\mc V)=
\{a\in\mc V\mid \ldb a_\lambda b\rdb=0\text{ for every }b\in\mc V\}
=Z(\mc V)\glslink{center-nc}{}
$$
is the center of the dPVA $\mc V$.
Next, let $X\in\widetilde{\Gamma}^1(\mc V)$ be a basic 1-cochain (cf. Example~\ref{exa:rdcd}). We say that $X$ is a \emph{Poisson 2-fold right conformal derivation} if
$$
X_{\lambda+\mu}^{(1)}(\ldb a_\lambda b\rdb)
=\ldb a_\lambda X_\mu(b)\rdb_L-\ldb b_\mu X_\lambda(a)\rdb_R
\,,
$$
for every $a,b\in \mc V$. We say that
$X$ is an \emph{inner 2-fold right conformal derivation}
if $X=\ldb -_\lambda a\rdb$, for some $a\in\mc V$.
Hence, from equations~\eqref{exa1} and~\eqref{exa2} we get that
$$
\dH^1_{\textrm{bas}}=
\frac{\{\text{Poisson 2-fold right conformal derivations}\}}{\{\text{inner 2-fold right conformal derivations}\}}
\,.
$$
\end{example}

%%%
\section[Reduced \MakeLowercase{d}PVA cohomology complex]{The reduced double Poisson vertex algebra cohomology}\label{sec:red-dPVA}

Let $\mc V$ be a dPVA and let $(\widetilde{\Gamma}(\mc V,\tilde\delta)$ be the complex of basic cochains introduced in Sections \ref{sec:omega1} and \ref{sec:omega2}.
Recall from Section~\ref{sec:1.1} the definition of the commutator space $[\widetilde{\Gamma}(\mc V),\widetilde{\Gamma}(\mc V)]$ and recall also from \eqref{20230802:eq1} that we have an action of the derivation $\partial$ of $\mc V$ on the space of basic cochains $\widetilde{\Gamma}(\mc V)$.
\begin{proposition}\label{20230802:prop1}
\begin{enumerate}[(a)]
\item
For homogeneous elements $X,Y\in\widetilde{\Gamma}(\mc V)$ we have
$$
\tilde \delta(XY)=\tilde\delta(X)Y+(-1)^{p(X)}X \tilde\delta(Y)\,.
$$
\item The (graded) commutator space $[\widetilde{\Gamma}(\mc V),\widetilde{\Gamma}(\mc V)]$ is compatible
with the $\mb Z_{\geq0}$-grading defined in~\eqref{eq:Omegatilde} and it is preserved by $\tilde\delta$.
\item The space $\partial\widetilde{\Gamma}(\mc V)$ is compatible
with the $\mb Z_{\geq0}$-grading~\eqref{eq:Omegatilde} and it is preserved by $\tilde\delta$.
\end{enumerate}
\end{proposition}
\begin{proof}
Let $X\in\widetilde{\Gamma}^m(\mc V)$ and $Y\in\widetilde{\Gamma}^n(\mc V)$.
For $1\leq s\leq m$, using Lemma~\ref{lem:PVAcoh-pre1}(i) and~\eqref{20240805:eq1},
we have
\begin{equation}
\begin{split}\label{20230802:eq3}
&\!\ldb {a_s}_{\lambda_s}
(XY)_{\lambda_1,\stackrel{s}{\check{\dots}},\lambda_{m+n+1}}
(a_1,\stackrel{s}{\check{\dots}},a_{m+n+1})\rdb_{(s)}
\\
&\!=\ldb {a_s}_{\lambda_s}
X_{\lambda_1,\stackrel{s}{\check{\dots}},\lambda_{m+1}}
(a_1,\stackrel{s}{\check{\dots}},a_{m+1})\rdb_{(s)}
Y_{\lambda_{m+2},\dots,\lambda_{m+n+1}}(a_{m+2},\dots,a_{m+n+1})
\,.
\end{split}
\end{equation}
Similarly,
for $m+2\leq s\leq m+n+1$, we have
\begin{equation}
\begin{split}\label{20230802:eq4}
&\ldb {a_s}_{\lambda_s}
(XY)_{\lambda_1,\stackrel{s}{\check{\dots}},\lambda_{m+n+1}}
(a_1,\stackrel{s}{\check{\dots}},a_{m+n+1})\rdb_{(s)}
\\
&=
X_{\lambda_{1},\dots,\lambda_{m}}(a_{1},\dots,a_{m})
\ldb {a_s}_{\lambda_s}
Y_{\lambda_{m+1},\stackrel{s}{\check{\dots}},\lambda_{m+n+1}}
(a_{m+1},\stackrel{s}{\check{\dots}},a_{m+n+1})\rdb_{(s-m)}
\,.
\end{split}
\end{equation}
Moreover, using Lemma~\ref{lem:PVAcoh-pre1}(i),~\eqref{20240805:eq1} and the Leibniz rule~\eqref{eq:lleibniz},
we have
\begin{equation}
\begin{split}\label{20230802:eq5}
&\ldb {a_{m+1}}_{\lambda_{m+1}}
(XY)_{\lambda_1,\stackrel{m+1}{\check{\dots}},\lambda_{m+n+1}}
(a_1,\stackrel{m+1}{\check{\dots}},a_{m+n+1})\rdb_{(m+1)}
\\
&=
\ldb {a_{m+1}}_{\lambda_{m+1}}
X_{\lambda_1,\dots,\lambda_{m}}(a_1,\dots,a_{m})
Y_{\lambda_{m+2},\dots,\lambda_{m+n+1}}(a_{m+2},\dots,a_{m+n+1})\rdb_{(m+1)}
\\
&=
\ldb {a_{m+1}}_{\lambda_{m+1}}
X_{\lambda_1,\dots,\lambda_{m}}(a_1,\dots,a_{m})
\rdb_{(m+1)}
Y_{\lambda_{m+2},\dots,\lambda_{m+n+1}}(a_{m+2},\dots,a_{m+n+1})\\
&+
X_{\lambda_1,\dots,\lambda_{m}}(a_1,\dots,a_{m})
\ldb {a_{m+1}}_{\lambda_{m+1}}
Y_{\lambda_{m+2},\dots,\lambda_{m+n+1}}(a_{m+2},\dots,a_{m+n+1})\rdb_{(1)}
\,.
\end{split}
\end{equation}
Part (a) then follows by a straightforward computation using the definition of $\tilde\delta$ given by
\eqref{eq:diff}, the identities~\eqref{20230802:eq3},
\eqref{20230802:eq4},~\eqref{20230802:eq5}
and Lemma~\ref{lem:PVAcoh-pre1}(ii).
Next, by definition of the associative product~\eqref{20230802:eq2} in $\widetilde{\Gamma}(\mc V)$ we clearly have $[\widetilde{\Gamma}^m(\mc V),\widetilde{\Gamma}^n(\mc V)]\subset\widetilde{\Gamma}^{m+n}(\mc V)$.
Hence, the commutator space $[\widetilde{\Gamma}(\mc V),\widetilde{\Gamma}(\mc V)]$ is compatible
with the $\mb Z_{\geq0}$-grading~\eqref{eq:Omegatilde}. By part (a),
$\tilde\delta$ is an odd derivation of the superalgebra $\widetilde{\Gamma}(\mc V)$ thus it preserves the commutator space $[\widetilde{\Gamma}(\mc V),\widetilde{\Gamma}(\mc V)]$. This proves part (b).

Let us prove part (c). From equation~\eqref{20230802:eq1} we have that $\partial$ preserves each homogeneous
component $\widetilde{\Gamma}^m(\mc V)$. Hence, the space $\partial\widetilde{\Gamma}(\mc V)$ is compatible
with the $\mb Z_{\geq0}$-grading~\eqref{eq:Omegatilde}. Furthermore, using~\eqref{20230802:eq1} and the second
sesquilinearity axiom~\eqref{eq:sesqui} it is straightforward to check that
$\tilde\delta\circ\partial=\partial\circ \tilde\delta$. The second part of the claim follows from this identity.
\end{proof}
By Proposition~\ref{20230802:prop1}, we can consider the $\mb Z_{\geq0}$-graded complex
$(\Gamma(\mc V),\delta)$, where
\begin{equation}\label{eq:Omega}
\Gamma(\mc V)
=\widetilde{\Gamma}(\mc V)/(\partial\widetilde{\Gamma}(\mc V)
+[\widetilde{\Gamma}(\mc V),\widetilde{\Gamma}(\mc V)])
=\bigoplus_{n\in\mb Z_{\geq0}}\Gamma^n(\mc V)\,,
\glslink{AsGammaRed}{}
\end{equation}
and $\delta:\widetilde{\Gamma}(\mc V)\to\widetilde{\Gamma}(\mc V)$
is induced by the differential $\tilde{\delta}$ of $\widetilde{\Gamma}(\mc V)$:
for $X\in\widetilde{\Gamma}(\mc V)$, the action of $\delta$ on the coset $[X]=X+\partial\widetilde{\Gamma}(\mc V)+[\widetilde{\Gamma}(\mc V),\widetilde{\Gamma}(\mc V)]\in\Gamma (V)$ is given by\glslink{deltadiff}{} 
$$
\delta([X])=[\tilde\delta(X)]
\,.
$$
For example, we have $\Gamma^0(\mc V)=\mc V/\partial \mc V+[\mc V,\mc V]=\mc V_\sharp$ and we have the canonical quotient map $\tint :\widetilde{\Gamma}^0(\mc V)\rightarrow \Gamma^0(\mc V)$. By an abuse of notation we denote by the same symbol $\tint: \widetilde{\Gamma}(\mc V)\to\Gamma(\mc V)$\glslink{tint}{}
the canonical quotient map with kernel $(\partial\widetilde{\Gamma}(\mc V)
+[\widetilde{\Gamma}(\mc V),\widetilde{\Gamma}(\mc V)])$.
Let $X\in\widetilde{\Gamma}^1(V)$ and define the linear map
$\widetilde{\proj}_1(X):\mc V\to\mc V$ by
($a\in V$)
\begin{equation}\label{eq:P1}
\widetilde{\proj}_1(X)(a)=(|_{x=\partial}\mult(X_{-x}(a))^\sigma)\in \mc V\,.
\end{equation}
It follows by a direct computation using the Leibniz rule \eqref{eq:Leibnizmaps} that $\widetilde{\proj}_1(X)$ is a derivation for the associative product of $\mc V$ which, by the sesquilinearity \eqref{eq:sesquimaps}, commutes with $\partial$. Hence, we get a map 
$\widetilde{\proj}_1:\widetilde{\Gamma}^1(\mc V)\to \Vect(\mc V)^{\partial}$,
where $\Vect(\mc V)^{\partial}$ denotes the space of derivations of the product of $\mc V$ commuting with $\partial$. It is not hard to check that
$\partial\widetilde{\Gamma}^1(V)+[\widetilde{\Gamma}^0(\mc V),\widetilde{\Gamma}^1(\mc V)]\subset\ker\tint|_{\widetilde{\Gamma}^1(\mc V)}$, so we get an induced
linear map
\begin{equation}\label{eq:20250702:eq1bis}
\proj_1:\Gamma^1(\mc V)\rightarrow \Vect(\mc V)^{\partial}\,.
\end{equation}
\begin{definition}
Let $\mc V$ be a dPVA. The complex $(\Gamma(\mc V),\delta)$
is called the \emph{reduced dPVA complex}
of $\mc V$ and its cohomology\glslink{dHredV}{}
$$
\dH_{\textrm{red}}(\mc V)=\coH(\Gamma(\mc V),\delta)
=\bigoplus_{n\in\mb Z_{\geq0}}\dH_{\textrm{red}}^n(\mc V) ,
\quad \dH_{\textrm{red}}^n(\mc V)=\ker(\delta|_{\Gamma^n(\mc V)})/\delta(\Gamma^{n-1}(\mc V)) ,
$$
is called the \emph{reduced dPVA cohomology} of $\mc V$.
\end{definition}
By Theorem~\ref{thm:PVAcoh} and Proposition~\ref{20230802:prop1} we have a short exact sequence of complexes
$$
0\rightarrow \partial\widetilde{\Gamma}(\mc V)
+[\widetilde{\Gamma}(\mc V),\widetilde{\Gamma}(\mc V)]
\rightarrow \widetilde{\Gamma}(\mc V)
\rightarrow \Gamma(\mc V)\rightarrow 0
$$
which leads to the following long exact sequence in cohomology (we denote $I(\mc V)=\partial\widetilde{\Gamma}(\mc V)
+[\widetilde{\Gamma}(\mc V),\widetilde{\Gamma}(\mc V)]$
)
\begin{equation}\label{eq:les2}
% \text{this diagram has been commented to speed up compiling}
% \begin{comment}
\hspace{-1.4cm}
\begin{tikzcd}
0 \arrow[r]& \coH^0(I(\mc V),\tilde\delta)\arrow[r]
& \dH^0_{\textrm{bas}}(\mc V)\arrow[r] &
\dH^0_{\textrm{red}}(\mc V) \arrow[dll,in=175,out=-5]&
\\
&\coH^1(I(\mc V),\tilde\delta)\arrow[r]&
\dH^1_{\textrm{bas}}(\mc V) \arrow[r]
&\dH^1_{\textrm{red}}(\mc V) \arrow[dll,in=175,out=-5]&
\\
&\coH^2(I(\mc V),\tilde\delta)\arrow[r]&
\dH^2_{\textrm{bas}}(\mc V) \arrow[r]
&\dH^2_{\textrm{red}}(\mc V)\arrow[r] &\dots
\end{tikzcd}
% \end{comment}
\end{equation}
By Proposition \ref{20230802:prop1}(c) we have that $(\partial \widetilde{\Gamma}(\mc V),\tilde\delta)\subset (\widetilde {\Gamma}(\mc V),\tilde\delta)$ is a subcomplex.
\begin{proposition}\label{20250812:prop1}
We have a morphism of complexes
$$
\widetilde{\Gamma}(\mc V)\to\partial\widetilde{\Gamma}(\mc V)\,,
\quad X\mapsto \partial X
\,,
$$
which induces isomorphisms $\dH_{\mathrm{bas}}^n(\mc V)\simeq \coH^n(\partial\widetilde{\Gamma}(\mc V),\tilde\delta)$, for every $n\geq1$. 
\end{proposition}
\begin{proof}
It is an immediate consequence of the fact that the action of $\partial$ on $\widetilde{\Gamma}^n(\mc V)$ given by \eqref{20230802:eq1} is injective for $n\geq1$ and of the identity $\tilde\delta\circ\partial=\partial\circ\tilde\delta$ which was proved in the proof of Proposition \ref{20230802:prop1}(c). 
\end{proof}
Since the commutator space
$[\widetilde{\Gamma}(\mc V),\widetilde{\Gamma}(\mc V)]$ is preserved by $\partial$,
by Proposition \ref{20230802:prop1}(b) and (c) we also have that $(\partial[\widetilde{\Gamma}(\mc V),\widetilde{\Gamma}(\mc V)],\tilde\delta)\subset ([\widetilde{\Gamma}(\mc V),\widetilde{\Gamma}(\mc V)],\tilde\delta)$ is a subcomplex.
The next result is proved similarly to Proposition \ref{20250812:prop1}.
\begin{proposition}\label{20250812:prop2}
We have a morphism of complexes
$$
[\widetilde{\Gamma}(\mc V),\widetilde{\Gamma}(\mc V)]
\to\partial[\widetilde{\Gamma}(\mc V),\widetilde{\Gamma}(\mc V)]\,,
\quad X\mapsto \partial X
\,,
$$
which induces isomorphisms $\coH^n([\widetilde{\Gamma}(\mc V),\widetilde{\Gamma}(\mc V)],\tilde\delta)\simeq \coH^n(\partial[\widetilde{\Gamma}(\mc V),\widetilde{\Gamma}(\mc V)],\tilde\delta)$, for every $n\geq1$. 
\end{proposition}

%%%%%%%%%%% NEW CHAPTER %%%%%%%%%%%%%%%
%%%%%%%%%%% NEW CHAPTER %%%%%%%%%%%%%%%
%%%%%%%%%%% NEW CHAPTER %%%%%%%%%%%%%%%
%%%%%%%%%%% NEW CHAPTER %%%%%%%%%%%%%%%

\chapter[Variational \MakeLowercase{d}PVA cohomology]{Variational double Poisson vertex algebra cohomology}\label{sec:compl-dPVA}
In this chapter we review the definition and properties of $n$-fold $\lambda$-brackets
following \cite{DSKV} and we define the variational double Poisson vertex algebra complex and its cohomology. This construction is the ``double" analogue of the one
presented in Section \ref{sec:var-PVA}.

%%%
\section{\texorpdfstring{$n$}{n}-fold \texorpdfstring{$\lambda$}{lambda}-brackets}\label{sec:compl-dPVA-1}
Let $\mc V$ be a differential algebra. In this section we generalize to the noncommutative case the 
construction of $n$-$\lambda$-brackets described in Section \ref{sec:var-PVA} by introducing the
$\kk$-vector space of $n$-fold $\lambda$-brackets
on $\mc V$, which we denote by $C^n(\mc V)$.

For $n=0$ we set $C^0(\mc V)=\mc V_\sharp=\mc V/(\partial \mc V+[\mc V,\mc V])$. For $n\geq1$ we
give the following definition.
\begin{definition}\label{def:nfoldlambda}
An $n$-\emph{fold} $\lambda$-\emph{bracket}\glslink{n-lambda-nc}{}
on $\mc V$ is a linear map
$$
\ldb-{}_{\lambda_1}-\dots-{}_{\lambda_{n-1}}-\rdb:\mc V^{\otimes n}\to\mc V^{\otimes n}[\lambda_1,\ldots,\lambda_{n-1}]
$$
which maps $a_1\otimes \dots\otimes a_n\mapsto \ldb a_1{}_{\lambda_1}a_2\dots a_{n-1}{}_{\lambda_{n-1}}a_n\rdb$
satisfying
\begin{enumerate}[(a)]
\item
sesquilinearity:
\begin{equation}\label{20140702:eq4}
\ldb a_1{}_{\lambda_1}\cdots{}_{\lambda_{i-1}}(\partial a_i)_{\lambda_i}\cdots a_{n-1}{}_{\lambda_{n-1}}a_n\rdb
=-\lambda_i\ldb a_1{}_{\lambda_1}\cdots a_{n-1}{}_{\lambda_{n-1}}a_n\rdb\,,
\end{equation}
for all $i=1,\ldots,n-1$, and
\begin{equation}\label{20140702:eq5}
\ldb a_1{}_{\lambda_1}\cdots a_{n-1}{}_{\lambda_{n-1}}(\partial a_n)\rdb
=-\lambda_n^\dagger
\ldb a_1{}_{\lambda_1}\cdots a_{n-1}{}_{\lambda_{n-1}}a_n\rdb\,,
\end{equation}
where we are using the notation \eqref{eq:dagger};
\item
skewsymmetry:
\begin{equation}\label{eq:nfold-skew}
\ldb{a_1}_{\lambda_1}\dots{a_{n-1}}_{\lambda_{n-1}}a_n\rdb
=(-1)^{n+1}\left(\big|_{\lambda_n=\lambda_n^\dagger}
\ldb{a_2}_{\lambda_2}\dots{a_{n}}_{\lambda_n}a_1\rdb^\sigma\right)
\,;
\end{equation}
\item Leibniz rules:
\begin{equation}\label{20140702:eq6}
\begin{array}{l}
\displaystyle{
\vphantom{\Big(}
\ldb a_1{}_{\lambda_1}\cdots bc_{\lambda_i}\dots a_{n-1}{}_{\lambda_{n-1}}a_n\rdb
=(|_{x=\partial} b)\star_{i}
\ldb a_1{}_{\lambda_1}\cdots c_{\lambda_i+x}\dots a_{n-1}{}_{\lambda_{n-1}}a_n\rdb
} \\
\displaystyle{
\vphantom{\Big(}
+\ldb a_1{}_{\lambda_1}\cdots b_{\lambda_i+x}\dots a_{n-1}{}_{\lambda_{n-1}}a_n\rdb
\star_{n-i} (|_{x=\partial}c)
\,,}
\end{array}
\end{equation}
for all $i=1,\ldots,n$, where we recall from Section \ref{sec:1.1} that $\star_n=\star_0$.
\end{enumerate}
\end{definition}
We denote by $C^n(\mc V)$ the space of $n$-fold $\lambda$-brackets 
and we let\glslink{AsPolyLamb}{}
$$
C(\mc V)=\bigoplus_{n\in\mb Z_{\geq0}}C^n(\mc V)
\,.
$$
For $n=1$, we have that $C^1(\mc V)=\Vect(\mc V)^\partial$\glslink{vect-partial}{} is the space of derivations of the associative product of $\mc V$ commuting with $\partial$, while, for $n=2$, $C^2(\mc V)$ consists of $2$-fold $\lambda$-brackets on $\mc V$ defined by \eqref{eq:sesqui}-\eqref{eq:rleibniz}.
\begin{remark}\label{rem:6.17}
The skewsymmetry condition \eqref{eq:nfold-skew} can be rewritten as
\begin{equation}\label{eq:nfold-skew2}
\ldb-{}_{\lambda_1}-\dots-{}_{\lambda_{n-1}}-\rdb
=(-1)^{n+1}|_{\lambda_n=\lambda_n^\dagger}\sigma\circ\ldb-{}_{\lambda_{\sigma(1)}}-\dots-{}_{\lambda_{\sigma(n-1)}}-\rdb\circ\sigma^{-1}
\,,
\end{equation}
and, iterating, we get
\begin{equation}\label{eq:nfold-skew2bis}
\ldb-{}_{\lambda_1}-\dots-{}_{\lambda_{n-1}}-\rdb
=(-1)^{s(n+1)}|_{\lambda_n=\lambda_n^\dagger}\sigma^s\circ\ldb-{}_{\lambda_{\sigma^s(1)}}-\dots-{}_{\lambda_{\sigma^s(n-1)}}-\rdb\circ\sigma^{-s}
\,.
\end{equation}
Moreover, we note that the Leibniz rules for $i=1,\dots,n-1$ can be obtained by repeatedly applying skewsymmetry to the Leibniz rule for $i=n$ (cf. Subsection~\ref{def:n-dbr} and Remark \ref{rem:skewPVA}).
\end{remark}
\begin{remark}\label{20250722:rem2}
Let $R\subset \mc V$ be a subset that generates $\mc V$ as a differential algebra. From sesquilinearity \eqref{20140702:eq4}, \eqref{20140702:eq5} and the Leibniz rules \eqref{20140702:eq6} it follows that any $n$-fold $\lambda$-bracket
$Q\in C^n(\mc V)$ is completely determined by its restriction on $R^{\otimes n}$.
\end{remark}

%%%
\section{Further notation and further properties of \texorpdfstring{$2$}{2}-fold \texorpdfstring{$\lambda$}{lambda}-brackets}\label{sec:further}
Let $\ldb-_{\lambda_1}-\dots-_{\lambda_{n-1}}-\rdb\in C^n(\mc V)$, $n\geq1$, and let
$a_1,\dots,a_{n-1}\in\mc V$. Then $\ldb{a_1}_{\lambda_1}\dots{a_{n-1}}_{\lambda_{n-1}}-\rdb:\mc V\to\mc V^{\otimes n}[\lambda_1,\dots,\lambda_{n-1}]$ is an $n$-fold derivation of the associative product of $\mc V$.
Recalling the notation \eqref{20240805:eq1b} we set ($B=b_1\otimes \dots\otimes b_{m}\in\mc V^{\otimes m})$
\begin{equation}\label{badnotation4}
\begin{split}
&\ldb a_1{}_{\lambda_1}\dots a_{n-1}{}_{\lambda_{n-1}}B\rdb_{(i)}
\\
&=b_1\otimes\dots\otimes b_{i-1}\otimes \ldb a_1{}_{\lambda_1}\dots a_{n-1}{}_{\lambda_{n-1}}b_i\rdb
\otimes b_{i+1}\otimes\dots\otimes b_m
\,,
\end{split}
\end{equation}
for $i=1,\dots,m$.
In particular, we denote
\begin{equation}\label{badnotation0}
\begin{split}
&\ldb a_1{}_{\lambda_1}\dots a_{n-1}{}_{\lambda_{n-1}}B\rdb_L
=\ldb a_1{}_{\lambda_1}\dots a_{n-1}{}_{\lambda_{n-1}}B\rdb_{(1)}
\\
&=\ldb a_1{}_{\lambda_1}\dots a_{n-1}{}_{\lambda_{n-1}}b_1\rdb
\otimes (b_2\otimes\dots\otimes b_m)
\,,
\end{split}
\end{equation}
and
\begin{equation}\label{badnotation3}
\begin{split}
&\ldb a_1{}_{\lambda_1}\dots a_{n-1}{}_{\lambda_{n-1}}B
\rdb_R
=
\ldb a_1{}_{\lambda_1}\dots a_{n-1}{}_{\lambda_{n-1}}B
\rdb_{(m)}
\\
&=(b_1\otimes \dots\otimes b_{m-1})\otimes
\ldb a_1{}_{\lambda_1}\dots a_{n-1}{}_{\lambda_{n-1}}b_m
\rdb
\,.
\end{split}
\end{equation}
Applying Lemma \ref{20140606:lem} to the linear map $\ldb a_1{}_{\lambda_1}\dots a_{n-1}{}_{\lambda_{n-1}}-\rdb$ we get the identity ($B\in\mc V^{\otimes m}$)
\begin{equation}\label{20240812:eq1}
\sigma^{m-1} \ldb a_1{}_{\lambda_1}\dots a_{n-1}{}_{\lambda_{n-1}}B\rdb_L
=\ldb a_1{}_{\lambda_1}\dots a_{n-1}{}_{\lambda_{n-1}}\sigma^{m-1}(B)
\rdb_R
\,.
\end{equation}
We generalize \eqref{badnotation0} by letting
\begin{equation}\label{badnotation2}
\begin{split}
&\ldb a_1{}_{\lambda_1}\dots a_{i-1}{}_{\lambda_{i-1}}B
{}_{\lambda_i}\dots a_n\rdb_L
\\
&=
\ldb a_1{}_{\lambda_1}\dots a_{i-1}{}_{\lambda_{i-1}}  b_1
{}_{\lambda_i+x}\dots a_n\rdb\otimes_{n-i} (|_{x=\partial}b_2\otimes\dots\otimes b_n)
\,,
\end{split}
\end{equation}
for every $i=1,\dots,n$. Equation \eqref{badnotation2} is motivated by the fact that,
using the skewsymmetry property \eqref{eq:nfold-skew}, we have
\begin{equation}\label{20240812:skew}
\ldb a_1{}_{\lambda_1}\dots a_{n-1}{}_{\lambda_{n-1}}B\rdb_L
=(-1)^{n+1}\left(
|_{\lambda_n=\lambda_n^\dagger} \ldb a_2{}_{\lambda_2}\dots a_{n-1}{}_{\lambda_{n-1}}B
{}_{\lambda_n}a_1\rdb_L^\sigma
\right)
\,.
\end{equation}
Let $\ldb-_\lambda-\rdb$ be a $2$-fold $\lambda$-bracket on $\mc V$. A special case of \eqref{badnotation2} which generalizes the third equation in \eqref{notation} and that we will use in the sequel is
\begin{equation}\label{badnotation1}
\ldb B_{\lambda} a\rdb_L=\ldb b_1{}_{\lambda+x}a\rdb\otimes_1(|_{x=\partial}b_2\otimes\dots\otimes b_m)
\,,
\end{equation}
where $B=b_1\otimes\dots\otimes b_m\in\mc V^{\otimes m}$ and $a\in\mc V$.

In the sequel, given a $2$-fold $\lambda$-bracket $\ldb-_{\lambda}-\rdb$, we will apply the  notations \eqref{badnotation0}, \eqref{badnotation3} and \eqref{badnotation1} repeatedly using the following convention. A notation of the form $\ldb-_{\lambda}-\rdb_{L,R}$ means that we should first apply
\eqref{badnotation3} to the second entry, and then apply \eqref{badnotation1} to the first entry, namely, for $a,b\in\mc V$, $A\in\mc V^{\otimes m}$ and $B\in\mc V^{\otimes n}$, we have
\begin{equation}\label{20240811:eq2}
\ldb a\otimes A_{\lambda}B\otimes b\rdb_{L,R}
=B\otimes \ldb a\otimes A_{\lambda}b\rdb_L
=B\otimes \left(\ldb a_{\lambda+x}b\rdb\otimes_1(|_{x=\partial}A)\right)
\,.
\end{equation}
Note that, in this particular case, the order of the operations is not important. Indeed
by \eqref{20230811:eq1a} we have
\begin{equation}\label{20240811:eq2bis}
\begin{split}
&B\otimes \left(\ldb a_{\lambda+x}b\rdb\otimes_1(|_{x=\partial}A)\right)
=(B\otimes \ldb a_{\lambda+x}b\rdb)\otimes_1(|_{x=\partial}A)
\\
&=\ldb a_{\lambda+x}B\otimes b\rdb_R\otimes_1(|_{x=\partial}A)
\,.
\end{split}
\end{equation}
However, this is not always the case. For example, we have
\begin{equation}
\begin{split}\label{20240811:eq3}
&\ldb a\otimes A_{\lambda}b\otimes B\rdb_{L,L}
=\ldb a\otimes A_{\lambda}b\rdb_{L}\otimes B
\\
&=\left(\ldb a_{\lambda+x}b\rdb\otimes_1(|_{x=\partial}A)\right)\otimes B
=\ldb a_{\lambda+x}b\otimes B\rdb_L\otimes_{n+1}(|_{x=\partial}A)
\,,
\end{split}
\end{equation}
but we obtain a different result (that is $\ldb a_{\lambda+x}b\otimes B\rdb_L\otimes 1(|_{x=\partial}A)$) if we use \eqref{badnotation1} in the first entry first, and \eqref{badnotation0} in the second entry after.
Note that this convention is consistent with the one used in the RHS of Jacobi identity in Lemma \ref{lem:PVAcoh-pre7}(ii) and in equation \eqref{20230811:eq4}.

In the case when we repeatedly use the notations \eqref{badnotation3} and \eqref{badnotation2} for an $n$-fold $\lambda$-bracket, $n\geq3$, we will specify in which position we should apply the notations as shown in the following examples: if $i<j$
and $B\in\mc V^{\otimes m}$, $C\in\mc V^{\otimes l}$, $a_1,\dots a_n\in\mc V$,
we have
\begin{equation}\label{20240717:eq3}
\begin{split}
&\ldb a_1{}_{\lambda_1}\dots a_{i-1}{}_{\lambda_{i-1}}a_i\otimes B
{}_{\lambda_i}\dots a_{j-1}{}_{\lambda_{j-1}}a_j\otimes C{}_{\lambda_j}\dots
a_{n-1}{}_{\lambda_{n-1}}a_n\rdb_{L_i,L_j}
\\
&=\Big(
\ldb a_1{}_{\lambda_1}\dots a_{i-1}{}_{\lambda_{i-1}}a_i
{}_{\lambda_i+x}\dots a_{j-1}{}_{\lambda_{j-1}}a_j{}_{\lambda_j+y}\dots
\\
&\quad\quad\quad\quad\quad\quad\quad\quad\quad\quad
\dots a_{n-1}{}_{\lambda_{n-1}}a_n\rdb
\otimes_{n-i}(|_{x=\partial}B)
\Big)
\otimes_{n+m-j}(|_{y=\partial} C)
\end{split}
\end{equation}
and
\begin{equation}\label{20240717:eq4}
\begin{split}
&\ldb a_1{}_{\lambda_1}\dots a_{i-1}{}_{\lambda_{i-1}}a_i\otimes B
{}_{\lambda_i}\dots a_{j-1}{}_{\lambda_{j-1}}a_j\otimes C{}_{\lambda_j}\dots
a_{n-1}{}_{\lambda_{n-1}}a_n\rdb_{L_j,L_i}
\\
&=\Big(
\ldb a_1{}_{\lambda_1}\dots a_{i-1}{}_{\lambda_{i-1}}a_i
{}_{\lambda_i+x}\dots a_{j-1}{}_{\lambda_{j-1}}a_j{}_{\lambda_j+y}\dots
\\
&\quad\quad\quad\quad\quad\quad\quad\quad\quad\quad
\dots a_{n-1}{}_{\lambda_{n-1}}a_n\rdb
\otimes_{n-j}(|_{y=\partial}C)
\Big)
\otimes_{n+m-i}(|_{x=\partial} B)
\,.
\end{split}
\end{equation}
The next result will be needed in the proof of Theorem \ref{Thm:g-dPVcoh1}.
\begin{lemma}\label{20240810:lem1}
Let $\ldb-_{\lambda}-\rdb$ be a $2$-fold $\lambda$-bracket on $\mc V$
and let $a,c\in\mc V$, $B=b_1\otimes\dots\otimes b_n \in\mc V^{\otimes n}$.
\begin{enumerate}[(i)]
\item For every $1\leq t\leq n+1$ we have
\begin{equation}\label{20240810:eq1}
\begin{split}
&\ldb a_{\lambda} \ldb B_\mu c\rdb_{L}\rdb_{(t)}
\\
&=\left\{
\begin{array}{ll}
\ldb a_{\lambda}\ldb b_1{}_{\mu+x}c\rdb\rdb_L\otimes_1
(|_{x=\partial}b_2\otimes\dots\otimes b_n)
\,,
&
t=1\,,
\\
\ldb b_1{}_{\lambda+\mu+x}c\rdb\otimes_1
(|_{x=\partial}b_2\otimes\dots\otimes\ldb a_\lambda b_t\rdb\otimes\dots\otimes b_n)
\,,
&
t\neq1,n+1\,,
\\
\ldb a_{\lambda}\ldb b_1{}_{\mu+x}c\rdb\rdb_R\otimes_2
(|_{x=\partial}b_2\otimes\dots\otimes b_n)
\,,
&
t=n+1
\,.
\end{array}
\right.
\end{split}
\end{equation}
\item For every $1\leq t\leq n$ we have
\begin{equation}\label{20240810:eq2}
\begin{split}
&\ldb \ldb a_{\lambda}  B\rdb_{(t)}{}_{\lambda+\mu} c\rdb_{L}
\\
&=\left\{
\begin{array}{ll}
\ldb \ldb a_{\lambda} b_1\rdb{}_{\lambda+\mu+x}c\rdb_L\otimes_1
(|_{x=\partial}b_2\otimes\dots\otimes b_n)
\,,
&
t=1\,,
\\
\ldb b_1{}_{\lambda+\mu+x}c\rdb\otimes_1
(|_{x=\partial}b_2\otimes\dots\otimes\ldb a_\lambda b_t\rdb\otimes\dots\otimes b_n)
\,,
&
t\neq1\,.
\end{array}
\right.
\end{split}
\end{equation}
\item Let us assume that the Jacobi identity \eqref{eq:jacobi2} holds. We have
\begin{equation}\label{20240810:eq3}
\ldb a_{\lambda} \ldb B_\mu c\rdb_{L}\rdb_{(1)}
-\ldb \ldb a_\lambda B\rdb_{(1)}{}_{\lambda+\mu} c\rdb_{L}
=\ldb B_{\mu} \ldb a_\lambda c\rdb \rdb_{L,R}
\end{equation}
and
\begin{equation}\label{20240810:eq4}
\ldb a_{\lambda} \ldb B_\mu c\rdb_{L}\rdb_{(n+1)}
+\ldb \ldb B_\mu a\rdb_{L}{}_{\lambda+\mu} c\rdb_{L}
=\ldb B_{\mu} \ldb a_\lambda c\rdb \rdb_{L,L}
\,.
\end{equation}
\end{enumerate}
\end{lemma}
\begin{proof}
Using \eqref{badnotation1} we have
\begin{equation}\label{20240811:eq1}
\ldb a_\lambda\ldb B_\mu c\rdb_L\rdb_{(t)}
=\ldb a_\lambda\ldb b_1{}_{\mu+x} c\rdb\otimes_1(|_{x=\partial}b_2\otimes \dots\otimes b_n)\rdb_{(t)}
\,.
\end{equation}
Equation \eqref{20240810:eq1} for $t\neq1,n+1$ follows from \eqref{20240811:eq1} by \eqref{eq:sesqui}, 
\eqref{20240805:eq1} for $D=\ldb a_\lambda-\rdb$, and \eqref{eq:tensor-i-notation}.
For $t=1$, equation \eqref{20240811:eq1} becomes
\begin{align*}
&\ldb a_\lambda \ldb b_1{}_{\mu+x}c\rdb'\rdb
\otimes (|_{x=\partial}b_2\otimes \dots\otimes b_n)
\otimes \ldb b_1{}_{\mu+x}c\rdb''
\\
&=\ldb a_\lambda \ldb b_1{}_{\mu+x}c\rdb\rdb_L
\otimes_1(|_{x=\partial}b_2\otimes \dots\otimes b_n)
\,,
\end{align*}
while for $t=n+1$ it becomes
\begin{align*}
&\ldb b_1{}_{\mu+x}c\rdb'
\otimes (|_{x=\partial}b_2\otimes \dots\otimes b_n)
\otimes\ldb a_\lambda \ldb b_1{}_{\mu+x}c\rdb''\rdb
\\
&=\ldb a_\lambda \ldb b_1{}_{\mu+x}c\rdb\rdb_R
\otimes_2(|_{x=\partial}b_2\otimes \dots\otimes b_n)
\,.
\end{align*}
This proves part (i). Part (ii) is proved similarly.
Finally, let us prove part (iii). For brevity we set $\tilde B=b_2\otimes\dots\otimes b_n$. From part (i) and (ii) for $t=1$ we have
\begin{align*}
&\ldb a_{\lambda} \ldb B_\mu c\rdb_{L}\rdb_{(1)}
-\ldb \ldb a_\lambda B\rdb_{(1)}{}_{\lambda+\mu} c\rdb_{L}
\\
&=\left(\ldb a_{\lambda}\ldb b_1{}_{\mu+x}c\rdb\rdb_L 
-\ldb \ldb a_{\lambda} b_1\rdb{}_{\lambda+\mu+x}c\rdb_L
\right)\otimes_1(|_{x=\partial} \tilde B)
\\
&=\ldb b_1{}_{\mu+x}\ldb a_\lambda c\rdb\rdb_R\otimes_1(|_{x=\partial} \tilde B)
=\ldb B_\mu\ldb a_\lambda c\rdb\rdb_{L,R}\,.
\end{align*}
In the second equality above we used the Jacobi identity \eqref{eq:jacobi2}, and in the 
third equality we used \eqref{20240811:eq2bis} and \eqref{20240811:eq2}. This proves equation \eqref{20240810:eq3}. To prove equation
\eqref{20240810:eq4} we first note that
\begin{equation}\label{20240811:eq4}
\begin{split}
&\ldb \ldb B_\mu a\rdb_L{}_{\lambda+\mu}c\rdb_L
=\ldb \ldb b_1{}_{\mu+x}a\rdb'_{\lambda+\mu+x+y}c\rdb\otimes_1
\left((|_{y=\partial}\ldb b_1{}_{\mu+x}a\rdb'')\otimes_1(|_{x=\partial}\tilde B)
\right)
\\
&=\ldb \ldb b_1{}_{\mu+x}a\rdb_{\lambda+\mu+x}c\rdb_L\otimes_2(|_{x=\partial}\tilde B)
\,.
\end{split}
\end{equation}
In the first equality above we used \eqref{badnotation1} twice, in the second equality we used 
Lemma \ref{lem:PVAcoh-pre6} and \eqref{badnotation1} again. Hence, equation \eqref{20240810:eq4}
follows by part (i) for $t=n+1$, equations \eqref{20240811:eq4} and \eqref{20240811:eq3},  and Jacobi identity \eqref{eq:jacobi2}.
\end{proof}

%%%
\section[The variational \MakeLowercase{d}PVA cohomology complex]{The variational double Poisson vertex algebra cohomology complex}\label{sec:var-dPVA}
Let $\mc V$ be a dPVA with $2$-fold $\lambda$-bracket $\llbracket-_\lambda-\rrbracket\in C^2(\mc V)$\glslink{lambda-b-nc}{} and let 
$C(\mc V)$ be the space of $n$-fold $\lambda$-brackets on $\mc V$ constructed in Section \ref{sec:compl-dPVA-1}. (We use square brackets to denote the fixed $2$-fold $\lambda$-bracket on $\mc V$ to distinguish it from an arbitrary  $n$-fold $\lambda$-bracket.)

For $\tint f\in\mc V_\sharp=C^0(\mc V)$ 
we set (cf. \eqref{20140707:eq3b-lie})\glslink{dd}{}
\begin{equation}\label{eq:dP0}
\dd(\tint f)=-[\tint f,-]=-\mult \circ \llbracket f_\lambda-\rrbracket|_{\lambda=0}:\mc V\to\mc V
\,.
\end{equation}
By Theorem \ref{20140707:thm}(b) we have that $[\tint f,-]\in C^1(\mc V)=\Vect(\mc V)^\partial$. Hence,
we have a well defined linear map $d:C^0(\mc V)\to C^1(\mc V)$.
Furthermore, let $Q=\ldb-_{\lambda_1}-\dots-_{\lambda_{n-1}}-\rdb\in C^n(\mc V)$. We define a linear map $\dd(Q):\mc V^{\otimes (n+1)}\to\mc V^{\otimes (n+1)}[\lambda_1,\dots,\lambda_n]$ by the formula
\begin{equation}\label{eq:dP-1}
\begin{split}
&\dd(Q)_{\lambda_1,\dots,\lambda_n}(a_1,\dots,a_{n+1})
\\
&=\sum_{s=1}^{n}(-1)^{n+s+1}
\llbracket a_s {}_{\lambda_s} \ldb a_1{}_{\lambda_1}\stackrel{s}{\check{\dots}}
a_n{}_{\lambda_n}a_{n+1}\rdb \rrbracket_{(s)}
\\
&-\llbracket\ldb a_1{}_{\lambda_1}\dots a_{n-1}{}_{\lambda_{n-1}}a_n\rdb_{\lambda_1+\dots+\lambda_n}a_{n+1}\rrbracket_{L}
\\
&+\sum_{s=1}^n(-1)^{n+s}
\ldb a_1{}_{\lambda_1}\dots a_{s-1}{}_{\lambda_{s-1}}
\llbracket a_{s}{}_{\lambda_s}a_{s+1}\rrbracket_{\lambda_{s}+\lambda_{s+1}}a_{s+2}\dots
{}_{\lambda_n}a_{n+1}\rdb_{L}
\\
&+(-1)^{n+1}\ldb a_2{}_{\lambda_2}\dots a_n{}_{\lambda_n}\llbracket a_1{}_{\lambda_1}a_{n+1}\rrbracket\rdb_R
\,.
\end{split}
\end{equation}
In \eqref{eq:dP-1} we are using the notation \eqref{20240805:eq1} for $D=\llbracket a_s{}_{\lambda_s}-\rrbracket$, and the notations \eqref{badnotation1}, \eqref{badnotation2} and \eqref{badnotation3}.
\begin{example}\label{exa:dP}
The linear map~\eqref{eq:dP-1} is the ``double" analogue of the linear map~\eqref{eq:dH-poly-n}. We show this analogy explicitly for $n=1$. Let $D\in C^1(\mc V)$ be a derivation commuting with $\partial$. From equation \eqref{eq:dP-1} we have
\begin{equation}\label{20240729:eq1}
\dd(D)_{\lambda}(a,b)
=D(\llbracket a_\lambda b \rrbracket)-\llbracket D(a)_\lambda b\rrbracket-
\llbracket a_\lambda D(b)\rrbracket
\,,
\end{equation}
where we are using the notation \eqref{mfold-ext}, namely $D(\llbracket a_\lambda b\rrbracket)
=D_{(1)}(\llbracket a_\lambda b\rrbracket)+D_{(2)}(\llbracket a_\lambda b\rrbracket)$,
which is the ``double" analogue of equation \eqref{20240828:eq1}.
Since $D$ is a derivation of $\mc V$ and, by \eqref{eq:lleibniz}, $\llbracket x_\lambda-\rrbracket:\mc V\to\mc V$ is also a derivation of $\mc V$, for every $x\in\mc V$, it can be checked directly that
$$
\dd (D)_{\lambda}(a,bc)=\dd(D)_{\lambda}(a,b)\, c+b\, \dd(D)_{\lambda}(a,c)
\,,
$$
for every $a,b,c\in\mc V$. Furthermore, since $\llbracket -_\lambda-\rrbracket$ satisfies the skewsymmetry axiom \eqref{eq:skew2}, we have
$$
\dd(D)_\lambda(a,b)=-|_{x=\partial} \dd(D)_{-\lambda-x}(b,a)^\sigma
\,,
$$
thus showing that $\dd(D)\in C^2(\mc V)$ (recall Remark \ref{rem:6.17}) and we have a well defined map
$\dd:C^1(\mc V)\to C^2(\mc V)$. 
Finally, we note that $\dd^2:C^0(\mc V)\to C^2(\mc V)$ is the trivial map. Indeed, combining equations \eqref{eq:dP0} and \eqref{20240729:eq1} we get
$$
\dd^2(\tint a)_{\lambda}(b,c)
=-[a,\llbracket b_\lambda c \rrbracket]+\llbracket[a,b]_\lambda c\rrbracket
+\llbracket b_\lambda [a,c]\rrbracket=0\,,
$$
which vanishes by \eqref{eq:quasijacobi2} since $\llbracket-_\lambda-\rrbracket$ satisfies the Jacobi identity \eqref{eq:jacobi2}.
\end{example}

\subsection{Main statement and its proof}

The next result is the analogue of Theorem \ref{Thm:g-dPcoh1} (and Theorem \ref{thm:PVAcoh})
for dPVA.
\begin{theorem} \label{Thm:g-dPVcoh1}
Let $\mc V$ be a dPVA and let $C(\mc V)$ be the space of $n$-fold $\lambda$-brackets on $\mc V$. Then
equations \eqref{eq:dP0} and \eqref{eq:dP-1} give a well defined map $\dd:C^n(\mc V)\to C^{n+1}(\mc V)$,
for every $n\in\mb Z_{\geq0}$, 
such that $\dd^2=0$.
\end{theorem}
\begin{proof}
We have already shown in Example \ref{exa:dP} that $\dd:C^n(\mc V)\to C^{n+1}(\mc V)$ is well defined for $n=0,1$, and that $\dd^2(\tint a)=0$ for every $\tint a\in C^0(\mc V)$.

Let $n\geq1$. It is straightforward, using \eqref{eq:sesqui},
\eqref{20140702:eq4} and \eqref{20140702:eq5} to verify that $\dd(Q)_{\lambda_1,\dots,\lambda_n}$ 
given in \eqref{eq:dP-1} satisfies the sesquilinearity axioms
\eqref{20140702:eq4} and \eqref{20140702:eq5}.
Furthermore, using \eqref{eq:dP-1} and the Leibniz rule \eqref{20140702:eq6} for $i=n$
we have $(a_1,\dots,a_n,b,c\in\mc V)$
\begin{equation}
\begin{split}\label{20240801:eq1}
&\dd(Q)_{\lambda_1,\dots,\lambda_n}(a_1,\dots,a_n,bc)
\\
&=
(-1)^n\llbracket a_1 {}_{\lambda_1} b \ldb a_2{}_{\lambda_2}\dots
a_n{}_{\lambda_n}c\rdb \rrbracket_{(1)}
+
\sum_{s=2}^{n}(-1)^{n+s+1}
b \llbracket a_s {}_{\lambda_s} \ldb a_1{}_{\lambda_1}\stackrel{s}{\check{\dots}}
a_n{}_{\lambda_n}c\rdb\rrbracket_{(s)}
\\
&+
\sum_{s=1}^{n-1}(-1)^{n+s+1}
\llbracket a_s {}_{\lambda_s} \ldb a_1{}_{\lambda_1}\stackrel{s}{\check{\dots}}
a_n{}_{\lambda_n}b\rdb  \rrbracket_{(s)} c
-
\llbracket a_{n} {}_{\lambda_{n}} \ldb a_1{}_{\lambda_1}\dots
a_{n-1}{}_{\lambda_{n-1}}b\rdb c \rrbracket_{(n)}
\\
&-\llbracket\ldb a_1{}_{\lambda_1}\dots a_{n-1}{}_{\lambda_{n-1}}a_n\rdb_{\lambda_1+\dots+\lambda_n}bc\rrbracket_{(1)}
\\
&+\sum_{s=1}^{n-1}(-1)^{n+s}
b\ldb a_1{}_{\lambda_1}\dots a_{s-1}{}_{\lambda_{s-1}}
\llbracket a_{s}{}_{\lambda_s}a_{s+1}\rrbracket{}_{\lambda_{s}+\lambda_{s+1}}a_{s+2}\dots
{}_{\lambda_n}c\rdb_{L}
\\
&+\sum_{s=1}^{n-1}(-1)^{n+s}
\ldb a_1{}_{\lambda_1}\dots a_{s-1}{}_{\lambda_{s-1}}
\llbracket a_{s}{}_{\lambda_s}a_{s+1}\rrbracket{}_{\lambda_{s}+\lambda_{s+1}}a_{s+2}\dots
{}_{\lambda_n}b\rdb_{L}
\\
&
+
\ldb a_1{}_{\lambda_1}\dots a_{n-1}{}_{\lambda_{n-1}}
\llbracket a_{n}{}_{\lambda_n}bc\rrbracket\rdb_{L}
+(-1)^{n+1}\ldb a_2{}_{\lambda_2}\dots a_n{}_{\lambda_n}
\llbracket a_1{}_{\lambda_1}bc\rrbracket\rdb_R
\,.
\end{split}
\end{equation}
Using the Leibniz rule \eqref{eq:lleibniz} for $\llbracket-_\lambda-\rrbracket$
and the definition of $\llbracket a_1{}_\lambda-\rrbracket_{(1)}$ (see \eqref{20240805:eq1}) we get
\begin{equation}\label{20240801:eq2}
\begin{split}
&\llbracket a_1 {}_{\lambda_1} b \ldb a_2{}_{\lambda_2}\dots
a_n{}_{\lambda_n}c\rdb\rrbracket_{(1)}\\
&=
b\llbracket a_1 {}_{\lambda_1} \ldb a_2{}_{\lambda_2}\dots
a_n{}_{\lambda_n}c\rdb\rrbracket_{(1)}
+
\llbracket a_1 {}_{\lambda_1} b\rrbracket\ldb a_2{}_{\lambda_2}\dots
a_n{}_{\lambda_n}c\rdb
\,.
\\
\end{split}
\end{equation}
Similarly, we have
\begin{equation}\label{20240801:eq3}
\begin{split}
&\llbracket a_{n} {}_{\lambda_{n}} \ldb a_1{}_{\lambda_1}\dots
a_{n-1}{}_{\lambda_{n-1}}b\rdb c \rrbracket_{(n)}
\\
&=
\llbracket a_{n} {}_{\lambda_{n}} \ldb a_1{}_{\lambda_1}\dots
a_{n-1}{}_{\lambda_{n-1}}b\rdb  \rrbracket_{(n)}c
+\ldb a_1{}_{\lambda_1}\dots
a_{n-1}{}_{\lambda_{n-1}}b\rdb
\llbracket a_{n} {}_{\lambda_{n}} c\rrbracket
\,.
\\
\end{split}
\end{equation}
Furthermore, using again the Leibniz rule \eqref{eq:lleibniz} for $\llbracket-_\lambda-\rrbracket$, the Leibniz rule
\eqref{20140702:eq6} for $i=n$ and the notations
\eqref{badnotation0} and \eqref{badnotation3} we get
\begin{equation}\label{20240801:eq4}
\begin{split}
&
\ldb a_1{}_{\lambda_1}\dots a_{n-1}{}_{\lambda_{n-1}}
\llbracket a_{n}{}_{\lambda_n}bc \rrbracket\rdb_{L}
=
\ldb a_1{}_{\lambda_1}\dots a_{n-1}{}_{\lambda_{n-1}}
\llbracket a_{n}{}_{\lambda_n}b \rrbracket\rdb_{L}c
\\
&
+\ldb a_1{}_{\lambda_1}\dots a_{n-1}{}_{\lambda_{n-1}}b\rdb
\llbracket a_{n}{}_{\lambda_n}c\rrbracket
+b\ldb a_1{}_{\lambda_1}\dots a_{n-1}{}_{\lambda_{n-1}}
\llbracket a_{n}{}_{\lambda_n}c\rrbracket\rdb_{L}
\\
\end{split}
\end{equation}
and
\begin{equation}\label{20240801:eq5}
\begin{split}
&
\ldb a_2{}_{\lambda_2}\dots a_n{}_{\lambda_n}
\llbracket a_1{}_{\lambda_1}bc\rrbracket\rdb_R
=
b\ldb a_2{}_{\lambda_2}\dots a_n{}_{\lambda_n}
\llbracket a_1{}_{\lambda_1}c \rrbracket\rdb_R
\\
&+
\llbracket a_1{}_{\lambda_1}b \rrbracket\ldb a_2{}_{\lambda_2}\dots a_n{}_{\lambda_n}c\rdb_R
+
\ldb a_2{}_{\lambda_2}\dots a_n{}_{\lambda_n}
\llbracket a_1{}_{\lambda_1}b \rrbracket\rdb_Rc
\,.
\\
\end{split}
\end{equation}
Using equations \eqref{20240801:eq2}, \eqref{20240801:eq3}, \eqref{20240801:eq4}
and \eqref{20240801:eq5} in \eqref{20240801:eq1} we have that $d(Q)$ satisfies
the Leibniz rule in its last entry ($a_1,\dots,a_n,b,c\in\mc V$):
\begin{align*}
&\dd(Q)_{\lambda_1,\dots,\lambda_n}(a_1,\dots,a_n,bc)
=b\, \dd(Q)_{\lambda_1,\dots,\lambda_n}(a_1,\dots,a_n,c)
\\
&+\dd(Q)_{\lambda_1,\dots,\lambda_n}(a_1,\dots,a_n,b)\,c
\,.
\end{align*}
Hence, by Remark \ref{rem:6.17}, we have that $\dd(Q)\in C^{n+1}(\mc V)$ provided that it satisfies the skewsymmetry axiom \eqref{eq:nfold-skew}. To this aim, we start by noticing that, using \eqref{20240812:skew} we have
\begin{equation}\label{skew1}
\begin{split}
&\llbracket a_1 {}_{\lambda_1} \ldb a_2{}_{\lambda_2}\dots
a_n{}_{\lambda_n}a_{n+1}\rdb \rrbracket_{L}
=-|_{x=\partial}\sigma\dsq{\ldb a_2{}_{\lambda_2}\dots
a_n{}_{\lambda_n}a_{n+1}\rdb_{-\lambda_1-x}a_1}_{L}
\\
&=-|_{\lambda_{n+1}=\lambda_{n+1}^\dagger}
\sigma \llbracket\ldb a_2{}_{\lambda_2}\dots a_n{}_{\lambda_n}a_{n+1}\rdb_{\lambda_2+\dots+\lambda_{n+1}}a_1\rrbracket_{L}
\,,
\end{split}
\end{equation}
where in the second equality we used \eqref{eq:dagger}. Furthermore, using skewsymmetry \eqref{eq:nfold-skew} and sesquilinearity \eqref{eq:sesqui} we have ($s=2,\dots,n$)
\begin{equation}\label{skew2}
\begin{split}
& \llbracket a_s {}_{\lambda_s} \ldb a_1{}_{\lambda_1}\stackrel{s}{\check{\dots}}
a_n{}_{\lambda_n}a_{n+1}\rdb \rrbracket_{(s)}
\\
&=(-1)^{n+1}|_{\lambda_{n+1}=\lambda_{n+1}^\dagger}
\llbracket a_s {}_{\lambda_s} \ldb a_2{}_{\lambda_2}\stackrel{s-1}{\check{\dots}}
a_{n+1}{}_{\lambda_{n+1}}a_{1}\rdb^\sigma \rrbracket_{(s)}
\\
&=(-1)^{n+1}|_{\lambda_{n+1}=\lambda_{n+1}^\dagger}\sigma
\llbracket a_s {}_{\lambda_s} \ldb a_2{}_{\lambda_2}\stackrel{s-1}{\check{\dots}}
a_{n+1}{}_{\lambda_{n+1}}a_{1}\rdb \rrbracket_{(s-1)}
\,,
\end{split}
\end{equation}
where in the last equality we used the first equation in \eqref{20140606:eq4}.
Similarly, using again skewsymmetry \eqref{eq:nfold-skew} and sesquilinearity \eqref{eq:sesqui} we have
\begin{equation}\label{skew3}
\begin{split}
&\llbracket \ldb a_1{}_{\lambda_1}\dots a_{n-1}{}_{\lambda_{n-1}}a_n\rdb_{\lambda_1+\dots+\lambda_n}a_{n+1}\rrbracket_{L}
\\
&=(-1)^{n+1}\llbracket\ldb a_2{}_{\lambda_2}\dots a_{n}{}_{\lambda_{n}}a_1\rdb^\sigma_{\lambda_1+\dots+\lambda_n}a_{n+1}\rrbracket_{L}
\\
&=(-1)^n|_{\lambda_{n+1}=\lambda_{n+1}^\dagger}\sigma^{-1}
\llbracket a_{n+1}{}_{\lambda_{n+1}}\ldb a_2{}_{\lambda_2}\dots a_{n}{}_{\lambda_{n}}a_1\rdb^\sigma\rrbracket_{(1)}
\\
&=(-1)^n|_{\lambda_{n+1}=\lambda_{n+1}^\dagger}\sigma
\llbracket a_{n+1}{}_{\lambda_{n+1}}\ldb a_2{}_{\lambda_2}\dots a_{n}{}_{\lambda_{n}}a_1\rdb\rrbracket_{(n)}
\,,
\end{split}
\end{equation}
where in the second equality we used \eqref{20240812:skew} and in the last equality we used the second equation in \eqref{20140606:eq4}. In a similar fashion, using
skewsymmetry, sesquilinearity and equations
\eqref{20240812:eq1} and \eqref{20240812:skew} we get the identities
\begin{equation}\label{skew4}
\begin{split}
&\ldb \llbracket a_{1}{}_{\lambda_1}a_{2}\rrbracket_{\lambda_{1}+\lambda_{2}}a_{3}{}_{\lambda_3}\dots
a_{n}{}_{\lambda_n}a_{n+1}\rdb_{L}
\\
&=(-1)^n|_{\lambda_{n+1}=\lambda_{n+1}^\dagger}\sigma
\ldb a_3{}_{\lambda_3}\dots a_{n+1}{}_{\lambda_{n+1}}\dsq{a_2{}_{\lambda_2}a_{1}}\rdb_R
\,,
\end{split}
\end{equation}
\begin{equation}\label{skew5}
\begin{split}
&\ldb a_1{}_{\lambda_1}\dots a_{s-1}{}_{\lambda_{s-1}}
\llbracket a_{s}{}_{\lambda_s}a_{s+1}\rrbracket_{\lambda_{s}+\lambda_{s+1}}a_{s+2}\dots
{}_{\lambda_n}a_{n+1}\rdb_{L}
\\
&
=(-1)^{n+1}|_{\lambda_{n+1}=\lambda_{n+1}^\dagger}
\ldb a_{2}{}_{\lambda_2}\dots a_{s-2}{}_{\lambda_{s-2}}
\llbracket a_{s-1}{}_{\lambda_{s-1}}a_{s}\rrbracket_{\lambda_{s-1}+\lambda_{s}}a_{s+1}\dots
a_{n+1}{}_{\lambda_{n+1}}a_1\rdb_{L}^\sigma
\,,
\end{split}
\end{equation}
for $s=2,\dots,n$, and
\begin{equation}\label{skew6}
\begin{split}
&\ldb a_2{}_{\lambda_2}\dots a_n{}_{\lambda_n}
\llbracket a_1{}_{\lambda_1}a_{n+1}\rrbracket\rdb_R
\\
&=-|_{\lambda_{n+1}=\lambda_{n+1}^\dagger}\sigma
\ldb a_2{}_{\lambda_2}\dots a_n{}_{\lambda_n}
\llbracket a_{n+1}{}_{\lambda_{n+1}}a_{1}\rrbracket\rdb_L
\,.
\end{split}
\end{equation}
Using equations \eqref{skew1}-\eqref{skew6} in \eqref{eq:dP-1} we get the identity
$$
\dd(Q)_{\lambda_1,\dots,\lambda_n}(a_1,\dots,a_{n+1})
=(-1)^n|_{\lambda_{n+1}=\lambda_{n+1}^\dagger}\sigma \left(
\dd(Q)_{\lambda_2,\dots,\lambda_{n+1}}(a_2,\dots,a_{n+1},a_1)
\right)
$$
thus showing that $\dd(Q)$ satisfies skewsymmetry \eqref{eq:nfold-skew} and concluding the proof that $\dd(Q)\in C^{n+1}(\mc V)$.

We are left to show that $\dd^2(Q)=0$. This will be done similarly to the analogous computation in the proof of Theorem \ref{thm:PVAcoh}.
To this aim, 
for $Q
=\ldb -{}_{\lambda_1}-\dots -{}_{\lambda_{n-1}}-\rdb\in C^n(\mc V)$
we denote
\begin{equation}\label{20240717:eq0}
Q_{\lambda_1,\dots,\lambda_{n-1}}(a_1,\dots,a_{n})=\ldb a_1{}_{\lambda_1}\dots a_{n-1}{}_{\lambda_{n-1}}a_n\rdb
\,.
\end{equation}
Hence, using the notation \eqref{eq:XL} we have
\begin{align*}
&Q^{(s)}_{\lambda_1,\dots,\lambda_{n-1}}(a_1,\dots,a_{s-1},a_s\otimes A,a_{s+1},\dots,a_n)
\\
&=Q_{\lambda_1,\dots,\lambda_s+x,\dots,\lambda_{n-1}}(a_1,\dots,a_s,\dots,a_n)\otimes_{n-s}(|_{x=\partial}A)\,,
\end{align*}
for $1\leq s\leq n$, $a_1,\dots,a_n\in\mc V$ and $A\in \mc V^{\otimes m}$,
so that
\begin{equation}\label{20240717:eq1}
\begin{split}
&Q^{(s)}_{\lambda_1,\dots,\lambda_{s-1},\lambda_{s}+\lambda_{s+1},\lambda_{s+2},\dots,\lambda_n}
(a_1,\dots, a_{s-1},\llbracket a_{s}{}_{\lambda_s}a_{s+1}\rrbracket,a_{s+2},\dots,a_{n+1})
\\
&=\ldb a_1{}_{\lambda_1}\dots a_{s-1}{}_{\lambda_{s-1}}
\llbracket a_{s}{}_{\lambda_s}a_{s+1}\rrbracket_{\lambda_{s}+\lambda_{s+1}}a_{s+2}\dots
{}_{\lambda_n}a_{n+1}\rdb_{L}
\,.
\end{split}
\end{equation}
Using the notation \eqref{20240805:eq1} for $D=Q_{\lambda_1,\dots,\lambda_{n-1}}(a_1,\dots,a_{n-1},-):\mc V\to \mc V^{\otimes n}[\lambda_1,\dots,\lambda_{n-1}]$ we also have
$$
(Q_R)_{\lambda_1,\dots,\lambda_{n-1}}(a_1,\dots, A\otimes a_n)
=
A\otimes
Q_{\lambda_1,\dots,\lambda_{n-1}}(a_1,\dots, a_n)
\,,
$$
for $a_1,\dots,a_n\in\mc V$ and $A\in\mc V^{\otimes m}$, so that
\begin{equation}\label{20240717:eq2}
(Q_{\lambda_2,\dots,\lambda_n})_R(a_2,\dots, a_n,\llbracket a_{1}{}_{\lambda_{1}}a_{n+1}\rrbracket)
=\ldb a_2{}_{\lambda_2}\dots a_n{}_{\lambda_n}\llbracket a_1{}_{\lambda_1}a_{n+1}\rrbracket\rdb_R
\,.
\end{equation}
Then, using \eqref{20240717:eq0}, \eqref{20240717:eq1} and
\eqref{20240717:eq2} we can rewrite equation \eqref{eq:dP-1}
as follows
\begin{equation}\label{eq:dP-1-Q}
\begin{split}
&\dd(Q)_{\lambda_1,\dots,\lambda_n}(a_1,\dots,a_{n+1})
\\
&=\sum_{s=1}^{n}(-1)^{n+s+1}
\llbracket a_s {}_{\lambda_s}
Q_{\lambda_1,\stackrel{s}{\check{\dots}},\lambda_n}(a_1{},\stackrel{s}{\check{\dots}},
a_n,a_{n+1})\rrbracket_{(s)}
\\
&-\llbracket Q_{\lambda_1,\dots,\lambda_{n-1}}(a_1,\dots,a_n)_{\lambda_1+\dots+\lambda_n}a_{n+1}\rrbracket_{L}
\\
&+\sum_{s=1}^n(-1)^{n+s}
Q^{(s)}_{\lambda_1,\dots,\lambda_{s-1},\lambda_{s}+\lambda_{s+1},\lambda_{s+2},\dots,\lambda_n}
(a_1,\dots, a_{s-1},\llbracket a_{s}{}_{\lambda_s}a_{s+1}\rrbracket,a_{s+2},\dots,a_{n+1})
\\
&+(-1)^{n+1}
(Q_R)_{\lambda_2,\dots,\lambda_n}(a_2,\dots, a_n,\llbracket a_{1}{}_{\lambda_{1}}a_{n+1}\rrbracket)
\,.
\end{split}
\end{equation}
In \eqref{eq:dP-1-Q} we are using the notation introduced in \eqref{Eq:NotCheck1}. From equation \eqref{eq:dP-1-Q} we get
\begin{equation}
\begin{split}\label{20240808:eq1}
&\dd^2(Q)_{\lambda_1,\dots,\lambda_{n+1}}(a_1,\dots,a_{n+2})\\
&=\sum_{t=1}^{n+1}(-1)^{n+t}
\llbracket a_t {}_{\lambda_t}
\dd (Q)_{\lambda_1,\stackrel{t}{\check{\dots}},\lambda_{n+1}}(a_1{},\stackrel{t}{\check{\dots}},
a_{n+1},a_{n+2})\rrbracket_{(t)}
\\
&-\llbracket \dd(Q)_{\lambda_1,\dots,\lambda_{n}}(a_1,\dots,a_{n+1})_{\lambda_1+\dots+\lambda_{n+1}}a_{n+2}\rrbracket_{L}
\\
&+\sum_{t=1}^{n+1}(-1)^{n+t+1}
\dd(Q)^{(t)}_{\lambda_1,\dots,\lambda_{t-1},\lambda_{t}+\lambda_{t+1},\lambda_{t+2},\dots,\lambda_{n+1}}
(a_1,\dots
\\
&\quad\quad\quad\quad\quad\quad\quad\quad\quad\quad\quad\quad\quad\quad\quad\quad
\dots, a_{t-1},\llbracket a_{t}{}_{\lambda_t}a_{t+1}\rrbracket,a_{t+2},\dots,a_{n+2})
\\
&+(-1)^{n}
(\dd(Q)_R)_{\lambda_2,\dots,\lambda_{n+1}}(a_2,\dots, a_{n+1},\llbracket a_1{}_{\lambda_1}a_{n+2}\rrbracket)
\,.
\end{split}
\end{equation}
Using again \eqref{eq:dP-1-Q} and the notations \eqref{20240717:eq0}-\eqref{20240717:eq2}, we have, for $1\leq t\leq n+1$:
\begin{equation}
\begin{split}\label{20240808:eq2}
&\dd(Q)_{\lambda_1,\stackrel{t}{\check{\dots}},\lambda_{n+1}}(a_1{},\stackrel{t}{\check{\dots}},
a_{n+1},a_{n+2})
\\
&=\sum_{s=1}^{t-1}(-1)^{n+s+1}
\llbracket a_s {}_{\lambda_s}
\ldb a_1{}_{\lambda_1}\stackrel{s}{\check{\dots}}\stackrel{t}{\check{\dots}}
a_{n+1}{}_{\lambda_{n+1}}a_{n+2}\rdb\rrbracket_{(s)}
\\
&+\sum_{s=t+1}^{n+1}(-1)^{n+s}
\llbracket a_s {}_{\lambda_s}
Q_{\lambda_1,\stackrel{t}{\check{\dots}}\stackrel{s}{\check{\dots}},\lambda_{n+1}}
\ldb a_1{}_{\lambda_1}\stackrel{t}{\check{\dots}}\stackrel{s}{\check{\dots}}
a_{n+1}{}_{\lambda_{n+1}}a_{n+2}\rdb\rrbracket_{(s-1)}
\\
&-\delta_{t\neq n+1}\llbracket
\ldb a_1{}_{\lambda_1}\stackrel{t}{\check{\dots}}a_n{}_{\lambda_n}a_{n+1}\rdb_{\lambda_1+\stackrel{t}{\check{\dots}}+\lambda_{n+1}}a_{n+2}\rrbracket_{L}
\\
&-\delta_{t,n+1}\llbracket
\ldb a_1{}_{\lambda_1}\dots a_{n-1}{}_{\lambda_{n-1}}a_{n}\rdb_{\lambda_1+\dots+\lambda_n}a_{n+2}\rrbracket_{L}
\\
&+\sum_{s=1}^{t-2}(-1)^{n+s}
\ldb a_1{}_{\lambda_1}\dots a_{s-1}{}_{\lambda_{s-1}}
\llbracket a_{s}{}_{\lambda_s}a_{s+1}\rrbracket_{\lambda_{s}+\lambda_{s+1}}a_{s+2}{}_{\lambda_{s+2}}\stackrel{t-1}{\check{\dots}}a_{n+1}{}_{\lambda_{n+1}}a_{n+2}\rdb_L
\\
&+(-1)^{n+t+1}
\ldb a_1{}_{\lambda_1}\dots a_{t-2}{}_{\lambda_{t-2}}
\llbracket a_{t-1}{}_{\lambda_{t-1}}a_{t+1}\rrbracket_{\lambda_{t-1}+\lambda_{t+1}}a_{t+2}{}_{\lambda_{t+2}}\dots a_{n+1}{}_{\lambda_{n+1}}a_{n+2}\rdb_L
\\
&+\sum_{s=t+1}^{n+1}(-1)^{n+s+1}
\ldb a_1{}_{\lambda_1}\stackrel{t}{\check{\dots}} a_{s-1}{}_{\lambda_{s-1}}
\llbracket a_{s}{}_{\lambda_s}a_{s+1}\rrbracket_{\lambda_{s}+\lambda_{s+1}}a_{s+2}{}_{\lambda_{s+2}}\dots a_{n+1}{}_{\lambda_{n+1}}a_{n+2}\rdb_L
\\
&+(-1)^{n+1}\delta_{t,1}
\ldb a_3{}_{\lambda_3}\dots a_{n+1}{}_{\lambda_{n+1}}
\llbracket a_{2}{}_{\lambda_{2}}a_{n+2}\rrbracket\rdb_R
\\
&+(-1)^{n+1}\delta_{t\neq1}
\ldb a_2{}_{\lambda_2}\stackrel{t-1}{\check{\dots}} a_{n+1}{}_{\lambda_{n+1}}
\llbracket a_{1}{}_{\lambda_{1}}a_{n+2}\rrbracket\rdb_R
\,.
\end{split}
\end{equation}
Moreover, with a similar computation to the one performed to derive equation \eqref{eq:dPsquare2}, and recalling the notations \eqref{20240717:eq0}-\eqref{20240717:eq2}
we get ($1\leq t\leq n+1$)
\begin{equation}
\begin{split}\label{20240808:eq3}
&\dd(Q)^{(t)}_{\lambda_1,\dots,\lambda_{t-1},\lambda_{t}+\lambda_{t+1},\lambda_{t+2},\dots,\lambda_{n+1}}
(a_1,\dots, a_{t-1},\llbracket a_{t}{}_{\lambda_t}a_{t+1} \rrbracket,a_{t+2},\dots,a_{n+2})
\\
&=\sum_{s=1}^{t-1}(-1)^{n+s+1}
\llbracket a_s{}_{\lambda_s}
\ldb a_1{}_{\lambda_1}\stackrel{s}{\check{\dots}} a_{t-1}{}_{\lambda_{t-1}}\dsq{ a_{t}{}_{\lambda_t}a_{t+1} }_{\lambda_{t}+\lambda_{t+1}}a_{t+2}{}_{\lambda_{t+2}}\dots
\\
&\quad\quad\quad\quad\quad\quad
\dots a_{n+1}{}_{\lambda_{n+1}}a_{n+2}\rdb_L
\rrbracket_{(s)}
\\
&+(-1)^{n+t+1}
\llbracket
\llbracket a_t{}_{\lambda_t}a_{t+1}\rrbracket_{\lambda_t+\lambda_{t+1}}
\ldb a_1{}_{\lambda_1}\dots a_{t-1}{}_{\lambda_{t-1}}a_{t+2}{}_{\lambda_{t+2}}\dots
a_{n+1}{}_{\lambda_{n+1}}a_{n+2} \rdb
\rrbracket_{L,(t)}
\\
&+\sum_{s=t+2}^{n+1}(-1)^{n+s}
\llbracket
a_s{}_{\lambda_s}
\ldb a_1{}_{\lambda_1}\dots a_{t-1}{}_{\lambda_{t-1}}
\llbracket a_{t}{}_{\lambda_t}a_{t+1}\rrbracket_{\lambda_{t}+\lambda_{t+1}}a_{t+2}{}_{\lambda_{t+2}}
\stackrel{s-1}{\check{\dots}}
\\
&\quad\quad\quad\quad\quad\quad
\dots a_{n+1}{}_{\lambda_{n+1}}a_{n+2}\rdb
\rrbracket_{(s)}
\\
&
-
\delta_{t\neq n+1}\llbracket
\ldb
a_1{}_{\lambda_1}\dots a_{t-1}{}_{\lambda_{t-1}}
\llbracket a_{t}{}_{\lambda_t}a_{t+1}\rrbracket_{\lambda_{t}+\lambda_{t+1}}
a_{t+2}{}_{\lambda_{t+2}}\dots
\\
&\quad\quad\quad\quad\quad\quad
\dots a_{n}{}_{\lambda_n}a_{n+1}
\rdb
_{\lambda_1+\dots+\lambda_{n+1}}a_{n+2}
\rrbracket_L
\\
&
-
\delta_{t,n+1}\llbracket
\ldb a_1{}_{\lambda_1}\dots a_{n-1}{}_{\lambda_{n-1}} a_{n}\rdb
_{\lambda_1+\dots+\lambda_{n}}\llbracket a_{n+1}{}_{\lambda_{n+1}}a_{n+2}\rrbracket
\rrbracket_{L,L}
\\
&
+\sum_{s=1}^{t-2}(-1)^{n+s}
\ldb a_1{}_{\lambda_1}\dots a_{s-1}{}_{\lambda_{s-1}}
\llbracket a_s{}_{\lambda_s} a_{s+1}\rrbracket_{\lambda_{s}+\lambda_{s+1}}a_{s+2}{}_{\lambda_{s+2}}\dots
\\
&\quad\quad\quad\quad\quad\quad
\dots a_{t-1}{}_{\lambda_{t-1}}\llbracket a_t{}_{\lambda_t}a_{t+1}\rrbracket_{\lambda_{t}+\lambda_{t+1}}
a_{t+2}{}_{\lambda_{t+2}}\ldots a_{n+1}{}_{\lambda_{n+1}}a_{n+2} \rdb_{L_s,L_{t-1}}
\\
&
+(-1)^{n+t+1}
\ldb a_1{}_{\lambda_1}\dots a_{t-2}{}_{\lambda_{t-2}}
\llbracket a_{t-1}{}_{\lambda_{t-1}}\llbracket a_t {}_{\lambda_t} a_{t+1}\rrbracket \rrbracket_L{}_{\lambda_{t-1}+\lambda_t+\lambda_{t+1}} a_{t+2}{}_{\lambda_{t+2}}\dots
\\
&\quad\quad\quad\quad\quad\quad
\dots a_{n+1}{}_{\lambda_{n+1}}a_{n+2}\rdb_L
\\
&
+(-1)^{n+t}
\ldb
a_1{}_{\lambda_1}\dots a_{t-1}{}_{\lambda_{t-1}}\llbracket \llbracket a_{t}{}_{\lambda_{t}} a_{t+1}\rrbracket_{\lambda_t+\lambda_{t+1}} a_{t+2}\rrbracket_L{}_{\lambda_{t}+\lambda_{t+1}+\lambda_{t+2}}
a_{t+3}{}_{\lambda_{t+3}}\dots
\\
&\quad\quad\quad\quad\quad\quad
\dots a_{n+1}{}_{\lambda_{n+1}}a_{n+2}
\rdb_L
\\
&
+\sum_{s=t+2}^{n+1}(-1)^{n+s+1}
\ldb a_1{}_{\lambda_1}\dots a_{t-1}{}_{\lambda_{t-1}}
\llbracket a_t{}_{\lambda_t} a_{t+1}\rrbracket_{\lambda_{t}+\lambda_{t+1}}a_{t+2}{}_{\lambda_{t+2}}\dots
\\
&\quad\quad\quad\quad\quad\quad
\dots a_{s-1}{}_{\lambda_{s-1}}\llbracket a_s{}_{\lambda_s}a_{s+1}\rrbracket_{\lambda_{s}+\lambda_{s+1}}
a_{s+2}{}_{\lambda_{s+2}}\ldots a_{n+1}{}_{\lambda_{n+1}}a_{n+2}\rdb_{L_{s-1},L_t}
\\
&+(-1)^{n+1}\delta_{t,1}
\ldb a_3{}_{\lambda_3}\dots a_{n+1}{}_{\lambda_{n+1}}\llbracket \llbracket a_1{}_{\lambda_1}a_2\rrbracket_{\lambda_1+\lambda_2}a_{n+2}\rrbracket_L\rdb_R
\\
&
+(-1)^{n+1}\delta_{t\neq1,n+1}
\ldb a_2{}_{\lambda_2}\dots a_{t-1}{}_{\lambda_{t-1}}
\llbracket a_{t}{}_{\lambda_t} a_{t+1}\rrbracket_{\lambda_t+\lambda_{t+1}}a_{t+2}{}_{\lambda_{t+2}}\dots \\
&\quad\quad\quad\quad\quad\quad
\dots a_{n+1}{}_{\lambda_{n+1}}\llbracket a_1{}_{\lambda_1} a_{n+2}\rrbracket\rdb_{L_{t-1},R}
\\
&
+(-1)^{n+1}
\ldb a_2{}_{\lambda_2}\dots a_{n}{}_{\lambda_{n}}\llbracket a_1{}_{\lambda_1}\llbracket a_{n+1}{}_{\lambda_{n+1}} a_{n+2}\rrbracket\rrbracket_L\rdb_{(2)}
\,.
\end{split}
\end{equation}
In \eqref{20240808:eq3} we are using the notations and conventions introduced in Section \ref{sec:further} (see, in particular, equations \eqref{badnotation4}, \eqref{20240811:eq3}, \eqref{20240717:eq3} and \eqref{20240717:eq4}).

Finally, using \eqref{eq:dP-1-Q} and some simple algebraic manipulations, we can also compute explicitly the second and last term in the
RHS of \eqref{20240808:eq1}. We get
\begin{equation}
\begin{split}\label{20240808:eq4}
&-\llbracket \dd(Q)_{\lambda_1,\dots,\lambda_{n}}(a_1,\dots,a_{n+1})_{\lambda_1+\dots+\lambda_{n+1}}a_{n+2}\rrbracket_{L}
%\\
%&
%=\sum_{t=1}^{n}(-1)^{n+t}
%\llbracket
%{\llbracket a_t {}_{\lambda_t}
%Q_{\lambda_1,\stackrel{t}{\check{\dots}},\lambda_n}(a_1{},\stackrel{t}{\check{\dots}},
%a_n,a_{n+1})\rrbracket_{(t)}}_{\lambda_1+\dots+\lambda_{n+1}}a_{n+2}
%\rrbracket_{L}
%\\
%&+\llbracket
%{\llbracket Q_{\lambda_1,\dots,\lambda_{n-1}}(a_1,\dots,a_n)_{\lambda_1+\dots+\lambda_n}a_{n+1}\rrbracket_{L}}_{\lambda_1+\dots+\lambda_{n+1}}a_{n+2}
%\rrbracket_{L}
%\\
%&+\sum_{t=1}^n(-1)^{n+t+1}
%\llbracket
%Q^{(t)}_{\lambda_1,\dots,\lambda_{t-1},\lambda_{t}+\lambda_{t+1},\lambda_{t+2},\dots,\lambda_n}
%(a_1,\dots, a_{t-1},\llbracket a_{t}{}_{\lambda_t}a_{t+1}\rrbracket,a_{t+2},\dots,a_{n+1})
%_{\lambda_1+\dots+\lambda_{n+1}}a_{n+2}
%\rrbracket_{L}
%\\
%&+(-1)^{n}
%\llbracket
%(Q_{\lambda_2,\dots,\lambda_n})_R(a_2,\dots, a_n,\llbracket a_{1}{}_{\lambda_{1}}a_{n+1}\rrbracket)
%_{\lambda_1+\dots+\lambda_{n+1}}a_{n+2}
%\rrbracket_{L}
\\
&
=\sum_{t=1}^{n}(-1)^{n+t}
\llbracket
{\llbracket a_t {}_{\lambda_t}
\ldb a_1{}_{\lambda_1}\stackrel{t}{\check{\dots}}
a_n{}_{\lambda_n}a_{n+1}\rdb\rrbracket_{(t)}}_{\lambda_1+\dots+\lambda_{n+1}}a_{n+2}
\rrbracket_{L}
\\
&+\llbracket
{\llbracket
\ldb a_1{}_{\lambda_1}\dots a_{n-1}{}_{\lambda_{n-1}}a_n\rdb_{\lambda_1+\dots+\lambda_n}a_{n+1}\rrbracket_{L}}_{\lambda_1+\dots+\lambda_{n+1}}a_{n+2}
\rrbracket_{L}
\\
&+\sum_{t=1}^n(-1)^{n+t+1}
\llbracket
\ldb a_1{}_{\lambda_1}\dots a_{t-1}{}_{\lambda_{t-1}}\llbracket a_{t}{}_{\lambda_t}a_{t+1}\rrbracket_{\lambda_{t}+\lambda_{t+1}}a_{t+2}{}_{\lambda_{t+2}}\dots
\\
&\quad\quad\quad\quad\quad\quad\quad\quad\quad\quad\quad\quad\quad
{\dots a_n{}_{\lambda_n}a_{n+1}\rdb_L}
_{\lambda_1+\dots+\lambda_{n+1}}a_{n+2}
\rrbracket_{L}
\\
&+(-1)^{n}
\llbracket
{\ldb a_2{}_{\lambda_2}\dots a_n{}_{\lambda_n}\llbracket a_{1}{}_{\lambda_{1}}a_{n+1}\rrbracket
\rdb_R}
_{\lambda_1+\dots+\lambda_{n+1}}a_{n+2}
\rrbracket_{L}
\,,
\end{split}
\end{equation}
and
\begin{equation}
\begin{split}\label{20240808:eq5}
&(-1)^{n}
(\dd(Q)_{\lambda_2,\dots,\lambda_{n+1}})_R(a_2,\dots, a_{n+1},\llbracket a_1{}_{\lambda_1}a_{n+2}\rrbracket)
%\\
%&
%=\sum_{t=2}^{n+1}(-1)^t\llbracket a_{t}{}_{\lambda_t}
%(Q_{\lambda_2,\stackrel{t-1}{\check{\dots}},\lambda_{n+1}})_R
%(a_2,\stackrel{t-1}{\check{\dots}},a_{n+1},\llbracket a_1{}_{\lambda_1}a_{n+2}\rrbracket)
%\rrbracket_{(t)}
%\\
%&
%+(-1)^{n+1}
%\llbracket Q_{\lambda_2,\dots,\lambda_n}(a_2,\dots,a_{n+1})_{\lambda_2+\dots+\lambda_{n+1}}
%\llbracket a_1{}_{\lambda_1}a_{n+2}\rrbracket\rrbracket_{L,R}
%\\
%&+\sum_{t=2}^{n}(-1)^{t+1}
%(Q^{(t-1)}_{\lambda_2,\dots,\lambda_t+\lambda_{t+1},\dots,\lambda_{n+1}})_R
%(a_2,\dots,\llbracket a_{t}{}_{\lambda_t} a_{t+1}\rrbracket,\dots,a_{n+1},\llbracket a_1{}_{\lambda_1} a_{n+2}\rrbracket)
%\\
%&+(-1)^{n}
%((Q_{\lambda_2,\dots,\lambda_{n}})_R)^{(n)}
%(a_2,\dots,\llbracket a_{n+1}{}_{\lambda_{n+1}}\llbracket a_1{}_{\lambda_1} a_{n+2}\rrbracket\rrbracket_R)
%\\
%&
%-(Q_{\lambda_3,\dots,\lambda_{n+1}})_R(a_3,\dots,a_{n+1},\llbracket a_2{}_{\lambda_2}\llbracket a_1{}_{\lambda_1}a_{n+2}\rrbracket \rrbracket_R)
\\
&
=\sum_{t=2}^{n+1}(-1)^t\llbracket a_{t}{}_{\lambda_t}
\ldb a_2{}_{\lambda_2}\stackrel{t-1}{\check{\dots}}a_{n+1}{}_{\lambda_{n+1}}\llbracket a_1{}_{\lambda_1}a_{n+2}\rrbracket\rdb_R
\rrbracket_{(t)}
\\
&
+(-1)^{n+1}
\llbracket
\ldb a_2{}_{\lambda_2}\dots a_n{}_{\lambda_n}a_{n+1}\rdb_{\lambda_2+\dots+\lambda_{n+1}}
\llbracket a_1{}_{\lambda_1}a_{n+2}\rrbracket\rrbracket_{L,R}
\\
&+\sum_{t=2}^{n}(-1)^{t+1}
\ldb a_2{}_{\lambda_2}\dots a_{t-1}{}_{\lambda_{t-1}}\llbracket a_{t}{}_{\lambda_t} a_{t+1}\rrbracket_{\lambda_t+\lambda_{t+1}}a_{t+2}{}_{\lambda_{t+2}}\dots
\\
&\quad\quad\quad\quad\quad\quad\quad\quad\quad\quad\quad\quad
\dots a_{n+1}{}_{\lambda_{n+1}}\llbracket a_1{}_{\lambda_1} a_{n+2}\rrbracket\rdb_{L_{t-1},R}
\\
&+(-1)^{n}
\ldb a_2{}_{\lambda_2}\dots a_n{\lambda_n}\llbracket a_{n+1}{}_{\lambda_{n+1}}\llbracket a_1{}_{\lambda_1} a_{n+2}\rrbracket\rrbracket_R\rdb_{(2)}
\\
&-\ldb a_3 {}_{\lambda_3}\dots a_{n+1}{}_{\lambda_{n+1}}\llbracket a_2{}_{\lambda_2}\llbracket a_1{}_{\lambda_1}a_{n+2}\rrbracket \rrbracket_R\rdb_R
\,.
\end{split}
\end{equation}
Plugging \eqref{20240808:eq2},\eqref{20240808:eq3}, \eqref{20240808:eq4} and \eqref{20240808:eq5} into \eqref{20240808:eq1}, and performing some algebraic manipulations, we get
\begin{subequations}
\begin{align}
&\dd^2(Q)_{\lambda_1,\dots,\lambda_{n+1}}(a_1,\dots,a_{n+2})
\notag
\\
&=
\sum_{t=1}^{n+1}
\sum_{s=1}^{t-1}(-1)^{t+s+1}
\llbracket a_t {}_{\lambda_t}
\llbracket a_s {}_{\lambda_s}
\ldb a_1{}_{\lambda_1}\stackrel{s}{\check{\dots}}\stackrel{t}{\check{\dots}}
a_{n+1}{}_{\lambda_{n+1}}a_{n+2}\rdb\rrbracket_{(s)}
\rrbracket_{(t)}
\label{20240808:A1}
\\
&+\sum_{t=1}^{n+1}
\sum_{s=t+1}^{n+1}(-1)^{t+s}
\llbracket a_t {}_{\lambda_t}
\llbracket a_s {}_{\lambda_s}
\ldb a_1{}_{\lambda_1}\stackrel{t}{\check{\dots}}\stackrel{s}{\check{\dots}}
a_{n+1}{}_{\lambda_{n+1}}a_{n+2}\rdb\rrbracket_{(s-1)}
\rrbracket_{(t)}
\label{20240808:A2}
\\
&-
\sum_{t=1}^{n}(-1)^{n+t}
\llbracket a_t {}_{\lambda_t}
\llbracket
\ldb a_1{}_{\lambda_1}\stackrel{t}{\check{\dots}}a_n {}_{\lambda_n}a_{n+1}\rdb_{\lambda_1+\stackrel{t}{\check{\dots}}+\lambda_{n+1}}a_{n+2}\rrbracket_{L}
\rrbracket_{(t)}
\label{20240808:M1}
\\
&+\llbracket a_{n+1} {}_{\lambda_{n+1}}
\llbracket
\ldb a_1{}_{\lambda_1}\dots a_{n-1}{}_{\lambda_{n-1}}a_{n}\rdb_{\lambda_1+\dots+\lambda_n}a_{n+2}\rrbracket_{L}
\rrbracket_{(n+1)}
\label{20240808:N1}\\
\begin{split}\label{20240808:B1}
&+\sum_{t=1}^{n+1}
\sum_{s=1}^{t-2}(-1)^{t+s}
\llbracket a_t {}_{\lambda_t}
\ldb a_1{}_{\lambda_1}\dots a_{s-1}{}_{\lambda_{s-1}}\llbracket a_{s}{}_{\lambda_s}a_{s+1}\rrbracket_{\lambda_{s}+\lambda_{s+1}}a_{s+2}{}_{\lambda_{s+2}}\stackrel{t-1}{\check{\dots}}
\\
&\quad\quad
\dots a_{n+1}{}_{\lambda_{n+1}}a_{n+2}\rdb_L
\rrbracket_{(t)}
\end{split}
\\
\begin{split}\label{20240808:D1}
&-\sum_{t=1}^{n+1}
\llbracket a_t {}_{\lambda_t}
\ldb a_1{}_{\lambda_1}\dots a_{t-2}{}_{\lambda_{t-2}}
\llbracket a_{t-1}{}_{\lambda_{t-1}}a_{t+1}\rrbracket_{\lambda_{t-1}+\lambda_{t+1}}
a_{t+2}{}_{\lambda_{t+2}}\dots
\\
&\quad\quad
\dots a_{n+1}{}_{\lambda_{n+1}}a_{n+2}\rdb_L
\rrbracket_{(t)}
\end{split}
\\
\begin{split}\label{20240808:C1}
&+\sum_{t=1}^{n+1}\sum_{s=t+1}^{n+1}(-1)^{t+s+1}
\llbracket a_t {}_{\lambda_t}
\ldb a_1{}_{\lambda_1}\stackrel{t}{\check{\dots}} a_{s-1}{}_{\lambda_{s-1}}\llbracket a_{s}{}_{\lambda_s}a_{s+1}\rrbracket_{\lambda_{s}+\lambda_{s+1}}a_{s+2}{}_{\lambda_{s+1}}\dots \\
&\quad\quad
\dots a_{n+1}{}_{\lambda_{n+1}}a_{n+2}\rdb_L
\rrbracket_{(t)}
\end{split}
\\
&+
\llbracket a_1{}_{\lambda_1}
\ldb a_3{}_{\lambda_3}\dots a_{n+1}{}_{\lambda_{n+1}} \llbracket a_{2}{}_{\lambda_{2}}a_{n+2}\rrbracket\rdb_R\rrbracket_{(1)}
\label{20240808:G1}
\\
&+\sum_{t=2}^{n+1}(-1)^{t+1}
\llbracket a_t {}_{\lambda_t}
\ldb a_2{}_{\lambda_2}\stackrel{t-1}{\check{\dots}} a_{n+1}{}_{\lambda_{n+1}}\llbracket a_{1}{}_{\lambda_{1}}a_{n+2}\rrbracket\rdb_R
\rrbracket_{(t)}
\label{20240808:I1}
\\
&
+\sum_{t=1}^{n}(-1)^{n+t}
\llbracket
{\llbracket a_t {}_{\lambda_t}
\ldb a_1{}_{\lambda_1}\stackrel{t}{\check{\dots}}
a_n{}_{\lambda_n}a_{n+1}\rdb\rrbracket_{(t)}}_{\lambda_1+\dots+\lambda_{n+1}}a_{n+2}
\rrbracket_{L}
\label{20240808:M2}
\\
&+\llbracket
{\llbracket
\ldb a_1{}_{\lambda_1}\dots a_{n-1}{}_{\lambda_{n-1}}a_n\rdb_{\lambda_1+\dots+\lambda_n}a_{n+1}\rrbracket_{L}}_{\lambda_1+\dots+\lambda_{n+1}}a_{n+2}
\rrbracket_{L}
\label{20240808:N2}\\
\begin{split}\label{20240808:H1}
&+\sum_{t=1}^n(-1)^{n+t+1}
\llbracket
\ldb a_1{}_{\lambda_1}\dots a_{t-1}{}_{\lambda_{t-1}}\llbracket a_{t}{}_{\lambda_t}a_{t+1}\rrbracket_{\lambda_{t}+\lambda_{t+1}}a_{t+2}{}_{\lambda_{t+2}}\dots
\\
&\quad\quad
{\dots a_n{}_{\lambda_n}a_{n+1}\rdb_L}
_{\lambda_1+\dots+\lambda_{n+1}}a_{n+2}
\rrbracket_{L}
\end{split}
\\
&+(-1)^{n}
\llbracket
{\ldb a_2{}_{\lambda_2}\dots a_n{}_{\lambda_n}\llbracket a_{1}{}_{\lambda_{1}}a_{n+1}\rrbracket\rdb_R}
_{\lambda_1+\dots+\lambda_{n+1}}a_{n+2}
\rrbracket_{L}
\label{20240808:L1}
\\
\begin{split}\label{20240808:C2}
&
+\sum_{t=1}^{n+1}\sum_{s=1}^{t-1}(-1)^{t+s}
\llbracket a_s{}_{\lambda_s}
\ldb a_1{}_{\lambda_1}\stackrel{s}{\check{\dots}} a_{t-1}{}_{\lambda_{t-1}}
\llbracket a_{t}{}_{\lambda_t}a_{t+1}\rrbracket_{\lambda_{t}+\lambda_{t+1}}a_{t+2}{}_{\lambda_{t+2}}\dots
\\
&\quad\quad
\dots a_{n+1}{}_{\lambda_{n+1}}a_{n+2}\rdb_L
\rrbracket_{(s)}
\end{split}
\\
&+\sum_{t=1}^{n+1}
\llbracket
\llbracket a_t{}_{\lambda_t}a_{t+1}\rrbracket_{\lambda_t+\lambda_{t+1}}
\ldb a_1{}_{\lambda_1}\dots a_{t-1}{}_{\lambda_{t-1}}a_{t+2}{}_{\lambda_{t+2}}\dots a_{n+1}{}_{\lambda_{n+1}}a_{n+2}\rdb
\rrbracket_{L,(t)}
\label{20240808:A3}\\
\begin{split}\label{20240808:B2}
&+\sum_{t=1}^{n+1}\sum_{s=t+2}^{n+1}(-1)^{t+s+1}
\llbracket a_s{}_{\lambda_s}
\ldb a_1{}_{\lambda_1}\dots a_{t-1}{}_{\lambda_{t-1}}
\llbracket a_{t}{}_{\lambda_t}a_{t+1}\rrbracket_{\lambda_{t}+\lambda_{t+1}}
a_{t+2}{}_{\lambda_{t+2}}\stackrel{s-1}{\check{\dots}}
\\
&\quad\quad
\dots a_{n+1}{}_{\lambda_{n+1}}a_{n+2}\rdb_L
\rrbracket_{(s)}
\end{split}
\\
\begin{split}
&
+\sum_{t=1}^{n}(-1)^{n+t}
\llbracket
\ldb a_1{}_{\lambda_1}\dots a_{t-1}{}_{\lambda_{t-1}}
\llbracket a_{t}{}_{\lambda_t}a_{t+1}\rrbracket_{\lambda_{t}+\lambda_{t+1}}a_{t+2}{}_{\lambda_{t+2}}\dots
\\
&\quad\quad
{\dots a_n{}_{\lambda_n}a_{n+1}\rdb_L}
_{\lambda_1+\dots+\lambda_{n+1}}a_{n+2}\rrbracket_L
\label{20240808:H2}
\end{split}
\\
&
-
\llbracket
\ldb a_1{}_{\lambda_1}\dots a_{n-1}{}_{\lambda_{n-1}} a_{n}\rdb
_{\lambda_1+\dots+\lambda_{n}}\llbracket a_{n+1}{}_{\lambda_{n+1}}a_{n+2}\rrbracket \rrbracket_{L,L}
\label{20240808:N3}\\
\begin{split}
&
+\sum_{t=1}^{n+1}\sum_{s=1}^{t-2}(-1)^{t+s+1}
\ldb a_1{}_{\lambda_1}\dots a_{s-1}{}_{\lambda_{s-1}}
\llbracket a_s{}_{\lambda_s} a_{s+1}\rrbracket_{\lambda_{s}+\lambda_{s+1}}
a_{s+2}{}_{\lambda_{s+2}}\dots
\\
&\quad\quad
\dots a_{t-1}{}_{\lambda_{t-1}}
\llbracket a_t{}_{\lambda_t}a_{t+1}\rrbracket_{\lambda_{t}+\lambda_{t+1}}a_{t+2}{}_{\lambda_{t+2}}\ldots a_{n+1}{}_{\lambda_{n+1}} a_{n+2} \rdb_{L_s,L_{t-1}}
\label{20240808:E1}
\end{split}
\\
\begin{split}
&
+\sum_{t=1}^{n+1}
\ldb a_1{}_{\lambda_1}\dots a_{t-2}{}_{\lambda_{t-2}}
{\llbracket a_{t-1}{}_{\lambda_{t-1}}\llbracket a_t {}_{\lambda_t} a_{t+1}\rrbracket\rrbracket _L}_{\lambda_{t-1}+\lambda_t+\lambda_{t+1}}
a_{t+2}{}_{\lambda_{t+2}}\dots
\\
&\quad\quad
\dots a_{n+1}{}_{\lambda_{n+1}}a_{n+2}\rdb_L
\label{20240808:D2}
\end{split}
\\
\begin{split}
&
-\sum_{t=1}^{n+1}
\ldb a_1{}_{\lambda_1}\dots a_{t-1}{}_{\lambda_{t-1}}
{\llbracket \llbracket a_{t}{}_{\lambda_{t}} a_{t+1}\rrbracket {}_{\lambda_t+\lambda_{t+1}} a_{t+2}\rrbracket_L}_{\lambda_{t}+\lambda_{t+1}+\lambda_{t+2}}a_{t+3}{}_{\lambda_{t+3}} \dots
\\
&\quad\quad
\dots a_{n+1}{}_{\lambda_{n+1}}a_{n+2}\rdb_L
\label{20240808:D3}
\end{split}
\\
\begin{split}
&
+\sum_{t=1}^{n+1}\sum_{s=t+2}^{n+1}(-1)^{t+s}
\ldb a_1{}_{\lambda_1}\dots a_{t-1}{}_{\lambda_{t-1}}
\llbracket a_t{}_{\lambda_t} a_{t+1}\rrbracket_{\lambda_{t}+\lambda_{t+1}}
a_{t+2}{}_{\lambda_{t+2}}\dots
\\
&\quad\quad
\dots a_{s-1}{}_{\lambda_{s-1}}
\llbracket a_s{}_{\lambda_s}a_{s+1}\rrbracket_{\lambda_{s}+\lambda_{s+1}}
a_{s+2}{}_{\lambda_{s+2}}\ldots a_{n+1}{}_{\lambda_{n-1}}a_{n+2} \rdb_{L_{s-1},L_t}
\label{20240808:E2}
\end{split}
\\
&-
\ldb a_3{}_{\lambda_3}\dots a_{n+1}{}_{\lambda_{n+1}}\llbracket \llbracket a_1{}_{\lambda_1}a_2\rrbracket_{\lambda_1+\lambda_2}a_{n+2}\rrbracket_L\rdb_R
\label{20240808:G3}\\
\begin{split}
&
+\sum_{t=2}^{n}(-1)^{t}
\ldb a_2{}_{\lambda_2}\dots a_{t-1}{}_{\lambda_{t-1}}
\llbracket a_{t}{}_{\lambda_t} a_{t+1}\rrbracket_{\lambda_t+\lambda_{t+1}}a_{t+2}{}_{\lambda_{t+2}}\dots
\\
&\quad\quad
\dots a_{n+1}{}_{\lambda_{n+1}}\llbracket a_1{}_{\lambda_1} a_{n+2}\rrbracket\rdb_{L_{t-1},R}
\label{20240808:F1}
\end{split}
\\
&
+(-1)^{n+1}
\ldb a_2{}_{\lambda_2}\dots a_{n}{}_{\lambda_n}\llbracket a_1{}_{\lambda_1}\llbracket a_{n+1}{}_{\lambda_{n+1}} a_{n+2}\rrbracket\rrbracket_L\rdb_{(2)}
\label{20240808:L3}
\\
&+\sum_{t=2}^{n+1}(-1)^t\llbracket a_{t}{}_{\lambda_t}
\ldb a_2{}_{\lambda_2}\stackrel{t-1}{\check{\dots}} a_{n+1}{}_{\lambda_{n+1}}\llbracket a_1{}_{\lambda_1}a_{n+2}\rrbracket\rdb_R
\rrbracket_{(t)}
\label{20240808:I2}
\\
&
+(-1)^{n+1}
\llbracket
\ldb a_2{}_{\lambda}\dots a_{n}{}_{\lambda_n}a_{n+1}\rdb_{\lambda_2+\dots+\lambda_{n+1}}
\llbracket a_1{}_{\lambda_1}a_{n+2}\rrbracket\rrbracket _{L,R}
\label{20240808:M3}\tag{10.48A}
\\
\tag{10.48B}
\begin{split}
&
+\sum_{t=2}^{n}(-1)^{t+1}
\ldb a_2{}_{\lambda_2}\dots a_{t-1}{}_{\lambda_{t-1}}
\llbracket a_{t}{}_{\lambda_t} a_{t+1}\rrbracket_{\lambda_t+\lambda_{t+1}}
a_{t+2}{}_{\lambda_{t+2}}\dots
\\
&\quad\quad
\dots a_{n+1}{}_{\lambda_{n+1}}\llbracket a_1{}_{\lambda_1} a_{n+2}\rrbracket\rdb_{L_{t-1},R}
\label{20240808:F2}
\end{split}
\\
&+(-1)^{n}
\ldb a_2{}_{\lambda_2}\dots a_n{}_{\lambda_n}
\llbracket a_{n+1}{}_{\lambda_{n+1}}\llbracket a_1{}_{\lambda_1} a_{n+2}\rrbracket \rrbracket_R\rdb_{(2)}
\label{20240808:L2}\tag{10.48C}
\\
&
-
\ldb a_3{}_{\lambda_3}\dots a_{n+1}{}_{\lambda_{n+1}}
\llbracket a_2{}_{\lambda_2}\llbracket a_1{}_{\lambda_1}a_{n+2}\rrbracket \rrbracket_R\rdb_R
\label{20240808:G2}\tag{10.48D}\,.
\end{align}
\end{subequations}
To conclude the proof we show that the RHS above vanishes. Note that
$$
\eqref{20240808:I1}+\eqref{20240808:I2}=0
\,,
\qquad
\eqref{20240808:H1}+\eqref{20240808:H2}=0
\,,
\qquad
\eqref{20240808:F1}+\eqref{20240808:F2}=0
\,.
$$
Next, swapping $s$ and $t$
and changing the order of summation in \eqref{20240808:A1},
and using parts (i) and (ii) of Lemma \ref{lem:PVAcoh-pre7} we get
\begin{align*}
&\eqref{20240808:A2}+\eqref{20240808:A1}
+
\eqref{20240808:A3}
\\
&
=\sum_{t=1}^{n+1}
\sum_{s=t+1}^{n+1}(-1)^{t+s}
\left(
\llbracket a_t {}_{\lambda_t}
\llbracket a_s {}_{\lambda_s}
\ldb a_1{}_{\lambda_1}\stackrel{t}{\check{\dots}}\stackrel{s}{\check{\dots}}
a_{n+1}{}_{\lambda_{n+1}}a_{n+2}\rdb\rrbracket_{(s-1)}
\rrbracket_{(t)}
\right.
\\
&\left.
-\llbracket a_s {}_{\lambda_s}\llbracket a_t {}_{\lambda_t}
\ldb a_1{}_{\lambda_1}\stackrel{t}{\check{\dots}}\stackrel{s}{\check{\dots}}
a_{n+1}{}_{\lambda_{n+1}}a_{n+2}\rdb\rrbracket_{(t)}
\rrbracket_{(s)}
\right)\\
&+
\sum_{t=1}^{n+1}
\llbracket
\llbracket a_t{}_{\lambda_t}a_{t+1}\rrbracket_{\lambda_t+\lambda_{t+1}}
\ldb a_1{}_{\lambda_1}\dots a_{t-1}{}_{\lambda_{t-1}}a_{t+2}{}_{\lambda_{t+2}}\dots a_{n+1}{}_{\lambda_{n+1}}a_{n+2}\rdb\rrbracket_{L,(t)}
\\
&
=-\sum_{t=1}^{n+1}
\Big(
\llbracket a_t {}_{\lambda_t}
\llbracket a_{t+1} {}_{\lambda_{t+1}}
C\rrbracket_{(t)}\rrbracket_{(t)}
-\llbracket a_{t+1} {}_{\lambda_{t+1}}\llbracket a_t {}_{\lambda_t}C\rrbracket_{(t)}
\rrbracket_{(t+1)}
\\
&-\llbracket\llbracket a_t{}_{\lambda_t}a_{t+1}\rrbracket_{\lambda_t+\lambda_{t+1}}C\rrbracket_{L,(t)} 
\Big)=0
\,,
\end{align*}
where $C=\ldb a_1{}_{\lambda_1}\dots a_{t-1}{}_{\lambda_{t-1}}a_{t+2}{}_{\lambda_{t+2}}\dots a_{n+1}{}_{\lambda_{n+1}}a_{n+2}\rdb$.
By swapping $s$ and $t$ in both equations \eqref{20240808:C2} and \eqref{20240808:B2} and changing the order of summation, it follows that we get the cancellations$$
\eqref{20240808:C1}+\eqref{20240808:C2}=0
\,,
\qquad
\eqref{20240808:B1}+\eqref{20240808:B2}=0
\,.
$$
Furthermore, we use \eqref{badnotation2}, the second identity in \eqref{20230811:eq1}
(with $s=t$ and $i=n+1-t$) and sesquilinearity \eqref{eq:sesqui} to get 
\begin{align*}
&\eqref{20240808:D1}
=-\sum_{t=1}^{n+1}
\ldb a_1{}_{\lambda_1}\dots a_{t-2}{}_{\lambda_{t-2}}
\llbracket a_{t-1}{}_{\lambda_{t-1}}a_{t+1}\rrbracket'_{\lambda_{t-1}+\lambda_t+\lambda_{t+1}+x}a_{t+2}{}_{\lambda_{t+2}}\dots
\\
&\quad\quad\quad\quad\quad\quad
\dots a_{n+1}{}_{\lambda_{n+1}}a_{n+2}
\rdb
\otimes_{n+1-t}\left(|_{x=\partial}\llbracket a_t {}_{\lambda_t}\llbracket a_{t-1}{}_{\lambda_{t-1}}a_{t+1}\rrbracket''\rrbracket\right)
\\
&
=-\sum_{t=1}^{n+1}
\ldb a_1{}_{\lambda_1}\dots a_{t-2}{}_{\lambda_{t-2}}
{\llbracket a_t {}_{\lambda_t}\llbracket a_{t-1}{}_{\lambda_{t-1}}a_{t+1}\rrbracket\rrbracket_R}_{\lambda_{t-1}+\lambda_t+\lambda_{t+1}} a_{t+2}{}_{\lambda_{t+2}}\dots
\\
&\quad\quad\quad\quad\quad\quad
\dots a_{n+1}{}_{\lambda_{n+1}}a_{n+2}\rdb_L
\,.
\end{align*}
Similarly, we get the identities
\begin{align*}
&\eqref{20240808:G1}
=
\ldb a_3{}_{\lambda_3}\dots a_{n+1}{}_{\lambda_{n+1}}\llbracket a_1{}_{\lambda_1}\llbracket a_2{}_{\lambda_2} a_{n+2}\rrbracket \rrbracket_L\rdb_R\,,
\\
&\eqref{20240808:L1}
=(-1)^n
\ldb a_2{}_{\lambda_2}\dots a_n{}_{\lambda_n}
\llbracket\llbracket a_1{}_{\lambda_{1}}a_{n+1}\rrbracket _{\lambda_1+\lambda_{n+1}}a_{n+2}\rrbracket_L\rdb_{(2)}
\,.
\end{align*}
%
%\pecetta{Maxime: In fact \eqref{20240808:G1} is already written in this simplified form. We should either remove the equation here, or put back the "non-simplified" form as \eqref{20240808:G1}.}
%
Hence, using the above computations and simple algebraic manipulations we get
\begin{align*}
&\eqref{20240808:D2}+\eqref{20240808:D1}+\eqref{20240808:D3}
\\
&=\sum_{t=1}^{n+1}
\ldb a_1{}_{\lambda_1}\dots a_{t-2}{}_{\lambda_{t-2}}
{J_t}_{\lambda_{t-1}+\lambda_t+\lambda_{t+1}} a_{t+2}{}_{\lambda_{t+2}}\dots
a_{n+1}{}_{\lambda_{n+1}}a_{n+2}\rdb_L
\,,
\\
&\eqref{20240808:G1}+\eqref{20240808:G2}+\eqref{20240808:G3}
=\ldb a_3{}_{\lambda_3}\dots a_{n+1}{}_{\lambda_{n+1}}J\rdb_R=0\,,
\\
&
\eqref{20240808:L3}+\eqref{20240808:L2}+\eqref{20240808:L1}
=
(-1)^{n+1}
\ldb a_2{}_{\lambda_2}\dots a_{n}{}_{\lambda_n}\tilde J\rdb_{(2)}\,,
\end{align*}
where
\begin{align*}
&
J_t:=\llbracket {a_{t-1}}_{\lambda_{t-1}}\llbracket {a_{t}}_{\lambda_t} a_{t+1}\rrbracket\rrbracket_L
-\llbracket {a_t}_{\lambda_t}\llbracket {a_{t-1}}_{\lambda_{t-1}} a_{t+1}\rrbracket\rrbracket_R
-\llbracket \llbracket {a_{t-1}}_{\lambda_{t-1}} a_{t}\rrbracket_{\lambda_{t-1}+\lambda_{t}} a_{t+1}\rrbracket_L\,,
\\
&
J:=\llbracket {a_{1}}_{\lambda_{1}}\llbracket {a_{2}}_{\lambda_2} a_{n+2}\rrbracket\rrbracket_L
-\llbracket {a_2}_{\lambda_2}\llbracket {a_{1}}_{\lambda_{1}} a_{n+2}\rrbracket\rrbracket_R
-\llbracket \llbracket {a_{1}}_{\lambda_{1}} a_{2}\rrbracket_{\lambda_{1}+\lambda_{2}} a_{n+2}\rrbracket_L\,,
\\
&
\tilde J:=\llbracket a_1{}_{\lambda_1}\llbracket a_{n+1}{}_{\lambda_{n+1}} a_{n+2}\rrbracket\rrbracket_L
-\llbracket a_{n+1}{}_{\lambda_{n+1}}\llbracket a_{1}{}_{\lambda_{1}} a_{n+2}\rrbracket\rrbracket_R
-\llbracket \llbracket a_1{}_{\lambda_1}a_{n+1}\rrbracket_{\lambda_1+\lambda_{n+1}}a_{n+2}\rrbracket_L
\end{align*}
The cancellations occur by linearity since $J_t=J=\tilde J=0$ by the Jacobi identity \eqref{eq:jacobi2}.
Next,  we note that,
a similar computations as in the proof of Theorem \ref{thm:PVAcoh}(c) shows that
$$
\eqref{20240808:E1}+\eqref{20240808:E2}
=0\,.
$$
Using Lemma \ref{20240810:lem1}(i) and (ii) (with $a=a_t$,
$B=\ldb a_1{}_{\lambda_1}\stackrel{t}{\check{\dots}}a_n {}_{\lambda_n}a_{n+1}\rdb$ and $c=a_{n+2}$) we have
\begin{align*}
&\eqref{20240808:M1}+\eqref{20240808:M2}
=(-1)^{n}
\Big(
\llbracket a_1 {}_{\lambda_1}
\llbracket \ldb a_2{}_{\lambda}\dots a_{n}{}_{\lambda_n}a_{n+1}\rdb_{\lambda_2+\dots+\lambda_{n+1}}a_{n+2}\rrbracket_{L}
\rrbracket_{(1)}
\\
&
-\llbracket
{\llbracket a_1 {}_{\lambda_1}
\ldb a_2{}_{\lambda}\dots a_{n}{}_{\lambda_n}a_{n+1}\rdb\rrbracket_{(1)}}_{\lambda_1+\dots+\lambda_{n+1}}a_{n+2}
\rrbracket_{L}\Big)
\\
&=
(-1)^{n}
\llbracket \ldb a_2{}_{\lambda}\dots a_{n}{}_{\lambda_n}a_{n+1}\rdb_{\lambda_2+\dots\lambda_{n+1}}
\llbracket a_1{}_{\lambda_1}a_{n+2}\rrbracket\rrbracket_{L,R}
=-\eqref{20240808:M3}\,,
\end{align*}
where in the second equality we used \eqref{20240810:eq3}. At last, using \eqref{20240810:eq4}, we also have the cancellation
$$
\eqref{20240808:N1}+\eqref{20240808:N2}+\eqref{20240808:N3}=0
\,,
$$
thus showing that $\dd^2(Q)=0$ and concluding the proof.
\end{proof}
By Theorem \ref{Thm:g-dPVcoh1} we have a complex $(C(\mc V),\dd)$.
\begin{definition}\label{def:var-dPVA}
Let $\mc V$ be a dPVA. The complex $(C(\mc V),\dd)$ is called the
\emph{variational dPVA complex} of $\mc V$. The cohomology\glslink{dPVHA}{}
$$
\dPVH(\mc V)=\coH(C(\mc V),\dd)=\bigoplus_{n\in\mb Z_{\geq0}} \dPVH^n(\mc V)
\,,
$$
where $\dPVH^n(\mc V)=\ker(\dd|_{C^n(\mc V)})/\dd(C^{n-1}(\mc V))$, 
is called the \emph{variational dPVA cohomology} of $\mc V$.
\end{definition}

\subsection{Cohomology for double Lie conformal algebras}

A \emph{double Lie conformal algebra} (dLCA) is a $\kk[\partial]$-module $\mc V$ equipped with a linear map
\[
\dsq{-_\lambda-}: \mc V\otimes \mc V \to (\mc V \otimes \mc V)[\lambda]
\]
satisfying sesquilinearity \eqref{eq:sesqui}, cyclic skewsymmetry \eqref{eq:skew2} and
the Jacobi identity~\eqref{eq:jacobi2}.
Put $C_{\operatorname{dLCA}}^0(\mc V)=\mc V$,
and for $n\geq 1$ let $C_{\operatorname{dLCA}}^n(\mc V)$ be the vector space spanned by linear maps
$$\ldb-{}_{\lambda_1}-\dots-{}_{\lambda_{n-1}}-\rdb:\mc V^{\otimes n}\to\mc V^{\otimes n}[\lambda_1,\ldots,\lambda_{n-1}]
$$
satisfying sesquilinearity \eqref{20140702:eq4}--\eqref{20140702:eq5} and skewsymmetry \eqref{eq:nfold-skew}.
Furthermore, we set
$C_{\operatorname{dLCA}}(\mc V)=\bigoplus_{n\geq 0} C_{\operatorname{dLCA}}^n(\mc V)$.
If one forgets about the Leibniz rules of $n$-fold $\lambda$-brackets in the proof of Theorem \ref{Thm:g-dPcoh1},
one readily sees that the operation
\[
 \dd : C_{\operatorname{dLCA}}(\mc V) \to C_{\operatorname{dLCA}}(\mc V)
\]
given by \eqref{eq:dP0} and \eqref{eq:dP-1} is well-defined, and it squares to zero.
Hence we can define the \emph{dLCA cohomology} of $\mc V$
as $\operatorname{dL_cH}(\mc V)= \coH(C_{\operatorname{dLCA}}(\mc V), \dd)$.

%%%
\section{Explicit description of the first few cohomology spaces}\label{sec:explicit}

\subsection{\texorpdfstring{$n=0$}{n=0}}\label{explicit:0}
By \eqref{eq:dP0} we have that
$$
\dPVH^0(\mc V)=\{\tint a\in\mc V_{\sharp}\mid \dsq{\tint a,b}=0\,,\text{ for every }b\in\mc V\}
=\Cas(\mc V)\glslink{casimir}{}\,,
$$
is the space of Casimir elements of $\mc V$, namely elements of the Lie algebra
$\mc V_\sharp$ which act trivially on $\mc V$ with respect to the Lie algebra action defined by \eqref{20140707:eq3b-lie}.

\subsection{\texorpdfstring{$n=1$}{n=1}}\label{explicit:1}
Recall that $C^1(\mc V)=\Vect(\mc V)^\partial$ consists of derivations of $\mc V$ commuting with $\partial$.
We call an element $D\in C^1(\mc V)$ a \emph{dPVA-derivation}
of $\mc V$ if
$$
D(\dsq{a_\lambda b})=\dsq{D(a)_\lambda b}+\dsq{a_\lambda D(b)}
\,,
$$
for every $a,b\in\mc V$ (see \eqref{20240729:eq1}),
and we call it an \emph{inner dPVA-derivation} of $\mc V$ if
$D=\llbracket\tint f,-\rrbracket$ for some $f\in\mc V$. It follows from equations
\eqref{eq:dP0} and \eqref{20240729:eq1} that
$$
\dPVH^1(\mc V)=
\frac{\{\text{dPVA-derivations of $\mc V$}\}}{\{\text{inner dPVA-derivations of $\mc V$}\}}
\,.
$$

\subsection{\texorpdfstring{$n=2$}{n=2}}\label{explicit:2}
We relate the cohomology space $\dPVH^2(\mc V)$ to equivalence classes of first-order deformations of the dPVA $\mc V$ (that preserve the differential algebra structure of $\mc V$) in the spirit of the paper \cite{NR}.

Let $\mc V$ be a dPVA with $2$-fold $\lambda$-bracket $\dsq{-_\lambda-}$. A \emph{formal deformation} of $\mc V$ is a $2$-fold $\lambda$-bracket on $\mc V[[\epsilon]]$
\begin{equation}\label{eq:for-def}
\ldb-_\lambda-\rdb_\epsilon=\dsq{-_\lambda-}+\sum_{n\geq1}\epsilon^n Q_n\in
C^2(\mc V)[[\epsilon]]
\,,
\end{equation}
satisfying the Jacobi identity \eqref{eq:jacobi2}.
Recall the map $\ldb-_\lambda-_\mu-\rdb_\epsilon:\mc V^{\otimes3}\to\mc V^{\otimes 3}[\lambda,\mu]$
defined in \eqref{eq:triple1} using the RHS of the Jacobi identity
\eqref{eq:jacobi2}. It is shown in \cite[Lem.~3.4]{DSKV} that
$\ldb-_\lambda-_\mu-\rdb_\epsilon$ is a $3$-fold $\lambda$-bracket on $\mc V$.
Using the Jacobi identity and \eqref{eq:for-def} we get ($a,b,c\in\mc V$)
\begin{equation}\label{20250715:eq1}
\begin{split}
&\ldb a_\lambda b_\mu c\rdb_\epsilon=\epsilon\Big(\dsq{a_\lambda\ldb b_\mu c\rdb}_L-\dsq{b_\mu\ldb a_\lambda c\rdb}_R
-\dsq{\ldb a_{\lambda}b\rdb_{\lambda+\mu}c}_L
\\
&-\ldb \dsq{a_\lambda b}_{\lambda+\mu}c\rdb_L+\ldb a_\lambda\dsq{b_\mu c}\rdb_L
-\ldb b_\mu\dsq{a_\lambda b} \rdb_R\Big)
+\epsilon^2 \mc V^{\otimes 3}[[\epsilon]]
\,,
\end{split}
\end{equation}
where we are denoting $(Q_1)_\lambda(a,b)=\ldb a_\lambda b\rdb$. A $2$-fold $\lambda$-bracket
$\ldb -_\lambda -\rdb\in C^2(\mc V)$ is called a \emph{first-order deformation}
of $\mc V$ if the coefficient of $\epsilon$ in the RHS of \eqref{20250715:eq1} vanishes.
On the other hand, by letting $Q=\ldb-_\lambda-\rdb\in C^2(\mc V)$ in 
\eqref{eq:dP-1} we get
\begin{equation}\label{20250715:eq1-b}
\begin{split}
\dd(Q)_{\lambda,\mu}(a,b,c)&=\dsq{a_\lambda\ldb b_\mu c\rdb}_L-\dsq{b_\mu\ldb a_\lambda c\rdb}_R
-\dsq{\ldb a_{\lambda}b\rdb_{\lambda+\mu}c}_L
\\
&-\ldb \dsq{a_\lambda b}_{\lambda+\mu}c\rdb_L+\ldb a_\lambda\dsq{b_\mu c}\rdb_L
-\ldb b_\mu\dsq{a_\lambda b}\rdb_R
\,.
\end{split}
\end{equation}
Comparing equations \eqref{20250715:eq1} and \eqref{20250715:eq1-b} we have that
$\ker(\dd|_{C^2(\mc V)})$ consists of the first-order deformations of the dPVA $\mc V$.

Next, let $(\mc V[[\epsilon]],\ldb-_\lambda-\rdb_\epsilon^1)$ and $(\mc V[[\epsilon]],\ldb-_\lambda-\rdb_{\epsilon}^2)$ be two formal
deformations of $\mc V$. We say that they are \emph{equivalent} if there exists a
$\kk[[\epsilon]]$-linear map
$\phi:\mc V[[\epsilon]]\rightarrow\mc V[[\epsilon]]$ such that ($a,b\in\mc V$):
\begin{enumerate}[(a)]
\item
$\phi(a)=a+\epsilon \mc V[[\epsilon]]$;
\item
$\phi(ab)=\phi(a)\phi(b)$;
\item
$\phi(\partial a)=\partial \phi(a)$;
\item
$(\phi\otimes\phi)(\ldb a_\lambda b\rdb_{\epsilon}^1)=\ldb \phi(a)_\lambda \phi(b)\rdb_\epsilon^2$.
\end{enumerate}
Condition (a) implies that $\phi$ is invertible, then (b) and (c) mean that it is a differential algebra automorphism, and (d) means that it is a dPVA isomorphism.
A formal deformation is called \emph{trivial} if it is equivalent to $(\mc V,\dsq{-_\lambda-})$.

Let $\ldb-_\lambda-\rdb_\epsilon^i=\dsq{-_\lambda-}+\epsilon\ldb-_\lambda-\rdb^i+\epsilon^2(\dots)$, $i=1,2$, be two equivalent formal
deformations of $\mc V$ and let us write
$\phi(a)=a+\epsilon D(a)+\epsilon^2(\dots)$, where $D:\mc V\to\mc V$ is a linear map.
Comparing the coefficient of $\epsilon$ in both sides of conditions (b) and (c) we get that $D\in C^1(\mc V)$. Furthermore, by comparing the coefficient of $\epsilon$ in both sides of condition
(d) we get the following identity ($a,b\in\mc V$)
\begin{equation}\label{20250716:eq1}
\ldb a_\lambda b\rdb^2-\ldb a_\lambda b\rdb^1
=D(\dsq{a_\lambda b})-\dsq{D(a)_\lambda b}-\dsq{a_\lambda D(b)}
\,.
\end{equation}
Comparing equations \eqref{20240729:eq1} and \eqref{20250716:eq1} we have that two
equivalent first-order deformations $\ldb-_\lambda-\rdb^i\in C^2(\mc V)$, $i=1,2$, differ by
an element in $\dd(C^1(\mc V))$.
Finally, note that an element in $\dd(C^1(\mc V))$ is a trivial first-order deformation. Indeed
for $D\in C^1(\mc V)$, let us define
$$
\ldb a_\lambda b\rdb_\epsilon^D=e^{\epsilon D}\left(\dsq{e^{-\epsilon D}(a)_\lambda
e^{-\epsilon D}(b)}\right)=\dsq{a_\lambda b}+\epsilon\, \dd(D)_\lambda(a,b)
+\epsilon^2(\dots)
\,,
$$
where we used equation \eqref{20240729:eq1}. It is immediate to verify that $\phi=e^{\epsilon D}:\mc V,\to\mc V[[\epsilon]]$ satisfies conditions (a)-(d) above showing that
$(\mc V,\dsq{-_\lambda-})$ and $(\mc V[[\epsilon]],\ldb -_\lambda -\rdb_\epsilon^D)$ are equivalent.

In conclusion, the space
$$
\dPVH^2(\mc V)
=\frac{\{\text{first-order deformations of }\mc V\}}{\{\text{trivial first-order deformations of }\mc V\}}
$$
parametrizes the equivalence classes of first-order deformations of the dPVA $\mc V$.

%%%
\section[Relation with equation (4.3)]{Relation with equation (\ref{Eq:dP-gen})}
Let $\mc V$ be a dPVA with $2$-fold $\lambda$ bracket $H=\llbracket -_{\lambda}-\rrbracket$, and let $(C(\mc V),\dd)$ be the variational dPVA complex introduced 
in Section \ref{sec:var-dPVA}. Given an $n$-fold $\lambda$-bracket
$Q\in C^{n}(\mc V)$, 
$n\geq1$, we can write the expression for the $(n+1)$-fold $\lambda$-bracket 
$\dd(Q)\in C^{n+1}(\mc V)$ given by \eqref{eq:dP-1} in a more covariant form similar to equation
\eqref{Eq:dP-gen} which describes the differential of the completed dPA complex.
\begin{proposition} \label{Thm:g-dPVcoh1-compare}
Let $Q\in C^{n}(\mc V)$. Then, $\dd(Q)\in C^{n+1}(\mc V)$ given by \eqref{eq:dP-1} can be equivalently defined
by the formula
\begin{equation}
\begin{split}\label{eq:dP}
&\dd(Q)
=\sum_{s=0}^n (-1)^{ns}
|_{\lambda_{n+1}=\lambda_{n+1}^\dagger}
\sigma^s \circ  
\left(Q_{\lambda_{\sigma^s(1)},\dots,\lambda_{\sigma^s(n-1)}}\otimes \Id_{\mc V}\right)
\\
&\quad\quad\quad\quad\quad\quad\quad\quad\quad\quad\quad\quad\quad\quad\quad
\circ \left(\Id_{\mc V}^{\otimes(n-1)}\otimes H_{\lambda_{\sigma^s(n)}}\right)
\circ \sigma^{-s} \\
&+
%\frac{n}{n+1}
\sum_{s=0}^{n} (-1)^{n(s+1)}
|_{\lambda_{n+1}=\lambda_{n+1}^\dagger}
\sigma^s \circ  
\left(H_{\lambda_{\sigma^s(1)}}\otimes \Id_{\mc V}^{\otimes(n-1)}\right)
\\
&\quad\quad\quad\quad\quad\quad\quad\quad\quad\quad\quad\quad\quad\quad\quad
\circ \left(\Id_{\mc V}\otimes Q_{\lambda_{\sigma^s(2)},\dots,\lambda_{\sigma^s(n)}}\right)
\circ \sigma^{-s} \,,
\end{split}
\end{equation}
where we are using \eqref{eq:dagger}.
\end{proposition}
We note that, in the RHS of \eqref{eq:dP}, the cyclic permutation $\sigma$ is acting on $\mc V^{\otimes (n+1)}$, hence it denotes the cycle $(1\dots n+1$).
\begin{proof}[Proof of Proposition \ref{Thm:g-dPVcoh1-compare}]
We denote the action of $Q\in C^n(\mc V)$ on $a_1\otimes a_{n}\in\mc V^{\otimes n}$ by $Q_{\lambda_1,\dots,\lambda_{n-1}}(a_1,\dots,a_n)
=\ldb a_1{}_{\lambda_1}\dots a_{n-1}{}_{\lambda_{n-1}}a_n\rdb$.
Let us also denote the RHS of \eqref{eq:dP} applied to
$a_1\otimes \dots\otimes a_{n+1}\in \mc V^{\otimes (n+1)}$ as
\begin{equation}\label{20250718:eq1}
\dd(Q)_{\lambda_1,\dots,\lambda_n}(a_1,\dots,a_{n+1})
=\sum_{s=0}^n(-1)^{ns}\left(A_s+(-1)^nB_n\right)\,,
\end{equation}
where we set
\begin{equation}
\begin{split}\label{eq:As}
&A_s=|_{\lambda_{n+1}=\lambda_{n+1}^\dagger}
\sigma^s \circ  
\left(Q_{\lambda_{\sigma^s(1)},\dots,\lambda_{\sigma^s(n-1)}}\otimes \Id_{\mc V}\right)
\\
&\quad\quad\quad\quad
\circ \left(\Id_{\mc V}^{\otimes(n-1)}\otimes H(\lambda_{\sigma^s(n)})\right)
\circ \sigma^{-s}(a_1\otimes\dots\otimes a_{n+1})
\\
&=|_{\lambda_{n+1}=\lambda_{n+1}^\dagger}
\sigma^s
\left(
\ldb a_{s+1}{}_{\lambda_{s+1}}\dots a_{n+1}{}_{\lambda_{n+1}} a_1{}_{\lambda_1}\dots\llbracket a_{s-1}{}_{\lambda_{s-1}}a_s\rrbracket\rdb_L\right)
\,,
\end{split}
\end{equation}
using \eqref{badnotation0} in the second equality,
and
\begin{equation}
\begin{split}\label{eq:Bs}
&B_s=|_{\lambda_{n+1}=\lambda_{n+1}^\dagger}
\sigma^s \circ  
\left(H(\lambda_{\sigma^s(1)})\otimes \Id_{\mc V}^{\otimes(n-1)}\right)
\\
&\quad\quad\quad\quad\quad
\circ \left(\Id_{\mc V}\otimes Q_{\lambda_{\sigma^s(2)},\dots,\lambda_{\sigma^s(n)}}\right)
\circ \sigma^{-s} (a_1\otimes\dots\otimes a_{n+1})
\\
&=|_{\lambda_{n+1}=\lambda_{n+1}^\dagger}
\sigma^s
\left(
\llbracket a_{s+1}{}_{\lambda_{s+1}}
\ldb a_{s+2}{}_{\lambda_{s+2}}\dots a_{n+1}{}_{\lambda_{n+1}} a_1{}_{\lambda_1}\dots a_{s-1}{}_{\lambda_{s-1}}a_s\rdb 
\rrbracket_{(1)}\right)
\,,
\end{split}
\end{equation}
using \eqref{badnotation4} in the second equality.
In equations \eqref{eq:As} and \eqref{eq:Bs} the index $k$ in $a_k$ is understood modulo $n+1$.
Setting $s=0$ in \eqref{eq:As} we have
\begin{equation}\label{eq:A0}
A_0
=\ldb a_{1}{}_{\lambda_{1}}\dots a_{n-1}{}_{\lambda_{n-1}}
\llbracket a_{n}{}_{\lambda_{n}}a_{n+1}\rrbracket\rdb_L
\,,
\end{equation}
while setting $s=1$ we have, using 
\eqref{20240812:eq1},
\begin{align*}
&A_1=
|_{\lambda_{n+1}=\lambda_{n+1}^\dagger}
\sigma
\left(
\ldb a_{2}{}_{\lambda_{2}}\dots a_{n}{}_{\lambda_{n}} \llbracket a_{n+1}{}_{\lambda_{n+1}}a_1\rrbracket\rdb_L\right)
\\
&=|_{\lambda_{n+1}=\lambda_{n+1}^\dagger}
\ldb a_{2}{}_{\lambda_{2}}\dots a_{n}{}_{\lambda_{n}}
\llbracket a_{n+1}{}_{\lambda_{n+1}}a_1\rrbracket^\sigma\rdb_R
\,.
\end{align*}
Using skewsymmetry \eqref{eq:skew2} and sesquilinearity \eqref{20140702:eq5} we rewrite the above equation as
\begin{equation}\label{eq:A1}
A_1
=-\ldb a_{2}{}_{\lambda_{2}}\dots a_{n}{}_{\lambda_{n}}
\llbracket a_{1}{}_{\lambda_{1}}a_{n+1}\rrbracket \rdb_R
\,.
\end{equation}
Similarly, by repeatedly applying \eqref{20240812:skew} and using sesquilinearity \eqref{20140702:eq4}
we have
\begin{equation}\label{eq:Asexplicit}
A_s
=(-1)^{(s+1)(n+1)}
\ldb a_1{}_{\lambda_1}\dots a_{s-2}{}_{\lambda_{s-2}}
\llbracket a_{s-1}{}_{\lambda_{s-1}}a_{s}\rrbracket_{\lambda_{s-1}+\lambda_{s}}a_{s+1}\dots
{}_{\lambda_n}a_{n+1}\rdb_{L}
\,,
\end{equation}
for $s=2,\dots,n$.
Moreover, using Lemma \ref{20140606:lem} and skewsymmetry \eqref{eq:nfold-skew} we have
\begin{equation}\label{eq:Bsexplicit}
B_s=(-1)^{s(n+1)}\llbracket a_{s+1} {}_{\lambda_{s+1}} \ldb a_1{}_{\lambda_1}\stackrel{s+1}{\check{\dots}}
a_n{}_{\lambda_n}a_{n+1}\rdb\rrbracket_{(s+1)}
\,,
\end{equation}
for $s=0,\dots,n-1$, while using 
skewsymmetry \eqref{eq:skew2} we can rewrite
\begin{equation}\label{eq:Bn}
B_n=
-\llbracket \ldb a_{1}{}_{\lambda_{1}}\dots a_{n-1}{}_{\lambda_{n-1}} a_n\rdb _{\lambda_1+\dots+\lambda_n}
a_{n+1}\rrbracket_L
\,.
\end{equation}
Plugging equations \eqref{eq:A0}, \eqref{eq:A1},
\eqref{eq:Asexplicit}, \eqref{eq:Bsexplicit} and \eqref{eq:Bn} in the RHS of \eqref{20250718:eq1} by a straightforward computation we get the RHS of \eqref{eq:dP-1}.
This concludes the proof.
\end{proof}
\begin{remark}\label{rem:chemla}
Using equations \eqref{20250718:eq1}, \eqref{eq:As} and \eqref{eq:Bs} (where we shift $s$ by $s-1$) we rewrite \eqref{eq:dP-1} as
\begin{equation}\label{20251017:eq1}
\begin{split}
&\dd(Q)_{\lambda_1,\dots,\lambda_n}(a_1,\dots,a_{n+1})
\\
&=\sum_{s=0}^n(-1)^{ns}|_{\lambda_{n+1}=\lambda_{n+1}^\dagger}
\sigma^s
\ldb a_{s+1}{}_{\lambda_{s+1}}\dots a_{n+1}{}_{\lambda_{n+1}} a_1{}_{\lambda_1}\dots\llbracket a_{s-1}{}_{\lambda_{s-1}}a_s\rrbracket\rdb_L\\
&+\sum_{s=0}^n(-1)^{ns}
|_{\lambda_{n+1}=\lambda_{n+1}^\dagger}
\sigma^{s-1}\llbracket a_{s}{}_{\lambda_{s}}
\ldb a_{s+1}{}_{\lambda_{s+1}}\dots
 a_{n+1}{}_{\lambda_{n+1}} a_1{}_{\lambda_1}\dots a_{s-2}{}_{\lambda_{s-2}}a_{s-1}\rdb
\rrbracket_{L}
\end{split}
\end{equation}
which is the dPVA analogue of Chemla's formula \eqref{Eq:Chemla}.
\end{remark}

\begin{remark}
After the first version of this memoir appeared on arXiv, the authors of \cite{ZT} have shared with us their work in which they define a dPVA cohomology.
Their construction relies on constructing a graded Lie bracket $[-,-]_{\operatorname{DPV}}$ on the space $C(\mc V)$ of $n$-fold $\lambda$-brackets.
Then, $\mc P:=\dsq{-,-}\in C^2(\mc V)$ defines a dPVA structure if and only if
$[\mc P,\mc P]_{\operatorname{DPV}}=0$. Moreover, if
$[\mc P,\mc P]_{\operatorname{DPV}}=0$, then $\dd_{\operatorname{PVA}}=[\mc P,-]_{\operatorname{DPV}}$ is a square-zero differential on $C^n(\mc V)$. With their notation, the action of $\dd_{\operatorname{PVA}}$ on $Q\in C^n(\mc V)$ is given by the formula
\[
\dd_{\operatorname{PVA}}(Q)= [\mc P,Q]_{\operatorname{DPV}} :=\mc P\diamond Q - (-1)^{n-1} Q\diamond \mc P\,.
\]
Evaluating the RHS above on $a_1\otimes \dots\otimes a_{n+1}\in \mc V^{\otimes (n+1)}$ according to
\cite[Eq.~(3.4)]{ZT}, we get $\dd_{\operatorname{PVA}}(Q)=(-1)^n\dd(Q)$, where $\dd(Q)$ is given by \eqref{20251017:eq1}.
Thus, the two approaches lead to the same cohomology theory.
As an interesting result, one obtains from \cite[Cor.~5.3]{ZT} that there is no obstruction to extend deformations of a dPVA $2$-fold $\lambda$-bracket if $\dPVH^3(\mc V)=0$.
\end{remark}

%%%
\section[Relation between reduced \& variational \MakeLowercase{d}PVA cohomology]{Relation between reduced and variational double Poisson vertex algebra cohomologies}

In this section we construct a homomorphism of complexes between the
reduced dPVA complex constructed in \ref{sec:red-dPVA}
and the variational dPVA complex
constructed in Section \ref{sec:var-dPVA}.
In Chapter \ref{sec:PVAdiff}, we will show that the two complexes are isomorphic for algebras of (noncommutative) differential polynomials.

%%%
\subsection{The projection operator}
Let $\mc V$ be a dPVA with $2$-fold $\lambda$-bracket $\ldb-_\lambda-\rdb$. Recall the definition of 
the space of basic cochains $\widetilde{\Gamma}(\mc V)$ given in Section \ref{sec:omega1}, of the quotient space $\Gamma(\mc V)$ constructed in Section \ref{sec:red-dPVA} and of the space of $n$-fold $\lambda$-brackets $C(\mc V)$ constructed in Section \ref{sec:compl-dPVA-1}.
In particular, we have that
$\widetilde{\Gamma}^0(\mc V)=\mc V$, $\Gamma^0(\mc V)=C^0(\mc V)=\mc V_{\sharp}$
and $C^1(\mc V)=\Vect(\mc V)^\partial$. In \eqref{eq:P1} we constructed a linear map $\widetilde{\proj}_1:\widetilde{\Gamma}^1(\mc V)\rightarrow C^1(\mc V)$ which induces a linear map $\proj_1:\Gamma^1(\mc V)\rightarrow C^1(\mc V)$ (cf. \eqref{eq:20250702:eq1bis}).
Let us generalize this construction to arbitrary $n\in\mb Z_{\geq0}$ by defining projection operators $\widetilde{\proj}_n:\widetilde{\Gamma}^n(\mc V)\to C^n(\mc V)$ which factor through the quotient and induce linear maps $\proj_n:\Gamma^n(\mc V)\rightarrow C^n(\mc V)$.\glslink{ProjOp}{}

For $f\in\mc V=\widetilde{\Gamma}^0(\mc V)$ we let $\widetilde{\proj}_0(f)=\tint f$.
Clearly, the induced map $P_0:\Gamma^0(\mc V)\to C^0(\mc V)$ is the identity map.
For $X\in\widetilde{\Gamma}^n(\mc V)$, $n\geq1$, we define a linear map
$\widetilde{\proj}_n(X):\mc V^{\otimes n}\to \mc V^{\otimes n}[\lambda_1,\dots,\lambda_{n-1}]$ by the formula
\begin{equation}\label{eq:Pn}
\widetilde{\proj}_n(X)
=\frac1n \sum_{s=0}^{n-1}(-1)^{s(n-s)}|_{\lambda_n=\lambda_n^\dagger}
\mult_{(s+1,s+2)}\circ \sigma^{s+1}\circ X_{\lambda_{\sigma^s(1)},\dots,\lambda_{\sigma^s(n)}}\circ \sigma^{-s}\,,
\end{equation}
where we are using \eqref{eq:dagger} and \eqref{eq:mii+1}.
Recall that $X:\mc V^{\otimes n}\to \mc V^{\otimes (n+1)}[\lambda_1,\dots,\lambda_{n}]$, thus, in \eqref{eq:Pn} the rightmost
permutation $\sigma$ is the cyclic permutation of $\mc V^{\otimes n}$
which is also applied to the indices of the $\lambda_j$,
while the leftmost $\sigma$ is the cyclic permutation 
of $\mc V^{\otimes (n+1)}$ (following Section \ref{sec:1.1} we always denote 
by $\sigma$ the cyclic permutation $(12\dots k)$ acting on $\mc V^{\otimes k}$, and the number of factors $k$ will always be clear from the context). Note that, for $n=1$, the RHS of
\eqref{eq:Pn} applied to $a\in\mc V$  coincides with the RHS of \eqref{eq:P1}.
\begin{proposition}\label{20230814:prop1}
\begin{enumerate}[(a)]
\item Let $X\in\widetilde{\Gamma}^n(\mc V)$. Then $\widetilde{\proj}_n(X)\in C^n(\mc V)$.
\item Let $X\in\widetilde{\Gamma}^n(\mc V)$. Then $\widetilde{\proj}_n(\partial X)=0$.
\item Let $X\in\widetilde{\Gamma}^m(\mc V)$ and $Y\in\widetilde{\Gamma}^n(\mc V)$. Then $\widetilde{\proj}_{m+n}([X,Y])=0$.
\end{enumerate}
\end{proposition}
\begin{proof}
For $n=0,1$, the claim of part (a) has been already verified.
Let $n\geq2$ and 
let $X\in\widetilde{\Gamma}^n(\mc V)$. To prove the claim of part (a) we need to show $\widetilde{\proj}_n(X)$ satisfies sesquilinearity \eqref{20140702:eq4}
and \eqref{20140702:eq5},
the Leibniz rule \eqref{20140702:eq6} and  skewsymmetry \eqref{eq:nfold-skew}.
We start by noticing that the sesquilinearity \eqref{eq:sesquimaps} for $X\in\widetilde{\Gamma}^n(\mc V)$
can be rewritten as
$$
X_{\lambda_1,\dots,\lambda_n}\circ\partial_{(i)}
=-\lambda_iX_{\lambda_1,\dots,\lambda_n}\,,
$$
where we are using the notation \eqref{20240805:eq1} applied to $\partial$. Moreover, it is straightforward to check that
\begin{equation}\label{useful1}
\sigma^s\circ \partial_{(i)}=\partial_{(\sigma^s(i))}\circ \sigma^s
\,,
\qquad s\in\mb Z\,.
\end{equation}
Hence, for every $i=1,\dots,n$, we have
\begin{equation}\label{20240624:eq1}
X_{\lambda_{\sigma^s(1)},\dots,\lambda_{\sigma^s(n)}}
\circ \partial_{(\sigma^{-s}(i))}
=-\lambda_iX_{\lambda_{\sigma^s(1)},\dots,\lambda_{\sigma^s(n)}}
\,.
\end{equation}
Then, it follows immediately from \eqref{eq:Pn}, \eqref{useful1}
and \eqref{20240624:eq1} that $\widetilde{\proj}_n(X)$ satisfies the sesquilinearity axioms \eqref{20140702:eq4} and
\eqref{20140702:eq5}.
Next, let us prove that $\widetilde{\proj}_n(X)$ satisfies the Leibniz rules
\eqref{20140702:eq6} for $i=1,\dots,n$. Using the
definition of $\widetilde{\proj}_n(X)$ and the Leibniz rules \eqref{eq:Leibnizmaps} satisfied by $X$ we get
($a_1,\dots,a_n,b,c\in\mc V$)
\begin{align*}
&\widetilde{\proj}_n(X)_{\lambda_1,\dots,\lambda_{n-1}}(a_1,\dots, a_{i-1}, bc,a_{i+1},\dots,a_n)
\\
&
=\frac1n\sum_{s=0}^{n-1}(-1)^{s(n-s)}|_{\lambda_n\to\lambda_n^\dagger}
\mult_{(s+1,s+2)}\circ\sigma^{s+1}\circ
\\
&
\Big[(|_{x=\partial}b)\star_{\sigma^{-s}(i)}
\left(X_{\lambda_{\sigma^{s}(1)},\dots,\lambda_i+x,\dots \lambda_{\sigma^{s}(n)}}
\circ \sigma^{-s}
(a_1\otimes\dots\otimes a_{i-1}\otimes c\otimes a_{i+1}\otimes \dots\otimes a_n)\right)
\\
&+
\left(X_{\lambda_{\sigma^{s}(1)},\dots,\lambda_i+x,\dots,\lambda_{\sigma^{s}(n)}}
\circ \sigma^{-s}
(a_1\otimes\dots\otimes a_{i-1}\otimes b\otimes a_{i+1}\otimes \dots\otimes a_n)\right)
\\
&\quad\quad\quad\quad\quad\quad\quad\quad\quad\quad\quad\quad\quad\quad\quad\quad\quad\quad\quad\quad\quad\quad\quad\quad
\star_{n+1-\sigma^{-s}(i)}(|_{x=\partial}c)
\Big]
\\
&=
(|_{x=\partial} b)\star_i\widetilde{\proj}_n(X)_{\lambda_1,\dots,\lambda_{i-1},\lambda_i+x,\lambda_{i+1},\dots\lambda_{n-1}}
(a_1,\dots,a_{i-1},c,a_{i+1},\dots,a_n)
\\
&+\widetilde{\proj}_n(X)_{\lambda_1,\dots,\lambda_{i-1},\lambda_i+x,\lambda_{i+1},\dots\lambda_{n-1}}
(a_1,\dots,a_{i-1},b,a_{i+1},\dots,a_n) \star_{n-i} (|_{x=\partial} c)\,,
\end{align*}
where in the last equality we used the identities
\eqref{20240724:eq4} and \eqref{20240724:eq5}.
Finally, let us prove that $\widetilde{\proj}_n(X)$ satisfies the skewsymmetry condition
\eqref{eq:nfold-skew}.
From the definition of $\widetilde{\proj}_n(X)$ given by \eqref{eq:Pn} we have
\begin{align*}
&|_{\lambda_n=\lambda_n^\dagger}\sigma\circ \widetilde{\proj}_n(X)_{\lambda_{\sigma(1)},\dots,\lambda_{\sigma(n-1)}}
\circ \sigma^{-1}
\\
&=\frac{1}{n}\sum_{s=0}^{n-1}(-1)^{s(n-s)}|_{\lambda_n=\lambda_n^\dagger}\sigma\circ \mult_{(s+1,s+2)}
\circ\sigma^{s+1}\circ X_{\lambda_{\sigma^{s+1}(1)},\dots,\lambda_{\sigma^{s+1}(n)}}
\circ \sigma^{-s-1}
\\
&=\frac{(-1)^{n+1}}{n}\sum_{s=1}^{n}(-1)^{s(n-s)}|_{\lambda_n=\lambda_n^\dagger}\sigma\circ \mult_{(s,s+1)}
\circ\sigma^{s}\circ X_{\lambda_{\sigma^{s}(1)},\dots,\lambda_{\sigma^{s}(n)}}
\circ \sigma^{-s}
\\
&=
\frac{(-1)^{n+1}}{n}\sum_{s=1}^{n-1}(-1)^{s(n-s)}|_{\lambda_n=\lambda_n^\dagger}\sigma\circ \mult_{(s,s+1)}
\circ\sigma^{s}\circ X_{\lambda_{\sigma^{s}(1)},\dots,\lambda_{\sigma^{s}(n)}}
\circ \sigma^{-s}
\\
&
+\frac{(-1)^{n+1}}{n}|_{\lambda_n=\lambda_n^\dagger}\sigma\circ \mult_{(n,n+1)}
\circ\sigma^{-1}\circ X_{\lambda_1,\dots,\lambda_{n}}
\\
&
=
\frac{(-1)^{n+1}}{n}\sum_{s=1}^{n-1}(-1)^{s(n-s)}|_{\lambda_n=\lambda_n^\dagger} \mult_{(s+1,s+2)}
\circ\sigma^{s+1}\circ X_{\lambda_{\sigma^{s}(1)},\dots,\lambda_{\sigma^{s}(n)}}
\circ \sigma^{-s}
\\
&+
\frac{(-1)^{n+1}}{n}|_{\lambda_n=\lambda_n^\dagger}\mult_{(1,2)} \circ \sigma
\circ X_{\lambda_1,\dots,\lambda_{n}}
=(-1)^{n+1}\widetilde{\proj}_n(X)_{\lambda_1,\dots,\lambda_{n-1}}\,.
\end{align*}
In the second equality above we shifted the index of summation $s$, in the fourth equality we used equation \eqref{20240724:eq1} for $k=1$ and in the last equality we used \eqref{eq:Pn}. 
This shows that $\widetilde{\proj}_n(X)$ satisfies \eqref{eq:nfold-skew2}
and completes the proof of part (a).

Part (b) follows by definition for $n=0$, and, for $n\geq1$, it follows immediately from
the action of $\partial$ given by \eqref{20230802:eq1} and the
fact that $(\lambda_{\sigma^s(1)}+\dots+\lambda_{\sigma^s(n)}+\partial)|_{\lambda_n=\lambda_n^\dagger}=0$, for every $s=0,\dots,n-1$.

The claim of part (c) for $m=n=0$ follows by definition of the integral map $\tint$. Let $(m,n)\in\mb Z_{\geq0}^2\setminus\{(0,0)\}$, and let $X\in\widetilde{\Gamma}^m(\mc V)$
and $Y\in\widetilde{\Gamma}^n(\mc V)$. By definition of the projection
map $\widetilde{\proj}_{m+n}$ given by \eqref{eq:Pn} we have
\begin{align*}
&\widetilde{\proj}_{m+n}(XY)_{\lambda_1,\dots,\lambda_{m+n-1}}
\\
%&=
%\frac1{m+n} \sum_{s=0}^{m+n-1}(-1)^{s(m+n-s)}|_{\lambda_{m+n}=\lambda_{m+n}^\dagger}
%\mult_{(s+1,s+2)}\circ \sigma^{s+1}\circ (XY)_{\lambda_{\sigma^s(1)},\dots,\lambda_{\sigma^s(m+n)}}\circ \sigma^{-s}
%\\
&=
\frac1{m+n} \sum_{s=0}^{n-1}(-1)^{s(m+n-s)}|_{\lambda_{m+n}=\lambda_{m+n}^\dagger}
\mult_{(s+1,s+2)}\circ \sigma^{s+1}
\\
&\quad\quad\quad\quad\quad\quad\quad\quad\quad\quad\quad\quad\quad\quad\quad\quad
\circ (XY)_{\lambda_{\sigma^s(1)},\dots,\lambda_{\sigma^s(m+n)}}\circ \sigma^{-s}
\\
&+
\frac1{m+n} \sum_{s=n}^{m+n-1}(-1)^{s(m+n-s)}|_{\lambda_{m+n}=\lambda_{m+n}^\dagger}
\mult_{(s+1,s+2)}\circ \sigma^{s+1}
\\
&\quad\quad\quad\quad\quad\quad\quad\quad\quad\quad\quad\quad\quad\quad\quad\quad
\circ (XY)_{\lambda_{\sigma^s(1)},\dots,\lambda_{\sigma^s(m+n)}}\circ \sigma^{-s}
\\
&=
\frac1{m+n} \sum_{s=0}^{n-1}(-1)^{s(m+n-s)}|_{\lambda_{m+n}=\lambda_{m+n}^\dagger}
\mult_{(m+s+1,m+s+2)}\circ \sigma^{m+s+1}
\\
&\quad\quad\quad\quad\quad\quad\quad\quad\quad\quad\quad\quad\quad\quad\quad\quad
\circ (YX)_{\lambda_{\sigma^{s+m}(1)},\dots,\lambda_{\sigma^{s+m}(m+n)}}\circ \sigma^{-s-m}
\\
&+
\frac1{m+n} \sum_{s=n}^{m+n-1}(-1)^{s(m+n-s)}|_{\lambda_{m+n}=\lambda_{m+n}^\dagger}
\mult_{(s-n+1,s-n+2)}\circ \sigma^{s-n+1}
\\
&\quad\quad\quad\quad\quad\quad\quad\quad\quad\quad\quad\quad\quad\quad\quad\quad
\circ (YX)_{\lambda_{\sigma^{s-n}(1)},\dots,\lambda_{\sigma^{s-n}(m+n)}}\circ \sigma^{-s+n}
\\
&=(-1)^{mn}\widetilde{\proj}_{m+n}(YX)_{\lambda_1,\dots,\lambda_{m+n}}\,,
\end{align*}
where in the second equality we used equation \eqref{20240724:eq3},
and in the last equality we changed the order of summation and used the definition \eqref{eq:Pn} of $\widetilde{\proj}_{m+n}$ again. The above computation shows
that $\widetilde{\proj}_{m+n}([X,Y])=\widetilde{\proj}_{m+n}(XY-(-1)^{mn}YX)=
\widetilde{\proj}_{m+n}(XY)-(-1)^{mn}\widetilde{\proj}_{m+n}(YX)=0$ thus concluding the proof
of part (c).
\end{proof}
The next result follows immediately from Proposition \ref{20230814:prop1}(a), (b) and (c) 
and the definition \eqref{eq:Omega} of the quotient space $\Gamma(\mc V)$.
\begin{corollary}\label{20240819:cor2}
For every $n\in\mb Z_{\geq0}$ we have well defined linear maps
\begin{equation}\label{eq:projmaps}
\proj_n:\Gamma^n(\mc V)\to C^n(\mc V)
\end{equation}
given by
$$
\proj_n([X])=\widetilde{\proj}_n(X)
$$
for every $X\in\widetilde{\Gamma}^n(\mc V)$, where $[X]$ denotes the coset of $X$ in
$\Gamma^n(\mc V)$. Hence, by linearity, we have a map
$$
\proj:\Gamma(\mc V)\to C(\mc V)
\,.
$$
\end{corollary}
We prove in the subsequent Proposition \ref{20240819:prop2} that the map $\proj$ is an 
isomorphism provided that $\mc V$ is an algebra of noncommutative
differential polynomials. 
\begin{remark}
The map $\proj:\Gamma(\mc V)\to C(\mc V)$ is the ``double" analogue of the injective map \eqref{20250702:eq2}. However, we were
not able to determine injectivity
in this case and we leave it as an open question.
In general, the map $\proj$ is not surjective as shown in the following example. Let $\mc V=\kk\langle x,u^{(n)}\mid n\in\mb Z_{\geq0}\rangle$ be the algebra of polynomials in
infinitely many noncommutative variables $x,u^{(0)},u^{(1)},\dots$. We make it a differential algebra by defining $\partial(x)=0$, $\partial(u^{(n)})=u^{(n+1)}$, $n\in\mb Z_{\geq0}$, and extending it  to $\mc V$ by the Leibniz rule. Let $D\in C^1(\mc V)$ be the derivation of $\mc V$ (commuting with $\partial$) such that $D(x)=1$ and $D(u^{(n)})=0$, for every $n\in\mb Z_{\geq0}$ (cf. Remark \ref{20250722:rem2}). Let $X\in\widetilde{\Gamma}^1(\mc V)$. By sesquilinearity \eqref{eq:sesquimaps}
we have that $0=X_\lambda(\partial(x))=-\lambda X_\lambda(x)$, from which follows that $X_\lambda(x)=0$.
Hence, from \eqref{eq:P1} we see that there is no $X\in\widetilde{\Gamma}^1(\mc V)$ such that
$\widetilde{\proj}_1(X)=D$.
\end{remark}

%%%
\subsection[The projection operator \& cohomology]{The projection operator and double Poisson vertex algebra cohomologies}
Recall the definition of the differential $\tilde\delta$ of $\widetilde{\Gamma}(\mc V)$ given in \eqref{eq:diff} (note that, in this chapter we are denoting the $2$-fold dPVA $\lambda$-bracket on $\mc V$ by $\llbracket\cdot\,_{\lambda}\,\cdot\rrbracket$)
and of the differential $\dd$ of $C(\mc V)$ 
given in \eqref{eq:dP-1}. We have the following compatibility conditions between
the complexes $(\widetilde{\Gamma}(\mc V),\tilde\delta)$, $(C(\mc V),\dd)$ and the projection
operator $\widetilde{\proj}:\widetilde{\Gamma}(\mc V)\rightarrow C(\mc V)$.
\begin{proposition}\label{20250720.prop1}
For every $X\in\widetilde{\Gamma}^n(\mc V)$, $n\in\mb Z_{\geq0}$, we have
\begin{equation}\label{delta-comm}
(n + 1)\widetilde{\proj}_{n+1}(\tilde\delta(X))=(n+\delta_{n,0})(-1)^n \dd(\widetilde{\proj}_n(X))
\,.
\end{equation}
\end{proposition}
\begin{proof}
Equation \eqref{delta-comm} for $n=0$ can be checked directly using \eqref{eq:P1},
\eqref{exa1}, \eqref{eq:dP0} and skewsymmetry \eqref{eq:skew2}.
For $n\geq1$, using equations \eqref{eq:Pn} and \eqref{eq:diff}, the LHS of \eqref{delta-comm}
applied to $a_1\otimes\dots\otimes a_{n+1}\in\mc V^{\otimes (n+1)}$ is
\begin{equation}
\begin{split}\label{20240820:eq1}
&(n + 1)\widetilde{\proj}_{n+1}(\tilde\delta(X))_{\lambda_1,\dots,\lambda_n}(a_1,\dots,a_{n+1})
\\
&=
\sum_{s=0}^n\sum_{t=1}^{n+1}(-1)^{s(n-s)+s+t+1}|_{\lambda_{n+1}=\lambda_{n+1}^\dagger}
\mult_{(s+1,s+2)}\sigma^{s+1}
\\
&
\qquad\qquad\qquad\qquad\llbracket a_{\sigma^s(t)}{}_{\lambda_{\sigma^s(t)}}
X_{\lambda_{\sigma^s(1)},\stackrel{t}{\check{\dots}},\lambda_{\sigma^s(n+1)}}(a_{\sigma^s(1)},\stackrel{t}{\check{\dots}},a_{\sigma^s(n+1)})
\rrbracket_{(t)}
\\
&
+
\sum_{s=0}^n\sum_{t=1}^{n}(-1)^{s(n-s)+s+t}
A(s,t)
\,,
\end{split}
\end{equation}
where we set
\begin{equation}\label{Ast}
\begin{split}
&A(s,t)
=|_{\lambda_{n+1}=\lambda_{n+1}^\dagger}
\mult_{(s+1,s+2)}\sigma^{s+1} \\
&X^{(t)}_{\lambda_{\sigma^s(1)},\dots,\lambda_{\sigma^s(t-1)},\lambda_{\sigma^s(t)}+\lambda_{\sigma^s(t+1)},\lambda_{\sigma^s(t+2)},\dots,\lambda_{\sigma^s(n+1)}}
(a_{\sigma^s(1)},\dots
\\
&\quad\quad\quad\quad\quad
\dots,a_{\sigma^s(t-1)},\llbracket a_{\sigma^s(t)}{}_{\lambda_{\sigma^s(t)}}a_{\sigma^s(t+1)}\rrbracket,a_{\sigma^s(t+2)},\dots,a_{\sigma^s(n+1)})
\,.
\end{split}
\end{equation}
Using Lemma \ref{20140606:lem} we rewrite the first sum in \eqref{20240820:eq1} as
\begin{equation}\label{20240822:eq1}
\begin{split}
&
\sum_{t=1}^{n}\sum_{s=0}^{t-1}(-1)^{s(n-s)+t+1}|_{\lambda_{n+1}=\lambda_{n+1}^\dagger}
\mult_{(s+1,s+2)}
\\
&\quad\quad\quad\quad\quad\quad
\llbracket a_{t}{}_{\lambda_{t}}
\sigma^{s+1}X_{\lambda_{\sigma^s(1)},\stackrel{t-s}{\check{\dots}},\lambda_{\sigma^s(n+1)}}(a_{\sigma^s(1)},\stackrel{t-s}{\check{\dots}},a_{\sigma^s(n+1)})
\rrbracket_{(t+1)}
\\
&+
\sum_{t=1}^{n}\sum_{s=t}^n(-1)^{(s-1)(n-s+1)+t+1}|_{\lambda_{n+1}=\lambda_{n+1}^\dagger}
\mult_{(s+1,s+2)}
\\
&\quad\quad\quad\quad\quad\quad
\llbracket a_{t}{}_{\lambda_{t}}
\sigma^{s}X_{\lambda_{\sigma^s(1)},\stackrel{n+1+t-s}{\check{\dots\dots}},\lambda_{\sigma^s(n+1)}}(a_{\sigma^s(1)},\stackrel{n+1+t-s}{\check{\dots\dots}},a_{\sigma^s(n+1)})
\rrbracket_{(t)}
\\
&+(-1)^n\sum_{s=0}^n(-1)^{s(n-s)}
|_{\lambda_{n+1}=\lambda_{n+1}^\dagger}
\mult_{(s+1,s+2)}
\sigma^{-1}
\\
&\quad\quad\quad\quad\quad\quad
\llbracket a_{n+1}{}_{\lambda_{n+1}}
\sigma^{s+1}X_{\lambda_{\sigma^s(1)},\dots,\lambda_{\sigma^s(n)}}(a_{\sigma^s(1)},\dots,a_{\sigma^s(n)})
\rrbracket_{(1)}
\,.
\end{split}
\end{equation}
We use Lemma \ref{20240821:lem1}(a) and (b) to rewrite the first two sums in \eqref{20240822:eq1}
as
\small 
\begin{equation}\label{20240822:eq2}
\begin{split}
&
\sum_{t=1}^{n}\sum_{s=0}^{t-2}(-1)^{s(n-s)+t+1}|_{\lambda_{n+1}=\lambda_{n+1}^\dagger}
\mult_{(s+1,s+2)}
\\
&\quad\quad\quad
\llbracket a_{t}{}_{\lambda_{t}}
\sigma^{s+1}X_{\lambda_{\sigma^s(1)},\stackrel{t-s}{\check{\dots}},\lambda_{\sigma^s(n+1)}}(a_{\sigma^s(1)},\stackrel{t-s}{\check{\dots}},a_{\sigma^s(n+1)})
\rrbracket_{(t+1)}
\\
&+
\sum_{t=1}^{n}(-1)^{(t-1)(n+1-t)+t+1}|_{\lambda_{n+1}=\lambda_{n+1}^\dagger}
\\
&\quad\quad\quad
\left(
\mult_{(t,t+1)}
\llbracket a_{t}{}_{\lambda_{t}}
\sigma^{t}X_{\lambda_{\sigma^{t}(1)},\dots,\lambda_{\sigma^{t}(n)}}(a_{\sigma^{t}(1)},\dots,a_{\sigma^{t}(n)})
\rrbracket_{(t+1)}
\right.\\
&\quad\quad\quad
\left.+
\mult_{(t+1,t+2)}
\llbracket a_{t}{}_{\lambda_{t}}
\sigma^{t}X_{\lambda_{\sigma^t(1)},\dots,\lambda_{\sigma^t(n)}}(a_{\sigma^t(1)},\dots,a_{\sigma^t(n)})
\rrbracket_{(t)}
\right)
\\
&+
\sum_{t=1}^{n}\sum_{s=t+1}^n(-1)^{(s-1)(n-s+1)+t+1}|_{\lambda_{n+1}=\lambda_{n+1}^\dagger}
\mult_{(s+1,s+2)}
\\
&\quad\quad\quad
\llbracket a_{t}{}_{\lambda_{t}}
\sigma^{s}X_{\lambda_{\sigma^s(1)},\stackrel{n+1+t-s}{\check{\dots\dots}},\lambda_{\sigma^s(n+1)}}(a_{\sigma^s(1)},\stackrel{n+1+t-s}{\check{\dots\dots}},a_{\sigma^s(n+1)})
\rrbracket_{(t)}
\\
&
=\sum_{t=1}^{n}\sum_{s=0}^{t-1}(-1)^{s(n-s)+t+1}|_{\lambda_{n+1}=\lambda_{n+1}^\dagger}
\\
&\quad\quad\quad
\llbracket a_{t}{}_{\lambda_{t}}
\mult_{(s+1,s+2)}\sigma^{s+1}X_{\lambda_{\sigma^s(1)},\stackrel{t-s}{\check{\dots}},\lambda_{\sigma^s(n+1)}}(a_{\sigma^s(1)},\stackrel{t-s}{\check{\dots}},a_{\sigma^s(n+1)})
\rrbracket_{(t)}
\\
&+
\sum_{t=1}^{n}\sum_{s=t+1}^n(-1)^{(s-1)(n-s+1)+t+1}|_{\lambda_{n+1}=\lambda_{n+1}^\dagger}
\\
&\quad\quad\quad
\llbracket a_{t}{}_{\lambda_{t}}
\mult_{(s,s+1)}\sigma^{s}X_{\lambda_{\sigma^s(1)},\stackrel{n+1+t-s}{\check{\dots\dots}},\lambda_{\sigma^s(n+1)}}(a_{\sigma^s(1)},\stackrel{n+1+t-s}{\check{\dots\dots}},a_{\sigma^s(n+1)})
\rrbracket_{(t)}
\,.
\end{split}
\end{equation}
\normalsize 
In a similar way, using \eqref{20240724:eq1} and Lemma \ref{20240821:lem1},
and using skewsymmetry \eqref{eq:skew2} (as for the derivation of \eqref{eq:Bn}), we can rewrite the third summand in \eqref{20240822:eq1} as
\begin{equation}\label{20240822:eq3}
\begin{split}
&(-1)^{n+1}\sum_{s=0}^{n-1}(-1)^{s(n-s)}
\\
&\quad\quad
\llbracket 
\mult_{(s+1,s+2)}
\sigma^{s+1}X_{\lambda_{\sigma^s(1)},\dots,\lambda_{\sigma^s(n)}}(a_{\sigma^s(1)},\dots,a_{\sigma^s(n)})_{\lambda_1+\dots+\lambda_n}a_{n+1}
\rrbracket_{L}
\,.
\end{split}
\end{equation}
Combining equations \eqref{20240822:eq1}, \eqref{20240822:eq2} and \eqref{20240822:eq3} we rewrite \eqref{20240820:eq1}
as follows
\begin{equation}
\begin{split}\label{LHS}
&(n + 1)\widetilde{\proj}_{n+1}(\tilde\delta(X))_{\lambda_1,\dots,\lambda_n}(a_1,\dots,a_{n+1})
\\
&=
\sum_{t=1}^{n}\sum_{s=0}^{t-1}(-1)^{s(n-s)+t+1}|_{\lambda_{n+1}=\lambda_{n+1}^\dagger}
\\
&\quad\quad
\llbracket a_{t}{}_{\lambda_{t}}
\mult_{(s+1,s+2)}\sigma^{s+1}X_{\lambda_{\sigma^s(1)},\stackrel{t-s}{\check{\dots}},\lambda_{\sigma^s(n+1)}}(a_{\sigma^s(1)},\stackrel{t-s}{\check{\dots}},a_{\sigma^s(n+1)})
\rrbracket_{(t)}
\\
&+
\sum_{t=1}^{n}\sum_{s=t+1}^n(-1)^{s(n-s)+n+t}|_{\lambda_{n+1}=\lambda_{n+1}^\dagger}
\\
&\quad\quad
\llbracket a_{t}{}_{\lambda_{t}}
\mult_{(s,s+1)}\sigma^{s}X_{\lambda_{\sigma^s(1)},\stackrel{n+1+t-s}{\check{\dots\dots}},\lambda_{\sigma^s(n+1)}}(a_{\sigma^s(1)},\stackrel{n+1+t-s}{\check{\dots\dots}},a_{\sigma^s(n+1)})
\rrbracket_{(t)}
\\
&
+(-1)^{n+1}\sum_{s=0}^{n-1}(-1)^{s(n-s)}
\\
&\quad\quad
\llbracket 
\mult_{(s+1,s+2)}
\sigma^{s+1}X_{\lambda_{\sigma^s(1)},\dots,\lambda_{\sigma^s(n)}}(a_{\sigma^s(1)},\dots,a_{\sigma^s(n)})_{\lambda_1+\dots+\lambda_n}a_{n+1}
\rrbracket_{L}
\\
&
+ \sum_{s=0}^n\sum_{t=1}^{n}(-1)^{s(n-s)+s+t} \,  A(s,t)\,.
\end{split}
\end{equation}
On the other hand, by definition of the projection operators $\widetilde{\proj}_n$
given in \eqref{eq:Pn}, we have
\begin{align*}
&n \widetilde{\proj}_n(X)_{\lambda_1,\stackrel{t}{\check{\dots}},\lambda_n}
(a_1,\stackrel{t}{\check{\dots}},a_{n+1})
\\
&=\sum_{s=0}^{t-1}
(-1)^{s(n-s)}|_{\lambda_{n+1}=\lambda_{n+1}^\dagger+\lambda_t}
\mult_{(s+1,s+2)}\sigma^{s+1}
\\
&\quad\quad\quad\quad\quad\quad\quad\quad\quad
X_{\lambda_{\sigma^s(1)},\stackrel{t-s}{\check{\dots}},\lambda_{\sigma^s(n+1)}}
(a_{\sigma^s(1)},\stackrel{t-s}{\check{\dots}},a_{\sigma^s(n+1)})
\\
&+(-1)^{n+1}\sum_{s=t+1}^{n}
(-1)^{s(n-s)}|_{\lambda_{n+1}=\lambda_{n+1}^\dagger+\lambda_t}
\mult_{(s,s+1)}\sigma^{s}
\\
&\quad\quad\quad\quad\quad\quad\quad\quad\quad
X_{\lambda_{\sigma^s(1)},\stackrel{n+1+t-s}{\check{\dots\dots}},\lambda_{\sigma^s(n+1)}}
(a_{\sigma^s(1)},\stackrel{n+1+t-s}{\check{\dots\dots}},a_{\sigma^s(n+1)})
\,,
\end{align*}
for every $t=1,\dots,n$.
Hence, using sesquilinearity \eqref{eq:sesqui}, the definition \eqref{eq:dP-1-Q} of the differential $\dd$ and
the definition \eqref{eq:Pn} of $\widetilde{\proj}_n(X)$ we have that the RHS of \eqref{delta-comm}
applied to $a_1\otimes\dots\otimes a_{n+1}\in\mc V^{\otimes (n+1)}$ is
\begin{equation}\label{20240823:eq1}
\begin{split}
&n(-1)^n \dd(\widetilde{\proj}_n(X))_{\lambda_1,\dots,\lambda_n}(a_1,\dots,a_{n+1})
\\
&=\sum_{t=1}^{n}\sum_{s=0}^{t-1}(-1)^{s(n-s)+t+1}
|_{\lambda_{n+1}=\lambda_{n+1}^\dagger}
\\
&\quad\quad
\llbracket a_t {}_{\lambda_t}
\mult_{(s+1,s+2)}\sigma^{s+1}
X_{\lambda_{\sigma^s(1)},\stackrel{t-s}{\check{\dots}},\lambda_{\sigma^s(n+1)}}
(a_{\sigma^s(1)},\stackrel{t-s}{\check{\dots}},a_{\sigma^s(n+1)})
\rrbracket_{(t)}
\\
&+\sum_{t=1}^{n}\sum_{s=t+1}^{n}(-1)^{s(n-s)+n+t}
|_{\lambda_{n+1}=\lambda_{n+1}^\dagger}
\\
&\quad\quad
\llbracket a_t {}_{\lambda_t}
\mult_{(s,s+1)}\sigma^{s}
X_{\lambda_{\sigma^s(1)},\stackrel{n+1+t-s}{\check{\dots\dots}},\lambda_{\sigma^s(n+1)}}
(a_{\sigma^s(1)},\stackrel{n+1+t-s}{\check{\dots\dots}},a_{\sigma^s(n+1)})
\rrbracket_{(t)}
\\
&+(-1)^{n+1}
\sum_{s=0}^{n-1}(-1)^{s(n-s)}
\\
&\quad\quad
\llbracket 
\mult_{(s+1,s+2)}\sigma^{s+1}X_{\lambda_{\sigma^s(1)},\dots,\lambda_{\sigma^s(n)}}
(a_{\sigma^s(1)},\dots,a_{\sigma^s(n)})
_{\lambda_1+\dots+\lambda_n}a_{n+1}\rrbracket_{L}
\\
&+n\sum_{t=1}^n(-1)^{t}
\widetilde{\proj}_n(X)^{(t)}_{\lambda_1,\dots,\lambda_{t-1},\lambda_{t}+\lambda_{t+1},\lambda_{t+2},\dots,\lambda_n}
(a_1,\dots
\\
&\quad\quad
\dots,a_{t-1},\llbracket a_{t}{}_{\lambda_t}a_{t+1}\rrbracket,a_{t+2},\dots,a_{n+1})
\\
&+n
|_{\lambda_{n+1}=\lambda_{n+1}^\dagger}
\sigma \widetilde{\proj}_n(X)^{(n)}_{\lambda_2,\dots,\lambda_n}(a_2,\dots, a_n,\llbracket a_{n+1}{}_{\lambda_{n+1}}a_{1}\rrbracket)
\,.
\end{split}
\end{equation}
Comparing \eqref{LHS} and \eqref{20240823:eq1} we see that, to conclude the proof of the identity \eqref{delta-comm} we are left to show that
\begin{equation}\label{20240823:eq2}
\begin{split}
&\sum_{s=0}^n\sum_{t=1}^{n}(-1)^{s(n-s)+s+t} \, A(s,t)
\\
&=n\sum_{t=1}^n(-1)^{t}
\widetilde{\proj}_n(X)^{(t)}_{\lambda_1,\dots,\lambda_{t-1},\lambda_{t}+\lambda_{t+1},\lambda_{t+2},\dots,\lambda_n}
(a_1,\dots
\\
&\quad\quad\quad\quad\quad\quad\quad\quad
\dots,a_{t-1},\llbracket a_{t}{}_{\lambda_t}a_{t+1}\rrbracket,a_{t+2},\dots,a_{n+1})
\\
&+n
|_{\lambda_{n+1}=\lambda_{n+1}^\dagger}
\sigma \widetilde{\proj}_n(X)^{(n)}_{\lambda_2,\dots,\lambda_n}(a_2,\dots, a_n,\llbracket a_{n+1}{}_{\lambda_{n+1}}a_{1}\rrbracket)
\,,
\end{split}
\end{equation}
where $A(s,t)$ is defined in \eqref{Ast}. Using the definition of the projection operators
and the notation \eqref{eq:XL} we have
\begin{equation}\label{20240823:eq3}
\begin{split}
&n\widetilde{\proj}_n(X)^{(t)}_{\lambda_1,\dots,\lambda_{t-1},\lambda_{t}+\lambda_{t+1},\lambda_{t+2},\dots,\lambda_n}
(a_1,\dots,\llbracket a_{t}{}_{\lambda_t}a_{t+1}\rrbracket,\dots,a_{n+1})
\\
&
=n\widetilde{\proj}_n(X)_{\lambda_1,\dots,\lambda_{t-1},\lambda_{t}+\lambda_{t+1}+x,\lambda_{t+2},\dots,\lambda_n}
(a_1,\dots
\\
&\quad\quad\quad\quad
\dots,a_{t-1},\llbracket a_{t}{}_{\lambda_t}a_{t+1}\rrbracket',a_{t+2},\dots,a_{n+1})
\otimes_{n-t}(|_{x=\partial}\llbracket a_{t}{}_{\lambda_t}a_{t+1}\rrbracket'')
\\
&=\sum_{s=0}^{t-1}(-1)^{s(n-s)}|_{\lambda_{n+1}=\lambda_{n+1}^\dagger}
\Big(
\mult_{(s+1,s+2)}\sigma^{s+1}
\\
&\quad\quad\quad\quad
X_{\lambda_{\sigma^s(1)},\dots,\lambda_{\sigma^s(t-s)}+\lambda_{\sigma^s(t-s+1)}+x,\dots,\lambda_{\sigma^s(n+1)}}
(a_{\sigma^s(1)},\dots
\\
&\quad\quad\quad\quad
\dots,\llbracket a_{\sigma^s(t-s)}{}_{\lambda_{\sigma^s(t-s)}}a_{\sigma^s(t-s+1)}\rrbracket',\dots,a_{\sigma^s(n+1)})
\Big)
\\
&\quad\quad\quad\quad
\otimes_{n-t}(|_{x=\partial}\llbracket a_{\sigma^s(t-s)}{}_{\lambda_{\sigma^s(t-s)}}a_{\sigma^s(t-s+1)}\rrbracket'')
\,.
\\
&+
\sum_{s=t}^{n}(-1)^{s(n-s)}|_{\lambda_{n+1}=\lambda_{n+1}^\dagger}
\Big(
\mult_{(s+1,s+2)}\sigma^{s+1}
\\
&\quad\quad\quad\quad
X_{\lambda_{\sigma^{s+1}(1)},\dots,\lambda_{\sigma^{s+1}(n+t-s)}+\lambda_{\sigma^{s+1}(n+1+t-s)}+x,\dots,\lambda_{\sigma^{s+1}(n+1)}}
(a_{\sigma^{s+1}(1)},\dots
\\
&\quad\quad\quad\quad
\dots,\llbracket a_{\sigma^{s+1}(n+t-s)}
{}_{\lambda_{\sigma^{s+1}(n+t-s)}}a_{\sigma^{s+1}(n+1+t-s)}\rrbracket',\dots,a_{\sigma^{s+1}(n+1)})
\Big)
\\
&\quad\quad\quad\quad
\otimes_{n-t}(|_{x=\partial}\llbracket a_{\sigma^{s+1}(n+t-s)}{}_{\lambda_{\sigma^{s+1}(n+t-s)}}a_{\sigma^{s+1}(n+1+t-s)}\rrbracket'')
\\
&=\sum_{s=0}^{t-1}(-1)^{s(n-s)}|_{\lambda_{n+1}=\lambda_{n+1}^\dagger}
\mult_{(s+1,s+2)}\sigma^{s+1}
\\
&\quad\quad\quad\quad
X_{\lambda_{\sigma^s(1)},\dots,\lambda_{\sigma^s(t-s)}+\lambda_{\sigma^s(t-s+1)},\dots,\lambda_{\sigma^s(n+1)}}^{(t-s)}
(a_{\sigma^s(1)},\dots
\\
&\quad\quad\quad\quad
\dots,\llbracket a_{\sigma^s(t-s)}{}_{\lambda_{\sigma^s(t-s)}}a_{\sigma^s(t+1)}\rrbracket,\dots,a_{\sigma^s(n+1)})
\\
&
+\sum_{s=t}^{n-1}(-1)^{s(n-s)}|_{\lambda_{n+1}=\lambda_{n+1}^\dagger}
\mult_{(s+2,s+3)}\sigma^{s+2}
\\
&\quad\quad\quad\quad
X_{\lambda_{\sigma^{s+1}(1)},\dots,\lambda_{\sigma^{s+1}(n+t-s)}+\lambda_{\sigma^{s+1}(n+1+t-s)},\dots,\lambda_{\sigma^{s+1}(n+1)}}^{(n+t-s)}
(a_{\sigma^{s+1}(1)},\dots
\\
&\quad\quad\quad\quad
\dots,\llbracket a_{\sigma^{s+1}(n+t-s)}{}_{\lambda_{\sigma^{s+1}(n+t-s)}}a_{\sigma^{s+1}(n+1+t-s)}\rrbracket,\dots,a_{\sigma^{s+1}(n+1)})
\\
&
=\sum_{s=0}^{t-1}(-1)^{s(n-s)}A(s,t-s)
+\sum_{s=t+1}^{n}(-1)^{(s+1)(n+1-s)}A(s,n+1+t-s)
\,.
\end{split}
\end{equation}
In the second equality we used \eqref{20240823:eq4} and \eqref{eq:XL}
and in the last equality we used \eqref{Ast}.
Summing the RHS of \eqref{20240823:eq3}
over $t=1,\dots,n$ we get
\begin{equation}\label{20240823:eq5}
\begin{split}
&n\sum_{t=1}^n (-1)^t\widetilde{\proj}_n(X)^{(t)}_{\lambda_1,\dots,\lambda_{t-1},\lambda_{t}+\lambda_{t+1},\lambda_{t+2},\dots,\lambda_n}
(a_1,\dots
\\
&\quad\quad\quad\quad\quad\quad\quad\quad\quad
\dots, a_{t-1},\llbracket a_{t}{}_{\lambda_t}a_{t+1}\rrbracket,a_{t+2},\dots,a_{n+1})
\\
&=
\sum_{t=1}^n\sum_{s=0}^n
(-1)^{s(n-s)+s+t}\left(1-\delta_{s,n+1-t}
\right)A(s,t)
\,.
\end{split}
\end{equation}
Finally, using \eqref{eq:XL} and \eqref{eq:Pn} we have
\begin{align*}
&n
|_{\lambda_{n+1}=\lambda_{n+1}^\dagger}
\sigma \widetilde{\proj}_n(X)^{(n)}_{\lambda_2,\dots,\lambda_n}(a_2,\dots, a_n,\llbracket a_{n+1}{}_{\lambda_{n+1}}a_{1}\rrbracket)
\\
&=|_{\lambda_{n+1}=\lambda_{n+1}^\dagger}
\sigma\left(
\widetilde{\proj}_n(X)_{\lambda_2,\dots,\lambda_n}(a_2,\dots, a_n,\llbracket a_{n+1}{}_{\lambda_{n+1}}a_{1}\rrbracket')
\otimes\llbracket a_{n+1}{}_{\lambda_{n+1}}a_{1}\rrbracket''\right)
\\
&
=\sum_{s=0}^{n-1}(-1)^{s(n-s)}
|_{\lambda_{n+1}=\lambda_{n+1}^\dagger}
\sigma\Big(
\mult_{(s+1,s+2)}\sigma^{s+1}
\\
&\quad\quad\quad\quad\quad
X_{\lambda_{\sigma^{s+1}(1)},\dots,\lambda_1+\lambda_{n+1} +x,\dots,\lambda_{\sigma^{s+1}(n+1)}}
(a_{\sigma^{s+1}(1)},\dots
\\
&\quad\quad\quad\quad\quad
\dots,\llbracket a_{n+1}{}_{\lambda_{n+1}}a_1\rrbracket',\dots a_{\sigma^{s+1}(n+1)})
\otimes (|_{x=\partial} \llbracket a_{n+1}{}_{\lambda_{n+1}}a_{1}\rrbracket'')\,
\Big)\,.
\end{align*}
We note that, in the array $(\lambda_{\sigma^{s+1}(1)},\dots,\lambda_1+\lambda_{n+1}+x,\dots,\lambda_{\sigma^{s+1}(n+1)})$, the entry $\lambda_1+\lambda_{n+1}+x$ appears in position
$n-s$.
Using equation \eqref{20240823:eq4} we rewrite the above identity as
\begin{equation}\label{20240823:eq6}
\begin{split}
&|_{\lambda_{n+1}=\lambda_{n+1}^\dagger}
\sigma \widetilde{\proj}_n(X)^{(n)}_{\lambda_2,\dots,\lambda_n}(a_2,\dots, a_n,\llbracket a_{n+1}{}_{\lambda_{n+1}}a_{1}\rrbracket)
\\
&
=\sum_{s=0}^{n-1}(-1)^{s(n-s)}
|_{\lambda_{n+1}=\lambda_{n+1}^\dagger}
\sigma
\mult_{(s+1,s+2)}\sigma^{s+1}
\\
&
\quad\quad
X_{\lambda_{\sigma^{s+1}(1)},\dots,\lambda_1+\lambda_{n+1},\dots,\lambda_{\sigma^{s+1}(n+1)}}^{(n-s)}
(a_{\sigma^{s+1}(1)},\dots,\llbracket a_{n+1}{}_{\lambda_{n+1}}a_1\rrbracket,\dots a_{\sigma^{s+1}(n+1)})
\\
&
=\sum_{s=1}^{n}(-1)^{(s+1)(n+1-s)}
|_{\lambda_{n+1}=\lambda_{n+1}^\dagger}
\mult_{(s+1,s+2)}\sigma^{s+1}
\\
&
\quad\quad
X_{\lambda_{\sigma^{s}(1)},\dots,\lambda_1+\lambda_{n+1},\dots,\lambda_{\sigma^{s}(n+1)}}^{(n+1-s)}
(a_{\sigma^{s}(1)},\dots,\llbracket a_{n+1}{}_{\lambda_{n+1}}a_1\rrbracket,\dots a_{\sigma^{s}(n+1)})
\\
&
=
\sum_{s=1}^{n}(-1)^{(s+1)(n+1-s)}
A(s,n+1-s)
=
\sum_{s=0}^n
\sum_{t=1}^n
(-1)^{s(n-s)+s+t}
\delta_{t,n+1-s}A(s,t)
\,.
\end{split}
\end{equation}
In the second equality we used Lemma \eqref{20240724:lem1}(a), in the third equality we used the definition of $A(s,t)$ given by \eqref{Ast} and the last equality is straightforward. 
The identity \eqref{20240823:eq2} follows combining equations \eqref{20240823:eq5} and \eqref{20240823:eq6}. This concludes the proof of part (d).
\end{proof}
As a consequence of Proposition \ref{20250720.prop1} and Corollary \ref{20240819:cor2}
we get the following commutative diagrams:
%\pecetta{this diagram has been commented to speed up compiling}
%\begin{comment}
\begin{equation*}
\begin{tikzcd}
\Gamma^0(\mc V) \arrow[r,"\delta"] \arrow[d,equal]& \Gamma^1(\mc V)\arrow[d,"\proj_1"]
\\
C^0(\mc V) \arrow[r,"\dd"]& C^1(\mc V)
\end{tikzcd}
\,,
\quad
\begin{tikzcd}
\Gamma^n(\mc V) \arrow[r,"\delta"] \arrow[d,"n\proj_n"]& \Gamma^{n+1}(\mc V)\arrow[d,"(n+1)\proj_{n+1}"]
\\
C^n(\mc V) \arrow[r,"(-1)^n\dd"]& C^{n+1}(\mc V)
\end{tikzcd}
\,,
\quad n\geq1\,.
\end{equation*}
%\end{comment}
%
In conclusion, we proved the following result.
\begin{theorem}\label{20250720:cor1}
For every $n\in\mb Z_{\geq0}$ the linear maps $\proj_n:\Gamma^n(\mc V)\to C^n(\mc V)$ given
by \eqref{eq:projmaps} give linear maps in cohomology
$$
\dH_{\mathrm{red}}^n(\mc V)\rightarrow \dPVH^n(\mc V)
\,.
$$
In particular, $\dH^0_{\mathrm{red}}(\mc V)\simeq \dPVH^0(\mc V)$ provided that $\proj_1$ defined in \eqref{eq:20250702:eq1bis} is injective.
\end{theorem}

%%%%%%%%%%% NEW CHAPTER %%%%%%%%%%%%%%%
%%%%%%%%%%% NEW CHAPTER %%%%%%%%%%%%%%%
%%%%%%%%%%% NEW CHAPTER %%%%%%%%%%%%%%%
%%%%%%%%%%% NEW CHAPTER %%%%%%%%%%%%%%%

\chapter[Algebras of differential polynomials \& \MakeLowercase{d}PVA cohomology]{Algebras of differential polynomials and double Poisson vertex algebra cohomologies
}\label{sec:PVAdiff}

%%%
\section{Algebras of noncommutative differential polynomials}
The \emph{algebra of noncommutative differential polynomials} 
$\mc R_\ell$ in $\ell$ variables $u_i$, $i\in I=\{1,\dots,\ell\}$
is the algebra of noncommutative polynomials in the variables $u_i^{(n)}$,\glslink{Rell}{}
$$
\mc R_\ell=\kk\langle u_i^{(n)}\mid i\in I,n\in\mb Z_{\geq0}\rangle\,,
$$
endowed with a derivation $\partial$, defined on generators by $\partial u_i^{(n)}=u_i^{(n+1)}$, and extended to $\mc V$ by the Leibniz rule.

We introduce partial derivatives $\frac{\partial}{\partial u_i^{(n)}}:\mc R_\ell\to\mc R_\ell\otimes\mc R_\ell$, for every $i\in I$ and
$n\in\mb Z_{\geq0}$, defined on monomials by
\begin{equation}\label{20140627:eq1}\glslink{partialder-nc}{}
\frac{\partial}{\partial u_i^{(n)}} (u_{i_1}^{(n_1)}\dots u_{i_s}^{(n_s)})
=
\sum_{k=1}^s
\delta_{i_k,i}\delta_{n_k,n}
\,
u_{i_1}^{(n_1)}\dots u_{i_{k-1}}^{(n_{k-1})}
\otimes
u_{i_{k+1}}^{(n_{k+1})}\dots u_{i_s}^{(n_s)}
\,,
\end{equation}
which are commuting $2$-fold derivations of $\mc R_\ell$ such that
\begin{equation}\label{eq:comm}
\left[\frac{\partial}{\partial u_i^{(n)}},\partial\right]=\frac{\partial}{\partial u_i^{(n-1)}}\,,
\end{equation}
where the RHS is zero for $n=0$.
Note that
\begin{equation}\label{20250723:eq1}
\bigcap_{(i,n)\in I\times \mb Z_{\geq0}}\ker\frac{\partial}{\partial u_i^{(n)}}
=\ker\partial
=\kk
\,.
\end{equation}
%
%It is convenient to rewrite equation \eqref{eq:comm} in terms of its generating series
%as follows ($i\in I$):
%\begin{equation}\label{eq:comm_series}
%\sum_{n\in\mb Z}z^n\frac{\partial}{\partial u_i^{(n)}}\circ\partial
%=(z+\partial)\circ\sum_{n\in\mb Z}z^n\frac{\partial}{\partial u_i^{(n)}}
%\end{equation}
%
It is shown in \cite{DSKV} that the partial derivatives \emph{strongly commute}, that is
$$
\left(\frac{\partial}{\partial u_{i}^{(m)}}\right)_L\frac{\partial f}{\partial u_j^{(n)}}
=\left(\frac{\partial}{\partial u_j^{(n)}}\right)_R\frac{\partial f}{\partial u_{i}^{(m)}}
\,,
$$
for any noncommutative differential polynomial $f\in \mc R_\ell$, and $i,j\in I$ and
$n,m\in\mb Z_{\geq0}$ (we are using the notation in \eqref{20240805:eq1}).
%
%\begin{definition}\label{def:diff-func}
%An \emph{algebra of (noncommutative) differential functions} in $\ell$ variables
%is a unital associative differential algebra $\mc V$, with derivation $\partial$,
%endowed with strongly commuting $2$-fold derivations
%$\frac{\partial}{\partial u_i^{(n)}}:\,\mc V\to\mc V\otimes\mc V$, $i\in I=\{1,\dots,\ell\},\,n\in\mb Z_{\geq0}$,
%such that \eqref{eq:comm} holds
%and, for every $f\in\mc V$, we have $\frac{\partial f}{\partial u_i^{(n)}}=0$
%for all but finitely many choices of indices $(i,n)\in I\times\mb Z_{\geq0}$.
%\end{definition}
%
%An example of such an algebra is the algebra $\mc R_\ell$,
%endowed with the $2$-fold derivations defined in \eqref{20140627:eq1},
%or its localization by non-zero elements.

The next result characterizes $2$-fold $\lambda$-brackets and dPVA structures on $\mc R_\ell$.
\begin{theorem}[{\cite[Thm.~3.10]{DSKV}}]\label{20130921:prop1}
\begin{enumerate}[(a)]
\item
Any $2$-fold (not necessarily skewsymmetric) $\lambda$-bracket on $\mc R_\ell$ is given by the following
Master Formula ($f,g\in \mc R_\ell$):
\begin{equation}\label{master-infinite}
\ldb f_{\lambda}g\rdb
=\sum_{\substack{i,j\in I\\m,n\in\mb Z_{\geq0}}}
\frac{\partial g}{\partial u_j^{(n)}}
\bullet
(\lambda+x+y)^n
(|_{x=\partial}\ldb u_{i}{}_{\lambda+y}u_{j}\rdb)
(-\lambda-y)^m
\bullet
\left(\Big|_{y=\partial}\frac{\partial f}{\partial u_i^{(m)}}\right)^\sigma\,.
\end{equation}
where $\bullet$ is as in \eqref{20140609:eqc1}.
\item
Equation \eqref{master-infinite} defines a structure of dPVA on $\mc R_\ell$
if and only if the skewsymmetry axiom \eqref{eq:skew2} and the Jacobi identity
\eqref{eq:jacobi2} hold on the $u_i$'s,
namely, for all $i,j,k\in I$, we have
\begin{equation}\label{skew-gen-b}
\ldb {u_i}_\lambda{u_j}\rdb=-(|_{x=\partial}\ldb {u_j}_{-\lambda-x}{u_i}\rdb^\sigma)
\,,
\end{equation}
and
\begin{equation}\label{jacobi-gen-b}
\ldb {u_i}_\lambda \ldb {u_j}_\mu {u_k}\rdb\rdb_L
-
\ldb {u_j}_\mu \ldb {u_i}_\lambda {u_k}\rdb\rdb_R
=
\ldb{\ldb {u_i}_\lambda {u_j}\rdb}_{\lambda+\mu} {u_k}\rdb_L
\,,
\end{equation}
where each term is understood via \eqref{master-infinite}.
%for example
%\begin{equation}\label{20140709:eq3}
%\ldb {u_i}_\lambda \ldb {u_j}_\mu {u_k}\rdb\rdb_L
%=
%\sum_{i\in I,n\in\mb Z_{\geq0}}
%\Big(\Big(\frac{\partial}{\partial u_h^{(n)}}\Big)_LH_{kj}(\mu)\Big)
%\bullet_3 (\lambda+\partial)^n H_{hi}(\lambda)
%\,.
%\end{equation}
\end{enumerate}
\end{theorem}
\begin{remark}
By Theorem \ref{20130921:prop1}(a), $2$-fold $\lambda$-brackets on $\mc R_\ell$ are in one-to-one correspondence with $\ell\times\ell$ matrices $H(\lambda)=(\ldb{u_j}_\lambda{u_i}\rdb)_{i,j\in I}$ with entries in $(\mc R_\ell\otimes\mc R_\ell)[\lambda]$.
On the other hand $H(\lambda)$ is the symbol of the matrix differential operator 
$H(\partial)=\big(H_{ij}(\partial)\big)_{i,j=1}^\ell$,
with entries in $(\mc R_\ell\otimes\mc R_\ell)[\partial]$.
The adjoint operator of $H(\partial)$ is the operator $H^\dagger(\partial)=(H^\dagger_{ij}(\partial))$
where $H^\dagger_{ij}(\partial)=((H_{ji})^*(\partial))^\sigma$ and
$A^*(\partial)$ denotes the formal adjoint of the differential operator $A(\partial)$.
Then property \eqref{skew-gen-b}
means that the matrix differential operator $H(\partial)$ is skewadjoint.
A skewadjoint matrix differential operator 
$H(\partial)$ satisfying \eqref{jacobi-gen-b}
is called a \emph{Poisson structure} on the algebra of noncommutative differential polynomials $\mc R_\ell$, see \cite{DSKV}.
\end{remark}

%%%
\section{Basic cochains and arrays}\label{sec:array}
Let $\mc V=\mc R_\ell$ denote the algebra of noncommutative differential polynomials in $\ell$ variables $u_i$, $i\in I=\{1,\dots,\ell\}$.
In this section we provide an explicit description of the space of basic cochains
$\widetilde{\Gamma}(\mc V)$ introduced in  Section \ref{sec:omega1}.

Recall from Remark \ref{20250722:rem1} that a basic $n$-cochain $X\in\widetilde{\Gamma}^n(\mc V)$ is
uniquely determined by its value on elements $u_{i_1}\otimes\dots\otimes u_{i_n}$, for $i_1,\dots,i_n\in I$.
Indeed, we have the following \emph{Master Formula for basic cochains} which expresses the action of $X$ on $\mc V^{\otimes n}$ in terms of its action on an $n$-tuple of generators.
\begin{proposition}\label{20250723:prop1}
For every $X\in\widetilde{\Gamma}^{n}(\mc V)$, $n\geq1$, and $f_1,\dots,f_n\in\mc V$, we have
\begin{align}
\begin{split}\label{eq:XA}
&X_{\lambda_1,\dots,\lambda_n}(f_1,\dots,f_n)
\\
&=\sum_{\substack{i_1,\dots,i_n\in I\\m_1,\dots,m_n\in\mb Z_{\geq0}}}\Bigg(\dots\Bigg(
X_{\lambda_1+x_1,\dots,\lambda_n+x_n}(u_{i_1},\dots,u_{i_n})
\\
&\quad\quad\quad\quad\quad\quad\quad\quad
\bullet_{(1,2)} (-\lambda_1-x_1)^{m_1}\bigg(\Big|_{x_1=\partial}\frac{\partial f_1}{\partial u_{i_1}^{(m_1)}}\bigg)^{\sigma}\Bigg)
%\bullet_{(2,3)}(-\lambda_2-x_2)^{m_2}\bigg(\Big|_{x_2=\partial}\frac{\partial f_2}{\partial u_{i_2}^{(m_2)}}\bigg)^{\sigma}\Bigg)
\dots 
\Bigg)
\\
&\quad\quad\quad\quad\quad\quad\quad\quad
\bullet_{(n,n+1)}(-\lambda_n-x_n)^{m_n}\bigg(\Big|_{x_n=\partial}\frac{\partial f_n}{\partial u_{i_n}^{(m_n)}}\bigg)^{\sigma}
\,.
\end{split}
\end{align}
\end{proposition}
In \eqref{eq:XA} the notation $\bullet_{(i,i+1)}$,
$i=1,\dots,n$, denotes the right action of $\mc V^{\otimes 2}$ on $\mc V^{\otimes (n+1)}$ defined by equation \eqref{bullet-ii}.
For example we have
$$
X_{\lambda}(f)=\sum_{\substack{i\in I\\m\in\mb Z_{\geq0}}}X_{\lambda+x}(u_i)\bullet (-\lambda-x)^m\bigg(\Big|_{x=\partial}\frac{\partial f}{\partial u_i^{(m)}}\bigg)^{\sigma}
\,,
$$
for $X\in\widetilde{\Gamma}^1(\mc V)$, and
\begin{align*}
&X_{\lambda,\mu}(f,g)
=\sum_{\substack{i,j\in I\\m,n\in\mb Z_{\geq0}}}\Bigg(
X_{\lambda+x,\mu+y}(u_i,u_j)
\bullet_{(1,2)} (-\lambda-x)^m\bigg(\Big|_{x=\partial}\frac{\partial f}{\partial u_i^{(m)}}\bigg )^{\sigma}
\Bigg)
\\
&\quad\quad\quad\quad\quad\quad\quad\quad\quad\quad
\bullet_{(2,3)}(-\mu-y)^n\left(\Big|_{y=\partial}\frac{\partial g}{\partial u_j^{(n)}}\right)^{\sigma}
\,,
\end{align*}
for $X\in\widetilde{\Gamma}^2(\mc V)$, where
\begin{align*}
&(a_1\otimes a_2\otimes a_3)\bullet_{(1,2)}(b_1\otimes b_2)=a_1b_1\otimes b_2a_2\otimes a_3\,,
\\
&(a_1\otimes a_2\otimes a_3)\bullet_{(2,3)}(b_2\otimes b_3)=a_1\otimes a_2b_2\otimes b_3a_3\,,
\end{align*}
for every $a_i,b_i\in\mc V$. It is proved in \cite[Lem.~1.1]{DSKV} that the actions $\bullet_{(1,2)}$ and $\bullet_{(2,3)}$ of $\mc V^{\otimes 2}$ on $\mc V^{\otimes (n+1)}$ commute (although these are denoted with a different notation in the cited paper).
For every $n\geq3$, it can be checked that the right actions $\bullet_{(i,i+1)}$ and $\bullet_{(j,j+1)}$ commute if $i\neq j$, so their order is not important in the Master Formula \eqref{eq:XA}. 
\begin{proof}[Proof of Proposition \ref{20250723:prop1}]
It is immediate to verify, using \eqref{eq:comm}, that the RHS of \eqref{eq:XA} satisfies the sesquilinearity
condition \eqref{eq:sesquimaps} and it is an easy exercise to check that it satisfies the Leibniz rules
\eqref{eq:Leibnizmaps}. Hence, the RHS of \eqref{eq:XA} defines a basic $n$-cochain, which evaluates to $X_{\lambda_1,\dots,\lambda_n}(u_{i_1},\dots,u_{i_n})$ for $f_1=u_{i_1}, \dots, f_n=u_{i_n}$. The claim then follows from the observation in Remark \ref{20250722:rem1}.
\end{proof}
Next, we recall from \cite{DSKV}, the definition of the space of arrays\glslink{Omegarray}{}
$$
\widetilde{\Omega}(\mc V)=\bigoplus_{n\in\mb Z_{\geq0}}\widetilde{\Omega}^n(\mc V)
\,.
$$
For $n=0$ we set $\widetilde{\Omega}^0(\mc V)=\mc V$ and, for $n\geq1$, we set
$\widetilde{\Omega}^n(\mc V)$ to be the space of arrays $A(\lambda_1,\dots,\lambda_n)=(A_{i_1,\dots,i_n}(\lambda_1,\dots,\lambda_{n}))_{i_1,\dots,i_n\in I}$,
with entries $A_{i_1,\dots,i_n}(\lambda_1,\dots,\lambda_{n})\in\mc V^{\otimes (n+1)}[\lambda_1,\dots,\lambda_{n}]$. Note that $\widetilde{\Omega}^n(\mc V)=
(\mc V^{\otimes (n+1)}[\lambda_1,\dots,\lambda_n])^{\ell^n}$.

Given an $n$-cochain $X\in\widetilde{\Gamma}^n(\mc V)$, $n\geq1$, we can define an array
$\widetilde{\gamma}(X)\in\widetilde{\Omega}^n(\mc V)$ by evaluating $X$ on an $n$-tuple of generators
\begin{equation}\label{eq:idsigma}
\widetilde{\gamma}(X)(\lambda_1,\dots,\lambda_n)
=(X_{\lambda_1,\dots,\lambda_n}(u_{i_1},\dots,u_{i_n}))_{i_1,\dots i_n\in I}
\,.
\end{equation}
We also set $\widetilde\gamma(f)=f$, for every $f\in\mc V$. Then, we have a linear map\glslink{gammapTil}{}
\begin{equation}\label{eq:gammatilde}
\widetilde{\gamma}:\widetilde{\Gamma}(\mc V)\rightarrow \widetilde{\Omega}(\mc V)
\,.
\end{equation}
On the other hand, given an array $A(\lambda_1,\dots,\lambda_n)=(A_{i_1,\dots,i_n}(\lambda_1,\dots,\lambda_{n}))_{i_1,\dots,i_n\in I}\in\widetilde{\Omega}^n(\mc V)$,
by the Master Formula \eqref{eq:XA} there is a unique
$X_A\in\widetilde{\Gamma}^n(\mc V)$ such that $\widetilde{\gamma}(X_A)(\lambda_1,\dots,\lambda_n)=A(\lambda_1,\dots,\lambda_n)$.
Hence, $\widetilde{\gamma}$ in \eqref{eq:gammatilde} is an identification.

Under this identification the product \eqref{20230802:eq2} of basic cochains 
translates into the following product of arrays:
for $A(\lambda_1,\dots,\lambda_m)=(A_{i_1,\dots,i_m}(\lambda_1,\dots,\lambda_{m}))_{i_1,\dots,i_m\in I}\in \widetilde{\Omega}^m(\mc V)$ and
$B(\lambda_1,\dots,\lambda_n)=(B_{i_1,\dots,i_n}(\lambda_1,\dots,\lambda_{n}))_{i_1,\dots,i_n\in I}\in \widetilde{\Omega}^n(\mc V)$ their product is the array
$(AB)(\lambda_1,\dots,\lambda_{m+n})\in \widetilde{\Omega}^{m+n}(\mc V)$ with entries
\begin{align*}
&(AB)_{i_1,\dots,i_{m+n}}(\lambda_1,\dots,\lambda_{m+n})
\\
&=A_{i_1,\dots,i_m}(\lambda_1,\dots,\lambda_{m})B_{i_{m+1},\dots,i_{m+n}}(\lambda_{m+1},\dots,\lambda_{m+n})\in\mc V^{\otimes (m+n+1)}[\lambda_1,\dots,\lambda_{m+n}]\,,
\end{align*}
and computed using the associative product \eqref{prod-nm}. Furthermore, the action of the derivation $\partial$ on $\widetilde{\Gamma}(\mc V)$ translates into the following action
on $\widetilde{\Omega}(\mc V)$: for $A(\lambda_1,\dots,\lambda_n)\in\widetilde{\Omega}^n(\mc V)$
we have that $(\partial A)(\lambda_1,\dots,\lambda_n)\in\widetilde{\Omega}^n(\mc V)$ is the array with entries
$$
(\partial A)_{i_1,\dots,i_n}(\lambda_1,\dots,\lambda_n)
=(\lambda_1+\dots+\lambda_n+\partial) \, A_{i_1,\dots,i_n}(\lambda_1,\dots,\lambda_n)
\,,
$$
for every $i_1,\dots,i_n\in I$. Hence, $\widetilde{\Omega}(\mc V)$ is a differential algebra
and we consider it as a superspace by assigning the parity of $A(\lambda_1,\dots,\lambda_n)\in\widetilde{\Omega}^n(\mc V)$ being equal to the parity of $n$.
We can then consider the graded quotient\glslink{Redomegarray}{}
\begin{equation}\label{eq:omega}
\Omega (\mc V)=
\widetilde{\Omega}(\mc V)/(\partial \widetilde{\Omega}(\mc V)+[\widetilde{\Omega}(\mc V),\widetilde{\Omega}(\mc V)])
=\bigoplus_{n\in\mb Z_{\geq0}}\Omega^n(\mc V)
\,.
\end{equation}
We clearly have $\widetilde{\Omega}^0(\mc V)=\mc V_\sharp$. The identification $\widetilde{\gamma}$
in \eqref{eq:gammatilde} descends to an identification
\begin{equation}\label{eq:gamma}
\gamma:\Gamma(\mc V)\stackrel{\simeq}{\rightarrow} \Omega(\mc V)\,,
\quad \gamma([X])=[\widetilde{\gamma}(X)]
\,,
\end{equation}
where $[X]$ denotes the coset of $X\in\widetilde{\Gamma}(\mc V)$ in $\Gamma(\mc V)$
and $[\widetilde{\gamma}(X)]$ denotes the coset of $\widetilde{\gamma}(X)\in\widetilde{\Omega}(\mc V)$
in $\Omega(\mc V)$.

Finally, let us assume that $\mc V$ is a dPVA with $2$-fold $\lambda$-bracket 
$\ldb-_{\lambda}-\rdb$. Using the identification \eqref{eq:gammatilde}
the differential $\tilde \delta$ defined in \eqref{eq:diff} induces a differential (which we still denote by the same symbol) $\tilde \delta:\widetilde{\Omega}^n(\mc V)\to \widetilde{\Omega}^{n+1}(\mc V)$,
$n\in\mb Z_{\geq0}$, defined by\glslink{Tildelta}{}
\begin{equation}
\begin{split}\label{eq:diff-sigma}
&\tilde\delta(A)_{i_1,\dots,i_{n+1}}(\lambda_1,\dots,\lambda_{n+1})
\\
&=
\sum_{s=1}^{n+1}(-1)^{s+1}\ldb {u_{i_s}}_{\lambda_s}
A_{i_1,\stackrel{s}{\check{\dots}},i_{n+1}}(\lambda_1,\stackrel{s}{\check{\dots}},\lambda_{n+1})
\rdb_{(s)}
\\
&
-\sum_{s=1}^{n}(-1)^{s+1}\sum_{\substack{j\in I\\m\in\mb Z_{\geq0}}}
A_{i_1,\dots,i_{s-1},j,i_{s+2},\dots,i_{n+1}}
(\lambda_1,\dots,\lambda_{s-1},\lambda_s+\lambda_{s+1}+x,\lambda_{s+2},\dots
\\
&
\dots,\lambda_{n+1})
\cdot_s(-\lambda_s-\lambda_{s+1}-x)^m\bigg(\Big|_{x=\partial}\bigg(\frac{\partial}{\partial {u_j^{(m)}}}\bigg)_L
\ldb u_{i_s}{}_{\lambda_s}u_{i_{s+1}}\rdb
\bigg)^{\sigma^2}
\,,
\end{split}
\end{equation}
for every $A(\lambda_1,\dots,\lambda_{n})=(A_{i_1,\dots,i_n}(\lambda_1,\dots,\lambda_n))\in\widetilde{\Omega}^n(\mc V)$. In \eqref{eq:diff-sigma} we are using the abbreviation
\begin{equation}\label{notationneverends}
\begin{split}
&(f_1\otimes\dots \otimes f_{n+1})\cdot_s(a\otimes b\otimes c)
\\
&=(f_1\otimes\dots \otimes f_s a\otimes b\otimes c  f_{s+1}  \otimes \dots\otimes f_{n+1})\
\\
&=\left((f_1\otimes\dots \otimes f_{n+1})\bullet_{(s,s+1)}(a\otimes c)\right)\otimes_{n+1-s}b
\,,
\end{split}
\end{equation}
where $f_1,\dots,f_{n+1},a,b,c\in\mc V$.

In conclusion by the identification \eqref{eq:gammatilde} (respectively \eqref{eq:gamma}) we have the complex $(\widetilde{\Omega}(\mc V),\tilde\delta)$ (respectively
$(\Omega(\mc V),\delta)$\glslink{deltadiff}{}) whose cohomology is isomorphic to the basic (respectively reduced)
dPVA cohomology of $\mc V$:
\begin{equation}\label{20250723:eq2}
\dH_{\mathrm{bas}}(\mc V)\simeq \coH(\widetilde{\Omega}(\mc V),\tilde\delta)
\quad \text{(respectively }\dH_{\textrm{red}}(\mc V)\simeq \coH(\Omega(\mc V),\delta)\text{)}
\,.
\end{equation}
\begin{remark}\label{rem:beltrami1}
In analogy with \cite{BDSK,DSKV}, we introduce the \emph{Beltrami} $2$\emph{-fold} $\lambda$\emph{-bracket}
$\ldb-_{\lambda}-\rdb^B$ using equation \eqref{master-infinite} where we let
$\ldb u_i{}_\lambda u_j\rdb^B=\delta_{ij}(1\otimes1)$, for $i,j\in I$ (note that
this $2$-fold $\lambda$-bracket is commutative and does not define a dPVA structure on $\mc V$).
By the Master Formula \eqref{master-infinite} we have ($f\in\mc V, j\in I$)
\begin{equation}\label{eq:beltrami1}
\ldb f_\lambda u_j \rdb^B
=\sum_{m\in\mb Z_{\geq0}}(-\lambda-\partial)^m\bigg(\frac{\partial f}{\partial u_j^{(m)}}\bigg)^\sigma
\,.
\end{equation}
Using \eqref{eq:beltrami1}, the identity $(a\otimes b)^\sigma\otimes_1c=\sigma^2(a\otimes b\otimes c)$, which holds for every $a,b,c\in\mc V$, and the definition
of $\ldb f\otimes g_\lambda u_j\rdb^B_L$ given in \eqref{notation} we get
\begin{equation}\label{eq:beltrami2}
\ldb f\otimes g_\lambda u_j\rdb^B_L=\sum_{m\in\mb Z_{\geq0}}
(-\lambda-\partial)^m
\bigg(\bigg(\frac{\partial}{\partial {u_j^{(m)}}}\bigg)_L\left(f\otimes g\right)\bigg)^{\sigma^2}
\,,
\end{equation}
for every $f,g\in\mc V$ and $j\in I$. Then, using \eqref{eq:beltrami2}, equation \eqref{eq:diff-sigma}
can be rewritten as
\begin{equation}
\begin{split}\label{eq:diff-sigma-beltrami}
&\tilde\delta(A)_{i_1,\dots,i_{n+1}}(\lambda_1,\dots,\lambda_{n+1})
=
\sum_{s=1}^{n+1}(-1)^{s+1}\ldb {u_{i_s}}_{\lambda_s}
A_{i_1,\stackrel{s}{\check{\dots}},i_{n+1}}(\lambda_1,\stackrel{s}{\check{\dots}},\lambda_{n+1})
\rdb_{(s)}
\\
&
-\sum_{s=1}^{n}(-1)^{s+1}\sum_{j\in I}
A_{i_1,\dots,i_{s-1},j,i_{s+2},\dots,i_{n+1}}
(\lambda_1,\dots,\lambda_s+\lambda_{s+1}+x,\dots
\\
&\quad\quad\quad\quad\quad\quad\quad\quad
\dots,\lambda_{n+1})
\cdot_s\left(|_{x=\partial}\ldb \ldb u_{i_s}{}_{\lambda_s}u_{i_{s+1}}\rdb_{\lambda_s+\lambda_{s+1}}u_j\rdb^B_L\right)
\,.
\end{split}
\end{equation}
\end{remark}
\begin{remark}\label{rem:deltaHconst}
Let $\ldb -_\lambda-\rdb$ be a (not necessarily skewsymmetric) 2-fold $\lambda$-bracket on $\mc V$ such that on generators we have
($i,j\in I$)
\begin{equation}\label{ass-const}
\ldb {u_i}{}_\lambda u_j\rdb\in\kk[\lambda]\subset\mc (\mc V\otimes\mc V)[\lambda]
\,.
\end{equation}
In this case, the Jacobi identity \eqref{jacobi-gen-b} holds
since all three terms are zero, but the skewsymmetry condition
\eqref{skew-gen-b} may not hold. Thus $\mc V$ is not necessarily a
dPVA.
Moreover, due to the fact that partial derivatives strongly commute, the following identity holds ($i,j\in I$, $f\in\mc V$)
\begin{equation}\label{20250801:eq1}
\ldb {u_i}_\lambda \ldb {u_j}_\mu f\rdb\rdb_L
=\ldb {u_j}_\mu \ldb {u_i}_\lambda f\rdb\rdb_R
\,.
\end{equation}
Let us consider a linear map
$\tilde \delta:\widetilde{\Omega}^n(\mc V)\rightarrow\widetilde{\Omega}^{n+1}(\mc V)$ given only by the first sum in \eqref{eq:diff-sigma} 
\begin{equation}\label{eq:dP_const}
\tilde\delta(A)_{i_1,\dots\,i_{n+1}}(\lambda_1,\dots,\lambda_{n+1})
=\sum_{s=1}^{n+1}(-1)^{s+1}
\ldb u_{i_s}{}_{\lambda_s}
A_{i_1,\stackrel{s}{\check{\dots}},i_{n+1}}(\lambda_1,\stackrel{s}{\check{\dots}},\lambda_{n+1})
\rdb_{(s)}
\,,
\end{equation}
for every $A=(A_{i_1,\dots,i_n}(\lambda_1,\dots,\lambda_{n}))_{i_1,\dots,i_n\in I}\in\widetilde{\Omega}^n(\mc V)$. It follows from the proof of Theorem
\ref{thm:PVAcoh} that $\tilde\delta$ gives a well defined differential on 
$\widetilde{\Omega}(\mc V)=\bigoplus_{n\geq0}\widetilde{\Omega}^n(\mc V)$, for every choice of $\ldb u_i{}_\lambda u_j\rdb$ as in 
\eqref{ass-const}. Indeed, in this case we only need to check the vanishing of the analogue of \eqref{eq:dPsquareA1}+\eqref{eq:dPsquareA2}, which is
equivalent to show that 
$$
\sum_{t=1}^{n+1}
   \ldb {u_{i_t}}{}_{\lambda_t} \ldb {u_{i_{t+1} }}{}_{\lambda_{t+1}} C \rdb_{(t)} \rdb_{(t)} 
=\sum_{t=1}^{n+1}\ldb {u_{i_{t+1} }}{}_{\lambda_{t+1}}  \ldb {u_{i_t}}{}_{\lambda_t} C \rdb_{(t)} \rdb_{(t+1)}
\,,
$$
where $C=A_{i_1,\ldots,i_{t-1},i_{t+2},\ldots,i_n}(\lambda_1,\ldots,\lambda_{t-1},\lambda_{t+2},\ldots,\lambda_n)$. In fact, the above identity follows from \eqref{20250801:eq1}.
\end{remark}

%%%
\section{\texorpdfstring{\MakeLowercase{$n$}}{n}-fold \texorpdfstring{$\lambda$}{lambda}-brackets and skewsymmetric arrays}\label{sec:skewarray}
In this section we give an explicit description of $n$-fold $\lambda$-brackets on $\mc V=\mc R_\ell$ by providing a generalization of the Master Formula \eqref{master-infinite}.

Recall from Remark \ref{20250722:rem2} that an $n$-fold $\lambda$-bracket $\ldb-_{\lambda_1}-\dots-_{\lambda_{n-1}}-\rdb\in C^n(\mc V)$ is uniquely determined by its value on elements $u_{i_1}\otimes \dots \otimes u_{i_n}\in\mc V^{\otimes n}$.
Using this fact, it was shown in \cite{DSKV} that any $D\in C^1(\mc V)=\Vect(\mc V)^\partial$ has the form ($f\in\mc V$)
\begin{equation}\label{eq:QS-1}
D(f)=\sum_{\substack{i\in I\\m\in\mb Z_{\geq0}}}\mult\bigg(\partial^m(D(u_i))\star \frac{\partial f}{\partial u_{i}^{(m)}}\bigg)
=\sum_{\substack{i\in I\\m\in\mb Z_{\geq0}}}
\frac{\partial f}{\partial u_{i}^{(m)}}\bullet
\partial^m(D(u_i))
\,,
\end{equation}
where in the second equality we used the left action \eqref{eq:bullyaction} of $\mc V^{\otimes 2}$ on
$\mc V$. The generalization of equations \eqref{eq:QS-1} and \eqref{master-infinite} is given by the following \emph{Master Formula for $n$-fold $\lambda$-brackets}
which expresses the action of $\ldb-_{\lambda_1}-\dots-_{\lambda_{n-1}}-\rdb$ on $\mc V^{\otimes n}$ in terms of its action on an $n$-tuple of generators.
\begin{proposition}\label{20250723:prop2}
For every $\ldb-_{\lambda_1}-\dots-_{\lambda_{n-1}}-\rdb\in C^{n}(\mc V)$, $n\geq1$, and $f_1,\dots,f_n\in\mc V$, we have
\begin{align}
\begin{split}\label{eq:QS}
&\ldb f_1{}_{\lambda_1}\dots f_{n-1}{}_{\lambda_{n-1}}f_n\rdb
=\sum_{\substack{i_1,\dots,i_n\in I\\m_1,\dots,m_n\in\mb Z_{\geq0}}}
\frac{\partial f_n}{\partial u_{i_n}^{(m_n)}}
(\lambda_1+\dots+\lambda_{n-1}+\partial)^{m_n}
\\
&
\bullet\Bigg(\Bigg(\dots\Bigg(
\ldb u_{i_1}{}_{\lambda_1+x_1}\dots u_{i_{n-1}}{}_{\lambda_{n-1}+x_{n-1}}u_{i_n}\rdb\bullet_{(1,2)} (-\lambda_1-x_1)^{m_1}\bigg(\Big|_{x_1=\partial}\frac{\partial f_1}{\partial u_{i_1}^{(m_1)}}\bigg)^{\sigma}\Bigg)
\\
&%\left.
%\bullet_{(2,3)}(-\lambda_2-x_2)^{m_2}\left(\Bigg|_{x_2=\partial}\frac{\partial f_2}{\partial u_{i_2}^{(m_2)}}\right)^{\sigma}\right)
\dots 
\Bigg)\bullet_{(n-1,n)}(-\lambda_{n-1}-x_{n-1})^{m_{n-1}}\bigg(\Big|_{x_{n-1}=\partial}\frac{\partial f_{n-1}}{\partial u_{i_{n-1}}^{(m_{n-1})}}\bigg)^{\sigma}\Bigg)
\,.
\end{split}
\end{align}
\end{proposition}
In \eqref{eq:QS} we are using the left (respectively right) action of $\mc V\otimes\mc V$ of
$\mc V^{\otimes n}$ given by \eqref{eq:bullyaction} (respectively \eqref{bullet-ii}). Recall from Section \ref{sec:1.1} that these actions commute, so their order is not important in the Master Formula \eqref{eq:QS}.
\begin{proof}[Proof of Proposition \ref{20250723:prop2}]
Analogous to the proof of Proposition \ref{20250723:prop2}, we briefly outline it for completeness. It is immediate to verify, using \eqref{eq:comm}, that the RHS of \eqref{eq:QS} satisfies the sesquilinearity
conditions \eqref{20140702:eq4} and \eqref{20140702:eq5}, and we leave as an exercise to verify that it satisfies skewsymmetry \eqref{eq:nfold-skew} and the Leibniz rules
\eqref{20140702:eq6}. Hence, the RHS of \eqref{eq:QS} defines an $n$-fold $\lambda$-bracket, which evaluates to $\ldb u_{i_1}{}_{\lambda_1}\dots u_{i_{n-1}}{}_{\lambda_{n-1}}u_{i_n}\rdb$ for $f_1=u_{i_1}, \dots, f_n=u_{i_n}$. The claim then follows from Remark \ref{20250722:rem2}.
\end{proof}
We recall the definition of the space of skewsymmetric arrays\glslink{Sigmarray}{}
$$
\Sigma(\mc V)=\bigoplus_{n\in\mb Z_{\geq0}}\Sigma^n(\mc V)
$$
from \cite{DSKV}. For $n=0$ we set
$\Sigma^0(\mc V)=\mc V_\sharp$.
For $n\geq1$,
we define $\Sigma^n(\mc V)$ as the space of arrays
$S(\lambda_1,\dots,\lambda_{n-1})=\big(S_{i_1\dots i_n}(\lambda_1,\dots,\lambda_{n-1})\big)_{i_1,\dots,i_n\in I}$
with entries 
$S_{i_1\dots i_n}(\lambda_1,\dots,\lambda_{n-1})\in \mc V^{\otimes n}[\lambda_1,\dots,\lambda_{n-1}]$,
satisfying the following skewadjointness condition ($i_1,\dots,i_n\in I$):
\begin{equation}\label{20140626:eq1-b}
S_{i_1\dots i_n}(\lambda_1,\dots,\lambda_{n-1})
=
(-1)^{n+1}\big(|_{x=\partial}
S_{i_2\dots i_n i_1}(\lambda_2,\dots,\lambda_{n-1},
-\lambda_1-\dots-\lambda_{n-1}-x)
\big)^\sigma
\,,
\end{equation}
where, as usual, $\sigma$ denotes the action of the cyclic permutation on $\mc V^{\otimes n}$.

Clearly, given an $n$-fold $\lambda$-bracket $Q=\ldb-_{\lambda_1}-\dots-_{\lambda_{n-1}-}\rdb\in C^n(\mc V)$, $n\geq1$, we can define a skewsymmetric array
$\psi(Q)(\lambda_1,\dots,\lambda_{n-1})\in\Sigma^n(\mc V)$ by evaluating $Q$ on an $n$-tuple of generators
\begin{equation}\label{eq:idsigma2}
\psi(Q)
=(\ldb u_{i_1}{}_{\lambda_1}u_{i_2}\dots {u_{i_{n-1}}}_{\lambda_{n-1}} u_{i_n}\rdb)_{i_1,\dots i_n\in I}
\,.
\end{equation}
We also set $\psi(\tint f)=\tint f$, for every $\tint f\in \mc V_\sharp$.
Then, we have a linear map\glslink{psimap}{}
\begin{equation}\label{eq:psi}
\psi:C(\mc V)\rightarrow \Sigma(\mc V)
\,.
\end{equation}
On the other hand, let 
$S(\lambda_1,\dots,\lambda_{n-1})=(S_{i_1,\dots,i_n}(\lambda_1,\dots,\lambda_{n-1}))_{i_1,\dots,i_n\in I}$
be a skewsymmetric array.
By the Master Formula \eqref{eq:QS} there is a unique $Q_S\in C^n(\mc V)$ such that
$\psi(Q_S)(\lambda_1,\dots,\lambda_{n-1})=S(\lambda_1,\dots,\lambda_{n-1})$.
Hence, the map $\psi$ in \eqref{eq:psi} is an identification.

Let us assume that $\mc V$ is a dPVA with $2$-fold $\lambda$-bracket
$\ldb-_{\lambda}-\rdb$. Using the identification $C(\mc V)\simeq\Sigma(\mc V)$ given by \eqref{eq:psi},
the differential $\dd$ defined by \eqref{eq:dP0} and \eqref{eq:dP-1} induces a differential (which we still denote by the same symbol) $\dd:\Sigma^n(\mc V)\to \Sigma^{n+1}(\mc V)$,
$n\in\mb Z_{\geq0}$, defined by
\begin{equation}\label{dH-tint}
\dd(\tint f)
%=\left(
%\mult\ldb {u_{i}}_{-\partial} f\rdb^\sigma
%\right)_{i\in I}
=\left(
-\mult\ldb f{}_{\lambda} u_i\rdb|_{\lambda=0}
\right)_{i\in I}
\in\mc V^\ell=\Sigma^1(\mc V)
\,,
\end{equation}
for $\tint f\in\mc V_\sharp=\Sigma^0(\mc V)$, and by\glslink{dd}{}
\begin{equation}
\begin{split}\label{eq:diff-sigma-1}
&\dd(S)_{i_1,\dots,i_{n+1}}(\lambda_1,\dots,\lambda_{n})
=
\sum_{s=1}^{n}(-1)^{n+s+1}\ldb {u_{i_s}}_{\lambda_s}
S_{i_1,\stackrel{s}{\check{\dots}},i_{n+1}}(\lambda_1,\stackrel{s}{\check{\dots}},\lambda_{n})
\rdb_{(s)}
\\
&-\ldb S_{i_1,\dots,i_n}(\lambda_1,\dots,\lambda_{n-1})_{\lambda_1+\dots+\lambda_n}u_{i_{n+1}}\rdb_{L}
\\
&
+\sum_{s=1}^{n}(-1)^{n+s}\sum_{\substack{j\in I\\m\in\mb Z_{\geq0}}}
S_{i_1,\dots,i_{s-1},j,i_{s+2},\dots,i_{n+1}}
(\lambda_1,\dots,\lambda_{s-1},\lambda_s+\lambda_{s+1}+x,\lambda_{s+2},\dots
\\
&\quad\quad\quad\quad\quad\quad
\dots,\lambda_{n})
\cdot_s(-\lambda_s-\lambda_{s+1}-x)^m\bigg(\Big|_{x=\partial}\bigg(\frac{\partial}{\partial {u_j^{(m)}}}\bigg)_L
\ldb u_{i_s}{}_{\lambda_s}u_{i_{s+1}}\rdb\bigg)^{\sigma^2}
\\
&+(-1)^{n+1}\sum_{\substack{j\in I\\m\in\mb Z_{\geq0}}}
\bigg(\bigg(\frac{\partial }{\partial u_{j}^{(m)}}\bigg)_R
\ldb u_{i_1}{}_{\lambda_1}u_{i_{n+1}}\rdb\bigg)
(\Id\otimes \bullet)(\lambda_2+\dots
\\
&\quad\quad\quad\quad\quad\quad
\dots+\lambda_{n}+\partial)^mS_{i_2,\dots,i_n,j}(\lambda_2,\dots,\lambda_n)
\,,
\end{split}
\end{equation}
for every $S(\lambda_1,\dots,\lambda_{n-1})\in\Sigma^n(\mc V)$. In equation \eqref{eq:diff-sigma-1} we are using the notation \eqref{notationneverends} and we denote
$$
(a\otimes b\otimes c)(\Id \otimes \bullet)A=
a\otimes\left( (b\otimes c)\bullet A\right)
\,,
$$
for every $a,b,c\in\mc V$ and $A\in\mc V^{\otimes n}$.

In conclusion by the identification \eqref{eq:psi} we have the complex $(\Sigma(\mc V),\dd)$
whose cohomology is isomorphic to the variational
dPVA cohomology of $\mc V$:
\begin{equation}\label{20250723:eq3}
\dPVH(\mc V)\simeq \coH(\Sigma(\mc V),\dd)
\,.
\end{equation}
\begin{remark}\label{rem:beltrami2}
Recall the Beltrami $\lambda$-bracket $\ldb-_\lambda-\rdb^B$ introduced in Remark \ref{rem:beltrami1}. By the Master-Formula \eqref{master-infinite} we have
\begin{equation}\label{beltramiright}
\ldb u_j{}_\lambda f\rdb^B=\sum_{m\in\mb Z_{\geq0}}\frac{\partial f}{\partial u_j^{(m)}}\lambda^m\,,
\end{equation}
for every $f\in\mc V$ and $i\in I$. Hence, using \eqref{notation}, we also have
$$
\ldb u_j{}_\lambda f\otimes g\rdb^B_R
=\sum_{m\in\mb Z_{\geq0}}\bigg(f\otimes \frac{\partial g}{\partial u_j^{(m)}}\bigg)\lambda^m
=\sum_{m\in\mb Z_{\geq0}}\bigg(\frac{\partial }{\partial u_j^{(m)}}\bigg)_R(f\otimes g)\lambda^m\,,
$$
for every $f,g\in\mc V$ and $i\in I$. Hence (cf. Remark \ref{rem:beltrami1}), we can rewrite   \eqref{eq:diff-sigma-1} as
\begin{equation}
\begin{split}\label{eq:diff-sigma-1-beltrami}
&\dd(S)_{i_1,\dots,i_{n+1}}(\lambda_1,\dots,\lambda_{n})=
\sum_{s=1}^{n}(-1)^{n+s+1}\ldb {u_{i_s}}_{\lambda_s}
S_{i_1,\stackrel{s}{\check{\dots}},i_{n+1}}(\lambda_1,\stackrel{s}{\check{\dots}},\lambda_{n})
\rdb_{(s)}
\\
&-\ldb S_{i_1,\dots,i_n}(\lambda_1,\dots,\lambda_{n-1})_{\lambda_1+\dots+\lambda_n}u_{i_{n+1}}\rdb_{L}
\\
&
+\sum_{s=1}^{n}(-1)^{n+s}\sum_{j\in I}
S_{i_1,\dots,i_{s-1},j,i_{s+2},\dots,i_{n+1}}
(\lambda_1,\dots,\lambda_s+\lambda_{s+1}+x,\dots
\\
&
\quad\quad\quad\quad\quad\quad\quad\quad
\dots,\lambda_{n})
\cdot_s\big(|_{x=\partial}\ldb 
\ldb u_{i_s}{}_{\lambda_s} u_{i_{s+1}}\rdb_{\lambda_s+\lambda_{s+1}}u_j\rdb^B_L \big)
\\
&+(-1)^{n+1}\sum_{j\in I}
\ldb u_j{}_{\lambda_2+\dots+\lambda_{n}+x}
\ldb u_{i_1}{}_{\lambda_1} u_{i_{n+1}}\rdb\rdb^B_R
(\Id\otimes \bullet)\big(|_{x=\partial}S_{i_2,\dots,i_n,j}(\lambda_2,\dots,\lambda_n)\big)
\,.
\end{split}
\end{equation}
\end{remark}
\begin{remark}\label{rem:deltaHconst2}
Let $\ldb -_\lambda-\rdb^H$ be a 2-fold $\lambda$-bracket on $\mc V$ which is constant on generators, that is \eqref{ass-const} holds.
Then, the same argument as in Remark \ref{rem:deltaHconst}, shows that the map $\dd$ defined by $\dd(\tint f)$ as in \eqref{dH-tint} for any $\tint f\in\Sigma^0(\mc V)$, and by 
\begin{equation}
\begin{split}\label{eq:diff-sigma-1-const}
&\dd(S)_{i_1,\dots,i_{n+1}}(\lambda_1,\dots,\lambda_{n})
=
\sum_{s=1}^{n}(-1)^{n+s+1}\ldb {u_{i_s}}_{\lambda_s}
S_{i_1,\stackrel{s}{\check{\dots}},i_{n+1}}(\lambda_1,\stackrel{s}{\check{\dots}},\lambda_{n})
\rdb_{(s)}
\\
&+(|_{x=\partial}\ldb u_{i_{n+1}}{}_{-\lambda_1-\dots-\lambda_n-x}S_{i_1,\dots,i_n}(\lambda_1,\dots,\lambda_{n-1})\rdb_{(n+1)})
\,,
\end{split}
\end{equation}
for every $S\in\Sigma^n(\mc V)$, is a well-defined differential on $\Sigma(\mc V)$.
\end{remark}

%%%
\section[Relation between dPVA cohomologies revisited]{Relation between reduced and variational double Poisson vertex algebra cohomologies revisited}\label{sec:rel}

In this section we show that $\Gamma(\mc V)\simeq C(\mc V)$ when $\mc V=\mc R_\ell$ is the algebra of noncommutative 
differential polynomials in $\ell$ variables $u_i$, $i\in I=\{1,\dots,\ell\}$.

For $n\geq1$, we define a linear map\glslink{phimap}{}
\begin{equation}\label{eq:Phi_n}
\widetilde{\phi}_n:\widetilde{\Omega}^n(\mc V)\rightarrow\Sigma^n(\mc V)
\end{equation}
as follows. Let $i_1,\dots,i_n\in I$, $m_1,\dots,m_n\in\mb Z_{\geq 0}$ and $f_1,\dots f_{n+1}\in\mc V$
be fixed and let
$$
A=A(\lambda_1,\dots,\lambda_n)=\left(\delta_{j_1,i_1}\dots\delta_{j_n,i_n}f_1\otimes\dots\otimes f_{n+1} \lambda_1^{m_1}\dots\lambda_n^{m_n}\right)_{j_1,\dots,j_n\in I}\in\widetilde{\Omega}^n(\mc V)
\,.
$$
We map it to the skewsymmetric array $\widetilde{\phi}_n(A)=\left(S_{j_1,\dots,j_n}(\lambda_1,\dots,\lambda_{n-1})\right)_{j_1,\dots,j_n\in I}\in\Sigma^n(\mc V)$
with entries 
$S_{j_1\dots j_n}(\lambda_1,\dots,\lambda_{n-1})=0$
unless $(j_1,\dots,j_n)$ is a cyclic permutation of $(i_1,\dots,i_n)$,
and
\begin{equation}\label{20140626:eq2-aff}
\begin{array}{l}
\displaystyle{
\vphantom{\Big)}
S_{j_1\dots j_n}(\lambda_1,\dots,\lambda_{n-1})
=
\frac1n (-1)^{s(n-s)}
\lambda_1^{m_{s+1}}\dots\lambda_{n-s}^{m_n}
\lambda_{n-s+1}^{m_1}\dots\lambda_{n-1}^{m_{s-1}}
} \\
\displaystyle{
\vphantom{\Big)}
(-\lambda_1-\dots-\lambda_{n-1}-\partial)^{m_s}
\big(
f_{s+1}\otimes\dots\otimes f_n\otimes f_{n+1}f_1\otimes f_2\otimes\dots\otimes f_s
\big)
\,,}
\end{array}
\end{equation}
for $(j_1,\dots,j_n)=(i_{\sigma^s(1)},\dots,i_{\sigma^s(n)})$, where $s=0,\dots,n-1$.
(For $s=0$, the tensor product in the RHS of \eqref{20140626:eq2-aff} is $f_{n+1} f_1\otimes \dots \otimes f_n$.)
We extend $\widetilde{\phi}_n$ by linearity to any element of $\widetilde{\Omega}(\mc V)$. Let us recall the following result.
\begin{proposition}[\cite{DSKV}]\label{20240819:prop1}
For every $n\geq0$ we have an isomorphism $\Omega^n(\mc V)\simeq\Sigma^n(\mc V)$.
\end{proposition}
\begin{proof}
By definition $\Omega^0(\mc V)=\Sigma^0(\mc V)=\mc V_\sharp$. For $n\geq1$, let
$\phi_n:\Omega^n(\mc V)\rightarrow\Sigma^n(\mc V)$ be the linear map defined by
$\phi_n\left([A(\lambda_1,\dots,\lambda_n)]\right)
=\widetilde{\phi}_n\left(A(\lambda_1,\dots,\lambda_n)\right)$,
for every $A(\lambda_1,\dots,\lambda_n)\in\widetilde{\Omega}^n(\mc V)$. This map is well defined 
(it does not depend on the choice of the representative in the coset $[A(\lambda_1,\dots,\lambda_n)]$) and the inverse map $\Sigma^n(\mc V)\rightarrow\Omega^n(\mc V)$ is given by
\begin{equation}\label{20140626:eq3-aff}
\begin{split}
&\left(S_{i_1,\dots,i_n}(\lambda_1,\dots,\lambda_{n-1})\right)_{i_1,\dots,i_n\in I}\in\Sigma^n(\mc V)
\\
&\mapsto
[\left(S_{i_1,\dots,i_n}(\lambda_1,\dots,\lambda_{n-1})\otimes 1\right)_{i_1,\dots,i_n\in I}]\in \Omega^n(\mc V)
\,.
\end{split}
\end{equation}
\end{proof}
Let us denote by $\phi:\Omega(\mc V)\rightarrow\Sigma(\mc V)$ the isomorphism of Proposition \ref{20240819:prop1}.
We have shown in Section \ref{sec:array} that we have an identification
$\gamma:\Gamma(\mc V)\stackrel{\simeq}{\rightarrow}\Omega(\mc V)$, cf.~\eqref{eq:gamma}, and in Section \ref{sec:skewarray}
that we have an identification $\psi:C(\mc V)\stackrel{\simeq}{\rightarrow}\Sigma(\mc V)$, cf.~\eqref{eq:psi}.
Hence, we have the following corollary of Proposition \ref{20240819:prop1}.
\begin{corollary}\label{20240819:cor1}
Let $\mc V=\mc R_\ell$ be the algebra of noncommutative differential polynomials in $\ell$ variables $u_i$. Then, we have an identification
\begin{equation}\label{eq:gammaBRA}
\Gamma(\mc V)=\widetilde{\Gamma}(\mc V)/(\partial\widetilde{\Gamma}(\mc V)+[\widetilde{\Gamma}(\mc V),\widetilde{\Gamma}(\mc V)])\simeq C(\mc V)
\,,
\end{equation}
given by the map $\psi^{-1}\circ\phi\circ\gamma$.
\end{corollary}
In particular, we have an isomorphism
$\Gamma^n(\mc V)\simeq C^n(\mc V)$, $n\in\mb Z_{\geq0}$, on each graded component
of $\Gamma(\mc V)$ and $C(\mc V)$.
\begin{proposition}\label{20240819:prop2}
The projection map $\proj:\Gamma(\mc V)\to C(\mc V)$, given on each graded component by \eqref{eq:projmaps}, is an isomorphism.
\end{proposition}
\begin{proof}
By definition $\Gamma^0(\mc V)=C^0(\mc V)=\mc V_\sharp$ and $\proj_0$ is the identity map.
Let $n\geq1$ and $X\in\widetilde{\Gamma}^n(\mc V)$. By definition of the projection
map $\proj_n$ (see Corollary \ref{20240819:cor2}) and by Corollary \ref{20240819:cor1} it suffices to show that
\begin{equation}\label{eq:toprove}
\widetilde{\proj}_n(X)=(\psi_n^{-1}\circ \widetilde{\phi}_n\circ \widetilde\gamma_n)(X)
\,.
\end{equation}
Let $i_1,\dots,i_n\in I$, $m_1,\dots,m_n\in\mb Z_{\geq 0}$ and $f_1,\dots f_{n+1}\in\mc V$
be fixed. By linearity, it suffices to check the identity \eqref{eq:toprove} for the unique
$n$-cochain $X$ such that
$$
\widetilde{\gamma}_n(X)=
\left(\delta_{j_1,i_1}\dots\delta_{j_n,i_n}f_1\otimes\dots\otimes f_{n+1} \lambda_1^{m_1}\dots\lambda_n^{m_n}\right)_{j_1,\dots,j_n\in I}
\,.
$$
Then, by the definitions of $\widetilde{\phi}_n$ given by \eqref{20140626:eq2-aff} and of $\psi_n^{-1}$ 
given in Section \ref{sec:skewarray} we have that 
$(\psi_n^{-1}\circ \widetilde{\phi}_n\circ \widetilde\gamma_n)(X)=\ldb-_{\lambda_1}-\dots-_{\lambda_{n-1}}-\rdb
\in\Sigma^n(\mc V)$ is the unique $n$-fold $\lambda$-bracket such that
$\ldb u_{j_1}{}_{\lambda_1}\dots u_{j_{n-1}}{}_{\lambda_{n-1}}u_{j_n}\rdb$ is given by the RHS of
\eqref{20140626:eq2-aff} for $(j_1,\dots,j_n)=(i_{\sigma^s(1)},\dots,i_{\sigma^s(n)})$, $s=0,\dots,n-1$,
and it is zero otherwise. It is immediate to check, using \eqref{eq:Pn}, that
$$
\widetilde{\proj}_n(X)_{\lambda_1,\dots,\lambda_{n-1}}(u_{j_1},\dots, u_{j_n})=\ldb u_{j_1}{}_{\lambda_1}\dots u_{j_{n-1}}{}_{\lambda_{n-1}}u_{j_n}\rdb
\,,
$$
for every $j_1,\dots,j_n\in I$. Hence, the identity \eqref{eq:toprove}
follows by the uniqueness of $\ldb-_{\lambda_1}-\dots-_{\lambda_{n-1}}-\rdb$ (cf. Remark \ref{20250722:rem2}).
\end{proof}
%
%\begin{remark}\label{diffconst}
By equation \eqref{delta-comm} and Proposition \ref{20240819:prop2} we have that the differentials
$\delta$ and $\dd$ on $\Gamma^n(\mc V)\simeq C^n(\mc V)$ differ by a constant factor:
$$
\delta=\frac{n+\delta_{n,0}}{n+1}(-1)^n \dd
\,,
\qquad n\in\mb Z_{\geq0}
\,.
$$
Hence, we have an isomorphism
\begin{equation}\label{20250723:eq4}
\dH_{\text{red}}^n(\mc V)\simeq \dPVH^n(\mc V)\,,
\quad n\in\mb Z_{\geq0}
\,.
\end{equation}
%\end{remark}

\medskip

We conclude this section stating a result that will be used in Chapter \ref{CH:PVAexa}.
Let us define the space\glslink{Sigmabar}{}
$$
\overline{\Sigma}(\mc V)=\bigoplus_{n\in\mb Z_{\geq0}}\overline{\Sigma}^n(\mc V)
\,.
$$
where we set $\overline{\Sigma}^0(\mc V)=\mc V/[\mc V,\mc V]$, and, for $n\geq1$, we define $\overline{\Sigma}^n(\mc V)$ as the space of arrays 
$\overline{S}(\lambda_1,\dots,\lambda_{n})=\big(\overline S_{i_1\dots i_n}(\lambda_1,\dots,\lambda_{n})\big)_{i_1,\dots,i_n\in I}$
with entries 
$\overline{S}_{i_1\dots i_n}(\lambda_1,\dots,\lambda_{n})\in \mc V^{\otimes n}[\lambda_1,\dots,\lambda_{n}]$,
satisfying the following skewadjointness condition ($i_1,\dots,i_n\in I$):
$$
\overline{S}_{i_1\dots i_n}(\lambda_1,\dots,\lambda_{n})
=
(-1)^{n+1}
\overline{S}_{i_2\dots i_n i_1}(\lambda_2,\dots,\lambda_{n},
\lambda_1)^\sigma
\,.
$$
Since the commutator space $[\mc V,\mc V]\subset \mc V$ is preserved by the derivation $\partial$ we have an induced action of $\partial$ on $\overline{\Sigma}^0(\mc V)$. We extend this action to an action on $\overline{\Sigma}^n(\mc V)$, $n\geq1$, by letting
$\partial \overline{S}(\lambda_1,\dots,\lambda_n)$ be the array with entries 
$$
(\partial \overline{S})_{i_1,\dots,i_n}(\lambda_1,\dots,\lambda_n)=(\lambda_1+\dots+\lambda_n+\partial)\overline{S}_{i_1,\dots,i_n}(\lambda_1,\dots,\lambda_n)
\,.
$$
Let $\ev|_{\lambda_n=\lambda_n^\dagger}:\overline{\Sigma}^n(\mc V)\to\Sigma(\mc V)$ be the evaluation map obtained by replacing $\lambda_n$ with $\lambda_n^\dagger$ as in \eqref{eq:dagger}, that is $\ev|_{\lambda_n=\lambda_n^\dagger}(\overline{S})=S\in\Sigma^n(\mc V)$ is the array with entries $ S_{i_1,\dots,i_n}(\lambda_1,\dots,\lambda_{n-1})=(|_{x=\partial}\overline S(\lambda_1,\dots,\lambda_{n-1},-\lambda_1-\dots-\lambda_{n-1}-x))$.
This map is surjective and its kernel is $\partial\overline{\Sigma}^n(\mc V)$. Hence, it induces an isomorphism\glslink{evmap}{}
$$
\ev:\overline{\Sigma}(\mc V)/\partial\overline{\Sigma}(\mc V)\stackrel{\simeq}{\longrightarrow}
\Sigma(\mc V)\,.
$$
Next, we let $\overline\phi_0:\widetilde{\Omega}^0(\mc V)=\mc V\to \overline{\Sigma}^0(\mc V)=\mc V/[\mc V,\mc V]$ be the canonical projection map, and,
for $n\geq1$, we define a linear map $\overline\phi_n:\widetilde{\Omega}^n(\mc V)\to\overline{\Sigma}^n(\mc V)$ as follows (cf. \eqref{20140626:eq2-aff}).
Let $i_1,\dots,i_n\in I$, $m_1,\dots,m_n\in\mb Z_{\geq0}$ and $f_1,\dots f_{n+1}\in\mc V$
be fixed and let
$$
A=A(\lambda_1,\dots,\lambda_n)=\left(\delta_{j_1,i_1}\dots\delta_{j_n,i_n}f_1\otimes\dots\otimes f_{n+1} \lambda_1^{m_1}\dots\lambda_n^{m_n}\right)_{j_1,\dots,j_n\in I}\in\widetilde{\Omega}^n(\mc V)
\,.
$$
We map it to the array $\overline{\phi}_n(A)=\left(\overline{S}_{j_1,\dots,j_n}(\lambda_1,\dots,\lambda_{n})\right)_{j_1,\dots,j_n\in I}\in\overline{\Sigma}^n(\mc V)$
with entries 
$\overline{S}_{j_1\dots j_n}(\lambda_1,\dots,\lambda_{n})=0$
unless $(j_1,\dots,j_n)$ is a cyclic permutation of $(i_1,\dots,i_n)$,
and
\begin{align*}
\overline{S}_{j_1\dots j_n}(\lambda_1,\dots,\lambda_{n})
&=\frac1n (-1)^{s(n-s)}
\big(
f_{s+1}\otimes\dots\otimes f_n\otimes f_{n+1}f_1\otimes f_2\otimes\dots\otimes f_s
\big)\times
\\
&
\times \lambda_1^{m_{s+1}}\dots\lambda_{n-s}^{m_n}
\lambda_{n-s+1}^{m_1}\dots\lambda_{n}^{m_{s}}
\,,
\end{align*}
for $(j_1,\dots,j_n)=(i_{\sigma^s(1)},\dots,i_{\sigma^s(n)})$. We extend $\overline{\phi}_n$ to $\widetilde{\Omega}(\mc V)$ by linearity.
Note that $\widetilde{\phi}_n=\ev|_{\lambda_n=\lambda_n^\dagger}\circ\overline{\phi}_n$, and we have the following diagram
%\pecetta{this diagram has been commented to speed up compiling}
%\begin{comment}
$$
\begin{tikzcd}
\widetilde{\Omega}(\mc V)\arrow[r,"\overline{\phi}"]\arrow[rr,bend left=30,"\widetilde{\phi}"]
&
\overline{\Sigma}(\mc V)\arrow[r,"\ev"]
&
\Sigma(\mc V)
\end{tikzcd}
$$
%\end{comment}
Since both $\ev$ and $\widetilde{\phi}$ are surjective (cf. Proposition \ref{20240819:prop1}) we have that $\overline{\phi}$ is surjective as well. We leave it as an exercise for the reader to check that $\ker\overline{\phi}=[\widetilde{\Omega}(\mc V),\widetilde{\Omega}(\mc V)]$. Hence, we have the isomorphism
$$
\overline{\Sigma}(\mc V)\simeq\widetilde{\Omega}(\mc V)/[\widetilde{\Omega}(\mc V),\widetilde{\Omega}(\mc V)]
\,.
$$
We can thus choose a subspace $\widetilde{\Sigma}(\mc V)\subset\widetilde{\Omega}(\mc V)$
such that we have the direct sum decomposition\glslink{Sigmatilde}{}
\begin{equation}\label{20250901:eq1}
\widetilde{\Omega}(\mc V)=\widetilde{\Sigma}(\mc V)\oplus [\widetilde{\Omega}(\mc V),\widetilde{\Omega}(\mc V)]
\end{equation}
and the restriction $\overline{\phi}|_{\widetilde{\Sigma}(\mc V)}:\widetilde{\Sigma}(\mc V)
\to\overline{\Sigma}(\mc V)$ is an isomorphism. Since $\partial$ and $\overline{\phi}$ commute we have that $\partial$ preserves the direct sum decomposition \eqref{20250901:eq1}. The direct sum \eqref{20250901:eq1} will play an important role in the computation of the dPVA cohomology discussed in Chapter \ref{CH:PVAexa}.

%%%%%%%%%%% NEW CHAPTER %%%%%%%%%%%%%%%
%%%%%%%%%%% NEW CHAPTER %%%%%%%%%%%%%%%
%%%%%%%%%%% NEW CHAPTER %%%%%%%%%%%%%%%
%%%%%%%%%%% NEW CHAPTER %%%%%%%%%%%%%%%

\chapter[Computation of the \MakeLowercase{d}PVA cohomology]{Computation of the double variational Poisson cohomology for a constant
\texorpdfstring{$2$}{2}-fold \texorpdfstring{$\lambda$}{lambda}-bracket}\label{CH:PVAexa}

Let $\mc V:=\mc R_\ell=\kk\langle u^{(p)}_i \mid i\in I, p\in\mb Z_{\geq 0} \rangle$
be the algebra of noncommutative differential polynomials in $\ell$ variables $u_1,\dots,u_\ell$ (here $I=\{1,\dots,\ell\}$).
Let $M\in\mb Z_{\geq0}$ and let $K=(K_{ij})_{i,j=1}^\ell\in \Gl_{\ell}(\kk)$ be a symmetric matrix, if $M$ is odd, or a skew-symmetric matrix, if $M$ is even, respectively. We define
\begin{equation} \label{Eq:uu-odd-gen}
    \dgal{u_i{}_\lambda u_j}= K_{ji}(1\otimes 1)\, \lambda^M\,,
\end{equation}
and we extend it to a $2$-fold $\lambda$-bracket on $\mc V$ using the Master Formula \eqref{master-infinite}, namely
\begin{equation} \label{Eq:uu-odd-full}
\ldb f_\lambda g\rdb=\sum_{\substack{i,j\in I\\p,q\in\mb Z_{\geq0}}}K_{ji}(-1)^p
\frac{\partial g}{\partial u^{(q)}_j}\bullet (\lambda+\partial)^{p+q+M}
\bigg(\frac{\partial f}{\partial u^{(p)}_i}\bigg)^\sigma\,,
\end{equation}
for every $f,g\in\mc V$. 
Note that skewsymmetry \eqref{skew-gen-b} holds on the $u_i$'s since $K_{ij}=(-1)^{M+1} K_{ji}$ by assumption on $K$, and the Jacobi identity \eqref{jacobi-gen-b} trivially holds since any of the three terms of it are zero due to the fact that the $\lambda$-bracket
\eqref{Eq:uu-odd-gen} has constant coefficients. Hence, by Theorem \ref{20130921:prop1}(c), the $2$-fold $\lambda$-bracket
\eqref{Eq:uu-odd-full} defines a dPVA structure on $\mc V$.

In this chapter we describe the variational dPVA cohomology $\dPVH(\mc V)$ defined in Section \ref{sec:var-dPVA}.
Using the isomorphism \eqref{20250723:eq3} and Proposition \ref{20240819:prop1} it suffices to describe the cohomology $\coH(\Omega(\mc V),\delta)$ defined in Section \ref{sec:skewarray}.
Our strategy is to compute first the cohomology $\coH(\widetilde{\Omega}(\mc V),\tilde\delta)$ defined in Section \ref{sec:array} (which is isomorphic to the 
basic dPVA cohomology $\dH_{\textrm{bas}}(\mc V)$ defined in Section \ref{sec:omega2}) and then use the long exact sequence \eqref{eq:les2} to compute the cohomology $\coH(\Omega(\mc V),\delta)$ (which is isomorphic to the reduced dPVA cohomology $\dH_{\textrm{red}}(\mc V)$). Eventually, we use the isomorphism \eqref{20250723:eq4}.
\begin{remark}
We point out that the symmetric/skewsymmetric assumption on the matrix $K$ is not essential and we will drop it. Indeed, as explained in Remarks \ref{rem:deltaHconst} and \ref{rem:deltaHconst2}, since the $\lambda$-bracket \eqref{Eq:uu-odd-full} satisfies \eqref{ass-const}, we have cohomology complexes $(\widetilde{\Omega}(\mc V),\tilde\delta)$ and
$(\Sigma(\mc V),d)$ for every $K\in\Gl_{\ell}(\kk)$.
\end{remark}
%

%%%
\section[Computation of the basic \MakeLowercase{d}PVA cohomology]{Computation of the basic double Poisson vertex algebra cohomology}\label{sec:PVAexa1}

Let $n\in\mb Z_{\geq0}$. We write $A(\lambda_1,\dots,\lambda_n)=(A_{\underline{i}}(\lambda_1,\dots,\lambda_n))_{\underline{i}\in I^n}\in\widetilde{\Omega}^n(\mc V)$
where, for $\underline{i}=(i_1,\ldots,i_n)\in I^n$,  
$A_{\underline{i}}(\lambda_1,\ldots,\lambda_n)\in \mc V^{\otimes (n+1)}[\lambda_1,\dots,\lambda_n]$ is given by  
\begin{equation} \label{Eq:cochX}
   A_{\underline{i}}(\lambda_1,\ldots,\lambda_n)
= \sum_{\underline{a}\in\Z_{\geq 0}^n} f_{\underline{i},\underline{a}}\,\lambda^{\underline a}
= \sum_{\underline{a}\in\Z_{\geq 0}^n} f_{\underline{i},\underline{a}}^1\otimes \ldots \otimes f_{\underline{i},\underline{a}}^{n+1}\, \lambda_1^{a_1} \cdots \lambda_n^{a_n}\,.
\end{equation}
In \eqref{Eq:cochX} we sum over multi-indices $\underline{a}=(a_1,\ldots,a_n) \in \Z_{\geq 0}^n$
and we are using the notations
$f_{\underline{i},\underline{a}}=f_{\underline{i},\underline{a}}^1\otimes \ldots \otimes f_{\underline{i},\underline{a}}^{n+1} \in \mc V^{\otimes (n+1)}$ and $\lambda^{\underline a}=\lambda_1^{a_1} \cdots \lambda_n^{a_n}$.
Using \eqref{Eq:uu-odd-full}, we write the action \eqref{eq:diff-sigma} of the differential $\tilde\delta$ on an array $A(\lambda_1,\dots,\lambda_n)\in\widetilde{\Omega}^n(\mc V)$ explicitly as
 (note that the second sum in the equation for $\tilde\delta$ vanishes since the $\lambda$-bracket
\eqref{Eq:uu-odd-gen} has constant coefficients)
\begin{align}
    &\tilde\delta(A)_{\underline i}(\lambda_1,\ldots,\lambda_{n+1})
\nonumber    \\
   &=\sum_{s=1}^{n+1} \sum_{\substack{j\in I\\p\in\mb Z_{\geq 0}}} (-1)^{s+1}\, K_{ji_s}
   \bigg(\frac{\partial}{\partial u_{j}^{(p)}}\bigg)_{(s)}  A_{i_1,\stackrel{s}{\check{\dots}},i_{n+1}}(\lambda_1,\stackrel{s}{\check{\dots}},\lambda_{n+1}) \lambda_s^{M+p} \,,
    \label{Eq:coch-dPX-1}
  \\
  &=\sum_{s=1}^{n+1} \sum_{\substack{j\in I\\p\in\mb Z_{\geq 0}}}\sum_{\underline a\in\mb Z_{\geq0}^n} (-1)^{s+1}\, K_{ji_s}
   \bigg(\frac{\partial}{\partial u_{j}^{(p)}}\bigg)_{(s)}  f_{\underline i{\stackrel{s}{\check{}}},\underline{a}} 
 \lambda_1^{a_1} \cdots \lambda_{s-1}^{a_{s-1}}   \lambda_s^{M+p} \lambda_{s+1}^{a_{s}} \cdots \lambda_{n+1}^{a_{n}}\,,
    \label{Eq:coch-dPX}
\end{align}
where $\underline i=(i_1,\ldots,i_{n+1})\in I^{n+1}$
and $\underline i{\stackrel{s}{\check{}}}=(i_1,\dots,i_{s-1},i_{s+1},\dots,i_{n+1})\in I^n$.

Let $\Delta\in\Vect^\partial(\mc V)$ be the \emph{degree operator} on $\mc V$. It is the derivation of $\mc V$ commuting with $\partial$  defined by (cf. \eqref{eq:QS-1})
\begin{equation}\label{eq:degreeop}
\Delta=\sum_{\substack{h\in I\\p\in\mb Z_{\geq0}}}
\mult\left(u^{(p)}_h\star\frac{\partial}{\partial u^{(p)}_h}\right)\,.
\end{equation}
The degree operator is diagonalizable on $\mc V$ and we have the direct sum decomposition
$$
\mc V=\bigoplus_{k\in\mb Z_{\geq0}}\mc V_k\,,
\qquad
\mc V_k=\{f\in\mc V\mid \Delta(f)=k f\}
\,.
$$
If $f\in\mc V$ is homogeneous of degree ${\Delta_f}$ and $g\in\mc V$ is homogeneous of degree $\Delta_g$, then clearly $fg$ is homogeneous
of degree $\Delta_{fg}=\Delta_f+\Delta_g$.
We extend $\Delta:\mc V\to\mc V$ to a derivation $\Delta:\mc V^{\otimes n}\to\mc V^{\otimes n}$
using \eqref{mfold-ext}. Explicitly, recalling \eqref{eq:degreeop}, we have
\begin{equation}\label{eq:Delta-times}
\begin{split}
&\Delta(f_1\otimes \dots\otimes f_{n})
\\
&=\sum_{i=1}^n\sum_{\substack{h\in I\\p\in\mb Z_{\geq0}}}f_1\otimes\dots\otimes f_{i-1}\otimes\bigg(\frac{\partial f_i}{\partial u_h^{(p)}}\bigg)'
u_h^{(p)}\bigg(\frac{\partial f_i}{\partial u_h^{(p)}}\bigg)''\otimes f_{i+1}\otimes\dots\otimes f_n
\\
&=\sum_{i=1}^n\sum_{\substack{h\in I\\p\in\mb Z_{\geq0}}}
\mult_{(i,i+1)}\bigg(u_h^{(p)}\star_i\bigg(\frac{\partial }{\partial u_h^{(p)}}
\bigg)_{(i)}
\bigg)
(f_1\otimes\dots\otimes f_n)
\,,
\end{split}
\end{equation}
for every $f_1,\dots,f_n\in \mc V$.
This operator is still diagonalizable and we have the direct sum decomposition in $\Delta$-eigenspaces of $\mc V^{\otimes (n+1)}$:
\begin{equation}\label{Delta-dec}
(\mc V^{\otimes (n+1)})=\bigoplus_{k\in\mb Z_{\geq0}}(\mc V^{\otimes (n+1)})_k
\,,
\qquad
(\mc V^{\otimes (n+1)})_k=\{f\in\mc V^{\otimes (n+1)}\mid \Delta(f)=k f\}
\,.
\end{equation}
Given $\underline{a}=(a_1,\dots,a_n)\in \Z_{\geq 0}^n$, we let 
\begin{equation}\label{eq:dM}
    d_M(\underline a) := \#\, \{a_j \mid a_j\geq M\}\,,
\end{equation}
which returns the number of entries of $\underline a$ which are greater or equal to $M$.
Note that $d_M(\underline a)\in\{0,1,\dots,n\}$ and that we have the direct sum decomposition
\begin{equation}\label{dm-dec}
\kk[\lambda_1,\dots,\lambda_n]=\bigoplus_{d=0}^n U_d^n
\,,
\qquad
U_d^n=\Span_\kk\{\lambda^{\underline a}\mid \underline a\in I^n\text{ s.t. }d_M(\underline a)=d\}
\,.
\end{equation}
Combining \eqref{Delta-dec} and \eqref{dm-dec} 
we obtain a bi-graded decomposition 
$$\mc V^{\otimes (n+1)}[\lambda_1,\ldots,\lambda_n]=\bigoplus_{k\in\mb Z_{\geq 0}} \bigoplus_{d=0}^n \mc U_{k,d}^n, 
\qquad \mc U_{k,d}^n:=(\mc V^{\otimes (n+1)})_k\otimes U_d^n\,,$$
which, in turn, induces the direct sum decomposition 
\begin{equation} \label{Eq:Biwei}
    \widetilde\Omega^n(\mc V) 
    = \bigoplus_{k\in\mb Z_{\geq 0}} \bigoplus_{d=0}^n \widetilde \Omega^n_{k,d}(\mc V)\,,
\end{equation} 
where $\widetilde{\Omega}_{k,d}^n$ is the linear span of elements $A(\lambda_1,\dots,\lambda_n)=(A_{\underline i}(\lambda_1,\dots,\lambda_n))\in\widetilde{\Omega}^n(\mc V)$ such that $A_{\underline{i}}(\lambda_1,\dots,\lambda_n)\in \mc U_{k,d}^n$ for every $\underline{i}=(i_1,\dots,i_n)\in I^n$. 
%In that case, $A_{\underline{i}}= \sum_{\underline a\in\mb Z_{\geq0}^n}f_{\underline{i},\underline{a}} \lambda^{\underline a}$ is such that 
%$\Delta(f_{\underline{i},\underline{a}} )=kf_{\underline{i},\underline{a}}$ 
%and $d_M(\underline a)=d$.
In particular, $\widetilde{\Omega}^n_{0,0}\subset\widetilde{\Omega}^n(\mc V)$ consists of arrays
whose entries $(A_{\underline i}(\lambda_1,\dots,\lambda_n))$ are polynomials of degree at most $N-1$ in each variable $\lambda_1,\dots,\lambda_n$, with constant coefficients.
Let
\begin{equation}\label{omega00}
\widetilde{\Omega}_{0,0}^\bullet
=\bigoplus_{n\in\mb Z_{\geq0}}\widetilde{\Omega}_{0,0}^n
\,.
\end{equation}
It is clear from \eqref{Eq:coch-dPX-1} and \eqref{20250723:eq1} that $\tilde\delta$ acts trivially on
$\widetilde{\Omega}^\bullet_{0,0}$ so that $(\widetilde{\Omega}^\bullet_{0,0},0)\subset(\widetilde{\Omega}(\mc V),\tilde\delta)$ is a subcomplex.

In this section we prove the following result.
\begin{theorem}  \label{Thm:DPVcoh-Const}
Consider $\mc V=\kk\langle u^{(p)} \mid p\in\mb Z_{\geq 0} \rangle$ with the dPVA structure defined by \eqref{Eq:uu-odd-full}. 
Then, the inclusion
$(\widetilde{\Omega}^\bullet_{0,0},0)\subset(\widetilde{\Omega}(\mc V),\tilde\delta)$
induces an isomorphism in cohomology
$$
\coH^n(\widetilde{\Omega}(\mc V),\tilde\delta)\simeq \widetilde{\Omega}_{0,0}^n\,,
$$
for every $n\in\mb Z_{\geq0}$.
Hence, we have    $\dim\left(\coH^n(\widetilde{\Omega}(\mc V),\tilde\delta) \right) = (\ell M)^n$,
for every $n\in\mb Z_{\geq0}$.
\end{theorem}
Theorem
\ref{Thm:DPVcoh-Const} is a (partial) noncommutative analogue of Theorem 11.2 in  \cite{DSK13}, see also \cite{CCS,DSK12}. For $M=0$ and $C=\Id_\ell$ it reduces to the computation of the noncommutative de Rham cohomology in \cite{DSKV}. Hence, we refer to the complex $(\widetilde{\Omega}(\mc V),\tilde\delta)$ as the \emph{generalized noncommutative de Rham complex}.
%
%\begin{remark}
%Recall from Remark \ref{rem:deltaHconst} that we have a well defined differential $\delta_H$
%on $\widetilde{\Sigma}(\mc V)$ also when $m$ is even in \eqref{Eq:uu-odd}.
%%
%By replacing, in the even case, the skewsymmetry axiom \eqref{eq:skew2} of the $\lambda$-bracket by its symmetric version, the proof, thus the statement, of Theorem \ref{Thm:DPVcoh-Const} is valid for any integer $m\geq 0$.
%The case $m=0$ boils down to the Beltrami $2$-fold $\lambda$-bracket and one recovers the acyclicity of the de Rham complex \eqref{20140623:eq6a-b} presented in \cite[\S3.5]{DSKV}. 
%\end{remark}

The remainder of the section is devoted to the proof of Theorem \ref{Thm:DPVcoh-Const}.
%First, note that, using the identification between $\widetilde{\Gamma}(\mc V)\simeq\widetilde{\Omega}(\mc V)$ given by the map $\widetilde{\gamma}$ defined by \eqref{eq:idsigma} and the isomorphism \eqref{20250723:eq2}, we can equivalently prove the statement for the complex $(\widetilde{\Omega}(\mc V),\tilde\delta)$, where $\tilde\delta$ is given by
%\eqref{eq:diff-sigma}.
First we prove the results that will be needed in its proof. In Lemma \ref{Lem:DPVCoh} we provide a necessary condition for an array $A(\lambda_1,\dots,\lambda_n)\in\widetilde{\Omega}^n(\mc V)$ to lie in $\ker\tilde\delta$. Then, we prove a global homotopy condition in Proposition \ref{Pr:Homot}.

For every $n,M\in\mb Z_{\geq0}$ we let $T_n(M)\subset\mb Z_{\geq0}^n$ be the subset of multi-indices defined by
$$
T_n(M)=\{\underline a=(a_1,\dots a_n)\mid a_t<M\text{ for every }t=1,\dots,n\}
\,.
$$

\begin{lemma} \label{Lem:DPVCoh}
    Let $A(\lambda_1,\dots,\lambda_n)\in\widetilde{\Omega}^n(\mc V)$ be as in \eqref{Eq:cochX}. 
    If $\tilde\delta(A)=0$, 
then $f_{\underline{i},\underline a}\in \kk$ for all $\underline i \in I^n$ and $\underline a\in T_n(M)$.
\end{lemma}
\begin{proof}
Let us assume that $\tilde\delta(A)=0$ and let us fix $\underline a\in T_n(M)$ and $s\in \{1,\ldots,n+1\}$.
%let $(\underline{i},\underline a)\in I^n \times T_n(M)$ be such that $f_{\underline i,\underline a}\neq0$. 
%\pecetta{I think we should request $f_{\underline a}\neq0$.}
For $1\leq t \leq n$, we successively apply to the RHS of \eqref{Eq:coch-dPX} the operator 
\begin{equation*}
    \left\{ 
\begin{array}{ll}
\frac{1}{a_t!} \,(\partial_{\lambda_t})^{a_t}\,, &\text{if }t<s\,,\\
\frac{1}{a_t!} \,(\partial_{\lambda_{t+1}})^{a_t}\,,&\text{if }t\geq s\,,
\end{array}
\right.
\end{equation*} 
where $\partial_{\lambda_i}$, $i=1,\dots,n$, denotes the usual partial derivative with respect to $\lambda_i$ of 
$\mc V^{\otimes (n+1)}[\lambda_1,\dots,\lambda_n]$ (it acts trivially on $\mc V^{\otimes (n+1)}$).  
We then evaluate at $\lambda_1=\ldots=\lambda_{s-1}=\lambda_{s+1}=\ldots=\lambda_{n+1}=0$ to get
$$
\sum_{\substack{j\in I\\p\in\mb Z_{\geq 0}}}\, K_{ji_s}
\bigg(\frac{\partial}{\partial u^{(p)}_j}\bigg)_{(s)}  f_{i_1,\stackrel{s}{\check{\ldots}},i_{n+1},\underline{a}} \, \lambda_s^{M+p} =0\,.
$$
%(This is where we used the assumption that all $a_i<m$, otherwise there would be extra terms associated with multi-indices $\underline{b} \neq \underline{a}$.) 
This equality is true for any multi-index $(i_1,\ldots,i_{n+1})\in I^{n+1}$, and therefore we deduce 
$$
\sum_{j\in I}\, K_{ji'}
\bigg(\frac{\partial}{\partial u^{(p)}_j}\bigg)_{(s)}  f_{\underline{i},\underline{a}} =0\,,
$$
for any fixed $\underline{i}\in I^n$, $i'\in I$ and $p\in\mb Z_{\geq0}$.  
Hence, invertibility of $K$ guarantees 
\begin{equation*}
 f_{\underline{i},\underline{a}}^1\otimes \ldots \otimes f_{\underline{i},\underline{a}}^{s-1}\otimes 
 \frac{\partial f_{\underline{i},\underline{a}}^s}{\partial u_j^{(p)}} 
 \otimes 
 f_{\underline{i},\underline{a}}^{s+1}\otimes \ldots 
 \otimes f_{\underline{i},\underline{a}}^{n+1} = 0\,,
\end{equation*}
for all $j\in I$ and $p\in\mb Z_{\geq0}$. By equation \eqref{20250723:eq1} it follows that $f_{\underline{i},\underline{a}}^s\in\kk$.
This is true for all $s=1,\dots,n+1$, therefore $f_{\underline{i},\underline{a}}\in \kk^{\otimes (n+1)} \simeq \kk$, as desired.   
\end{proof}
We introduce the further notation 
\begin{equation} \label{Eq:BiweiB}
\widetilde \Omega^n_{\bullet,\geq 1}(\mc V) := \bigoplus_{k\geq 0} \bigoplus_{d=1}^n \widetilde\Omega^n_{k,d}(\mc V)\,.
\end{equation} 
The Lie derivative $L_\Delta:\widetilde{\Omega}(\mc V)\to\widetilde{\Omega}(\mc V)$ (cf. \cite[Sec.~3.5]{DSKV}) of the degree operator $\Delta$ is called the \emph{Euler operator}
and it is an infinite-variable analogue of the derivation $\mathrm{E}$ defined in \eqref{Eq:Euler}, see \cite{PV}.
%
%Using equation \eqref{eq:LP} and the identification $\widetilde{\Omega}(\mc V)\simeq\widetilde{\Sigma}(\mc V)$ we get the explicit expression for the induced action of $L_\Delta$ on $\widetilde{\Sigma}(\mc V)$. 
%
For $A=(A_{\underline{i}}(\lambda_1,\dots,\lambda_n))_{\underline{i}\in I^n}$ as in \eqref{Eq:cochX}
the Euler operator is defined as the array $L_\Delta (A)\in \widetilde{\Omega}^n(\mc V)$ with entries 
\begin{equation}\label{eq:LDelta}
(L_\Delta (A))_{\underline{i}}=\sum_{\underline a\in\mb Z_{\geq0}}(\Delta +n)(f_{\underline{i},\underline a})\lambda^{\underline a}\,.
\end{equation}
Note that $L_{\Delta}$ is an even derivation of $\widetilde{\Omega}(\mc V)$. There is also a contraction operator $\iota_\Delta:\widetilde{\Omega}(\mc V)\to\widetilde{\Omega}(\mc V)$ associated to the degree operator $\Delta$
(cf. \cite[Sec.~3.5]{DSKV}) which is defined as follows: for $A=(A_{\underline{i}}(\lambda_1,\dots,\lambda_n))_{\underline{i}\in I^n}$ as in \eqref{Eq:cochX} we have that $\iota_\Delta(A)\in\widetilde{\Omega}^{n-1}(\mc V)$ is the array with entries ($\underline i=(i_1,\dots,i_{n-1})\in I^{n-1}$)
\begin{align}
&(\iota_\Delta (A))_{\underline i}(\lambda_1,\dots,\lambda_{n-1})\nonumber
\\
\begin{split}\label{eq:iotaDelta}
&=\sum_{s=1}^n\sum_{j\in I}\sum_{\underline a\in\mb Z_{\geq0}^n}(-1)^{s+1}
f_{\underline{i}_{\hat{s}},\underline a}^1\otimes \dots \otimes f_{\underline{i}_{\hat{s}},\underline a}^{s-1}\otimes f_{\underline{i}_{\hat{s}},\underline a}^{s} u_{j}^{(a_s)}
f_{\underline{i}_{\hat{s}},\underline a}^{s+1}\otimes 
f_{\underline{i}_{\hat{s}},\underline a}^{s+2}\otimes \dots
\\
&\quad\quad\quad\quad\quad\quad\quad\quad
\dots\otimes f_{\underline{i}_{\hat{s}},\underline a}^{n+1} \, 
\lambda_1^{a_1}\dots\lambda_{s-1}^{a_{s-1}}\lambda_{s}^{a_{s+1}}\dots \,\lambda_{n-1}^{a_n}\,,
\end{split}
\\
\begin{split}
&=\sum_{s=1}^n\sum_{j\in I}(-1)^{s+1}
\mult_{(s,s+1)}
\big(
(|_{x=\partial}u_{j})
\star_s
A_{i_1,\dots,i_{s-1},j,i_s,\dots,i_{n-1}}
(\lambda_1,\dots
\\
&\quad\quad\quad\quad\quad\quad\quad\quad
\dots,\lambda_{s-1},x,\lambda_{s},\dots,\lambda_{n-1})
\big)
\,,
\end{split}
\end{align}
where $\underline{i}_{\hat{s}}:=(i_1,\ldots,i_{s-1},j,i_{s},\ldots,i_{n-1})\in I^n$ means that we insert the index of summation $j$ in position $s$ (that is between $i_{s-1}$ and $i_s$). Note that $\iota_{\Delta}$ is an odd derivation of $\widetilde{\Omega}(\mc V)$.

We introduce the following modifications of the operators $L_\Delta$ and $\iota_\Delta$
which depend on $M\in\mb Z$, and on $K$ for the contraction operator.
We let $L_{\Delta,M},\iota_{\Delta,M}:\widetilde{\Omega}(\mc V)\to\widetilde{\Omega}(\mc V)$ be the 
linear operators defined, for $A=A(\lambda_1,\dots,\lambda_n)$ as in \eqref{Eq:cochX}, by\glslink{PVAL}{}
\begin{equation}\label{Eq:Lm}
(L_{\Delta,M} (A))_{\underline{i}}=\sum_{\underline a\in\mb Z_{\geq0}}(\Delta + d_M(\underline a))(f_{\underline{i},\underline a})\lambda^{\underline a}\,.
\end{equation}
and\glslink{PVAiota}{}
\begin{align}
&(\iota_{\Delta,M} (A))_{i_1,\dots,i_{n-1}}(\lambda_1,\dots,\lambda_{n-1})
\nonumber
\\
\begin{split}
\label{Eq:iotam}
&=\sum_{s=1}^n\sum_{i,j\in I}\sum_{\underline a\in\mb Z_{\geq0}^n}(-1)^{s+1}
(K^{-1})_{jk}
f_{\underline{i}_{\hat{s}},\underline a}^1\otimes \dots
\otimes f_{\underline{i}_{\hat{s}},\underline a}^{s-1}\otimes f_{\underline{i}_{\hat{s}},\underline a}^{s} u_{k}^{(a_s-M)}
f_{\underline{i}_{\hat{s}},\underline a}^{s+1}
\\
&\quad\quad\quad\quad\quad\quad\quad\quad
\otimes 
f_{\underline{i}_{\hat{s}},\underline a}^{s+2}\otimes \dots\otimes f_{\underline{i}_{\hat{s}},\underline a}^{n+1} \, 
\lambda_1^{a_1}\dots\lambda_{s-1}^{a_{s-1}}\lambda_{s}^{a_{s+1}}\dots \,\lambda_{n-1}^{a_n}
\end{split}
\\
\begin{split}\label{Eq:iotam-1}
&=\sum_{s=1}^n\sum_{j,k\in I}(-1)^{s+1}
(K^{-1})_{jk}
\mult_{(s,s+1)}
\big(
(|_{u=\partial}x^{-M}u_{k})
\star_s
A_{i_1,\dots,i_{s-1},j,i_s,\dots,i_{n-1}}
(\lambda_1,\dots
\\
&\quad\quad\quad\quad\quad\quad\quad\quad
\dots,\lambda_{s-1},x,\lambda_{s},\dots,\lambda_{n-1})
\big)
\,,
\end{split}
\end{align}
where, as before, $\underline{i}_{\hat{s}}:=(i_1,\ldots,i_{s-1},j,i_{s},\ldots,i_{n-1})\in I^n$.
In \eqref{Eq:iotam} (and \eqref{Eq:iotam-1})
we assume that $u^{(p)}=0$ for $p<0$. Clearly, $L_{\Delta,0}=L_{\Delta}$ and $\iota_{\Delta,0}=\iota_{\Delta}$ when $C=\Id_\ell$. Moreover, $L_{\Delta,M}$ is an even derivation of $\widetilde{\Omega}(\mc V)$, while $\iota_{\Delta,M}$ is an odd derivation
of $\widetilde{\Omega}(\mc V)$ thus they both preserve the commutator subspace $[\widetilde{\Omega}(\mc V),\widetilde{\Omega}(\mc V)]\subset\widetilde{\Omega}(\mc V)$.
Note that $L_{\Delta,M}$ is diagonalizable on $\widetilde{\Omega}(\mc V)$,
it preserves the direct sum decomposition \eqref{Eq:Biwei}, and it is invertible on
$$
W:=\bigoplus_{\substack{n,k\in\mb Z_{\geq0}\\ 0\leq d\leq n\\ (k,d)\neq(0,0) }}
\widetilde{\Omega}^n_{k,d}\subset\widetilde{\Omega}(\mc V)
\,.
$$
In particular, on any direct summand $\widetilde{\Omega}^n_{k,d}$ of $W$, the inverse of $L_{\Delta,M}$ is the operator of multiplication by $(k+d)^{-1}$.

From \eqref{Eq:iotam} we have that
$\iota_{\Delta,M}(\widetilde{\Omega}^n_{k,d})\subset\widetilde{\Omega}^{n-1}_{k+1,d-1}$. 
Hence,
$$
\iota_{\Delta,M}\bigg(\bigoplus_{n\geq1}\widetilde{\Omega}^n_{\bullet,\geq1}\bigg)\subset W
\,,
$$
and therefore we can define the \emph{homotopy operator}\glslink{PVAhm}{} 
\begin{equation}\label{eq:hm}
 h_M = (L_{\Delta,M})^{-1} \circ \iota_{\Delta,M}\big|_{\bigoplus_{n\geq1}\widetilde{\Omega}^n_{\bullet,\geq1}} : \bigoplus_{n\geq 1} \widetilde \Omega^n_{\bullet,\geq 1}(\mc V) \to \widetilde \Omega(\mc V)\,.
\end{equation}
The next result is motivated by \cite[Sec.~3.5]{DSKV}. 
\begin{proposition} \label{Pr:Homot}
    The following identity holds for every $A\in\bigoplus_{n\geq 1}\widetilde \Omega^n_{\bullet,\geq 1}(\mc V)$: 
\begin{equation}\label{eq:homotopy}
    (\tilde\delta \circ h_M + h_M \circ \tilde\delta)(A) = A\,.
\end{equation}
\end{proposition}
\begin{proof}
By inspection of \eqref{Eq:coch-dPX} we have $\tilde\delta(\widetilde{\Omega}^n_{k,d})\subset\widetilde{\Omega}^{n+1}_{k-1,d+1}$. Then,
\begin{equation} \label{Eq:ImdeltP}
    \im \tilde\delta \subseteq \bigoplus_{n\geq 1} \widetilde \Omega^n_{\bullet,\geq 1}(\mc V)
\,,
\end{equation}
and therefore the LHS of \eqref{eq:homotopy} is well-defined.
Furthermore, we already noted that $\iota_{\Delta,M}(\widetilde{\Omega}^n_{k,d})\subset\widetilde{\Omega}^{n-1}_{k+1,d-1}$. Then, it is immediate to check that
$$
(\tilde\delta \circ h_M + h_M \circ \tilde\delta)(\widetilde{\Omega}^n_{k,d})
\subset\widetilde{\Omega}^n_{k,d}
\,.
$$
Recall also that $(L_{\Delta,M})^{-1}$ is the multiplication by $(k+d)^{-1}$ on $\widetilde{\Omega}^n_{k,d}$. Hence, for $A\in\widetilde{\Omega}^n_{k,d}$, the identity
\eqref{eq:homotopy} is equivalent to
\begin{equation}\label{eq:homotopy2}
    (\tilde\delta \circ \iota_{\Delta,M} + \iota_{\Delta,M} \circ \tilde\delta)(A) =(k+d) A\,.
\end{equation}
(This is a generalization of Cartan's formula from \cite{DSKV}.)
Let $A(\lambda_1,\dots,\lambda_n)\in\widetilde{\Omega}^{n}_{k,d}$. Using \eqref{Eq:coch-dPX-1} and \eqref{Eq:iotam-1}, we find
\begin{subequations}
\begin{align}
&(\tilde\delta\circ \iota_{\Delta,M}(A))_{i_1,\dots,i_n}(\lambda_1,\dots,\lambda_{n})\nonumber
\\
\begin{split}\label{A1}
&=\sum_{t=1}^{n} \sum_{s=1}^{t-1}
\sum_{\substack{i,j,h\in I\\p\in\mb Z_{\geq0}}}(-1)^{s+t}K_{hi_t}(K^{-1})_{ij}
\bigg(\frac{\partial}{\partial u_{h}^{(p)}}\bigg)_{(t)}\mult_{(s,s+1)}
\\
&\left(
(|_{x=\partial}x^{-M}u_{j})
\star_s
A_{i_1,\dots,i_{s-1},i,i_s,\stackrel{t+1}{\check{\dots}},i_{n}}
(\lambda_1,\dots,\lambda_{s-1},x,\lambda_{s},\stackrel{t+1}{\check{\dots}},\lambda_{n})
\right)\lambda_t^{p+M}
\end{split}
\\
\begin{split}\label{A2}
&+\sum_{t=1}^{n}
\sum_{\substack{i,j,k\in I\\p\in\mb Z_{\geq0}}}K_{hi_t}(K^{-1})_{ij}
\bigg(\frac{\partial}{\partial u_{h}^{(p)}}\bigg)_{(t)}\mult_{(t,t+1)}
\\
&\left(
(|_{x=\partial}x^{-M}u_{j})
\star_t
A_{i_1,\dots,i_{t-1},i,i_{t+1},\dots,i_{n}}
(\lambda_1,\dots,\lambda_{t-1},x,\lambda_{t+1},\dots,\lambda_{n})
\right)\lambda_t^{p+M}
\end{split}
\\
\begin{split}\label{A3}
&+\sum_{t=1}^{n} \sum_{s=t+1}^{n}
\sum_{\substack{i,j,h\in I\\p\in\mb Z_{\geq0}}}(-1)^{s+t}K_{hi_t}(K^{-1})_{ij}
\bigg(\frac{\partial}{\partial u_{h}^{(p)}}\bigg)_{(t)}\mult_{(s,s+1)}
\\
&\left(
(|_{x=\partial}x^{-M}u_{j})
\star_s
A_{i_1,\stackrel{t}{\check{\dots}},i_{s},i,i_{s+1},\dots,i_{n}}
(\lambda_1,\stackrel{t}{\check{\dots}},\lambda_{s},x,\lambda_{s+1},\dots,\lambda_{n})
\right)\lambda_t^{p+M}
\,.
\end{split}
\end{align}
\end{subequations}
Using the identity \eqref{Eq:Dmult1}, 
which holds for every $2$-fold derivation $D:\mc V\to\mc V^{\otimes 2}$, $a\in \mc V$, $X\in\mc V^{\otimes (n+1)}$, 
%$t=1,\dots,n$ and $1\leq s\leq t-1$,
we rewrite \eqref{A1} as
\begin{equation}\label{A1-1}
\begin{split}
&\sum_{t=1}^{n} \sum_{s=1}^{t-1}
\sum_{\substack{i,j,h\in I\\p\in\mb Z_{\geq0}}}(-1)^{s+t}K_{hi_t}(K^{-1})_{ij}
\mult_{(s,s+1)}\bigg((|_{x=\partial}x^{-M}u_{j})
\star_s
\\
&
\bigg(\frac{\partial}{\partial u_{h}^{(p)}}\bigg)_{(t+1)}
A_{i_1,\dots,i_{s-1},i,i_s,\stackrel{t+1}{\check{\dots}},i_{n}}
(\lambda_1,\dots,\lambda_{s-1},x,\lambda_{s},\stackrel{t+1}{\check{\dots}},\lambda_{n})
\bigg)
\lambda_t^{p+M}
\end{split}
\end{equation}
Using the identity \eqref{Eq:Dmult2}, 
%which holds for every $2$-fold derivation $D:\mc V\to\mc V^{\otimes 2}$, $a\in \mc V$, $X\in\mc V^{\otimes n+1}$ and $t=1,\dots,n$, 
and performing some algebraic manipulations,
we rewrite \eqref{A2} as
\begin{equation}\label{A2-1}
\begin{split}
&\sum_{t=1}^{n}
\sum_{\substack{i,j,h\in I\\p\in\mb Z_{\geq0}}}K_{hi_t}(K^{-1})_{ij}
\mult_{(t+1,t+2)}\bigg(
(|_{x=\partial}x^{-M}u_{j})
\star_{t+1}
\\
&
\bigg(\frac{\partial}{\partial u_{h}^{(p)}}\bigg)_{(t)}
A_{i_1,\dots,i_{t-1},i,i_{t+1},\dots,i_{n}}
(\lambda_1,\dots,\lambda_{t-1},x,\lambda_{t+1},\dots,\lambda_{n})
\bigg)\lambda_t^{p+M}
\\
&+\sum_{t=1}^{n}
\#\{\lambda_t^a\mid a-M\geq0\}
A_{i_1,\dots,i_{n}}
(\lambda_1,\dots,\lambda_{n})
\\
&+\sum_{t=1}^{n}
\sum_{\substack{i,j,h\in I\\p\in\mb Z_{\geq0}}}K_{hi_t}(K^{-1})_{ij}
\mult_{(t,t+1)}\bigg(
(|_{x=\partial}x^{-M}u_{j})
\star_t
\\
&
\bigg(\frac{\partial}{\partial u_{h}^{(p)}}\bigg)_{(t+1)}
A_{i_1,\dots,i_{t-1},i,i_{t+1},\dots,i_{n}}
(\lambda_1,\dots,\lambda_{t-1},x,\lambda_{t+1},\dots,\lambda_{n})
\bigg)\lambda_t^{p+M}
\,.
\end{split}
\end{equation}
Moreover, using the identity \eqref{Eq:Dmult3}, 
%which holds for every $2$-fold derivation $D:\mc V\to\mc V^{\otimes 2}$, $a\in \mc V$, $X\in\mc V^{\otimes n+1}$, $t=1,\dots,n$ and $t+1\leq s\leq n$,
we rewrite \eqref{A3} as
\begin{equation}\label{A3-1}
\begin{split}
&\sum_{t=1}^{n} \sum_{s=t+1}^{n}
\sum_{\substack{i,j,h\in I\\p\in\mb Z_{\geq0}}}(-1)^{s+t}K_{hi_t}(K^{-1})_{i,j}
\mult_{(s+1,s+2)}
\bigg((|_{x=\partial}x^{-M}u_{j})\star_{s+1}
\\
&
\bigg(\frac{\partial}{\partial u_{h}^{(p)}}\bigg)_{(t)}
A_{i_1,\stackrel{t}{\check{\dots}},i_{s-1},i,i_s,\dots,i_{n}}
(\lambda_1,\stackrel{t}{\check{\dots}},\lambda_{s-1},x,\lambda_{s},\dots,\lambda_{n})
\bigg)\lambda_t^{p+M}
\end{split}
\end{equation}
Combining \eqref{A1-1} and the third sum in \eqref{A2-1}, combining \eqref{A3-1} and the first sum
in \eqref{A2-1}, and using the definition of $d_M$ given in \eqref{eq:dM} and the fact that $A_{i_1,\dots,i_n}\in\mc U_{k,d}^n$ we then have
\begin{equation}\label{eq:AA1}
\begin{split}
&(\tilde\delta\circ \iota_{\Delta,M}(A))_{i_1,\dots,i_n}(\lambda_1,\dots,\lambda_{n})
\\
&=\sum_{t=1}^{n} \sum_{s=1}^{t}
\sum_{\substack{i,j,h\in I\\p\in\mb Z_{\geq0}}}(-1)^{s+t}K_{hi_t}(K^{-1})_{ij}
\mult_{(s,s+1)}\bigg((|_{x=\partial}x^{-M}u_{j})
\star_s
\\
&
\bigg(\frac{\partial}{\partial u_{h}^{(p)}}\bigg)_{(t+1)}
A_{i_1,\dots,i_{s-1},i,i_s,\stackrel{t+1}{\check{\dots}},i_{n}}
(\lambda_1,\dots,\lambda_{s-1},x,\lambda_{s},\stackrel{t+1}{\check{\dots}},\lambda_{n})
\bigg)
\lambda_t^{p+M}
\\
&+d\,
A_{i_1,\dots,i_{n}}
(\lambda_1,\dots,\lambda_{n})
\\
&+\sum_{t=1}^{n} \sum_{s=t}^{n}
\sum_{\substack{i,j,h\in I\\p\in\mb Z_{\geq0}}}(-1)^{s+t}K_{hi_t}(K^{-1})_{ij}
\mult_{(s+1,s+2)}
\bigg((|_{x=\partial}x^{-M}u_{j})\star_{s+1}
\\
&
\bigg(\frac{\partial}{\partial u_{h}^{(p)}}\bigg)_{(t)}
A_{i_1,\stackrel{t}{\check{\dots}},i_{s},i,i_{s+1},\dots,i_{n}}
(\lambda_1,\stackrel{t}{\check{\dots}},\lambda_{s},x,\lambda_{s+1},\dots,\lambda_{n})
\bigg)\lambda_t^{p+M}
\,.
\end{split}
\end{equation}
On the other hand, using \eqref{Eq:iotam-1} and \eqref{Eq:coch-dPX-1}, we have
\begin{equation}\label{B1}
\begin{split}
&(\iota_{\Delta,m}\circ \tilde\delta)(A)_{i_1,\dots,i_n}(\lambda_1,\dots,\lambda_{n})
\\
&=\sum_{s=1}^{n+1}\sum_{t=1}^{s-1}
\sum_{\substack{i,j,h\in I\\p\in\mb Z_{\geq0}}}(-1)^{s+t}
(K^{-1})_{ij}K_{hi_t}\mult_{(s,s+1)}\bigg(
(|_{x=\partial}x^{-M}u_j)\star_s
\\
&\bigg(\frac{\partial}{\partial u_h^{(p)}}
\bigg)_{(t)}A_{i_1,\stackrel{t}{\check{\dots}},i_{s-1},i,i_s,\dots,i_n}
(\lambda_1,\stackrel{t}{\check{\dots}},\lambda_{s-1},x,\lambda_s,\dots,\lambda_n)\lambda_t^{p+M}\bigg)
\\
&+\sum_{s=1}^{n+1}
\sum_{\substack{h\in I\\p\in\mb Z_{\geq0}}}
\mult_{(s,s+1)}\bigg(
u_h^{(p)}\star_s
\bigg(\frac{\partial}{\partial u_h^{(p)}}
\bigg)_{(s)}A_{i_1,\dots,i_n}
(\lambda_1,\dots,\lambda_n)\bigg)
\\
&+\sum_{s=1}^{n+1}\sum_{t=s+1}^{n+1}
\sum_{\substack{i,j,h\in I\\p\in\mb Z_{\geq0}}}(-1)^{s+t}
(K^{-1})_{ij}K_{hi_{t-1}}\mult_{(s,s+1)}\bigg(
(|_{x=\partial}x^{-M}u_j)\star_s
\\
&\bigg(\frac{\partial}{\partial u_h^{(p)}}
\bigg)_{(t)}A_{i_1,\dots,i_{s-1},i,i_s,\stackrel{t}{\check{\dots}},i_n}
(\lambda_1,\dots,\lambda_{s-1},x,\lambda_s,\stackrel{t}{\check{\dots}},\lambda_n)\lambda_{t-1}^{p+M}\bigg)
\\
&=-\sum_{t=1}^{n}\sum_{s=t}^{n}
\sum_{\substack{i,j,h\in I\\p\in\mb Z_{\geq0}}}(-1)^{s+t}
(K^{-1})_{ij}K_{hi_t}\mult_{(s+1,s+2)}\bigg(
(|_{x=\partial}x^{-M}u_j)\star_{s+1}
\\
&\bigg(\frac{\partial}{\partial u_h^{(p)}}
\bigg)_{(t)}A_{i_1,\stackrel{t}{\check{\dots}},i_{s},i,i_{s+1},\dots,i_n}
(\lambda_1,\stackrel{t}{\check{\dots}},\lambda_{s},x,\lambda_{s+1},\dots,\lambda_n)\lambda_t^{p+M}\bigg)
\\
&+k\,A_{i_1,\dots,i_n}
(\lambda_1,\dots,\lambda_n)
\\
&-\sum_{s=1}^{n+1}\sum_{t=s}^{n}
\sum_{\substack{i,j,h\in I\\p\in\mb Z_{\geq0}}}(-1)^{s+t}
(K^{-1})_{ij}K_{hi_{t}}\mult_{(s,s+1)}\bigg(
(|_{x=\partial}x^{-M}u_j)\star_s
\\
&\bigg(\frac{\partial}{\partial u_h^{(p)}}
\bigg)_{(t+1)}A_{i_1,\dots,i_{s-1},i,i_s,\stackrel{t+1}{\check{\dots}},i_n}
(\lambda_1,\dots,\lambda_{s-1},x,\lambda_s,\stackrel{t+1}{\check{\dots}},\lambda_n)\lambda_{t}^{p+M}\bigg)
\,.
\end{split}
\end{equation}
To derive the second equality above, we changed the order of summation in the first sum and shifted the index of summation $s$, we used \eqref{eq:Delta-times} and the fact that $A(\lambda_1,\dots,\lambda_n)\in\widetilde{\Omega}^n_{k,d}(\mc V)$ to rewrite the second sum,
and we changed the order of summation and shifted the index of summation $t$ in the third sum.
Equation \eqref{eq:homotopy2} follows immediately from \eqref{eq:AA1} and \eqref{B1}.
\end{proof}
\begin{remark}\label{20250913:rem1}
Since $L_{\Delta,M}$ is an even derivation of $\widetilde{\Omega}(\mc V)$ and $\iota_{\Delta,M}$ is an odd derivation of $\widetilde{\Omega}(\mc V)$ we have that the homotopy operator $h_M$ defined in \eqref{eq:hm} preserves the space $[\widetilde{\Omega}(\mc V),\widetilde{\Omega}(\mc V)]\cap\bigoplus_{n\geq1}\widetilde{\Omega}_{\bullet,\geq1}^n(\mc V)$.
\end{remark}
We denote by $B^{n}=\ker \tilde\delta|_{\widetilde{\Omega}^n(\mc V)}$, $n\in\mb Z_{\geq0}$, and by
$Z^n=\im \tilde\delta|_{\widetilde{\Omega}^{n-1}(\mc V)}$, $n\geq1$, so that we have
$\coH^0(\widetilde{\Omega}(\mc V),\tilde\delta)=B^0$ and $\coH^n(\widetilde{\Omega}(\mc V),\tilde\delta)=B^n/Z^n$, for $n\geq1$.
We are now ready to prove Theorem \ref{Thm:DPVcoh-Const}. 
\begin{proof}[Proof of Theorem \ref{Thm:DPVcoh-Const}]
Recall that $\widetilde{\Omega}^0(\mc V)=\mc V$. Let $f\in\mc V$. By Lemma \ref{Lem:DPVCoh},
$\tilde\delta(f)=0$ implies that $f\in\kk$. Hence, $\coH^0(\widetilde{\Omega}(\mc V),\tilde\delta)=B^0=\kk$.

Next, let us fix $n \geq 1$. We are going to show that $B^n=\widetilde{\Omega}^n_{0,0}\oplus Z^n$, where
\begin{equation} \label{Eq:DPVcons1}
\widetilde{\Omega}^n_{0,0}
=\{ (C_{\underline i}(\lambda_1,\dots,\lambda_n))_{\underline i\in I^n}\in\widetilde{\Omega}^n(\mc V)
\mid C_{\underline i}(\lambda_1,\dots,\lambda_n)\in U^{n}_{0}\text{ for every }\underline i\in I^n
\}
\,,
\end{equation}
is defined in \eqref{Eq:Biwei} and $U_0^n=\bigoplus_{\underline a\in T_n(M)} \kk \, \lambda^{\underline a}$ is defined in \eqref{dm-dec}.

Let $A=A(\lambda_1,\dots,\lambda_n)$ be as in \eqref{Eq:cochX}, and let us write
(in a unique way)
$A=X+Y$, where $X=X(\lambda_1,\dots,\lambda_n)
=(X_{\underline i}(\lambda_1,\dots,\lambda_n))_{\underline i\in I^n}$ with entries
$$
X_{\underline i}(\lambda_1,\dots,\lambda_n)=\sum_{\underline a\in T_{n}(M)}f_{\underline a}\lambda^{\underline a}
\,,
$$
and
$Y=Y(\lambda_1,\dots,\lambda_n)=(Y_{\underline i}(\lambda_1,\dots,\lambda_n))_{\underline i\in I^n}$ with entries
$$
Y_{\underline i}(\lambda_1,\dots,\lambda_n)=\sum_{\underline a\not\in T_{n}(M)}f_{\underline a}\lambda^{\underline a}
\,.
$$
Let us assume that $A\in B^n$. Then, by Lemma \ref{Lem:DPVCoh} we have that $X\in\widetilde{\Omega}^n_{0,0}$.
Moreover, since $\tilde\delta$ acts trivially on $\widetilde{\Omega}^n_{0,0}$, we have that
$$
\tilde\delta( Y)=\tilde\delta(A)-\tilde\delta(C)=0
\,.
$$
On the other hand, $Y\in\widetilde{\Omega}^n_{\bullet,\geq1}(\mc V)$. Then, by
\eqref{eq:homotopy} we have
$$
Y=\tilde\delta (h_M(Y))+h_M(\tilde\delta(Y))
=\tilde \delta(h_M(Y))
\in Z^n
\,,
$$
from which follows the direct sum decomposition $B^n=\widetilde{\Omega}^n_{0,0}\oplus Z^n$.
We have thus shown that
$$
\coH^n(\widetilde{\Omega}(\mc V),\tilde\delta)\simeq
\bigoplus_{\substack{\underline i\in I^n\\\underline a\in T_n(M)}}\kk\, \eta_{\underline i,\underline a}
\,\qquad n\in\mb Z_{\geq0}
\,,
$$
where $\eta_{\underline i,\underline a}=[(\delta_{\underline j,\underline i}\lambda^{\underline a})_{\underline j\in I^n}]\in B^n/Z^n$ denotes the cohomology class of
the array $(\delta_{\underline j,\underline i}\lambda^{\underline a})_{\underline j\in I^n}\in\widetilde{\Omega}^n_{0,0}$.
%
%The statement of the theorem follows from the isomorphism
%\eqref{20250723:eq2}.
\end{proof}
%

%%%
\section[Computation of the reduced \MakeLowercase{d}PVA cohomology]{Computation of the reduced double Poisson vertex algebra cohomology}\label{sec:PVAexa2}

Using the isomorphism \eqref{20250723:eq2}
and the results in Section \ref{sec:red-dPVA} we recall that we have a short exact sequence
of complexes
$$
0\rightarrow (\partial\widetilde{\Omega}(\mc V)
+[\widetilde{\Omega}(\mc V),\widetilde{\Omega}(\mc V)],\tilde\delta)
\stackrel{\alpha}{\longrightarrow} (\widetilde{\Omega}(\mc V),\tilde\delta)
\stackrel{\beta}{\longrightarrow} (\Omega(\mc V),\delta)\rightarrow 0
\,,
$$
where $\alpha$ denotes the inclusion map and $\beta=\tint$\glslink{tint}{} denotes the quotient map, see \eqref{eq:omega} (we are using the same notation for the quotient map as in Section \ref{sec:red-dPVA} since we can identify $\Gamma(\mc V)$ with $\Omega(\mc V)$ using \eqref{eq:gamma}; if we identify $\Omega(\mc V)$ with $\Sigma(\mc V)$ using Proposition \ref{20240819:prop1}, then $\beta=\widetilde\phi$ given in \eqref{eq:Phi_n}),
which leads to the following long exact sequence in cohomology (we denote $I(\mc V)=\partial\widetilde{\Omega}(\mc V)
+[\widetilde{\Omega}(\mc V),\widetilde{\Omega}(\mc V)]$)
\begin{equation}\label{eq:les2Omega}
% \text{this diagram has been commented to speed up compiling}
% \begin{comment}
\hspace{-1.5cm}
\begin{tikzcd}
0 \arrow[r]& \coH^0(I(\mc V),\tilde\delta)\arrow[r,"\alpha_0"]
& \coH^0(\widetilde{\Omega}(\mc V),\tilde\delta)\arrow[r,"\beta_0"] &
\coH^0(\Omega(\mc V),\delta) \arrow[dll,in=175,out=-5,"\gamma_0"']&
\\
&\coH^1(I(\mc V),\tilde\delta)\arrow[r,"\alpha_1"]&
\coH^1(\widetilde{\Omega}(\mc V),\tilde\delta) \arrow[r,"\beta_1"]
&\coH^1(\Omega(\mc V),\delta) \arrow[dll,in=175,out=-5,"\gamma_1"']&
\\
&\coH^2(I(\mc V),\tilde\delta)\arrow[r,"\alpha_2"]&
\coH^2(\widetilde{\Omega}(\mc V),\tilde\delta) \arrow[r,"\beta_2"]
&\coH^2(\Omega(\mc V),\delta)\arrow[r,"\gamma_2"] &\dots
\end{tikzcd}
% \end{comment}
\end{equation}
We want to use the long exact sequence \eqref{eq:les2Omega} to get information on
$\coH(\Omega(\mc V),\delta)$. Borrowing the terminology from the commutative case (cf. \cite{DSK13}), we refer to the complex $(\Omega(\mc V),\delta)$ as the \emph{generalized noncommutative variational complex} (for $M=0$ and $C=\Id_\ell$ it reduces to the noncommutative variational complex defined in \cite{DSKV}).

%%%
\subsection{Description of $\coH^n(I(\mc V),\tilde\delta)$}\label{sec:I(V)}
Recall that, by Theorem \ref{Thm:DPVcoh-Const}, for every $n\in\mb Z_{\geq0}$ we have a (canonical)
isomorphism $\coH^{n}(\widetilde{\Omega}(\mc V),\tilde\delta)\simeq \widetilde{\Omega}^n_{0,0}$.
Next, we describe $\coH^{n}(I(\mc V),\tilde\delta)$.
To this aim, we first need the results in Lemma \ref{20250812:lem1}.
Before stating and proving it let us introduce some notation. We set $C=[\widetilde{\Omega}(\mc V),\widetilde{\Omega}(\mc V)]$. Since $L_{\Delta,M}$ preserves $C$ we have a direct sum decomposition
(cf. \eqref{Eq:Biwei})
$$
C=\bigoplus_{k,n\in\mb Z_{\geq0}}\bigoplus_{d=0}^nC^n_{k,d}\,,
\quad
C^n_{k,d}=[\widetilde{\Omega}(\mc V),\widetilde{\Omega}(\mc V)]\cap\widetilde{\Omega}^n_{k,d}
\,.
$$
We also denote $C^\bullet_{0,0}=\bigoplus_{n\in\mb Z_{\geq0}}C^n_{0,0}\subset\widetilde{\Omega}^{\bullet}_{0,0}$.
It is not hard to check, using the definition of the commutator space $C$, that $C^0_{0,0}=0$ and that $C^{n}_{0,0}$, $n\geq1$, consists of the arrays $A(\lambda_1,\dots,\lambda_n)\in\widetilde{\Omega}^n_{0,0}$ whose entries satisfy
\begin{equation}\label{20250901:eqC}
\sum_{s=0}^{n-1}(-1)^{s(n-s)}A_{i_{\sigma^s(1)},\dots,i_{\sigma^s(n)}}(\lambda_{\sigma^s(1)},\dots,\lambda_{\sigma^s(n)})=0\,,
\end{equation}
for every $i_1,\dots,i_n\in I$. In particular, $C^1_{0,0}=0$.
Let us also fix, as in Section \ref{sec:rel}, a subspace $\widetilde{\Sigma}(\mc V)\subset\widetilde{\Omega}{(\mc V)}$
complementary to $C$, namely the direct sum decomposition \eqref{20250901:eq1} holds, such that $\overline{\phi}|_{\widetilde{\Sigma}(\mc V)}:\widetilde{\Sigma}(\mc V)\rightarrow\overline{\Sigma}(\mc V)$ is an isomorphism. Recall also that $\widetilde{\Sigma}(\mc V)$ is preserved by $\partial$. Moreover, it is straightforward to verify that $\widetilde{\phi}$ and $L_{\Delta,M}$ commute. Hence, we have the direct sum decomposition
$$
\widetilde{\Sigma}(\mc V)=\bigoplus_{k,n\in\mb Z_{\geq0}}\bigoplus_{d=0}^n\widetilde{\Sigma}^n_{k,d}
\,,
\quad
\widetilde{\Sigma}^n_{k,d}=\widetilde{\Sigma}(\mc V)\cap\widetilde{\Omega}^n_{k,d}
\,.
$$
Let us set $\widetilde{\Sigma}^\bullet_{0,0}=\bigoplus_{n\in\mb Z_{\geq0}}\widetilde{\Sigma}^n_{0,0}$. Then, we have the direct sum decomposition
\begin{equation}\label{20250901:eq4}
\widetilde{\Omega}^\bullet_{0,0}
=\widetilde{\Sigma}^\bullet_{0,0}\oplus C^\bullet_{0,0}
\,.
\end{equation}
Without loss of generality we may assume that $\widetilde{\Sigma}^0_{0,0}=\kk$,
and that $\widetilde{\Sigma}^{n}_{0,0}$, $n\geq1$, consists of the arrays $A(\lambda_1,\dots,\lambda_n)\in\widetilde{\Omega}^n_{0,0}$ whose entries satisfy
\begin{equation}\label{20250901:eqSigma}
A_{i_1,i_2,\dots,i_n}(\lambda_1,\lambda_2,\dots,\lambda_n)=
(-1)^{n+1}A_{i_2,\dots,i_n,i_1}(\lambda_2,\dots,\lambda_n,\lambda_1)\,,
\end{equation}
for every $i_1,\dots,i_n\in I$. In particular, $\widetilde{\Sigma}^1_{0,0}=\widetilde{\Omega}^1_{0,0}$.
\begin{lemma}\label{20250812:lem1}
\begin{enumerate}[(a)]
\item We have canonical isomorphisms
$$
\coH^n(\partial\widetilde{\Omega}(\mc V),\tilde\delta)\simeq\partial\widetilde{\Omega}^n_{0,0}\simeq
\left\{
\begin{array}{ll}
0\,, & n=0\,,
\\
\widetilde{\Omega}_{00}^n\,, & n\geq1\,.
\end{array}
\right.
$$
\item 
The inclusion $(C^\bullet_{0,0},0)\subset ([\widetilde{\Omega}(\mc V),\widetilde{\Omega}(\mc V)],\tilde\delta)$ induces a canonical isomorphism in cohomology
$$
\coH^n([\widetilde{\Omega}(\mc V),\widetilde{\Omega}(\mc V)],\tilde\delta)\simeq C^n_{0,0}\,,
$$
for every $n\geq0$.
\item
We have canonical isomorphisms
$$
\coH^n(\partial\widetilde{\Omega}(\mc V)\cap[\widetilde{\Omega}(\mc V),\widetilde{\Omega}(\mc V)],\tilde\delta)\simeq\partial C^n_{0,0}\simeq
\left\{
\begin{array}{ll}
0\,, & n=0\,,
\\
C_{00}^n\,, & n\geq1\,.
\end{array}
\right.
$$
\end{enumerate}
\end{lemma}
\begin{proof}
Since $\partial$ and $\tilde\delta$ commute we have
$$
H^{0}(\partial\widetilde{\Omega}(\mc V),\tilde\delta)
=\ker\tilde\delta|_{\partial\widetilde{\Omega}^0(\mc V)}
=\partial \ker\tilde\delta|_{\widetilde{\Omega}^0(\mc V)}
=\partial(\kk)=0
\,.
$$
For $n\geq1$, the isomorphism $\coH^n(\partial\widetilde{\Omega}(\mc V),\tilde\delta)\simeq\widetilde{\Omega}_{00}^n$, follows from Proposition \ref{20250812:prop1}
and the identification \eqref{20250723:eq2}. This proves part (a). For part (b) note that the inclusion $C^n_{0,0}\subset [\widetilde{\Omega}(\mc V),\widetilde{\Omega}(\mc V)]$ induces a linear map $C^n_{0,0}\rightarrow \coH^n([\widetilde{\Omega}(\mc V),\widetilde{\Omega}(\mc V)],\tilde\delta)$ sending  $A\in C_{00}^n$ to its cohomology class in $\coH^n([\widetilde{\Omega}(\mc V),\widetilde{\Omega}(\mc V)],\tilde\delta)$. This map is injective by Theorem \ref{Thm:DPVcoh-Const}. We are left to prove it is also surjective. Let us introduce the following notation: for $X=X(\lambda_1,\dots,\lambda_n)\in\widetilde{\Omega}^n(\mc V)$ we write
$$
X=\sum_{k\in\mb Z_{\geq0}}\sum_{d=0}^n X_{k,d}
\,,
\quad
X_{k,d}\in\widetilde{\Omega}_{k,d}^n
\,,
$$
for the decompositions of $X$ with respect to the direct sum \eqref{Eq:Biwei}.
Let $A=A(\lambda_1,\dots,\lambda_n)\in C\cap\widetilde{\Omega}^n(\mc V)$ and let us assume that
$\tilde\delta(A)=0$. Let $X=A-A_{00}$. Clearly $\tilde \delta X=0$ and $X_{00}=0$.
Hence, by Lemma \ref{Lem:DPVCoh} we have that $X\in C\cap \widetilde{\Omega}_{\bullet,\geq1}^n(\mc V)$. As in the proof of Theorem \ref{Thm:DPVcoh-Const} it follows that
$$
A-A_{00}=X=\tilde\delta(h_M(X))\in \tilde\delta([\widetilde{\Omega}(\mc V),\widetilde{\Omega}(\mc V)])
\,,
$$
which belongs to $\tilde\delta([\widetilde{\Omega}(\mc V),\widetilde{\Omega}(\mc V)])$
since $h_M(X)\in[\widetilde{\Omega}(\mc V),\widetilde{\Omega}(\mc V)]$
(cf. Remark \ref{20250913:rem1})
and $\tilde\delta$ preserves the commutator space $[\widetilde{\Omega}(\mc V),\widetilde{\Omega}(\mc V)]$ by Proposition \ref{20230802:prop1}(b).
This proves that the map $C^n_{0,0}\rightarrow \coH^n([\widetilde{\Omega}(\mc V),\widetilde{\Omega}(\mc V)],\tilde\delta)$ is surjective thus concluding the proof of part (b). For part (c),
since $\partial$ preserves the direct sum decomposition \eqref{20250901:eq1} we have
$$
\partial\widetilde{\Omega}(\mc V)\cap[\widetilde{\Omega}(\mc V),\widetilde{\Omega}(\mc V)]
=\partial[\widetilde{\Omega}(\mc V),\widetilde{\Omega}(\mc V)]
\,.
$$
Since $\partial$ and $\tilde\delta$ commute, it follows that
$$
H^{0}(\partial[\widetilde{\Omega}(\mc V),\widetilde{\Omega}(\mc V)],\tilde\delta)
=\ker\tilde\delta|_{\partial[\mc V,\mc V]}
=\partial \ker\tilde\delta|_{[\mc V,\mc V]}=0\,,
$$
where in the last equality we used part (b) and the fact that $C^0_{0,0}=0$.
For $n\geq1$, the isomorphism $\coH^n(\partial[\widetilde{\Omega}(\mc V),\widetilde{\Omega}(\mc V)],\tilde\delta)\simeq C_{00}^n$, follows from Proposition \ref{20250812:prop2}
and the identification \eqref{20250723:eq2}. This proves part (c).
\end{proof}
We are now able to provide a description of the cohomology of the subcomplex $(I(\mc V),\tilde\delta)\subset(\widetilde{\Omega}(\mc V),\tilde\delta)$.
\begin{proposition}\label{prop:I(V)}
We have canonical isomorphisms
$$
\coH^n(I(\mc V),\tilde\delta)\simeq\partial\widetilde{\Sigma}^n_{0,0}\oplus C^n_{0,0}\simeq
\left\{
\begin{array}{ll}
0\,, & n=0\,,
\\
\widetilde{\Omega}_{00}^n\,, & n\geq1\,.
\end{array}
\right.
$$
\end{proposition}
\begin{proof}
For brevity let us set $B=\partial\widetilde{\Omega}(\mc V)$ (recall that $C=[\widetilde{\Omega}(\mc V),\widetilde{\Omega}(\mc V)]$, so $\mc I(\mc V)=B+C$).
We consider the short exact sequence of complexes
\begin{equation}\label{20250901:eq2}
0\rightarrow (B\cap C,\tilde\delta)
\stackrel{\mu}{\longrightarrow} (B,\tilde\delta)\oplus(C,\tilde\delta)
\stackrel{\nu}{\longrightarrow} (B+C,\tilde\delta)\rightarrow 0
\,,
\end{equation}
where $\mu (X)=(X,-X)$, for $X\in B\cap C$ and $\nu(X,Y)=X+Y$, for $X\in B$ and $Y\in C$.
Recalling that $B\cap C=\partial C$, the short exact sequence \eqref{20250901:eq2} leads to the 
following long exact sequence in cohomology
\begin{equation}\label{eq:les4Omega}
% \text{this diagram has been commented to speed up compiling}
% \begin{comment}
\hspace{-1.5cm}
\begin{tikzcd}
0 \arrow[r]& \coH^0(\partial C,\tilde\delta)\arrow[r,"\mu_0"]
& \coH^0(B,\tilde\delta)\oplus\coH^0(C,\tilde\delta)\arrow[r,"\nu_0"] &
\coH^0(I(\mc V),\tilde \delta) \arrow[dll,in=175,out=-5,"\xi_0"']&
\\
&\coH^1(\partial C,\tilde\delta)\arrow[r,"\mu_1"]&
\coH^1(B,\tilde\delta)\oplus\coH^1(C,\tilde\delta)\arrow[r,"\nu_1"]
&\coH^1(I(\mc V),\tilde \delta) \arrow[dll,in=175,out=-5,"\xi_1"']&
\\
&\coH^2(\partial C,\tilde\delta)\arrow[r,"\mu_2"]&
\coH^2(B,\tilde\delta)\oplus\coH^2(C,\tilde\delta) \arrow[r,"\nu_2"]
&\coH^2(I(\mc V),\tilde \delta)\arrow[r,"\xi_2"] &\dots
\end{tikzcd}
% \end{comment}
\end{equation}
%By Lemma \ref{20250812:lem1} we have that $\coH^0(\partial C,\tilde\delta)=\coH^0(B,\tilde\delta)=
%\coH^0(C,\tilde\delta)=0$. Plugging this information in the long exact sequence \eqref{eq:les4Omega}
%yields $\coH^0(I(\mc V),\tilde\delta)=0$.
The maps $\mu_n$ in \eqref{eq:les4Omega}
induced by the map $\mu$ in \eqref{20250901:eq2} are as follows: for $X\in C^n$, the cohomology class $\partial X+\tilde\delta(\partial C)\in \coH^n(\partial C,\tilde\delta)$ is mapped to the pair of cohomology classes
$(\partial X+\tilde\delta(B),-\partial X+\tilde\delta (C))\in\coH^n(B,\tilde\delta)\oplus\coH^n(C,\tilde\delta)$. We use Lemma \ref{20250812:lem1}(c) to identify $\coH^n(\partial C,\tilde\delta)\simeq \partial C^n_{0,0}$, Lemma \ref{20250812:lem1}(a) to identify $\coH^n(B,\tilde\delta)\simeq \partial \widetilde{\Omega}^n_{0,0}$ and Lemma \ref{20250812:lem1}(b) to identify $\coH^n( C,\tilde\delta)\simeq C^n_{0,0}$. With these identifications the maps $\mu_n$ become
\begin{equation}\label{20250901:eq3}
\mu_n:\partial C^{n}_{0,0}\rightarrow \partial \widetilde{\Omega}^{n}_{0,0}\oplus C^n_{0,0}\,,
\quad
\partial X\mapsto (\partial X,-(\partial X)_{00})\,,
\end{equation}
where we are using the notation introduced in the proof of Lemma \ref{20250812:lem1}(b):
$(\partial X)_{00}$ is the unique element in $C^n_{0,0}$ such that $\partial X-(\partial X)_{00}\in\tilde\delta(C)$.
In particular, the maps $\mu_n$ are injective, from which follows that $\xi_n$ is trivial for every $n\geq0$. Using the exactness of the sequence \eqref{eq:les4Omega}, we then have that
$$
\nu_n:\partial\widetilde{\Omega}^{n}_{0,0}\oplus C^n_{0,0}\to \coH^n(I(\mc V),\tilde\delta)\,,
$$
is surjective for every $n\geq0$. Hence, using again the exactness of the long exact sequence we get
$$
\coH^n(I(\mc V),\tilde\delta)\simeq \partial\widetilde{\Omega}^{n}_{0,0}\oplus C^n_{0,0}/\ker \nu_n
=\partial\widetilde{\Omega}^{n}_{0,0}\oplus C^n_{0,0}/\im \mu_n
\,.
$$
Recalling the direct sum decomposition \eqref{20250901:eq4} and the definition of the maps $\mu_n$ given in \eqref{20250901:eq3} we can then factor the map $\nu_n$ through the surjective map $\pi_n$
(with $\ker\pi_n=\im\mu_n$)
$$
\begin{tikzcd}
\partial\widetilde{\Omega}^n_{0,0}\oplus C^n_{0,0} \arrow[r,"\nu_n"]\arrow[d,"\pi_n"']&
\coH^n(I(\mc V),\tilde\delta)\\
\partial\widetilde{\Sigma}^n_{0,0}\oplus C^n_{0,0}
\arrow[ur,bend right,"\exists!"']
\end{tikzcd}
$$
which sends $(\partial U,\partial V,Z)\in\partial\widetilde{\Sigma}^n_{0,0}\oplus\partial C^n_{0,0}\oplus C^n_{0,0}$ to $(\partial U,Z+(\partial V)_{0,0})\in\partial\widetilde{\Sigma}^n_{0,0}\oplus C^n_{0,0}$. This yields the
isomorphism 
$\coH^n(I(\mc V),\tilde\delta)\simeq \partial\widetilde{\Sigma}^{n}_{0,0}\oplus C^n_{0,0}$, for every $n\geq0$,
and concludes the proof.
\end{proof}
%
%%%
\subsection{Description of $\coH^n(\Omega(\mc V),\delta)$}\label{sec:Omega(V)}
In this section we want to use the long exact sequence \eqref{eq:les2Omega}, Theorem \ref{Thm:DPVcoh-Const} and Proposition \ref{prop:I(V)} to describe the cohomology of the complex $(\Omega(\mc V),\delta)$.
To this aim, we first provide an explicit description of the maps $\alpha_n$ in the long exact sequence \eqref{eq:les2Omega}.
\begin{lemma}\label{20250902:lem1}
\begin{enumerate}[(a)]
\item The map $\alpha_0$ is trivial.
\item For $n\geq1$, identifying $\coH^n(I(\mc V),\tilde\delta)\simeq\partial\widetilde{\Sigma}^n_{0,0}\oplus C^n_{0,0}\simeq\widetilde{\Omega}_{0,0}^n\simeq \coH^n(\widetilde{\Omega}(\mc V),\tilde\delta)$ as in Proposition \ref{prop:I(V)} and in Theorem \ref{Thm:DPVcoh-Const}, the map $\alpha_n$ induces
$\alpha_n\in\End(\widetilde{\Omega}_{0,0}^n)$ defined as follows.
For $P\in\Sigma^n_{0,0}$ and $Q\in C^n_{0,0}$, there exist $Y\in \widetilde{\Omega}^{n-1}(\mc V)$ and a unique $X\in\widetilde{\Omega}_{0,0}^n$ such that
$\partial P+Q=X+\tilde\delta(Y)$. Then,
\begin{equation}\label{eq:alphan}
\alpha_n(P+Q)=X\,.
\end{equation}
In particular, $\alpha_n$ restricts to the identity map on $C^{n}_{0,0}$.
\end{enumerate}
\end{lemma}
\begin{proof}
Part (a) follows from the fact that $\coH^0(I(\mc V),\tilde\delta)=0$.
Let us prove part (b). For $n\geq1$, the isomorphism $H^n(\widetilde{\Omega}(\mc V),\tilde\delta)\simeq\widetilde{\Omega}^n_{0,0}$ given by Theorem \ref{Thm:DPVcoh-Const} maps
$R+\tilde\delta(\widetilde{\Omega}^{n-1}(\mc V))$ to the unique element $X\in\widetilde{\Omega}_{0,0}^n$ such that $R-X\in\tilde\delta(\widetilde{\Omega}^{n-1}(\mc V))$ and the inverse map sends $X\in
\widetilde{\Omega}_{0,0}^n$ to $X+\tilde\delta(\widetilde{\Omega}^{n-1}(\mc V))$.
Similarly, the composition of isomorphisms $\coH^n(I(\mc V),\tilde\delta)\simeq\partial\widetilde{\Sigma}^n_{0,0}\oplus C^n_{0,0}\simeq\widetilde{\Omega}^n_{0,0}$
given by Proposition \ref{prop:I(V)} maps $R+\tilde\delta(I(\mc V)\cap \widetilde{\Omega}^{n-1}(\mc V))$
to the unique $X\in\widetilde{\Omega}^n_{0,0}$ such that $R-\partial X_1-X_2\in\tilde\delta(I(\mc V)\cap \widetilde{\Omega}^{n-1}(\mc V))$,
where $X=X_1+X_2$, $X_1\in\widetilde{\Sigma}^n_{0,0}$, $X_2\in C^n_{0,0}$, is the decomposition of $X$ 
with respect to the direct sum \eqref{20250901:eq4}, and the inverse map sends 
$P+Q\in\widetilde{\Omega}^n_{0,0}$, where $P\in\widetilde{\Sigma}^{n}_{0,0}$ and 
$Q\in C^{n}_{0,0}$ to $\partial P+Q+\tilde\delta(I(\mc V)\cap \widetilde{\Omega}^{n-1}(\mc V))$.
Equation \eqref{eq:alphan} follows by composing these isomorphisms with the map
$\alpha_n:\coH^n(I(\mc V),\tilde\delta)\rightarrow \coH^n(\widetilde{\Omega},\tilde\delta)$, induced
by the inclusion $I(\mc V)\subset \widetilde{\Omega}(\mc V)$, which sends 
$R+\tilde\delta(I(\mc V)\cap \widetilde{\Omega}^{n-1}(\mc V))\in\coH^n(I(\mc V),\tilde\delta)$ to 
$R+\tilde\delta(\widetilde{\Omega}^{n-1}(\mc V))\in\coH^n(\widetilde{\Omega}(\mc V),\tilde\delta)$.
\end{proof}
Using Lemma \ref{20250902:lem1} the long exact sequence \eqref{eq:les2Omega1} becomes
\begin{equation}\label{eq:les2Omega1}
%\text{this diagram has been commented to speed up compiling}
%\begin{comment}
\begin{tikzcd}
& 0\arrow[r,"\alpha_0"]
& \widetilde{\Omega}^0_{0,0}\arrow[r,"\beta_0"] &
\coH^0(\Omega(\mc V),\delta) \arrow[dll,in=175,out=-5,"\gamma_0"']&
\\
&\widetilde{\Omega}^1_{0,0}\arrow[r,"\alpha_1"]&
\widetilde{\Omega}^1_{0,0} \arrow[r,"\beta_1"]
&\coH^1(\Omega(\mc V),\delta) \arrow[dll,in=175,out=-5,"\gamma_1"']&
\\
&\widetilde{\Omega}^2_{0,0}\arrow[r,"\alpha_2"]&
\widetilde{\Omega}^2_{0,0} \arrow[r,"\beta_2"]
&\coH^2(\Omega(\mc V),\delta)\arrow[r,"\gamma_2"] &\dots
\end{tikzcd}
%\end{comment}
\end{equation}
As an immediate consequence we get the following result which agrees with the description of the noncommutative variational complex in \cite{DSKV}.
\begin{corollary}\label{M=0}
For $M=0$, we have 
$$
\coH^n(\Omega(\mc V),\delta)=\delta_{n0}\kk\,.
$$
\end{corollary}
\begin{proof}
For $M=0$ we have that $\widetilde{\Omega}^n_{0,0}=\delta_{n0}\kk$. The claim then follows from the long exact sequence \eqref{eq:les2Omega1}. 
\end{proof}
From now on we assume that $M\geq1$.
Recall that $C^1_{0,0}=0$ and $\widetilde{\Sigma}^1_{0,0}=\widetilde{\Omega}^1_{0,0}$. Let
$P=(P_{i}(\lambda))_{i\in I}\in\widetilde{\Omega}^1_{0,0}$ and write its entries as
$$
P_{i}(\lambda)=\sum_{k=0}^{M-1}p_{i,k}\lambda^k\in\kk[\lambda]
\,.
$$
Then,
\begin{equation}\label{20250902:eq1}
\lambda P_{i}(\lambda)=\sum_{k=1}^{M-1}p_{i,k-1}\lambda^k+p_{i,M-1}\lambda^M\in\kk[\lambda]
\,.
\end{equation}
From Theorem \ref{Thm:DPVcoh-Const} we have that the array $R=(p_{i,M-1}\lambda^M)_{i\in I}\in\tilde \delta(\mc V)$. In fact $R=\tilde\delta f$, where $f=\sum_{j,h\in I}p_{h,M-1}(K^{-1})_{hj}u_j$.
Hence, by \eqref{eq:alphan} we have that $\alpha_1(P)=X$, where
$X_i(\lambda)=\sum_{k=1}^{M-1}p_{i,k-1}\lambda^k$. Clearly, $\dim\ker\alpha_1=\ell$.
\begin{proposition}\label{cor:dim}
We have
$$
\dim H^{0}(\Omega(\mc V),\delta)=1+\ell\,,
$$
and
$$
\dim \coH^n(\Omega(\mc V),\delta)=\dim\ker\alpha_n+\dim\ker\alpha_{n+1}
\,,
$$
for every $n\geq1$. 
\end{proposition}
\begin{proof}
For every $n\geq0$ we have $\dim\coH^n(\Omega(\mc V),\delta)=\dim\ker\gamma_n+\dim\im\gamma_{n}$,
where $\gamma_n$ is the boundary map in the long exact sequence \eqref{eq:les2Omega}. By exactness of the sequence \eqref{eq:les2Omega} we have
that $\dim\im\gamma_n=\dim\ker\alpha_{n+1}$ and $\dim\ker\gamma_n=\dim\im\beta_n$.
From Lemma \ref{20250902:lem1}(a) and Theorem \ref{Thm:DPVcoh-Const} we have that
$\beta_0:\kk\rightarrow H^0(\Omega(\mc V),\delta)$ is injective. Hence,
we have
$$
\dim H^{0}(\Omega(\mc V),\delta)=1+\dim \ker\alpha_1=1+\ell\,.
$$
Again, by exactness of the sequence \eqref{eq:les2Omega}, we have
$\dim\im\beta_n=\dim\coH^n(\widetilde{\Omega}(\mc V),\tilde\delta)-\dim\ker\beta_n
=\dim\coH^n(\widetilde{\Omega}(\mc V),\tilde\delta)-\dim\im\alpha_n=\dim\ker\alpha_n
$, for every $n\geq1$.
Hence,
we conclude that
$$
\dim \coH^n(\Omega(\mc V),\delta)=\dim\ker\alpha_n+\dim\ker\alpha_{n+1}
\,. 
$$
\end{proof}
It is immediate to verify using \eqref{Eq:coch-dPX-1} that $\delta \tint 1=\tint \tilde\delta(1)=0$ and, similarly, that $\delta \tint u_i=\tint\tilde \delta(u_i)=0$, for every $i\in I$. Since these elements are linearly independent in $\mc V_{\sharp}$, by Proposition
\ref{cor:dim} we have that (cf. Subsection \ref{explicit:0})
$$
\coH^0(\Omega(\mc V),\delta)=\Cas(\mc V)
=\kk\tint 1\oplus\left(\bigoplus_{i\in I}\kk\tint u_i
\right)
\,.
$$
We want to find an analogous description for $\coH^n(\Omega(\mc V),\delta)$ for every $n\geq1$. To do so, we start by generalizing equation \eqref{20250902:eq1} to arbitrary $P\in\widetilde{\Sigma}^n_{0,0}$.
\begin{lemma}\label{20250902:lem2}
Let $P=(P_{\underline i}(\lambda_1,\dots,\lambda_n))_{\underline i\in I^n}\in\widetilde{\Sigma}^n_{0,0}$, $n\geq1$. Then, there exist $X=(X_{\underline i}(\lambda_1,\dots,\lambda_n))_{\underline i\in I^n}\in\widetilde{\Sigma}^n_{0,0}$ and $Y=(Y_{\underline i}(\lambda_1,\dots,\lambda_{n-1}))_{\underline i\in I^n}\in\widetilde{\Omega}^n_{0,0}$ such that
\begin{equation}\label{20250902:eq2}
\begin{split}
&(\lambda_1+\dots+\lambda_n) P_{\underline i}(\lambda_1,\dots,\lambda_n)
=X_{\underline i}(\lambda_1,\dots,\lambda_n)
\\
&+
\sum_{s=0}^{n-1}(-1)^{s(n-s)}Y_{i_{\sigma^s(1)},\dots,i_{\sigma^s(n)}}
(\lambda_{\sigma^s(2)},\dots,\lambda_{\sigma^s(n)})\lambda_{\sigma^s(1)}^M
\,.
\end{split}
\end{equation}
\end{lemma}
\begin{proof}
Recall that
$\widetilde{\Sigma}^n_{00}$ consists of arrays $P(\lambda_1,\dots,\lambda_n)=(P_{\underline i}(\lambda_1,\dots,\lambda_n))_{\underline i\in I^n}$ whose entries
$P_{\underline i}(\lambda_1,\dots,\lambda_n)\in\kk[\lambda_1,\dots,\lambda_n]$ are polynomials of degree at most
$N-1$ in each variable and satisfy \eqref{20250901:eqSigma}.
In particular, by repeatedly applying \eqref{20250901:eqSigma}, these entries satisfy
\begin{equation}\label{20250901:eqSigma-bis}
P_{i_1,\dots,i_n}(\lambda_1,\dots,\lambda_n)=
(-1)^{s(n-s)}P_{i_{\sigma^s(1)},\dots,i_{\sigma^s(n)}}(\lambda_{\sigma^s(1)},\dots,\lambda_{\sigma^s(n)})\,,
\end{equation}
for every $s=0,\dots,n-1$.
Let $P(\lambda_1,\dots,\lambda_n)=(P_{\underline i}(\lambda_1,\dots,\lambda_n))_{\underline i\in I^n}\in\widetilde{\Sigma}^n_{00}$. For every $\underline i\in I^n$ there exist a unique $X_{\underline i}(\lambda_1,\dots,\lambda_n)\in\kk[\lambda_1,\dots,\lambda_n]$ and unique $Y_{\underline i}^{(s)}(\lambda_1,\dots,\lambda_{n-1})\in\kk[\lambda_1,\dots,\lambda_{n-1}]$,
$s=0,\dots,n-1$, 
of degree at most
$N-1$ in each variable such that
\begin{equation}\label{20250812:eq1}
(\lambda_1+\dots+\lambda_n)P_{\underline i}(\lambda_1,\dots,\lambda_n)
=X_{\underline i}(\lambda_1,\dots,\lambda_n)
+\sum_{s=0}^n Y_{\underline i}^{(s)}(\lambda_1,\stackrel{s+1}{\check{\dots}},\lambda_n)\lambda_{s+1}^M
\,.
\end{equation}
Note that
\begin{equation}\label{20250903:eq1}
Y_{\underline i}^{(s)}(\lambda_1,\stackrel{s+1}{\check{\dots}},\lambda_n)
=\res_{\lambda_{s+1}}P_{\underline i}(\lambda_1,\dots,\lambda_n)\lambda_{s+1}^{-M}\,,
\end{equation}
where $\res_{x}p(x)$ denotes the coefficient of $x^{-1}$ in $p(x)\in\kk[x,x^{-1}]$.
Hence, using \eqref{20250812:eq1} and \eqref{20250903:eq1}, we have
\begin{align*}
Y_{\underline i}^{(s)}(\lambda_1,\stackrel{s+1}{\check{\dots}},\lambda_n)
&=(-1)^{s(n-s)}\res_{\lambda_{\sigma^s(1)}}P_{i_{\sigma^s(1)},\dots,i_{\sigma^s(n)}}(\lambda_{\sigma^s(1)},\dots,\lambda_{\sigma^s(n)})\lambda_{\sigma^s(1)}^{-M}
\\
&=(-1)^{s(n-s)}Y_{i_{\sigma^s(1)},\dots,i_{\sigma^s(n)}}^{(0)}(\lambda_{\sigma^s(2)},\dots,\lambda_{\sigma^s(n)})
\,,
\end{align*}
for every $s=1,\dots,n-1$. Since the array
$$
\left(
\sum_{s=0}^{n-1}(-1)^{s(n-s)}Y_{i_{\sigma^s(1)},\dots,i_{\sigma^s(n)}}^{(0)}
(\lambda_{\sigma^s(2)},\dots,\lambda_{\sigma^s(n)})\lambda_{\sigma^s(1)}^M
\right)_{i_1,\dots,i_n\in I}
$$
satisfies the skewsymmetry condition \eqref{20250901:eqSigma}, from \eqref{20250812:eq1}
we have that $X=(X_{\underline i}(\lambda_1,\dots,\lambda_n))\in\widetilde{\Sigma}^n_{0,0}$.
\end{proof}
Let $P=(P_{\underline i}(\lambda_1,\dots,\lambda_n))_{\underline i\in I^n}\in\widetilde{\Sigma}^n_{0,0}$ and decompose it as in \eqref{20250902:eq2}.
It follows using the same argument as in the proof of Theorem \ref{Thm:DPVcoh-Const} that
$$\left(
\sum_{s=0}^{n-1}(-1)^{s(n-s)}Y_{i_{\sigma^s(1)},\dots,i_{\sigma^s(n)}}
(\lambda_{\sigma^s(2)},\dots,\lambda_{\sigma^s(n)})\lambda_{\sigma^s(1)}^M
\right)_{i_1,\dots,i_n\in I}\in\tilde\delta(\widetilde{\Omega}^{n-1}(\mc V))\,.
$$
Hence, recalling the definition of the linear map $\alpha_n$ given in Lemma \ref{20250902:lem1}(b) we have that
$\alpha_n (P)=X$, where $X\in\widetilde{\Sigma}^n_{0,0}$ is the unique element defined by
\eqref{20250902:eq2}.

For every $n\geq1$, we introduce the subspace $V_n\subset\widetilde{\Omega}^n_{0,0}$ consisting of arrays $Y=(Y_{\underline i}(\lambda_1,\dots,\lambda_{n-1}))_{\underline i\in I^n}$ such that
\begin{equation}\label{20250902:eq3}
\sum_{s=0}^{n-1}(-1)^{s(n-s)}Y_{i_{\sigma^s(1)},\dots,i_{\sigma^s(n)}}
(\lambda_{\sigma^s(2)},\dots,\lambda_{\sigma^s(n)})\lambda_{\sigma^s(1)}^M
\in(\lambda_1+\dots+\lambda_n)\kk[\lambda_1,\dots,\lambda_n]
\,,
\end{equation}
for all indices $i_1,\dots,i_n\in I$. By Lemma \ref{20250902:lem2} and the fact that $\alpha_n$ restricts to the identity on $C^{n}_{0,0}$ we have a canonical isomorphism
$$
\phi_n:V_n\rightarrow \ker\alpha_n\,,
$$
mapping $Y\in V_n$ to the unique element $\phi_n(Y)=P\in \ker\alpha_n\subset\widetilde{\Sigma}^n_{0,0}$ such that
\begin{equation}\label{20250905:eq2}
\begin{split}
&\sum_{s=0}^{n-1}(-1)^{s(n-s)}Y_{i_{\sigma^s(1)},\dots,i_{\sigma^s(n)}}
(\lambda_{\sigma^s(2)},\dots,\lambda_{\sigma^s(n)})\lambda_{\sigma^s(1)}^M
\\
&=
(\lambda_1+\dots+\lambda_n) P_{\underline i}(\lambda_1,\dots,\lambda_n)
\,,
\end{split}
\end{equation}
for every $\underline i=(i_1,\dots,i_n)\in I^n$.
In the case $n=1$ illustrated by equation \eqref{20250902:eq1} we have that $\phi_1:V_1=\kk^\ell\rightarrow\ker\alpha_1$ maps the array $a=(a_1,\dots,a_\ell)\in\kk^\ell$
to the array $\phi_1(a)=(a_1\lambda_1^{M-1},\dots,a_\ell \lambda_1^{M-1})\in\ker\alpha_1$.

Furthermore we define the linear map $\chi_n:V_{n+1}\rightarrow \coH^n(\Omega(\mc V),\delta)$ by sending the array $Y=(Y_{i_0,i_1,\dots,i_n}(\lambda_1,\dots,\lambda_n))_{\underline i\in I^{n+1}}\in V_{n+1}$ to the cohomology class $\chi_n(Y)$ with
representative
$\widetilde{Y}=(\widetilde{Y}_{\underline i}(\lambda_1,\dots,\lambda_n))_{\underline i\in I^n}\in\widetilde{\Omega}^n(\mc V)$, where
\begin{equation}\label{20250905:eq1}
\begin{split}
&
\widetilde{Y}_{\underline i}(\lambda_1,\dots,\lambda_n)
\\
&=
\sum_{s=0}^{n}\sum_{j\in I}(-1)^{s(n-s)} (1^{\otimes s}\otimes w_j\otimes1^{\otimes (n-s)})Y_{j,i_{\sigma^{s}(1)},\dots,i_{\sigma^s(n)}}(\lambda_{\sigma^s(1)},\dots,\lambda_{\sigma^s(n)})
\,,
\end{split}
\end{equation}
and we set
\begin{equation}\label{wj}
w_j=\sum_{h\in I}(K^{-1})_{jh}u_h\in\mc V\,.
\end{equation}
The map $\chi_n$ is well defined since $\tint\widetilde{Y}\in\Omega^{n}(\mc V)$ is a closed element thus $\chi_n(Y)=[\tint \widetilde Y]\in\coH^{n}(\Omega(\mc V),\delta)$ is indeed a cohomology class.
In fact, by an explicit computation, using \eqref{Eq:coch-dPX-1}, we get that
$\tilde\delta(\widetilde Y)_{i_1,\dots,i_{n+1}}(\lambda_1,\dots,\lambda_{n+1})$
is given by the LHS of \eqref{20250905:eq2} with $n$ replaced by $n+1$. Hence, $\tilde\delta{\widetilde{Y}}=\partial P$, for some $P\in\widetilde{\Omega}^{n+1}(\mc V)$ and we have that
$\delta(\tint \widetilde{Y})=\tint\tilde\delta(\widetilde{Y})=0$.
\begin{proposition}\label{20250903:prop1}
For every $n\in\mb Z_{\geq0}$, we have that
$$
\phi_{n+1}=\gamma_n\circ \chi_n
\,.
$$
\end{proposition}
\begin{proof}
Let us recall the standard definition of the boundary map $\gamma_n$ in the long exact sequence
\eqref{eq:les2Omega}. Let $Y\in\Omega^n(\mc V)$ be such that $\delta (Y)=0$ and let us denote by $[Y]\in\coH^n(\Omega(\mc V),\delta)$ its cohomology class. Choose an element $\widetilde Y\in\widetilde{\Omega}^n(\mc V)$ such that $\tint \widetilde Y=Y$. Since $\delta (Y)=\delta(\tint \widetilde Y)=\tint \tilde \delta(\widetilde Y)=0$, we have that
$\tilde \delta (\widetilde Y)\in \widetilde{\Omega}(\mc V)+[\widetilde{\Omega}(\mc V),\widetilde{\Omega}(\mc V)]$. Hence, there exits $Q\in I(\mc V)\cap\widetilde{\Omega}^{n+1}(\mc V)$ such that $\tilde\delta(\widetilde Y)=Q$ (in particular, $\tilde\delta(Q)=0$). Then, $\gamma_n([Y])=[Q]\in\coH^{n+1}(I(\mc V),\tilde\delta)$. Using the identification
$\coH^n(I(\mc V),\tilde\delta)\simeq \partial\widetilde{\Sigma}^n_{0,0}\oplus C^n_{0,0}$ given in Proposition \ref{prop:I(V)} and the fact that $\im\gamma_n=\ker\alpha_{n+1}\subset\partial\widetilde{\Sigma}^{n+1}_{0,0}$, there exists a unique $P\in\widetilde{\Sigma}^{n+1}_{0,0}$ such that $Q-\partial P\in\tilde\delta(I(\mc V)\cap\widetilde{\Omega}^{n-1}(\mc V))$. 

Let $Y\in V_{n+1}$. Then, the array $\phi_{n+1}(Y)=P\in\ker\alpha_{n+1}$ is defined by \eqref{20250905:eq2}  (where $n$ is replaced by $n+1$).
On the other hand, a representative for the cohomology class of $\chi_n(Y)\in\coH^n(\Omega(\mc V),\delta)$ is the array $X\in\widetilde{\Omega}^n(\mc V)$ given by \eqref{20250905:eq1}. Let
$\gamma_n(\chi_n(Y))=P_1\in\ker\alpha_{n+1}$. By the previous considerations $P_1$ is the unique array such that $\tilde\delta(X)-\partial P_1\in\tilde\delta(I(\mc V)\cap\widetilde{\Omega}^{n-1}(\mc V))$. Using \eqref{Eq:coch-dPX-1} and the definition of $P$ given by \eqref{20250905:eq2} (where $n$ is replaced by $n+1$) it is straightforward to verify that
$\tilde\delta(X)-\partial P=0$. Hence, by uniqueness, we have $P=P_1$.
\end{proof}
\begin{theorem} \label{Thm:HOmega-Const}
For every $n\geq0$ we have the direct sum decomposition
$$
\coH^n(\Omega(\mc V),\delta)=\im\beta_n\oplus \chi_n(V_{n+1})
\,.
$$
In particular, $\coH^n(\Omega(\mc V),\delta)$ is spanned by the cohomology classes with representative $P\in\widetilde{\Sigma}^n_{00}\subset\widetilde{\Omega}(\mc V)$ and
the cohomology classes $\chi_n(X)$, for $X\in V_{n+1}$.
\end{theorem}
\begin{proof}
Since $\phi_{n+1}:V_{n+1}\rightarrow\ker\alpha_{n+1}$ is an isomorphism, by Proposition \ref{20250903:prop1} we have that
$\chi_n:V_{n+1}\to \coH^n(\Omega(\mc V),\delta)$ is injective and $\chi_n(V_{n+1})\cap\ker\gamma_n=0$. The claim then follows by a dimensionality argument based on Proposition \ref{cor:dim}: in fact, $\dim\chi_{n}(V_{n+1})=\dim V_{n+1}=\dim\ker\alpha_{n+1}$ by injectivity and $\dim\im\beta_n=\dim\ker\alpha_n$ by exactness of the sequence \eqref{eq:les2Omega}.   
\end{proof}
We conclude this section by outlining a strategy to compute the dimension of the kernel of the maps $\alpha_n$ in the long exact sequence \eqref{eq:les2Omega}, hence the dimension of $\coH^n(\Omega(\mc V),\delta)$ (recall Proposition \ref{cor:dim}). As an upshot it turns out that the cohomology spaces $\coH^{n}(\Omega(\mc V),\delta)$ are always nontrivial (while in the analogous situation in the commutative case, studied in \cite{DSK13}, these space are trivial for $n$ large enough).

For $n,M\geq1$, let
$$A_{n,M}=\kk[\lambda_1,\dots,\lambda_n]/(\lambda_1^M,\dots,\lambda_n^M)\,.$$
Recall that the space
$\widetilde{\Omega}^n_{0,0}$ consists of arrays $P(\lambda_1,\dots,\lambda_n)=(P_{\underline i}(\lambda_1,\dots,\lambda_n))_{\underline i\in I^n}$ whose entries
$P_{\underline i}(\lambda_1,\dots,\lambda_n)\in\kk[\lambda_1,\dots,\lambda_n]$ are polynomials of degree at most
$N-1$ in each variable. We can then identify
each entry with the corresponding coset in $A_{n,M}$. Hence, $\widetilde{\Omega}^n_{0,0}\simeq A_{n,M}^{\ell^n}$.
Let $P(\lambda_1,\dots,\lambda_n)=(A_{\underline i}(\lambda_1,\dots,\lambda_n))_{\underline i\in I^n}\in\widetilde{\Sigma}^n_{0,0}$. Recall from Lemma \ref{20250902:lem2} that $\alpha_n(P)=X\in\widetilde{\Sigma}_{0,0}^n$ is the array defined by equation \eqref{20250902:eq2}.
In terms of the identification 
$\widetilde{\Omega}^n_{0,0}\simeq A_{n,M}^{\ell^n}$, equation \eqref{20250902:eq2} allows us to rewrite (for $P\in\widetilde{\Sigma}_{0,0}^n$)
\begin{equation}\label{eq:alphan2}
\alpha_n(P)=(E_{n,M}(P_{\underline i}))_{\underline i\in I^n}\,,
\end{equation}
where
$$
E_{n,M}:A_{n,M}\rightarrow A_{n,M}\,,
\quad p(\lambda_1,\dots,\lambda_n)\mapsto (\lambda_1+\dots+ \lambda_n)p(\lambda_1,\dots,\lambda_n)
\,.
$$
We denote by $h_{n,M}=\dim\ker E_{n,M}$, for every $n,M\geq1$.
\begin{proposition}\label{prop:sl2}
We have that
$$
h_{n,M}
=\text{coefficient of }t^{[\frac{n(M-1)}{2}]}\text{ in } (1+t+\dots+t^{M-1})^n
\,.\
$$
\end{proposition}
\begin{proof}
Let us denote by $V=A_{n,M}$ and write
$$
V=\bigoplus_{k=0}^{n(M-1)} V[k]
\,,
$$
where $V[k]$ is the subspace of homogeneous polynomials of degree $k$. The Hilbert-Poincaré function of $V$ (with respect to the polynomial grading) is
$$
p(t)=\left(\frac{1-t^M}{1-t}\right)^n=(1+t+\dots+t^{M-1})^n
\,.
$$
To prove the claim it suffices to show that $h_{n,M}=\dim V[[\frac{n(M-1)}{2}]]$.
Since $p(t)=t^{n(M-1)}p(\frac1t)$, we have that
\begin{equation}\label{20250808:eq1}
\dim V[k]=\dim V[n(M-1)-k]
\,,\quad
k=0,\dots,n(M-1)\,.
\end{equation}
Let $E=E_{n,M}\in\End(V)$
and
$$
H=2\sum_{i=1}^n\lambda_i\frac{d}{d\lambda_i}-n(M-1)\Id_V\in \End(V)
\,.
$$
Note that $E$ is nilpotent, so it belongs to $\mf{sl}(V)$.
The operator $H$ is semisimple and
we have the $H$-eigenspace decomposition
$$
V=\bigoplus_{k=-n(M-1)}^{n(M-1)}V_k
\,,
\quad V_k=\{v\in V\mid Hv=kv\}
\,,
$$
where $V_k=V[\frac{k+n(M-1)}{2}]\subset V$ is the subspace of homogeneous polynomials of degree $\frac{k+n(M-1)}{2}$ (which could be zero if $\frac{k+n(M-1)}{2}\not\in\mb Z_{\geq0}$). By \eqref{20250808:eq1} we have that $H\in \mf{sl}(V)$ and by a direct computation we get
$[H,E]=2E$. Hence, by the Jacobson-Morozov Theorem, we have an $\mf{sl}(2)$-triple $\{E,H,F\}\subset \mf{sl}(V)$. By representation theory of $\mf{sl}(2)$ it follows that
$$
h_{n,M}=\dim\ker E=\dim V_0+\dim V_1
=\dim V[[\frac{n(M-1)}{2}]]\,.
$$
\end{proof}
For example, we have $h_{1,M}=1$, $h_{2,M}=M$ ,
$h_{3,M}=\frac{3M^2}{4}$, for $M$ even, and $h_{3,M}=\frac{3M^2+1}{4}$, for $M$ odd, and $h_{4,M}=\frac{M(2M^2+1)}{3}$. However, we do not know a closed formula for $h_{n,M}$ for arbitrary $n,M\geq1$.

Let us denote by $B_{n,M}\subset A_{n,M}$ the subspace consisting of cosets whose representative is
a polynomial $p(\lambda_1,\dots,\lambda_n)\in\kk[\lambda_1,\dots,\lambda_{n}]$ of degree at most
$M-1$ in each variable satisfying the condition
\begin{equation}\label{20250905:eq3}
p(\lambda_1,\lambda_2,\dots,\lambda_n)=(-1)^{n+1}p(\lambda_2,\dots,\lambda_n,\lambda_1)
\,.
\end{equation}
Then, clearly $E_{n,M}(B_{n,M})\subset B_{n,M}$. Moreover, the operator $H$ in the proof of
Proposition \ref{prop:sl2} still acts diagonally on $B_{n,M}$, hence $B_{n,M}\subset A_{n,M}$ is a subrepresentation of $\mf{sl}_2$. We can then apply the same argument
in the proof of Proposition \ref{prop:sl2} to derive that
$$
s_{n,M}=\dim\ker E_{n,M}|_{B_{n,M}}=
\dim \left(A_{n,M}[[\frac{n(M-1)}{2}]]\cap B_{n,M}\right)\,,
$$
namely it is the number of homogeneous polynomials of degree $[\frac{n(M-1)}{2}]$ in $A_{n,M}$
satisfying the skewsymmetry condition \eqref{20250905:eq3}. We were not able to determine the 
Hilbert-Poincar\'e function of $B_{n,M}$, however it is straightforward to verify that the first few values of $s_{n,M}$ are $s_{1,M}=\ell$, $s_{2,M}=[\frac{h_{2,M}}{2}]$ and $s_{3,M}=[\frac{h_{3,M}+2}{3}]$.

Let us also denote by $o_n=o_n(\ell)$ the number of orbits of the action of the cyclic group $G_n=\langle \sigma\rangle$ of order $n$ on $I^n$ given by ($(i_1,\dots,i_n)\in I^n$)
\begin{equation}\label{eq:cyclic-action}
\sigma\cdot (i_1,i_2,\dots, i_n)=(i_{2},\dots,i_n,i_1)
\,.
\end{equation}
Using Burnside's Lemma one finds
$$
o_n
=\frac{1}{n}\sum_{d|n}\varphi(d)\ell^{\frac n d}\,,
$$
where $\varphi(d)$ is the Euler totient function. The polynomials (in the variable $\ell$) $o_n$ are known as \emph{necklace polynomials} \cite{M}.
\begin{proposition}\label{20250905:prop1}
For every $n,M\geq1$ we have
$$
\dim\ker\alpha_{n}=(o_n-\ell)h_{n,M}+\ell s_{n,M}
\,.
$$
\end{proposition}
\begin{proof}
Let $P\in\widetilde{\Sigma}^{n}_{0,0}$. Note that 
$P_{i,i,\dots,i}\in B_{n,M}$ (cf. \eqref{20250901:eqC}) for every $i\in I$, and there are $s_{n,M}$ polynomials
such that $E_{n,M}(P_{i,i,\dots,i})=0$. This gives a contribution of $\ell s_{n,M}$ linearly independent elements in $\ker\alpha_n$. If $(i_1,\dots,i_n)\neq(i,i,\dots,i)$, $i\in I$, then its orbit
under the action of the group $G_n$ contains more than one element and if $E_{n,M}(P_{i_1,\dots,i_n})=0$, then we have also $E_{n,M}(P_{i_{\sigma^s(1)},\dots,i_\sigma^s(n)})=0$, for every $s=0,\dots,n-1$, by the skewsymmetry property \eqref{20250901:eqC}. Since $P_{i_1,\dots,i_n}\in A_{n,M}$, we have also a contribution of $(o_n-\ell)h_{n,M}$ linearly independent polynomials in $\ker\alpha_n$. This proves the claim.
\end{proof}

%%%%%

\subsection{Explicit example for $M=1$}
Let us describe explicitly the case $M=1$. First, in this case we have $\widetilde{\Omega}^n_{0,0}=\kk^{\ell^n}$, $n\geq0$, and
$$
\widetilde{\Sigma}^{n}_{0,0}=\{(a_{\underline i})_{\underline i\in I^n}\in\kk^{\ell^n}
\mid a_{i_1,i_2,\dots,i_n}=(-1)^{n+1}a_{i_2,\dots,i_n,i_1}\text{ for every } i_1,\dots,i_n\in I\}
\,.
$$
Moreover, from equations \eqref{20250902:eq2} and \eqref{20250905:eq2} we immediately get that $\widetilde{\Sigma}^n_{0,0}=\ker\alpha_n=V_n$. For brevity we use $V_n$ to denote this space.
From Proposition \ref{prop:sl2} we have that $h_{n,1}=1$. Moreover, we have that
$s_{n,1}=1$ if $n$ is odd, otherwise $s_{n,1}=0$. Hence, from Proposition \ref{20250905:prop1} we get
\begin{equation}\label{20250905:eq4}
\dim V_n=o_n-\ell+\delta_{n\equiv1(2)}\ell\,,
\end{equation}
for every $n\geq1$.
By Proposition \ref{cor:dim} and equation \eqref{20250905:eq4} we then get
$$
\dim\coH^{n}(\Omega(\mc V),\delta)=o_n+o_{n+1}-\ell\,,
\quad
n\geq1\,.
$$
Recall that the maps $\beta_n$ in the long exact sequence \eqref{eq:les2Omega} are induced by the quotient map $\tint$. Identifying $\widetilde{\Omega}^n_{0,0}$ and $\coH^n(\widetilde{\Omega}(\mc V),\tilde\delta)$ using Theorem \ref{Thm:DPVcoh-Const}, we have that $\beta_n$ maps an array
$a\in\widetilde{\Omega}^n_{0,0}$ to the cohomology class in $\coH^{n}(\Omega(\mc V),\delta)$
whose representative is the projection of $a$ onto $\widetilde{\Sigma}^n_{0,0}$ with respect to the direct sum decomposition \eqref{20250901:eq4}. From Theorem \ref{Thm:HOmega-Const} we then get that
$\coH^n(\Omega(\mc V),\delta)$ is linearly spanned by the cohomology classes with representative
$a\in V_n$ and by the cohomology classes with representative (see \eqref{20250905:eq1})
\begin{align*}
&X_b(\lambda_1,\dots,\lambda_n)
\\
&=\big(
\sum_{s=0}^{n}\sum_{j\in I}(-1)^{s(n-s)} b_{j,i_{\sigma^{s}(1)},\dots,i_{\sigma^s(n)}}
(1^{\otimes s}\otimes w_j\otimes1^{\otimes (n-s)})\big)_{i_1,\dots,i_n\in I}\in\widetilde{\Omega}^n(\mc V)
\,,
\end{align*}
for every $b\in V_{n+1}$.

Recall from Subsection \ref{explicit:2} (using the isomorphisms \eqref{20250723:eq4} and \eqref{20250723:eq2}) that $\coH^1(\Omega(\mc V),\delta)$ parametrizes the equivalence classes
of dPVA-derivations of $\mc V$ up to inner dPVA-derivations. Using the above description of $\coH^n(\Omega(\mc V),\delta)$, and the identification of $\Omega(\mc V)$ with the space $\Sigma(\mc V)$ given by Proposition \ref{20240819:prop1} (recall also Proposition \ref{20240819:prop2}), we can pick as representatives for these equivalence classes the derivations $D_a\in\Vect(\mc V)^\partial$, $a\in\kk^\ell$, defined on generators by
$$
D_a(u_i)=a_i\,,
\quad i\in I
\,,
$$
and the derivations $D_b\in\Vect(\mc V)^\partial$, $b\in V_{2}\simeq \mf{so}_\ell$, defined on generators by
$$
D_b(u_i)=\sum_{j\in I}b_{ji}w_j
\,,
\quad
i\in I\,,
$$
where $w_j$ is given by \eqref{wj}
Hence, $\dPVH^1(\mc V)\simeq \kk^\ell\oplus \mf{so}_\ell$. In a similar way we can describe non-equivalent first-order formal deformations of the dPVA $\mc V$ (with $2$-fold $\lambda$-bracket given by \eqref{Eq:uu-odd-full}) since they are parametrized by $\coH^2(\Omega(\mc V),\delta)$, see Subsection \ref{explicit:2}. For $a\in V_2$ and $b\in V_3$, the associated first-order formal deformation is the $2$-fold $\lambda$-bracket $\ldb-_\lambda-\rdb_{\epsilon}^{a,b}$ defined on generators by
\begin{equation}\label{20250906:eq1}
\ldb{u_i}_\lambda u_j\rdb_{\epsilon}^{a,b}
=K_{ji}(1\otimes1)\lambda+\epsilon a_{ji}(1\otimes1)
+\epsilon\sum_{h\in I}\left(
b_{hji}w_h\otimes 1-1\otimes b_{hij}w_h\right)
\,,
\end{equation}
for every $i,j\in I$ (recall from Subsection \ref{explicit:2} that the $2$-fold $\lambda$-bracket \eqref{20250906:eq1} is skewsymmetric but it only satisfies the Jacobi identity up to order $\epsilon$). Let $\mc A=\bigoplus_{i\in I}\kk u_i\subset\mc V$ be the generating space
(as a differential algebra) of $\mc V$. We define a symmetric bilinear form
$(\cdot|\cdot)$ on $\mc A$ by setting $(u_i|u_j)=K_{ji}$, a skewsymmetric bilinear form
$\langle\cdot|\cdot\rangle_a$ on $\mc A$ by setting
$\langle u_i|u_j\rangle_a =a_{ji}$, and we define a bilinear product $\circ_b$ on $\mc A$ by setting
$$
u_i\circ_b u_j=\sum_{h\in I}b_{hji}w_h
\,,
\quad i,j\in I\,.
$$
By linearity, we rewrite \eqref{20250906:eq1} for every $x,y\in\mc A$ as follows
\begin{equation}\label{20250906:eq2}
\ldb x_\lambda y\rdb_{\epsilon}^{a,b}
=(x|y)(1\otimes1)\lambda+\epsilon \left(
\langle x|y\rangle_a(1\otimes1)
+x\circ_b y\otimes 1-1\otimes y\circ_b x\right)
\end{equation}
Note that using \eqref{wj} and the fact that $b\in V_3$ we have 
$$
(u_i\circ_b u_j|u_k)
=b_{kji}
=b_{ikj}
=(u_i|u_j\circ_b u_k)
\,,
$$
for every $i,j,k\in I$. Hence, the bilinear form $(\cdot|\cdot)$ on $\mc A$ is a trace form for the bilinear product $\circ_b$ on $\mc A$. Let $\Bil_K(\mc A)$ denote the space of bilinear products on $\mc A$ for which the symmetric bilinear form on $\mc A$ defined by $K$ is a trace form. From \eqref{20250906:eq1}
we have that non-equivalent first-order formal deformations of $\mc V$ are parametrized by the space
$\dPVH^2(\mc V)\simeq \mf{so}_\ell\oplus\Bil_K(\mc A)$.
\begin{remark}
Let $\mc A=\bigoplus_{i\in I}\kk u_i\subset\mc V$ be the generating space
of the algebra of differential polynomials in $\ell$ variables $u_1,\dots,u_\ell$. Let $(\cdot|\cdot)$ and $\langle\cdot|\cdot\rangle$ be bilinear forms on $\mc A$, and let $\circ$ be a bilinear product on $\mc A$.
We define the following $2$-fold $\lambda$-bracket for elements $x,y\in\mc A$
(see \eqref{20250906:eq2})
\begin{equation}\label{20250906:eq3}
\ldb x_\lambda y\rdb
=\langle x|y\rangle(1\otimes1)+ (x\circ y)\otimes 1-1\otimes (y\circ x)
+(x|y)(1\otimes 1)\lambda
\,,
\end{equation}
and we extend it to $\mc V$ using the Master Formula \eqref{master-infinite}.
Since
$$
\ldb y_{-\lambda-\partial} x\rdb^\sigma
=\langle y|x\rangle(1\otimes1)
+1\otimes (y\circ x)-(x\circ y)\otimes 1
-(y|x)(1\otimes 1)\lambda
\,,
$$
skewsymmetry \eqref{skew-gen-b} holds if and only if  $\langle\cdot|\cdot\rangle$ is skewsymmetric and $(\cdot|\cdot)$ is symmetric.
Moreover, by a straightforward computation we get
\begin{align*}
\ldb x_\lambda \ldb y_\mu z\rdb\rdb_L
&=\langle x|y\circ z\rangle(1\otimes1\otimes 1)
+
(x\circ (y\circ z))\otimes1\otimes 1
\\
&-1\otimes ((y\circ z)\circ x)\otimes 1
+ (x|y\circ z)(1\otimes1\otimes1)\lambda
\,;
\\
\ldb y_\mu \ldb x_\lambda z\rdb\rdb_R
&=
-\langle y|z\circ x\rangle(1\otimes1\otimes1)
-1\otimes (y\circ (z\circ x))\otimes1
\\
&+1\otimes1\otimes((z\circ x)\circ y)
-(y|z\circ x)(1\otimes1\otimes1)\mu
\,;
\\
\ldb \ldb x_\lambda y\rdb _{\lambda+\mu}z\rdb_L
&=\langle x\circ y|z\rangle(1\otimes1\otimes1)+
((x\circ y)\circ z)\otimes1\otimes 1
\\
&-1\otimes 1\otimes (z\circ (x\circ y))
+(x\circ y|z)(1\otimes1\otimes1)(\lambda+\mu) 
\,,
\end{align*}
for every $x,y,z\in\mc A$. Hence, assuming that 
$\langle\cdot|\cdot\rangle$ is skewsymmetric and 
$(\cdot|\cdot)$ is symmetric, the $2$-fold $\lambda$-bracket \eqref{20250906:eq3} satisfies the Jacobi identity \eqref{jacobi-gen-b} if and only if 
\begin{equation}\label{eq:cond}
x\circ (y\circ z)=(x\circ y)\circ z
\,,
\quad
(x|y\circ z)=(x\circ y|z)
\,,
\quad
\langle x| y\circ z\rangle
+\langle y| z\circ x\rangle
+\langle z| x\circ y\rangle=0
\,.
\end{equation}
for every $x,y,z\in\mc A$. We can thus associate, using \eqref{20250906:eq3}, a dPVA $\mc V=\mc V(\mc A)$ to any associative algebra $\mc A$ with a trace form and a skewsymmetric bilinear form satisfying the third condition in  \eqref{eq:cond}. For example, \eqref{20250906:eq3} extends the dPA associated to the path algebra of a quiver in Section \ref{Sec:Quiv}. This dPVA is a noncommutative analogue of the \emph{affine PVA} $\mc V(\mf g)$ associated to a Lie algebra $\mf g$ and used to define classical affine $W$-algebras, see \cite{DSKV13}. We will study the dPVA cohomology of $\mc V(\mc A)$ in a subsequent work.

\end{remark}

%

%%%%%%%%%%% NEW CHAPTER %%%%%%%%%%%%%%%
%%%%%%%%%%% NEW CHAPTER %%%%%%%%%%%%%%%
%%%%%%%%%%% NEW CHAPTER %%%%%%%%%%%%%%%
%%%%%%%%%%% NEW CHAPTER %%%%%%%%%%%%%%%

\chapter[Variational \MakeLowercase{d}PVA cohomology \& representation spaces]{Variational double Poisson vertex algebra cohomology and representation spaces}
\label{Ch:repVardPVA}

In this chapter we want to generalize the results of Chapter \ref{CH:rep-dPA} to the framework
of dPVA. In particular, we want to relate the variational dPVA cohomology defined in Section \ref{sec:var-dPVA}
with the variational Poisson cohomology defined in Section \ref{sec:var-PVA}.

%%%
\section{Double Poisson vertex algebras and representation spaces}\label{sec:dPVAtoPVA-1}
Let $\mc V$ be a differential algebra. Recall from Section \ref{ss:Rep-Not} the 
commutative algebra $\mc V_N$, $N\geq1$. It is the commutative algebra generated by symbols $a_{ij}$, $a\in\mc V$ and $1\leq i,j\leq N$, subject to the relations
($\alpha\in\kk$, $a,b\in\mc V$):
\begin{equation}\label{20130917:eq2}
(\alpha a)_{ij}=\alpha a_{ij}\,,
\qquad
(a+b)_{ij}=a_{ij}+b_{ij}\,,
\qquad
(ab)_{ij}=\sum_{k=1}^Na_{ik}b_{kj}\,.
\end{equation}
We make $\mc V_N$ a commutative differential algebra by defining the derivation (which we still denote by $\partial$) 
\begin{equation}\label{eq:partial-VN}
\partial( a_{ij})=(\partial a)_{ij}\,.
\end{equation}
If $\mc V$ is a dPVA with $2$-fold $\lambda$-bracket
$\ldb -_{\lambda}-\rdb$, it is shown in \cite{DSKV}, following the seminal work \cite{VdB1}, that $\mc V_N$ is a PVA
with $\lambda$-bracket $\{-_{\lambda}-\}$ defined by
(using Sweedler's notation)
\begin{equation}\label{20240829:eq1}
\{a_{ij}{}_\lambda b_{hk}\}=(\ldb a_{\lambda}b\rdb')_{hj}(\ldb a_{\lambda} b\rdb'')_{ik}
\,,
\end{equation}
for every $a_{ij},b_{hk}\in\mc V_N$, and extended to $\mc V_N$ by sesquilinearity \eqref{sesquiLCA},
and the Leibniz rules \eqref{lleibnizPVA}, \eqref{rleibnizPVA}. The notation in \eqref{20240829:eq1}
has the following meaning: if we write
$$
\ldb a_\lambda b\rdb=\sum_{n\in\mb Z_{\geq0}} (a_nb)'\otimes (a_nb)''\lambda^n
\,,
$$
then
$$
\{a_{ij}{}_\lambda b_{hk}\}=\sum_{n\in\mb Z_{\geq0}}(a_n b)'_{hj}(a_n b)''_{ik}\lambda^n
\,.
$$

%%%
\section{Motivational interlude}
\label{sec:dPVAtoPVA-mot}
Out first goal is to define a linear map
$\tr:C(\mc V)\to C(\mc V_N)$,
where $C(\mc V)$ denotes the space of $n$-fold $\lambda$-brackets on $\mc V$ defined in Section \ref{sec:compl-dPVA-1} and $C(\mc V_N)$ denotes the space of poly-$\lambda$-brackets on $\mc V_N$ defined in Section \ref{sec:var-PVA}, 
such that $\tr(C^n(\mc V))\subset C^n(\mc V)$ for every $n\in\mb Z_{\geq0}$. As a motivation for the construction of this map, let us illustrate the cases $n=0,1,2$ first. The Bourbakist reader can decide to skip this section without any harm.

We start by considering the case $n=0$: recall the definition of the trace map $\tr:\mc V\to\mc V_N$ from Section \ref{ss:Rep-Not}
given by
\begin{equation}\label{eq:trace0}
\tr (f)=\sum_{i=1}^Nf_{ii}\in\mc V_N
\,,\quad f\in\mc V
\,.
\end{equation}
\begin{lemma}\label{20240829:lem1}
\begin{enumerate}[(a)]
\item For every $f\in\mc V$ we have $\tr(\partial f)=\partial(\tr (f))$ (note that the $\partial$ symbol in the LHS denotes the derivation of $\mc V$, while the $\partial$ symbol in the RHS denotes the derivation of $\mc V_N$).
\item For every $f,g\in\mc V$ we have $\tr(fg)=\tr(gf)$.
\item
We have a well defined linear map
$$
\tr:C^0(\mc V)=\mc V_\sharp\to C^0(\mc V_N)=\mc V_N/\partial \mc V_N
$$
given by $\tr(\tint f)=\tint(\tr f)$ (note that the $\tint$ symbol in the LHS denotes the projection map $\mc V\to\mc V_\sharp$, while the $\tint$ symbol in the RHS denotes the projection map $\mc V_N\to\mc V_N/\partial\mc V_N$).
\end{enumerate}
\end{lemma}
\begin{proof}
Part (a) follows from \eqref{eq:partial-VN} and part (b) follows from the third relation in \eqref{20130917:eq2} and the fact that the product in $\mc V_N$ is commutative. Part (c) is an immediate consequence of parts (a) and (b).
\end{proof}
Next, let us illustrate the case $n=1$. Recall that $C^1(\mc V)=\Vect^\partial(\mc V)$
consists of all the derivations of the (noncommutative) associative product of $\mc V$ which commute with $\partial:\mc V\to\mc V$, while
$C^1(\mc V_N)=\Vect^\partial(\mc V_N)$ consists of all the derivations of the (commutative) associative product of $\mc V_N$ which commute with
$\partial:\mc V_N\to\mc V_N$.
Given $D\in C^1(\mc V)$ we define the map $\tr(D):\mc V_N\to\mc V_N$ by
\begin{equation}\label{eq:trace1}
\tr(D)(a_{ij})=D(a)_{ij}
\,,
\end{equation}
for every $a\in\mc V$ and $1\leq i,j\leq N$, and we extend it to $\mc V_N$ using the Leibniz rule.
\begin{lemma}\label{20240902:lem1}
For every $D\in C^1(\mc V)$ we have that $\tr(D)\in C^1(\mc V_N)$.
\end{lemma}
\begin{proof}
First, we need to verify that
the map $\tr(D)$ is well defined, namely, it is compatible with the defining relations \eqref{20130917:eq2} of $\mc V_N$. Clearly, $\tr(D)$ is linear, since $D$ is linear. Hence, we are left to show that
$$
\tr(D)((ab)_{ij})=\sum_{k=1}^N\tr(D)(a_{ik}b_{kj})
\,.
$$
This follows by construction. In fact, we have
\begin{align*}
&\tr(D)((ab)_{ij})=D(ab)_{ij}=\sum_{k=1}^N \left(D(a)_{ik}b_{kj}+a_{ik} D(b)_{kj}\right)
\\
&=\sum_{k=1}^N\left(\tr(D)(a_{ik})b_{kj}+a_{ik}\tr(D)(b_{kj})\right)
=\sum_{k=1}^N\tr(D)(a_{ik}b_{kj})
\,,
\end{align*}
where in the first equality we used \eqref{eq:trace1}, in the second equality we used the fact that $D$ is a derivation of $\mc V$, in the third equality we used again \eqref{eq:trace1}, and finally we used the fact that, by construction, $\tr(D)$ is extended to $\mc V_N$ by the Leibniz rule. Hence, we have that
$\tr(D)\in\Vect(\mc V_N)$. To conclude the proof we need to show that $\tr(D)$ commutes with $\partial:\mc V_N\to\mc V_N$. For this, it suffices to check that
$\partial\tr(D)(a_{ij})=\tr(D)((\partial a)_{ij})$, for every $a\in\mc V$ and $1\leq i,j\leq N$.
This is an immediate consequence of the definition of $\tr(D)$ given in \eqref{eq:trace1},
the fact that $D$ commutes with $\partial$, and the action of $\partial$ on $\mc V_N$ given in
\eqref{eq:partial-VN}.
\end{proof}
For $n=2$, given a skewsymmetric $2$-fold $\lambda$-bracket $\ldb-_{\lambda}-\rdb\in C^2(\mc V)$, we let $\tr(\ldb-_{\lambda}-\rdb):\mc V_N\otimes\mc V_N\to\mc V_N[\lambda]$ be the map defined by the RHS of \eqref{20240829:eq1} on the pair $(a_{ij},b_{hk})$,
$a,b\in \mc V$, $1\leq i,j,h,k\leq N$, and extended on  $\mc V_N\otimes\mc V_N$ using the left and right Leibniz rules \eqref{lleibnizPVA}, \eqref{rleibnizPVA}.
As previously mentioned, it is shown in \cite{DSKV} that $\tr(\ldb-_{\lambda}-\rdb)$ is a well defined linear map, satisfying
sesquilinearity \eqref{sesquiLCA} and skewsymmetry
\eqref{skewLCA} (and, of course, the Leibniz rules \eqref{lleibnizPVA}, \eqref{rleibnizPVA} by construction). Hence, $\tr(\ldb-_{\lambda}-\rdb)\in C^2(\mc V_N)$.

%%%
\section{From \texorpdfstring{\MakeLowercase{$n$}}{n}-fold \texorpdfstring{$\lambda$}{lambda}-brackets to poly-\texorpdfstring{$\lambda$}{lambda}-brackets}\label{sec:13.3}

The next result gives the generalization of equations \eqref{eq:trace0}, \eqref{eq:trace1}
and \eqref{20240829:eq1} to arbitrary $n$-fold skewsymmetric $\lambda$-brackets,
cf. Theorem \ref{Thm:IndBr}.
\begin{theorem}\label{Thm:IndBr-PVA}
Let $\mc V$ be a differential algebra, and let $\ldb-_{\lambda_1}-\dots-_{\lambda_{n-1}}-\rdb\in C^n(\mc V)$ be an $n$-fold $\lambda$-bracket on $\mc V$, $n\geq1$.
Then, for every $N\geq1$,
we have a well-defined $n$-$\lambda$-bracket $\tr(\ldb-_{\lambda_1}-\dots-_{\lambda_{n-1}}-\rdb)\in C^n(\mc V_N)$, given on generators $a^1_{i_1j_1},\dots a^n_{i_nj_n}\in\mc V_N$
by
\begin{equation}\label{Eq:TrBrRep-PVA}
\begin{split}
&\tr(\ldb-_{\lambda_1}-\dots-_{\lambda_{n-1}}-\rdb)(a^1_{i_1j_1},a^2_{i_2j_2},\dots,a^n_{i_nj_n})
\\
&=\sum_{\tau\in S_{n-1}}\sgn(\tau)\ldb {a^{\tau(1)}}_{\lambda_{\tau(1)}}{a^{\tau(2)}}_{\lambda_{\tau(2)}}\dots {a^{\tau(n-1)}}_{\lambda_{\tau(n-1)}}a^n\rdb_{\tau(\underline i,\underline j)}
\,,
\end{split}
\end{equation}
where we are using the notation \eqref{Eq:Not-RepIndex} and we are setting
$$
\tau(\underline i,\underline j)=(i_nj_{\tau(1)},i_{\tau(1)}j_{\tau(2)},\dots,i_{\tau(n-1)}j_n)
\,,
$$
and extended to $\mc V_N^{\otimes n}$
by the Leibniz rules \eqref{eq:leibnizpoly2} and \eqref{eq:leibnizpoly}.
\end{theorem}
The notation \eqref{Eq:Not-RepIndex} takes the following form in this context: write
$$
\ldb a^1{}_{\lambda_1}\dots a^{n-1}{}_{\lambda_{n-1}}a^n\rdb
%=\sum_{\underline m\in\mb Z_{\geq0}^{n-1}}
%f_{\underline m}\lambda^{\underline m}
=\sum_{\underline m\in\mb Z_{\geq0}^{n-1}}
f_{\underline m,1}\otimes \dots\otimes f_{\underline m,n-1}\otimes f_{\underline m,n}\lambda_1^{m_1}\lambda_2^{m_2}\dots\lambda_{n-1}^{m_{n-1}}
\,,
$$
where we are denoting $\underline m=(m_1,\dots,m_{n-1})$. Then
\begin{align*}
&\ldb a^1{}_{\lambda_1}\dots a^{n-1}{}_{\lambda_{n-1}}a^n\rdb_{i_1j_1,i_2j_2,\dots, i_nj_n}
\\
&=
\sum_{\underline m\in\mb Z_{\geq0}^{n-1}}
(f_{\underline m,1})_{i_1j_1}\dots(f_{\underline m,n-1})_{i_{n-1}j_{n-1}} (f_{\underline m,n})_{i_n,j_n}\lambda_1^{m_1}\lambda_2^{m_2}\dots\lambda_{n-1}^{m_{n-1}}
\,.
\end{align*}
In analogy with the notation used in \eqref{20240829:eq1}
we denote the LHS of \eqref{Eq:TrBrRep-PVA}
by $\{a^1_{i_1j_1}{}_{\lambda_1}\dots a^{n-1}_{i_{n-1}j_{n-1}}{}_{\lambda_{n-1}}a^n_{i_nj_n}\}$. Moreover, using the subgroup $S_n^{(k)}=\{\tau\in S_n\mid \tau(k)=k\}\subset S_n$, for $k=1,\dots,n$, we rewrite
\eqref{Eq:TrBrRep-PVA} as
\begin{equation}\label{Eq:TrBrRep-PVA-2}
\begin{split}
&\{a^1_{i_1j_1}{}_{\lambda_1}\dots a^{n-1}_{i_{n-1}j_{n-1}}{}_{\lambda_{n-1}}a^n_{i_nj_n}\}
\\
&=\sum_{\tau\in S_{n}^{(n)}}\sgn(\tau)\ldb {a^{\tau(1)}}_{\lambda_{\tau(1)}}\dots {a^{\tau(n-1)}}_{\lambda_{\tau(n-1)}}a^{\tau(n)}\rdb_{\tau(\underline i,\underline j)}
\,.
\end{split}
\end{equation}
\begin{proof}[Proof of Theorem \ref{Thm:IndBr-PVA}]
First we need to verify that the map
$\{-_{\lambda_1}-\dots-_{\lambda_{n-1}}-\}:
\mc V_N^{\otimes n}\to\mc V_N[\lambda_1,\dots,\lambda_{n-1}]$
is well defined, that is, it is compatible with the defining relations \eqref{20130917:eq2} of $\mc V_N$. Clearly, the RHS of \eqref{Eq:TrBrRep-PVA-2} is linear in $a^1,\dots,a^n$ since
$\ldb-_{\lambda_1}-\dots-_{\lambda_{n-1}}-\rdb$ is a linear map.
Hence, we are left to show that
\begin{equation}\label{20240902:toprove1}
\{a^1_{i_1j_1}{}_{\lambda_1}\dots a^{n-1}_{i_{n-1}j_{n-1}}{}_{\lambda_{n-1}}(bc)_{ij}\}
=\sum_{k=1}^N\{a^1_{i_1j_1}{}_{\lambda_1}\dots a^{n-1}_{i_{n-1}j_{n-1}}{}_{\lambda_{n-1}}b_{ik}c_{kj}\}
\end{equation}
and
\begin{equation}\label{20240902:toprove2}
\begin{split}
&\{a^1_{i_1j_1}{}_{\lambda_1}\dots (bc)_{ij}{}_{\lambda_s}\dots
a^{n-1}_{i_{n-1}j_{n-1}}{}_{\lambda_{n-1}}a^n_{i_nj_n}\}
\\
&=\sum_{k=1}^N\{a^1_{i_1j_1}{}_{\lambda_1}\dots b_{ik}c_{kj}{}_{\lambda_s}\dots
a^{n-1}_{i_{n-1}j_{n-1}}{}_{\lambda_{n-1}}a^n_{i_nj_n}\}
\,,
\end{split}
\end{equation}
for every $s=1,\dots,n-1$. Equation \eqref{20240902:toprove1}
follows using the same computation as in Lemma \ref{20240902:lem1}
applied to $D=\{a^1_{\lambda_1}\dots a^{n-1}_{\lambda_{n-1}}-\}$.
Before proving \eqref{20240902:toprove2}, let us prove that the
skewsymmetry property
\begin{equation}\label{20240902:toprove3}
\begin{split}
&\{a^1_{i_1j_1}{}_{\lambda_1}\cdots a^{n-1}_{i_{n-1}j_{n-1}}{}_{\lambda_{n-1}}a^n_{i_nj_n}\}
\\
&=\sgn(\tau)
|_{\lambda_n=\lambda_n^\dagger}\{ a^{\tau(1)}_{i_{\tau(1)}j_{\tau(1)}}{}_{\lambda_{\tau(1)}}\cdots a^{\tau(n-1)}_{i_{\tau(n-1)}j_{\tau(n-1)}}{}_{\lambda_{\tau(n-1)}}a^{\tau(n)}_{i_{\tau(n)}j_{\tau(n)}}\}
\end{split}
\end{equation}
holds for every $\tau\in S_n$ and $a^{1}_{i_1j_1},\dots, a^n_{i_nj_n}\in\mc V_N$. It actually suffices to prove \eqref{20240902:toprove3}
for $\tau=(12)$ and $\tau=\sigma=(12\dots n)$ since they generate the symmetric group
$S_n$. Equation \eqref{20240902:toprove3} for $\tau=(12)$ follows 
immediately since the RHS of \eqref{Eq:TrBrRep-PVA-2} changes sign if 
we swap $a^{1}_{i_1j_1}$ by $a^{2}_{i_2j_2}$ and $\lambda_1$ by $\lambda_2$.
On the other hand, using the cyclic skewsymmetry \eqref{eq:nfold-skew2bis} we have
\begin{equation}\label{20240903:eq1}
\begin{split}
&\{a^1_{i_1j_1}{}_{\lambda_1}\dots
a^{n-1}_{i_{n-1}j_{n-1}}{}_{\lambda_{n-1}}a^n_{i_nj_n}\}
=\sgn(\sigma)\Big|_{\lambda_n=\lambda_n^\dagger}
\sum_{k=1}^{n-1}\sum_{\tau\in S_{n}^{(n)}\text{s.t.}\tau(k)=1}
\sgn(\tau\sigma^{k-1})
\\
&\left(
\ldb{a^{\tau\sigma^{k-1}(2)}}_{\lambda_{\tau\sigma^{k-1}(2)}}
\dots
{a^{\tau\sigma^{k-1}(n)}}_{\lambda_{\tau\sigma^{k-1}(n)}}
a^{\tau\sigma^{k-1}(1)}\rdb^{\sigma^k}
\right)
_{\tau(\underline i,\underline j)}
\,.
\end{split}
\end{equation}
Note that, using \eqref{Eq:RepIndex1} we have
$$
(A^{\sigma^k})_{\tau(\underline i,\underline j)}
=A_{\tau\sigma^{k-1}(\sigma(\underline i),\sigma(\underline j))}
\,,
$$
where $\sigma(\underline i)=(i_{\sigma(1)},\dots,i_{\sigma(n)})=(i_2,\dots,i_n,i_1)$ (similarly for $\sigma(\underline j)$).
Clearly, if $\tau\in S_n^{(n)}$ is such that $\tau(k)=1$, then
$\tau\sigma^{k-1}\in S_n^{(1)}$ and $\tau\sigma^{k-1}(n+1-k)=n$. Hence, we have
$$
\bigcup_{k=1}^{n-1}\{\tau\sigma^{k-1}\mid \tau\in S_n^{(n)}\text{ s.t. }
\tau(k)=1\}
=\bigcup_{k=2}^{n}\{\tilde \tau\mid \tilde\tau\in S_n^{(1)}\text{ s.t. }
\tilde \tau(k)=n\}=S_n^{(1)}
\,.
$$
Using these observations, we rewrite \eqref{20240903:eq1} as
\begin{align*}
&\{a^1_{i_1j_1}{}_{\lambda_1}a^2_{i_2j_2}{}_{\lambda_2}\dots_{\lambda_{n-1}}a^n_{i_nj_n}\}
\\
&=\sgn(\sigma)|_{\lambda_n=\lambda_n^\dagger}
\sum_{\tau\in S_{n}^{(1)}}
\sgn(\tau)
\ldb{a^{\tau(2)}}_{\lambda_{\tau(2)}}
{a^{\tau(3)}}_{\lambda_{\tau(3)}}
\dots
{a^{\tau(n)}}_{\lambda_{\tau(n)}}
a^{\tau(1)}\rdb
_{\tau(\sigma(\underline i),\sigma(\underline j))}
\\
&=(-1)^{n+1}|_{\lambda_n=\lambda_n^\dagger}
\{a^2_{i_2j_2}{}_{\lambda_2}a^3_{i_3j_3}{}_{\lambda_3}\dots_{\lambda_{n}}a^1_{i_1j_1}\}
\,.
\end{align*}
In the last identity above we used \eqref{Eq:TrBrRep-PVA-2}. 
This completes the proof of \eqref{20240902:toprove3}.
Using the fact that $\{-_{\lambda_1}-\dots-_{\lambda_{n-1}}-\}$ is 
extended using the Leibniz rules \eqref{eq:leibnizpoly2} and \eqref{eq:leibnizpoly} we have that equations \eqref{20240902:toprove1} and \eqref{20240902:toprove3} imply \eqref{20240902:toprove2}
for all $s=1,\dots,n-1$, thus showing that $\{-_{\lambda_1}-\dots-_{\lambda_{n-1}}-\}$ is well defined.
The sesquilinearity axioms \eqref{eq:sesquipoly1} and \eqref{eq:sesquipoly2} for $\{-_{\lambda_1}-\dots-_{\lambda_{n-1}}-\}$ follow immediately from the definition \eqref{Eq:TrBrRep-PVA-2}.
We then have that $\{-_{\lambda_1}-\dots-_{\lambda_{n-1}}-\}=\tr(\ldb -_{\lambda_1}-\dots-_{\lambda_{n-1}}-\rdb)\in C^{n}(\mc V_N)$, concluding the proof.
\end{proof}
%

%%%
\section[From \MakeLowercase{d}PVA cohomology to PVA cohomology]{From variational double Poisson vertex algebra cohomology to variational Poisson vertex algebra cohomology}
\label{sec:dPVAtoPVA}

As a consequence of Lemma \ref{20240829:lem1} and Theorem \ref{Thm:IndBr-PVA} we have a well defined linear map $\tr:C(\mc V)\to C(\mc V_N)$, for every $N\geq1$, which we call the \emph{trace map}, such that
$\tr(C^n(\mc V))\subset C^n(\mc V_N)$, for every $n\in\mb Z_{\geq0}$.
\begin{theorem} \label{Thm:dP-rep2-PVA} 
Let $\mc V$ be a dPVA with $2$-fold $\lambda$-bracket $\llbracket-_\lambda-\rrbracket$
and let $\mc V_N$ be the corresponding PVA with $\lambda$-bracket $[-_\lambda-]$ defined using \eqref{20240829:eq1}, $N\geq1$.
Let $\dd: C(\mc V)\to\mc V$ be the differential of the variational dPVA complex
given by \eqref{eq:dP0} and \eqref{eq:dP-1}, and let $\dd_N: C(\mc V_N)\rightarrow C(\mc V_N)$ be the differential of the variational PVA complex given by \eqref{eq:dH-poly-0}
and \eqref{eq:dH-poly-n}. The trace map is a morphism of complexes 
\begin{equation} \label{Eq:dP-rep2-PVA}
    \tr : (C(\mc V),\dd) \longrightarrow 
    (C(\mc V_N), (-1)^{\bullet}\dd_{N}),  
\end{equation}
which descends to a linear map $\dPVH(\mc V)\to \PVH(\mc V_N)$ in cohomology. 
\end{theorem} 
\begin{proof}
We need to show that the following diagram
\begin{equation}\label{diag1}
% \text{this diagram has been commented to speed up compiling}
% \begin{comment}
\begin{tikzcd}
C^n(\mc V) \arrow[r,"\dd"] \arrow[d,"\tr"] & C^{n+1}(\mc V)\arrow[d,"\tr"]
\\
C^n(\mc V_N) \arrow[r,"(-1)^n \dd_N"] & C^{n+1}(\mc V_N)
\end{tikzcd}
% \end{comment}
\end{equation}
is commutative for every $n\in\mb Z_{\geq0}$.
Let us consider the case $n=0$ first. For $\tint f\in C^0(\mc V)=\mc V_\sharp$,
using \eqref{eq:dP0} and \eqref{eq:trace1}, we have ($a_{ij}\in\mc V_N$)
\begin{align*}
&\tr(\dd(\tint f))(a_{ij})=-\tr(\mult \ldb f_\lambda - \rdb|_{\lambda=0})(a_{ij})
\\
&=-(\mult\llbracket f_{\lambda}a\rrbracket)_{ij}|_{\lambda=0}
=-\sum_{k=1}^N\llbracket f_{\lambda}a\rrbracket'_{ik}\llbracket f_{\lambda}a\rrbracket''_{kj}|_{\lambda=0}
\,.
\end{align*}
On the other hand, using Lemma \ref{20240902:lem1}(c), \eqref{eq:trace0} and \eqref{eq:dH-poly-0}
we have
$$
\dd_{N}(\tr(\tint f))(a_{ij})=-\sum_{k=1}^N[f_{kk}{}_\lambda a_{ij}]|_{\lambda=0}
=-\sum_{k=1}^N\llbracket f_\lambda a\rrbracket'_{ik}\llbracket f_\lambda a\rrbracket''_{kj}
|_{\lambda=0}
\,,
$$
where in the last equality we used \eqref{20240829:eq1}. Hence, the diagram \eqref{diag1} commutes for $n=0$.

We consider now the case $n=1$. Let $D\in C^1(\mc V)$. Using \eqref{20240729:eq1} and the definition of the trace map \eqref{Eq:TrBrRep-PVA} (see also \eqref{20240829:eq1}) we have
($a_{ij},b_{hk}\in\mc V_N$)
\begin{equation}\label{20240903:eq2}
\tr(\dd(D)_\lambda)(a_{ij},b_{hk})
=\left(
D(\llbracket a_\lambda b\rrbracket)-\llbracket D(a)_\lambda b\rrbracket-\llbracket a_\lambda D(b)\rrbracket
\right)_{hj,ik}
\,.
\end{equation}
On the other hand, using \eqref{eq:trace1} and \eqref{20240828:eq1}
we have
\begin{equation}\label{20240903:eq3}
\dd_{N}(\tr(D))_\lambda(a_{ij},b_{hk})
=[D(a)_{ij}{}_\lambda b_{hk}] +[a_{ij}{}_\lambda D(b)_{hk}]
-\tr(D)([a_{ij}{}_\lambda b_{hk})])
\end{equation}
By Lemma \ref{20240902:lem1} and equations \eqref{20240829:eq1} and
\eqref{eq:trace1} we have
\begin{equation}\label{20240903:eq4}
\begin{split}
&\tr(D)([a_{ij}{}_\lambda b_{hk})])
=\tr(D)(\llbracket a_\lambda b\rrbracket'_{hj}\llbracket a_\lambda b\rrbracket''_{ik})
\\
&=\tr(D)(\llbracket a_\lambda b\rrbracket'_{hj})\llbracket a_\lambda b\rrbracket''_{ik}
+\llbracket a_\lambda b\rrbracket'_{hj}\tr(D)(\llbracket a_\lambda b\rrbracket''_{ik})
\\
&=D(\llbracket a_\lambda b\rrbracket')_{hj}\, \llbracket a_\lambda b\rrbracket''_{ik}
+\llbracket a_\lambda b\rrbracket'_{hj}\, D(\llbracket a_\lambda b\rrbracket'')_{ik}
=D(\llbracket a_\lambda b\rrbracket)_{hj,ik}
\,.
\end{split}
\end{equation}
From \eqref{20240829:eq1} and \eqref{20240903:eq4} we see that the RHS of \eqref{20240903:eq2} and \eqref{20240903:eq3}
are opposite to each other, thus proving that the diagram \eqref{diag1}
is commutative for $n=1$.

For $n\geq2$, let $Q=\ldb-_{\lambda_1}-\dots-_{\lambda_{n-1}}-\rdb\in C^{n}(\mc V)$
and $\tr(Q)=\{-_{\lambda_1}-\dots_{\lambda_{n-1}}-\}\in C^n(\mc V_N)$.
Using \eqref{Eq:TrBrRep-PVA-2} and \eqref{eq:dP-1} we have
\begin{equation}\label{20240905:LHS}
\begin{split}
&(-1)^n\tr(\dd (Q)_{\lambda_1,\dots,\lambda_n})(a^1_{i_1j_1},\dots, a^{n+1}_{i_{n+1}j_{n+1}})
\\
&=(-1)^n\sum_{\tau\in S_{n+1}^{(n+1)}}\sgn(\tau)
\Big(
d(Q)_{\lambda_{\tau(1)},\dots,\lambda_{\tau(n)}}(a^{\tau(1)},\dots,a^{\tau(n+1)})
\Big)_{\tau(\underline i,\underline j)}
\\
&=\sum_{\tau\in S_{n+1}^{(n+1)}}\sgn(\tau)
\bigg(
\sum_{s=1}^{n}(-1)^{s+1}
\llbracket a^{\tau(s)} {}_{\lambda_{\tau(s)}} \ldb a^{\tau(1)}{}_{\lambda_{\tau(1)}}\stackrel{s}{\check{\dots}}
a^{\tau(n)}{}_{\lambda_{\tau(n)}}a^{\tau(n+1)}\rdb\rrbracket_{(s)}
\\
&+(-1)^{n+1}\llbracket \ldb a^{\tau(1)}{}_{\lambda_{\tau(1)}}\dots a^{\tau(n-1)}{}_{\lambda_{\tau(n-1)}}a^{\tau(n)}\rdb_{\lambda_{\tau(1)}+\dots+\lambda_{\tau(n)}}a^{\tau(n+1)}\rrbracket_{L}
\\
&+\sum_{s=1}^n(-1)^{s}
\ldb a^{\tau(1)}{}_{\lambda_{\tau(1)}}\dots a^{\tau(s-1)}{}_{\lambda_{\tau(s-1)}}
\llbracket a^{\tau(s)}{}_{\lambda_{\tau(s)}}a^{\tau(s+1)}\rrbracket_{\lambda_{\tau(s)}+\lambda_{\tau(s+1)}}a^{\tau(s+2)}_{\lambda_{\tau(s+2)}}\dots
\\
&\quad\quad\quad\quad\quad\quad\quad\quad\quad\quad\quad\quad\quad\quad\quad\quad
\dots a^{\tau(n)}{}_{\lambda_{\tau(n)}}a^{\tau(n+1)}\rdb_{L}
\\
&
-\ldb a^{\tau(2)}{}_{\lambda_{\tau(2)}}\dots a^{\tau(n)}{}_{\lambda_{\tau(n)}}
\llbracket a^{\tau(1)}{}_{\lambda_{\tau(1)}}a^{\tau(n+1)}\rrbracket\rdb_R
\bigg)_{\tau(\underline i,\underline j)}
\,.
\end{split}
\end{equation}
On the other hand, using \eqref{20250720:eq1} with $c=\tr(Q)$ we have
\begin{equation}\label{20240905:RHS}
\begin{split}
&\dd_{N}(\tr(Q))_{\lambda_1,\dots,\lambda_{n}}(a^1_{i_1j_1},\dots,a^{n+1}_{i_{n+1}j_{n+1}})
\\
&=
\sum_{s=1}^n(-1)^{s+1}[ a^s_{i_sj_s}{}_{\lambda_s}\tr(Q)_{\lambda_1,\stackrel{s}{\check{\dots}},\lambda_n}
(a^1_{i_1j_1},\stackrel{s}{\check{\dots}},a^n_{i_nj_n},a^{n+1}_{i_{n+1}j_{n+1}})]
\\
&+(-1)^{n+1}
[\tr(Q)_{\lambda_1,\dots,\lambda_{n-1}}(a^1_{i_1,j_1},\dots,a^n_{i_nj_n})_{\lambda_1+\dots+\lambda_n}a^{n+1}_{i_{n+1}j_{n+1}}]
\\
&+\sum_{1\leq s< t\leq n}(-1)^{s+t}
\tr(Q)_{\lambda_s+\lambda_t,\lambda_1,\stackrel{s}{\check{\dots}} \stackrel{t}{\check{\dots}},\lambda_n}
([a^s_{i_sj_s}{}_{\lambda_s}a^t_{i_tj_t}],a^1_{i_1j_1},\stackrel{s}{\check{\dots}} \stackrel{t}{\check{\dots}}, a^{n}_{i_nj_n},a^{n+1}_{i_{n+1}j_{n+1}})
\\
&+\sum_{s=1}^n(-1)^{s}
\tr(Q)_{\lambda_1,\stackrel{s}{\check{\dots}},\lambda_n}(a^1_{i_1j_1},\stackrel{s}{\check{\dots}},a^n_{i_nj_n},[a^s_{i_sj_s}{}_{\lambda_s}a^{n+1}_{i_{n+1}j_{n+1}}])
\,.
\end{split}
\end{equation}
To conclude the proof we need to show that the RHS of \eqref{20240905:LHS}
and the RHS of \eqref{20240905:RHS} are equal. This follows from the following identities that we are going to prove
\begin{subequations}
\begin{align}
\begin{split}\label{20240905:toprove1}
&[\tr(Q)_{\lambda_1,\dots,\lambda_{n-1}}(a^1_{i_1,j_1},\dots,a^n_{i_nj_n})_{\lambda_1+\dots+\lambda_n}a^{n+1}_{i_{n+1}j_{n+1}}]
\\
&=\!\!\sum_{\tau\in S_{n+1}^{(n+1)}}\!\!\!\!\!\sgn(\tau)
\Big(\llbracket \ldb a^{\tau(1)}{}_{\lambda_{\tau(1)}}\dots
a^{\tau(n-1)}{}_{\lambda_{\tau(n-1)}}a^{\tau(n)}\rdb_{\lambda_{\tau(1)}+\dots+\lambda_{\tau(n)}}a^{\tau(n+1)}\rrbracket_{L}
\Big)_{\tau(\underline i,\underline j)}
\,;
\end{split}
\\
\begin{split}\label{20240905:toprove2}
&\sum_{s=1}^n(-1)^{s}
\tr(Q)_{\lambda_1,\stackrel{s}{\check{\dots}},\lambda_n}
(a^1_{i_1j_1},\stackrel{s}{\check{\dots}},a^n_{i_nj_n},[a^s_{i_sj_s}{}_{\lambda_s}a^{n+1}_{i_{n+1}j_{n+1}}])
\\
&=\!\!\sum_{\tau\in S_{n+1}^{(n+1)}}\!\!\!\!\!\sgn(\tau)
\Big(
(-1)^{n}
\ldb a^{\tau(1)}{}_{\lambda_{\tau(1)}}\dots a^{\tau(n-1)}{}_{\lambda_{\tau(n-1)}}
\llbracket a^{\tau(n)}{}_{\lambda_{\tau(n)}}a^{\tau(n+1)}\rrbracket\rdb_{L}
\\
&
-\ldb a^{\tau(2)}{}_{\lambda_{\tau(2)}}\dots a^{\tau(n)}{}_{\lambda_{\tau(n)}}
\llbracket a^{\tau(1)}{}_{\lambda_{\tau(1)}}a^{\tau(n+1)}\rrbracket\rdb_R
\Big)_{\tau(\underline i,\underline j)}
\,;
\end{split}
\\
\begin{split}\label{20240905:toprove3}
&
\sum_{s=1}^n(-1)^{s+1}[ a^s_{i_sj_s}{}_{\lambda_s}\tr(Q)_{\lambda_1,\stackrel{s}{\check{\dots}},\lambda_n}
(a^1_{i_1j_1},\stackrel{s}{\check{\dots}},a^n_{i_nj_n},a^{n+1}_{i_{n+1}j_{n+1}})]
\\
&=
\sum_{s=1}^{n}
\sum_{\tau\in S_{n+1}^{(n+1)}}(-1)^{s+1}\sgn(\tau)
\Big(
\llbracket a^{\tau(s)} {}_{\lambda_{\tau(s)}} \ldb a^{\tau(1)}{}_{\lambda_{\tau(1)}}\stackrel{s}{\check{\dots}}
a^{\tau(n)}{}_{\lambda_{\tau(n)}}a^{\tau(n+1)}\rdb\rrbracket_{(s)}
\Big)_{\tau(\underline i,\underline j)}\,;
\end{split}
\\
\begin{split}\label{20240905:toprove4}
&\sum_{1\leq s< t\leq n}(-1)^{s+t}
\tr(Q)_{\lambda_s+\lambda_t,\lambda_1,\stackrel{s}{\check{\dots}} \stackrel{t}{\check{\dots}},\lambda_n}
([a^s_{i_sj_s}{}_{\lambda_s}a^t_{i_tj_t}],a^1_{i_1j_1},\stackrel{s}{\check{\dots}} \stackrel{t}{\check{\dots}}, a^{n}_{i_nj_n},a^{n+1}_{i_{n+1}j_{n+1}})
\\
&=\sum_{s=1}^{n-1}
\sum_{\tau\in S_{n+1}^{(n+1)}}(-1)^{s}\sgn(\tau)
\Big(
\ldb a^{\tau(1)}{}_{\lambda_{\tau(1)}}\dots
\\
&
\dots a^{\tau(s-1)}{}_{\lambda_{\tau(s-1)}}
\llbracket a^{\tau(s)}{}_{\lambda_{\tau(s)}}a^{\tau(s+1)}\rrbracket_{\lambda_{\tau(s)}+\lambda_{\tau(s+1)}}a^{\tau(s+2)}\dots
{}_{\lambda_{\tau(n)}}a^{\tau(n+1)}\rdb_{L}
\Big)_{\tau(\underline i,\underline j)}
\,.
\end{split}
\end{align}
\end{subequations}
We start by proving \eqref{20240905:toprove1}.
Using the definition
\eqref{Eq:TrBrRep-PVA-2} of the trace map we have that the LHS of \eqref{20240905:toprove1} is equal to
\begin{equation}\label{20240905:eq1}
\begin{split}
&\sum_{\tau\in S_n^{(n)}}\sgn(\tau)[
{\ldb {a^{\tau(1)}}_{\lambda_{\tau(1)}}\dots {a^{\tau(n-1)}}_{\lambda_{\tau(n-1)}}a^{\tau(n)}\rdb_{\tau(\underline i,\underline j)}}_{\lambda_1+\dots+\lambda_n}a^{n+1}_{i_{n+1}j_{n+1}}
]
\\
&=
\sum_{\tau\in S_n^{(n)}}\sum_{k=1}^n\sgn(\tau)
\Big(\llbracket \sigma^{1-k}\ldb {a^{\tau(1)}}_{\lambda_{\tau(1)}}\dots
\\
&\quad\quad\quad\quad
\dots{a^{\tau(n-1)}}_{\lambda_{\tau(n-1)}}a^{\tau(n)}\rdb_{\lambda_1+\dots+\lambda_n}a^{n+1}\rrbracket_L
\Big)_{\tau\sigma^{k-1}(\underline i,\underline j)}
\\
&=
\sum_{\tau\in S_n^{(n)}}\sum_{k=1}^n\sgn(\tau\sigma^{k-1})
\Big(\llbracket \ldb {a^{\tau\sigma^{k-1}(1)}}_{\lambda_{\tau\sigma^{k-1}(1)}}\dots
\\
&\quad\quad\quad\quad
\dots{a^{\tau\sigma^{k-1}(n-1)}}_{\lambda_{\tau\sigma^{k-1}(n-1)}}a^{\tau\sigma^{k-1}(n)}\rdb_{\lambda_1+\dots+\lambda_n}a^{n+1}\rrbracket_L
\Big)_{\tau\sigma^{k-1}(\underline i,\underline j)}
\,.
\end{split}
\end{equation}
In the first equality above we used the identity
$$
[{A_{(\underline i,\underline j)}}_\lambda a^{n+1}_{i_{n+1},j_{n+1}}]
=\sum_{k=1}^n\Big(\llbracket \sigma^{1-k}(A)_{\lambda}a^{n+1}\rrbracket_L
\Big)_{\sigma^{k-1}(\underline i,\underline j)}
\,,
$$
(here $\sigma=(12\dots n)$) which can be easily checked for $A\in\mc V^{\otimes n}$,
and $\underline i=(i_1,\dots,i_n)$, $\underline j=(j_1,\dots,j_n)$,
where $1\leq i_k,j_k\leq N$, and in the second equality we used the cyclic skewsymmetry \eqref{eq:nfold-skew2bis} and sesquilinearity \eqref{20140702:eq4}.
Clearly, $\tau\sigma^{k-1}\in S_{n+1}^{(n+1)}$. Moreover,
$\tau\sigma^{k-1}(n+1-k)=n$, for every $k=1,\dots,n$. Hence,
$$
\bigcup_{k=1}^{n}\{\tau\sigma^{k-1}\mid \tau\in S_n^{(n)}\}
=\bigcup_{k=1}^{n}\{\tilde \tau\mid \tilde \tau\in S_{n+1}^{(n+1)}\text{ s.t. }
\tilde \tau(k)=n\}=S_{n+1}^{(n+1)}
\,,
$$
from which follows that the last term in \eqref{20240905:eq1} is equal to the RHS of \eqref{20240905:toprove1}.

Next, using \eqref{20240829:eq1} and the Leibniz rule \eqref{eq:leibnizpoly} we rewrite the LHS of\eqref{20240905:toprove2} as
\begin{align}
\begin{split}
&\sum_{s=1}^n(-1)^{s}
\tr(Q)_{\lambda_1,\stackrel{s}{\check{\dots}},\lambda_n} \big( a^1_{i_1j_1},\stackrel{s}{\check{\dots}},a^n_{i_nj_n},
(\llbracket a^s{}_{\lambda_s}a^{n+1}\rrbracket)'_{i_{n+1} j_s}\big)\, 
\times
\\
&\quad\quad\quad\quad\quad\quad\quad\quad
\times(\llbracket a^s{}_{\lambda_s}a^{n+1}\rrbracket)''_{i_s j_{n+1}}
\label{20240905:eq2}
\end{split}
\\
\begin{split}
&+\sum_{s=1}^n(-1)^{s}
\tr(Q)_{\lambda_1,\stackrel{s}{\check{\dots}},\lambda_n} \big(a^1_{i_1j_1},\stackrel{s}{\check{\dots}},a^n_{i_nj_n},
(\llbracket a^s{}_{\lambda_s}a^{n+1}\rrbracket)''_{i_s j_{n+1}}\big)
\times
\\
&\quad\quad\quad\quad\quad\quad\quad\quad
\times (\llbracket a^s{}_{\lambda_s}a^{n+1}\rrbracket)'_{i_{n+1} j_s}
\,.
\label{20240905:eq3}
\end{split}
\end{align}
Using the definition of the trace map \eqref{Eq:TrBrRep-PVA-2} 
and the notation \eqref{badnotation0}
we have
(recall from \S\ref{ss:Proof-dP-rep2} the notation $S_n^{(i,j)}=S_n^{(i)}\cap S_n^{(j)}$)
\begin{align*}
&\eqref{20240905:eq2}
=\sum_{s=1}^n\sum_{\tau\in S_{n+1}^{(s,n+1)}}
(-1)^s\sgn(\tau)
\Big(
\ldb {a^{\tau(1)}}_{\lambda_{\tau(1)}}\stackrel{s}{\check{\dots}}
\\
&\quad\quad\quad\quad\quad\quad\quad\quad
\dots{a^{\tau(n)}}_{\lambda_{\tau(n)}}\llbracket {a^s}_{\lambda_s} a^{n+1}\rrbracket\rdb_L
\Big)_{\tau(s s+1\dots n)(\underline i,\underline j)}
\,.
\end{align*}
Note that we have a bijective correspondence between the sets
$S_{n+1}^{(s,n+1)}$ and $\{\tilde\tau\in S_{n+1}^{(n+1)}\mid \tilde{\tau}(n)=s\}$, for every $s=1,\dots,n$, given by
$$
\tau\in S_{n+1}^{(s,n+1)}\mapsto\tilde\tau=\tau(s s+1\dots n)
\,.
$$
Moreover, $\sgn(\tilde\tau)=(-1)^{n+s}\sgn(\tau)$ and $S_{n+1}^{(n+1)}=\bigcup_{s=1}^{n}\{\tilde\tau\in S_{n+1}^{(n+1)}\mid \tilde{\tau}(n)=s\}$. Hence, by replacing $\tau$ with $\tilde \tau(n\, n-1\dots s)$ we have
\begin{align*}
&\eqref{20240905:eq2}
=(-1)^n\sum_{\tau\in S_{n+1}^{(n+1)}}
\sgn(\tau)
\Big(
\ldb {a^{\tau(1)}}_{\lambda_{\tau(1)}}\dots
\\
&\quad\quad\quad\quad\quad\quad\quad\quad
\dots{a^{\tau(n-1)}}_{\lambda_{\tau(n-1)}}\llbracket {a^{\tau(n)}}_{\lambda_{\tau(n)}} a^{\tau(n+1)}\rrbracket\rdb_L
\Big)_{\tau(\underline i,\underline j)}
\,,
\end{align*}
which gives the first sum in the RHS of \eqref{20240905:toprove2}.
Similarly, using now the notation \eqref{badnotation3}, it follows that
\begin{align*}
&\eqref{20240905:eq3}
=-\sum_{\tau\in S_{n+1}^{(n+1)}}
\sgn(\tau)
\Big(
\ldb {a^{\tau(2)}}_{\lambda_{\tau(2)}}\dots
{a^{\tau(n)}}_{\lambda_{\tau(n)}}\llbracket {a^{\tau(1)}}_{\lambda_{\tau(1)}} a^{\tau(n+1)}\rrbracket\rdb_R
\Big)_{\tau(\underline i,\underline j)}
\,,
\end{align*}
which gives the second sum in the RHS of \eqref{20240905:toprove2}.

Let us now prove \eqref{20240905:toprove3}. To ease notation, for $s=1,\dots,n$, $\tau\in S_{n+1}^{s,n+1}$ and $a^1,\dots,a^{n+1}\in\mc V$, we simply denote
$$
\ldb {a^{\tau(1)}}_{\lambda_{\tau(1)}}\stackrel{s}{\check{\dots}}{a^{\tau(n)}}_{\lambda_{\tau(n)}}a^{\tau(n+1)}\rdb
=x^1\otimes\dots\otimes x^n\in\mc V^{\otimes n}[\lambda_1,\dots,\lambda_{n-1}]\,,
$$
and, for $i_1,\dots,i_{n+1},j_1,\dots,j_{n+1}$, we let
$$
\ldb {a^{\tau(1)}}_{\lambda_{\tau(1)}}\stackrel{s}{\check{\dots}}{a^{\tau(n)}}_{\lambda_{\tau(n)}}a^{\tau(n+1)}\rdb_{i_{n+1}j_{\tau(1)},\dots,i_{\tau(s-1)}j_{\tau(s+1)},\dots,i_{\tau(n)},j_{n+1}}
=y_1\dots y_n\in\mc V_N\,,
$$
where
$$
y_k=\left\{
\begin{array}{ll}
(x^k)_{i_{\tau(k-1)}j_{\tau(k)}}\,,
&
k<s
\\
(x^s)_{i_{\tau(s-1)}j_{\tau(s+1)}}\,,
&
k=s
\\
(x^k)_{i_{\tau(k)}j_{\tau(k+1)}}\,,
&
k>s\,.
\end{array}
\right.
$$
Using the definition of the trace map, the above shorthand and the Leibniz rule \eqref{lleibnizPVA},
we can rewrite the LHS of \eqref{20240905:toprove3} as
\begin{align*}
&
\sum_{s=1}^n(-1)^{s+1}
\sum_{k=1}^n\sum_{\tau\in S_{n+1}^{(s,n+1)}}\sgn(\tau)
y_1\dots y_{k-1}[ a^s_{i_sj_s}{}_{\lambda_s}y_k]y_{k+1}\dots y_n
\\
&=
\sum_{s=1}^n(-1)^{s+1}\Bigg(
\sum_{k=1}^{s-1}\sum_{\tau\in S_{n+1}^{(s,n+1)}}\sgn(\tau)
y_1\dots y_{k-1}\llbracket a^s{}_{\lambda_s}x^k\rrbracket_{i_{\tau(k-1)}j_s,i_sj_{\tau(k)}}
y_{k+1}\dots y_n
\\
&+\sum_{\tau\in S_{n+1}^{(s,n+1)}}\sgn(\tau)
y_1\dots y_{s-1}\llbracket a^s{}_{\lambda_s}x^s\rrbracket_{i_{\tau(s-1)}j_s,i_sj_{\tau(s+1)}}
y_{s+1}\dots y_n
\\
&+\sum_{k=s+1}^n\sum_{\tau\in S_{n+1}^{(s,n+1)}}\sgn(\tau)
y_1\dots y_{k-1}\llbracket a^s{}_{\lambda_s}x^k\rrbracket_{i_{\tau(k)}j_s,i_sj_{\tau(k+1)}}
y_{k+1}\dots y_n
\Bigg)
\\
&=
\sum_{s=1}^n(-1)^{s+1}\Bigg(
\sum_{k=1}^{s-1}\sum_{\tau\in S_{n+1}^{(s,n+1)}}\sgn(\tau)
\Big(\llbracket a^s{}_{\lambda_s}
\ldb {a^{\tau(1)}}_{\lambda_{\tau(1)}}\stackrel{s}{\check{\dots}}
\\
&\quad\quad\quad\quad\quad\quad\quad\quad\quad\quad\quad\quad\quad
\dots{a^{\tau(n)}}_{\lambda_{\tau(n)}}a^{\tau(n+1)}\rdb\rrbracket_{(k)}
\Big)_{\tau(s s-1\dots k)(\underline i,\underline j)}
\\
&+\sum_{\tau\in S_{n+1}^{(s,n+1)}}\sgn(\tau)
\Big(\llbracket a^s{}_{\lambda_s}
\ldb {a^{\tau(1)}}_{\lambda_{\tau(1)}}\stackrel{s}{\check{\dots}}{a^{\tau(n)}}_{\lambda_{\tau(n)}}a^{\tau(n+1)}\rdb\rrbracket_{(s)}
\Big)_{\tau(\underline i,\underline j)}
\\
&+\sum_{k=s+1}^n\sum_{\tau\in S_{n+1}^{(s,n+1)}}\sgn(\tau)
\Big(\llbracket a^s{}_{\lambda_s}\ldb {a^{\tau(1)}}_{\lambda_{\tau(1)}}\stackrel{s}{\check{\dots}}
\\
&\quad\quad\quad\quad\quad\quad\quad\quad\quad\quad\quad\quad\quad
\dots{a^{\tau(n)}}_{\lambda_{\tau(n)}}a^{\tau(n+1)}\rdb\rrbracket_{(k)}
\Big)_{\tau(s s+1\dots k)(\underline i,\underline j)}
\Bigg)
\\
&=
\sum_{k,s=1}^n
\sum_{\tilde \tau\in S_{n+1}^{(n+1)}\text{s.t.}\tilde \tau(k)=s}
(-1)^{k+1}\sgn(\tilde \tau)
\Big(\llbracket a^{\tilde \tau(k)}{}_{\lambda_{\tilde \tau(k)}}
\ldb {a^{\tilde \tau(1)}}_{\lambda_{\tilde \tau(1)}}\stackrel{k}{\check{\dots}}
\\
&\quad\quad\quad\quad\quad\quad\quad\quad\quad\quad\quad\quad\quad
\dots{a^{\tilde \tau(n)}}_{\lambda_{\tilde \tau(n)}}a^{\tilde \tau(n+1)}\rdb\rrbracket_{(k)}
\Big)_{\tilde \tau(\underline i,\underline j)}
\,.
\end{align*}
In the second equality above we used the notation \eqref{20240805:eq1} for $D=\llbracket a^s_{\lambda_s}-\rrbracket$, and in the last equality we replaced $\tau$ by $\tilde \tau(k\, k+1\dots s)$ for $k<s$ and by $\tilde \tau(k\, k-1\dots s)$ for $k>s$, with $\tilde\tau\in S_{n+1}^{(n+1)}$ such that $\tilde \tau(k)=s$. The last sum above is clearly equal to the RHS of \eqref{20240905:toprove3}.

We are left to prove \eqref{20240905:toprove4}. Using the Leibniz rule \eqref{eq:leibnizpoly2} for $i=1$ and \eqref{20240829:eq1} we rewrite the LHS of \eqref{20240905:toprove4} as
\begin{equation}\label{20240906:eq1}
\begin{split}
&\sum_{1\leq s< t\leq n}(-1)^{s+t}
\tr(Q)_{\lambda_s+\lambda_t+x,\lambda_1,\stackrel{s}{\check{\dots}} \stackrel{t}{\check{\dots}},\lambda_n}
((\llbracket a^s{}_{\lambda_s}a^t\rrbracket')_{i_tj_s},a^1_{i_1j_1},\stackrel{s}{\check{\dots}} \stackrel{t}{\check{\dots}}
\\
&\quad\quad\quad\quad\quad\quad\quad\quad
\dots,a^{n+1}_{i_{n+1}j_{n+1}})
\big(|_{x=\partial}(\llbracket a^s{}_{\lambda_s}a^t\rrbracket'')_{i_sj_t}\big)
\\
&+\sum_{1\leq s< t\leq n}(-1)^{s+t}
\tr(Q)_{\lambda_s+\lambda_t+x,\lambda_1,\stackrel{s}{\check{\dots}} \stackrel{t}{\check{\dots}},\lambda_n}
((\llbracket a^s{}_{\lambda_s}a^t\rrbracket'')_{i_sj_t},a^1_{i_1j_1},\stackrel{s}{\check{\dots}} \stackrel{t}{\check{\dots}}
\\
&\quad\quad\quad\quad\quad\quad\quad\quad
\dots,a^{n+1}_{i_{n+1}j_{n+1}})
\big(|_{x=\partial}(\llbracket a^s{}_{\lambda_s}a^t\rrbracket')_{i_tj_s}\big)
\,.
\end{split}
\end{equation}
Using the skewsymmetry property \eqref{eq:skewpoly} of the trace map, the first sum in \eqref{20240906:eq1} can be rewritten as
\begin{equation}\label{20240906:eq2}
\begin{split}
&\sum_{1\leq s< t\leq n}(-1)^{s}
\tr(Q)_{\lambda_1,\stackrel{s}{\check{\dots}},\lambda_{t-1},\lambda_s+\lambda_t+x,\lambda_{t+1},\dots,\lambda_n}
(a^1_{i_1j_1},\stackrel{s}{\check{\dots}}
\\
&\quad\quad\quad\quad
\dots,(\llbracket a^s{}_{\lambda_s}a^t\rrbracket')_{i_tj_s},\dots
,a^{n+1}_{i_{n+1}j_{n+1}})
\big(|_{x=\partial}(\llbracket a^s{}_{\lambda_s}a^t\rrbracket'')_{i_sj_t}\big)
\end{split}
\end{equation}
Similarly, using also the sesquilinearity \eqref{eq:sesquipoly1} and the skewsymmetry property \eqref{eq:skew2} of the $2$-fold $\lambda$-bracket $\llbracket-_{\lambda}-\rrbracket$, we rewrite the second sum in \eqref{20240906:eq1} as
\begin{equation}\label{20240906:eq3}
\begin{split}
&\sum_{1\leq t< s\leq n}(-1)^{s}
\tr(Q)_{\lambda_1,\dots,\lambda_{t-1},\lambda_s+\lambda_t+x,\lambda_{t+1},\stackrel{s}{\check{\dots}},\lambda_n}
(a^1_{i_1j_1},\dots
\\
&\quad\quad\quad\quad
\dots,(\llbracket a^s{}_{\lambda_s}a^t\rrbracket')_{i_tj_s},\stackrel{s}{\check{\dots}},a^{n+1}_{i_{n+1}j_{n+1}})
\big(|_{x=\partial}(\llbracket a^s{}_{\lambda_s}a^t\rrbracket'')_{i_sj_t}\big)
\end{split}
\end{equation}
Next, we use the definition of the trace map \eqref{Eq:TrBrRep-PVA-2} and equation
\eqref{badnotation2}
to get
\begin{align*}
&\eqref{20240906:eq2}
=\sum_{1\leq s<t\leq n}\sum_{k=1}^{s-1}\sum_{\substack{\tau\in S_{n+1}^{(s,n+1)}\\\text{s.t. }\tau(k)=t}}
(-1)^s\sgn(\tau)
\\
&\quad\quad
\Big(\ldb {a^{\tau(1)}}_{\lambda_{\tau(1)}}\dots
\llbracket {a^s}_{\lambda s}a^t\rrbracket_{\lambda_s+\lambda_t}
\stackrel{s}{\check{\dots}}{a^{\tau(n)}}_{\lambda_{\tau(n)}}a^{\tau(n+1)}\rdb_L
\Big)_{\tau(s s-1\dots k)(\underline i,\underline j)}
\\
&+\sum_{1\leq s<t\leq n}\sum_{k=s+1}^{n}\sum_{\substack{\tau\in S_{n+1}^{(s,n+1)}\\\text{s.t. }\tau(k)=t}}
(-1)^s\sgn(\tau)
\\
&\quad\quad
\Big(\ldb {a^{\tau(1)}}_{\lambda_{\tau(1)}} \stackrel{s}{\check{\dots}}
\llbracket {a^s}_{\lambda_s}a^t\rrbracket_{\lambda_s+\lambda_t}
 \dots {a^{\tau(n)}}_{\lambda_{\tau(n)}}a^{\tau(n+1)}\rdb_L
\Big)_{\tau(s s+1\dots k-1)(\underline i,\underline j)}
\\
&
=\sum_{k=1}^{n-1}\sum_{1\leq s<t\leq n}
\sum_{\substack{\tau\in S_{n+1}^{(s,n+1)}\\ \text{s.t.} \tau(k)=s,\tau(k+1)=t}}
(-1)^k\sgn(\tau)
\\
&\quad\quad
\Big(\ldb {a^{\tau(1)}}_{\lambda_{\tau(1)}}\dots
\llbracket {a^{\tau(k)}}_{\lambda_{\tau(k)}}a^{\tau(k+1)}\rrbracket_{\lambda_{\tau(k)}+\lambda_{\tau(k+1)}}
\dots {a^{\tau(n)}}_{\lambda_{\tau(n)}}a^{\tau(n+1)}\rdb_L
\Big)_{\tau(\underline i,\underline j)}
\,.
\end{align*}
To derive the second equality we substituted, in the first sum, $\tau$ by $\tilde \tau (k\, k+1\dots s)$ with 
$\tilde\tau\in S_{n+1}^{(s,n+1)}$ such that $\tilde\tau(k)=s$ and $\tilde\tau(k+1)=t$, and, in the second sum, $\tau$ by $\tilde \tau(k-1\, k-2\dots s)$, where $\tilde\tau\in S_{n+1}^{(s,n+1)}$ is such that $\tilde\tau(k-1)=s$ and $\tilde\tau(k)=t$ and shifted the index of summation $k$.
Similarly, we get
\begin{align*}
&\eqref{20240906:eq3}
=\sum_{k=1}^{n-1}\sum_{1\leq t<s\leq n}
\sum_{\substack{\tau\in S_{n+1}^{(s,n+1)}\\ \text{s.t.} \tau(k)=s,\tau(k+1)=t}}
(-1)^k\sgn(\tau)
\\
&\quad\quad
\Big(\ldb {a^{\tau(1)}}_{\lambda_{\tau(1)}}\dots \ldb {a^{\tau(k)}}_{\lambda_{\tau(k)}}a^{\tau(k+1)}\rdb^H_{\lambda_{\tau(k)}+\lambda_{\tau(k+1)}}
\dots {a^{\tau(n)}}_{\lambda_{\tau(n)}}a^{\tau(n+1)}\rdb_L
\Big)_{\tau(\underline i,\underline j)}
\,.
\end{align*}
Combining \eqref{20240906:eq1}, \eqref{20240906:eq2} and \eqref{20240906:eq3} the LHS of \eqref{20240905:toprove4} is
\begin{align*}
&\sum_{k=1}^{n-1}\sum_{1\leq s\neq t\leq n}
\sum_{\substack{\tau\in S_{n+1}^{(s,n+1)}\\ \text{s.t.} \tau(k)=s,\tau(k+1)=t}}
(-1)^k\sgn(\tau)
\\
&\quad\quad
\Big(\ldb {a^{\tau(1)}}_{\lambda_{\tau(1)}}\dots
\llbracket {a^{\tau(k)}}_{\lambda_{\tau(k)}}a^{\tau(k+1)}\rrbracket_{\lambda_{\tau(k)}+\lambda_{\tau(k+1)}}
\dots {a^{\tau(n)}}_{\lambda_{\tau(n)}}a^{\tau(n+1)}\rdb_L
\Big)_{\tau(\underline i,\underline j)}
\,,
\end{align*}
which coincides with the RHS of \eqref{20240905:toprove4}.
\end{proof}
\begin{remark}
It is also possible to define a linear map $\tr:\widetilde{\Gamma}(\mc V)\to\widetilde{\Gamma}(\mc V_N)$ from the space of basic $n$-cochains over $\mc V$ to the space of basic $n$-cochains
over $\mc V_N$. For $f\in\mc V=\widetilde{\Gamma}^0(\mc V)$ we set $\tr(f)\in\mc V_N=\widetilde{\Gamma}^0(\mc V_N)$ as in \eqref{eq:trace0} and for $X=X_{\lambda_1,\dots,\lambda_n}\in\widetilde{\Gamma}^n(\mc V)$, $n\geq1$, let (cf. \eqref{eq:Pn})
$$
\bar\proj_n(X)_{\lambda_1,\dots,\lambda_n}
=\frac1n \sum_{s=0}^{n-1}(-1)^{s(n-s)}
\mult_{(s+1,s+2)}\circ \sigma^{s+1}\circ X_{\lambda_{\sigma^s(1)},\dots,\lambda_{\sigma^s(n)}}\circ \sigma^{-s}\,.
$$
Then
\begin{align*}
&\tr(X)_{\lambda_1,\dots,\lambda_n}(a^1_{i_1j_i},\dots,a^{n}_{i_nj_n})
\\
& =\sum_{\tau\in S_n^{(n)}}\sgn(\tau)
\bar\proj(X)_{\lambda_{\tau(1)},\dots,\lambda_{\tau(n)}}(a^{\tau(1)},\dots,a^{\tau(n)})_{\tau(\underline i,\underline j)}
\,,
\end{align*}
and we extend it to $\mc V_N^{\otimes n}$ by the Leibniz rules \eqref{eq:Leibnizmaps-comm}.
The map $\tr:\widetilde{\Gamma}(\mc V)\to\widetilde{\Gamma}(\mc V_N)$ induces a map between the basic dPVA cohomology defined in Definition \ref{def:bas-dPVA} and the basic PVA cohomology defined in Section \ref{sec:basic}. It also descends to a map between the corresponding reduced cohomologies.
\end{remark}

%%%
\section{\texorpdfstring{$\Gl_N$}{GLN}-invariance}\label{sec:GLN-PVA}
Let $\mc V$ be a differential algebra and let $\mc V_N$, $N\geq1$, be the associated commutative differential algebra defined in Section \ref{sec:dPVAtoPVA-1}.
In analogy with Chapter \ref{CH:rep-dPA} we have an action of $G=\Gl_N(\kk)$ on $\mc V_N$ defined on generators $a_{ij}\in\mc V_N$, $a\in\mc V$, $1\leq i,j\leq N$, as follows ($g\in G$)
\begin{equation}\label{eq:G-action}
g\cdot a_{ij}=\sum_{h,k=1}^N(g^{-1})_{ih}g_{kj}a_{hk}
\,,
\end{equation}
and extended to $\mc V_N$ by letting $g$ act as a differential algebra automorphism of $\mc V_N$, namely
$$
g\cdot( pq )=(g\cdot p)(g\cdot q)\,,
\qquad
g\cdot(\partial p)=\partial(g\cdot p)\,,
\qquad
p,q\in\mc V_N
\,.
$$
It is straightforward to check that this action is compatible with the defining relations \eqref{20130917:eq2} of $\mc V_N$.
\begin{remark}
To any $a\in\mc V$ we associate the matrix $M_a=(a_{ij})_{i,j=1}^N\in\Mat_{N\times N}(\mc V_N)$. Let
$g\cdot M_a=(g\cdot a_{ij})_{i,j=1}^N\in\Mat_{N\times N}(\mc V_N)$ denote the matrix whose entry $(i,j)$ is the image of the action of $g$ on $a_{ij}$.
Then, the action \eqref{eq:G-action} can be rewritten concisely in matrix form as
$$
g\cdot M_a=g^{-1}M_a g\,,
$$
that is, it corresponds to conjugation by the matrix element $g^{-1}$.
\end{remark}
Recall from Section \ref{sec:dPVAtoPVA} that we have a map
$\tr:\mc V_\sharp\to\mc V_N/\partial\mc V_N$ defined by
$\tr(\tint f)=\tint \tr f=\sum_{i=1}^N\tint f_{ii}=\tint \tr(M_f)$.
Hence, for every $g\in G$ we have that $g\cdot \tr(\tint f)=\tint \tr(g^{-1}M_f g)=\tint\tr(f)=\tr(\tint f)$, so that
$\tr(\mc V_\sharp)\subset\left(\mc V_N/\partial\mc V_N\right)^G$. Hence,
$\tr(C^0(\mc V))\subset C^0(\mc V_N)^G$.
\begin{proposition}\label{20250627:prop1}
Let $\mc V$ be a differential algebra and let $\mc V_N$, $N\geq1$, be the associated commutative differential algebra. Let $\tr:C(\mc V)\to C(\mc V_N)$ be the trace map defined in \eqref{Eq:TrBrRep-PVA}. Then,
$\tr(C(\mc V))\subset C(\mc V_N)^{G}$.
\end{proposition}
\begin{proof}
We show that $\tr(C^n(\mc V))\subset C^n(\mc V_N)^{G}$,
for every $n\in\mb Z_{\geq0}$. The case $n=0$ has already been proved, so we are left to consider the case $n\geq1$.
Let $\ldb-_{\lambda_1}-\dots-_{\lambda_{n-1}}-\rdb\in C^n(\mc V)$ be an $n$-fold $\lambda$-bracket on $\mc V$. For every $g\in G$
and $a^{1}_{i_1j_1},\dots, a^{n}_{i_nj_n}\in\mc V_N$, using \eqref{Eq:TrBrRep-PVA},
the definition \eqref{eq:G-action} of the action of $G$ on $\mc V_N$
and the definition \eqref{eq:G-action-C} of the action of $G$ on $C^n(\mc V_N)$, we have
\begin{align*}
&(g\cdot\tr(\ldb-_{\lambda_1}-\dots-_{\lambda_{n-1}}-\rdb))(a^1_{i_1j_1},\dots,a^n_{i_nj_n})
\\
&=\sum_{\tau\in S_{n-1}}
\sum_{h_1,k_1,\dots, h_n,k_n=1}^N\left(\prod_{l=1}^ng_{i_lh_l}(g^{-1})_{k_lj_l}\right)\sgn(\tau)
\\
&\quad\quad\quad\quad\quad\quad\quad\quad
g\cdot \ldb {a^{\tau(1)}}_{\lambda_{\tau(1)}}\dots {a^{\tau(n-1)}}_{\lambda_{\tau(n-1)}}a^n\rdb_{\tau(\underline h,\underline k)}
\,.
\end{align*}
Let us simply denote
$$
\ldb {a^{\tau(1)}}_{\lambda_{\tau(1)}}\dots {a^{\tau(n-1)}}_{\lambda_{\tau(n-1)}}a^n\rdb
=x^1\otimes\cdots\otimes x^n\in\mc V^{\otimes n}[\lambda_1,\dots,\lambda_{n-1}]
\,.
$$
Then, using the notation \eqref{Eq:Not-RepIndex}
and the properties of the action of $G$ on $\mc V_N$,
the above equation becomes 
\begin{align*}
&(g\cdot\tr(\ldb-_{\lambda_1}-\dots-_{\lambda_{n-1}}-\rdb))(a^1_{i_1j_1},\dots,a^n_{i_nj_n})
\\
&=\sum_{\tau\in S_{n-1}}\sgn(\tau)
\sum_{h_1,k_1,\dots, h_n,k_n=1}^N\Bigg(\prod_{l=1}^ng_{i_lh_l}(g^{-1})_{k_lj_l}\Bigg)
x^1_{h_nk_{\tau(1)}}x^2_{h_{\tau(1)}k_{\tau(2)}}\dots x^n_{h_{\tau(n-1)}k_n}
\\
&=\sum_{\tau\in S_{n-1}}\sgn(\tau)
\sum_{h_1,k_1,\dots, h_n,k_n=1}^N
g\cdot\big(g_{i_nh_n}(g^{-1})_{k_{\tau(1)}j_{\tau(1)}}x^1_{h_nk_{\tau(1)}}\times
\\
&\quad\quad\quad
\times
g_{i_{\tau(1)}h_{\tau(1)}}(g^{-1})_{k_{\tau(2)}j_{\tau(2)}}
x^2_{h_{\tau(1)}k_{\tau(2)}}\dots
g_{i_{\tau(n-1)}h_{\tau(n-1)}}(g^{-1})_{k_nj_n} x^n_{h_{\tau(n-1)}k_n}
\big)
\\
&=\sum_{\tau\in S_{n-1}}\sgn(\tau)
g\cdot\big((g^{-1}\cdot x^1_{i_nj_{\tau(1)}})
(g^{-1}\cdot x^2_{i_{\tau(1)}j_{\tau(2)}})\dots
(g^{-1}\cdot x^n_{i_{\tau(n-1)}j_n})
\big)
\\
&=\sum_{\tau\in S_{n-1}}\sgn(\tau)
x^1_{i_nj_{\tau(1)}}x^2_{i_{\tau(1)}j_{\tau(2)}}\dots
x^n_{i_{\tau(n-1)}j_n}
\\
&=
\tr(\ldb-_{\lambda_1}-\dots-_{\lambda_{n-1}}-\rdb)(a^1_{i_1j_1},\dots,a^n_{i_nj_n})
\,,
\end{align*}
from which follows that $\tr(\ldb-_{\lambda_1}-\dots-_{\lambda_{n-1}}-\rdb)\in C^n(\mc V_N)^G$.
\end{proof}
Let us assume that $\mc V$ is a dPVA
with $2$-fold $\lambda$-bracket $\llbracket-_{\lambda}-\rrbracket$,
then the commutative differential algebra $\mc V_N$ is a PVA $\lambda$-bracket $[-_{\lambda}-]$ given by \eqref{20240829:eq1}. Then, the action of $G$ on $\mc V_N$ is given by PVA automorphisms. Indeed, we have ($g\in G,a_{ij},b_{hk}\in\mc V_N$)
\begin{align*}
&g\cdot \{a_{ij}{}_{\lambda} b_{hk}\}_{H_N}
=\left(g\cdot(\llbracket a_\lambda b\rrbracket')_{hj}\right)\left(g\cdot(\llbracket a_\lambda b\rrbracket'')_{ik}\right)
\\
&=\sum_{\alpha,\beta,\gamma,\delta=1}^N
(g^{-1})_{h\alpha}g_{\beta j}(\llbracket a_{\lambda}b\rrbracket')_{\alpha\beta}
(g^{-1})_{i\gamma}g_{\delta k}(\llbracket a_{\lambda}b\rrbracket')_{\gamma\delta}
\\
&=\sum_{\alpha,\beta,\gamma,\delta=1}^N
(g^{-1})_{h\alpha}g_{\beta j}
(g^{-1})_{i\gamma}g_{\delta k}
[a_{\gamma\beta}{}_{\lambda}b_{\alpha\delta}]
=[g\cdot a_{ij}{}_\lambda g\cdot b_{hk}]
\,,
\end{align*}
thus $g\cdot [p_\lambda q]=[g\cdot p_{\lambda}g\cdot q]$ for every
$p,q\in\mc V_N$ by  sesquilinearity, Leibniz rules and the fact that $G$ acts 
by automorphisms of the differential algebra structure of $\mc V_N$.

By Proposition \ref{20250627:prop1} and the results of Section \ref{sec:PVAcoh-inv}
we get the following.
\begin{theorem}
Let $\mc V$ be a dPVA
and let $\mc V_N$, $N\geq1$, be the corresponding PVA defined using \eqref{20240829:eq1}.
Let $d$ be the differential of $C(\mc V)$ defined by \eqref{eq:dP0} and \eqref{eq:dP}, and let
$d_{N}$ be the differential of $ C(\mc V_N)$ defined by \eqref{eq:dH-poly-0}
and \eqref{eq:dH-poly-n}. Then, the trace map \eqref{Eq:dP-rep2-PVA} is a morphism of complexes 
\begin{equation} \label{Eq:dP-rep2-PVA-bis}
    \tr : (C(\mc V),\dd) \longrightarrow 
    (C(\mc V_N)^G, (-1)^{\bullet}\dd_{N}|_{C(\mc V_N)^G}),  
\end{equation}
which descends to a linear map $\dPVH(\mc V)\to \PVH_G(\mc V_N)$ in cohomology. 
Furthermore, we have a restriction map
$C^n(\mc V_N)^G\to C^n(\mc V_N^G)$, for every $n\geq1$, which induces the linear map in cohomology
$\PVH^n_G(\mc V_N)\to \PVH^n(\mc V_N^G)$
(see \eqref{diag0}).
\end{theorem}

%%%%%%%%%%% NEW PART %%%%%%%%%%%%%%%
%%%%%%%%%%% NEW PART %%%%%%%%%%%%%%%
%%%%%%%%%%% NEW PART %%%%%%%%%%%%%%%
%%%%%%%%%%% NEW PART %%%%%%%%%%%%%%%
%%%%%%%%%%% NEW PART %%%%%%%%%%%%%%%
%%%%%%%%%%% NEW PART %%%%%%%%%%%%%%%
%%%%%%%%%%% NEW PART %%%%%%%%%%%%%%%

\addtocontents{toc}{\protect\newpage}%it puts part 3 at the beginning of page 5 rather than at the end of page 4 in the table of contents
\part{Relations under the jet and quotient functors}

%%%%%%%%%%% NEW CHAPTER %%%%%%%%%%%%%%%
%%%%%%%%%%% NEW CHAPTER %%%%%%%%%%%%%%%
%%%%%%%%%%% NEW CHAPTER %%%%%%%%%%%%%%%
%%%%%%%%%%% NEW CHAPTER %%%%%%%%%%%%%%%

\chapter[From PA to PVA cohomology, and back]{From Poisson algebra to Poisson vertex algebra cohomology, and back}
\label{Ch:PA-PVA}

We show how to relate the Poisson cohomology (using the Chevalley-Eilenberg complex) and the variational PVA cohomology using the jet and quotient functors. 
While this relation is natural and certainly known to experts, we include the presentation for completeness since we have not found references handling this topic.  

\section{Constructions with the jet functor}
We consider the jet functor from the category of finitely generated commutative algebras to the one of differential commutative algebras introduced in Example \ref{exa:jetPVA}. 
It associates to any algebra $A$ a universal representative denoted $J_\infty A$, called the jet algebra, cf. \cite[\S1.1]{AM}. 
To be explicit, let $R_\ell^{(0)}=\kk[x_1,\dots,x_\ell]$ be the commutative algebra of polynomials in $\ell$ variables $x_1,\dots,x_\ell$. 
Then  $J_\infty R_\ell^{(0)}=R_\ell$ for 
$$R_\ell=\kk[x_1^{(k)},\dots,x_\ell^{(k)}\mid k\geq 0]$$ 
equipped with the differential satisfying $\partial(x_j^{(k)})=x_j^{(k+1)}$, and we have an inclusion of commutative algebras $R_\ell^{(0)} \hookrightarrow R_\ell$ given by $x_j \mapsto x_j^{(0)}$. 
Note that $R_\ell=J_\infty R_\ell^{(0)}$ is the commutative analogue of $\mc R_\ell$ from Chapter \ref{sec:PVAdiff}. 

Any finitely generated commutative algebra $A$ can be realized as a quotient $A=R_\ell^{(0)}/I$. 
In this case the jet algebra can be realized explicitly as the quotient 
$$
J_\infty A= R_\ell/\langle I\rangle_{\partial}\,, 
\qquad 
\langle I\rangle_{\partial}=\sum_{n=0}^{\infty}R_\ell\, \partial^nI\,,
$$
with $\langle I\rangle_{\partial}$ the differential ideal generated by $I$.  
As in the case of $R_\ell^{(0)}$, there is a natural inclusion $\iota: A \hookrightarrow J_\infty A$. 
Arakawa showed that any Poisson bracket on $A$ extends to a PVA $\lambda$-bracket on $J_\infty A$ \cite[\S2.3]{Ar12}. The next result is an analogous extension for multiderivations. 
\begin{lemma}\label{lem:jetP}
Let $A$ be a finitely generated commutative algebra. For every skewsymmetric $n$-linear derivation
$Q\in\mf X^n(A)$, there exists a unique $n$-$\lambda$-bracket $J_\infty Q=\{-_{\lambda_1}-\dots-_{\lambda_{n-1}}-\}\in C^n(J_{\infty} A)$ making the following diagram commute:
\begin{equation}\label{diag:jetP}
% \text{this diagram has been commented to speed up compiling}
% \begin{comment}
\begin{tikzcd}
A^{n}\arrow[r,"Q"] \arrow[d,hook,"\iota"] & A\arrow[d,hook,"\iota"]
\\
(J_{\infty}A)^n \arrow[r,"J_\infty Q"] & J_{\infty}A[\lambda_1,\dots,\lambda_{n-1}]
\end{tikzcd}
% \end{comment}
\end{equation}
\end{lemma}
\begin{proof} 
For $n\geq 1$, the $n$-$\lambda$-bracket $J_\infty Q$ is explicitly defined as follows. Let $[x_j]_I\in A$
be the classes of $x_j\in R_\ell^{(0)}$, $j=1,\dots, \ell$, in the quotient algebra. 
Let
$$Q([x_{j_1}]_I,\dots,[x_{j_n}]_I)=[f_{j_1,\dots,j_n}]_I, \qquad j_1,\dots,j_n=1,\dots,\ell$$
be the images on an $n$-uple of these classes via the $n$-derivation $Q$, and fix representatives $f_{j_1,\dots,j_n}\in R_{\ell}^{(0)}$. 
If $p_1,\dots,p_n\in R_\ell$ and if $c_1=[p_1]_{\langle I\rangle_\partial},\dots,c_n=[p_n]_{\langle I\rangle_\partial}$ are the corresponding classes in $J_\infty A$, then
\begin{align}
&\{c_1{}_{\lambda_1}\dots c_{n-1}
{}_{\lambda_{n-1}}c_n\}
\nonumber 
\\
&=\Bigg[
\sum_{\substack{j_1,\dots,j_n=1,\dots,\ell\\m_1,\dots,m_n\in\mb Z_{\geq0}}}
\frac{\partial p_n}{\partial x_{j_n}^{(m_n)}}(\lambda_1+\dots+\lambda_{n-1}+\partial)^{m_n}
f_{j_1,\dots,j_n}\\
&\quad\quad\qquad
\bigg((-\lambda_1-\partial)^{m_1}\frac{\partial p_1}{\partial x_{j_1}^{(m_1)}}\bigg)
\dots \bigg((-\lambda_{n-1}-\partial)^{m_{n-1}}\frac{\partial p_{n-1}}{\partial x_{j_{n-1}}^{(m_{n-1})}}\bigg)
\Bigg]_{\langle I\rangle_\partial}
. \nonumber 
\end{align}
It is easy to check that $J_\infty Q$ is well defined by the above formula, that it is indeed an element of $C^n(J_\infty A)$ and that $J_\infty Q\circ \iota=\iota\circ Q$. The uniqueness of such
$n$-$\lambda$-bracket is obvious.
\end{proof}
\begin{remark}
We identify $A\simeq\iota(A)\subset J_{\infty}A$. Then, since the diagram \eqref{diag:jetP} is commutative,
the $n$-$\lambda$-bracket provided by Lemma \ref{lem:jetP} satisfies  
\begin{equation} \label{Eq:PAPVA-ext1}
  \{a_1 {}_{\lambda_1}\dots{}_{\lambda_{n-1}} a_n\}
  =Q(a_1,\dots,a_n),  
\end{equation}
when evaluated on $a_1,\ldots, a_n \in A$.
\end{remark}
\begin{remark}
Any finitely generated differential commutative algebra $V$ can be realized as a quotient 
$V=R_\ell/I_V$ for some $\ell \geq 1$ and some differential ideal $I_V\subset R_\ell$. 
Based on \cite[Eq.~(9.3)]{DSK13}, any $n$-$\lambda$-bracket on $V$ is also determined by the value that it takes on the generators of $V$ given by the classes $[x_j]_{I_V} \in V$ of the generators of $R_\ell^{(0)}$. 
Indeed, if $c_1,\ldots,c_n\in V$ admit representatives $p_1,\ldots,p_n \in R_\ell$,  
\begin{align}
&\{c_1{}_{\lambda_1}\dots c_{n-1}
{}_{\lambda_{n-1}}c_n\}
\nonumber 
\\
&=\Bigg[
\sum_{\substack{j_1,\dots,j_n=1,\dots,\ell\\m_1,\dots,m_n\in\mb Z_{\geq0}}}
\frac{\partial p_n}{\partial x_{j_n}^{(m_n)}}(\lambda_1+\dots+\lambda_{n-1}+\partial)^{m_n}
\nonumber  \\
&\quad\qquad
f_{j_1,\dots,j_n}(\lambda_1+x_1,\ldots,\lambda_{n-1}+x_{n-1})
\bigg((-\lambda_1-x_1)^{m_1}\big|_{x_1=\partial}\frac{\partial p_1}{\partial x_{j_1}^{(m_1)}}\bigg)
\dots  \nonumber \\
&\qquad
\dots \bigg((-\lambda_{n-1}-x_{n-1})^{m_{n-1}} \big|_{x_{n-1}=\partial}\frac{\partial p_{n-1}}{\partial x_{j_{n-1}}^{(m_{n-1})}}\bigg)
\Bigg]_{I_V}
. \label{Eq:PVAnMaster}
\end{align}
where $f_{j_1,\dots,j_n}(\lambda_1,\ldots,\lambda_{n-1})\in R_\ell[\lambda_1,\ldots,\lambda_{n-1}]$ is a representative of 
$$\{ [x_{j_1}]_{I_V}{}_{\lambda_1}\dots_{\lambda_{n-1}}[x_{j_n}]_{I_V} \} \in V[\lambda_1,\ldots,\lambda_{n-1}], \qquad j_1,\dots,j_n=1,\dots,\ell\,.$$
\end{remark}

Let us now assume that $(A,\br{-,-})$ is a (still finitely generated) Poisson algebra.  
Thanks to Lemma \ref{lem:jetP}, we get a morphism $J_\infty: \mf X(A) \to C(J_\infty A)$ after completing it in degree $0$ by 
\begin{equation} \label{Eq:PAPVA-ext0}
  A \to C^0(J_\infty A)= (J_\infty A)_\sharp, \quad a\mapsto \tint \iota(a)\,.  
\end{equation}
Now, we can view $\mf X(A)$ as a complex for the Chevalley-Eilenberg differential \eqref{Eq:Diff-Pcoh1}, which we denote as $\delta_A$. 
As already mentioned, Arakawa \cite{Ar12} proved that $J_\infty \br{-,-}$ defines a PVA structure on $J_\infty A$ (this is now a  consequence of Lemma \ref{lem:jetP}). 
Hence we can view $C(J_\infty A)$ as a complex for the variational differential \eqref{eq:dH-poly-n}, which we denote as $\delta_{J_\infty A}$.   
To ease notation, we denote $J_\infty \br{-,-}$ as $\br{-_\lambda -}$, which satisfies due to \eqref{Eq:PAPVA-ext1}
\begin{equation} \label{Eq:PAPVA-br}
    \br{a_\lambda b} = \br{a,b} \,, \quad \text{for all } a,b\in A\simeq\iota(A)\subset J_\infty A\,.
\end{equation}

\begin{proposition}\label{Pr:PAPVA1}
For every $Q\in \mf X^n(A)$, $n\in\mb Z_{\geq0}$, we have
\begin{equation}\label{PAPVA-comm}
J_\infty(\delta_A(Q))= \delta_{J_\infty A}(J_\infty(Q))
\,.
\end{equation}
In particular, we get a morphism of complexes 
$$J_\infty: (\mf X(A) , \delta_A) \to (C(J_\infty A),\delta_{J_\infty A})\,.$$
\end{proposition}
\begin{proof}
    In view of the Master Formula \eqref{Eq:PVAnMaster} for $n$-$\lambda$-brackets, it suffices to check the equality \eqref{PAPVA-comm} on the generators $x_j$, $1\leq j \leq \ell$, viewed as elements of $J_\infty A$. 
(We use the same notation for these elements in $A$, which we previously denoted as $[x_j]_I\in A$ and $[x_j]_{\langle I\rangle_\partial}\in J_\infty A$.) 

For $n=0$, $Q=a\in A$ and we find 
\begin{align*}
 J_\infty(\delta_A(a))(x_j)=   
 \br{x_j,a},
\end{align*}  
where we used \eqref{Eq:Diff-Pcoh1} with \eqref{Eq:PAPVA-ext1}, then
\begin{align*} 
\delta_{J_\infty A}(J_\infty(a)) (x_j)
&=\br{x_j\, {}_{-\lambda-\partial} a} \big|_{\lambda=0} = \br{x_j,a}, 
\end{align*}  
where we used \eqref{eq:dH-poly-0} with \eqref{Eq:PAPVA-ext0} and skewsymmetry for the first equality, 
and \eqref{Eq:PAPVA-br} for the second equality.

Next, we look at $n\geq 1$. Fix $j_1,\ldots,j_{n+1}\in \{1,\ldots,\ell\}$ and evaluate the LHS of \eqref{PAPVA-comm} as follows: 
\begin{align*}
J_\infty(\delta_A(Q))(x_{j_1} ,\ldots&, x_{j_{n+1}} )  
=
\sum_{1\leq i \leq n+1} (-1)^{i+1} \br{x_{j_i}, Q(x_{j_1},\stackrel{i}{\check{\dots}}, x_{j_{n+1}}) } \\
&+  \sum_{1\leq i<k \leq n+1} (-1)^{i+k}  Q(\br{x_{j_i},x_{j_k}},x_{j_1},\stackrel{i}{\check{\dots}},\stackrel{k}{\check{\dots}}, x_{j_{n+1}})\,.
\end{align*}
where we used \eqref{Eq:PAPVA-ext1} and \eqref{Eq:Diff-Pcoh1} (with $a_i=x_{j_i}$). 
For the RHS, we use \eqref{eq:dH-poly-n} to obtain 
\begin{align*}
&\delta_{J_\infty A}(J_\infty(Q))(x_{j_1} ,\ldots, x_{j_{n+1}} )  \\
&=
\sum_{1\leq i \leq n+1} (-1)^{i+1} \big|_{\lambda_{n+1}=\lambda_{n+1}^\dagger} 
\big\{x_{j_i}{}_{\lambda_i} (J_\infty Q)(x_{j_1}{}_{\lambda_1}\stackrel{i}{\check{\dots}}
x_{j_n} {}_{\lambda_n} x_{j_{n+1}}) \big\} \\
&+  \sum_{1\leq i<k \leq n+1} (-1)^{i+k} \big|_{\lambda_{n+1}=\lambda_{n+1}^\dagger}   
(J_\infty Q)(\br{x_{j_i} {}_{\lambda_i} x_{j_k}} {}_{\lambda_i+\lambda_k} 
x_{j_1} {}_{\lambda_1}\stackrel{i}{\check{\dots}}\stackrel{k}{\check{\dots}}
x_{j_n} {}_{\lambda_n} x_{j_{n+1}})\,.
\end{align*}
Making use of \eqref{Eq:PAPVA-ext1} and \eqref{Eq:PAPVA-br}, we get the same expression as for the LHS. 
\end{proof}

In conclusion, we proved the following result.
\begin{corollary}\label{Cor:PAPVA1}
For every $n\in\mb Z_{\geq0}$ the linear map $J_\infty: \mf X^n(A) \to C^n(J_\infty A)$ given
by Lemma \ref{lem:jetP} (or \eqref{Eq:PAPVA-ext0} for $n=0$) gives a morphism in cohomology 
$$
\coH_{CE}^n(A)\rightarrow \PVH^n(J_\infty A)
\,.
$$ 
\end{corollary}

\section{Constructions with the quotient functor}

We consider the quotient functor %$\QQ$
from the category of (finitely generated) differential commutative algebras to the one of commutative algebras as in Example \ref{exa:PVAtoPA}.   
Recall that it is given by $V\mapsto q(V)=V/\langle \partial V \rangle$ with the projection $\pi:V \to q(V)$, and that we can produce a Poisson bracket from a PVA $\lambda$-bracket. The next result is a straightforward generalization. 
\begin{lemma}\label{lem:QuotPV}
Let $V$ be a differential commutative algebra. For every 
$n$-$\lambda$-bracket $c:=\{-_{\lambda_1}-\dots-_{\lambda_{n-1}}-\}\in C^n(V)$, there exists a unique 
skewsymmetric $n$-linear derivation 
$q(c)\in\mf X^n(q(V))$ making the following diagram commute:
\begin{equation}\label{diag:QuotPV}
% \text{this diagram has been commented to speed up compiling}
% \begin{comment}
\begin{tikzcd}
V^{n}\arrow[r,"c"] \arrow[d,twoheadrightarrow,"\pi"] & V[\lambda_1,\dots,\lambda_{n-1}] \arrow[d,twoheadrightarrow,"\pi"]
\\
q(V)^n \arrow[r,"q(c)"] & q(V)
\end{tikzcd}
% \end{comment}
\end{equation}
\end{lemma}
\begin{proof}
We introduce $q(c)$ by setting 
\begin{equation} \label{Eq:PAPVA-qc}
 q(c)(\pi(a_1),\ldots,\pi(a_n))
 = \pi( \{a_1{}_{\lambda_1}\dots a_{n-1}{}_{\lambda_{n-1}} a_n\} 
 |_{\lambda_1=\ldots=\lambda_{n-1}=0} )\, ,
\end{equation}
for arbitrary $a_1,\ldots,a_n \in V$. One can then check that this gives a well-defined element of $\mf X^n(q(V))$, which is the unique one making the diagram commute.
\end{proof}
For $n=0$, we use the map $C^0(V) \to q(V)$, $\tint a \mapsto \pi(a)$, obtained by factoring 
$\pi:V\to V/\langle \partial V \rangle$ through the linear map $\tint : V \mapsto V_\sharp=V/\partial V$. 
Assuming from now on that $(V,\br{-_\lambda-})$ is a PVA inducing the Poisson algebra structure 
$(q(V),\br{-,-})$, 
we can view the $C(V)$ and $\mf X(q(V))$ as complexes equipped with the differentials 
\eqref{eq:dH-poly-n} and \eqref{Eq:Diff-Pcoh1}, which we denote as $\delta_V$ and $\delta_{q(V)}$, respectively. 
\begin{proposition}\label{Pr:PAPVA2}
For every $c\in C^n(V)$, $n\in\mb Z_{\geq0}$, we have
\begin{equation}\label{PAPVA-comm2}
q(\delta_V(c))= \delta_{q(V)}(q(c)) \,.
\end{equation}
In particular, we get a morphism of complexes 
$$q: (C(V), \delta_V) \to (\mf X(q(V)) ,\delta_{q(V)})\,.$$
\end{proposition}
\begin{proof}
For $n=0$ with $c=\tint a\in V_\sharp$, evaluating both sides of \eqref{PAPVA-comm2} on $\pi(b)\in q(V)$ gives 
\begin{align*}
q(\delta_V(\tint a))(\pi(b))&= - \pi( \br{a_\lambda b} |_{\lambda =0}) \,, \\
\delta_{q(V)}(q(\tint a))(\pi(b))&=  \br{\pi(b),\pi(a)}\,,
\end{align*}
where we used \eqref{Eq:PAPVA-qc} with \eqref{eq:dH-poly-0} in the first equality with arbitrary lifts $a,b\in V$, then \eqref{Eq:Diff-Pcoh1} with $q(\tint a)=\pi(a)$ for the second equality. 
Both expressions are equal by \eqref{Eq:br-PVAquotient}. 
Since $\pi(b)$ is an arbitrary element in $q(V)$, we get $q(\delta_V(\tint a))=\delta_{q(V)}(q(\tint a))$. 

For $n\geq 1$, given $c_{\lambda_1,\dots,\lambda_{n-1}}\in C^n(V)$ and 
$a_1,\ldots,a_{n+1}\in V$, we find 
\begin{align*}
&q(\delta_V(c))(\pi(a_1),\ldots,\pi(a_{n+1})) \\
&=
\sum_{1\leq i \leq n+1} (-1)^{i+1} \pi\bigg( 
\big\{a_{i}{}_{\lambda_i} c_{\lambda_1,\stackrel{i}{\check{\dots}},\lambda_n}(a_1,\stackrel{i}{\check{\dots}}
a_{n+1}) \big\} |_{\lambda_1=\ldots=\lambda_{n+1} =0}\bigg) \\
&+  \sum_{1\leq i<k \leq n+1} (-1)^{i+k} \pi\bigg( 
c_{\lambda_i+\lambda_k,\lambda_1,\stackrel{i,k}{\check{\dots}},\lambda_n}(\br{ a_i {}_{\lambda_i} a_k},
a_1,\stackrel{i,k}{\check{\dots}}a_{n+1}) \big\} 
|_{\lambda_1=\ldots=\lambda_{n+1} =0}  \bigg)\,, 
\end{align*}
where we used \eqref{Eq:PAPVA-qc} and \eqref{eq:dH-poly-n}. Then, we obtain 
\begin{align*}
&\delta_{q(V)}(q(c))(\pi(a_1),\ldots,\pi(a_{n+1})) \\
&=
\sum_{1\leq i \leq n+1} (-1)^{i+1}  
\big\{\pi(a_{i}) , q(c)(\pi(a_1),\stackrel{i}{\check{\dots}} \pi(a_{n+1})) \big\} \\
&+  \sum_{1\leq i<k \leq n+1} (-1)^{i+k} \, 
q(c)\big(\br{ \pi(a_i) , \pi(a_k)} ,
\pi(a_1),\stackrel{i}{\check{\dots}}\stackrel{k}{\check{\dots}},\pi(a_{n+1}) \big)\,, 
\end{align*}
using \eqref{Eq:Diff-Pcoh1}. 
Both expressions perfectly coincide if one uses \eqref{Eq:br-PVAquotient} and \eqref{Eq:PAPVA-qc}. 
This entails the desired equality in \eqref{PAPVA-comm2}. 
\end{proof}
This time, we deduce the following result.
\begin{corollary}\label{Cor:PAPVA2}
For every $n\in\mb Z_{\geq0}$ the linear map $q: C^n(V) \to \mf X^n(q(V))$ given
by Lemma \ref{lem:QuotPV} (or $\tint a \mapsto \pi(a)$ for $n=0$) gives a morphism in cohomology 
$$
\PVH^n(V) \to \coH_{CE}^n(q(V)) \,.
$$ 
\end{corollary}

In Example \ref{exa:jetPVA}, we have seen that composing the quotient functor after the jet functor gives back the original Poisson algebra, i.e. $q(J_\infty(A))=A$ (up to isomorphism).  
Similarly, the composition 
$$
\mf X(A) \stackrel{J_\infty}{\longrightarrow} C(J_\infty A) \stackrel{q}{\longrightarrow} 
\mf X(q(J_\infty(A))) = \mf X(A)
$$
obtained from Lemmas \ref{lem:jetP} and \ref{lem:QuotPV} is easily seen to be an isomorphism. 
Combining Corollaries \ref{Cor:PAPVA1} and \ref{Cor:PAPVA2}, we get the following result in cohomology. 
\begin{theorem}\label{20250822:thm1}
 The composition 
 $$
\coH_{CE}(A) \longrightarrow \PVH(J_\infty A) \longrightarrow \coH_{CE}(A)
$$   
is an isomorphism. In particular, the linear map $\coH_{CE}(A) \to \PVH(J_\infty A)$ is injective,  
and $\PVH(J_\infty A) \to \coH_{CE}(A)$ is surjective. 
\end{theorem}

%%%%%%%%%%% NEW CHAPTER %%%%%%%%%%%%%%%
%%%%%%%%%%% NEW CHAPTER %%%%%%%%%%%%%%%
%%%%%%%%%%% NEW CHAPTER %%%%%%%%%%%%%%%
%%%%%%%%%%% NEW CHAPTER %%%%%%%%%%%%%%%

\chapter[From \MakeLowercase{d}PA to \MakeLowercase{d}PVA cohomology, and back]{From double Poisson algebra to double Poisson vertex algebra cohomology, and back}
\label{Ch:dPA-dPVA}

We describe a relation between the completed double Poisson cohomology  and the variational double Poisson vertex algebra cohomology using the jet and quotient functors. 
We also show that this is compatible with the constructions in the commutative case described in Chapter \ref{Ch:PA-PVA}. 

\section{Constructions with the jet functor}
We consider the \emph{associative} jet functor from the category of finitely generated algebras to the one of differential algebras introduced in Example \ref{exa:jetdPVA}. 
It associates to any algebra $\cA$ a universal representative denoted $J_\infty \cA$, called the jet algebra, cf. \cite[\S4.2]{BFM}. 
To be explicit, let $\mc R_\ell^{(0)}=\kk\langle u_1,\dots,u_\ell\rangle$ be the algebra of polynomials in $\ell$ (noncommutative) variables $u_1,\dots,u_\ell$.  
Then  
$$J_\infty \mc R_\ell^{(0)}=\kk\big\langle u_1^{(k)},\dots,u_\ell^{(k)}\big| k\geq 0 \big\rangle$$ 
with the derivation defined by $\partial(u_j^{(k)})=u_j^{(k+1)}$ and extended to $\mc V$ by the Leibniz rule. 
This is precisely the algebra $\mc R_\ell$ from Chapter \ref{sec:PVAdiff}. 
We have an inclusion of algebras $\mc R_\ell^{(0)} \hookrightarrow \mc R_\ell$ given by $u_j \mapsto u_j^{(0)}$. 
Any finitely generated algebra $\cA$ can be realized as a quotient $\cA=\mc R_\ell^{(0)}/\mc I$ for some two-sided ideal $\mc I$.  
In this case the jet algebra can be realized explicitly as the quotient 
$$
J_\infty \cA =\mc R_\ell/\langle \mc I\rangle_{\partial}\,, 
\qquad 
\langle \mc I\rangle_{\partial}=\sum_{n=0}^{\infty}(\mc R_\ell)\,\partial^n\mc I\, (\mc R_\ell)\,,
$$
with $\langle \mc I\rangle_{\partial}$ the two-sided differential ideal generated by $\mc I$.  
Again, there is a natural inclusion $\iota: \cA \hookrightarrow J_\infty \cA$. 
It is shown in \cite[Lem.~4.5]{BFM} that any double Poisson bracket on $\cA$ extends to a dPVA $2$-fold $\lambda$-bracket on $J_\infty \cA$. The next result is an analogous extension for $n$-brackets.

\begin{lemma}\label{lem:jetdP}
Let $\mc A$ be a finitely generated algebra. For every $Q=\ldb-\rdb\in\BRA(\mc A)_n$, 
there exists a unique $n$-fold $\lambda$-bracket $J_\infty Q=\ldb-_{\lambda_1}-\dots-_{\lambda_{n-1}}-\rdb\in C^n(J_{\infty} \mc A)$ making the following diagram commute:
\begin{equation}\label{diag:jetdP}
% \text{this diagram has been commented to speed up compiling}
% \begin{comment}
\begin{tikzcd}
\mc A\arrow[r,"Q"] \arrow[d,hook,"\iota"] & \mc A\arrow[d,hook,"\iota"]
\\
J_{\infty}\mc A \arrow[r,"J_\infty Q"] & (J_{\infty}\mc A)^{\otimes n}[\lambda_1,\dots,\lambda_{n-1}]
\end{tikzcd}
% \end{comment}
\end{equation}
\end{lemma}
\begin{proof}
Similarly to the proof of Lemma \ref{lem:jetP} we define the $n$-fold $\lambda$-bracket $J_\infty Q$ explicitly as follows. Let $[x_j]_{\mc I}\in \mc A$
be the classes of $x_j\in \mc R_\ell^{(0)}$, $j=1,\dots, \ell$, in the quotient algebra, let
$$\ldb [x_{j_1}]_{\mc I},\dots,[x_{j_n}]_{\mc I}\rdb=[f_{j_1,\dots,j_n}]_{\mc I}, \qquad j_1,\dots,j_n=1,\dots,\ell$$
be the images on an $n$-uple of these classes via the $n$-bracket $Q$, with fixed representative $f_{j_1,\dots,j_n}\in (\mc R_{\ell}^{(0)})^{\otimes n}$. If $p_1,\dots,p_n\in J_{\infty}\mc R_\ell$ and if $c_1=[p_1]_{\langle \mc I\rangle_\partial},\dots,c_n=[p_n]_{\langle \mc I\rangle_\partial}$ are the corresponding classes in $J_\infty \mc A$, then (cf. equation \eqref{eq:QS})
\begin{equation}
 \begin{aligned}
\ldb c_1{}_{\lambda_1}\dots c_{n-1}&{}_{\lambda_{n-1}}c_n\rdb
=\Bigg[
\sum_{\substack{j_1,\dots,j_n=1,\dots,\ell\\m_1,\dots,m_n\in\mb Z_{\geq0}}}
\frac{\partial p_n}{\partial u_{j_n}^{(m_n)}}
(\lambda_1+\dots+\lambda_{n-1}+\partial)^{m_n}
\\
&
\bullet\Bigg(\Bigg(\dots\Bigg(
f_{j_1,\dots,j_n}\bullet_{(1,2)} (-\lambda_1-\partial)^{m_1}\bigg(\frac{\partial p_1}{\partial u_{i_1}^{(m_1)}}\bigg)^{\sigma}\Bigg)
\\
&\dots 
\Bigg)\bullet_{(n-1,n)}(-\lambda_{n-1}-\partial)^{m_{n-1}}\bigg(\frac{\partial p_{n-1}}{\partial u_{i_{n-1}}^{(m_{n-1})}}\bigg)^{\sigma}\Bigg)
\Bigg]_{\langle \mc I\rangle_\partial}
\end{aligned}   
\end{equation}
As in the proof of Lemma \ref{lem:jetP}, 
one can check that $J_\infty Q$ is well defined by the above formula, that it is an element of $C^n(J_\infty \mc A)$, and that $J_\infty Q\circ \iota=\iota\circ Q$. The uniqueness of such
$n$-fold $\lambda$-bracket is clear.
\end{proof}
\begin{remark}
We identify $\mc A\simeq\iota(\mc A)\subset J_{\infty}\mc A$. Then, since the diagram \eqref{diag:jetdP} is commutative,
the $n$-fold $\lambda$-bracket provided by Lemma \ref{lem:jetdP} satisfies  
\begin{equation} \label{Eq:dPAdPVA-ext1}
  \ldb a_1 {}_{\lambda_1}\dots{}_{\lambda_{n-1}} a_n\rdb
  =Q(a_1,\dots,a_n),  
\end{equation}
when evaluated on $a_1,\ldots, a_n \in\mc A$. 
\end{remark}
\begin{remark}
Any finitely generated differential algebra $\mc V$ can be realized as a quotient 
$\mc V=\mc R_\ell/\mc I_{\mc V}$ for some $\ell \geq 1$ and some two-sided differential ideal $\mc I_{\mc V}\subset \mc R_\ell$. 
Based on Remark \ref{20250722:rem2} any $n$-fold $\lambda$-bracket on $\mc V$ is also determined by the value that it takes on the generators of $\mc V$ given by the classes $[u_j]_{\mc I_{\mc V}} \in \mc V$ of the generators of $\mc R_\ell^{(0)}$. 
Hence,  if $c_1,\ldots,c_n\in \mc V$ admit representatives $p_1,\ldots,p_n \in \mc R_\ell$,
using equation \eqref{eq:QS}, we have 
\begin{align}
&\ldb c_1{}_{\lambda_1}\dots c_{n-1}
{}_{\lambda_{n-1}}c_n\rdb
\nonumber 
\\
&=\Bigg[
\sum_{\substack{j_1,\dots,j_n=1,\dots,\ell\\m_1,\dots,m_n\in\mb Z_{\geq0}}}
\frac{\partial p_n}{\partial u_{j_n}^{(m_n)}}
(\lambda_1+\dots+\lambda_{n-1}+\partial)^{m_n}
\nonumber\\
&
\bullet\Bigg(\Bigg(\dots\Bigg(
f_{j_1,\dots,j_n}(\lambda_1+x_1,\dots,\lambda_{n-1}+x_{n-1})\bullet_{(1,2)} (-\lambda_1-\partial)^{m_1}\bigg(\frac{\partial p_1}{\partial u_{i_1}^{(m_1)}}\bigg)^{\sigma}\Bigg)
\nonumber
\\
&\dots 
\Bigg)\bullet_{(n-1,n)}(-\lambda_{n-1}-\partial)^{m_{n-1}}\bigg(\frac{\partial p_{n-1}}{\partial u_{i_{n-1}}^{(m_{n-1})}}\bigg)^{\sigma}\Bigg)
\Bigg]_{\mc I_{\mc V}}
. \label{Eq:dPVAnMaster}
\end{align}
where $f_{j_1,\dots,j_n}(\lambda_1,\ldots,\lambda_{n-1})\in \mc R_\ell^{\otimes n}[\lambda_1,\ldots,\lambda_{n-1}]$ is a representative of 
$$\ldb [u_{j_1}]_{\mc I_{\mc V}}{}_{\lambda_1}\dots_{\lambda_{n-1}}[u_{j_n}]_{\mc I_{\mc V}} \rdb \in \mc V^{\otimes n}[\lambda_1,\ldots,\lambda_{n-1}], \qquad j_1,\dots,j_n=1,\dots,\ell\,.$$
\end{remark}
Let us now assume that $(\mc A,\ldb-,-\rdb)$ is a finitely generated dPA.  
Thanks to Lemma \ref{lem:jetdP}, we get a morphism $J_\infty: \wBRA(\cA)\to C(J_\infty \mc A)$ after completing it in degree $0$ by 
\begin{equation} \label{Eq:dPAdPVA-ext0}
  \wBRA(\cA)_0=\mc A_\sharp \to C^0(J_\infty \mc A)= (J_\infty \mc A)_{\sharp}
  , \quad \bar a\mapsto \tint \iota(a)\,.  
\end{equation}
Let us denote by $\dd_{\mc A}$ the differential on $\wBRA(\mc A)$ given by \eqref{Eq:dP-gen-0} and \eqref{Eq:dP-gen}.  
As already mentioned, it is proved in \cite{BFM} that $J_\infty \ldb-,-\rdb$ defines a dPVA structure on $J_\infty \mc A$ (this follows also from Lemma \ref{lem:jetdP}). 
Hence we can consider the variational dPVA differential on $C(J_{\infty}\mc A)$
given by \eqref{eq:dP0} and \eqref{eq:dP-1}, which we denote as $\dd_{J_\infty \mc A}$.   
To ease notation, we denote $J_\infty \ldb-,-\rdb$ as $\ldb-_\lambda -\rdb$, which satisfies due to \eqref{Eq:dPAdPVA-ext1}
\begin{equation} \label{Eq:dPAdPVA-br}
    \ldb a_\lambda b\rdb = \ldb a,b\rdb \in \mc A^{\otimes 2}\hookrightarrow (J_\infty \mc A)^{\otimes 2}, \quad \forall  a,b\in \mc A \hookrightarrow J_\infty \mc A\,.
\end{equation}
\begin{proposition}\label{Pr:dPAdPVA1}
For every $Q\in \wBRA(\mc A)_n$, $n\in\mb Z_{\geq0}$, we have
\begin{equation}\label{dPAdPVA-comm}
J_\infty(\dd_{\mc A}(Q))= \dd_{J_\infty \mc A}(J_\infty(Q))
\,.
\end{equation}
In particular, we get a morphism of complexes 
$$J_\infty: (\wBRA(\mc A) , \dd_{\mc A}) \to 
(C(J_\infty \mc A),\dd_{J_\infty \mc A})\,.$$
\end{proposition}
\begin{proof}
    In view of the master formula \eqref{Eq:dPVAnMaster} for $n$-fold $\lambda$-brackets on
    $J_{\infty}\mc A$, it suffices to check the equality \eqref{dPAdPVA-comm} on the generators $u_j$, $1\leq j \leq \ell$, viewed as elements of $J_\infty A$. 
(We use the same notation for these elements in $\mc A$, which we previously denoted as $[u_j]_{\mc I}\in \mc A$ and $[u_j]_{\langle \mc I\rangle_\partial}\in J_\infty \mc A$.) 

For $n=0$, $Q=\bar a\in \mc A_\sharp$ and we find 
\begin{align*}
 J_\infty(\dd_{\mc A}(\bar a))(u_j)=-\mult \ldb a,u_j\rdb,
\end{align*}  
where we used \eqref{Eq:dP-gen-0} with \eqref{Eq:dPAdPVA-ext1}. On the other hand, we have
\begin{align*} 
\dd_{J_\infty \mc A}(J_\infty(\bar a)) (u_j)
&=-\mult \ldb a_{\lambda} u_j\rdb|_{\lambda=0}= -\mult \ldb a,u_j\rdb, 
\end{align*}  
where we used \eqref{eq:dP0} with \eqref{Eq:dPAdPVA-ext0}. This proves \eqref{dPAdPVA-comm} for $n=0$.

Next, we look at $n\geq 1$. Fix $j_1,\ldots,j_{n+1}\in \{1,\ldots,\ell\}$. From Lemma \ref{lem:jetdP}
(cf. \eqref{Eq:dPAdPVA-ext1}) we have that
$$
J_\infty(\dd_{\mc A}(Q))(u_{j_1} ,\ldots, u_{j_{n+1}} )  
= \dd_{\mc A}(Q)(u_{j_1},\ldots,u_{j_{n+1}})
\,.
$$
On the other hand, we use \eqref{eq:dP} to obtain 
\begin{align*}
&\dd_{J_\infty \mc A}(J_\infty(Q))(u_{j_1} ,\ldots, u_{j_{n+1}} )  \\
&
=\sum_{s=0}^n (-1)^{ns}
|_{\lambda_{n+1}=\lambda_{n+1}^\dagger}
\sigma^s \circ  
\left(J_{\infty}(Q)_{\lambda_{\sigma^s(1)},\dots,\lambda_{\sigma^s(n-1)}}\otimes \Id_{\mc V}\right)
\\
&
\quad\quad\quad\quad\quad\quad\quad\quad\quad\quad\quad\quad
\circ \left(\Id_{\mc V}^{\otimes(n-1)}\otimes \ldb-_{\lambda_{\sigma^s(n)}}-\rdb\right)
\circ \sigma^{-s}
(u_{j_1} ,\ldots, u_{j_{n+1}} )\\
&+
%\frac{n}{n+1}
\sum_{s=0}^{n} (-1)^{n(s+1)}
|_{\lambda_{n+1}=\lambda_{n+1}^\dagger}
\sigma^s \circ  
\left(\ldb-_{\lambda_{\sigma^s(1)}}-\rdb\otimes \Id_{\mc V}^{\otimes(n-1)}\right)
\\
&
\quad\quad\quad\quad\quad\quad\quad\quad\quad\quad\quad\quad
\circ \left(\Id_{\mc V}\otimes J_{\infty} (Q)_{\lambda_{\sigma^s(2)},\dots,\lambda_{\sigma^s(n)}}\right)
\circ \sigma^{-s}
(u_{j_1} ,\ldots, u_{j_{n+1}} )
\\
&=
\sum_{s=0}^n (-1)^{ns}
\sigma^s \circ  
\left(Q\otimes \Id_{\mc V}\right)
\circ \left(\Id_{\mc V}^{\otimes(n-1)}\otimes \ldb-,-\rdb\right)
\circ \sigma^{-s}
(u_{j_1} ,\ldots, u_{j_{n+1}} )\\
&+
%\frac{n}{n+1}
\sum_{s=0}^{n} (-1)^{n(s+1)}
\sigma^s \circ  
\left(\ldb-,-\rdb\otimes \Id_{\mc V}^{\otimes(n-1)}\right)
\circ \left(\Id_{\mc V}\otimes Q\right)
\circ \sigma^{-s}
(u_{j_1} ,\ldots, u_{j_{n+1}} )
\\
&= \dd_{\mc A}(Q)(u_{j_1},\ldots,u_{j_{n+1}})
\,,
\end{align*}
where in the second equality we used \eqref{Eq:dPAdPVA-ext1} to replace $J_{\infty }Q$ by $Q$ and $\ldb-_{\lambda}-\rdb$ by $\ldb-,-\rdb$ and in the third equality we used equation \eqref{Eq:dP-gen}. Hence, we get that the RHS and LHS of \eqref{dPAdPVA-comm} coincide when evaluated at $u_{j_1},\dots,u_{j_{n+1}}$ thus concluding the proof. 
\end{proof}
In conclusion, we proved the following result.
\begin{corollary}\label{Cor:dPAdPVA1}
For every $n\in\mb Z_{\geq0}$ the linear map $J_\infty: \wBRA(\mc A)_n \to C^n(J_\infty \mc A)$ given
by Lemma \ref{lem:jetdP} (or \eqref{Eq:dPAdPVA-ext0} for $n=0$) gives a morphism in cohomology 
$$
\widehat{\dPH}^n(\cA)\rightarrow \dPVH^n(J_\infty \mc A)
\,.
$$ 
\end{corollary}

\section{Constructions with the quotient functor}

We consider the quotient functor %$\QQ$
from the category of differential algebras to the one of algebras as in Example \ref{exa:dPVAtodPA}.
Recall that it is given by the projection map $\pi:\mc V\mapsto q(\mc V)=\mc V/\langle \partial \mc V \rangle$. Moreover, if $\mc V$ is a dPVA, then $q(\mc V)$ is a dPA with double Poisson bracket defined by \eqref{20250822:eq1}. The next result is a straightforward generalization of this construction to any $n$-fold $\lambda$-bracket. 
\begin{lemma}\label{lem:QuotdPV}
Let $\mc V$ be a differential algebra. For every 
$n$-fold $\lambda$-bracket $c:=\ldb-_{\lambda_1}-\dots-_{\lambda_{n-1}}-\rdb\in C^n(\mc V)$, $n\geq1$, there exists a unique $n$-bracket 
$q(c)\in\wBRA(q(\mc V))_n$ making the following diagram commute:
\begin{equation}\label{diag:QuotdPV}
% \text{the diagram has been commented to speed up compiling}
% \begin{comment}
\begin{tikzcd}
\mc V^{\otimes n}\arrow[r,"c"] \arrow[d,twoheadrightarrow,"\pi^{\otimes n}"] & \mc V^{\otimes n}[\lambda_1,\dots,\lambda_{n-1}] \arrow[d,twoheadrightarrow,"\pi^{\otimes n}"]
\\
q(\mc V)^{\otimes n} \arrow[r,"q(c)"] & q(\mc V)^{\otimes n}
\end{tikzcd}
% \end{comment}
\end{equation}
\end{lemma}
\begin{proof}
For $a_1,\dots,a_n\in\mc V$ we define 
\begin{equation} \label{Eq:dPAdPVA-qc}
 q(c)(\pi(a_1),\ldots,\pi(a_n))
 = \pi^{\otimes n}( \ldb a_1{}_{\lambda_1}\dots a_{n-1}{}_{\lambda_{n-1}} a_n\rdb 
 |_{\lambda_1=\ldots=\lambda_{n-1}=0} )\, .
\end{equation}
One can then check that this gives a well-defined element of $\wBRA(\mc V)_n$, which is the unique one making the diagram \eqref{diag:QuotdPV} commute.
\end{proof}
For $n=0$, we define the map $C^0(\mc V)=\mc V_{\sharp} \to \wBRA(\mc V)_0=q(\mc V)_\sharp$ by mapping $\tint a \mapsto q(\tint a)=\overline{\pi(a)}$, which is obtained by factoring 
the map $\mc V\to q(\mc V)_{\sharp}=q(\mc V)/[q(\mc V),q(\mc V)]$, $a\mapsto \overline{\pi(a)}$, through the linear map $\tint : \mc V \mapsto \mc V_\sharp=\mc V/(\partial \mc V+[\mc V,\mc V])$. 
From now on let us assume $(\mc V,\ldb-_\lambda-\rdb)$ is a dPVA inducing the dPA structure 
$(q(\mc V),\ldb-,-\rdb)$ using \eqref{20250822:eq1}.
We denote by $\dd_{\mc V}$ the differential on $C(\mc V)$ given by \eqref{eq:dP0} and \eqref{eq:dP-1}, and we denote by $\dd_{q(\mc V)}$ the differential on $\wBRA(q(\mc V))$ given by \eqref{Eq:dP-gen-0} and \eqref{Eq:dP-gen}.
\begin{proposition}\label{Pr:dPAdPVA2}
For every $c\in C^n(\mc V)$, $n\in\mb Z_{\geq0}$, we have
\begin{equation}\label{dPAdPVA-comm2}
q(\dd_{\mc V}(c))= \dd_{q(\mc V)}(q(c)) \,.
\end{equation}
In particular, we get a morphism of complexes 
$$q: (C(\mc V), \dd_{\mc V}) \to (\wBRA(q(\mc V)) ,\dd_{q(\mc V)})\,.$$
\end{proposition}
\begin{proof}
Let $n=0$ and let $c=\tint a\in \mc V_\sharp$, evaluating both sides of \eqref{dPAdPVA-comm2} on $\pi(b)\in q(\mc V)$ gives 
\begin{align*}
q(\dd_\mc V(\tint a))(\pi(b))&= - \pi\circ\mult( \ldb a_\lambda b\rdb |_{\lambda =0}) \,, \\
\dd_{q(\mc V)}(q(\tint a))(\pi(b))&=-\mult  \ldb \pi(a),\pi(b)\rdb
=-\mult\circ(\pi\otimes\pi)(\ldb a_\lambda b\rdb|_{\lambda=0})\,,
\end{align*}
where we used \eqref{Eq:dPAdPVA-qc} with \eqref{eq:dP0} in the first equality, and
\eqref{Eq:dP-gen-0} and \eqref{20250822:eq1} in the second equality. 
Both expressions are equal since the following diagram is commutative:
%\pecetta{this diagram has been commented to speed up compiling}
%\begin{comment}
\begin{equation*}
\begin{tikzcd}
\mc V\otimes\mc V\arrow[r,"\mult"] \arrow[d,,"\pi\otimes\pi"]
& \mc V\arrow[d,"\pi"]
\\
q(\mc V)\otimes q(\mc V)\arrow[r,"\mult"] & q(\mc V)
\end{tikzcd}
\end{equation*}
%\end{comment}
%
This proves that $q(\dd_{\mc V}(\tint a))=\dd_{q(\mc V)}(q(\tint a))$. 

For $n\geq 1$, given $c=c_{\lambda_1,\dots,\lambda_{n-1}}\in C^n(\mc V)$ and 
$a_1,\ldots,a_{n+1}\in \mc V$, using \eqref{dPAdPVA-comm2} we find
$$
q(\dd_{\mc V}(c))(\pi(a_1),\ldots,\pi(a_{n+1}))
=
\pi^{\otimes (n+1)}d_{\mc V}(c)(a_1,\dots,a_{n+1})|_{\lambda_1=\dots=\lambda_n=0}
\,.
$$
On the other hand, using \eqref{Eq:dP-gen} we have
\begin{align*}
&\dd_{q(V\mc )}(q(c))(\pi(a_1),\ldots,\pi(a_{n+1})) \\
&=
\sum_{s=0}^{n} (-1)^{ns}\sigma^s \circ  
(q(c)\otimes \Id_{\mc V})\circ (\Id_{\mc V}^{\otimes(n-1)}\otimes \ldb-,-\rdb)
\circ \sigma^{-s} (\pi(a_1),\ldots,\pi(a_{n+1})) \\
&+\sum_{s=0}^n (-1)^{n(s+1)}\sigma^s \circ  
(\ldb-,-\rdb\otimes \Id_{\mc V}^{\otimes(n-1)}) \circ (\Id_{\mc V}\otimes q(c))
\circ \sigma^{-s}
(\pi(a_1),\ldots,\pi(a_{n+1}))
\\
&
=\sum_{s=0}^n (-1)^{ns}
\pi^{\otimes (n+1)}\circ |_{\lambda_{n+1}=\lambda_{n+1}^\dagger}
\sigma^s \circ  
\left(c_{\lambda_{\sigma^s(1)},\dots,\lambda_{\sigma^s(n-1)}}\otimes \Id_{\mc V}\right)
\\
&
\quad\quad\quad\quad\quad\quad\quad\quad\quad
\circ \left(\Id_{\mc V}^{\otimes(n-1)}\otimes \ldb-_{\lambda_{\sigma^s(n)}}-\rdb\right)
\circ \sigma^{-s}
(a_{1} ,\ldots, a_{n+1} )|_{\lambda_1=\dots=\lambda_n=0}\\
&+
%\frac{n}{n+1}
\sum_{s=0}^{n} (-1)^{n(s+1)}
\pi^{\otimes (n+1)}\circ |_{\lambda_{n+1}=\lambda_{n+1}^\dagger}\sigma^s \circ  
\left(\ldb-_{\lambda_{\sigma^s(1)}}-\rdb\otimes \Id_{\mc V}^{\otimes(n-1)}\right)
\\
&
\quad\quad\quad\quad\quad\quad\quad\quad\quad
\circ \left(\Id_{\mc V}\otimes c_{\lambda_{\sigma^s(2)},\dots,\lambda_{\sigma^s(n)}}\right)
\circ \sigma^{-s}
(a_{1} ,\ldots, a_{n+1} )|_{\lambda_1=\dots=\lambda_n=0}
\\
&
=\pi^{\otimes (n+1)} \dd_{\mc V}(c)(a_1,\dots,a_{n+1})|_{\lambda_1=\dots=\lambda_n=0}
\,,
\end{align*}
where in the second equality we used the fact that $\sigma$ and $\pi^{\otimes n}$ commute,
we used equation \eqref{Eq:dPAdPVA-qc} and we used the fact that $\pi(\lambda_{n+1}^\dagger X)|_{\lambda_1=\dots=\lambda_n=0}=0$ for
every $X\in\mc V^{\otimes n}$.
This proves \eqref{dPAdPVA-comm2} for $n\geq1$. 
\end{proof}
As a consequence, we get the following result.
\begin{corollary}\label{Cor:dPAdPVA2}
For every $n\in\mb Z_{\geq0}$ the linear map $q: C^n(\mc V) \to \wBRA(q(\mc V))_n$ given
by Lemma \ref{lem:QuotdPV} (or $\tint a \mapsto \overline{\pi(a)}$ for $n=0$) gives a morphism in cohomology 
$$
\dPVH^n(\mc V) \to \widehat{\dPH}^n(q(\mc V)) \,.
$$ 
\end{corollary}
In Example \ref{exa:jetdPVA}, we have seen that composing the quotient functor after the jet functor gives back the original dPA, i.e. $q(J_\infty(\mc A))=\mc A$ (up to isomorphism).  
Similarly, under the assumption that $\mc A$ is finitely generated, the composition 
$$
\wBRA(\mc A) \stackrel{J_\infty}{\longrightarrow} C(J_\infty \mc A) \stackrel{q}{\longrightarrow} 
\wBRA(q(J_\infty(\mc A))) = \wBRA(\mc A)
$$
obtained from Lemmas \ref{lem:jetdP} and \ref{lem:QuotdPV} is easily seen to be an isomorphism. 
Combining Corollaries \ref{Cor:dPAdPVA1} and \ref{Cor:dPAdPVA2}, we get the following result in cohomology. 
\begin{theorem}\label{20250822:thm2}
 The composition 
 $$
\widehat{\dPH}(\mc A) \longrightarrow \dPVH(J_\infty \mc A) \longrightarrow \widehat{\dPH}(\mc A)
$$   
is an isomorphism. In particular, the linear map $\widehat{\dPH}(\mc A) \to \dPVH(J_\infty \mc A)$ is injective,  
and $\dPVH(J_\infty \mc A) \to \widehat{\dPH}(\mc A)$ is surjective. 
\end{theorem}

\begin{example}
    Consider the case of a quiver $Q$ made of $\ell \in 2\Z_{>0}$ loops all based at a single vertex, so that $\cA:=\kk Q \simeq \mc R_\ell^{(0)}$. As in Section \ref{Sec:Quiv}, fix a skewsymmetric matrix $C\in \Gl_{\ell}(\kk)$. Then, consider the corresponding dPA structure given by \eqref{Eq:Quiv-gen}, where all $e_s=1$ since $B=\kk$. 
    By Theorem \ref{Thm:DPcoh-Quiv}, $\widehat{\dPH}^n(\cA)\simeq \delta_{n,0} \kk$. 
Under the map $J_\infty$, the double Poisson bracket on $\cA$ gives rise to a $2$-fold $\lambda$-bracket on $J_\infty \cA \simeq \mc R_\ell$ that is given by \eqref{Eq:uu-odd-gen} with $K=C$ and $M=0$, cf. \eqref{Eq:dPAdPVA-ext1}. 
Using this last dPVA structure on $J_\infty \cA$, 
we get in view of Corollary \ref{M=0} the corresponding cohomology $\coH^n(\Omega(J_\infty \cA),\delta)\simeq \delta_{n,0} \kk$. 
Since we have an isomorphism $\dPVH(J_\infty \cA)\simeq \coH(\Omega(J_\infty \cA),\delta)$ by the identifications of Section \ref{sec:rel}, 
we conclude that the two maps in Theorem \ref{20250822:thm2} are isomorphisms. 
\end{example}

\begin{example}
    In general, the two maps in Theorem \ref{20250822:thm2} can fail to be isomorphisms. 
To see this, take $\cA=\mc R_\ell^{(0)}$ ($\ell \geq 1)$ with the zero double Poisson bracket, so that 
$J_\infty \cA=\mc R_\ell$ has the zero dPVA structure. 
Then, 
\[
\widehat{\dPH}(\cA)=\wBRA_\kk(\mc R_\ell^{(0)}) , \qquad 
\dPVH(J_\infty \cA)=C(\mc R_\ell)\,.
\]
Since the image of any $2$-fold $\lambda$-bracket of the form \eqref{Eq:uu-odd-gen} with $M>0$ via the map $C^2(\mc R_\ell)\rightarrow\wBRA(\mc R_\ell^{(0)})_2$ is the trivial double Poisson bracket on $\mc R_\ell^{(0)}$, we conclude that
the map $C(\mc R_\ell)\rightarrow\wBRA(\mc R_\ell^{(0)})$ is not injective (similarly, one checks that the map $\wBRA(\mc R_{\ell}^{(0)})\rightarrow C(\mc R_\ell)$ is not surjective).
\end{example}

%%%
\section{Going to representation algebras}

In this section we show the compatibility with the results of Chapter \ref{Ch:PA-PVA}
and the induced cohomologies on representation spaces defined in Chapters \ref{CH:rep-dPA} and \ref{Ch:repVardPVA}.

Let $\mc A$ be an associative algebra and let $\mc A_N$, $N\geq1$, be the commutative algebra defined in Section \ref{ss:Rep-Not}. Recall from Section \ref{sec:7.1} that we have a map
$$
\tr:\wBRA(\mc A)_n\to \mf X^n(\mc A_N)\,,\quad
n\geq0\,,
$$
explicitly given by \eqref{Eq:TrBrRep}. Recall also from \cite{BFM} that
$(J_{\infty}\mc A)_N\simeq J_{\infty}\mc A_N$. Hence, from Section \ref{sec:13.3} we get a map
$$
\tr:C^n(J_{\infty}\mc A)\rightarrow C^n(J_{\infty}\mc A_N)
\,,
\quad
n\geq0\,,
$$
explicitly given by \eqref{eq:trace0} for $n=0$ and by \eqref{Eq:TrBrRep-PVA} for $n\geq1$.
\begin{proposition}\label{20250822:prop1}
Let $\mc A$ be a finitely generated associative algebra. For every $n\geq0$ the following diagram is commutative
%\pecetta{this diagram has been commented to speed up compiling}
%\begin{comment}
\begin{equation*}
\begin{tikzcd}
\wBRA(\mc A)_n \arrow[r,"J_{\infty}"] \arrow[d,"\tr"]& C^n(J_{\infty}\mc A)\arrow[d,"\tr"]
\\
\mf X^n(\mc A_N) \arrow[r,"J_{\infty}"] & C^n(J_{\infty}\mc A_N)
\end{tikzcd}
\end{equation*}
%\end{comment}
where $J_{\infty}:\wBRA(\mc A)_n\rightarrow C^n(J_{\infty}\mc A)$ is the map given by Lemma \ref{lem:jetdP} and $J_{\infty}:\mf X^n(\mc A_N)\rightarrow C^n(J_{\infty}\mc A_N)$ is the map given by Lemma \ref{lem:jetP}.
\end{proposition}
\begin{proof}
The case $n=0$ follows immediately from the definitions (see Lemma \ref{20240829:lem1}(c), \eqref{Eq:Tr-morph}, \eqref{Eq:dPAdPVA-ext0} and \eqref{Eq:PAPVA-ext0}).
Let $Q=\ldb-\rdb\in\wBRA(\mc A)_n$, $n\geq1$, and let $a^1,\dots,a^n\in\mc A$. In view of the Master Formula \eqref{Eq:PVAnMaster} for $n$-$\lambda$-brackets, to prove the claim it suffices to show that
$$
\tr(J_{\infty}Q)_{\lambda_1,\dots,\lambda_n}(a^{1}_{i_1j_1},\dots,a^n_{i_nj_n})
=J_{\infty}(\tr Q)_{\lambda_1,\dots,\lambda_n}(a^{1}_{i_1j_1},\dots,a^n_{i_nj_n})
\,,
$$
for every $i_1,j_1,\dots,i_n,j_n=1,\dots N$.
We have
\begin{equation}\label{20250826:eq1}
\begin{split}
&\tr(J_{\infty}Q)_{\lambda_1,\dots,\lambda_n}(a^{1}_{i_1j_1},\dots,a^n_{i_nj_n})
\\
&=
\sum_{\tau\in S_{n-1}}\sgn(\tau)\left(
(J_{\infty }Q)_{\lambda_{\tau(1)},\dots,\lambda_{\tau(n-1)}}(a^{\tau(1)},\dots,a^{\tau(n-1)},a^n)
\right)_{\tau(\underline i,\underline j)}
\\
&=
\sum_{\tau\in S_{n-1}}\sgn(\tau)
\ldb a^{\tau(1)},\dots,a^{\tau(n-1)},a^n\rdb
_{\tau(\underline i,\underline j)}
\,,
\end{split}
\end{equation}
where in the first equality we used the definition of the trace map \eqref{Eq:TrBrRep-PVA}, and in the second equality we used \eqref{Eq:dPAdPVA-ext1}.
On the other hand, using \eqref{Eq:PAPVA-ext1} and \eqref{Eq:TrBrRep} we have
\begin{equation}\label{20250826:eq2}
\begin{split}
&J_\infty(\tr Q)_{\lambda_1,\dots,\lambda_n}(a^1_{i_1j_1},\dots,a^n_{i_nj_n})
=\tr Q(a^1_{i_1j_1},\dots,a^n_{i_nj_n})
\\
&=\sum_{\tilde\sigma\in S_{n}^{(1)}}\sgn(\tilde\sigma)\ldb a^1,a^{\tilde\sigma(2)},\dots a^{\tilde\sigma(n)}\rdb_{\tilde\sigma (\underline i,\underline j)}
\,.
\end{split}
\end{equation}
The RHS of equations \eqref{20250826:eq1} and \eqref{20250826:eq2} coincide since we have the bijection
\begin{equation}\label{20250826:bij}
S_n^{(1)}\rightarrow S_{n-1}\,,
\quad
\tilde\sigma\mapsto \sigma^{-1}\tilde\sigma\sigma
\,,
\end{equation}
where $\sigma=(1 2\dots n)\in S_n$.
\end{proof}
Let us assume that $\mc A$ is a dPA. Then, we have the following corollary of Proposition \ref{20250822:prop1}, of Theorems \ref{Thm:dP-rep2} and \ref{Thm:dP-rep2-PVA}, and of Propositions
\ref{Pr:PAPVA1} and \ref{Pr:dPAdPVA1}.
\begin{corollary}\label{20250822:cor1}
We have the following commutative diagram of complexes
%\pecetta{this diagram has been commented to speed up compiling}
%\begin{comment}
\begin{equation*}
\begin{tikzcd}
(\wBRA(\mc A),\dd_{\mc A}) \arrow[r,"J_{\infty}"] \arrow[d,"\tr"]& (C(J_{\infty}\mc A),\dd_{J_{\infty}\mc A})\arrow[d,"\tr"]
\\
(\mf X(\mc A_N),(-1)^\bullet\delta_{\mc A_N}) \arrow[r,"J_{\infty}"] & (C(J_{\infty}\mc A_N),(-1)^\bullet\delta_{J_{\infty\mc A_N}})
\end{tikzcd}
\end{equation*}
%\end{comment}
\end{corollary}
Let now $\mc V$ be a differential algebra, and let $\mc V_N$, $N\geq1$, be the commutative differential algebra defined in Section \ref{sec:dPVAtoPVA-1}. The map $\pi(a)_{ij}\mapsto \pi(a_{ij})$, $i,j=1,\dots,N$, $a\in\mc V$, gives an isomorphism $q(\mc V)_N\simeq q(\mc V_N)$,
and, as before, we have maps
$$
\tr:C^n(\mc V)\rightarrow C^n(\mc V_N)\,,
\quad
\tr:\wBRA(q(\mc V))_n\rightarrow\mf X^n(q(\mc V_N))\,,
\quad
n\geq0\,.
$$
\begin{proposition}\label{20250822:prop2}
Let $\mc V$ be a differential algebra. For every $n\geq0$ the following diagram is commutative
% \begin{comment}
\begin{equation*}
\begin{tikzcd}
C^n(\mc V)\arrow[r,"q"] \arrow[d,"\tr"]& \wBRA(q(\mc V))_n \arrow[d,"\tr"]
\\
C^n(\mc V_N) \arrow[r,"q"] & \mf X^n(q(\mc V_N))
\end{tikzcd}
\end{equation*}
% \end{comment}
where $q:C^n(\mc V)\rightarrow \wBRA(q(\mc V))_n$ is the map given by Lemma \ref{lem:QuotdPV} and $q:C^n(\mc V_N)\rightarrow \mf X^n(q(\mc V_N))$ is the map given by Lemma \ref{lem:QuotPV}.
\end{proposition}
\begin{proof}
The case $n=0$ follows immediately from the definitions. Let $c=\ldb-_{\lambda_1}-\dots-_{\lambda_{n-1}}-\rdb\in C^n(\mc V)$, $n\geq1$, and let $a^1,\dots, a^n\in\mc V$. To prove the claim it suffices to show that
$$
\tr(q(c))(\pi(a^1_{i_1,j_1}),\dots ,\pi(a^n_{i_nj_n}))
=q(\tr(c))(\pi(a^1_{i_1,j_1}),\dots ,\pi(a^n_{i_nj_n}))
\,,
$$
for every $i_1,j_1,\dots,i_n,j_n=1,\dots,N$.
We have
\begin{equation}\label{20250826:eq3}
\begin{split}
&\tr(q(c))(\pi(a^{1}_{i_1j_1}),\dots,\pi(a^n_{i_nj_n}))
\\
&=\sum_{\tilde\sigma\in S_n^{(1)}}\sgn(\tilde\sigma)\left(
q(c)(\pi(a^1),\pi(a^{\tilde\sigma(2)}),\dots,\pi(a^{\tilde\sigma(n)}))
\right)_{\tilde\sigma(\underline i,\underline j)}
\\
&=
\sum_{\tilde\sigma\in S_{n}^{(1)}}\sgn(\tilde\sigma)
\pi^{\otimes n}
\left(
\ldb a^{1}_{\lambda_1}\dots a^{\tilde\sigma(n-1)}_{\lambda_{\tilde\sigma(n-1)}}a^{\tilde\sigma(n)}\rdb
_{\tilde\sigma(\underline i,\underline j)}
\right)|_{\lambda_1=\lambda_2=\dots=\lambda_{n}=0}\,,
\end{split}
\end{equation}
where in the first equality we used the definition of the trace map \eqref{Eq:TrBrRep}, and in the second equality we used \eqref{Eq:dPAdPVA-qc}.
On the other hand, using \eqref{Eq:dPAdPVA-qc} and \eqref{Eq:TrBrRep-PVA} we have
\begin{equation}\label{20250826:eq4}
\begin{split}
&q(\tr c)(\pi(a^1_{i_1j_1}),\dots,\pi(a^n_{i_nj_n}))
\\
&=\pi^{\otimes n}\left(
(\tr c)_{\lambda_1,\dots,\lambda_{n-1}}(a^1_{i_1j_1},\dots,a^n_{i_nj_n})
\right)|_{\lambda_1=\lambda_2=\dots\lambda_n=0}
\\
&=\sum_{\tau\in S_{n-1}}\sgn(\tau)\pi^{\otimes n}
\left(
\ldb a^{\tau(1)}_{\lambda_{\tau(1)}}\dots a^{\tau(n-1)}_{\lambda_{\tau(n-1)}} a^n\rdb_{\tau(\underline i,\underline j)}
\right)|_{\lambda_1=\lambda_2=\dots\lambda_n=0}
\,.
\end{split}
\end{equation}
The RHS of equations \eqref{20250826:eq3} and \eqref{20250826:eq4} coincide in view of the bijection~\eqref{20250826:bij}.
\end{proof}
Let us assume that $\mc V$ is a dPVA. Then, we have the following corollary of Proposition \ref{20250822:prop2}, of Theorems \ref{Thm:dP-rep2} and \ref{Thm:dP-rep2-PVA}, and of Propositions
\ref{Pr:PAPVA2} and \ref{Pr:dPAdPVA2}.
\begin{corollary}\label{20250822:cor2}
We have the following commutative diagram of complexes
%\pecetta{the diagram has been commented to speed up compiling}
%\begin{comment}
\begin{equation*}
\begin{tikzcd}
(C(\mc V),\dd_{\mc V})\arrow[r,"q"] \arrow[d,"\tr"]& (\wBRA(q(\mc V)),\dd_{q(\mc V)}) \arrow[d,"\tr"]
\\
(C(\mc V_N),\delta_{\mc V_N}) \arrow[r,"q"] & (\mf X(q(\mc V_N)),\delta_{q(\mc V_N)})
\end{tikzcd}
\end{equation*}
%\end{comment}
\end{corollary}
In conclusion, combining Corollaries \ref{20250822:cor1} and \ref{20250822:cor2}, and Theorems
\ref{20250822:thm1} and \ref{20250822:thm2} we get the following result.
\begin{theorem} \label{Thm:RepJetQuot}
Let $\mc A$ be a finitely generated algebra. For every $N\geq1$ we have the following commutative diagram
%\pecetta{the diagram has been commented to speed up compiling}
%\begin{comment}
\begin{equation*}
\begin{tikzcd}
\widehat{\dPH}(\mc A)\arrow[r]\arrow[d]
&
\dPVH(J_\infty \mc A)\arrow[r]\arrow[d]
&
\widehat{\dPH}(\mc A)\arrow[d]
\\
\coH_{CE}(\mc A_N) \arrow[r]
&
\PVH(J_\infty \mc A_N) \arrow[r]
&
\coH_{CE}(\mc A_N)
\end{tikzcd}
\end{equation*}
%\end{comment}
where the composition of horizontal arrows is an isomorphism.
\end{theorem}

%%%%%%%%%%%%%%%%%%%%%%%%%%%%%%%%%%%
%%%%%%%%%%%%%%%%%%%%%%%%%%%%%%%%%%%
%%%%%%%%%%% INDEX   %%%%%%%%%%%%%%%
%%%%%%%%%%%%%%%%%%%%%%%%%%%%%%%%%%%
%%%%%%%%%%%%%%%%%%%%%%%%%%%%%%%%%%%

\small 
\printglossary[title={Index of Notation}, toctitle={Index of Notation}]
\normalsize

%%%%%%%%%%%%%%%%%%%%%%%%%%%%%%%%%%%
%%%%%%%%%%%%%%%%%%%%%%%%%%%%%%%%%%%
%%%%%%%%%%%%%%%%%%%%%%%%%%%%%%%%%%%
%%%%%%%%%%% BIBLIO  %%%%%%%%%%%%%%%
%%%%%%%%%%%%%%%%%%%%%%%%%%%%%%%%%%%
%%%%%%%%%%%%%%%%%%%%%%%%%%%%%%%%%%%
%%%%%%%%%%%%%%%%%%%%%%%%%%%%%%%%%%%

%\chapter{Ghost chapter -- end of the text}

\end{document}